**А. В. Ивако**

# ВВЕДЕНИЕ В ТЕОРИЮ МОЛЕКУЛЯРНЫХ ПРОСТРАНСТВ

A. V. Evako

# INTRODUCTION TO THE THEORY OF MOLECULAR SPACES

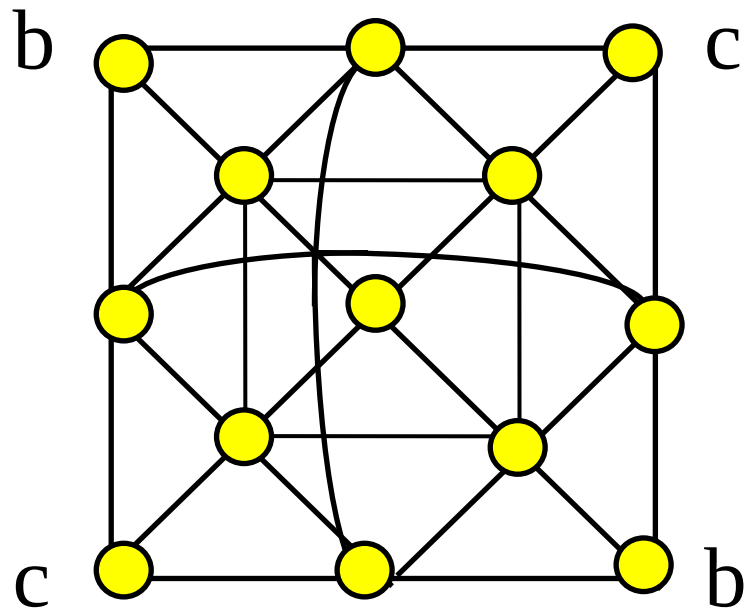

Книга базируется на предыдущей книге автора “Теория молекулярных пространств и ее приложения к компьютерам, физике и другим областям”, опубликованной в 1997 г. В ней подробно рассказывается о молекулярных пространствах, являющихся дискретными моделями многомерных непрерывных объектов, и предназначенных для использования в компьютерах. Книга может быть полезной математикам, программистам, создателям компьютерной графики и всем специалистам, работающим с многомерными компьютерными образами. Предлагаемые модели представляют интерес для физиков, а также могут применяться в биологии, химии и многих других областях. Книга рассчитана на широкие круги читателей, использующих компьютеры в своей деятельности.

The book is based on the previous book of the author “The Theory of Molecular Spaces and Its Applications to Computers, Physics, And Other Fields”, published in 1997. This book describes molecular spaces which are digital models of continuous many-dimensional objects. The book can be useful to mathematicians, programmers, designers of computer graphics and all professionals working with many-dimensional computer images. These models are of interest to physicists and can be utilized in biology, chemistry and many other fields. The book is intended to a wide circle of readers using computers in their work.

ВВЕДЕНИЕ В ТЕОРИЮ МОЛЕКУЛЯРНЫХ ПРОСТРАНСТВ.

Теория молекулярных пространств (ТМП) является одним из направлений дигитальной топологии. Она изучает образование, видоизменение и свойства 2-х, 3-х, и многомерных образов при работе компьютера и в его памяти. Необходимость создания математической теории многомерных пространств, построенных из конечного числа элементов, возникла сравнительно недавно в связи с использованием многомерных образов различного назначения в компьютерных программах.
Области применения ТМП весьма разнообразны.
В первую очередь, методы ТМП могут быть использованы в компьютерной графике, распознавании образов, моделировании объектов, описании динамически движущихся поверхностей, научных исследованиях и научной визуализации и других областях, связанных с компьютерами и их применением.
В математике подход ТМП найдет применение во многих сферах, включая топологию, геометрию и теорию дифференциальных уравнений.
Дигитальные многомерные объекты, исследуемые в ТМП, предлагают физикам “зоопарк” разнообразных дискретных моделей для описания физического вакуума и реального физического пространства.
Модели и методы ТМП могут успешно работать в таких областях, как биология, химия, медицина, промышленность, военное дело, а также во многих других.
Книга рассчитана на самый широкий круг читателей, начиная от математиков, программистов, физиков, химиков, биологов, специалистов в области медицины, промышленности, военного дела и общественных наук, использующих компьютеры в своей работе, и кончая школьниками старших классов.

Александр В. Ивако.

INTRODUCTION TO THE THEORY OF MOLECULAR SPACES.

The theory of molecular spaces is one of alternative branches of digital topology that studies constructing and modifying 2, 3 and n-dimensional digital image arrays in a computer and its memory. Demands for mathematical theory of multidimensional spaces built of finite number of points appeared only recently in connection with the use of many dimensional images in computer programs.
In the first place methods of TMS can be used in computer graphics, pattern analysis and recognition, object modeling, detection of dynamically moving surfaces, scientific visualization and other areas connecting to computers and its applications.
In mathematics methods of TMS can be used in topology, geometry, differential equations, etc.
Digital n-dimensional images offer to physicists a “zoo” of various discrete models of a physical vacuum and the real physical space.
Methods of TMS can work successfully in such areas as biology, chemistry, industry, medicine, etc. This book may concern a wide circle of readers, from programmers, mathematicians, physicists, biologists, medicine, industry, military and other specialists to high school students, in other words, anyone using computers in his practice.

Alexander V. Evako.

Сведения об авторе.

Александр Вадимович Ивако окончил ф-т теоретической и экспериментальной физики Московского инженерно-физического института и механико-математический факультет Московского государственного университета. Работал на кафедрах математики и физики в ряде московских вузов. В течение нескольких лет работал в университетах США. В настоящий момент его основные интересы лежат в области теории молекулярных пространств и ее приложений к компьютерам, математике, теоретическому и экспериментальному исследованию физического вакуума и реального физического пространства. Работы А. Ивако опубликованы в ряде журналов, включая “Discrete Mathematics”, “International Journal of Theoretical Physics”, “Journal of Mathematical Imaging and Vision”, а также в трудах ряда международных конференций. Автор монографии “Теория молекулярных пространств и ее приложения к компьютерам, физике и другим областям”, изданной в 1997 г.

Alexander V. Evako graduated from Moscow Engineering Physical Institute in Physics and Moscow University in Mathematics. After receiving Ph. D. he worked at Physics, Mathematics and Computer Science Departments for various Universities in Moscow. For several years he worked at Universities in the USA. His major interests include theory of molecular spaces and its applications to computers, mathematics, theoretical and experimental study of a physical vacuum and the real physical space. His works were published in a number of journals including “Discrete Mathematics”, “International Journal of Theoretical Physics”, “Journal of Mathematical Imaging and Vision” and proceedings of several international conferences. He is the author of the book “The Theory of Molecular Spaces and Its Applications to Computers, Physics, And Other Fields”, 1997.

# ОГЛАВЛЕНИЕ

**CONTENTS**

ОГЛАВЛЕНИЕ

# ПРЕДИСЛОВИЕ

Книга посвящена теории молекулярных пространств -- одному из направлений дигитальной топологии. Дигитальная топология является областью науки, возникшей сравнительно недавно на стыке таких дисциплин, как математика и компьютерные науки, и имеющей приложения во многих отраслях знаний. По одному из определений дигитальная топология (digital topology) есть наука о топологических свойствах дигитальных образов пространственных объектов, возникающих при работе компьютера (topological properties of digital image arrays). Компьютерные образы различных объектов в своей основе являются дигитальными, то есть построенными из одинаковых неделимых элементов, и конечными, состоящими всегда из конечного числа элементов, ограниченного объемом памяти машины. Хотя потребность в такой теории возникла в науке уже несколько десятилетий назад, только развитие компьютеров дало реальный толчок появлению дигитальной топологии.

Традиционно считается, что внесение математических методов в прикладные науки всегда приносит пользу последним и вызывает качественный скачок вперед в той или иной прикладной области. Следует ожидать, что применение методов дигитальной топологии в компьютерных науках, а также в чистой и прикладной математике явится плодотворным катализатором в их развитии.

Перечислим более подробно те области применения компьютеров, где, по мнению специалистов, использование дигитальной топологии и ТМП принесет реальную и немедленную пользу. Это -- компьютерная графика (computer graphics), компьютерное моделирование объектов (computer object modeling), диаграммный анализ (pattern analysis), научное представление объектов в виде зрительных образов -- научная визуализация (scientific visualization).

Сходство молекулярных пространств с химическими и биологическими объектами дает основание надеяться, что книга представит интерес и может быть плодотворно использована специалистами в этих областях науки.

Однако, все более вырисовывается перспектива, что, помимо научного и прикладного применения, дигитальная топология может претендовать на некоторое более значительное место на шкале научных теорий, чем только как чисто математическое направление. Дело в том, что методы этой науки дают реальную возможность применить дигитальный подход как ко всей математике, которая всегда базировалась на понятиях непрерывности и бесконечности, так и к окружающему нас материальному миру в целом, включая реальное физическое пространство, в котором мы существуем.

Если это произойдет, то будет затронут существенный слой нашего мировоззрения и восприятия окружающего нас мира.

При этом в матеиатике исчезают многочисленные особенности и парадоксы, связанные с несоизмеримостью, бесконечно малыми и бесконечно большими величинами, отпадает надобность в иррациональных числах и тому подобное.

В представлениях об окружающем нас мире дискретный взгляд на природу и физическое пространство более подходит нашему опыту, чем непрерывный. Из существования квантовых явлений следует, что свойства и взаимодействие материи имеют дискретную основу и меняются с уменьшением масштабов арены действия. Интуиция и здравый смысл подсказывают, что аналогичными особенностями должно обладать также и пространство. Дигитальная топология в состоянии обеспечить необходимый теоретический фундамент на пути перехода от непрерывной картины мира к дискретной.

Эта книга базируется на работах, выполненных автором самостоятельно и в соавторстве с коллегами в 1983 -- 1998 гг., а также на монографии автора "Теория молекулярных пространств и ее приложения к компьютерам, физике и другим областям", опубликованной в 1997 году. Следует отметить, что монография являлась первой книгой в научной литературе, посвященной теории молекулярных пространств. Автор стремился привлечь внимание к дигитальной топологии со стороны как можно более широкого круга читателей, и это обусловило стиль изложения материала и характер построения монографии. Монография носила обзорный характер, технические детали опускались, доказательства теорем и утверждений отсутствовали. Содержание этой книги существенно отличается от содержания монографии в следующих основных пунктах:

В отличие от монографии, где вообще не приводилось доказательств теорем, в этой книге даются подробные доказательства каждого утверждения или теоремы.

В книгу добавлены новые результаты, не вошедшие в монографию, и составляющие значительную часть объема книги.

Изменен порядок изложения материала с целью улучшить логическое построение книги.

Значительно увеличено количество рисунков и иллюстраций, позволяющих лучше понять материал.

Книга имеет практически независимый характер и доступна для понимания без обращения к ссылкам.

При первом знакомстве с книгой совершенно не обязательно читать подряд весь текст. Достаточно прочитать приложение и главы 1, 3, 12, 14. В главе 15 рассказано о компьютерных экспериментах, приведших к пониманию идей и принципов, лежащих в основе теории молекулярных пространств. Тем, кто интересуется физикой, стоит обратить внимание на главу 18, где рассматриваются некоторые очевидные следствия применения теории к структуре дискретного пространства-времени, модели замкнутой вселенной, а также анализируются физические обоснования трехмерности физического пространства. Специалисты в области дифференциальных уравнений могут прочитать главу 16, где предлагается общая структура простейших уравнений математической физики на молекулярном пространстве.

Следует подчеркнуть тот факт, что как теория молекулярных пространств, так и дигитальная топология в целом, возникли как научные теории сравнительно недавно и находятся, по сути дела, лишь в начальной стадии своего становления. Как и во всякой новой теории здесь имеется большое количество специфических проблем, начиная от отсутствия своей собственной терминологии и до отсутствия научных журналов, публикующих работы по этому направлению. С другой стороны, именно благодаря начальному этапу развития здесь существует множество новых направлений для исследования применительно как к фундаментальным, так и к прикладным областям.

Замечание.

Название “дигитальная топология” является, в сущности, нонсенсом, по крайней мере, в отношении к молекулярным пространствам. Действительно, топология считается заданной, если соответствующим образом определено семейство множеств, называемых открытыми. На молекулярном пространстве, невозможно разумным образом определить такую систему, и по этой причине понятие открытого множества заменяется на понятие окаема и шара. Следовательно, молекулярное пространство не является топологическим пространством в традиционном смысле этого слова.

Для чтения и понимания основной части материала этой книги достаточно знания математики на уровне старших классов средней школы. Книга рассчитана на самый широкий круг читателей, желающих познакомиться с этим направлением в математике и компьютерных науках и использующих компьютеры в своей работе. В нескольких главах описываются алгоритмы, позволяющие непосредственно переводить результаты, полученные в книге, в компьютерные программы. В книге имеются повторения там, где это необходимо для облегчения чтения и понимания нового материала.

Пользуясь случаем автор благодарит всех тех, кто помогал ему в разное время в работе над материалом, который составил содержание этой книги и, особенно, своих соавторов по научным статьям: профессора Я.-Н. Йеха (Yeong Nan Yeh) (Тайпейский институт математики, Тайвань), профессора Р. Коппермана (R. Kopperman) (Сити колледж, Нью-Йоркский университет, США), профессора Р. Мелтэра (R. Melter) (Лонг Айленд университет, Саузэмптон, Нью-Йорк, США), доктора Ю. Мухина (Y. Mukhin) (Сиракузский университет, США), Автор искренне признателен за полезные советы и замечания профессору Л. Левитину (L. Levitin) (Бостонский университет, США), профессору Д. Клейтману (D. Kleitman) (Массачусетсский технологический институт, США), профессору Т. Конгу (T. Kong) (Квинс колледж, Нью-Йоркский университет, США), профессору П. Саулсону (P. Saulson) (Сиракузский университет, США), профессору Ж.-К. Рота (G.-C. Rota) (Массачусетсский технологический институт, США). Автор благодарит профессора С. Агафонова (Московский технический университет, Россия), профессора Р. Соркина (R. Sorkin) (Сиракузский университет, США), доктора В. Стэна (W. Stan) (США) за обсуждения, способствовавшие разъяснению многих неясных и спорных вопросов.

Всю ответственность за возможные неточности автор целиком и полностью принимает на себя.

Александр Ивако.

# ЧАСТЬ 1
# МОЛЕКУЛЯРНЫЕ ПРОСТРАНСТВА И ИХ СВОЙСТВА

## 1. МОЛЕКУЛЯРНЫЕ ПРОСТРАНСТВА. ОПРЕДЕЛЕНИЯ И ОПЕРАЦИИ НАД МОЛЕКУЛЯРНЫМИ ПРОСТРАНСТВАМИ

We introduce and define molecular spaces which are digital models of continuos spaces and operations on molecular spaces.

Содержание этой главы частично изложено в работах [3,28,43,45,49,51].

В дигитальной топологии существует несколько подходов к построению дискретных пространств на конечном или счетном множестве точек. Традиционный путь, взятый из классической математики, предполагает задание системы подмножеств, называемых открытыми и

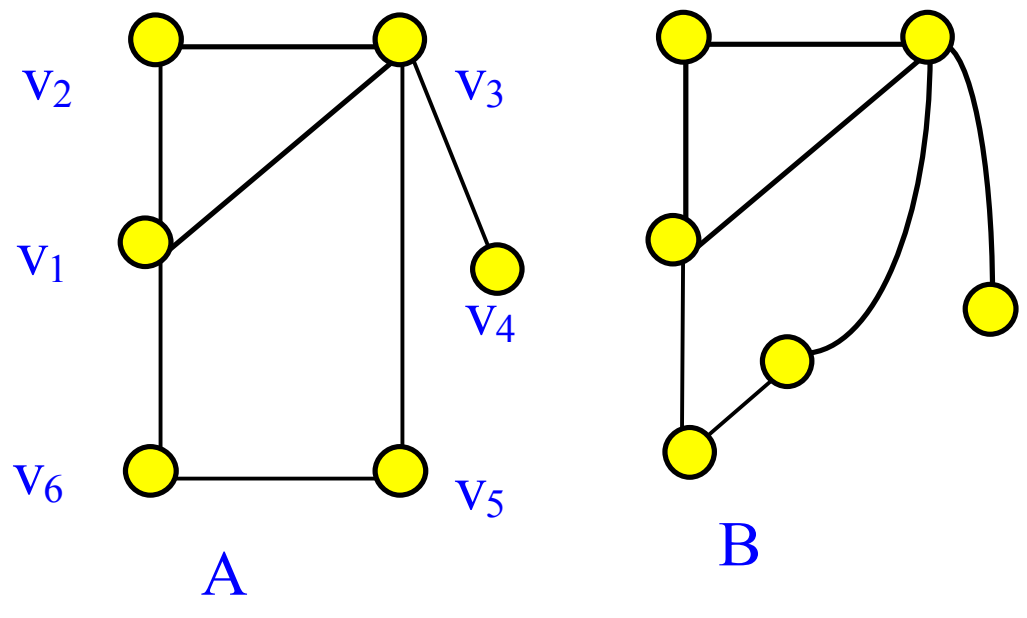

*Рис. А* Молекулярные пространства А и В изоморфны.

удовлетворяющих определенным соотношениям. Такой подход развивался рядом авторов [30,31,32,33,34,35,38,39]. Однако при переходе к пространствам размерности 2 и выше здесь возникают трудности, связанные с тем, что появляется несколько различных типов открытых множеств, причем, количество типов увеличивается с увеличением размерности пространства. В нашем подходе мы, вообще, отказываемся от понятия открытого множества. Для этого существует достаточно оснований. Например, в топологии широко используется понятие нерва покрытия топологического пространства. Каждому открытому множеству, принадлежащему покрытию, ставится в соответствие точка нерва, причем, если два множества имеют непустое пересечение, то точки, им соответствующие, соединены связью. При достаточно мелком покрытии нерв должен отражать все основные свойства топологического пространства. Именно таким образом строится молекулярное пространство. Нерв отражает наиболее фундаментальные топологические свойства молекулярных (и классических топологических, в частности, непрерывных) пространств. Однако следует отметить, что нерв является только одной из нескольких возможных форм описания молекулярных

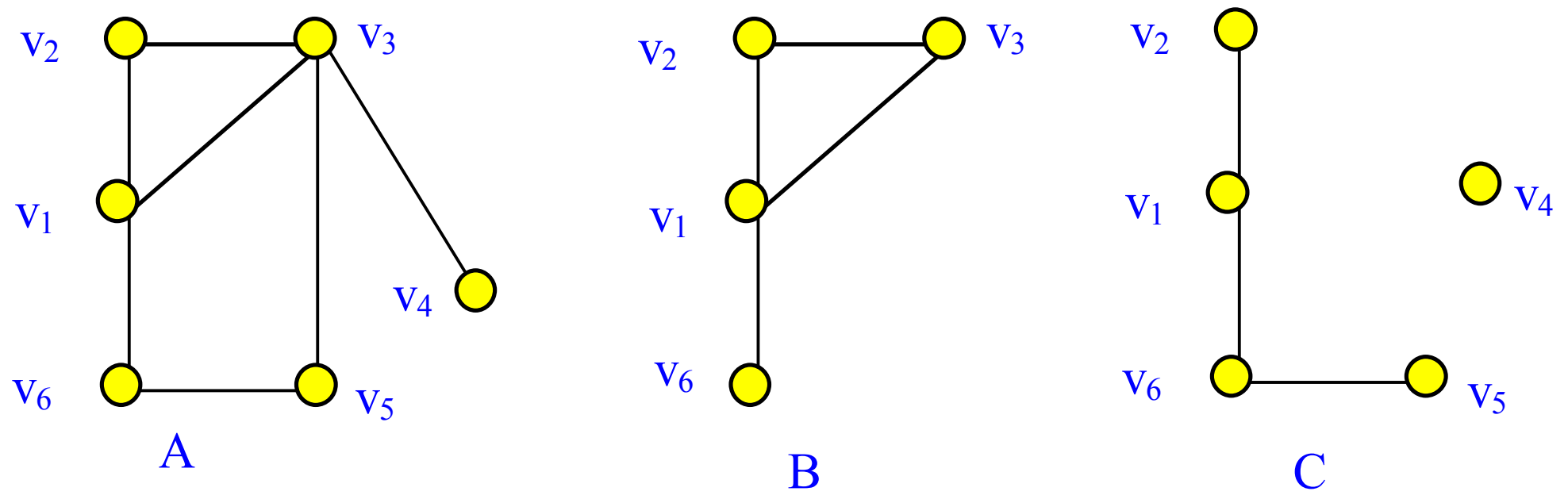

• Рис. В *Молекулярные пространства B и C являются подпространствами пространства A. МП B получено удалением из A точек* $v_4$*, и* $v_5$*, пространство C-удалением точки* $v_3$*. Молекулярные пространства A и B-связны, Молекулярное пространство C-несвязно.* $|A|=6$, $|B|=4$, $|C|=5$.

пространств, в некоторых случаях явно недостаточных. Используя только нерв, мы не в состоянии, например, описать внешнюю геометрию молекулярного пространства при его вложении в более общее молекулярное пространство или в непрерывное евклидово пространство. Нерв, по всей вероятности, удобен, когда речь идет только о базовых, топологических свойствах молекулярного пространства. Таким образом, молекулярное пространство есть пространство, построенное на конечном или счетном множестве точек с определенными топологическими и геометрическими свойствами, такими как размерность, связность, метрика, объем, группы гомологий и тому подобное [43]. Все эти свойства определяются структурой пространства. Возникает вопрос, как определить эту структуру?

Определение молекулярного пространства

Молекулярным пространством G, сокращенно МП, называется дигитальное пространство, состоящее из конечного или счетного множества точек V, на котором задана топологическая структура W, определяемая наличием или отсутствием связей между точками, $G=(V,W)$, где $V=(v_1,v_2,...v_n)$, $W = ((v_p v_q),....)$ (рис. a, рис. b).

Связь между точками $v_p$ и $v_q$ означает, что пара $(v_p,v_q)$ является неупорядоченной, то есть $(v_p v_q)$ и $(v_q v_p)$ есть одно и то же. Если между точками $v_p$ и $v_q$ имеется связь $(v_p v_q)$, то точки $v_p$ и $v_q$ называются

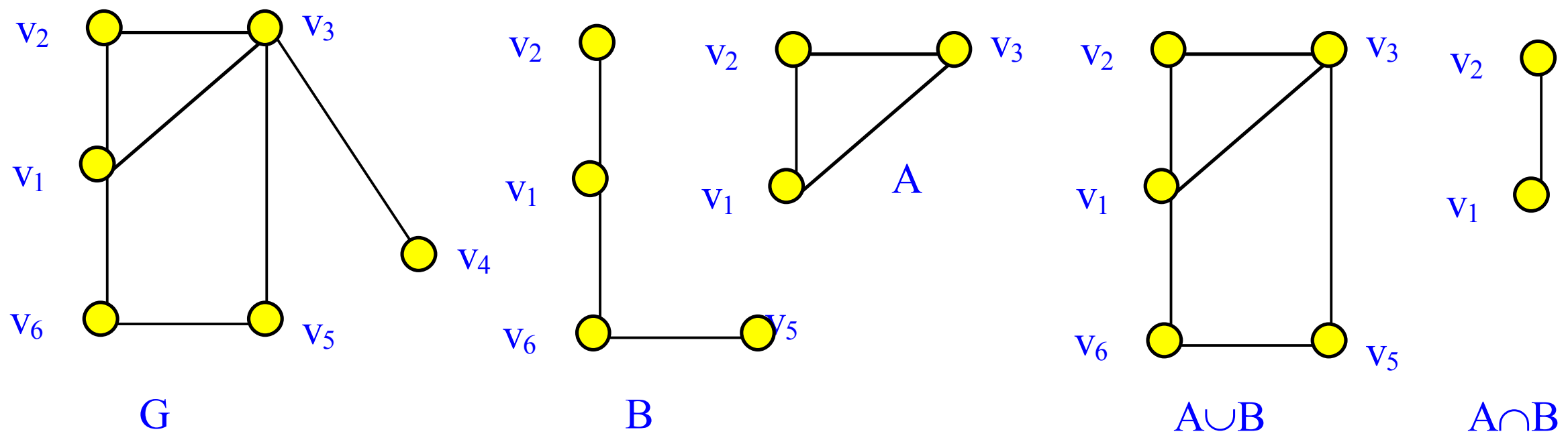

- Рис. С *Объединение A∪B и пересечение A∩B двух подпространств A и B пространства G.*

смежными, соседними или связанными; точки $v_p$ и $v_q$ называются границами связи $(v_p v_q)$; точки $v_p$ и $v_q$ и связь $(v_p v_q)$ называются инцидентными друг другу; связи $(v_p v_q)$ и $(v_q v_s)$ называются смежными, если они инцидентны одной вершине.

Связь между двумя точками предполагает, что они находятся в ближайшем соседстве, причем, все точки, находящиеся в ближайшем соседстве к данной точке, являются равноправными по отношению к этой точке. Таким образом, в молекулярном пространстве каждая точка окружена семейством точек, являющихся ее ближайшими соседями. Эта точка образует связь со всеми точками семейства. Следует подчеркнуть, что метрика, то есть расстояние между двумя точками пространства, определяется естественным образом. Таким образом, для полного определения молекулярного пространства достаточно задать два множества: множество точек V и множество связей W, определяющее топологию пространства.

О б о з н а ч е н и е .

Для молекулярного пространства G, состоящего из конечного или счетного множества точек $V=(v_1,v_2,...v_n)$, будет использоваться обозначение $G=(v_1,v_2,...v_n)$.

Определение. объема МП

Количество точек МП G называется объемом молекулярного пространства G и обозначается $|G|$ (рис. b).

Определение веса МП

Весом пространства G называется число ребер P(G)=|W| этого пространства.

Введем определение изоморфных пространств. Изоморфные пространства есть одинаковые копии, отличаться могут только обозначения.

Определение изоморфизма МП.

Два МП G и H называются изоморфными, если между множествами их точек существует взаимно-однозначное соответствие h: G→H, сохраняющее связь точек. $v_1$ и $v_2$ смежны в G тогда и только тогда, когда $h(v_1)$и $h(v_2)$ смежны в H.

МП A и B, изображенные на рис. a, являются изоморфными. Вообще говоря, между двумя пространствами может быть установлено несколько

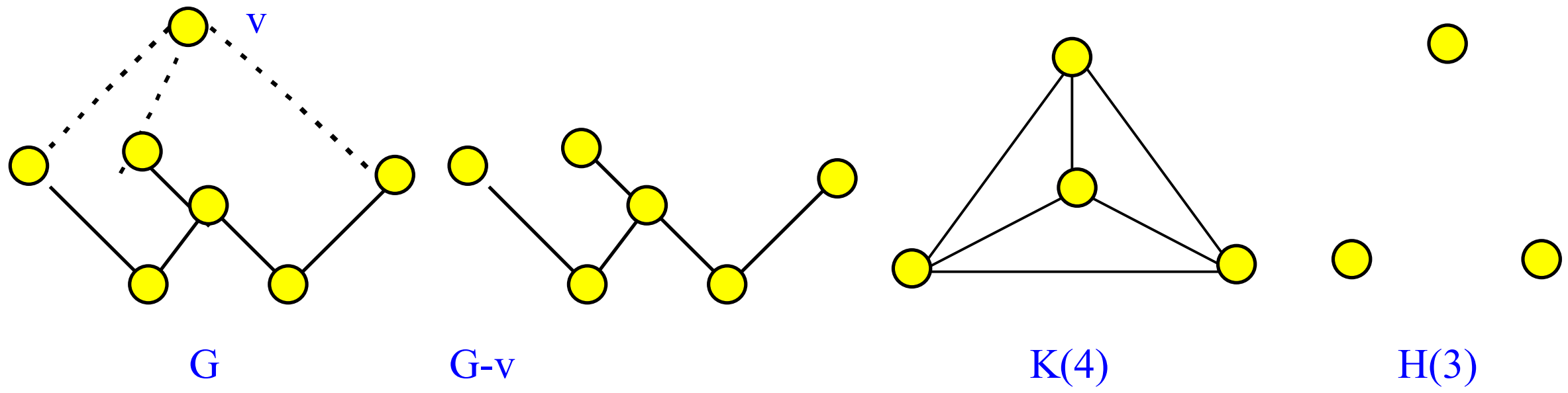

• Рис. D *Пространство G-v получено из пространства G удалением точки v и всех инцидентных ей связей. Пространство K(4) есть связка, пространство H(3) является вполне несвязным.*

отличающихся один от другого изоморфизмов.

Важным является понятие связного пространства. Удобно определить связность через цепь.

Определение цепи.

Цепью в МП G между точками $v_1$ и $v_n$ называется последовательность различных точек $(v_1, v_2,... v_n)$, где соседние точки попарно смежны. Если $v_1$ и $v_n$ в цепи также смежны, то цепь называется циклом.

На рис. b $(v_1,v_2,v_3,v_5,v_6)$ является цепью, и $(v_1,v_2,v_3,v_5,v_6,v_1)$-циклом.

Определение связного МП.

МП называется связным, если любая пара точек соединена цепью. Две точки принадлежат разным компонентам связности МП, если не существует цепи, соединяющей их.

Как и в случае непрерывных пространств нам необходимо иметь определение подпространства. Это связано с определенными сложностями. Дело в том, что подпространство определяется как точками, так и топологией, то есть связями между точками. При выделении подпространства из пространства связи между точками подпространства и остальными точками, определяющие внешнюю топологию подпространства, уже не учитываются. Всегда необходимо иметь это ввиду. Иными словами, связи между точками подпространства определяют его топологию саму по себе, а связи между точками подпространства и остальными точками определяют топологию размещения подпространства во внешнем пространстве.

Определение подпространства.

Индуцированным или порожденным подпространством в МП G называется МП H, если:

1. множество точек H является подмножеством множества точек МП G.

2. связь в МП H существует тогда и только тогда, когда она существует в G.

Следовательно, порожденным подпространством называется максимальное подпространство на данном множестве точек. Это означает, что если точки смежны в МП, то они всегда смежны в подпространстве. В

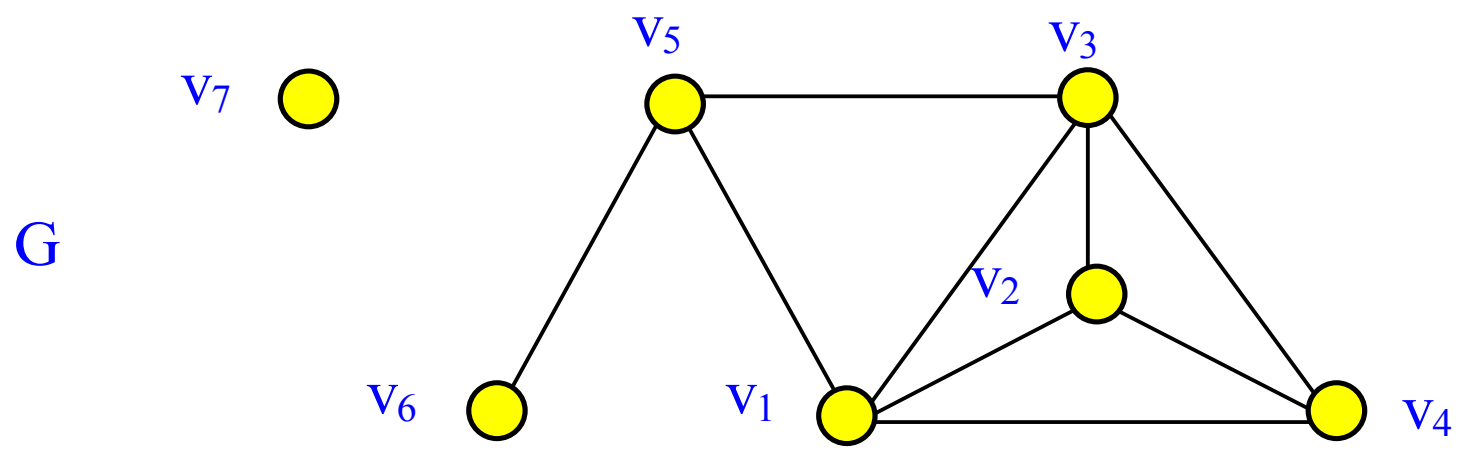

- Рис. E *Пространство G имеет следующие связки: одна связка на одной точке K(1) = ($v_7$); одна связка на двух точках K(2) = ($v_5$, $v_6$); одна связка на трех точках K(3) = ($v_1$, $v_3$, $v_5$); одна связка на четырех точках K(4) = ($v_1$, $v_2$, $v_3$, $v_4$).*

теории МП рассматриваются и другие виды подпространств, кроме порожденных. Однако, так как мы используем только порожденные подпространства (кроме особо оговоренных случаев), мы будем использовать название подпространство для порожденного подпространства. Подпространство также можно получить путем удаления из пространства всех точек (и инцидентным им связей), не принадлежащих данному подпространству (рис. b).

На рис. b пространство A и его подпространства B и C определяются следующими множествами:

A:V=($v_1$,$v_2$,$v_3$,$v_4$,$v_5$,$v_6$),

W=(($v_1$,$v_2$),($v_1$,$v_3$),($v_1$,$v_6$),($v_2$,$v_3$),($v_3$,$v_4$),($v_3$,$v_5$),($v_5$,$v_6$)).

B:V=($v_1$,$v_2$,$v_3$,$v_6$), W=(($v_1$,$v_2$),($v_1$,$v_3$),($v_1$,$v_6$),($v_2$,$v_3$)).

C:$V=(v_1,v_2,v_4,v_5,v_6)$, $W=((v_1,v_2),(v_1,v_6),(v_5,v_6))$.
В дальнейшем нам также понадобятся операции объединения и пересечения двух или нескольких подпространств данного пространства (рис. с). Эти операции мы определим стандартным образом на множестве точек подпространств. Сформулируем эти операции в виде определений.

Определение объединения подпространств.

Объединением двух подпространств A с множеством точек $V_p = (v_{p1}, v_{p2},\dots v_{pk})$ и B с множеством точек $V_q = (v_{q1}, v_{q2},\dots v_{qm})$ данного пространства G с множеством точек $V = (v_1, v_2,\dots v_k)$ называется подпространство $A\cup B$ с множеством точек $V_p\cup V_q = (v_{p1}, v_{p2},\dots v_{pk})\cup(v_{q1}, v_{q2},\dots v_{qm})$, являющимся теоретико-множественным объединением подмножеств $V_p\subseteq V$ и $V_q\subseteq V$.

Аналогично вводится определение пересечения подпространств.

Определение пересечения подпространств.

Пересечением двух подпространств A с множеством точек $V_p = (v_{p1}, v_{p2},\dots v_{pk})$ и B с множеством точек $V_q = (v_{q1},v_{q2},\dots v_{qm})$ данного пространства G с множеством точек $V = (v_1, v_2,\dots v_k)$ называется подпространство $A\cap B$ с множеством точек $V_p\cap V_q = (v_{p1}, v_{p2},\dots v_{pk})\cap(v_{q1}, v_{q2},\dots v_{qm})$, являющимся теоретико-множественным пересечением подмножеств $V_p\subseteq V$ и $V_q\subseteq V$.

Для пересечения $A\cap B$ подпространств A и B мы также будем использовать обозначение A|B.
Объединение подпространств является операцией, более чувствительной к объемлющему пространству, чем их пересечение. Пересечение подпространств работает на любом пространстве и однозначно для любых подпространств. Объединение подпространств зависит от объемлющего пространства в следующем смысле: в объединении могут возникнуть связи между точками двух подпространств, которых нет в каждом из подпространств в отдельности. Например, на рис. с подпространства на точках $(v_1, v_2)$ и $(v_2,v_3)$ содержат каждое в отдельности по одной связи, тогда как их обединение-подпространство на точках $(v_1,v_2,v_3)$ содержит уже три связи. Таким образом, появилась связь, не принадлежащая ни одному из подпространств. Объединение $A\cup B$ содержит связь $(v_3,v_5)$, которая отсутствует как в A, так и в B. То есть множество связей не только объединяется, но и могут появиться дополнительные связи. При пересечении двух подпространств этого не происходит, множество связей пересекается. Это всегда необходимо учитывать, работая с молекулярными пространствами. Еще раз подчеркнем, что в дальнейшем, если не оговорено особо, под объединением подпространств всегда понимается подпространство, порожденное объединением множеств точек.

Определение разности G-v.

Подпространством G-v называется подпространство, полученное из пространства G вычеркиванием точки v, (рис. d).

Это определение можно расширить на подпространство H пространства G.

Определение разности G-H.

Подпространством G-H называется подпространство, полученное из пространства G вычеркиванием всех точек подпространства H.

Очевидно, что $G=(G-v)\cup v$, $G=(G-H)\cup H$.

Определение связки K(n) и вполне несвязного пространства H(n).

Пространство K(n) на n точках называется связкой (полным связным пространством), если каждая пара его точек смежна. Пространство H(n) на n точках называется вполне несвязным, если оно не содержит связей. (рис. d).

Определение клики.

Любая максимальная полная связка как подпространство называется кликой (рис. e).

В комбинаторной топологии любая связка на p точках называется (p-1)-симплексом. Такой симплекс имеет по определению размерность (p-1).

Определение прямой суммы пространств G⊕H.

Пусть имеется два пространства: G с набором точек $V= (v_{1,}, v_{2,}, v_{3,}.... v_{n,})$ и набором связей W, и H с набором точек $R= (r_{1,}, r_{2,}, r_{3,}.... r_m)$ и набором связей S. Прямой суммой (суммой) двух пространств G и H назовем пространство G⊕H, состоящее из точек V и R, связей W и S, и всех связей, соединяющих точки V и R (рис. f).

Как видно, при этой операции число точек пространства G⊕H равно сумме точек пространств G и H, $|G\oplus H|=|G|\oplus|H|$, число связей $P)G)=|W|+|S|+nm$.

Определение конуса G⊕v.

Если одно из пространств, например, второе в предыдущем определении является точкой v, то сумма G⊕v называется конусом (cone) пространства G (рис. f).

Определение прямой суммы подпространств G⊕H.

Пусть в пространстве A имеется два подпространства: G с набором точек $V=(v_1,v_2,v_3,....v_n)$ и H с набором точек $R=(r_1,r_2,r_3,....r_m)$, причем множества V и H не пересекаются. Если каждая точка из V смежна с каждой точкой из R, то подпространство $G\cup H=G\oplus H$, состоящее из точек V и R, называется прямой суммой (суммой) двух подпространств G и H.

<u>Определение р-дольного пространства.</u>

Пространство $K(n_1,n_2,...n_p)$ называется р-дольным, если его можно представить как сумму р вполне несвязных пространств $H(n_1)$, $H(n_2)...H(n_p)$, $K(n_1,n_2,...n_p)=H(n_1)\oplus H(n_2)\oplus...\oplus H(n_p)$.

На рис. f изображены двудольные пространства К(1,3) и К(2,3).
Одним из наиболее важных понятий в теории молекулярных пространств являются окаем и шар точки и подпространства. Эти два понятия выполняют в ТМП ту же функцию, что и открытые и замкнутые окрестности в топологии. На дигитальном пространстве затруднительно,

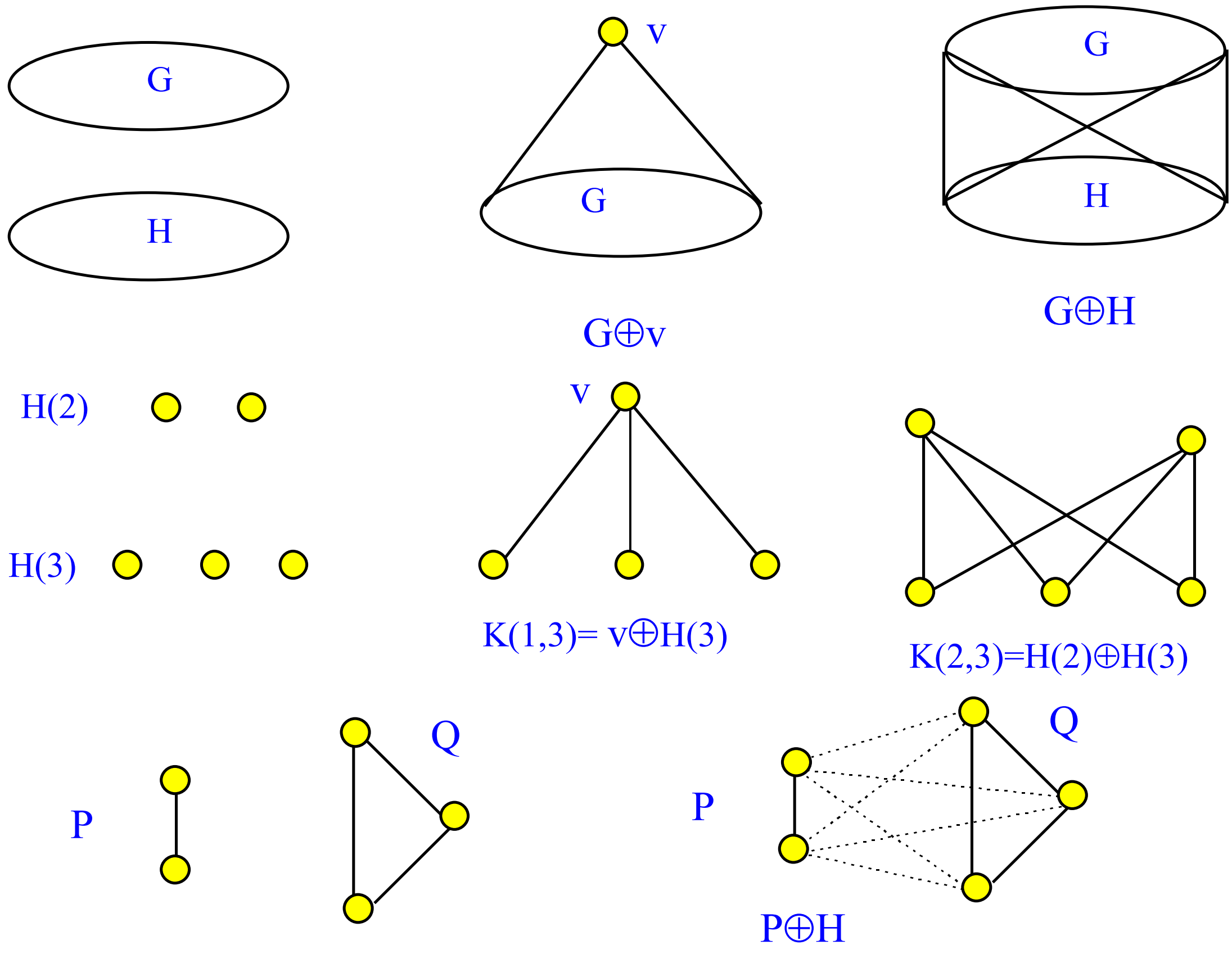

• Рис. F *Прямые суммы пространств v⊕H, G⊕H, v⊕H(3).*

если не невозможно разумным образом определить понятия открытого и замкнутого множеств, аналогичных классическим. Это вполне объяснимо, поскольку понятие открытого множества основано на непрерывности и бесконечности. В ТМП как науке, возникшей из потребностей компьютерной практики, непрерывность и бесконечность не используются.

<u>Определение шара и окаема точки.</u>

Пусть G и v являются пространством и его точкой. Назовем шаром (ball) точки v подпространство U(v), содержащее как точку v, так и все точки пространства G, смежные с точкой v (рис. g, рис. h).

Назовем окаемом (rim) точки v подпространство O(v), содержащее все точки пространства G, смежные с v, но не содержащее самой точки v.

Очевидно, что O(v) получается из U(v) путем вычеркивания точки v, O(v) = U(v)-v (рис. g, рис. h). Связь между шаром и окаемом определяется соотношениями: шар точки v является конусом и объединением точки v и ее окаема. $U(v) = v \oplus O(v) = v \cup O(v)$, $O(v) = U(v) - v$.

Определение общего окаема точек.

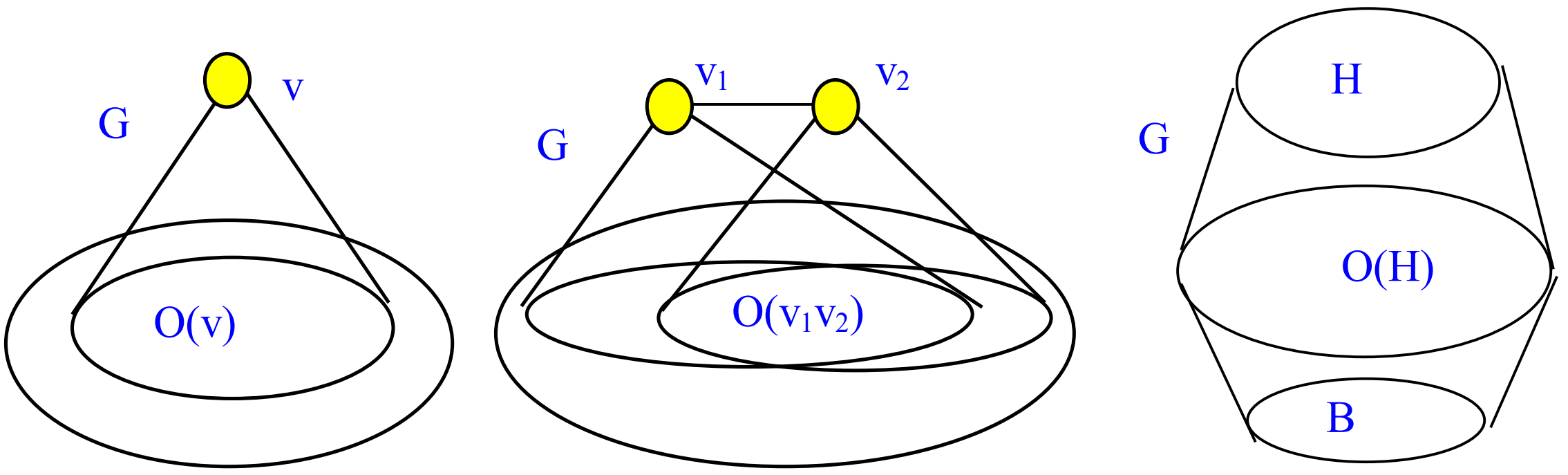

- Рис. G *Окаем O(v) точки v и общий окаем O($v_1v_2$) двух точек $v_1$ и $v_2$ в молекулярном пространстве G. Шар точки v есть подпространство, содержащее O(v) и точку v. Окаем*

Подпространство O($v_1,v_2,v_3,....v_P$), образуемое пересечением O($v_1$), O($v_2$), O($v_3$),...O($v_P$), называется общим окаемом (joint rim) точек $v_1$, $v_2$, $v_3$,.... $v_P$ (рис. g, рис. h).

$O(v_1, v_2, v_3, .... v_P) = O(v_1) \cap O(v_2) \cap O(v_3) \cap ... \cap O(v_P)$.

Иными словами, каждая точка пространства O($v_1,v_2,v_3,....v_P$) должна быть смежна всем точкам $v_1$, $v_2$, $v_3$,.... $v_P$. Операция пересечения рассматривается на множестве точек пространства G. При этом точки $v_1$, $v_2$, $v_3$,.... $v_P$ не обязаны быть смежными. На рис. h, серым цветом изображен окаем двух точек O(vu). Аналогично дается определение шара и окаема подпространства H пространства G. Окаем подпространства H играет роль границы подпространства, множества точек, составляющих ближайшее окружение H.

Определение шара и окаема подпространства.

Пусть G и H являются пространством и его подпространством. Назовем шаром (ball) или полным шаром подпространства H подпространство U(H), содержащее как точки подпространства H, так и все точки пространства G, смежные с хотя бы огдной точкой подпространства H (рис. h, рис. i). Назовем окаемом (rim) или полным окаемом подпространства H подпространство O(H), содержащее все точки пространства G, смежные с хотя бы одной точкой подпространства H, кроме точек подпространства H. Пусть H содержит точки $v_{1,}$, $v_{2,}$, $v_{3,}$.... $v_P$. Тогда U(H) есть объединение шаров всех точек, принадлежащих H.

Следующие свойства полного окаема и полного шара очевидны из определений.

С в о й с т в а .

Пусть G и H являются пространством и его подпространством. Пусть H содержит точки $v_{1,}$, $v_{2,}$, $v_{3,}$.... $v_P$. Тогда:

$U(H)=U(v_1)\cup U(v_2)\cup U(v_3)\cup ...\cup U(v_P)$.

$O(H)=O(v_1)\cup O(v_2)\cup O(v_3)\cup ...\cup O(v_P)-H$.

$O(H)= U(H)-H=U(v_1)\cup U(v_2)\cup U(v_3)\cup ...\cup U(v_P)-H$.

$U(H)=O(H)\cup H$.

На рис. h точка v имеет черный цвет, ее окаем O(v) изображен точками серого цвета и, следовательно, черная и серые точки вместе обозначают шар U(v) точки v.

Очевидна разница между полным окаемом O(H) подпространства H, состоящего из точек $(v_1,v_2,...v_p)$, и общим окаемом точек, его составляющих. Еще раз напомним, что

$O(H)=O(v_1)\cup O(v_2)\cup O(v_3)\cup ...\cup O(v_P)-H$

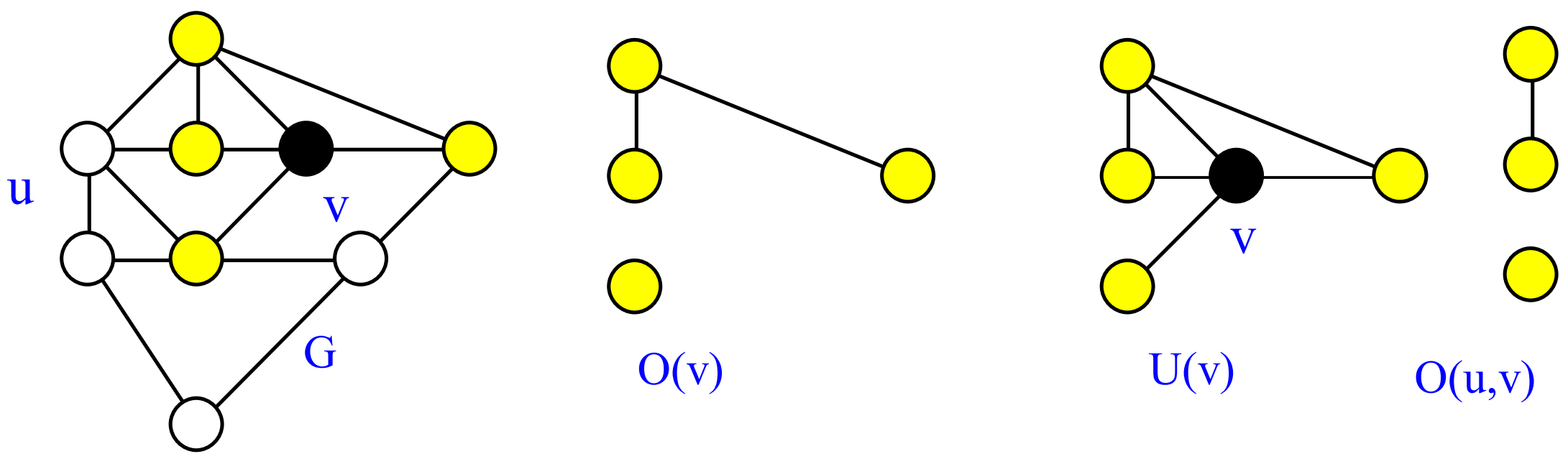

• Рис. H *Окаем O(v) и шар U(v) точки v в пространстве G. Крайний справа-общий окаем точек u и v.*

$O(v_1,v_2,v_3,....v_P)=O(v_1)\cap O(v_2)\cap O(v_3)\cap ...\cap O(v_P)$.

Важным понятием является метрика, заданная на молекулярном пространстве. Обычно метрикой называется неотрицательная функция d(u,v) двух точек u и v, заданная на всем пространстве и удовлетворяющая трем аксиомам.

$d(u,v)\geq 0$ и $d(u,v)=0$ тогда и только тогда, когда $u=v$.

$d(u,v) = d(v,u)$

$d(u,v)+d(v,w)\geq d(u,w)$ для любых точек u, v, w.

Метрику можно определить несколькими способами. Мы введем определение, естественным образом вытекающее из топологии

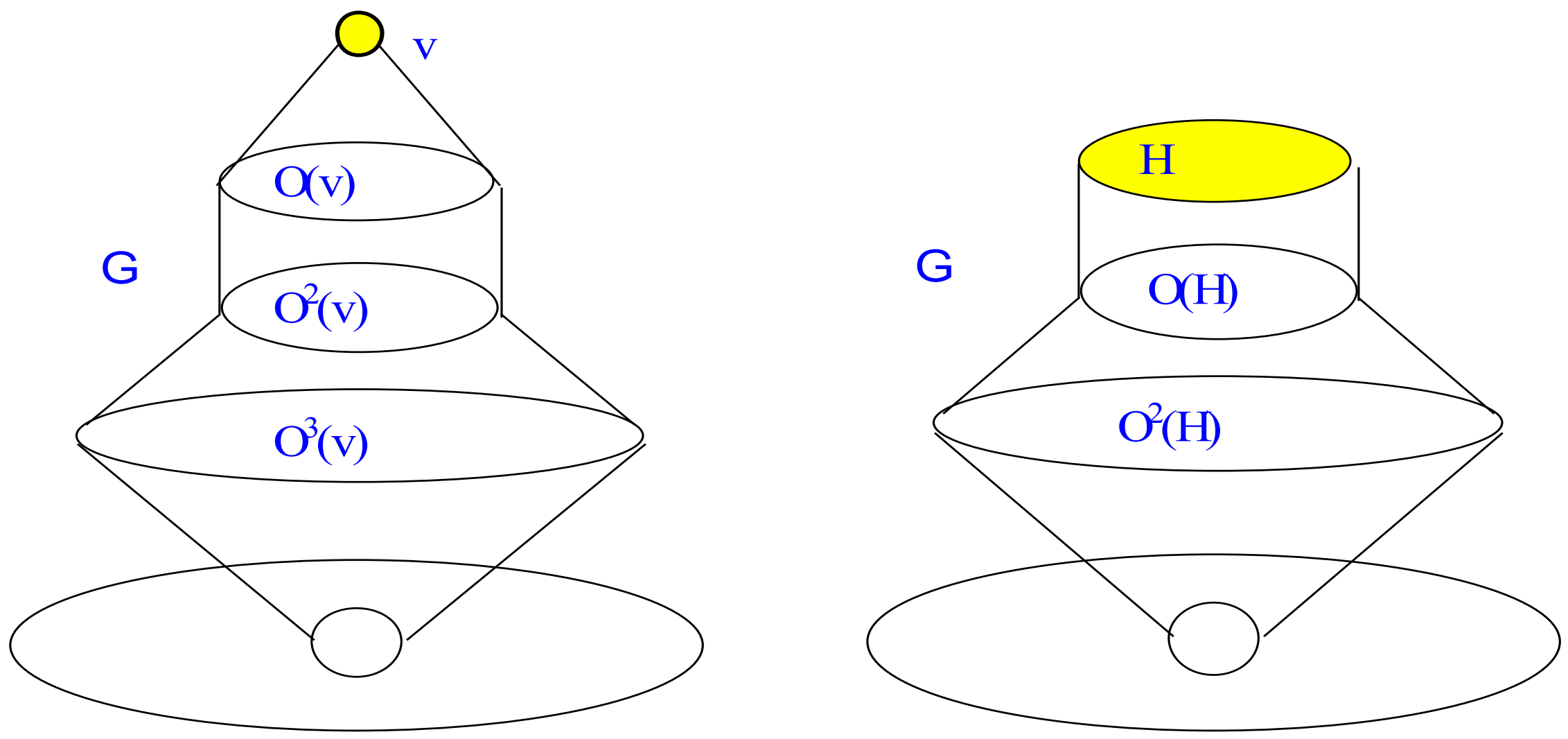

- Рис. I *$O(v)=O^1(v)$, $O^2(v)$ и $O^3(v)$ есть первый, второй и третий полные окаемы точки v. $O(H)=O^1(H)$ и $O^2(H)$ являются первым и вторым полными окаемами подпространства H.*

молекулярного пространства. Вначале определим последовательную серию полных окаемов и полных шаров точки и подпространства.

<u>Определение полного n-окаема и n-шара точки.</u>

Пусть v есть некоторая точка молекулярного пространства G.

Полным 2-окаемом $O^2(v)$ точки v называется полный окаем ее шара $U(v)$, $O^2(v)=O(U(v))$ (рис. i).

Полным 2-шаром $U^2(v)$ точки v называется полный шар подпространства $U(v)$, $U^2(v)=U(U(v))$ (рис. i).

Полным n-окаемом $O^n(v)$ точки v называется полный окаем полного (n-1)-шара $U^{n-1}(v)$ точки v, $O^n(v)=O(U^{n-1}(v))$.

Полным n-шаром $U^n(v)$ точки v называется полный шар полного (n-1)-шара $U^{n-1}(v)$ точки v, $U^n(v)=U(U^{n-1}(v))$.

Очевидна связь между n-окаемом и n-шаром точки v.

С в о й с т в а .

$U(v)=O(v)\cup v$

$U^2(v)=O^2(v)\cup U(v)= O^2(v)\cup O(v)\cup v.$

.......................

$U^n(v)=O^n(v)\cup U^{n-1}(v)= O^n(v)\cup O^{n-1}(v)\cup ...O(v)\cup v.$

Нам также понадобится понятие полных n-окаема и n-шара подпространства H пространства G. Введем эти понятия через обобщение

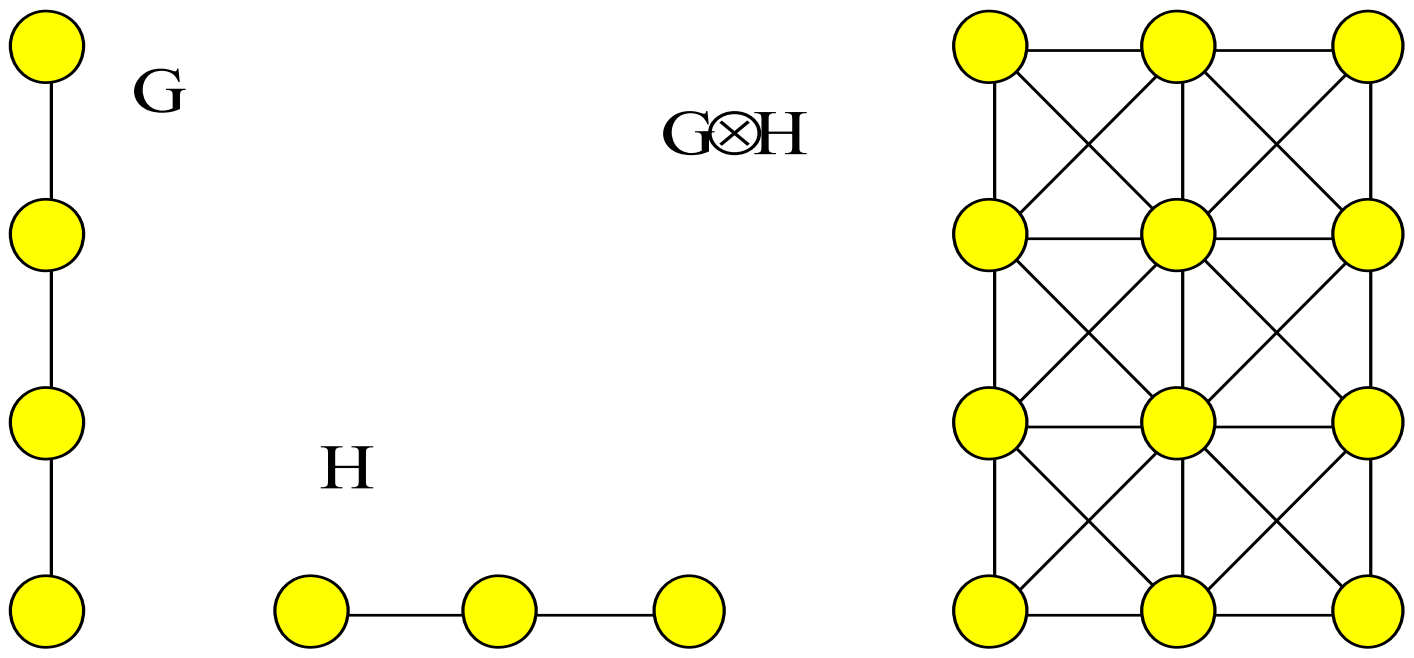

• Рис. J *Прямое произведение пространств G и H. Каждая точка пространства G заменяется на пространство H. Между смежными и одноименными точками смежных копий H устанавливаются связи.*

предыдущих определений на подпространство H пространства G.

<u>Определение полного n-окаема и n-шара подпространства.</u>

Пусть H есть некоторое подпространство молекулярного пространства G. Полным n-окаемом $O^n(H)$ подпространства H называется полный окаем подпространства $U^{n-1}(H)$, $O^n(H)=O(U^{n-1}(H))$. Полным n-шаром $U^n(H)$ называется полный шар подпространства $U^{n-1}(H)$, $U^n(H)=U(U^{n-1}(H))$.

Очевидна связь между n-окаемом и n-шаром подпространства H.

Свойства.

$U(H)=O(H)\cup H$

$U^2(H)=O^2(H)\cup U(H)= O^2(H)\cup O(H)\cup H.$

......................

$U^n(H)=O^n(H)\cup U^{n-1}(H)= O^n(H)\cup O^{n-1}(H)\cup ...O(H)\cup H.$

На рис. i изображены различные окаемы точки v и подпространства H в пространстве G. Так как полный шар U(v) точки v есть $v\oplus O(v)$, $U(v)=v\oplus O(v)$, то второй полный окаем точки v есть полный первый окаем подпространства U(v), $O^2(v)=O(U(v))$. Очевидно также, что $U^2(v)=O^2(v)\cup U(v)$.

<u>Определение расстояния (метрики).</u>

Пусть v есть некоторая точка молекулярного пространства G. Расстояние d(v,u) от v до u равно n, d(v,u)=n, если точка u принадлежит полному n-окаему $O^n(v)$ точки v, $u\in O^n(v)$. Если v и u принадлежат разным компонентам связности МП, то $d(u,v)=\infty$.

Легко видеть, что определенное таким образом расстояние удовлетворяет всем метрическим аксиомам и равно числу связей между точками в самой короткой по числу точек цепи, соединяющей точки v и u.
Расстояние между смежными точками всегда равно 1. Аналогично вводится расстояние между подпространством H пространства G и точкой, а также между двумя подпространствами H и F.

$d(H,v)=n \Leftrightarrow v\in O^{n}(H)$

$d(H,F)=n \Leftrightarrow O^{n-1}(H)\cap F=\varnothing,\ O^{n}(H)\cap F\neq\varnothing.$

<u>Определение прямого произведения G⊗H.</u>

Пусть имеется два пространства: G с набором точек ($v_{1,}$, $v_{2,}$, $v_{3,}$,.... $v_m$...) и набором связей W, и H с набором точек ($u_{1,}$, $u_{2,}$, $u_{3,}$,.... $u_n$...) и набором связей S. Прямым произведением (произведением) двух пространств G и H назовем пространство G⊗H, состоящее из точек $w_{kp}=v_k\otimes u_p$, в котором шар $U(v_k\otimes u_p)$ каждой точки $w_{kp}=v_k\otimes u_p$ состоит из точек $w_{st}=v_s\otimes u_t$, где $v_s$ принадлежит $U(v_k)$, $u_t$ принадлежит $U(u_p)$,

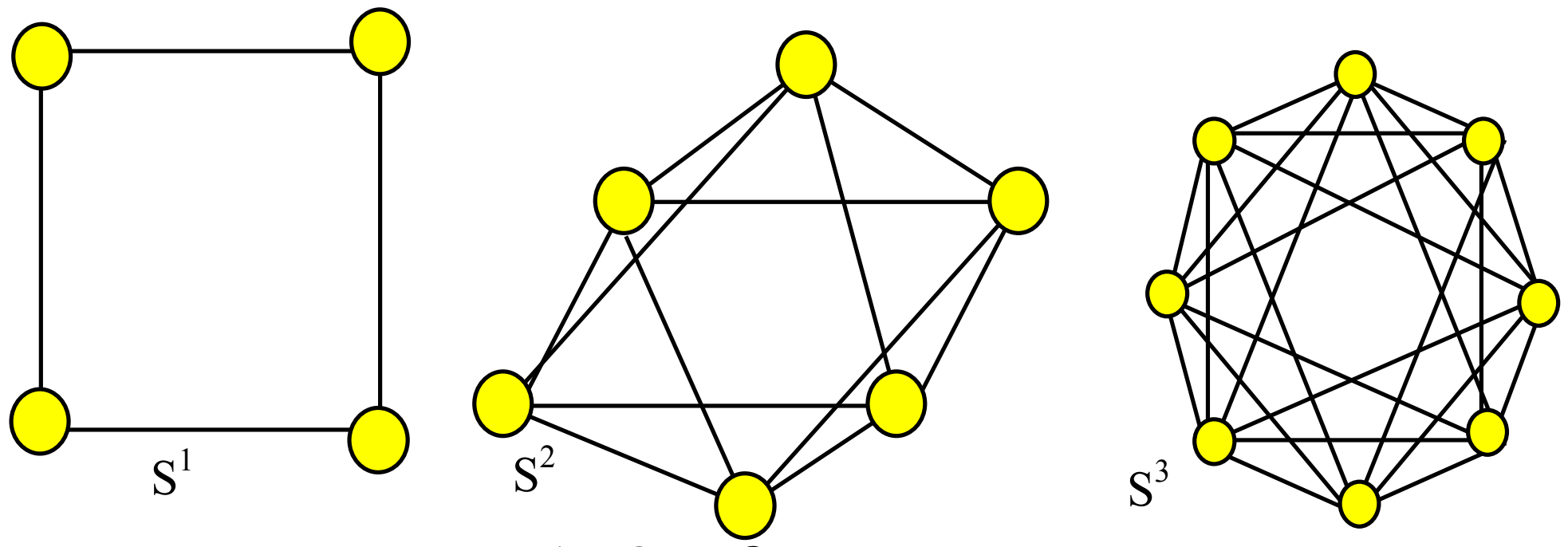

- Рис. K Пространства $S^1$, $S^2$ и $S^3$ являются однородными пространствами.

$w_{st}=v_s\otimes u_t\in U(v_k\otimes u_p) \Leftrightarrow v_s\in U(v_k),\ u_t,\in U(u_p).$

Прямое произведение молекулярных пространств является аналогом топологического произведения пространств, или прямой суммы линейных векторных пространств [53]. В определении отчетливо видна связь между классическим определением топологического произведения и прямым произведением молекулярных пространств. Наглядно прямое произведение можно представить следующим образом. Каждая точка одного из пространств, например G, заменяется пространством H. Полученное пространство имеет структуру G, в котором точками служат копии пространства H. Затем, если две копии H соответствуют двум смежным точкам в G, то устанавливаются связи как между одноименными точками, так и между смежными точками этих копий независимо от того, какой копии они принадлежат (рис. j).
Введем также понятие однородного молекулярного пространства. Однородность является важным свойством пространства. Стандартное евклидово пространство, например, является однородным. Пространства, используемые в физике, также являются, в большинстве, однородными. Молекулярные пространства могут существовать как в однородной, так и в неоднородной формах. Однородная форма во многих случаях является предпочтительной, поскольку пространство описывается более просто. Тем не менее, будет показано, что некоторые, хорошо известные в классической математике пространства не имеют однородной формы, по крайней мере в своем минимальном виде. К таким пространствам относится, например, проективная плоскость, получаемая заклеиванием лентой Мебиуса дырки на сфере.

Определение однородного пространства.

> Молекулярное пространство называется однородным, если окаемы (или шары) двух любых точек пространства являются изоморфными (рис. k).

Иными словами, окаемы являются одинаковыми копиями и при правильном совмещении совпадают один с другим.
В заключение этого раздела следует отметить, что молекулярное пространство можно рассматривать также как граф. Нерв покрытия топологического пространства есть не что иное, как граф. В сущности любой граф является молекулярным пространством с определенными математическими характеристиками. Поэтому все, что применимо к графам, с полным правом может быть применимо и к молекулярным пространствам. Для того, чтобы понять связь между молекулярным пространством и графом, дадим определение графа [21,22].

Определение графа.

> Графом G называется пара множеств (V,W), где V является множеством точек, $V=(v_1,v_2,...v_n)$, W является множеством неупорядоченных пар различных точек из V, $W = ((v_p, v_q),....)$.

В литературе точки графа называются также вершинами и узлами, а пары-связями, ребрами, дугами или линиями. Если две точки соединены связью, то есть образуют пару, то они называются смежными или связными. Связь (u,v) называется инцидентной точкам u и v. Две связи называются смежными, если они инцидентны одной и той же точке.
Следует сказать несколько слов по поводу определений и терминологии. Даже если брать только математику, одни и те же термины в различных областях этой науки используются для обозначения, часто, совершенно различных понятий и объектов. Когда возникает новое научное направление, неизбежно возникают трудности каким образом использовать уже существующие термины для обозначения новых объектов. Наиболее приемлемый путь, как нам кажется, это введение новых обозначений и терминологии, не существующей в тех областях, на которых базируется данное направление. Однако, часто это вызывает активные возражения коллег, работающих в смежных областях. Поэтому, наиболее безболезненным будет некоторый средний путь, когда часть терминов берется из традиционных областей науки по принципу схожести; новые же термины и обозначения вводятся только тогда, когда они совершенно необходимы. С другой стороны, часто бывает, что невозможно напрямую использовать привычную и подходящую по смыслу терминологию, взятую из какого-либо раздела математики, например, теории графов, в дигитальной топологии. Дело здесь в том, что одни и те же названия часто используются для определения различных понятий в различных разделах математики. Например, с точки зрения специалиста в теории графов полный граф, рассматриваемый как молекулярное пространство, логично было бы назвать полным молекулярным пространством. Однако, такое название было бы неестественным с точки зрения топологии, в которой используется понятие полного пространства. Поэтому, чтобы избежать недоразумений в таких случаях мы будем вводить определения, не допускающие множественного толкования. Это позволит, хотя бы частично, избежать путаницы и противоречивости. Хотя эта глава посвящена именно определениям, в ней мы ввели только базовые понятия, совершенно необходимые для дальнейшего продвижения в выбранном направлении. Это не значит, что мы уже определили все, что нам необходимо. Однако, новые определения будут появляться в тексте по мере необходимости там, где без них нельзя совершенно нельзя обойтись.

## Список литературы к главе 1.

# ТОЧЕЧНЫЕ ПРЕОБРАЗОВАНИЯ МОЛЕКУЛЯРНЫХ ПРОСТРАНСТВ

We introduce and study contractible molecular spaces and special class of transformations of molecular spaces that are called contractible. These transformations are much similar to homotopic transformations of classical topology spaces. In this connection we define homotopic molecular spaces and study their properties.

Результаты, полученные в этой главе, частично изложены в работах [15,26,43,49,52].

## *ТОЧЕЧНЫЕ ПРОСТРАНСТВА. ТОЧЕЧНЫЕ ПРЕОБРАЗОВАНИЯ МОЛЕКУЛЯРНЫХ ПРОСТРАНСТВ И ИХ СВОЙСТВА*

В дальнейшем нам понадобятся преобразования, позволяющие видоизменять молекулярное пространство, например, увеличивать или уменьшать количество его точек, устанавливать или удалять связи между различными точками. В определенной мере такие преобразования моделируют преобразования непрерывных пространств. Рассмотрим простейший пример. Предположим, что в каком-либо процессе

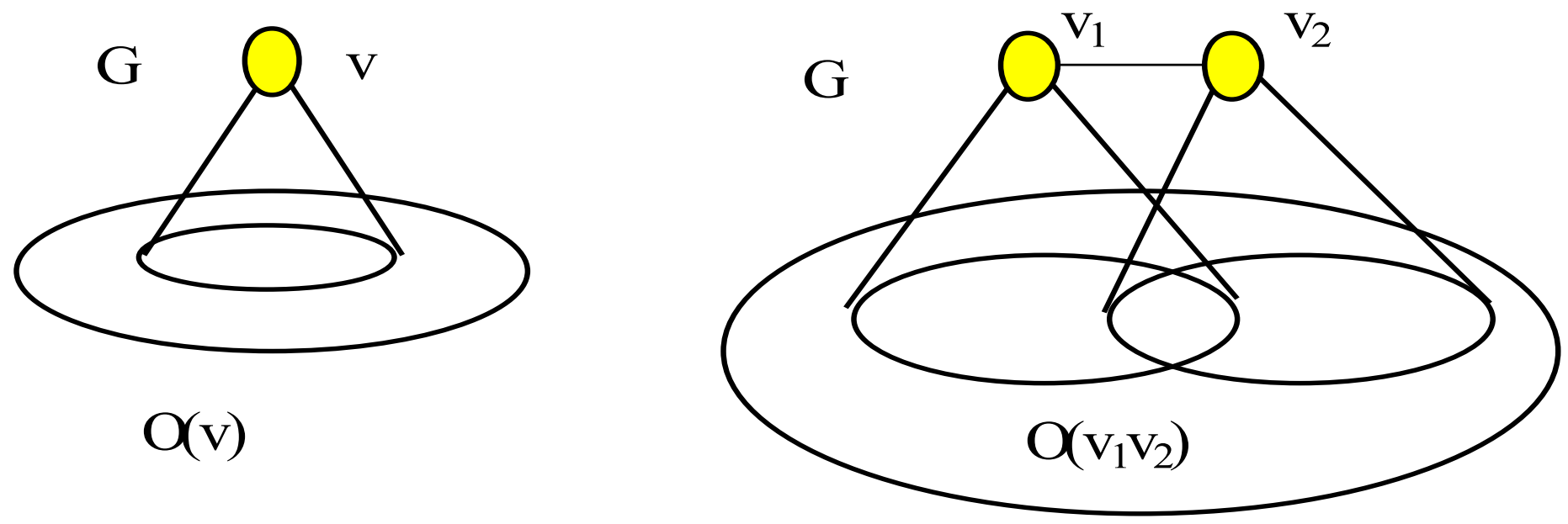

*Рис. 12* Окаем O(v) точки v в МП G является точечным, поэтому v может быть удалена (приклеена) из МП. Общий окаем $O(v_1v_2)$ двух точек $v_1$ и $v_2$ является точечным и поэтому связь $(v_1v_2)$ может быть удалена (приклеено) из МП G.

материальный шар увеличивает свой объем. Молекулярный шар, являющийся дигитальным образом этого непрерывного шара, также должен увеличивать свой объем, например, за счет увеличения числа точек, его образующих. Иными словами, к молекулярному шару добавляются новые точки, но при этом форма и другие инвариантные характеристики этого шара не меняются.

Возникает вопрос, каким образом мы должны добавлять точки к уже имеющимся, чтобы на каждом этапе пространственный объект сохранял все свои свойства. При этом, имеются ввиду строгие математические

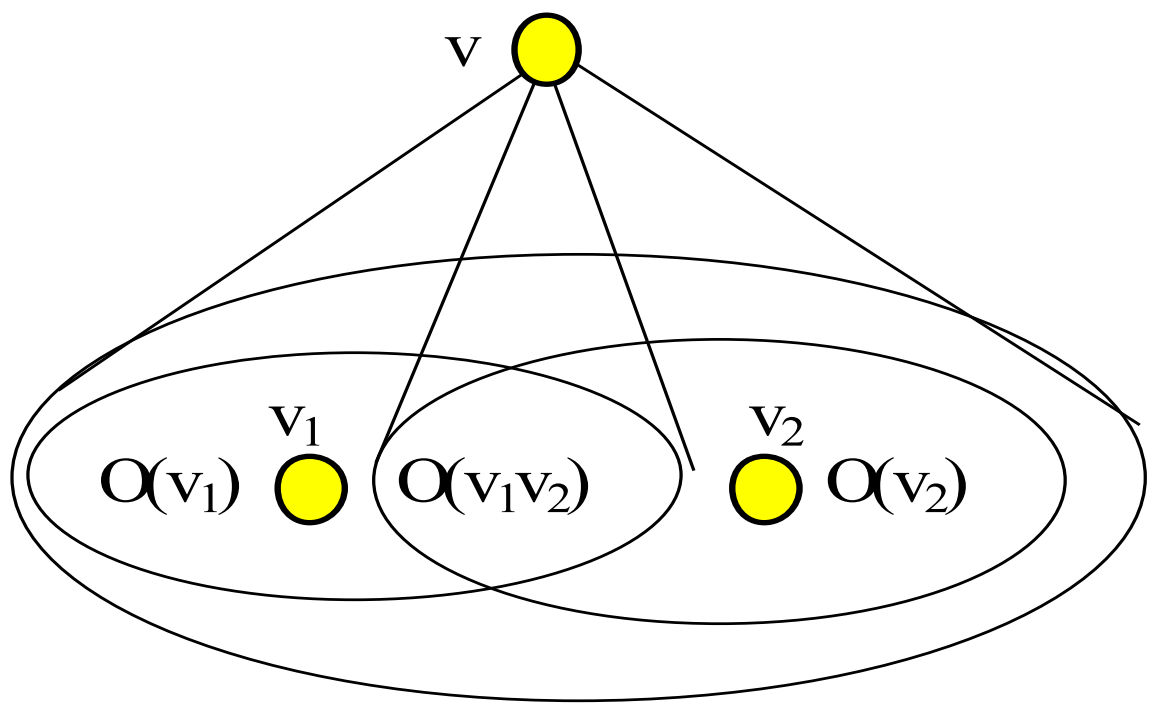

*Рис. 13* В пространстве $v\oplus G$ общий окаем $O(v_1v_2) = v\oplus(O(v_1v_2)|G)$ точек $v_1$ и $v_2$ является конусом, и поэтому точки $v_1$ и $v_2$ могут быть соединены ребром.

характеристики пространственного объекта. Преобразования, которые мы будем вводить ниже, как раз и призваны ответить на этот вопрос на наиболее общем уровне.. В предыдущих работах мы представили точечные преобразования молекулярных пространств состоящими из четырех операций. Предварительно дадим несколько необходимых определений. Переходим теперь к определению точечных (contractible) преобразований молекулярных пространств.

Определение точечных пространств

Семейство T молекулярных пространств $G_1$, $G_2$, $G_3$,...$G_n$,..., $T=(G_1,G_2,G_3,...G_n,...)$, называется точечным (contractible), если:

1. Тривиальное молекулярное пространство K(1), состоящее из одной точки, принадлежит T.
2. Любое молекулярное пространство семейства T может быть получено из тривиального молекулярного пространства с помощью точечных преобразований.

Все молекулярные пространства семейства T называются точечными (contractible).

Определение точечных преобразований

Следующие преобразования называются точечными (contractible):

1. Удаление точки v.

Точка v молекулярного пространства G может быть удалена, если окаем O(v) этой точки является точечным молекулярным пространством, то есть O(v) принадлежит T.

2. Приклеивание точки v.

Точка v может быть приклеена к молекулярному пространству G по молекулярному подпространству $G_1$, то есть соединена связями со всеми точками молекулярного пространства $G_1$, если $G_1$ является точечным молекулярным пространством то есть $G_1$ принадлежит T.

3. Удаление связи ($v_1v_2$).

Связь ($v_1v_2$), соединяющая две смежные точки $v_1$ и $v_2$ молекулярного пространства G, может быть удалена, если общий окаем O($v_1$,$v_2$) этих точек является точечным молекулярным пространством, то есть O($v_1$,$v_2$) принадлежит T.

4. Приклеивание связи ($v_1v_2$).

Две несмежные точки $v_1$ и $v_2$ молекулярного пространства G, могут быть соединены связью ($v_1v_2$), если общий окаем O($v_1$,$v_2$) этих точек является точечным молекулярным пространством, то есть O($v_1$,$v_2$) принадлежит T.

Очевидно, что семейство T всех точечных молекулярных пространств может быть получено из тривиального молекулярного пространства K(1) последовательным применением преобразований 1-4. На Рис. 12 схематично изображены точечные преобразования. Мы будем говорить, что точка v, или связь ($v_1v_2$) точечные, если O(v), или O($v_1v_2$) являются точечными молекулярными пространствами.

Точечные преобразования играют важную роль в теории молекулярных пространств. Они позволяют установить связь между двумя

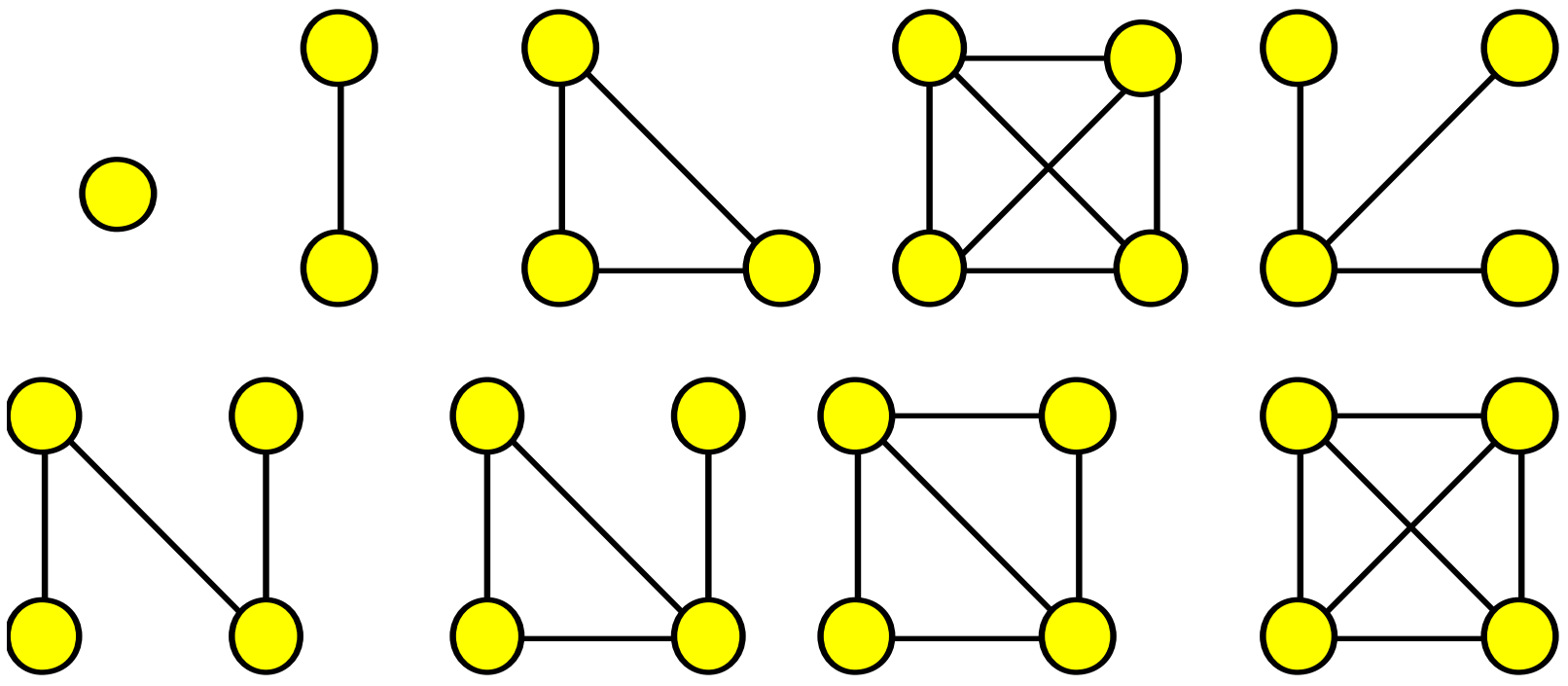

*Рис. 14* Изображены все точечные молекулярные пространства с числом точек менее пяти. Каждое из них может быть стянуто к одной точке, изображенной вверху справа точечными отбрасываниями точек.

молекулярными пространствами, являющимися дискретными образами двух топологически эквивалентных непрерывных поверхностей. А именно, если две поверхности топологически эквивалентны, то молекулярное пространство одной из них может быть приведено путем применения точечных преобразований к молекулярному пространству другой поверхности. Однако, класс точечных преобразований более широк и соответствует скорее преобразованиям между гомотопически эквивалентными пространствами. Иными словами, если две поверхности гомотопически эквивалентны, то их молекулярные пространства связаны серией точечных преобразований. Об этом мы будем писать дальше.
А сейчас введем еще одно определение, смысл которого более понятен с учетом вышесказанного.

Определение гомотопных пространств

> Два молекулярных пространства G и H называются гомотопически эквивалентными или гомотопными (G~H), если одно из них может быть получено из другого путем точечных преобразований.

Согласно этому определению, все точечные молекулярные пространства гомотопны тривиальному молекулярному пространству K(1). На Рис. 14 изображены все точечные молекулярные пространства с числом точек не более четырех. Докажем несколько важных для дальнейшего теорем, которые конкретизируют свойства точечных преобразований.

Теорема 1

> *Всякая связка K(n), имеющая n связных точек, является точечным МП.*

Доказательство.
Используем метод математической индукции. Для малых n точечность молекулярных пространств проверяется непосредственно. Предположим, что для произвольного n K(n) является точечным. Очевидно, что молекулярное пространство K(n+1) получается из молекулярного пространства K(n) при помощи приклеивания точки к v молекулярному пространству K(n)=O(v), что является точечным преобразованием 2. Следовательно, K(n+1) является точечным. Теорема доказана.□
Следующая теорема является обобщением теоремы 1 и важным инструментом при доказательстве многих других утверждений.

Теорема 2

> *Конус v⊕G любого молекулярного пространства G является точечным.*

Доказательство.
Для молекулярных пространств с малым числом точек теорема проверяется непосредственно. Предположим, что теорема верна для всех молекулярных пространств G с числом точек $|G| \le n$, то есть молекулярные пространства v⊕G являются точечными. Пусть G является молекулярным пространством с числом точек n+1, $|G|= n+1$. Выберем две несмежные точки $v_1$ и $v_2$ молекулярного пространства G. Тогда их окаем $O(v_1v_2) = v\oplus(O(v_1v_2)|G)$ в v⊕G является конусом. Так как $|O(v_1v_2)|G|<n$, то из индукционного предположения следует, что окаем $O(v_1v_2) = v\oplus(O(v_1v_2)|G)$ есть точечное молекулярное пространство, и мы можем соединить точки $v_1$ и $v_2$ связью. Повторяя эту операцию для всех несмежных точек молекулярного пространства мы получаем полное МП K(n+2) на n+2 точках, которое являются точечным МП согласно предыдущей теореме. Отбрасывая те же самые связи в обратном порядке, мы переходим путем точечных преобразований от точечного МП K(n+2) к исходному(Рис. 13). Следовательно, $v\oplus G$ является также точечным. Это завершает доказательство.□

Теорема 3

*Всякое точечное МП является связным.*

Д о к а з а т е л ь с т в о.

Для точечных МП с малым числом точек теорема проверяется

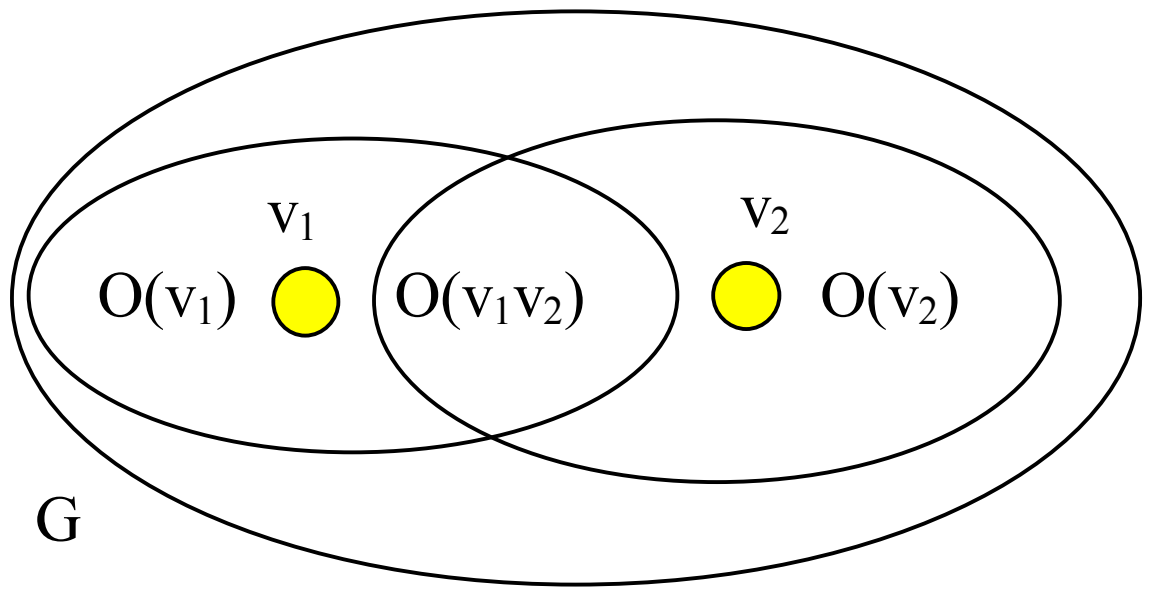

*Рис. 15* Если G точечное и точка $v_1$ несмежна с некоторыми точками, то среди них найдется точечная точка $v_2$.

непосредственно. Предположим, что теорема верна для всех молекулярных пространств G с числом точек $|G|<n$. Пусть G является точечным пространством с числом точек n. Приклеивание связи или точки, очевидно, оставляет пространство связным. При отбрасывании точки или связи по точечному подпространству H $|H|<n$, так как H связно, то связность G также не нарушается. Теорема доказана (Рис. 16). □

Введем аксиому, которая легко проверяется на МП с малым количеством точек.

А к с и о м а.

Пусть G есть точечное МП, и его точка v является несмежной с хотя бы одной из точек G. Тогда существует хотя бы одна точка u, несмежная v, такая, что их общий окаем O(vu) есть точечное подпространство (Рис. 15).

Эта аксиома может быть сформулирована несколькими различными способами. В качестве аксиомы может быть взята, например, Теорема 4.

На МП с малым количеством точек аксиома проверяется непосредственно. Возможно, что в дальнейшем эту аксиому удастся доказать как теорему.

Теорема 4

*Пусть G есть точечное МП с числом точек большим единицы. Тогда оно имеет, по крайней мере, две точечные точки*

Доказательство

Предположим, что точка v в МП G является смежной всем остальным точкам, то есть $G=v\oplus U$. Тогда любая точка МП U будет точечной, так как ее окаем есть конус. Предположим, что G не является конусом. Пусть $V = (v_1,v_2,...v_n)$ есть множество точек. Выберем $v_1$ и в соответствии с аксиомой установим связь с несмежной точкой $v_p$, так как $O(v_1v_p) \in T$. Повторяем

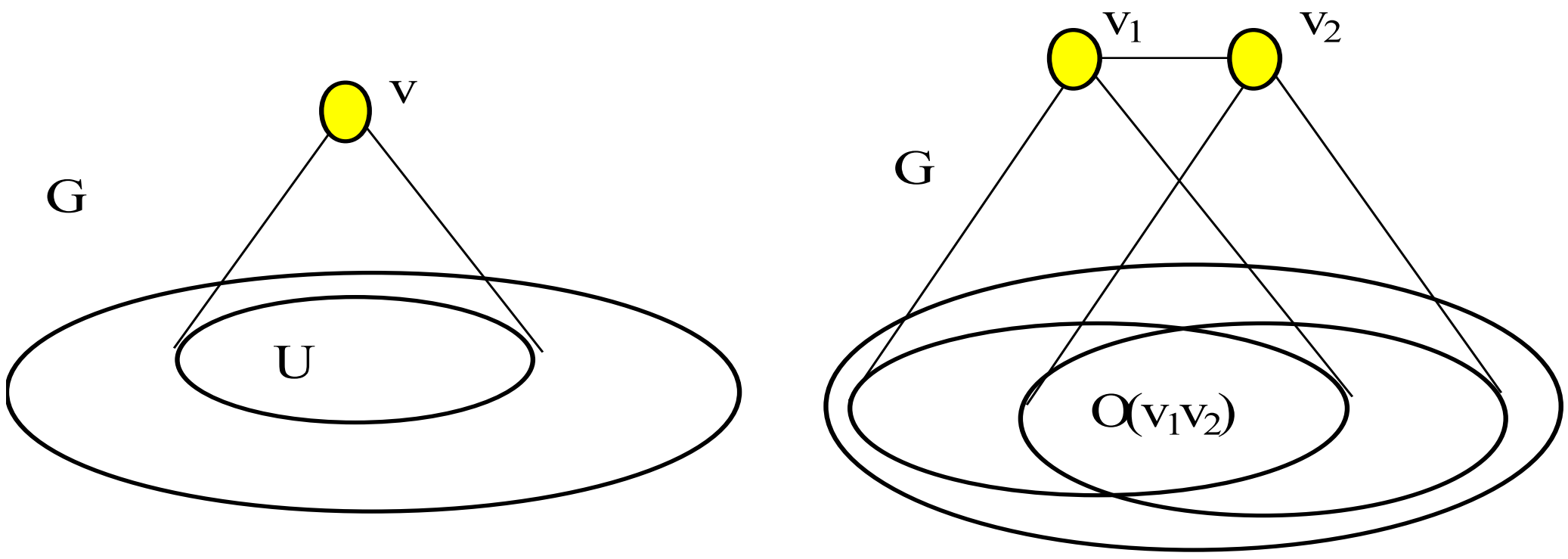

*Рис. 16* При приклеивании точки v к точечному МП G по точечному подпространству U связность сохраняется. Точно так же связность сохраняется при приклеивании связи $(v_1v_2)$ по точечному общему окаему $O(v_1v_2)$.

этот процесс, пока не соединим $v_1$ со всеми остальными точками МП G. Очевидно, что последняя точка $v_s$ является точечной, так как $O(v_1v_s) = O(v_s)$ после установления всех предыдущих связей. Получаем первую точечную точку в G. Тот же самый процесс повторяем начиная с точки $v_s$. Получаем вторую точечную точку $v_2$. Теорема доказана.□

Эта теорема имеет множество применений в этом разделе. В большинстве теорем ниже мы используем эту теорему.

Следствие

Как следствие теоремы 4 можно сформулировать следующее предложение: В точечном МП G для любой точки v найдется отличная от нее точечная точка.

Теорема 5

*Всякое точечное МП на n точках может быть получено из связки*

*K(n) на n точках путем точечных удалений связей и наоборот. При этом связи, инцидентные любой данной точк, могут быть удалены в последнюю очередь (Рис. 17).*

Доказательство.

Предположим, что G = $v_1 \oplus U$ есть конус, и $V=(v_1,v_2,...v_n)$ есть множество его точек. Мы можем приклеить любую связь ($v_p v_s$), так как O($v_p v_s$) есть конус и, следовательно, точечное подпространство. Приклеиваем все недостающие связи и получаем K(n). Отбрасывая связи в обратном порядке мы получаем G из K(n).

Пусть G ≠ $v_1 \oplus U$. Выберем несмежную с $v_1$ точку $v_s$ в соответствии с аксиомой такую, что O($v_1 v_s$) есть точечное подпространство. Соединим эти две точки связью и повторим этот процесс, пока не установив все связи, соединяющие точку $v_1$ с остальными точками. Все другие связи можно приклеить в произвольном порядке и получить связку K(n). Отбрасывая

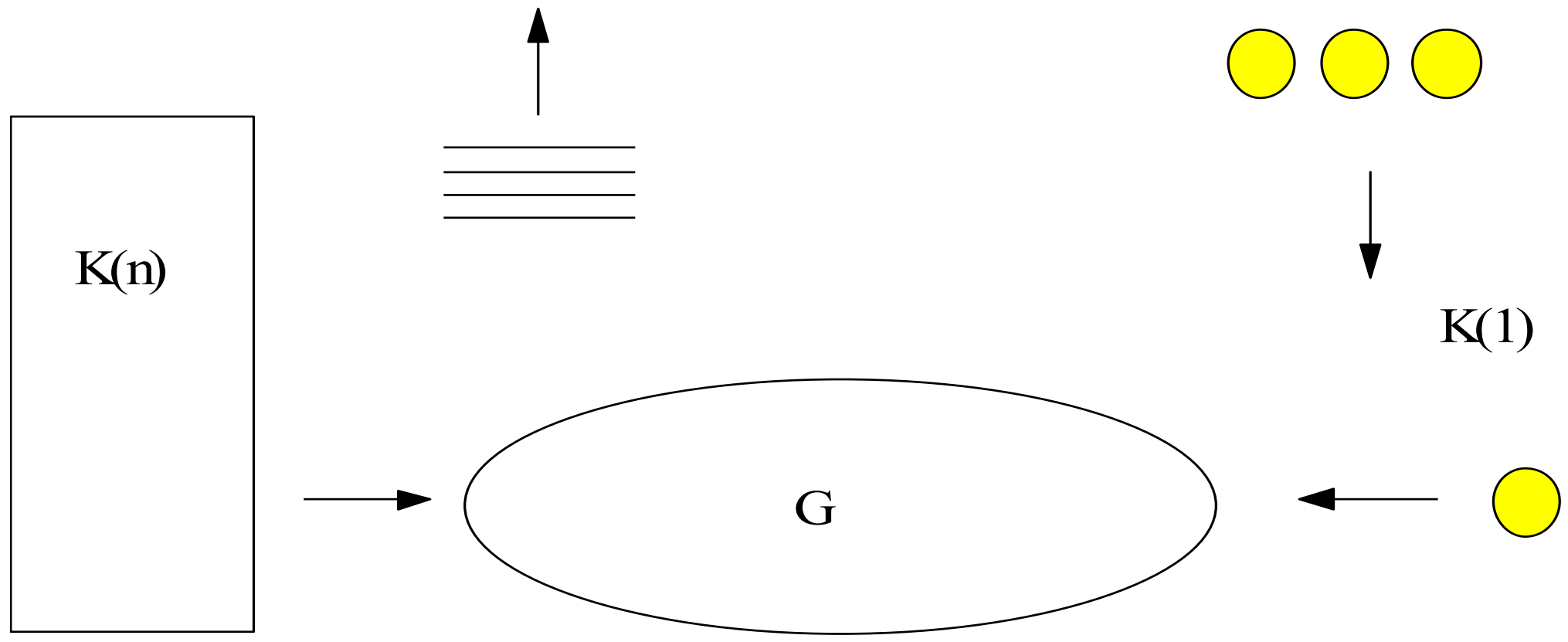

*Рис. 17* Слева. Любое точечное молекулярное пространство G, |G|= n, может быть получено из связки K(n) отбрасыванием только связей. Справа. Любое точечное пространство G, |G|= n, может быть получено из одноточечного МП K(1) приклеиванием только точек.

связи в обратном порядке, получаем G из K(n). Доказательство закончено.□

Нам хотелось бы еще раз подчеркнуть, что МП G получается из МП K(n) только при помощи операций отбрасывания связей без затрагивания точек МП.

Следующая теорема позволяет получить точечное МП G только при помощи приклеивания точек к одноточечному МП K(1).

Теорема 6

*Всякое точечное МП может быть получен из тривиального МП K(1) путем точечных приклеиваний точек и наоборот (Рис. 17).*

Доказательство.

Для точечных МП с небольшим числом точек теорема проверяется непосредственно.

Пусть G есть точечное МП. Тогда в соответствии с теоремой оно имеет по крайней мере две точечных точки, одна из которых может быть отброшена. Продолжая процесс отбрасывания, мы получаем одноточечное МП K(1). Приклеивая точки в обратном порядке, мы из K(1) получаем G. Доказательство закончено.□

Суть этой теоремы в том, что в сущности нет необходимости использовать

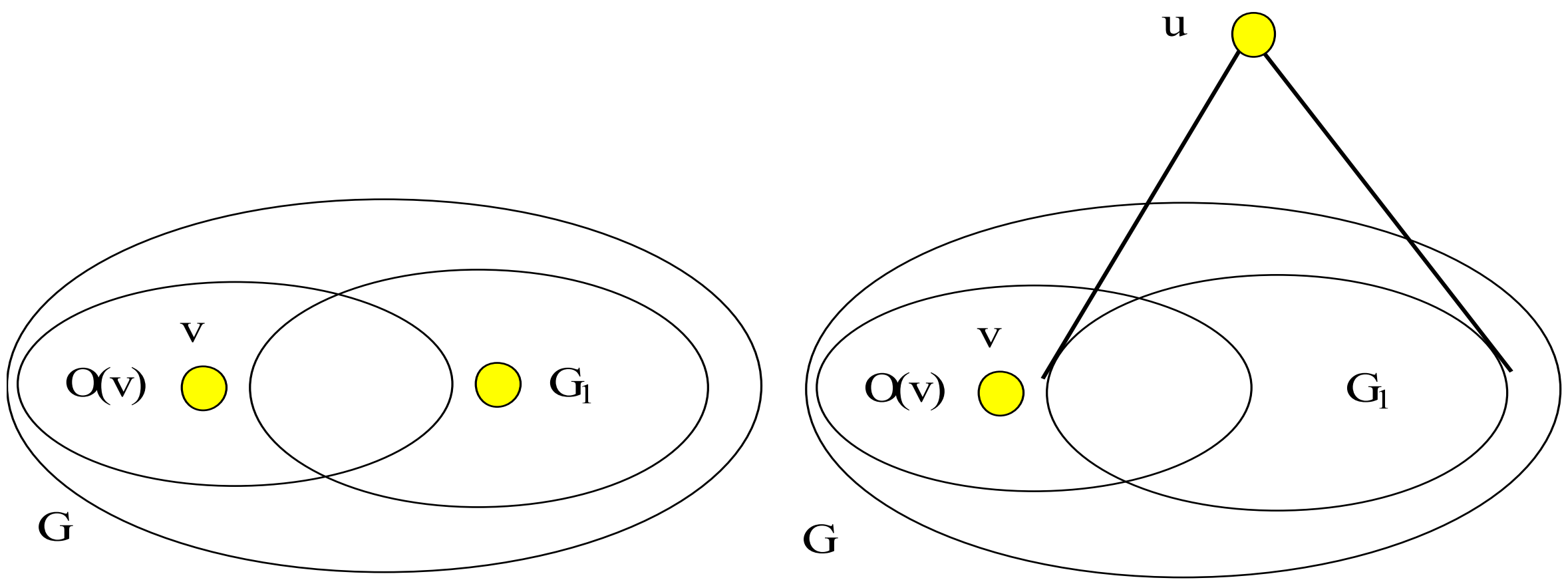

*Рис. 18* В точечном МП G для любого точечного подпространства $G_1$, не совпадающего с G, существует точечная вершина v, не принадлежащая G.

приклеивания и отбрасывания связей. Достаточно обойтись точками. Используя эту теорему, можно дать определение точечных преобразований на основе приклеиваний и отбрасываний только точек.

<u>Определение В точечных пространств</u>

Семейство T МП $G_1$, $G_2$, $G_3$,...$G_n$,..., T = ($G_1$,$G_2$,$G_3$,...$G_n$,...), называется точечным (contractible), если:

1. Тривиальное МП K(1), состоящее из одной точки, принадлежит T.
2. Любое МП семейства T может быть получено из тривиального МП с помощью точечных преобразований.

Все МП семейства T называются точечными (contractible).

<u>Определение В точечных преобразований</u>

Следующие преобразования называются точечными (contractible):

1. Удаление точки v.

Точка v МП G может быть удалена, если окаем O(v) этой точки является точечным МП, то есть O(v) принадлежит T.

2. Приклеивание точки v.

Точка v может быть приклеена к МП G по подпространству $G_1$, то есть соединена связями со всеми точками МП $G_1$, если $G_1$ является точечным МП то есть $G_1$ принадлежит T.

Теорема 7

*Пусть G и $G_1$ являются точечным МП и его точечным подпространством, отличным от G. Тогда в G существует точечная точка v, не принадлежащая $G_1$ (Рис. 18).*

Доказательство.

Приклеим точку u к G таким образом, что $O(u)=G_1$. Полученное МП $G\cup(u\oplus G_1)$ является точечным. Согласно следствию, для u в $G\cup(u\oplus G_1)$ существует несмежная с u точечная точка v. Следовательно, v не принадлежит $G_1$. Доказательство закончено.□

Теорема 8

*Пусть G и $G_1$ являются точечным МП и его точечным подпространством, отличным от G. Тогда $G_1$ может быть преобразован в G последовательным приклеиванием точек, оставаясь каждый раз подпространством G.*

Доказательство.

Приклеим точку $v_1$ к G таким образом, что $O(v_1) = G_1$ (Рис. 19). Полученное МП $G\cup(v_1\oplus G_1)$ является точечным. Согласно аксиоме существует несмежная с $v_1$ точка $v_2$ такая, что $O(v_1v_2) = O(v_2)\cap G_1$ есть точечное МП. Установление связи в $G\cup(u\oplus G_1)$ равносильно приклеиванию точки $v_2$ к МП $G_1$ по $O(v_2)\cap G_1$. Очевидно, точечное приклеивание связей, инцидентных точке $v_1$ можно продолжить, пока

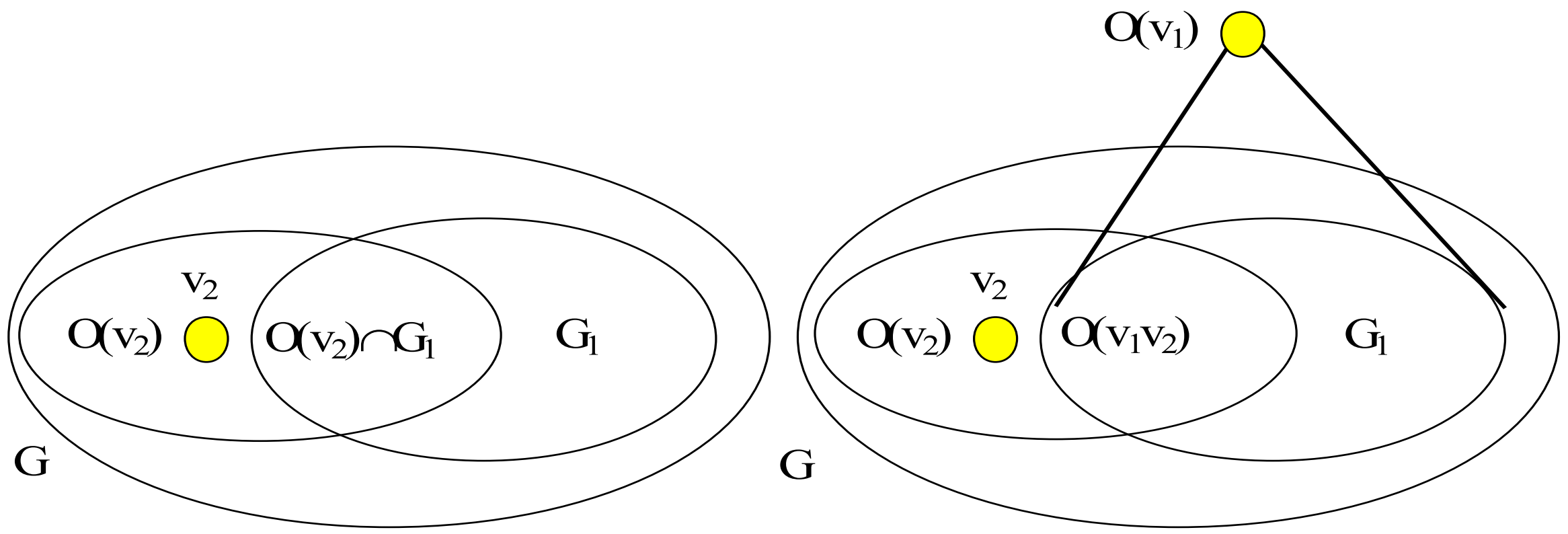

*Рис. 19* Приклеивание точки к подпространству $G_1$ по точечному подпространству $O(v_2)\cap G_1 = O(v_1v_2)$ равносильно приклеиванию связи между точками $v_1$ и $v_2$ в МП G.

точка $v_1$ не будет смежной всем остальным точкам МП G (Рис. 18, Рис. 19). Этот процесс равносилен приклеиванию точек к МП $G_1$. При этом $G_1$ переходит в G. Доказательство закончено.□

Теорема 9

$G \in T$, $G\text{-}v \in T \Rightarrow O(v) \in T$;
$G \in T$, $O(v) \in T \Rightarrow G\text{-}v \in T$;
$O(v) \in T$, $G\text{-}v \in T \Rightarrow G \in T$.

Доказательство.

Эти три утверждения следуют непосредственно из определения точечного приклеивания и отбрасывания точки. □

Теорема 10

*Пусть H есть произвольное МП и G-его точечное подпространство(Рис. 20). Приклеивание точки v к H по G может быть заменено приклеиванием v к любой точке из G и последующим последовательным приклеиванием связей между v и остальными точками из G.*

Доказательство.

Пусть точка $v_1$ принадлежит G, $v_1 \in G$. Соединим v и $v_1$ связью. Очевидно, что МП $G \cup v$ остается точечным (Рис. 20). Предположим, мы приклеили точечные связи $(vv_1)$, $(vv_2)$, $(vv_3)$,.... $(vv_{k-1})$. Так как $G \cup v$ остается точечным, то согласно аксиоме в G имеется точка $v_k$ такая, что $O(vv_k) \cap G$ является точечным подпространством, что означает, что v и $v_k$ могут быть соединены связью. Таким путем мы можем приклеить все связи между v и G. Доказательство закончено.□

## *ГОМОТОПНЫЕ ПРОСТРАНСТВА И ИХ СВОЙСТВА*

Как уже говорилось, два молекулярных пространства G и H называются гомотопически эквивалентными или гомотопными, если одно из них может быть получено из другого путем точечных преобразований.

Теорема 11

*Пусть G и v являются точечным МП и его произвольной точкой. Тогда G-v и O(v) гомотопны, и O(v) может быть преобразован точечным приклеиванием точек в G-v (Рис. 21).*

Доказательство.

В соответствии с аксиомой произвольная точка v точечного МП G может быть соединена со всеми несмежными ей точками последовательным приклеиванием точечных связей. Приклеивание точечной связи в G (Рис. 21) равносильно точечному приклеиванию точки к O(v). Следовательно, G-v гомотопно O(v) и может быть получено из O(v) точечным приклеиванием к последней точек. Доказательство закончено.□

Рассмотрим связь пространства и его подпространства для неточечных пространств, используя результаты, полученные выше.

Определение вложенного гомотопного пространства

Подпространство G пространства H называется вложенным и гомотопным неточечному пространству H, если H может быть стянуто в G точечными преобразованиями, не затрагивающими G.

(При этом очевидно, что G гомотопно H, G~H).

Это определение означает, что производятся точечные преобразования, при которых как точки подпространства G так и связи между ними не меняются. Точки могут отбрасываться только из H-G, связи могут отбрасываться или устанавливаться только между точками, не принадлежащими одновременно G. Вообще говоря, можно предположить, что МП H стягивается в G преобразованиями, вовлекающими точки и связи из G, то есть меняющими G в процессе стягивания. Так это или нет

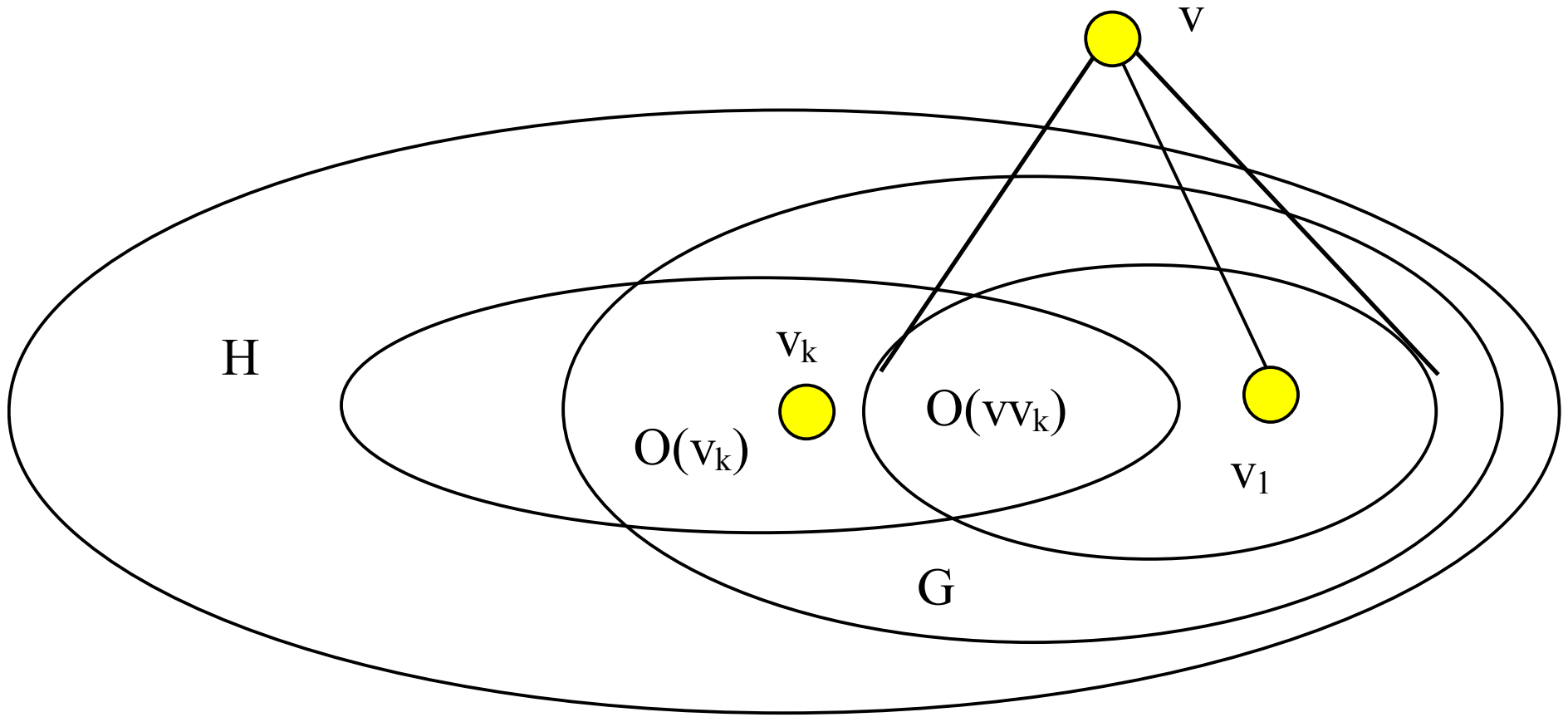

*Рис. 20* Точка v соединяется связью с точкой $v_1$ подпространства $G_1$. Затем устанавливаются связи между v и остальными точками подпространства $G_1$.

пока неясно. Сформулируем эту проблему в виде задачи.

З а д а ч а

Доказать или опровергнуть следующее:

Пусть неточечное подпространство G пространства H гомотопно H. Тогда G вложено и гомотопно H.

Теорема 12

*Пусть неточечное МП G вложено и гомотопно H. Тогда МП, полученное присоединением точки v к H по G, O(v)=G, является точечным пространством.*

Доказательство

Согласно определению подпространста G, вложенного в H, G может быть получено из H путем точечных преобразований, не затрагивающих G. Пометив эти преобразования, соединим точку v со всеми точками из H. Затем отбрасываем связь точки v в той же последовательности с теми точками из H, которые отбрасывались при преобразовании H в G (Рис. 22). Очевидно, что отбрасывания этих связей также будут точечными

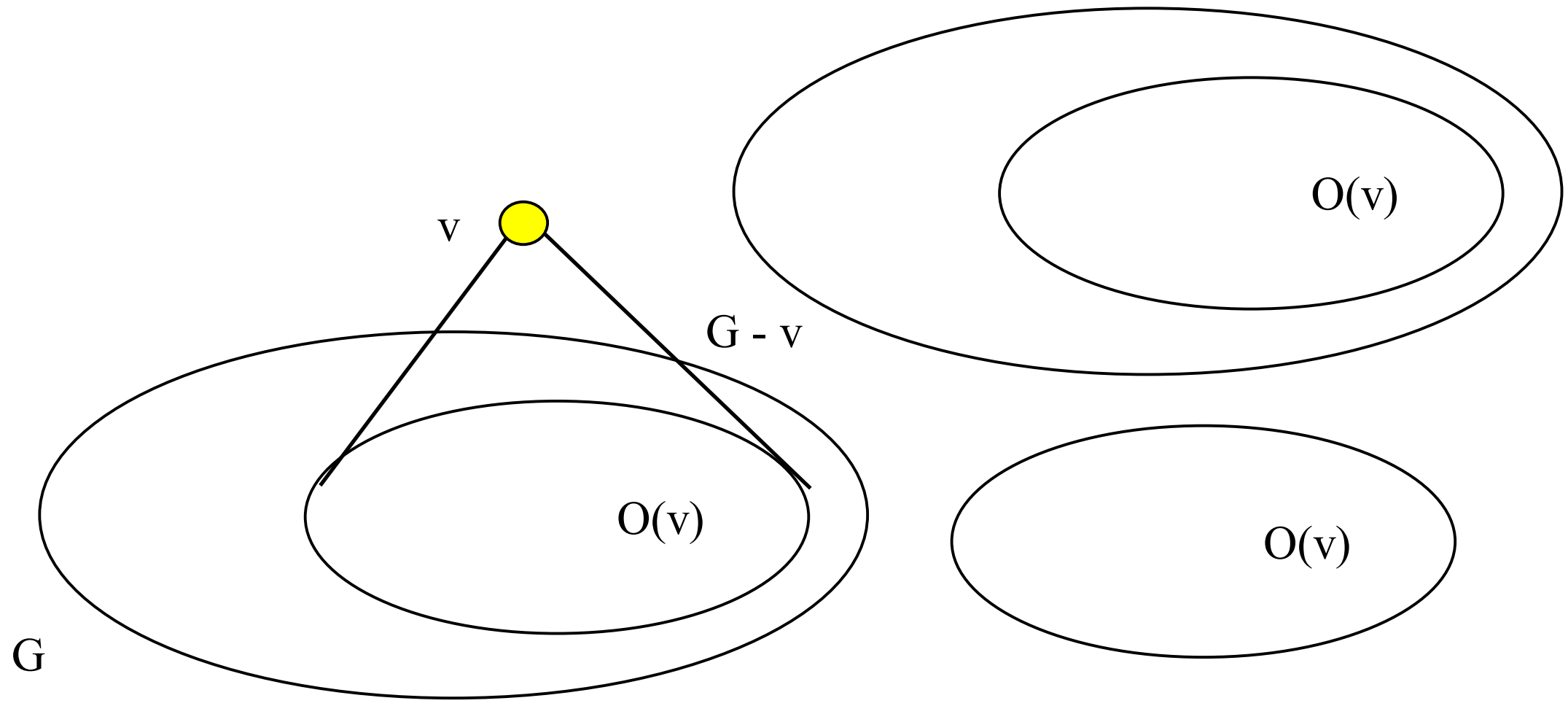

*Рис. 21* Если МП G является точечным, то окаем O(v) гомотопен подпространству G-v.

преобразованиями. Следовательно, точечное пространство H⊕v стягивается в H∪(G⊕v), которое также будет точечным. Теорема доказана.□

Следует обратить внимание на то обстоятельство, что, в отличие от точечного пространства, не накладывается никаких ограничений на связность МП H и его вложенного и гомотопного подпространства G.

Теорема 13

*Пусть подпространство G вложено и гомотопно H.. Если H не содержит точек с точечными окаемами, то H=G..*

Доказательство.

Предположим противное: H≠G. Соединим точку v со всеми точками из G в H. Полученное пространство H∪(G⊕v) является точечным по предыдущей теореме. Следовательно, для точки v должна существовать несмежная точка с точечным окаемом. Так как такой точки нет по условию, то наше предположение неверно. Это означает, что H=G. Доказательство закончено.□

Аналогичная теорема также будет доказана для нормальных молекулярных пространств.

Следующая теорема позволяет уточнить процедуру преобразования неточечного молекулярного пространства к его вложенному гомотопному подпространству. Оказывается, что последовательность преобразований можно свести только к отбрасываниям точечных точек, не привлекая установление или отбрасывание связей между точками.

Теорема 14

*Пусть H и G являются неточечными МП и его вложенным гомотопным подпространством и H≠G. Тогда H может быть преобразовано в G последовательным отбрасыванием точечных точек, не принадлежащих G.*

Доказательство.

Приклеим точку v к H таким образом, что $O(v_1) = G$ (Рис. 23). Полученное пространство $H\cup(G\oplus v)$ является точечным. Согласно аксиоме, существует несмежная с v точка $v_1$ такая, что $O(vv_1) = O(v_1)\cap G$ есть точечное МП. Установление связи в $H\cup(G\oplus v)$ равносильно приклеиванию точки $v_1$ к МП G по $O(v_1)\cap G_1$. Очевидно, точечное приклеивание связей, инцидентных точке v можно продолжить, пока точка v не будет смежной всем остальным точкам МП H. Этот процесс равносилен приклеиванию точек к МП G. При этом G переходит в H. Отклеивание точек от H в обратной последовательности переводит H в G. Доказательство закончено. □

На Рис. 24 пространство H является неточечным пространством, содержащим нормальную окружность G, состоящую из точек $v_1$, $v_2$, $v_3$, и $v_4$, и являющуюся вложенным гомотопным подпространством. В соответствии с теоремой, H стягивается в G только отбрасыванием ненумерованных точек, без привлечения связей. Последние теоремы

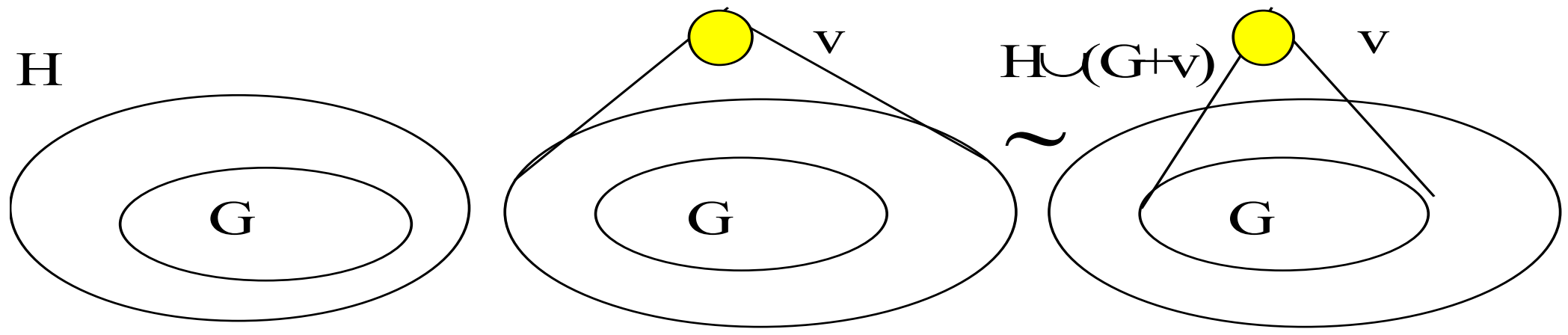

*Рис. 22* Пространство G вложено в H.. Присоединение точки v, O(v)=G, превращает H в точечное пространство.

представляют принципиальный интерес. Смысл этих теорем в том, что МП может быть получено из любого другого гомотопного ему одним лишь приклеиванием или отбрасыванием точек, не задействуя приклеивание или отбрасывание связей. В связи с этим, нам хотелось бы сформулировать задачу, которая еще не решена, но кажется очевидной.

Задача.

Если два МП G и H с одинаковым количеством точек гомотопны, то одно из них может быть стянуто в другое последовательностью

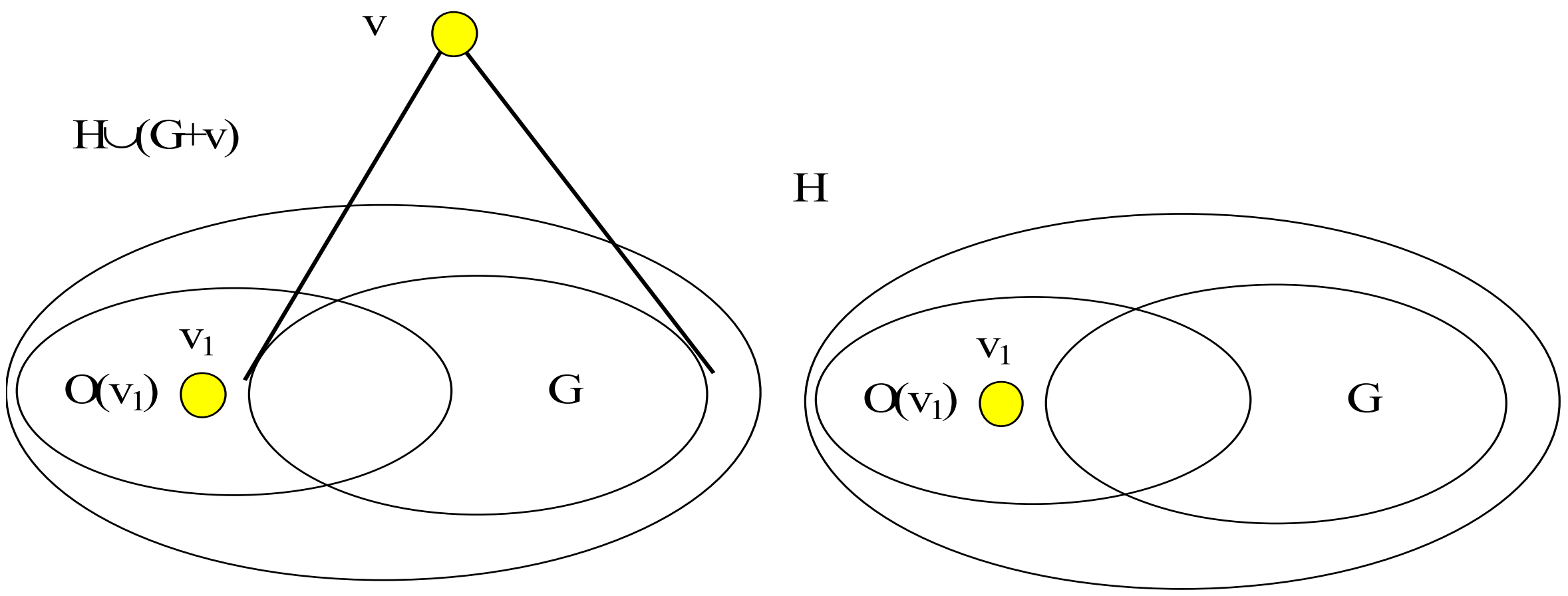

*Рис. 23* В неточечном МП H вложенное гомотопное подпространство G может быть из H отбрасыванием точек с точечными окаемами.

точечных отбрасываний и приклеиваний связей. (сравните с теоремой о точках).

## ПЕРЕСТРОЙКИ МОЛЕКУЛЯРНЫХ ПРОСТРАНСТВ

Мы имеем четыре точечных преобразования молекулярного пространства: приклеивание и отбрасывание точки по точечному окаему и приклеивание и отбрасывание связи между двумя точками, если общий окаем этих точек есть точечное подпространство. Пусть пространство G, состоящее из точек $(v_1,v_2,v_3,...v_d)$, приводится к пространству H последовательностью точечных преобразований $\Phi_1,\Phi_2,\Phi_3,...\Phi_n$. Мы хотим построить некоторое новое пространство, содержащее G и H как подпространства, и принадлежащее тому же классу пространств, что и G и H, то есть гомотопное каждому из них. Введем прямое произведение $G\otimes K(2)= G\otimes K(u_0,u_1)$, которое гомотопно G, $G\otimes K(2)\sim G$. Очевидно, что оно может быть представлено в виде двух копий пространства G. Обозначим $G_0=G\otimes u_0$ и $G_1=G\otimes u_1$. Все точки $v_k\otimes u_0$, подпространства $G_0$ имеют точечные окаемы вида $O(v_k\otimes u_0)= (v_k\otimes u_1\oplus O(v_k)\otimes u_0)\cup O(v_k)\otimes u_1$. Точно также все точки $v_k\otimes u_1$, подпространства $G_1$ (а также $G_0$) имеют точечные окаемы вида $O(v_k\otimes u_1)=(v_k\otimes u_0\oplus O(v_k)\otimes u_1)\cup O(v_k)\otimes u_0$. Применим к $G_1$ точечное преобразование Φ. Это преобразование будет точечным также в $G\otimes K(2)$.Следует подчеркнуть, что $\Phi(G_1)(G\otimes K(2))\neq(\Phi_1 G)\otimes K(2)$. Получается пространство, которое уже не является прямым произведением, поскольку над ним уже совершено преобразование Φ. Оно гомотопно G, $\Phi(G_1)(G\otimes K(2))\sim G\otimes K(2)\sim G$, и состоит из двух слоев: один слой $G_0$ есть пространство G, а другой слой $G_1$ является пространством $G_1=\Phi_1 G$.

Теорема 15

*Для любого пространства G прямое произведение G⊗K(2), где K(2) является связкой из двух точек, гомотопно G (Рис. 26).*

Доказательство.

Справедливость теоремы следует непосредственно из свойств прямого произведения. □

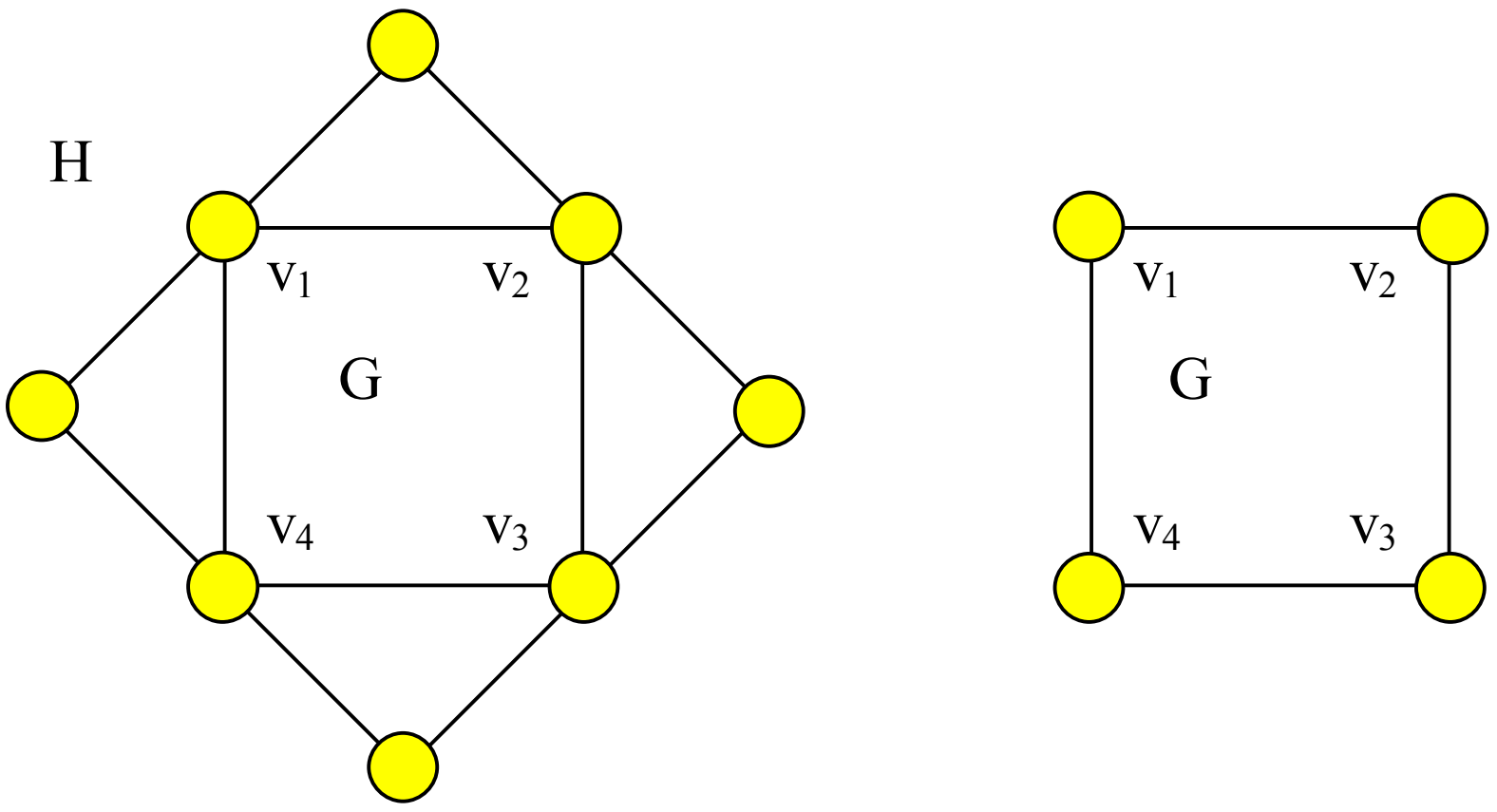

*Рис. 24* Пространство H является неточечным пространством, содержащим вложенное гомотопное пространство G. H стягивается в G отбрасыванием четырех точечных точек

Введем определение вышеописанной важной конструкции.

<u>Определение индуцированного точечного преобразования на G⊗K(2)</u>

Пусть Φ есть точечное отбрасывание либо точки $v_k$, либо связи $(v_kv_p)$ в $G=(v_1,v_2,v_3,...v_d)$, вовлекающее либо $O(v_k)$, либо $O(v_kv_p)$ соответственно. Тогда точечное отбрасывание Φ индуцирует отбрасывание $\Phi(G_1)$ в G⊗K(2) либо точки $v_k \otimes u_1$, либо связи $((v_k \otimes u_1)(v_p \otimes u_1))$ соответственно. Пусть Φ есть точечное приклеивание либо точки v, либо связи $(v_kv_p)$ в $G=(v_1,v_2,v_3,...v_d)$ по точечным подпространствам либо O(v), либо $O(v_kv_p)$ соответственно. Тогда точечное приклеивание Φ индуцирует приклеивание $\Phi(G_1)$ в G⊗K(2) либо точки $v \otimes u_1$, либо связи $((v_k \otimes u_1)(v_p \otimes u_1))$ по подпространствам либо O(v)⊗K(2), либо $O(v_kv_p) \otimes K(2)$ соответственно (Рис. 26, Рис. 25). Преобразования $\Phi(G_1)$ называются индуцированными точечными преобразованиями на G⊗K(2).

Теорема 16

*Индуцированное преобразование $\Phi(G_1)$ над G⊗K(2) всегда является точечным.*

Доказательство.

По сути дела это уже следует из предыдущего рассмотрения. Пусть $O(v_k)$ в G является точечным окаемом, и точку $v_k$ можно отбросить. Тогда $O(v_k \otimes u_1)=v_k \otimes u_0 \oplus O(v_k) \otimes K(2)$.

Так как точка $v_k \otimes u_0$ смежна со всеми остальными точками из $O(v_k \otimes u_1)$, то $O(v_k \otimes u_1)$ является конусом (вообще говоря, независимо от вида $O(v_k)$), и,

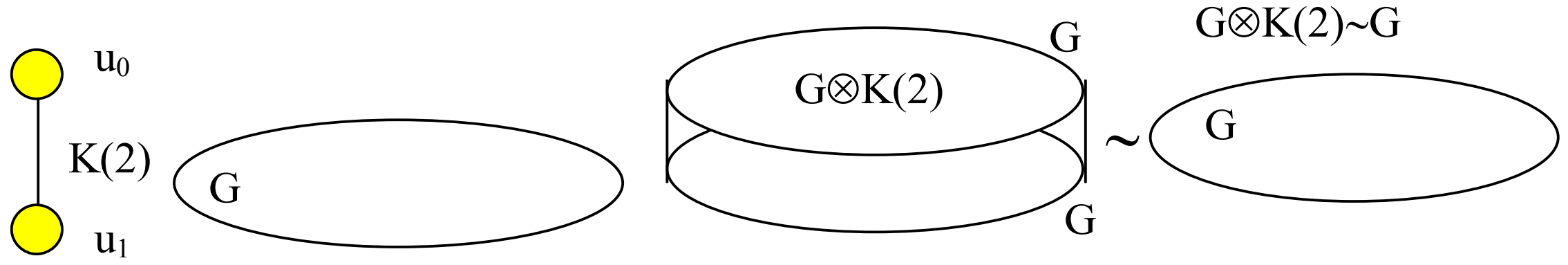

*Рис. 26* Прямое произведение пространства G и двухточечной связки K(2) содержит G как вложенное подпространство и гомотопно G.

следовательно, точечным пространством (Рис. 27). Пусть $O(v_k v_p)$ в G является точечным общим окаемом, и связь $(v_k v_p)$ можно отбросить. Тогда $O((v_k \otimes u_1)(v_p \otimes u_1))=O(v_k \otimes u_1) \cap O(v_p \otimes u_1)=$ $[(v_k \otimes u_0 \oplus O(v_k) \otimes u_1) \cup O(v_k) \otimes u_0] \cap [(v_p \otimes u_0 \oplus O(v_p) \otimes u_1) \cup O(v_p) \otimes u_0]=v_k \otimes u_0 \oplus$ $v_p \otimes u_0 \oplus O(v_k v_p) \otimes K(2)$. Так как точки $v_k \otimes u_0$ и $v_p \otimes u_0$ смежны со всеми остальными точками из $O((v_k \otimes u_1)(v_p \otimes u_1))$, то $O((v_k \otimes u_1)(v_p \otimes u_1))$ является конусом (вообще говоря независимо от вида $O(v_k v_p)$) и, следовательно, точечным пространством.

Таким образом, при отбрасывании точки или связи преобразование $\Phi(G_1)$ всегда точечное. Пусть Φ есть точечное приклеивание точки v в G по O(v). Рассмотрим $O(v \otimes u_1)$ в $G \otimes K(2)$. По определению $O(v \otimes u_1)= O(v) \otimes K(2)$ есть

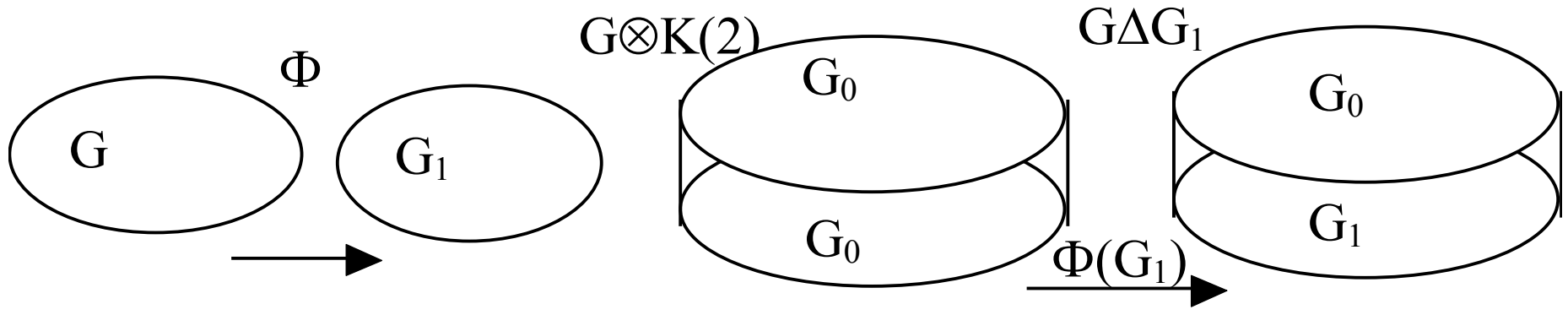

*Рис. 25* Точечное преобразование Φ пространства G в пространство $G_1$ путем приклеивания точки по точечному окаему. Оно индуцирует точечное преобразование $\Phi(G_1)$ пространства $G \otimes K(2)$ в пространство $G \Delta G_1$. Приклеенная точка v отнесена к $G_1$.

прямое произведение двух точечных пространств и, следовательно, является точечным. Пусть теперь Φ есть точечное приклеивание связи $(v_k v_p)$ в G по O $(v_k v_p)$. Рассмотрим $O((v_k \otimes u_1)(v_p \otimes u_1))$ в $G \otimes K(2)$. Очевидно, $O((v_k \otimes u_1)(v_p \otimes u_1))= O(v_k v_p) \otimes K(2)$ есть прямое произведение двух точечных пространств и, следовательно, является точечным. Таким образом, индуцированные приклеивания $\Phi(G_1)$ являются точечными. Теорема доказана. □

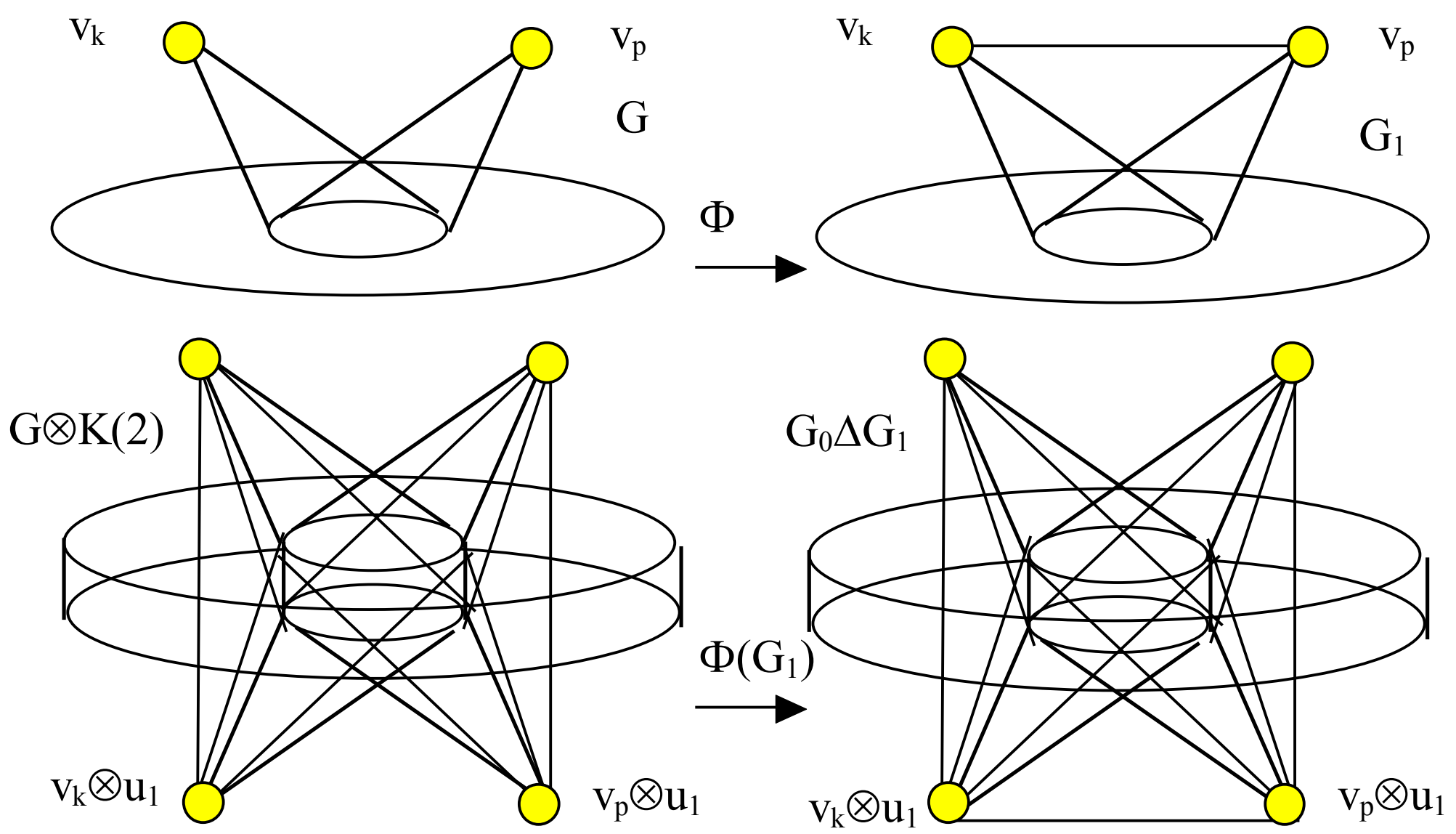

*Рис. 27* Точечное преобразование Φ пространства G в пространство $G_1$ путем приклеивания связи между точками по точечному окаему. Оно индуцирует точечное приклеивание связи $\Phi(G_1)$ между точками $v_k\otimes u_1$ и $v_p\otimes u_1$, то есть точечное преобразование пространства $G\otimes K(2)$ в пространство $G\Delta G_1$. Установленная связь $((v_k\otimes u_1)(v_p\otimes u_1))$отнесена к $G_1$.

Определение одношаговой Φ перестройки или башни над пространством G.

Одношаговой Φ перестройкой (Рис. 28) или башней над пространством $G=(v_1,v_2,v_3,...v_d)$, называется пространство $G_0\Delta G_1$, полученное из $G\otimes K(2)=G\otimes K(u_0,u_1)$, к которому применено индуцированное точечное преобразование $\Phi(G_1)$, $G_0\Delta G_1=\Phi(G_1)(G\otimes K(2))$. При этом $G_0=G$ и состоит из точек $G\otimes u_0$, $G_1\approx\Phi G$.

Определим свойства одношаговой перестройки в виде теоремы.

Теорема 17

*Одношаговая перестройка $G_0\Delta G_1$ над пространством G есть пространство, обладающее следующими свойствами:*

*$G_0\Delta G_1$ гомотопно как $G_0$, так и $G_1$, в котором $G_0= G\otimes u_0\approx G$, а $G_1=\Phi G$.*

*$G_0\Delta G_1$ может быть стянуто точечными отбрасываниями вершин либо к $G_0$, либо к $G_1$.*

*$G_0\Delta G_1$ можно рассматривать как одношаговую перестройку над $G_1$.*

Доказательство.

Докажем, что $G_0\Delta G_1$ гомотопно как $G_0$, так и $G_1$, где $G_0= G\otimes u_0\approx G$, а $G_1=\Phi G$. Так как $G_0\Delta G_1\sim G_0\otimes K(2)$, и $G_0\otimes K(2)\sim G$, то $G_0\Delta G_1\sim G\sim G_1=\Phi G$. Докажем, что $G_0\Delta G_1$ может быть стянуто точечными отбрасываниями вершин либо к $G_0$.

Предположим, что $\Phi(G_1)$ в $G_0\otimes K(2)$ есть отбрасывание точки $v_k\otimes u_1$ с точечным окаемом $O(v_k\otimes u_1)=v_k\otimes u_0\oplus O(v_k)\otimes K(2)$, где $O(v_k)$ является точечным пространством. Окаем любой точки вида $v_p\otimes u_1$ определяется выражением $O(v_p\otimes u_1)=v_p\otimes u_0\oplus(O(v_p)\otimes K(2))$, если $v_k\otimes u_1$ и $v_p\otimes u_1$ не смежны, $O(v_p\otimes u_1)=v_p\otimes u_0\oplus(O(v_p)\otimes K(2)-(v_k\otimes u_1))$, если $v_k\otimes u_1$ и $v_p\otimes u_1$

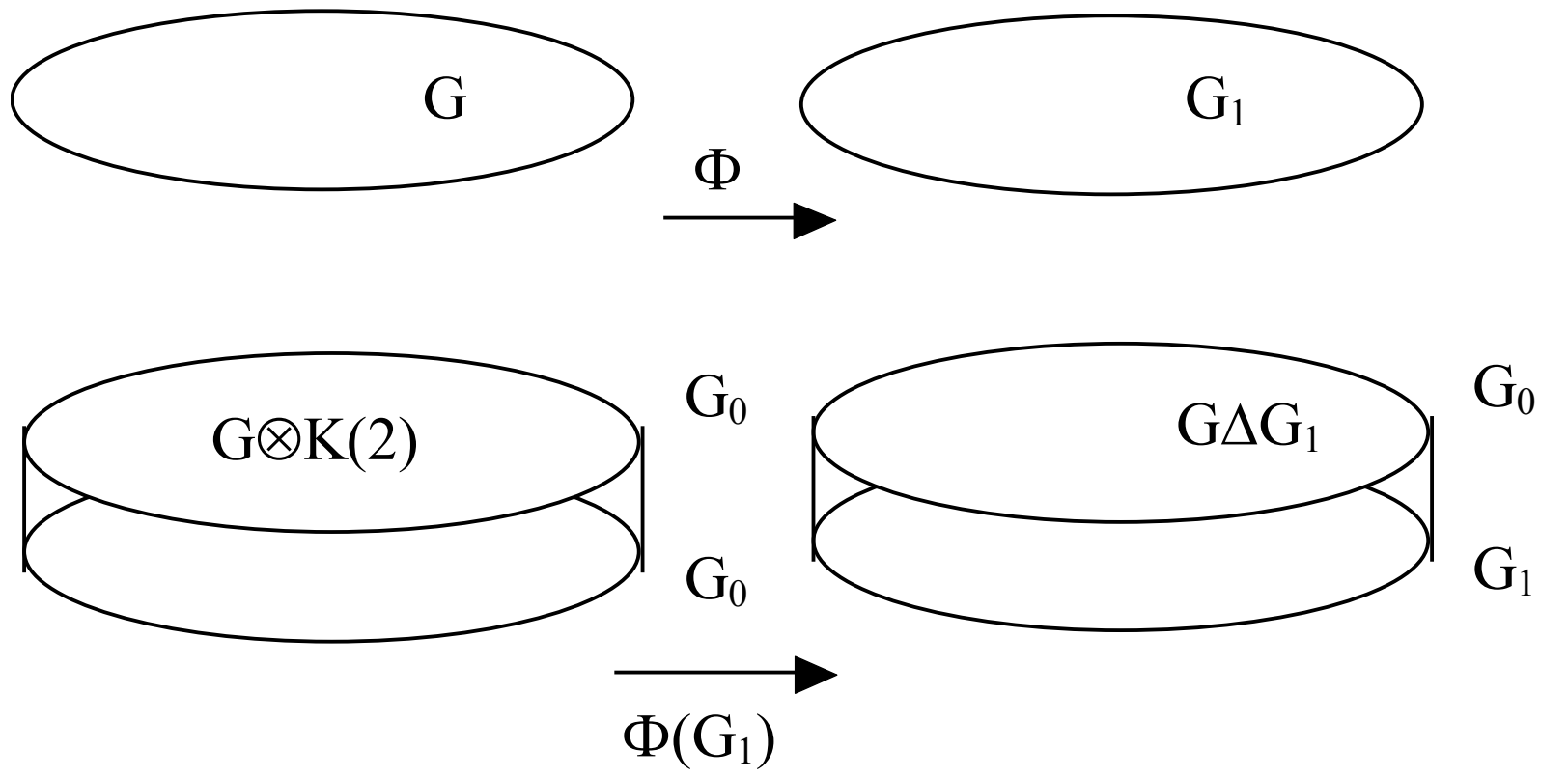

*Рис. 28* Точечное преобразование Φ пространства G в пространство $G_1$ индуцирует точечное преобразование $\Phi(G_1)$ пространства $G\otimes K(2)$ в пространство $G\Delta G_1$. $G\Delta G_1$. называется одношаговой Φ перестройкой пространства G.

смежны, и является конусом, то есть точечным пространством. Следовательно, все точки из $G_1=(G\otimes u_1-(v_k\otimes u_1))$ могут быть отброшены, и $G_0\Delta G_1$ может быть стянуто в $G_0$ точечными отбрасываниями точек.

Предположим, что $\Phi(G_1)$ в $G_0\otimes K(2)$ есть приклеивание точки $v\otimes u_1$ по точечному подпространству $O(v\otimes u_1)= O(v)\otimes K(2)$, где $O(v)$ является точечным пространством. Тогда окаем этой же точки $v\otimes u_1$ в $G_0\Delta G_1$ является точечным по построению, и эта точка может быть отброшена. Окаем любой другой точки $v_k\otimes u_1$ из $G_0\Delta G_1$ определяется выражением $O(v_k\otimes u_1)=v_k\otimes u_0\oplus O(v_k)\otimes K(2)$. Следовательно, $O(v_k\otimes u_1)=v_k\otimes u_0\oplus O(v_k)\otimes K(2)$ есть конус, то есть точечное пространство. Это означает, что все остальные точки из $G_1$ могут быть отброшены, и $G_0\Delta G_1$ может быть стянуто в $G_0$ этими точечными отбрасываниями точек.

Предположим, что $\Phi(G_1)$ в $G_0\otimes K(2)$ есть отбрасывание связи $((v_k\otimes u_1)(v_p\otimes u_1))$ между точками $(v_k\otimes u_1)$ и $(v_p\otimes u_1)$ по точечному подпространству $O((v_k \otimes u_1)(v_p \otimes u_1))=v_k \otimes u_0 \oplus v_p \otimes u_0 \oplus O(v_k v_p) \otimes K(2)$, где $O(v_k v_p)$ является точечным пространством.

Очевидно, что окаем любой точки $v_s \otimes u_1$, есть конус вида $O(v_s \otimes u_1)=v_s \otimes u_0 \oplus (O(v_s) \otimes K(2)-A)$, то есть точечное пространство (здесь вид

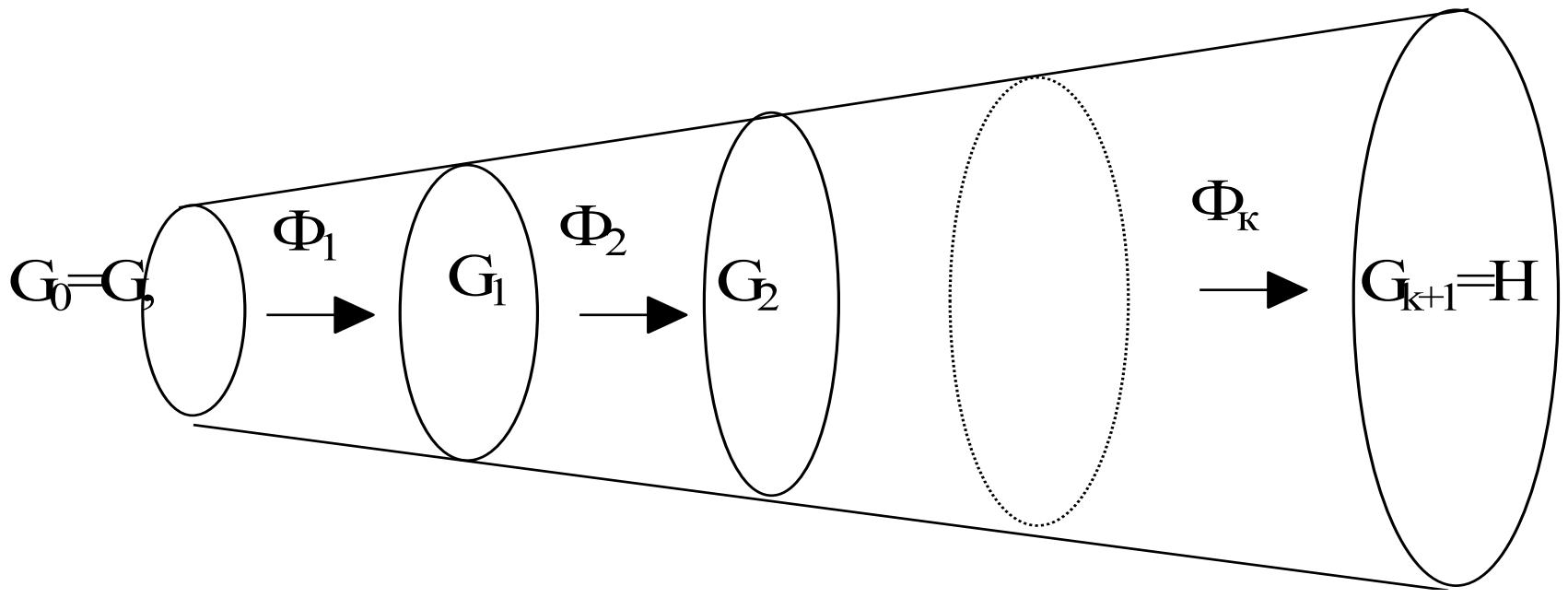

*Рис. 29* n-шаговая $\Phi_1,\Phi_2,\Phi_3,...\Phi_n$ перестройка или башня над пространством G есть пространство $W=G_0 \Delta G_1 \Delta G_2 \Delta G_3 \Delta ... G_n$, где $G_k \Delta G_{k+1}$, k=0,1,2,...n, являются одношаговыми перестройками, $G_{k+1}=\Phi_k G_k$, $G_0=G$, $G_{k+1}=H$. W гомотопно G и H.

A не имеет значения). Это означает, что все точки из $G_1$ могут быть отброшены, и $G_0 \Delta G_1$ может быть стянуто в $G_0$ этими точечными отбрасываниями точек. Рассмотрим две несмежные точки $(v_k \otimes u_1)$ и $(v_p \otimes u_1)$ в $G \otimes K(2)$, чей точечный общий окаем определяется выражением $O((v_k \otimes u_1)(v_p \otimes u_1))= O(v_k v_p) \otimes K(2)$, где $O(v_k v_p)$ является точечным пространством. Установим связь $((v_k \otimes u_1)(v_p \otimes u_1))$. При этом $G \otimes K(2)$ перейдет в $G_0 \Delta G_1$. Рассмотрим теперь точку $v_k \otimes u_1$ в $G_0 \Delta G_1$. $O(v_k \otimes u_1)=v_k \otimes u_0 \oplus O(v_k) \otimes K(2) \cup (O(v_k v_p) \otimes K(2) \oplus v_p \otimes u_1)$. Вид этого окаема означает, что к точечному пространству $v_k \otimes u_0 \oplus O(v_k) \otimes K(2)$ приклеена точка $v_p \otimes u_1$ по точечному подпространству $(O(v_k v_p) \otimes K(2) \oplus v_p \otimes u_1)$. Следовательно, $O(v_k \otimes u_1)$ является точечным, и точка $v_k \otimes u_1$ может быть отброшена То же самое рассуждение справедливо для точки $v_p \otimes u_1$. Рассмотрим точку $v_s \otimes u_1$, смежную точкам $v_k \otimes u_1$ и $v_p \otimes u_1$. Легко видеть, что $O(v_s \otimes u_1)=v_s \otimes u_0 \oplus (O(v_s) \otimes K(2) \cup (v_k v_p))$ является точечным окаемом, и точка $v_s \otimes u_1$ может быть отброшена Окаемы всех остальных точек $v_r \otimes u_1$ при установлении связи $((v_k \otimes u_1)(v_p \otimes u_1))$ не меняются и определяются выражением $O(v_r \otimes u_1)=v_r \otimes u_0 \oplus O(v_r) \otimes K(2)$. Так как они являются конусами, то есть точечными пространствами, то все эти точки могут быть отброшены. Это означает, что все точки из $G_1$ могут быть отброшены, и $G_0 \Delta G_1$ может быть стянуто в $G_0$ этими точечными отбрасываниями точек. Таким образом, мы доказали, что $G_0 \Delta G_1$ может быть стянуто в $G_0$ точечными отбрасываниями точек из $G_1$. Рассмотрим теперь то же самое пространство, но представленное в виде $G_1 \Delta G_0$, что достигается простой заменой 0 на 1 и наоборот в обозначении точек $v_s \otimes u_0$ и $v_s \otimes u_1$. При этом

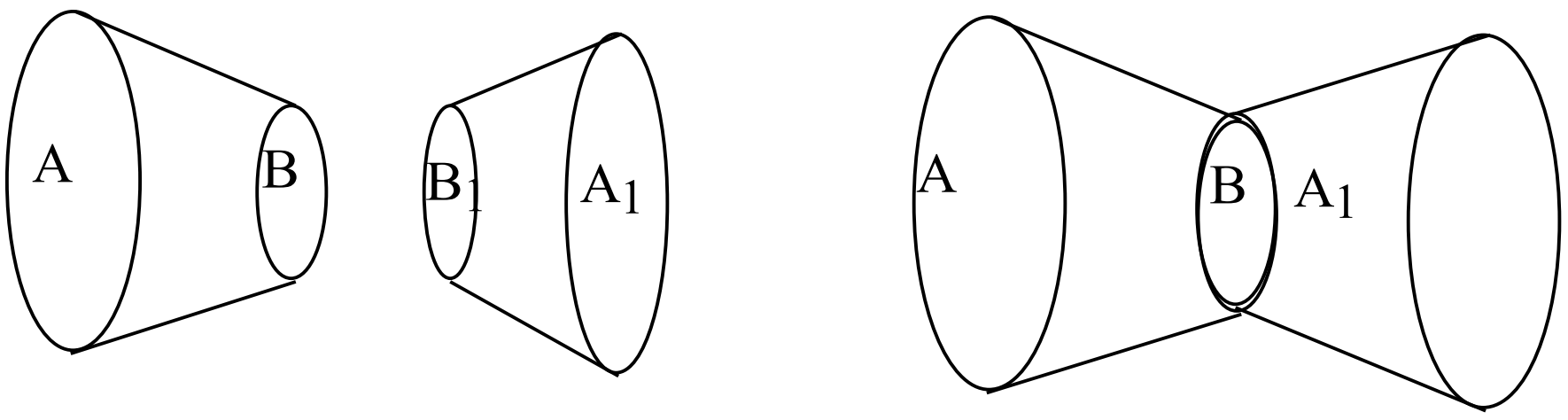

*Рис. 30* Склейка двух пространств по подпространству B, изоморфному $B_1$.

$G_0 = \Phi^{-1} G_1$. Легко видеть, что при такой форме записи предыдущее доказательство означает, что в $G_1 \Delta G_0$ все точки из $G_0$, могут быть отброшены, и $G_1 \Delta G_0$ может быть стянуто в $G_1$ этими точечными отбрасываниями. Это рассуждение доказывает также, что $G_0 \Delta G_1$ можно рассматривать как одношаговую перестройку над $G_1$. Теорема доказана.□

Предположим, что пространство G, приводится к пространству H последовательностью точечных преобразований $\Phi_1, \Phi_2, \Phi_3, \ldots \Phi_n$. На первом этапе строим одношаговую $\Phi_1$ перестройку $G_0 \Delta G_1$, где $G_0 = G$, $G_1 \approx \Phi_1 G$. На следующем этапе строим одношаговую $\Phi_2$ перестройку $G_1 \Delta G_2$ над пространством $G_1$, в которой $G_2 = \Phi_2 G_1 = \Phi_2 \Phi_1 G$. Склеивая две одношаговые перестройки $G_0 \Delta G_1$ и $G_1 \Delta G_2$ по $G_1$ получаем двухшаговую $\Phi_1 \Phi_2$ перестройку $G_0 \Delta G_1 \Delta G_2$. Продолжая этот процесс n раз, получаем пространство, на одном конце которого находится $G_0 = G$, а на другом конце пространство $G_n = H$.

<u>Определение n-шаговой перестройки или башни над пространством G</u>

n-шаговой $\Phi_1, \Phi_2, \Phi_3, \ldots \Phi_n$ перестройкой или башней над пространством G называется пространство $W = G_0 \Delta G_1 \Delta G_2 \Delta G_3 \Delta \ldots G_n$, где $G_k \Delta G_{k+1}$, $k = 0,1,2,\ldots n$, является одношаговой перестройкой, $G_{k+1} = \Phi_k G_k$, $G_0 = G$, $G_{k+1} = H$ (Рис. 29).

Из построения легко следуют свойства этой конструкции.

Теорема 18

*Пусть пространство G, стягивается к пространству H последовательностью точечных преобразований* $\Phi_1, \Phi_2, \Phi_3, \ldots \Phi_n$.

*Тогда существует пространство $W=G_0\Delta G_1\Delta G_2\Delta G_3\Delta ...G_n$, где $G_0=G$, $G_n=H$, обладающее следующими свойствами:*

- *W гомотопно G и H, W~G~H.*
- *W содержит G и H как вложенные гомотопные подпространства, $G\subseteq W$, $H\subseteq W$.*
- *Точечными отбрасываниями вершин W стягивается как к G, так и к H.*

Д о к а з а т е л ь с т в о.

Построим n-шаговую башню $W=G_0\Delta G_1\Delta G_2\Delta G_3\Delta ...G_n$, над пространством с точечными преобразованиями $\Phi_1,\Phi_2,\Phi_3,...\Phi_n$, где $G_{k+1}=\Phi_k G_k$, $G_0=G$, $G_{k+1}=H$. В соответствии с предыдущим, любую точку из слоя $G_n=H$ можно отбросить точечным образом, поскольку она принадлежит одноточечной перестройке $G_{n-1}\Delta G_n$. Отбрасывая послойно подпространства $G_n$, $G_{n-1}$, $G_{n-2}\Delta,...G_1$, мы приводим W к $G_0=G$. Начав тот же самый процесс отбрасывания точек с подпространства $G_0=G$, мы приводим W к $G_n=H$. Теорема доказана.□

Цель всех этих построений в том, что для любого конечного набора гомотопных пространств существует некоторое универсальное пространство, обладающее теми же свойствами, что и данные пространства, содержащее их как подпространства. Прежде всего напомним определение склейки двух пространств.

Определение склейки двух пространств

Пусть $G=A\cup B$ и $H=A_1\cup B_1$ есть два молекулярных пространства и подпространство B изоморфно $B_1$, h: $B\rightarrow B_1$. Тогда склейкой G и H по изоморфизму h называется пространство $F=G\Diamond H=(G\cup H)/B=A\cup B\cup A_1=(G,B)\#(B,H)$, содержащее точки из A, B и $A_1$, в котором любая точка $v\in B$ отождествлена с $h(v)\in B_1$.

Теорема 19

*Пусть имеется конечный набор гомотопных пространств $G_0,G_1,G_2,...G_n$,, Тогда существует пространство $W_0$, обладающее следующими свойствами:*

- *W0 гомотопно Gk, k=0,1,2,...n*
- *W0 содержит Gk, k=0,1,2,...n, как вложенные гомотопные подпространства, $G_k \subseteq W$.*
- *Точечными отбрасываниями вершин $W_0$ стягивается к любому $G_k$, k=0,1,2,...n.*

Доказательство.

Для каждой соседней пары пространств $G_{k-1}$, $G_k$, k=1,2,...n строим башню $W_k=(G_{k-1},G_k)$ в соответствии с предыдущей теоремой. Затем склеиваем $W_k=(G_{k-1},G_k)$ и $W_{k+1}=(G_k,G_{k+1})$ по $G_k$ (отождествляем точки, принадлежащие разным башням). Обозначая # операцию склейки получаем пространство $W_0=W_1\#W_2\#W_3\#...W_n=(G_0,G_1)\#(G_1,G_2)\#(G_2,G_3)\#...(G_{n-1},G_n)$. Легко видеть, что это пространство удовлетворяет всем указанным выше свойствам Теорема доказана. □

Список литературы к главе 2.

# НОРМАЛЬНОЕ МОЛЕКУЛЯРНОЕ ПРОСТРАНСТВО И ЕГО РАЗМЕРНОСТЬ

In this chapter we define the dimension of a normal molecular space by means of axioms that are based on obvious geometrical background. It also presents digital models of some well known spaces such as n-dimensional spheres a torus, projective plane, Klein bottle, n-dimensional Euclidean spaces.

Некоторые результаты этой главы имеются в работах [3,7,10, 42,43,50,51]

## *ГЕОМЕТРИЧЕСКИ-ИНТУИТИВНАЯ ОСНОВА ДЛЯ ОПРЕДЕЛЕНИЯ РАЗМЕРНОСТИ НА МОЛЕКУЛЯРНОМ ПРОСТРАНСТВЕ*

Наша задача-построить молекулярное пространство с определенными свойствами, позволяющими рассматривать его как комбинаторный образ непрерывного n-мерного пространства. Для этого мы выберем в элементарной геометрии те признаки, которые описывают понятия близости, непрерывности и размерности на наиболее фундаментальном уровне, и которые позволят нам в дальнейшем легко перейти к молекулярной модели непрерывного пространства. Картина, которую мы собираемся здесь предложить читателю, не является доказательством, но позволяет понять на интуитивном уровне, что является размерностью в теории молекулярных пространств.

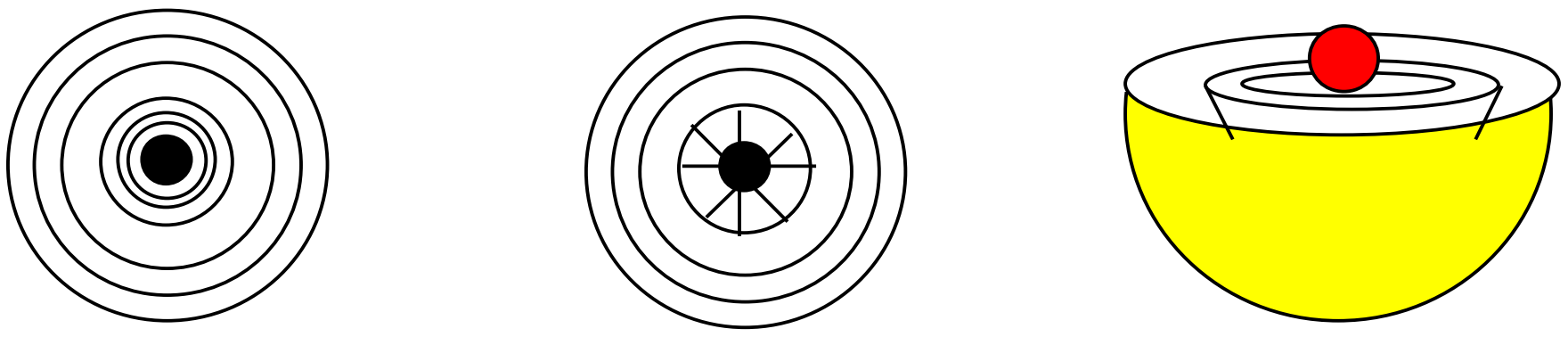

Рис. EE Слева. Непрерывный круг содержит внутри себя бесконечное множество концентрических вложенных окружностей уменьшающегося радиуса. В центре. Молекулярный круг содержит внутри себя всегда конечное число концентрических вложенных окружностей. Справа. Молекулярный шар содержит внутри себя, также как и круг, всегда конечное число концентрических вложенных сфер, окружающих центральную точку.

Рассмотрим 3-мерное (в общем случае n-мерное) евклидово пространство $E^3$, в котором выберем некоторую точку p. Пусть $D^3$ будет (трехмерным) сплошным шаром радиуса R с центром в точке p. Границей этого шара является, естественно, (двумерная) сфера $S^2$. Как известно [1,40] множество всех таких шаров со всевозможными центрами и радиусами образует базу топологии для $E^3$. Шар $D^3$ содержит внутри себя бесконечное множество концентрических вложенных шаров меньших радиусов (рис. ee).

$D^3 \supseteq D^3_1 \supseteq D^3_2 \supseteq D^3_3 \supseteq ... D^3_n \supseteq ....,$

$R > R_1 > R_2 > R_3 >..... > R_n >....,$

где $R_n$-радиус шара $D_n$.

Например, мы можем выбрать радиусы по формуле $R_n=R/(1+n)$.На рис. ee для простоты мы изображаем двумерную, а не трехмерную картинку.

Теперь предположим, что пространство молекулярно. Что изменится в той картине? Так как в молекулярном пространстве любой конечный объем содержит конечное число элементов, следовательно, шар D состоит из конечного числа точек. Это означает, что последовательность вложенных в D шаров будет всегда конечной, а не бесконечной, как в непрерывном случае.

Наглядной иллюстрацией этого случая является русская игрушка-матрешка, где внутри деревянной фигурки находится вторая меньшего размера, внутри которой находится третья и так далее до последней сплошной куклы, внутри которой уже ничего нет. Другой наглядный пример этого математического объекта-китайская игрушка, состоящая из костяных шаров, вырезанных один внутри другого, число которых может составлять несколько десятков, а порой и сотен. Можно придумать и растительные образы, например, рассматривать молекулярный шар как кочан капусты или луковицу.

Для нас здесь важно следующее: внешняя поверхность молекулярного шара является молекулярной сферой, которую можно ободрать как внешний слой кожуры с луковицы и получить молекулярный шар меньшего размера. Продолжая этот процесс обдирания, мы получим наименьший молекулярный шар, состоящий из последнего слоя, который еще можно ободрать, и неделимого элемента внутри этого слоя, играющего роль центральной точки. Таким образом, мы можем определить наименьший молекулярный шар как центральную точку, окруженную молекулярной сферой. Точки сферы являются ближайшими к центральной точке. Для выделения этой близости установим, что центральная точка соединена связью с каждой точкой ближайшей молекулярной сферы.

Мы уже почти получили определение: молекулярный минимальный трехмерный шар является точкой, окруженной молекулярной двумерной сферой, или, в общем случае, молекулярный минимальный n-мерный шар является точкой, окруженной молекулярной (n-1)-мерной сферой.

## *ОПРЕДЕЛЕНИЕ НОРМАЛЬНОГО МОЛЕКУЛЯРНОГО ПРОСТРАНСТВА И ЕГО РАЗМЕРНОСТИ*

Перейдем теперь к определению размерности на молекулярном пространстве. Размерность является важнейшей характеристикой пространства. Эволюция этого понятия прошла долгий путь, начиная с античности. Однако только в начале 20 века сформировались современные представления об этой характеристике. В математике существует несколько подходов к определению размерности, начиная с самого простого параметрического, используемого также в физике, и до топологических подходов. Однако все они неприменимы к дигитальному

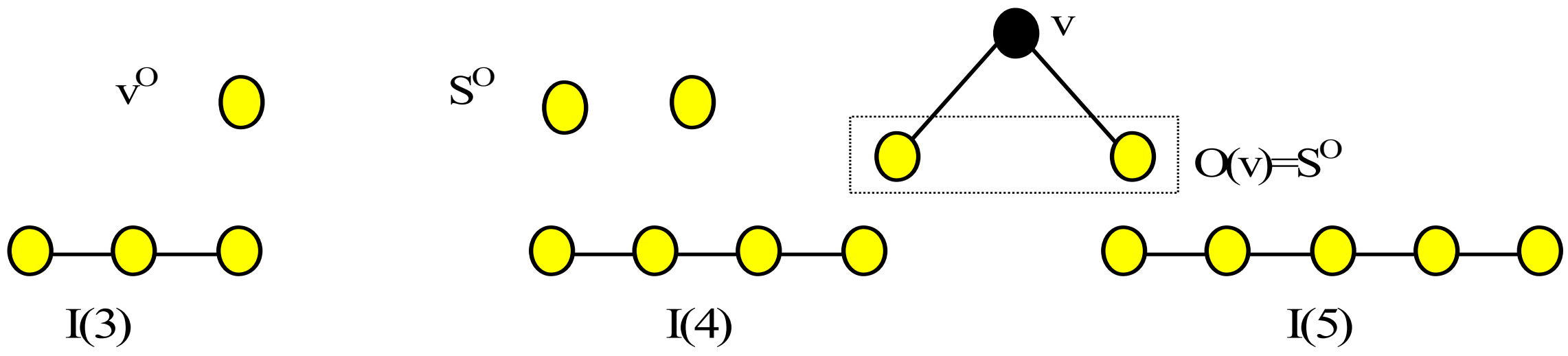

Рис. FF $S^0$ является нуль-мерной сферой. Точка v является одномерной, поскольку ее окаем есть нульмерная сфера. I(3), I(4) и I(5) есть одномерные отрезки, содержащие одну, две и три одномерные точки. Одномерные точки являются внутренними точками отрезков.

пространству, построенному на конечном множестве точек. Эта неприменимость в своей основе кроется в том, что подходы классической математики предполагают непрерывность и бесконечность, тогда как дигитальная топология, как практическая наука, вообще не использует таких понятий. Поэтому необходимо или же расширить классические понятия, сделав их применимыми к дигитальному пространству, или же заново определить размерность дигитального пространства. Мы пойдем вторым путем, определяя размерность молекулярного пространства аксиоматически с учетом компьютерных экспериментов, о которых мы расскажем ниже. Вначале перечислим 0-мерные пространства. Мы имеем два 0-мерных пространства, 0-мерную точку и 0-мерную сферу.

<u>Определение нормальной 0-мерной точки</u>

Уединенная (изолированная) точка называется нормальной 0-мерной точкой (рис. ff).

Согласно классическому определению 0-мерная сфера состоит из двух изолированных точек. Используем это стандартное определение нуль-мерной сферы в удобной для нас форме.

Определение 0-мерной нормальной сферы

Нормальной нуль-мерной сферой $S^0$ является МП, состоящее из двух изолированных точек (рис. ff).

Этими двумя пространствами ограничивается множество всех нормальных

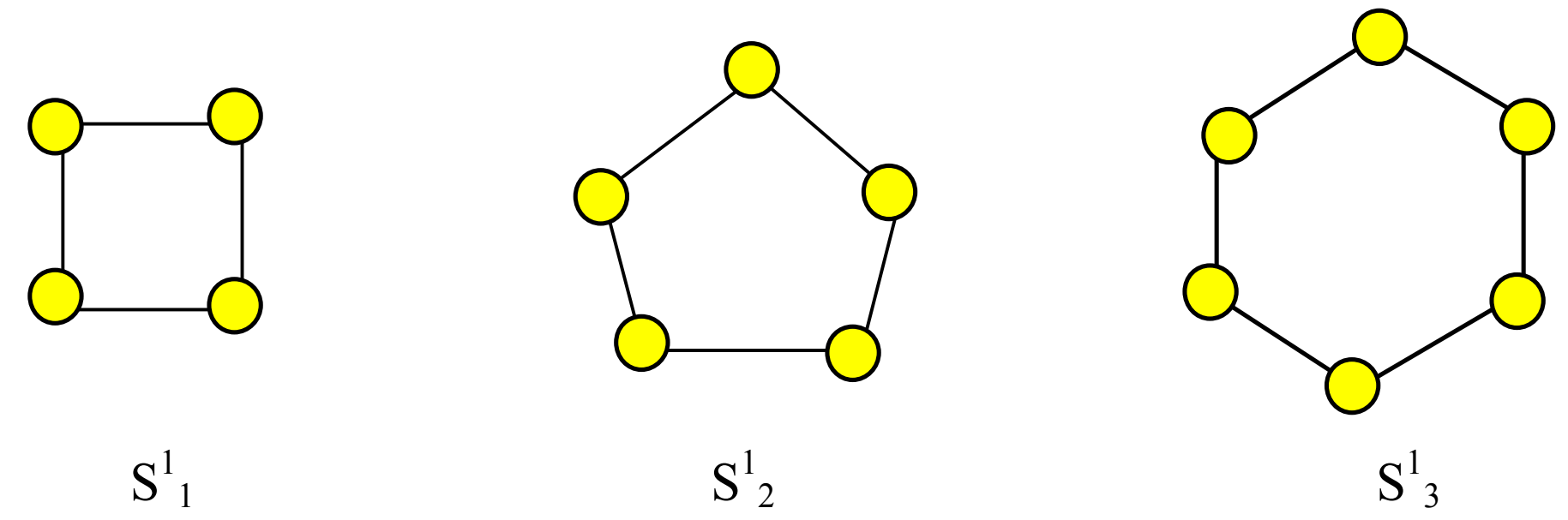

Рис. GG Изображены три одномерные нормальные сферы $S^1_1$, $S^1_2$, и $S^1_3$, состоящие из 4, 5 и 6 точек. Сфера $S^1_1$ , состоящая из 4-х точек, минимальна.

0-мерных пространств. Для сохранения последовательности определений будем считать по определению, что уединенная точка есть 0-мерное нормальное пространство с краем, а нормальная 0-мерная сфера есть единственное 0-мерное замкнутое пространство. Теперь, используя вышеописанный геометрический подход, мы можем дать определение одномерной точки, которая есть не что иное, как центральная точка одномерного минимального молекулярного шара.

Определение 1-мерной нормальной точки (рис. ff).

Точка v молекулярного пространства G называется одномерной нормальной точкой , если окаем O(v) этой точки является нормальной нуль-мерной сферой $S^0$.

Теперь мы можем определить одномерные пространства-окружности, отрезки и бесконечную прямую. В дальнейшем мы будем рассматривать

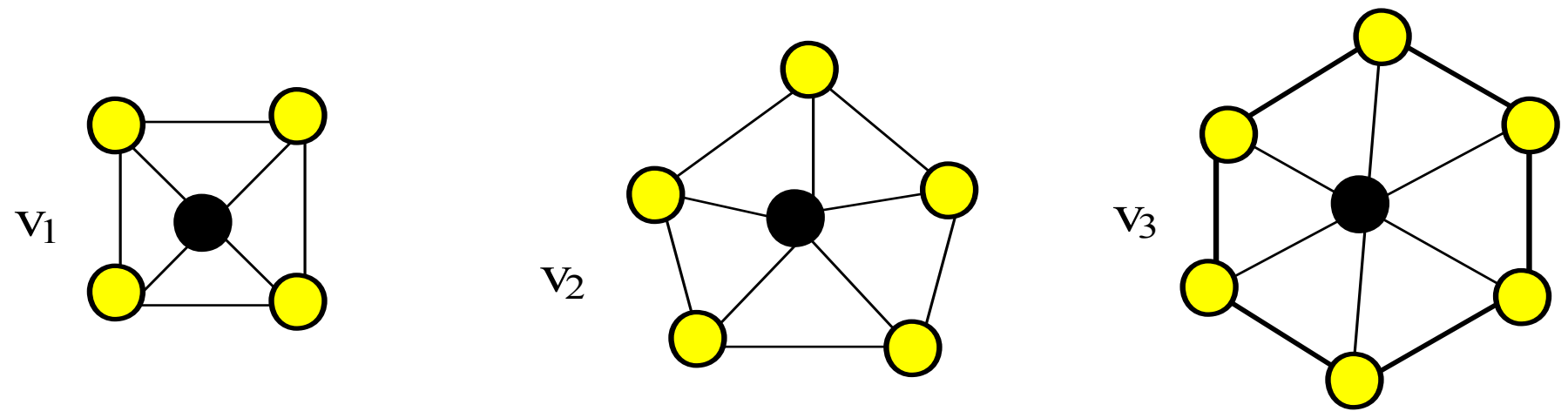

Рис. HH Точки $v_1$, $v_2$ и $v_3$ двумерны, так как их окаемы-одномерные сферы.

преимущественно связные пространства, в которых каждая пара точек может быть соединена цепью.

Определение нормальной прямой и отрезка

Молекулярное связное пространство, состоящее из бесконечного множества одномерных точек, будет нормальной бесконечной прямой. Нормальным отрезком прямой или нормальным одномерным пространством с краем будет связное пространство, состоящее из не менее чем трех точек, все точки которого одномерны, кроме двух крайних точек, окаемы которых содержат по одной точке.

На рис. ff I(3), I(4) и I(5) являются одномерными отрезками, состоящими из

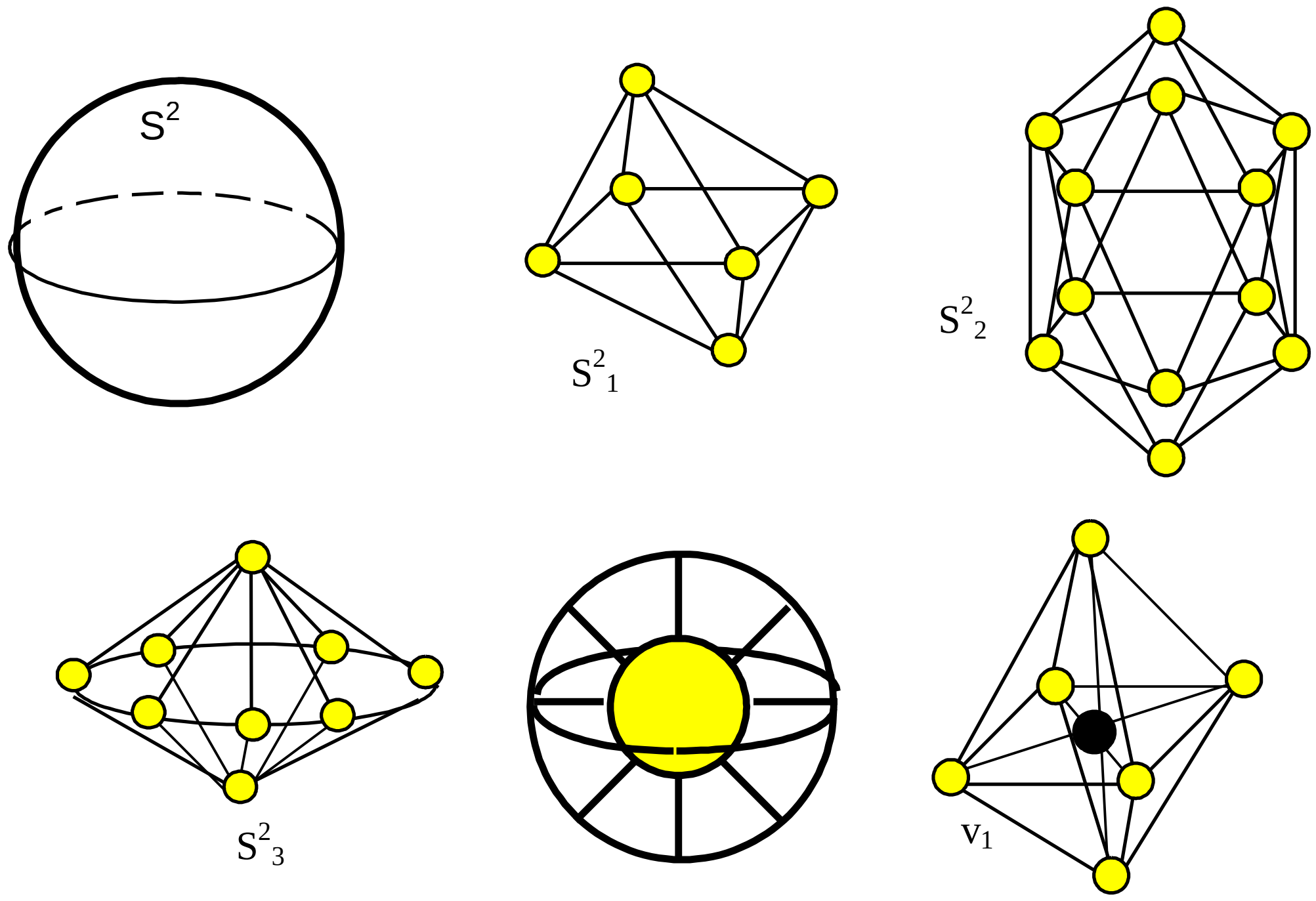

Рис. II Двумерные сферы и трехмерная точка. Окаем каждой точки двумерной сферы есть окружность. Двумерная сфера $S^2{}_1$ минимальна, она состоит из 6 точек-наименьшего числа точек, необходимых для образования двумерной сферы. Сфера $S^2{}_2$ состоит из 12 точек, а сфера $S^2{}_3$ из 9 точек. Точка $v_1$ трехмерна, как имеющая своим окаемом двумерную сферу.

3-х, 4-х и 5-и точек. В отрезке I(3) содержится только одна одномерная точка. В отрезках I(4) и I(5) содержится по 2 и 3 одномерных точки. Крайние точки всех этих отрезков смежны только с одной точкой и не являются одномерными. Иными словами одномерные точки являются внутренними точками отрезков.

Для дальнейшего наиболее важной является одномерная сфера-окружность, мы дадим ее определение.

Определение нормальной окружности

Связное пространство $S^1$, состоящее из конечного числа точек, называется нормальной замкнутой одномерной сферой-окружностью, если все его точки нормальны и одномерны (рис. gg).

Легко видеть, что наименьшее число точек окружности равно 4. На рис. gg показаны три окружности с числом точек 3, 4 и 5. Окружность является единственным замкнутым нормальным одномерным пространством. Теперь дадим общее определение, используя индукцию.

Определение n-мерного нормального замкнутого пространства

Рассмотрим целые числа n, n > 0. Связное пространство G на конечном множестве точек называется n-мерным нормальным замкнутым пространством, если окаем любой его точки есть замкнутое нормальное (n-1)-мерное пространство.

Определение n-мерной нормальной точки

Точка v пространства G называется n-мерной нормальной точкой, если окаем O(v) этой точки является нормальным (n-1)-мерным замкнутым пространством.

Подчеркнем, что отрезок прямой (рис. ff) не является замкнутым нормальным пространством, так как не все его точки являются одномерными; это-пара концевых точек отрезка. Такое пространство

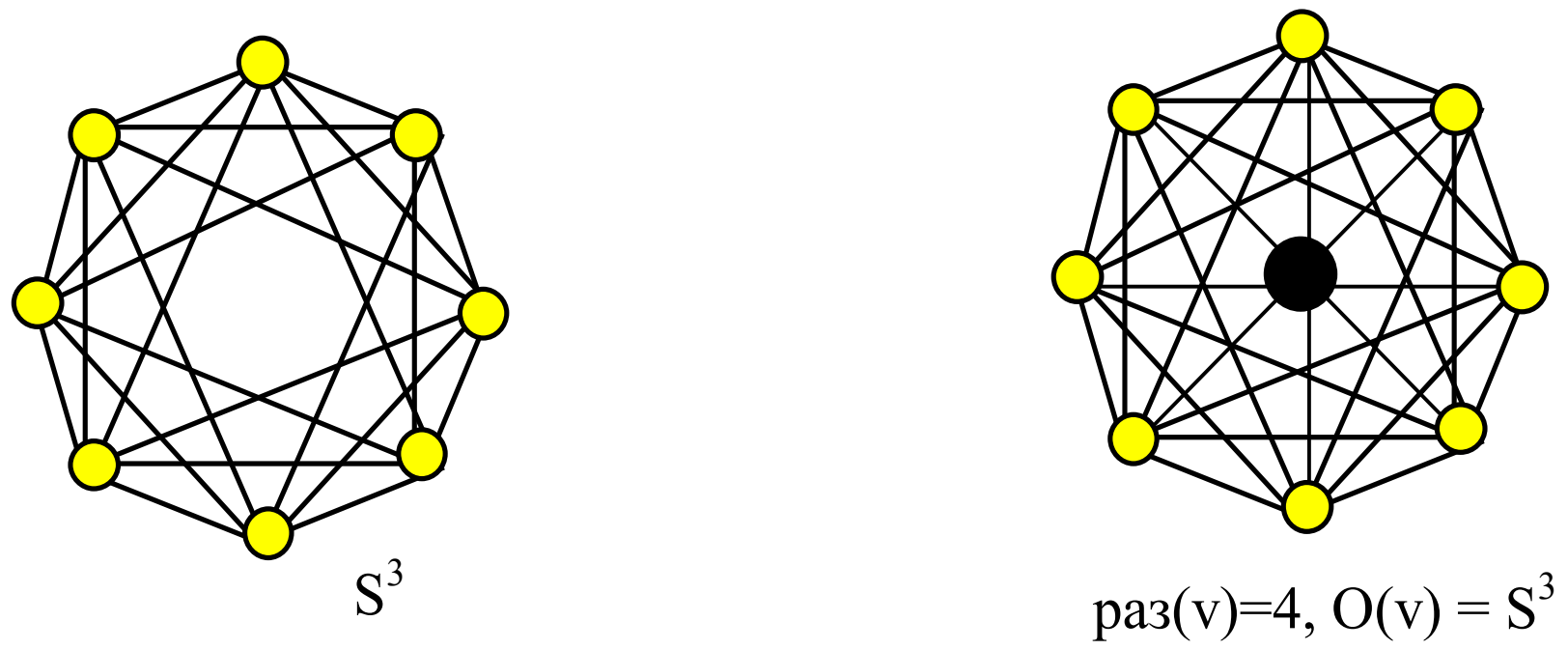

Рис. JJ Слева изображена трехмерная сфера $S^3$, состоящая из 8 точек-наименьшего числа точек, необходимых для образования трехмерной сферы. Справа находится четырехмерная точка, называется пространством с краем, позднее мы будем их изучать. Следует отметить, что в классической теории размерности имеется определение размерности топологического пространства, аналогичное вышеданному [40,41]. Это так называемая малая индуктивная размерность ind X. Кроме нее имеются размерности Ind X и dim X, а также другие. Для молекулярных пространств также можно ввести размерности по аналогии с Ind X и dim X. Однако, различия между ними невелики, поэтому мы будем использовать на протяжении всей книги данное определение.

Обозначение.

Если точка v или пространство X имеют размерность n, мы будем использовать обозначение раз(v) или раз(X).

Рассмотрим примеры молекулярных пространств различной размерности. рис. hh показывает, как образуется 2-мерная точка. Точки $v_2$, $v_3$ и $v_4$ двумерны, так как их окаемы-одномерные сферы. На рис. ii мы изобразили двумерные сферы и трехмерную точку. Окаем каждой точки двумерной сферы есть окружность. Двумерная сфера $S^2_1$ минимальна, она состоит из 6 точек-наименьшего числа точек, необходимых для образования двумерной сферы. Сфера $S^2_2$ состоит из 12 точек, а сфера $S^2_3$ из 9 точек.

$E^2_1$

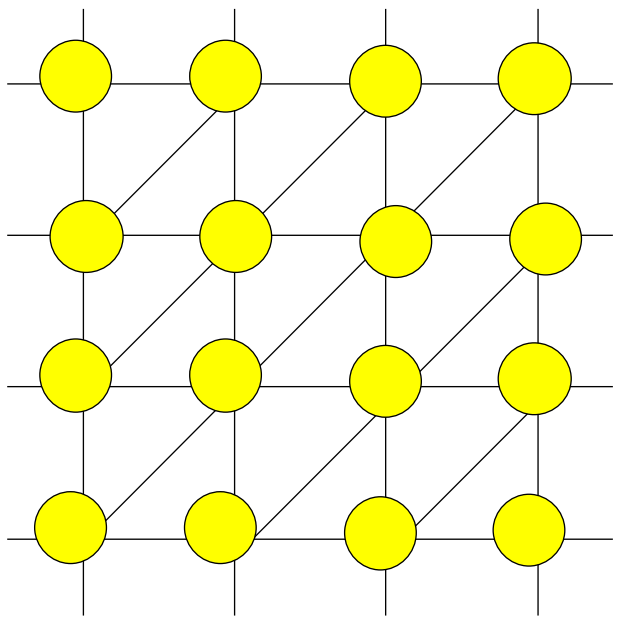

$E^2_2$

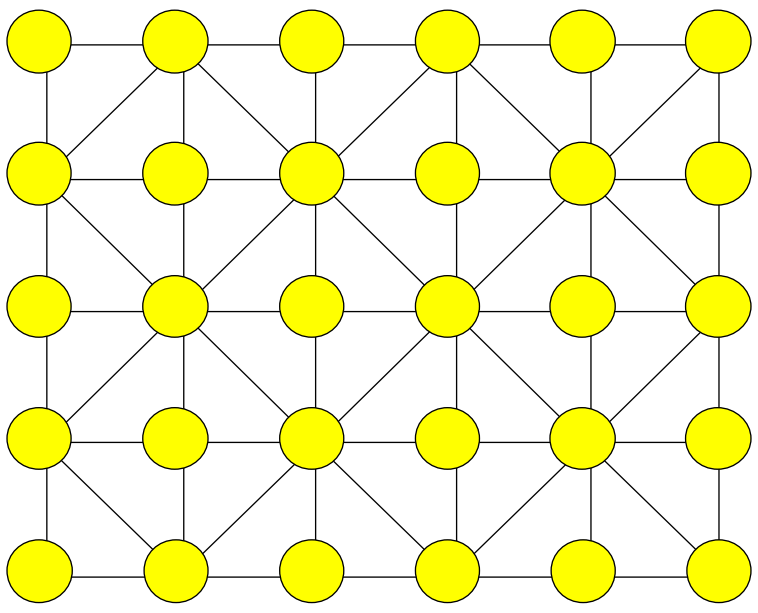

$E^2_3$

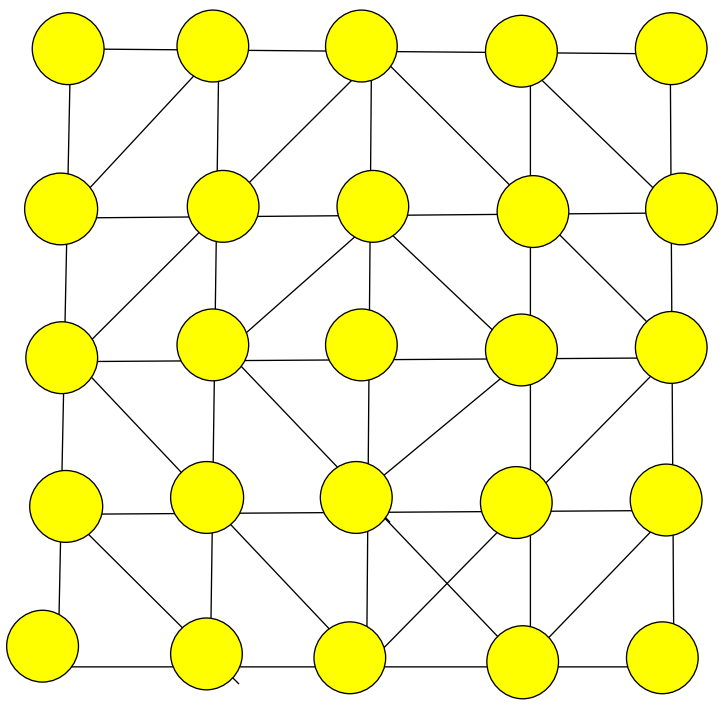

Рис. КК Различные воплощения молекулярной двумерной плоскости. $E^2_1$ состоит из точек двух видов, $E^2_2$ представляет собой однородную двумерную плоскость, $E^2_3$ есть молекулярная плоскость с центральной симметрией.

Легко видеть, что $S^2_1$ и $S^2_2$ являются однородными в том смысле, что окаемы всех точек в каждой из них одинаковы, или, как говорят, изоморфны. Было доказано, что не существует других однородных сфер, кроме $S^2_1$ и $S^2_2$. Точка $v_1$ трехмерна, как имеющая своим окаемом двумерную сферу.

На рис. jj изображена трехмерная сфера $S^3$, состоящая из 8 точек-наименьшего числа точек, необходимых для образования трехмерной сферы. Окаем каждой точки этой сферы есть минимальная двумерная сфера $S^2_1$. Сфера $S^3$, как легко видеть, однородна. До сих пор неизвестно, существует ли еще какая-либо нормальная однородная сфера с числом точек, отличным от 8. Точка v, изображенная на том же рисунке, является четырехмерной, поскольку имеет своим окаемом трехмерную сферу.

З а д а ч а

Найти однородную трехмерную сферу с числом точек,

превышающим 8, или доказать, что ее не существует.

Известно, что в физике рассматривается замкнутая модель вселенной, являющаяся 3-мерной расширяющейся сферой. В момент своего образования эта сфера, очевидно, была наименьшей, состоящей из 8 точек. Количество точек пространства в момент рождения должно было как то повлиять на всю последующую историю вселенной, вплоть до настоящего момента. Действительно, число 8 встречается, например, в физике элементарных частиц,

З а д а ч а

Выяснить, имеется ли связь между физическими объектами, характеризуемыми числом 8 и трехмерной наименьшей сферой, состоящей из 8 точек.

Помимо замкнутых нормальных молекулярных n-мерных пространств, содержащих всегда конечное множество точек существуют пространства на бесконечном счетном множестве точек, такие как модели цилиндрических поверхностей, плоскостей и тому подобное.

Определение n-мерного пространства

Рассмотрим целые числа n, n≥1. Пространство G на бесконечном счетном множестве точек называется n-мерным нормальным (не замкнутым) пространством, если окаем любой его точки есть замкнутое нормальное (n-1)-мерное пространство.

Общее определение пространств, моделирующих евклидовы пространства различной размерности n>1 мы введем ниже, когда будем рассматривать

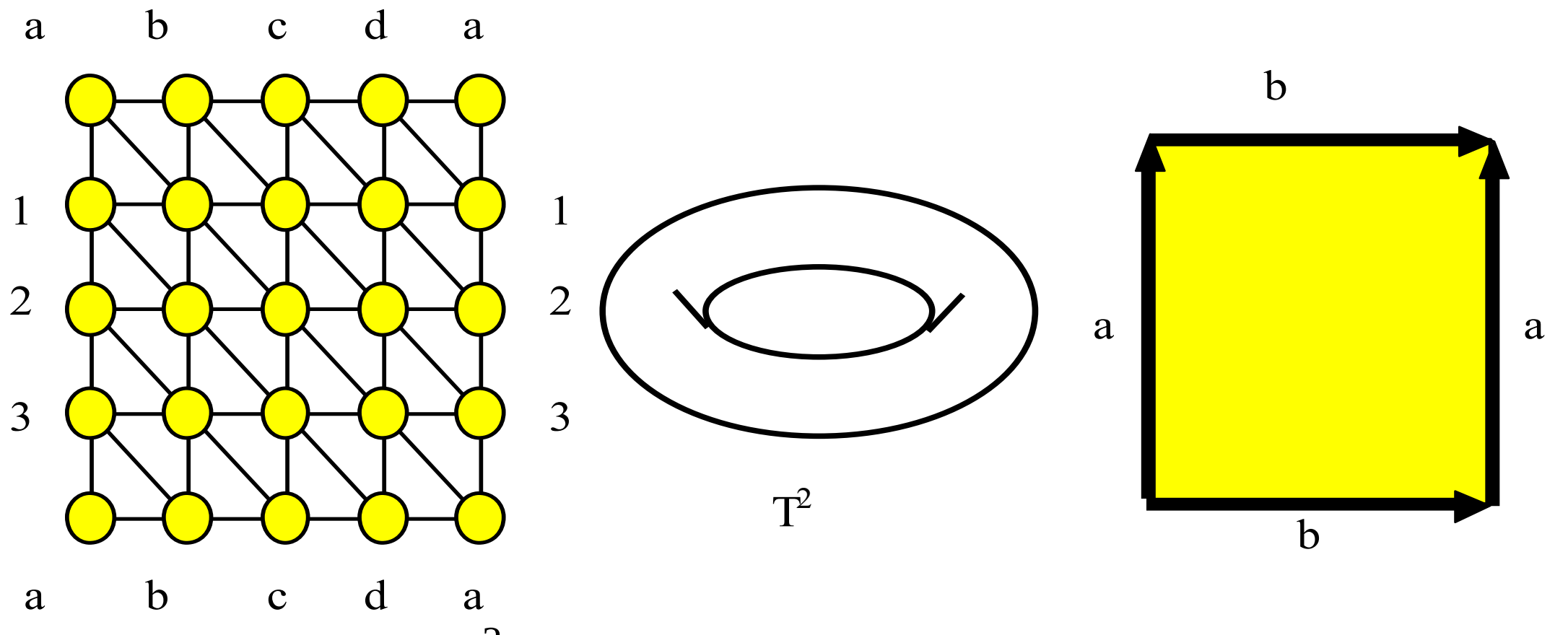

Рис. LL Двумерный тор $T^2$. В непрерывном случае тор может быть получен склеиванием сторон квадрата как это показано справа. Одинаковыми буквами или цифрами обозначены два или несколько образов одной и той же вершины регулярные пространства. Определение одномерного евклидова пространства-одномерной бесконечной прямой (в топологическом смысле) мы уже ввели. Здесь просто отметим, что евклидово нормальное молекулярное n-мерное пространство состоит из нормальных n-мерных точек и обладает основными свойствами непрерывного евклидова пространства.

рис. kk представляет различные молекулярные двумерные плоскости. $E^2_1$ является молекулярным аналогом дискретного двумерного частично-упорядоченного пространства $X^2$ [30,31,32,33,34]. Структура $X^2$ не позволяла рассматривать его как молекулярное пространство, в отличие от $E^2_1$. $E^2_2$ представляет собой однородную двумерную плоскость. $E^2_3$ есть молекулярная плоскость с центральной симметрией.

На рис. ll мы имеем двумерный тор $T^2$. В непрерывном случае эта поверхность может быть получена прямым склеиванием противоположных сторон квадрата как это показано на рисунке [1]. Молекулярный нормальный тор состоит из 16 точек. Прямая проверка показывает, что 16 является наименьшим числом точек, необходимым для построения этой поверхности. В этом смысле тор “сложнее” трехмерной сферы, для которой требуется как минимум 8 точек. Нормальный молекулярный

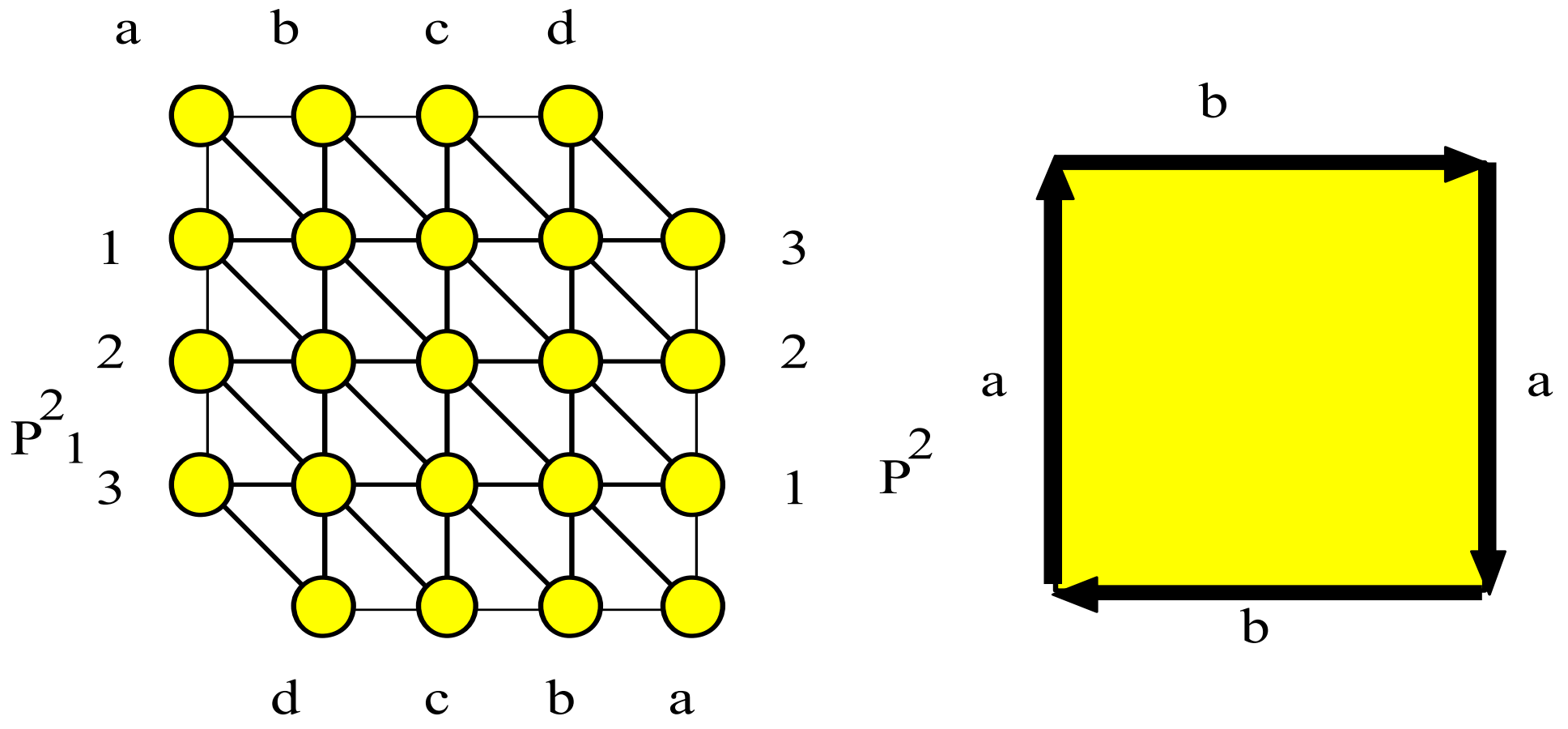

Рис. ММ Двумерная проективная плоскость $P^2$. Как в непрерывном, так и в молекулярном случае $P^2$ может быть получена склеиванием крест накрест сторон квадрата как это показано справа. Одинаковыми буквами или цифрами обозначены два или несколько образов одной и той же вершины.

минимальный тор является, как это легко видеть из рисунка, однородным пространством.

Следующая поверхность (рис. mm) является проективной плоскостью, поверхностью, возникающей при склеивании противоположных сторон квадрата крест накрест, как это показано на рисунке [1]. Точно такое же склеивание крест накрест точек, лежащих на противоположных сторонах нормального молекулярного квдрата, дает нормальную проективную плоскость. Проективную плоскость можно также получить склеивая попарно точки сферы, лежащие на противоположных концах диаметров, проходящих через центр сферы, рис. nn. Поэтому проективную плоскость также называют склеенной сферой.

Еще один путь получения проективной плоскости-заклеивание лентой Мебиуса дырки в сфере (рис. nn). Все эти возможности реализуются в молекулярном пространстве. Наименьшая нормальная проективная

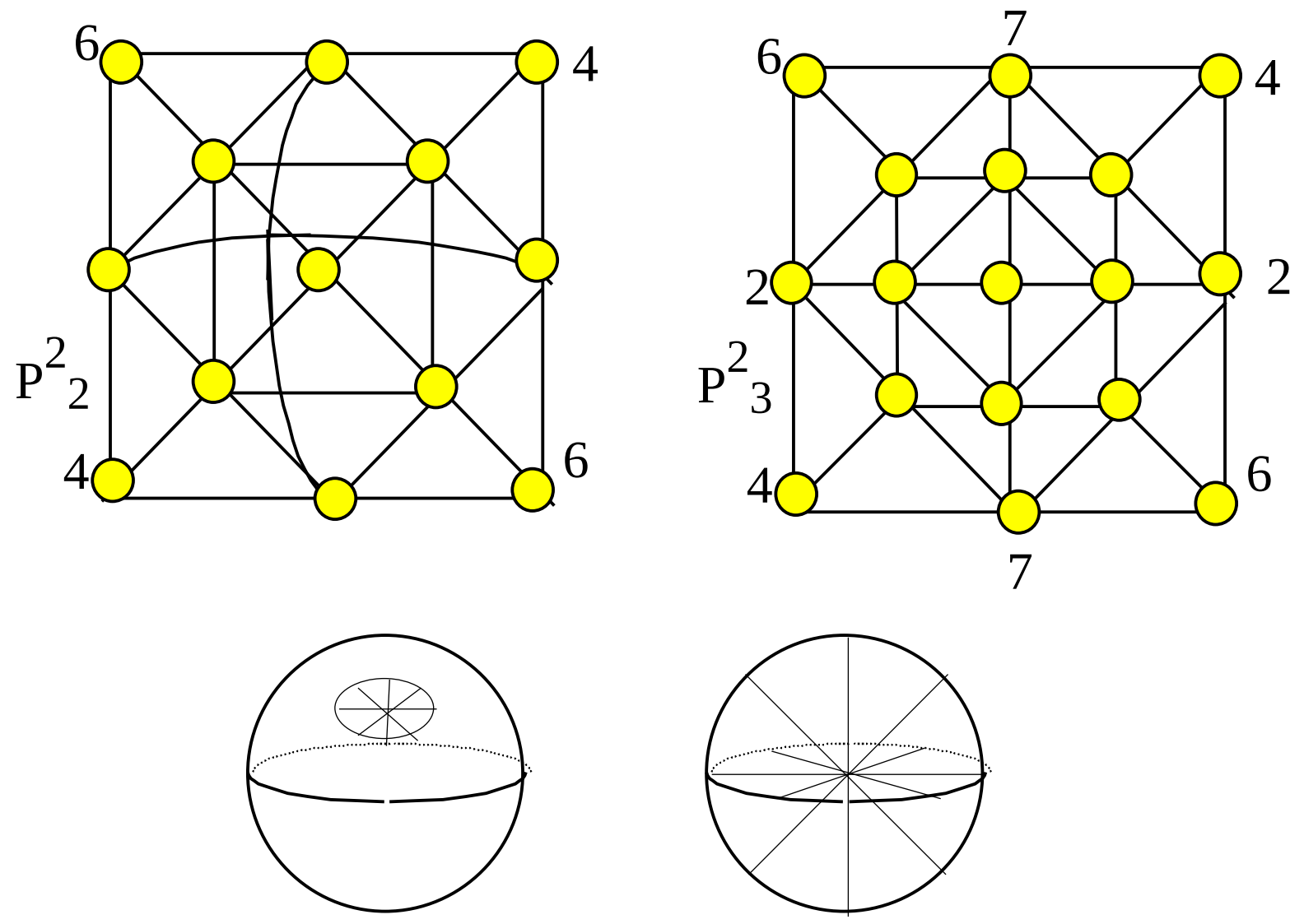

Рис. NN Две нормальные молекулярные проективные плоскости. $P^2_2$ состоит из 11 точек и является наименьшей возможной моделью. В непрерывном случае проективная плоскость может быть также получена заклеиванием дырки на сфере при помощи ленты Мебиуса, или же попарным склеиванием точек сферы, лежащих на противоположных концах диаметров. Одинаковыми буквами или цифрами обозначены склеенные точки.

плоскость состоит из 11 точек. Если рассматривать только замкнутые двумерные поверхности, то проективая плоскость является второй по объему поверхностью после двумерной сферы, для которой минимальный объем равен 6. Еще одна особенность этой поверхности состоит в том, что нормальная молекулярная проективная плоскость является всегда неоднородной поверхностью, независимо от количества точек, ее образующих (это будет показано позднее). На рис. nn нарисованы две модели проективной плоскости, одна из которых есть наименьшая проективная плоскость из 11 точек. Легко видеть, что все модели неоднородны. Каждая из трех моделей (рис. mm, рис. nn) может быть получена из любой другой путем точечных преобразований.

Следующее пространство (рис. oo) является нормальной молекулярной бутылкой Клейна. Как в классическом, так и в молекулярном случае, она возникает при склеивании двух противоположных сторон квадрата крест-накрест, а двух других-напрямую. В следующих главах мы подробно изучим свойства этих поверхностей и покажем, такие структуры, как эйлерова характеристика и гомологические группы молекулярных

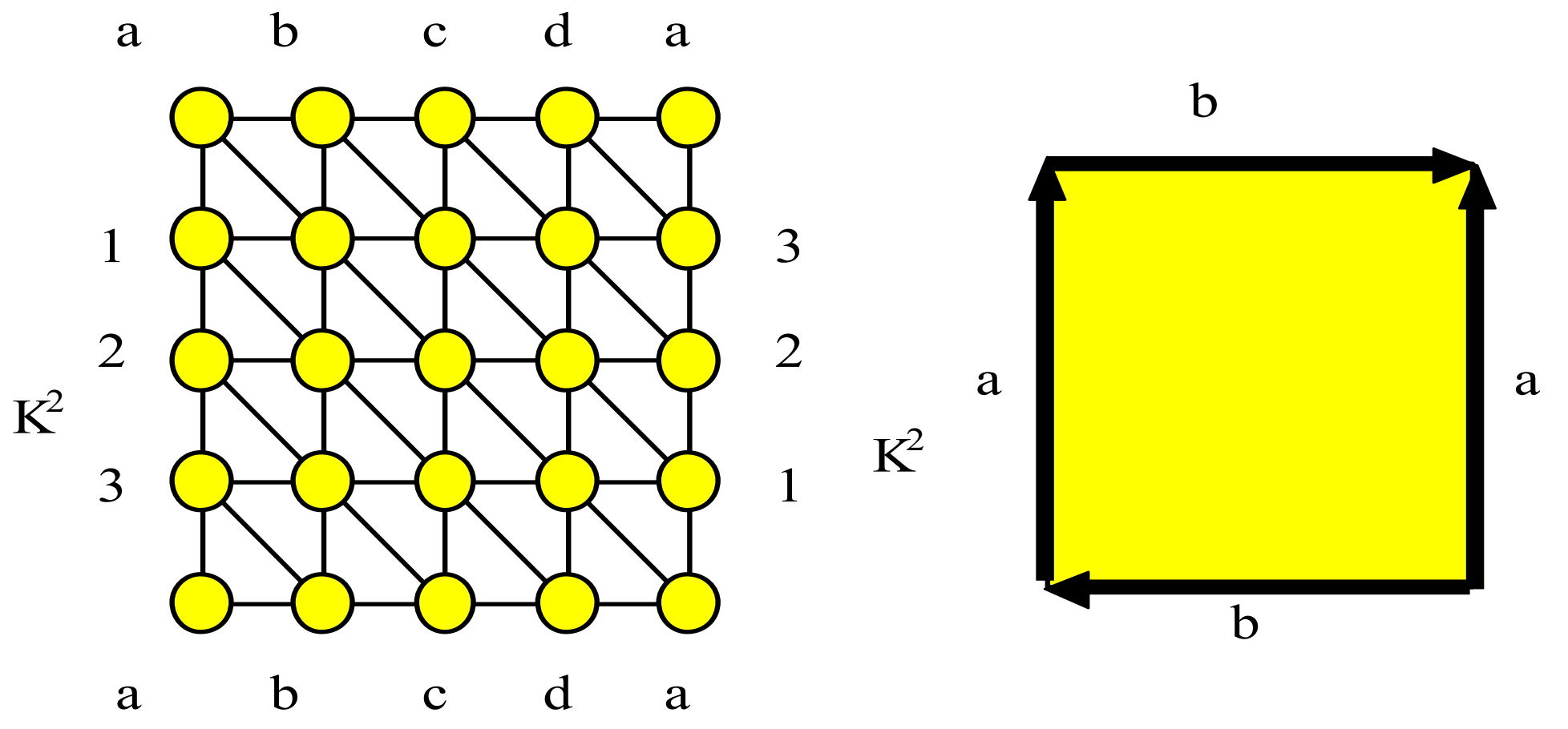

Рис. OO Двумерная молекулярная бутылка Клейна $K^2$. В непрерывном случае бутылка Клейна может быть получена склеиванием сторон квадрата как это показано справа. Одинаковыми буквами или цифрами обозначены два или несколько образов одной и той же точки.

пространств и их непрерывных прообразов совпадают. На всех рисунках одинаковыми буквами или цифрами обозначены два или несколько образов одной и той же точки. Это сделано, чтобы избежать излишней усложненности картинок.

Задача.

Доказать, что молекулярная бутылка Клейна $K^2$ минимальна, или найти минимальную.

## Список литературы к главе 3.

# СВОЙСТВА НОРМАЛЬНЫХ МОЛЕКУЛЯРНЫХ ПРОСТРАНСТВ

We investigate some properties of normal molecular spaces that are much similar to properties of continuous spaces, and define and study homeomorphic transformations on normal spaces.

Некоторые результаты этой главы имеются в работах [7,10,43,50,51]

## *СВОЙСТВА РАЗМЕРНОСТИ НОРМАЛЬНЫХ МОЛЕКУЛЯРНЫХ ПРОСТРАНСТВ*

Молекулярные пространства обладают свойствами, во многом аналогичными свойствам непрерывных пространств, образами которых они являются. К таким свойствам относится структура пространства в малом, особенности операций суммы и произведения нескольких нормальных пространств, операции, позволяющая получить новые

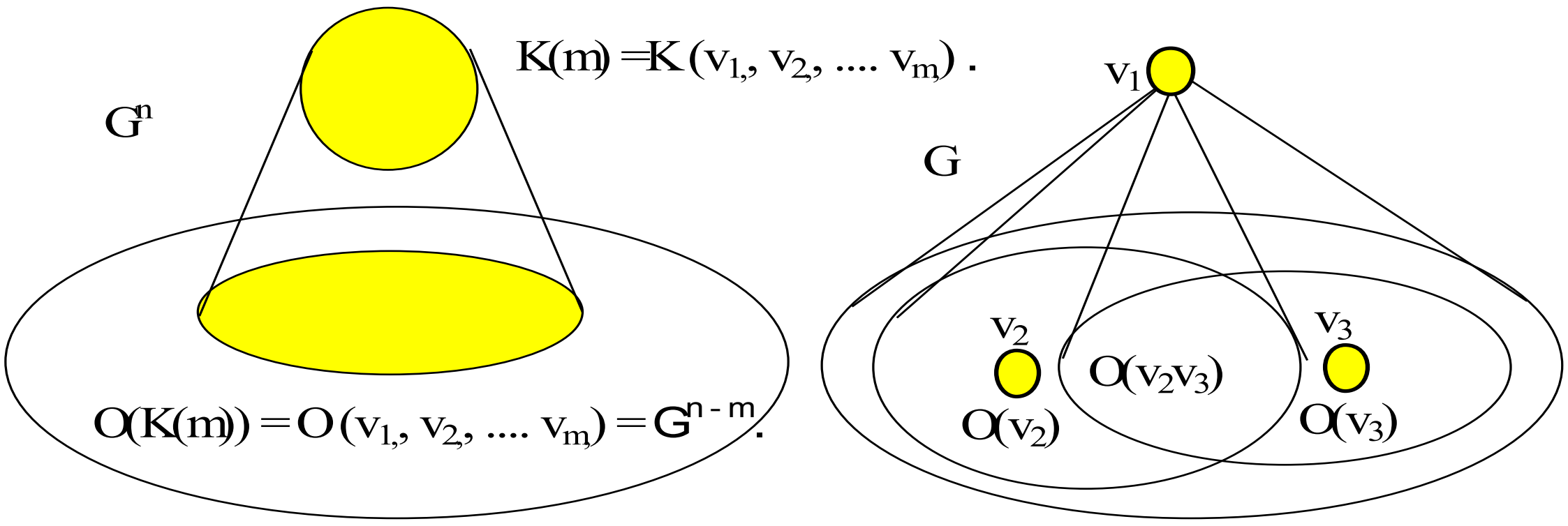

Рис. 42 Слева. Общий окаем подпространства K(m) из m попарно связных точек есть пространство $G^{n-m}$ размерности (n-m). Справа. Пространство G является конусом. Окаем $O(v_2)$ точки $v_2$, а также общие окаемы любого множества точек, принадлежащих $O(v_1)$ являются конусами.

пространства, группы гомологий, эйлерова характеристика и некоторые другие свойства. Мы рассмотрим основные свойства молекулярных пространств и методы конструирования, позволяющие естественным путем получить новые молекулярные пространства , в которых различные точки имеют различную топологию. Будет полезно еще раз подчеркнуть те базовые свойства молекулярного нормального пространства, которые следуют непосредственно из его определения.

Теорема 20

1. *Нормальное пространство $G^n$ ненулевой размерности всегда является связным.*

*2. Окаем любой точки нормального пространства $G^n$, n>1, всегда является связным.*

Утверждения 1 и 2 следуют непосредственно из определения нормального пространства.

Теорема 21

*Пусть Gn является молекулярным (не обязательно замкнутым) нормальным n-мерным пространством, где K(m) есть связка (подпространство из m попарно связных точек (v1,v2,v3,.. vm)). Тогда общий окаем этих точек есть молекулярное замкнутое нормальное (n-m)-мерное пространство.*

*$O(v1v2\ v3 \ldots vm)=O(v1)\cap O(v2) \cap O(v3)\cap\ldots\cap O(vm)=Gn-m$.*

Доказательство.

По определению нормального молекулярного n-мерного пространства окаем $O(v_1)=G^{n-1}$ произвольной точки пространства $v_1$ является молекулярным замкнутым нормальным (n-1) мерным пространством. Пусть $v_1$ и $v_2$ смежны, то есть между этими точками имеется связь. Это означает, что $v_2 \in G^{n-1}=O(v_1)$. Следовательно, $O(v_2)|G^{n-1}=O(v_2v_1)$ является молекулярным замкнутым нормальным (n-2)-мерным пространством. Пусть $v_3$ смежна с $v_2$ и $v_1$, то есть $v_3 \in O(v_2v_1)=G^{n-2}$. Следовательно,

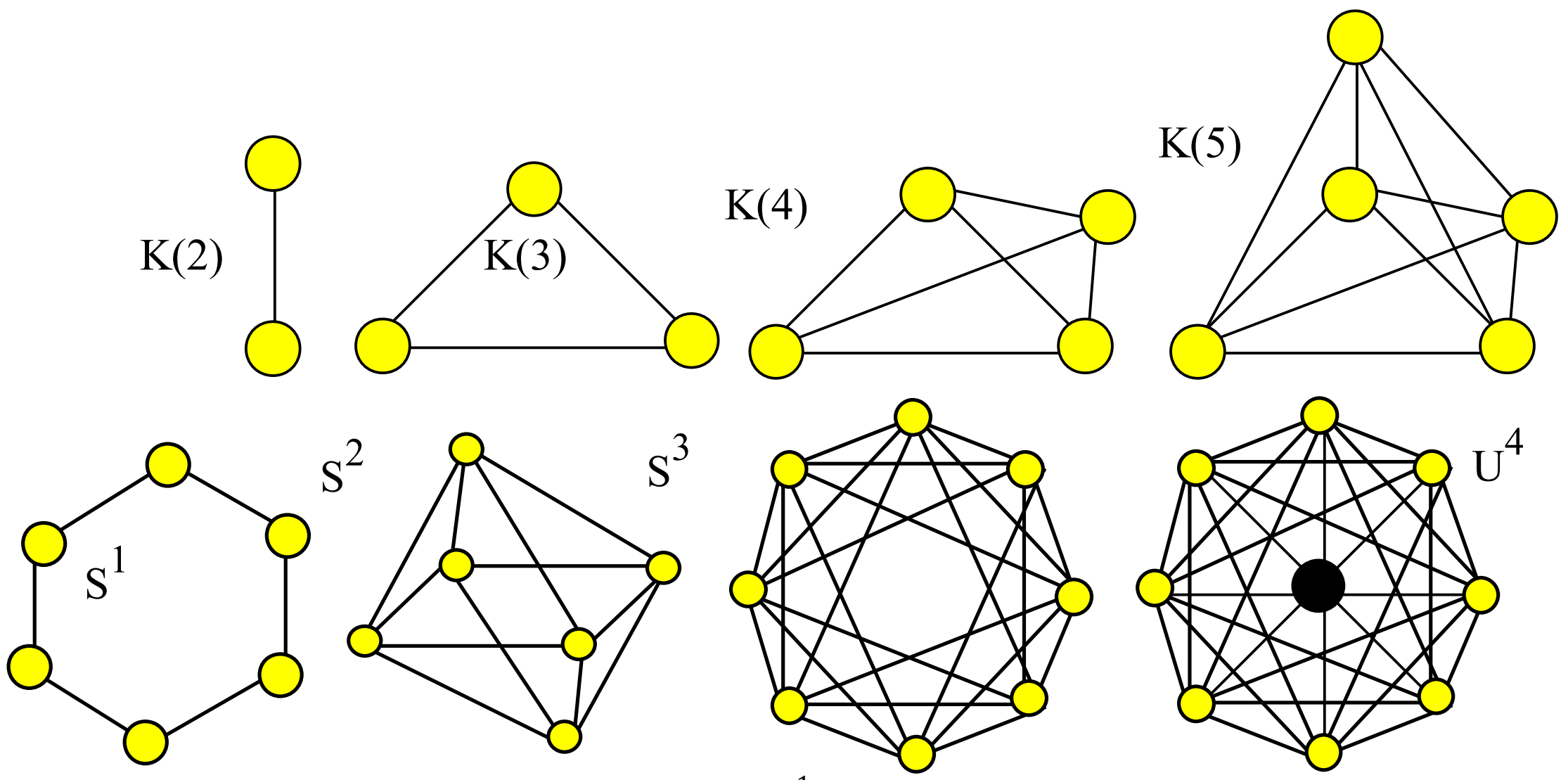

Рис. 43 Одномерная сфера-окружность $S^1$ состоит из связок K(2), двумерная сфера $S^2$ построена из связок K(3), трехмерная сфера $S^3$ образована связками K(4), элемент четырехмерного пространства-четырехмерный шар U(v) точки v состоит из связок K(5).

$O(v_3)|G^{n-2}=O(v_3v_2v_1)$ является молекулярным замкнутым нормальным (n-3) мерным пространством. Повторяя этот процесс для точек $(v_1,v_2,v_3,\ldots v_m)$ мы получаем, что $O(v_1v_2v_3\ldots v_m)=O(v_1)\cap O(v_2)\cap O(v_3)\cap\ldots\cap O(v_m)=G^{n-m}$.

Теорема доказана. □

Следующая теорема позволяет довольно точно определить локальную структуру нормального молекулярного пространства. Оказывается, что n-мерное нормальное пространство состоит из связок K(n+1) на (n+1) точках. Иными словами двумерные нормальные пространства состоят из точек,

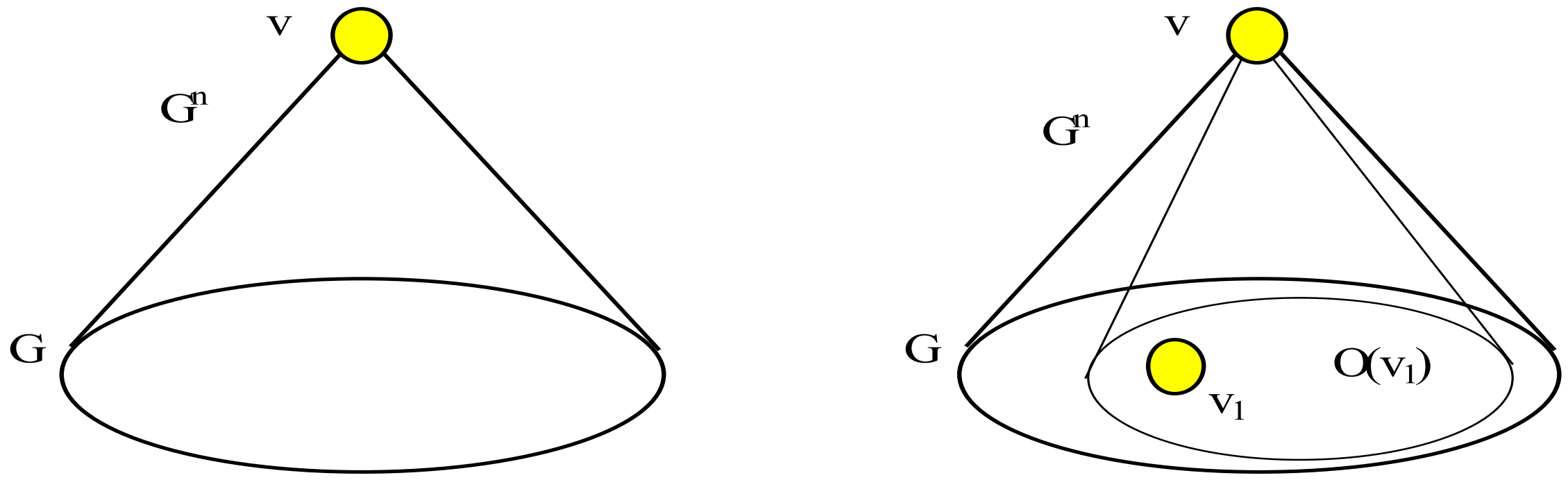

Рис. 44 Нормальное пространство $G^n$ не может быть конусом v⊕G. Если пространство есть конус, то окаем $O(v_1)$ любой точки $v_1$, принадлежащей G также есть конус.

образующих треугольники, трехмерные пространства состоят из тетраэдров и так далее. Как это видно из Рис. 43 одномерная сфера-окружность состоит из связок K(2), двумерная сфера построена из связок K(3), трехмерная сфера образована связками K(4), элемент четырехмерного пространства-четырехмерный шар точки v состоит из связок K(5).

Теорема 22

*Пусть $G^n$ есть нормальное молекулярное пространство размерности n. Тогда для любой связки $K(s)\subseteq G^n$, $s \le n+1$, существует связка K(n+1), содержащая K(s), $K(s)\subseteq K(n+1)\subseteq G^n$.*

Доказательство.

Используем построения предыдущей теоремы. Пусть $K(s)=K(v_1,v_2,v_3,...v_s)$. Согласно предыдущему общий окаем точек $v_1$, $v_2$, $v_3$,... $v_s$ есть нормальное молекулярное замкнутое пространство размерности (n-s), $O(v_1v_2\ v_3,...\ v_s)=O(v_1)\cap O(v_2)\cap O(v_3)\cap...\cap O(v_s)=G^{n-s}$. Выберем произвольную точку $v_{s+1}$ из $G^{n-1}$. Она является смежной со всеми предыдущими и $O(v_1v_2v_3....v_{s+1})=O(v_1)\cap O(v_2)\cap O(v_3)\cap...\ \cap O(v_{s+1})=G^{n-s-1}$. Повторяя этот процесс для точек $(v_1,v_2,v_3,...v_{n+1})$ мы получаем, что $O(v_1v_2v_3...v_{n+1})$ есть пустое множество, а все эти точки связны между собой. Связка $K(n+1)=K(v_1,v_2,v_3,..v_{n+1})$ состоит из n+1 точек и содержит исходную связку. Теорема доказана. □

Иллюстрацией этой теоремы для одномерных, двумерных и трехмерных случаев служат пространства, изображенные на ранее и на Рис. 43. Из этой теоремы следует, что n-мерное нормальное молекулярное пространство не содержит связок K(n+2).

С л е д с т в и е .

Нормальное n-мерное молекулярное пространство не содержит связки K(n+2) из (n+2) точек.

Одним из важных свойств замкнутого молекулярного пространства является то, что оно не может быть конусом, то есть пространством, в котором одна точка имеет связь со всеми остальными. Напомним, что конусом v⊕G называется молекулярное пространство G, в котором имеется точка v, соединенная со вмести остальными точками пространства.

Теорема 23

*Нормальное молекулярное пространство любой размерности не может быть конусом.*□

Д о к а з а т е л ь с т в о .

Используем доказательство от противного. Предположим, что нормальное молекулярное замкнутое пространство $G^n$ размерности n есть конус, и точка v связана со всеми остальными точками пространства, то есть $G=O(v)=G^n-v$ (Рис. 44). Окаем O(v) является нормальным молекулярным замкнутым пространством размерности (n-1). Возьмем любую точку из O(v), скажем, $v_1$. Окаем этой точки $O(v_1)|O(v)$ в пространстве O(v) есть нормальное молекулярное замкнутое пространство размерности (n-2), и также является конусом. Возьмем любую точку из $O(v_1)|O(v)$, скажем, $v_2$. Окаем этой точки $O(v_2)|O(v_1)|O(v)$ в пространстве $O(v_1)|O(v)$ есть нормальное молекулярное замкнутое пространство размерности (n-3) и также конус. Продолжая этот процесс n раз мы получаем, что окаем точки $v_{n-1}$ $O(v_{n-1})|O(v_n)|O(v_{n-1})|\ldots |O(v_1)|O(v)$ в пространстве $O(v_{n-2})|O(v_{n-3})|\ldots |O(v_1)|O(v)$ есть нормальное молекулярное замкнутое пространство размерности (n-n)=0 и также конус. Однако, этот окаем должен быть нуль-мерным нормальным замкнутым пространством и состоять, согласно определению, из двух несвязных точек, то есть не может быть конусом. Мы получили противоречие. Следовательно наше предположение неверно. Теорема доказана. □

Прямым и очевидным следствием этой теоремы является в нормальном пространстве наличие для любой точки пространства несмежной ей точки.

С л е д с т в и е .

Пусть G есть нормальное пространство. Тогда для любой точки a∈G найдется несмежная ей точка b∈G..

Докажем еще одну теорему, которая определяет минимальное число точек замкнутого нормального молекулярного пространства.

Теорема 24

*Нормальное молекулярное замкнутое пространство $G^n$ размерности n имеет не менее чем 2n+2 точек $|G^n| \geq 2n+2$.*

Д о к а з а т е л ь с т в о.

Используем доказательство от противного.

Предположим, что нормальное молекулярное замкнутое пространство $G^n$

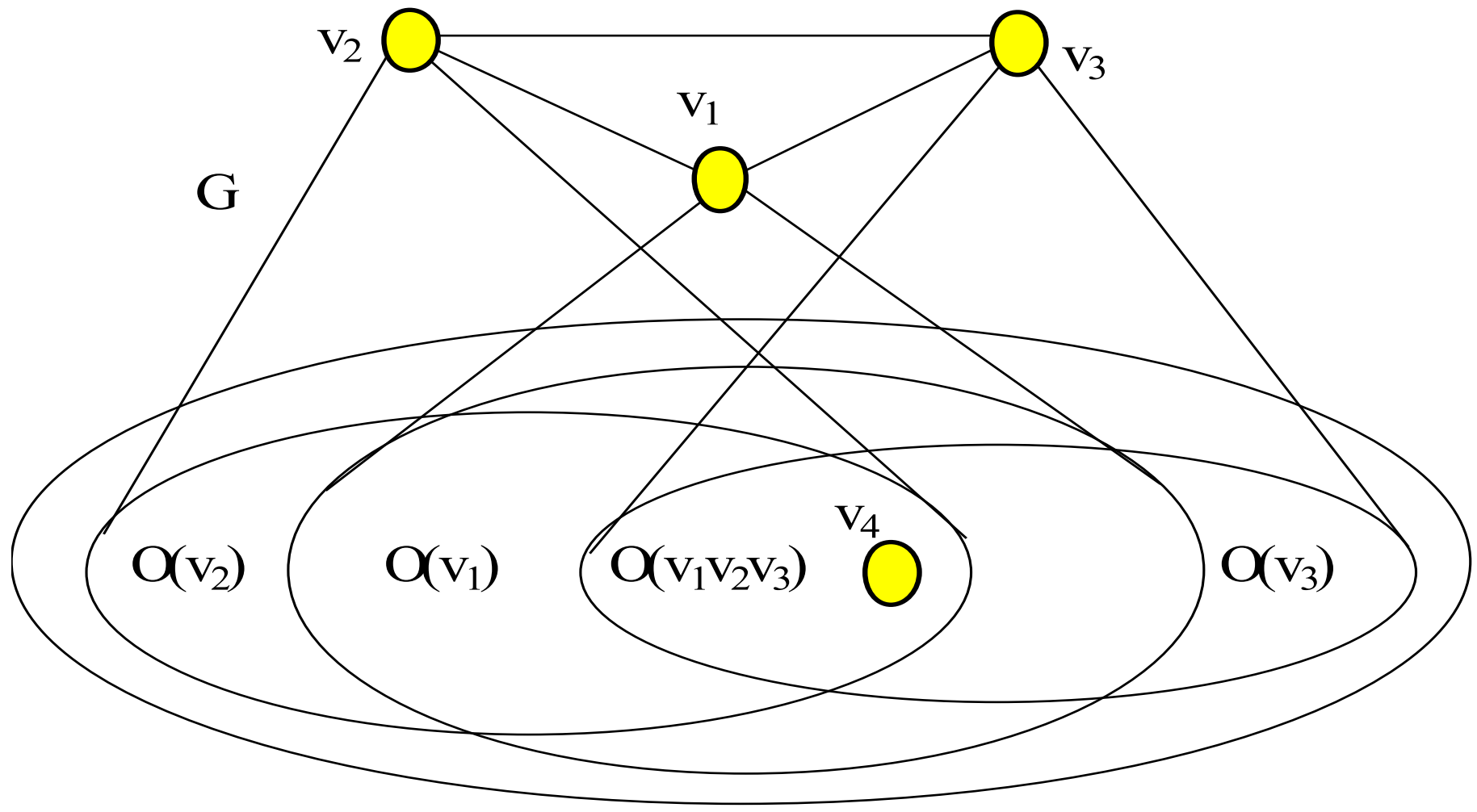

Рис. 45 Объемы общих окаемов уменьшаются не менее чем на 1 с увеличением числа точек в общем окаеме на 1. $|O(v_1)|-1 \geq |O(v_1v_2)|$, $|O(v_1v_2)|-1 \geq |O(v_1v_2v_3)|,..,|O(v_1)|-n \geq |O(v_1v_2...v_{n+1})|$.

размерности n содержит менее 2n+2 точек (Рис. 45). Выберем произвольную точку $v_1$. Окаем $O(v_1)$ является нормальным молекулярным замкнутым пространством размерности (n-1). Так как согласно предыдущей теореме $v_1$ не имеет связи с некоторыми точками из $G^n$, поскольку $G^n$ не может быть конусом, объем окаема $O(v_1)$ менее 2n, $|O(v_1)|<2n$. Возьмем любую точку из $O(v_1)$, скажем, $v_2$. Окаем этой точки $O(v_2)|O(v_1)$ в пространстве $O(v_1)$ есть нормальное молекулярное замкнутое пространство размерности (n-2) и также не является конусом. Следовательно $|O(v_2)|O(v_1)|<2n-2$. Возьмем любую точку из $O(v_2)|O(v_1)$, скажем, $v_3$. Окаем этой точки $O(v_3)|O(v_2)|O(v_1)$ в пространстве $O(v_2)|O(v_1)$ есть нормальное молекулярное замкнутое пространство размерности (n-3) и также не есть конус. Это означает, что $|O(v_3)|O(v_2)|O(v_1)|<2n-4$. Продолжая этот процесс n раз мы получаем, что окаем точки $v_n$ $O(v_n)|O(v_{n-1})|O(v_{n-2})|... |O(v_2)|O(v_1)$ в пространстве $O(v_{n-1})|O(v_{n-2})|... |O(v_2)|O(v_1)$ есть нормальное молекулярное замкнутое пространство размерности (n-n)=0, и также не конус. При этом объем этого

пространства согласно предположению, удовлетворяет неравенству $|O(v_n)|O(v_{n-1})|O(v_{n-2})|\dots|O(v_2)|O(v_1)|<2n-2n=0$. Однако этот окаем должен быть нуль-мерным нормальным замкнутым пространством и состоять согласно определению из 2-х несвязных точек. Мы получили противоречие. Следовательно наше предположение неверно. Теорема доказана. □

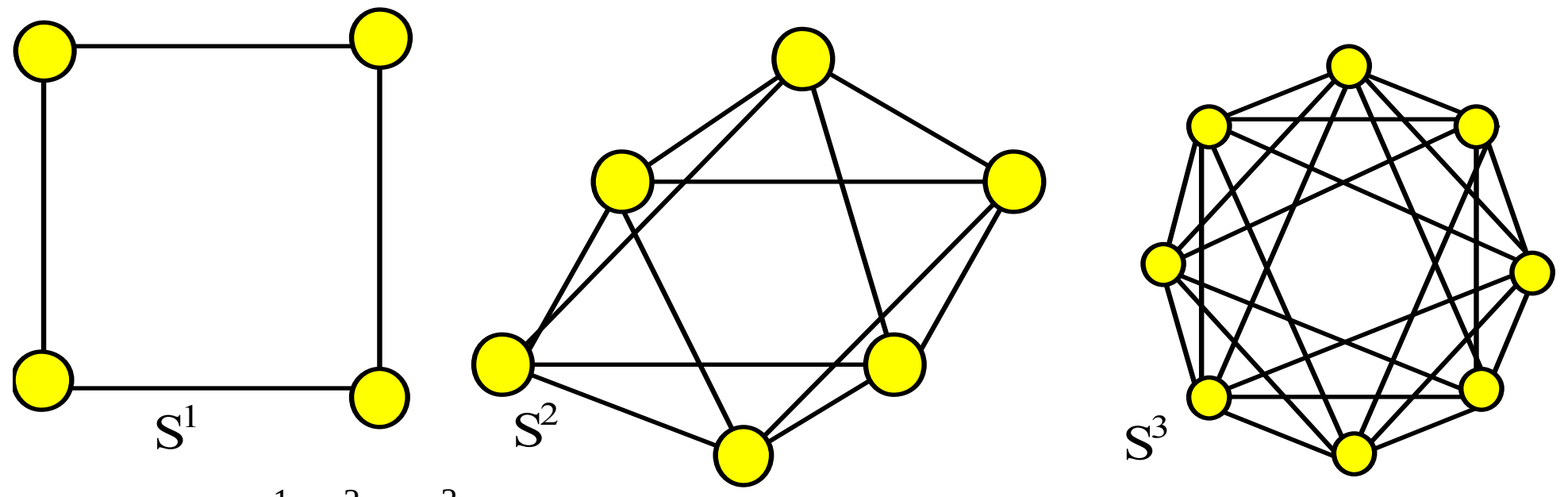

Рис. 46 Сферы $S^1$, $S^2$ и $S^3$ являются единственными замкнутыми n-мерными пространствами, для которых в размерностях 1, 2 и 3 выполняется соотношение $|G^n|=2n+2$, $S^2=S^0\oplus S^1$, $S^3=S^2\oplus S^1$.

Согласно этой теореме, одномерное замкнутое нормальное пространство состоит не менее чем из 4-х точек, двумерное-не менее чем из 6-и, трехмерное-не менее чем из 8-и и так далее. Как мы видели выше минимальная окружность содержит 4 точки, минимальная сфера-6 точек, минимальная трехмерная сфера-8 точек. Возникает вопрос, существует ли замкнутое нормальное n-мерное пространство, отличное от сферы, состоящее из (2n+2) точек? Прямая проверка для пространств с небольшим объемом показывает, что таких пространств не существует. Однако в общем виде это нуждается в доказательстве. На Рис. 46 изображены нормальные замкнутые пространства $G^n$ размерности n=1, 2 и 3, объемы которых составляют 2n+2.

Операция сложения молекулярных уже была определена выше. Эта операция позволяет получить молекулярные пространства с новыми топологическими свойствами. Размерное свойство такой операции аналогично свойству соответствующей операции в классической топологии.

Теорема 25

*Пусть $G^m$ и $S^0$-нормальное молекулярное замкнутое пространство размерности m и 0-мерная сфера с набором точек $V=(v_1,v_2,v_3,\dots v_q)$ и $R=(a,b)$ соответственно. Тогда их прямая сумма $G^m\oplus S^0$ является нормальным замкнутым молекулярным пространством размерности (m+1), $W^{m+1}=G^m\oplus S^0$.*

Доказательство.

Используем индукцию. Для m=0, 1 теорема проверяется непосредственно. Предположим, что теорема верна для $G^m$, где $m \le k$. Пусть m=k+1. Возьмем любую точку из $G^m$, скажем, $v_1$. Окаем этой точки в пространстве $W^{m+1}$ есть сумма окаема этой точки в пространстве $G^m$ и пространства $S^0$,

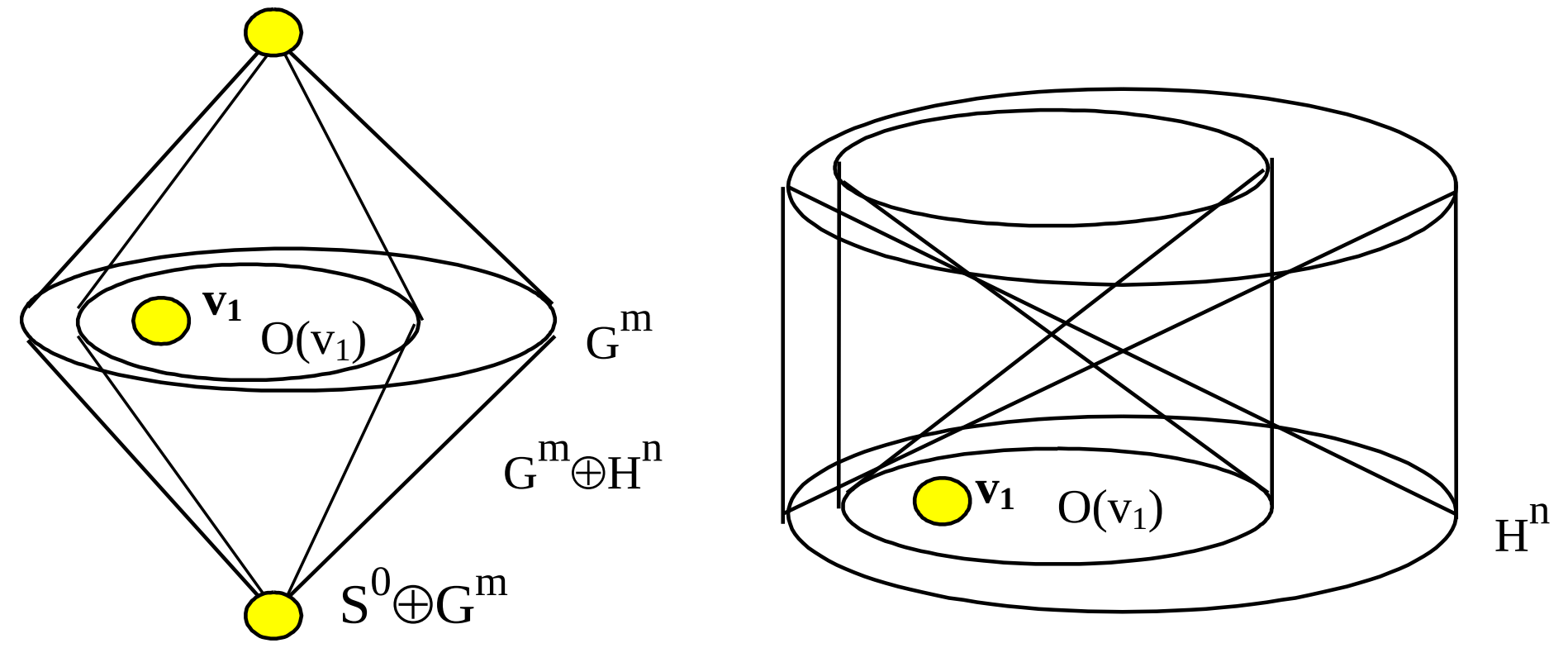

Рис. 47 Сумма 0-мерной сферы $S^0$ и замкнутого n-мерного пространства $G^n$ имеет размерность (n+1). Сумма замкнутого n-пространства $G^n$ и замкнутого m-пространства $H^m$ имеет размерность (n+m+1) то есть $O(v_1)|W^{m+1}=(O(v_1)|G^m)\oplus S^0$. Так как окаем этой точки в пространстве $G^m$ является к-мерным пространством, то, согласно индукции, $O(v_1)|W^{m+1}$ будет нормальным (к+1) мерным пространством. Окаемы точек a и b также являются нормальным (к+1) мерным пространством $G^m$. Следовательно, $W^{m+1}=G^m\oplus S^0$ есть нормальное (к+2) мерное пространство (Рис. 47, Рис. 48). Теорема доказана.□

На Рис. 48 иллюстрируется сумма двух 0-мерных сфер, дающее одномерную сферу, то есть окружность, $S^1=S^0\oplus S^0$.

Следующая теорема обобщает предыдущую на сумму произвольных нормальных замкнутых пространств.

Теорема 26

*Пусть $G^m$ и $H^n$-два нормальных молекулярных замкнутых пространства размерности m и n с наборами точек $V=(v_1,v_2,v_3,...v_q)$ и $R=(r_1,r_2,r_3,...r_t)$ соответственно. Тогда их сумма $G^m\oplus H^n$ является нормальным замкнутым молекулярным пространством размерности (m+n+1), $W^{m+n+1}=G^m\oplus H^n$.*

Доказательство.

Используем индукцию. Для n, m=0, 1 теорема проверяется непосредственно. Предположим, что теорема верна для $G^m$ и $H^n$, где $m+n \le p$. Пусть n+m=p+1. Так как теорема уже доказана для m=0 или n=0, рассмотрим случай, когда $m > 0$ или $n > 0$. Выберем произвольную точку $v_1 \in G^m$ и рассмотрим окаем этой точки в $W^{m+n+1}=G^m \oplus H^n$. Очевидно, что $O(v_1)|W^{m+n+1}=(O(v_1)|G^m)\oplus H^n$, тогда раз$(O(v_1)|G^m)=m-1$, раз$(H^n)=n$. Следовательно, сумма размерностей пространств $(O(v_1)|G^m)$ и $H^n$ равна $(m-1)+n=p$, и, по предположению, их сумма является нормальным замкнутым молекулярным пространством размерности $m+n$, раз$(O(v_1)|W^{m+n+1})$=раз$(O(v_1)|G^m)$+раз$(H^n)+1=(m-1)+n+1=m+n$. Это означает, что размерность пространства $W^{m+n+1}=G^m\oplus H^n$ равна $(m+n+1)$. Теорема доказана.□

Прежде чем переходить к примерам сформулируем еще одну теорему о свойствах суммы, которую доказывать не будем в силу ее очевидности.

Теорема 27

*Сумма нормальных молекулярных замкнутых пространств коммутативна и ассоциативна,* $W^{m+n+1}=G^m\oplus H^n=H^n\oplus G^m$, $W^{m+n+p+1}=(G^m\oplus H^n)\oplus F^p=G^m\oplus(H^n\oplus F^p)$.

Доказательство очевидно.□

Из Рис. 49 видно, что сумма одномерной и нульмерной сфер является двумерной сферой. Следует отметить то обстоятельство, что сумма минимальных сфер также дает минимальную сферу. В связи с этим

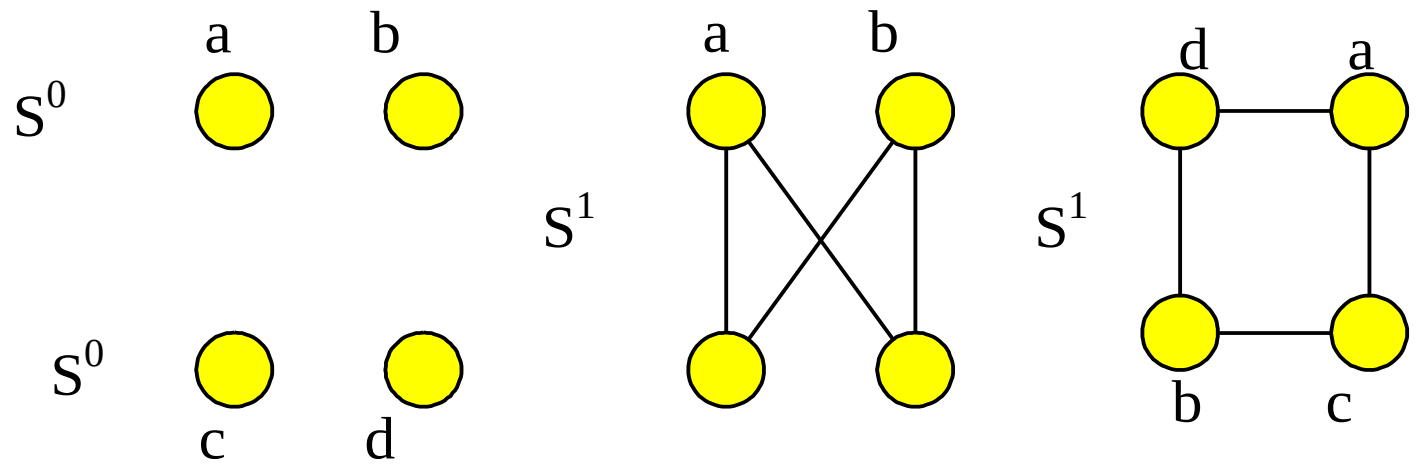

Рис. 48 Сумма двух 0-мерных сфер является одномерной окружностью. $S^1=S^0\oplus S^0$.

сформулируем в виде теоремы следующее предположение.

Теорема 28

*Пусть* $G^m$ *и* $H^n$*-два нормальных молекулярных замкнутых минимальных пространства размерности m и n. Тогда их сумма* $W^{m+n+1}=G^m\oplus H^n$ *также будет минимальным пространством.*

Самостоятельно доказать или опровергнуть.□

До сих пор мы рассматривали нормальные пространства общего вида, не вводя почти никаких ограничений на глобальную и локальную структуру пространства. Однако наибольший практический интерес представляют те пространства, которые моделируют многомерные объекты, возникающие в реальных процессах, обрабатываемых компьютером. Рассмотрим соотношение между нормальным пространством и его нормальным подпространством той же размерности.

Теорема 29

*Пусть $G^n$ и $H^n$-два нормальных молекулярных замкнутых пространства размерности n, и $G^n$ есть подпространство пространства $H^n$ , $G^n \subseteq H^n$,. Тогда $G^n$=$H^n$.*

Доказательство.

Используем индукцию. Для n, m=0, 1 теорема проверяется непосредственно. Предположим, что теорема верна для $G^m$ и $H^n$, где $n \leq p$. Пусть n=p+1. Выберем произвольную точку $v_1 \in G^n$ и рассмотрим окаем этой точки в $H^n$. Очевидно, что $O(v_1)|G^n \subseteq O(v_1)|H^n$, тогда $O(v_1)|G^n$=$O(v_1)|H^n$, так как размерность этих окаемов равна p, и $O(v_1)|H^n$ является связным пространством. Кроме того, так как $H^n$ связно, для любой точки $v_k \in G^n$ и для любой точки $v_p \in H^n$ существует цепь,

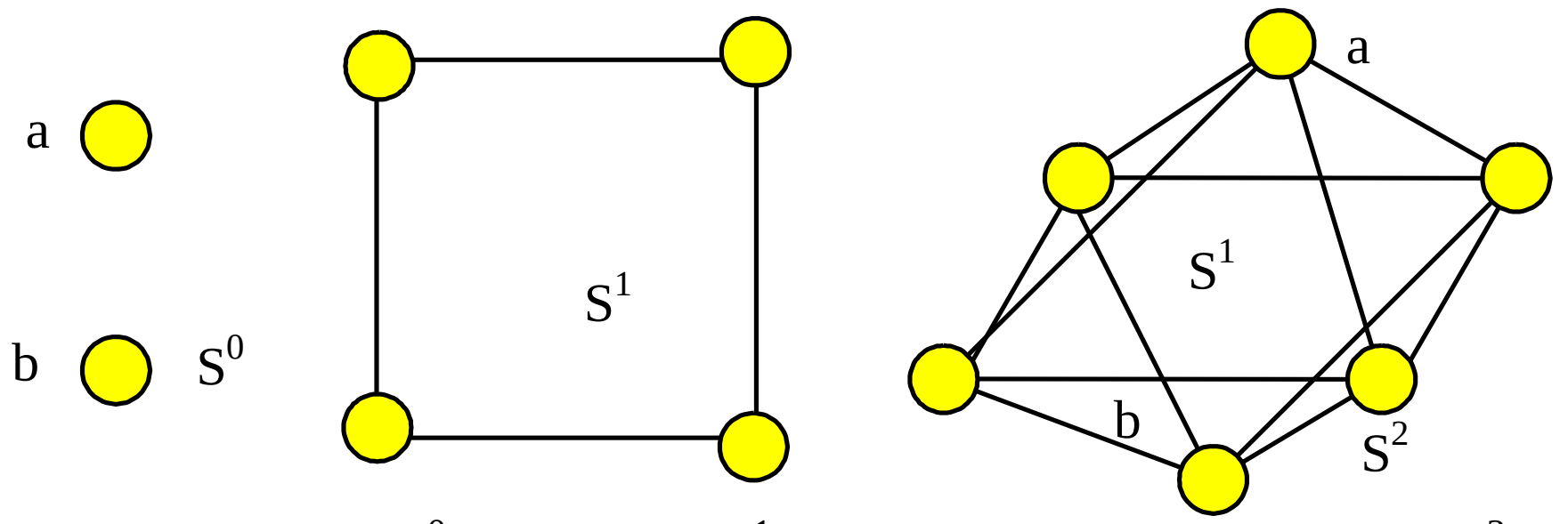

Рис. 49 Сумма 0-сферы $S^0$ и 1-сферы $S^1$ является 2-мерной сферой. $S^2$=$S^0 \oplus S^1$.=$S^1 \oplus S^0$.

связывающая их. Так как для каждой точки этой цепи $O(v_s)|G^n$=$O(v_s)|H^n$ , то $O(v_s)|G^n$=$O(v_s)|H^n$ для всех точек из $H^n$. Это означает, что $G^n$=$H^n$ . Теорема доказана. □

Следствия данной теоремы очевидны. Никакое замкнутое нормальное пространство размерности n не может содержать замкнутое нормальное подпространство той же размерности, кроме самое себя. Например, сфера не может содержать другую сферу, тор или бутылку Клейна. В то же время, если пространство $H^n$ не является нормальным, то оно может содержать подпространство той же размерности, но меньшего объема, как это видно из Рис. 50. Пространство A одномерно, замкнуто, но несвязно. Оно не является нормальным и содержит две одномерные окружности, как подпространства. Пространство B является двумерным, но также не является нормальным, так как окаемы точек a и b являются несвязными окружностями. B содержит два двумерных, замкнутых и нормальных пространства, являющихся двумерными сферами. Одна сфера состоит прямой суммы точек a и b и правой окружности, другая сфера состоит прямой суммы тех же точек a и b и левой окружности. Прямым следствием этой теоремы является следующая теорема, уточняющая структуру специального вида нормального пространства.

Теорема 30

*Пусть $G^n$ есть нормальное молекулярное замкнутое пространство размерности n, в котором существуют две несвязные точки a и b, общий окаем которых есть нормальное молекулярное замкнутое пространство O(ab) размерности (n-1). Тогда $G^n$ есть пряма сумма сферы $S^0$ и O(ab), $G^n=S^0\oplus O(ab)$, , где состоит $S^0$ из двух точек a и b.*

Доказательство.

Рассмотрим подпространство $H=a\cup b\cup O(ab)$, состоящее из точек a, b и O(ab) (Рис. 51). H есть прямая сумма 0-мерной сферы $S^0$, состоящей из точек a и b, и O(ab). Следовательно, H есть нормальное замкнутое подпространство размерности (n+1), раз(H)=раз($S^0$)+раз(O(ab)+1=0+n-1+1=n. По предыдущей теореме $G^n$=H. Теорема доказана. □

Докажем еще одну теорему, показывающую связь между окаемами различных точек нормального пространства. Предположим, что все точки некоторого пространства имеют окаемы, являющиеся нормальными змкнутыми пространствами. Возникает вопрос, имеют ли окаемы этих точек одну и ту же размерность? Оказывается, что имеют.

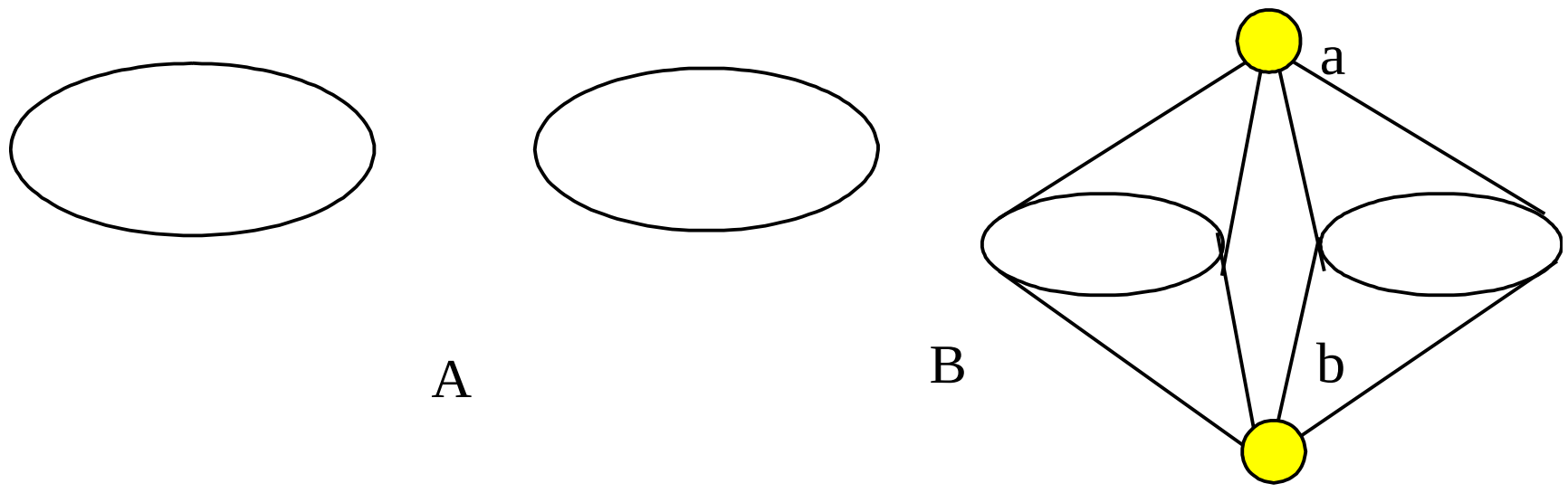

Рис. 50 Пространство A одномерно и несвязно, так как состоит из двух несвязных окружностей $S^1$. Пространство B есть прямая сумма нульмерной сферы $S^0$ и пространства A. Оба этих пространства не принадлежат к классу нормальных.

Теорема 31

*Пусть $G^n$ есть некоторое молекулярное связное пространство, в котором окаем каждой точки $v_k$ есть некоторое нормальное замкнутое пространство размерности $n_k$ . Тогда окаемы всех точек имеют одну и ту же размерность n, и пространство G является (n+1)-мерным нормальным пространством $G^{n+1}$.*

Доказательство.

Выберем две смежные точки $v_1\in G$ и $v_2\in G$ и рассмотрим $O(v_1v_2)$. Так как $O(v_1v_2)\subseteq O(v_1)$, то раз($O(v_1v_2)$)=$n_1$-1, Так как $O(v_1v_2)\subseteq O(v_2)$, то раз($O(v_1v_2)$)=$n_2$-1. Следовательно, $n_1=n_2$ для любых двух смежных точек. Так как пространство связно, то все точки имеют одну и ту же размерность n. пространство G является (n+1)-мерным нормальным пространством

$G^{n+1}$. Теорема доказана. □

С точки зрения теории еще один аспект представляет определенный интерес. Предположим, что имеется некоторое связное молекулярное пространство G, гомотопное некоторому замкнутому нормальному пространству $H^n$. Гомотопность означает, что точечными преобразованиями G может быть стянуто в $H^n$. Будет ли G нормальным n-

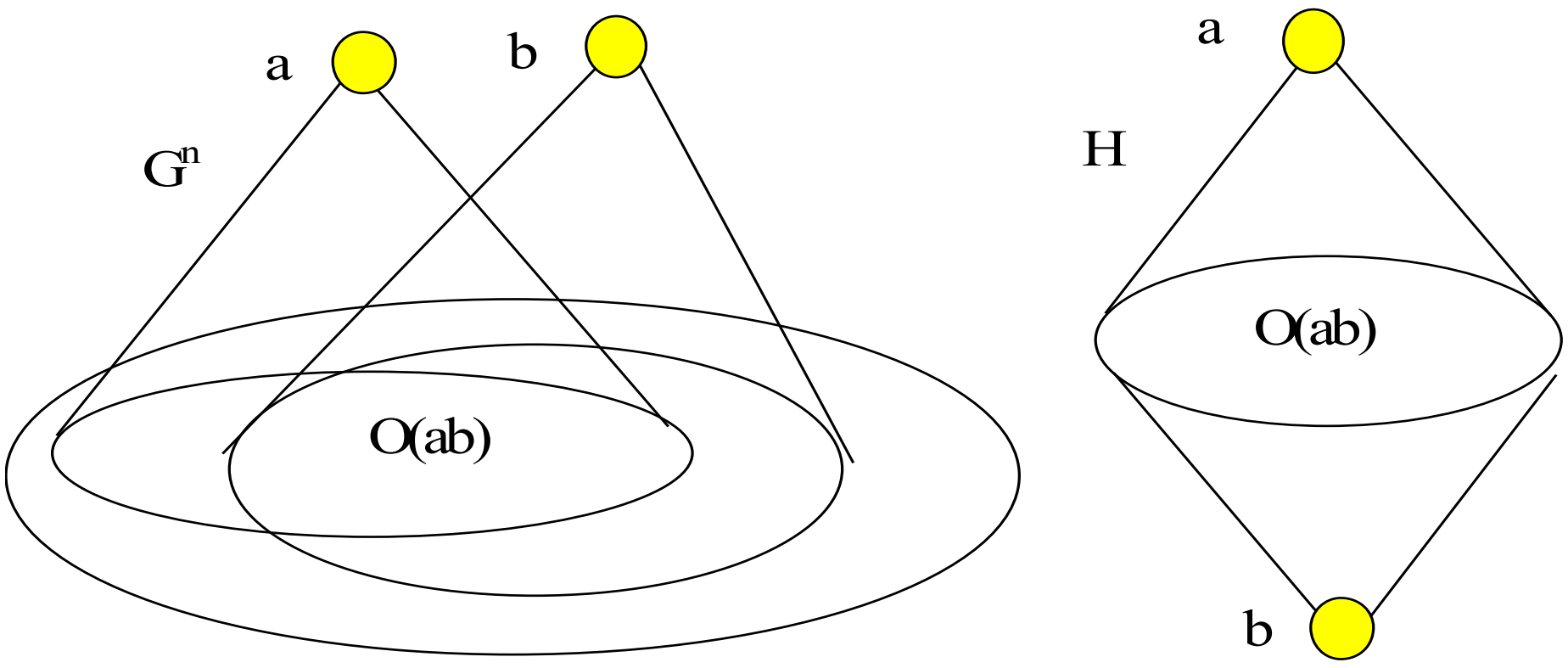

Рис. 51 В нормальном пространстве $G^n$ общий окаем несмежных точек a и d есть (n-1)-мерное нормальное пространство O(ab). Тогда $G^n=H=S^0\oplus O(ab)$ мерным пространством. Очевидно, что в общем случае это не так. Однако, что можно сказать в случае, если G не имеет точечных связей и точек, то есть мы не можем отбросить ни одну точку или связь, следовательно, процесс перехода может быть начат только с приклеивания точки или связи. Во всех рассмотренных конкретных случаях G также являлось замкнутым нормальным пространством. Однако нам не удалось подкрепить этот экспериментальный факт доказательством. Поэтому сформулируем это в виде аксиомы.

А к с и о м а   о   г о м о т о п и и   н о р м а л ь н о м у   п р о с т р а н с т в у

> Если пространство G не имеет точечных связей или точек и гомотопно некоторому n-мерному замкнутому нормальному пространству $H^n$, то G само является n-мерным замкнутым нормальном пространством.

Отметим еще одну любопытную особенность молекулярных пространств, которую мы назвали эвристическим принципом “отсутствие одной особенности” или принципом “ООО”. Было замечено, что в молекулярном пространстве не может быть так, чтобы одна точка отличалась от всех остальных. Иными словами этот принцип может быть сформулирован следующим образом:

Э в р и с т и ч е с к и й   п р и н ц и п   “О О О”

> В молекулярном пространстве с числом точек не меньше 2 свойством А или не обладает ни одна точка, или обладают не меньше двух точек.

С проявлениями такой особенности мы будем встречаться неоднократно на протяжении этой книги. Как проявление принципа отсутствие одной особенности следующая теорема утверждает, что если окаемы всех точек

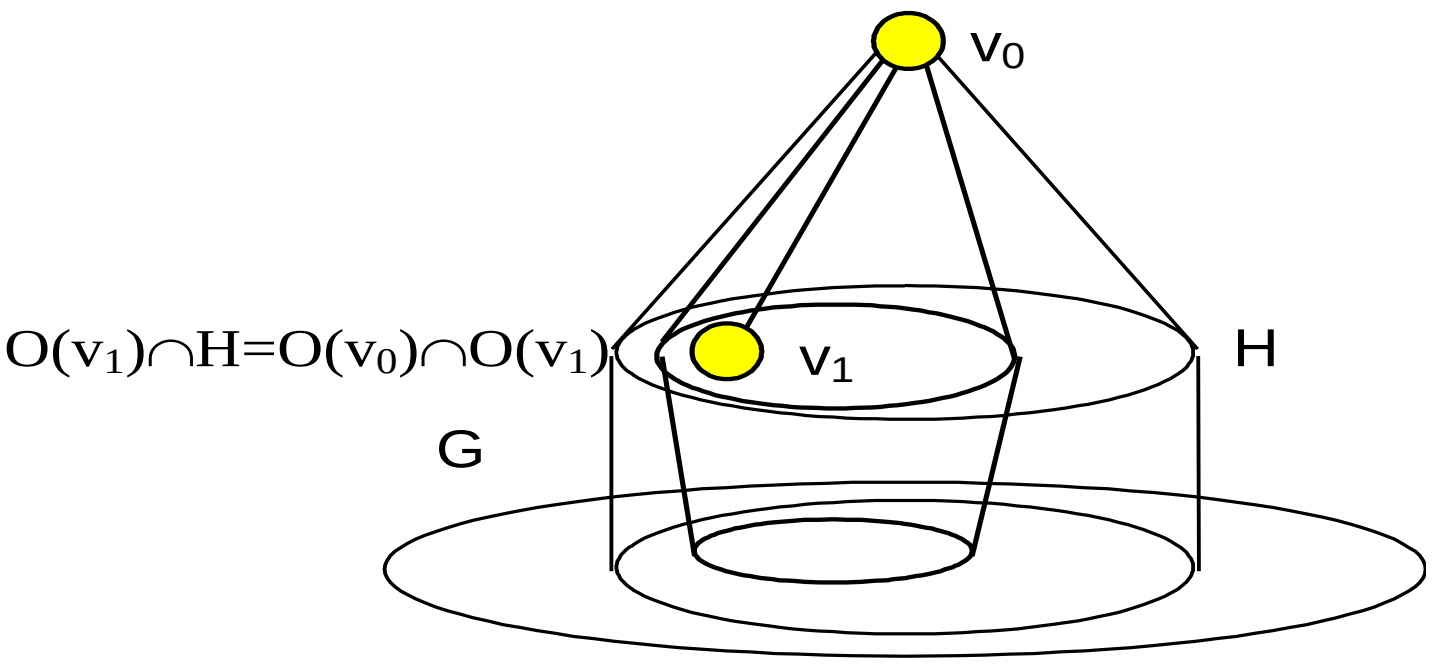

Рис. 52 Окаем $O(v_1)$ произвольной точки $v_1$ в H есть пересечение $O(v_0)\cap O(v_1)$ двух окаемов и, следовательно, (n-2)-мерное нормальное замкнутое пространство. Таким образом, окаем каждой точки в H есть (n-1)-мерное нормальное замкнутое пространство.

некоторого молекулярного пространства, кроме одной, являются n-мерными нормальными пространствами, то и эта последняя точка также n-мерна и молекулярное пространство является нормальным n-мерным.

Теорема 32

*Пусть G есть некоторое молекулярное замкнутое пространство на множестве точек $V=(v_0,v_1,v_2,...v_n)$. Пусть полные окаемы $O(v_i)$ точек $v_i$, $i=1,2,3...n$, есть (n-1)-мерные нормальные молекулярные замкнутые пространства. Тогда полный окаем $O(v_0)$ точки $v_0$ также есть (n-1)-мерное нормальное молекулярное замкнутое пространство, и G есть молекулярное нормальное n-мерное пространство. (Рис. 52).*

Доказательство.

Необходимо показать, что $O(v_0)$ точки $v_0$ есть (n-1)-мерное нормальное молекулярные замкнутое пространство. Обозначим $H=O(v_0)$. Возьмем любую точку из H, например, $v_1$. $O(v_1)$ есть по условию есть (n-1)-мерное нормальное молекулярное замкнутое пространство. Тогда $O(v_0)\cap O(v_1)=O(v_0v_1)=O(v_1)|H$ есть (n-2)-мерное нормальное молекулярные замкнутое пространство. Очевидно, что это справедливо для любой точки, принадлежащей H. Следовательно $H=O(v_0)$ по определению есть (n-1)-мерное нормальное молекулярное замкнутое пространство. Отсюда G есть молекулярное нормальное n-мерное пространство. Теорема доказана. □

## *ГОМЕОМОРФНЫЕ ПРЕОБРАЗОВАНИЯ НОРМАЛЬНЫХ МОЛЕКУЛЯРНЫХ ПРОСТРАНСТВ*

В классической топологии существует понятие гомеоморфизма топологических пространств--это, по определению, непрерывное взаимно-однозначное отображение пространства X на Y, такое что обратное отображение пространства Y на X также непрерывно. Пространства X и Y гомеоморфны, если между ними существует гомеоморфизм. Топологические особенности гомеоморфных пространств совпадают, хотя геометрические характеристики могут быть различны, такие пространства могут отличаться, например, размером или формой. С точки зрения геометрии поверхность шара и поверхность куба различаются, хотя топологически они неразличимы, поскольку поверхность куба может непрерывным образом (без разрывов и склеиваний) быть деформирована в поверхность шара. Нечто подобное необходимо иметь и для молекулярных пространств. Введем преобразования, которые мы назовем гомеоморфными преобразованиями молекулярного пространства, или гм-преобразованиями, и которые будут в некотором роде аналогами гомеоморфизмов.

Нам необходимо для дальнейшего дать несколько однотипных определений, расширяющих уже сделанные ранее определения полного окаема и полного шара точки и подпространства. Предположим, что окаем O(v) точки v имеет определенную топологическую структуру. Возьмем

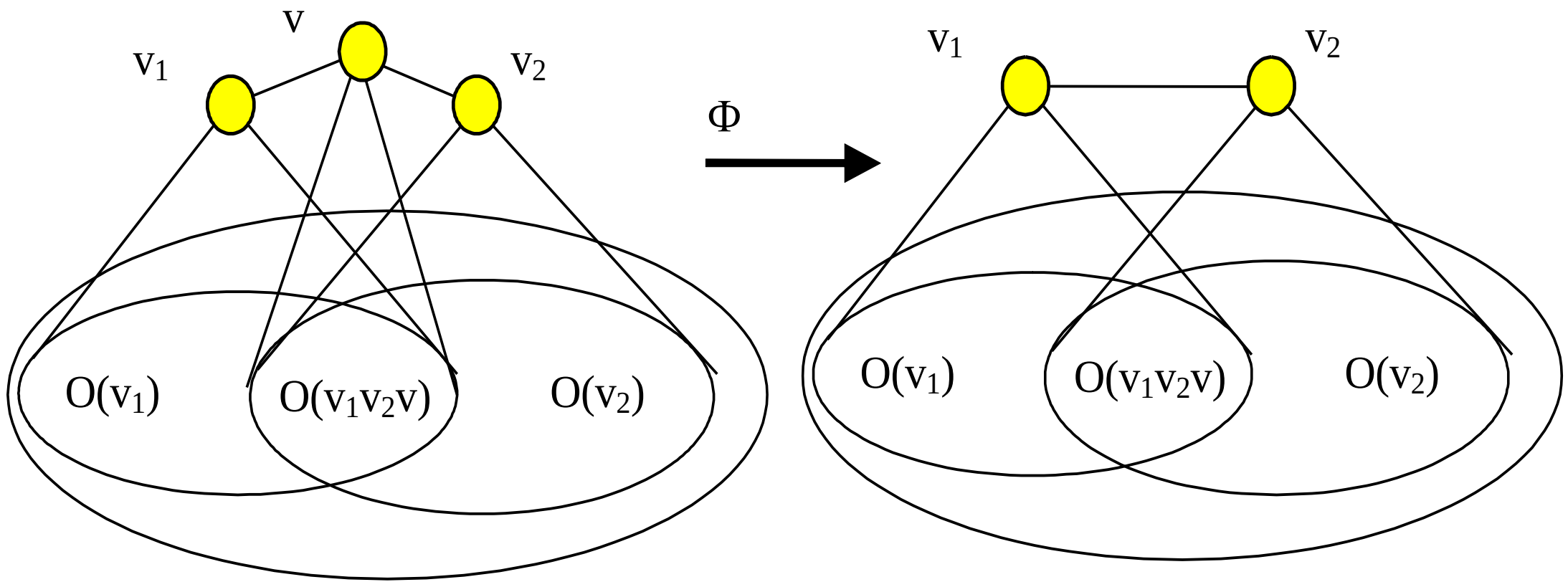

Рис. 53 Гомеоморфные преобразования молекулярного пространства. Точка v меняется на связь ($v_1v_2$). При этом должны соблюдаться условия: $O(v) = S^0(v_1,v_2) \oplus O(v_1v_2v)$, $O(v_1v_2) = v \oplus O(v_1v_2v)$.

окаем O(U(v)) шара U(v) точки. Если структура O(U(v)) такая же как O(v), то мы можем связать точку v с O(U(v)), отбросив O(v). С классической точки зрения это будет выглядеть как небольшое сжатие пространства, уменьшение его объема при сохранении всех основных топологических характеристик.

Следует отметить, что гомеоморфные преобразования молекулярных пространств определяются прежде всего из соображений их практического использования при компьютерной обработке молекулярных моделей.

Поэтому возможны другие, отличные от вводимых ниже, определения

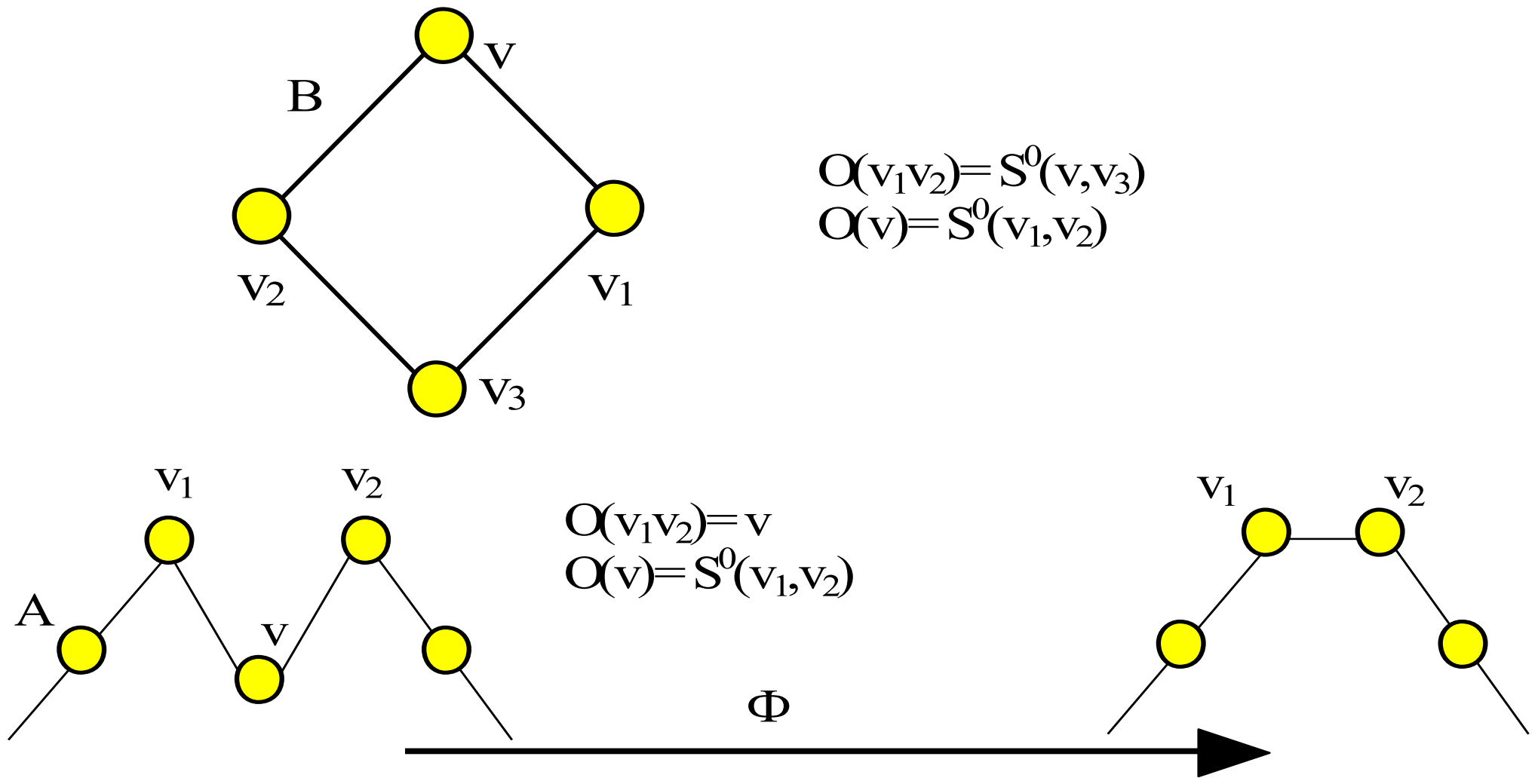

Рис. 54 Замена точки на связь. Условия $O(v) = S^0(v_1,v_2) \oplus O(v_1v_2v)$ и $O(v_1v_2) = v \oplus O(v_1v_2v)$ замены выполняются для одномерного случая внизу и не выполняется для одномерного случая вверху, где замена точки на связь запрещена.

таких преобразований. Однако в любом случае гм-преобразования должны быть точечными преобразованиями.

Определение замены связи на точку (ЗСТ)

Пусть G есть некоторое молекулярное пространство, $v_1$ и $v_2$ есть две точки этого пространства со связью $(v_1v_2)$. Преобразование, состоящее из отбрасывания связи $(v_1v_2)$ и присоединения точки v называется гомеоморфным, если окаем точки v O(v) содержит точки $v_1$, $v_2$ и $O(v_1v_2)$, то есть $O(v)=S^0(v_1,v_2) \oplus O(v_1v_2)$, где сфера $S^0(v_1,v_2)$ состоит из двух точек $v_1$, $v_2$ (Рис. 55).

Определение замены точки на связь (ЗТС)

Пусть G есть некоторое молекулярное нормальное n-мерное пространство, v есть некоторая точка этого пространства, смежная с точками $v_1$ и $v_2$, которые несмежны между собой. Преобразование, состоящее из отбрасывания точки v и установления связи $(v_1v_2)$ называется гомеоморфным, если $O(v)=S^0(v_1,v_2) \oplus O(v_1v_2v)$ и $O(v_1v_2)=v \oplus O(v_1v_2v)$, где сфера $S^0$ состоит из двух несвязных друг с другом точек $v_1$, $v_2$.

Легко видеть, что первое преобразование увеличивает число точек пространства на 1, второе преобразование уменьшает число точек на 1.

Эти преобразования являются в некотором смысле взаимообратными, отбрасывание связи и приклеивание точки означает, что в обратном направлении эту точку можно заменить на отброшенную связь. Гомеоморфные преобразования (ГМ) обладают важной особенностью; они сохраняют размерность нормальных пространств. С помощью гомеоморфных преобразований возможно изменение объема пространства (напомним, что под объемом понимается число точек, формирующих пространство).

З а м е ч а н и е.

Вообще говоря эти преобразования в общем виде применимы не только к нормальным пространствам, но к любым молекулярным пространствам. При этом замена связи на точку осуществима без всяких ограничений, замена точки на связь определяется окаемом точки v. Преобразования такого вида будут использоваться в регулярных пространствам, расширяющих понятие нормальных пространств, а также в точечных преобразованиях молекулярных

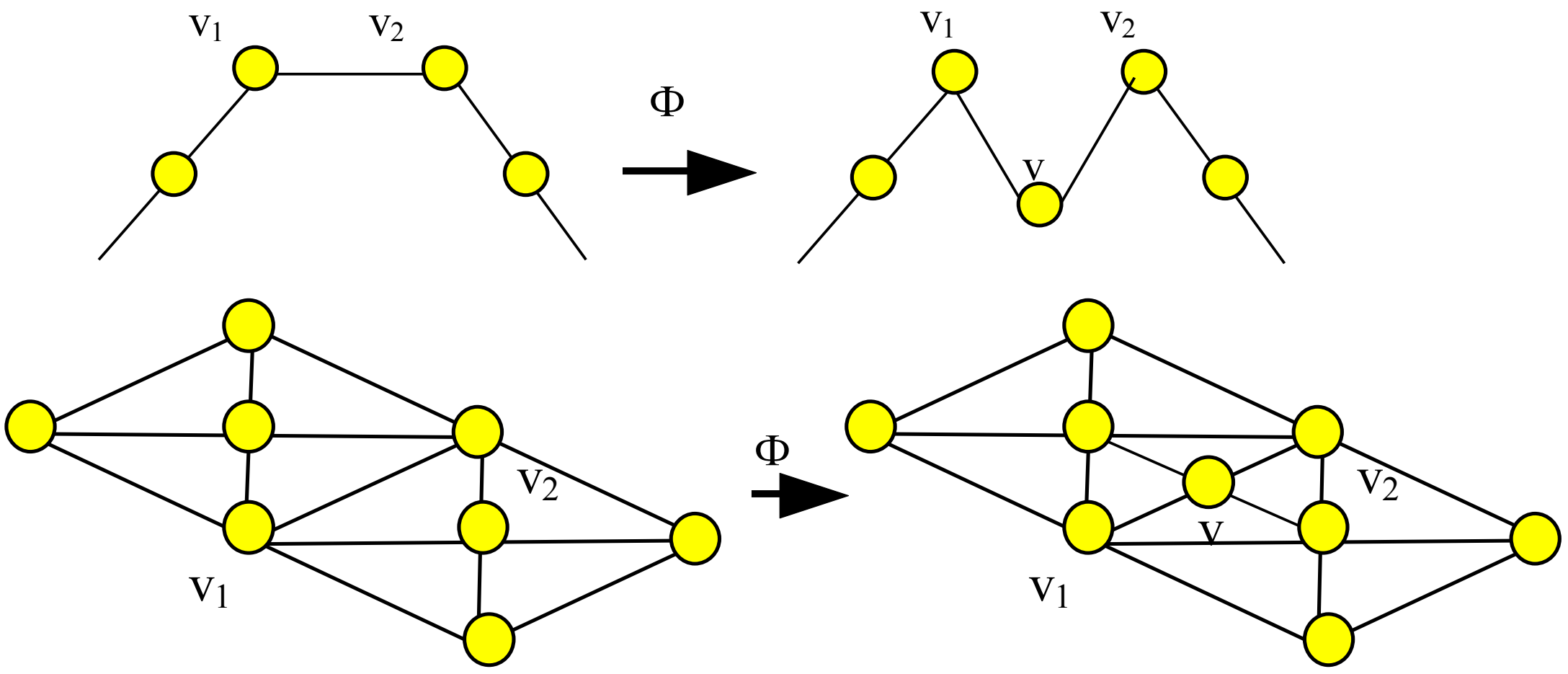

Рис. 55 Замена связи на точку для одномерного (вверху) и двумерного (внизу) случаев. Связь ($v_1v_2$) отбрасывается и заменяется на точку v.

пространств.

На Рис. 55 связь ($v_1v_2$) отбрасывается и заменяется на точку v. При этом точка v соединяется в одномерном случае с точками $v_1$ и $v_2$, поскольку $O(v_1v_2)$ есть пустое подпространство. В двумерном случае точка v соединяется с точками $v_1$ , $v_2$, а также с общим окаемом $O(v_1v_2)$ этих точек (внизу). Легко видеть, что окаемы всех точек после преобразования остаются или 0-мерными сферами для одномерного случая, или окружностями для двумерного случая.

Более сложным и не всегда возможным является замена точки на связь. В некоторых случаях эту замену невозможно произвести. Как скажем, в случае одномерной окружности, состоящей из четырех точек, мы не можем уменьшить число точек до трех, поскольку топология пространства, в частности, размерность, изменится. На Рис. 57 показана замена точки на связь при выполнении условий $O(v) = S^0(v_1,v_2)\oplus O(v_1v_2v)$ и $O(v_1v_2) = v\oplus O(v_1v_2v)$ для одномерного случая. Внизу точка v заменяется на связь $(v_1v_2)$, поскольку выполняются эти условия; вверху точка v не может быть заменена на связь $(v_1v_2)$, поскольку не выполняется условие $O(v_1v_2) = v\oplus O(v_1v_2v)$. В двумерном случае (Рис. 56) точка v заменяется на связь

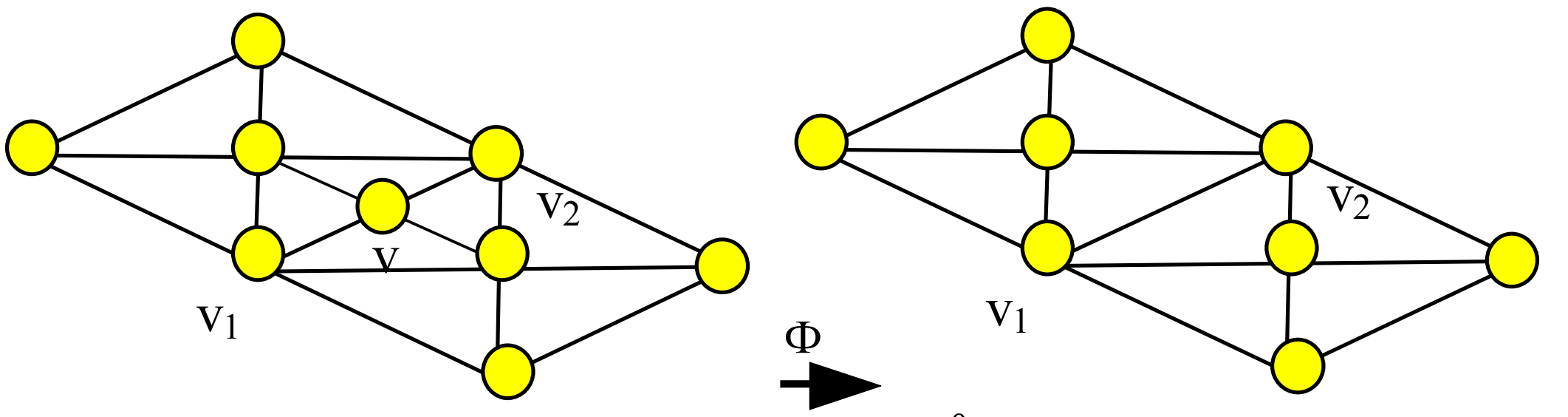

Рис. 56 Замена точки на связь. Условия $O(v) = S^0(v_1,v_2)\oplus O(v_1v_2v)$ и $O(v_1v_2) = v\oplus O(v_1v_2v)$ замены выполняются для двумерного пространства слева. Точка v точка отбрасывается, и устанавливается связь $(v_1v_2)$ в пространстве справа.

$(v_1v_2)$, поскольку оба условия выполняются.

Определение обратного гомеоморфного преобразования. XE "

Пусть G есть некоторое молекулярное пространство, $v_1$ и $v_2$ есть две точки этого пространства со связью $(v_1v_2)$. Гомеоморфное преобразование Ф замена связи на точку, состоящее из отбрасывания связи $(v_1v_2)$ и присоединения точки v называется обратным к гомеоморфному преобразованию $Ф^{-1}$ замена точки на связь, состоящему из отбрасывания точки v и присоединения связи $(v_1v_2)$. Очевидно, что для каждого гомеоморфного преобразования Ф существует обратное ему гомеоморфное преобразование $Ф^{-1}$ такое, что $Ф^{-1}ФG==G$.

Определение гомеоморфных пространств

Пространство H называется гомеоморфным пространству G, $H\approx G$, если H может быть получено из G последовательным применением конечного числа гомеоморфных преобразований, $H=Ф_n \ldots Ф_2Ф_1G$.

Все нормальные пространства, гомеоморфные данному, образуют класс гомеоморфных молекулярных пространств. Введем обозначение Г(G) для множества всех пространств, гомеоморфных G, тогда $H\in Г(G)$. Очевидно, что множество всех нормальных молекулярных пространств разбивается на непересекающиеся классы гомеоморфных пространств.

Из рисунков видно, что оба гомеоморфных преобразования-замена связи на точку и замена точки на связь, не меняют размерности пространства. Теперь докажем это сначала для первого преобразования-замена связи на точку.

Теорема 33

*Гомеоморфное преобразование ЗСТ не меняет размерность и нормальность нормального молекулярного пространства $G^n$.*

Доказательство.

Используем индукцию. Для малых размерностей теорема проверяется непосредственно. Предположим, что теорема верна для размерности $n < k$. Пусть $n = k$. Пусть $G^n$ есть некоторое молекулярное пространство, $v_1$ и $v_2$ есть две точки этого пространства со связью $(v_1v_2)$ (Рис. 55). Согласно свойствам нормального пространства $O(v_1v_2)$ есть нормальное замкнутое (n-2)-мерное пространство. Удалим связь $(v_1v_2)$ и соединим точку v с точками со всеми точками из $O(v_1v_2)$, а также с $v_1$ и $v_2$.

Легко видеть, что окаем $O(v_1)$ точки $v_1$, не изменился кроме того, что точка $v_2$ в этом окаеме заменилась на точку v. Сходным образом не изменился окаем $O(v_2)$, точки $v_2$. Следовательно, окаемы этих точек остались (n-1)-мерными замкнутыми нормальными пространствами после замены связи на точку.

Окаем O(v) точки v имеет структуру суммы 0--мерной сферы, состоящей из двух несвязных точек $v_1$ и $v_2$ и (n-2)-мерного пространства $O(v_1v_2)$, $O(v)=S^0\oplus O(v_1v_2)$. Следовательно, O(v) есть нормальное замкнутое (n-1) пространство (Рис. 55). Пусть точка $v_3$ лежит в $O(v_1v_2)$. $O(v_3)$ есть нормальное замкнутое подпространство размерностью (n-1) (Рис. 55). В нем происходит точно такое же преобразование замена связи на точку. Согласно предположению, после этого преобразования пространство остается нормальным, замкнутым и той же размерности. Все остальные точки пространства $G^n$ своих окаемов не меняют. Таким образом, окаем каждой точки остается замкнутым нормальным (n-1)-мерным пространством. Теорема доказана. □

Преобразование замена точки на связь есть преобразование, обратное замене связи на точку. Почти очевидно, что если в исходном пространстве $O(v) = S^0(v_1,v_2)\oplus O(v_1v_2v)$ и $O(v_1v_2) = v\oplus O(v_1v_2v)$, то точка v может быть заменена на связь $(v_1v_2)$. Тем не менее докажем, что это справедливо.

Теорема 34

*Гомеоморфное преобразование ЗТС не меняет размерность и нормальность нормального молекулярного пространства $G^n$.*

Доказательство.

Используем индукцию. Для размерности n=1 теорема проверяется непосредственно (Рис. 54). Предположим, что теорема верна для размерности $n < k$. Пусть $n = k$. Пусть $G^n$ есть некоторое молекулярное пространство, v есть какая-то точка этого пространства, для которой $O(v) = S^0(v_1,v_2)\oplus O(v_1v_2v)$ и $O(v_1v_2) = v\oplus O(v_1v_2v)$, где $S^0$ есть две несвязные точки $v_1$ и $v_2$ из U(v). Удалим точку v и установим связь $(v_1v_2)$. Легко видеть, что окаемы всех точек, не лежащих в $O(v_1v_2v)$, не изменились. Следовательно, окаемы этих точек остались (n-1)-мерными замкнутыми нормальными пространствами после замены точки на связь. Рассмотрим окаем $O(v_3)$

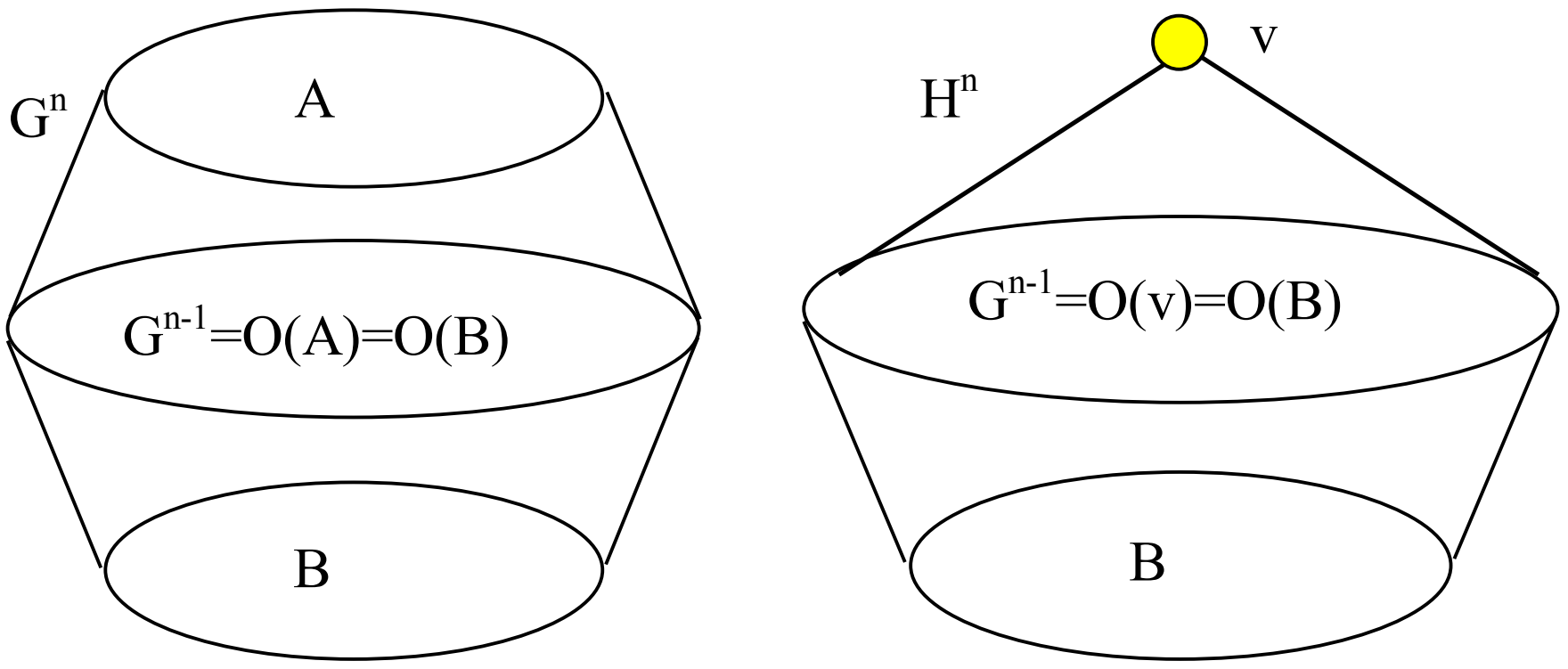

Рис. 57 Полные окаемы подпространств A и B есть замкнутое нормальное (n-1)-мерное подпространство $G^{n-1}$. Отбрасывание A и соединение точки v со всеми точками из $G^{n-1}$ не меняет размерность пространства.

точки $v_3$ лежащей в $O(v_1v_2v)$. Этот окаем является нормальным замкнутым (n-1)-мерным пространством и имеет структуру исходного пространства $G^n$, то есть $O(v)|O(v_3) = S^0(v_1,v_2)\oplus O(v_1v_2v)|O(v_3)$ и $O(v_1v_2)|O(v_3) = v\oplus O(v_1v_2v)|O(v_3)$, где $S^0$ есть две несвязные точки $v_1$ и $v_2$ из $O(v)|O(v_3)$. В $O(v_3)$ происходит точно такое же преобразование замена связи на точку. Согласно предположению индукции, после этого преобразования подпространство $O(v_3)$ остается нормальным, замкнутым и той же размерности (n-1). Таким образом окаем каждой точки остается замкнутым нормальным (n-1)-мерным пространством. Теорема доказана. □

Внимательный взгляд на эти два преобразования позволяет заметить, что, вообще говоря, окаемы точек могут менять свой тип. Например, при преобразовании замена связи на точку может появиться окаем, топологически не свойственный данному пространству. Однако, в одном случае топологические свойства окаемов при гомеоморфных преобразованиях не меняются. Это происходит, если все окаемы сферические. В самом деле, умножение сферы любой размерности на 0-сферу всегда дает сферу. Иными словами при гомеоморфных преобразованиях молекулярное многообразие остается многообразием.

Перейдем теперь к прямой сумме двух молекулярных пространств. Пусть $\Phi_G$ есть гомеоморфное преобразование замена связи на точку на $G^m$. Очевидно, что оно будет также гомеоморфным преобразованиям на $W^{m+n+1}=G^m \oplus H^n$, поскольку такое преобразование применимо к любой связи в нормальном пространстве. Несколько сложнее убедиться в том, что гомеоморфное преобразование замена точки на связь на $G^m$ будет также гомеоморфным на $W^{m+n+1} = G^m \oplus H^n$. Докажем это.

Теорема 35

*Пусть $G^m$ и $H^n$-два нормальных молекулярных замкнутых пространства размерности m и n. Пусть $\Phi_G$ есть некоторое гомеоморфное преобразование замена связи на точку на $G^m$. Тогда $\Phi_G$ есть также гомеоморфное преобразование замена связи на точку на $G^m \oplus H^n$.*

Доказательство.

Пусть $G^m$ есть замкнутое нормальное пространство, $v_1$ и $v_2$ есть две точки этого пространства без связи, и для точки v выполняются условия $O(v) = S^0(v_1,v_2) \oplus O(v_1v_2v)$ и $O(v_1v_2) = v \oplus O(v_1v_2v)$, где $A^{m-2} = O(v_1v_2v)$ есть нормальное замкнутое молекулярное пространство размерности (m-2). Заменим теперь точку v на связь $(v_1v_2)$ гомеоморфным образом. При этом $G^m$ переходит в $\Phi g^m$. Рассмотрим $G^m \oplus H^n$. Выберем ту же точку v те же две точки $v_1$ и $v_2$ без связи из подпространства $G^m$. При этом $O(v) = S^0(v_1,v_2) \oplus (A^{m-2} \oplus H^n)$ и $O(v_1v_2) = v \oplus (A^{m-2} \oplus H^n)$, где $(A^{m-2} \oplus H^n) = O(v_1v_2v)$ в $G^m \oplus H^n$. Следовательно в $G^m \oplus H^n$ можно гомеоморфно отбросить точку v и установить связь между точками $v_1$ и $v_2$. Это завершает доказательство теоремы. □

Рассмотрим случай, когда в прямой сумме двух молекулярных пространств одно из этих пространств подвергается гомеоморфному преобразованию Ф. Пусть $W^{m+n+1} = G^m \oplus H^n$ и $W_1^{m+n+1} = \Phi_G(G^m) \oplus H^n$. Возникает вопрос, будут ли гомеоморфны $W^{m+n+1}$ и $W_1^{m+n+1}$?

Теорема 36. О гомеоморфизме прямой суммы пространств.

*Пусть $G^m$ и $H^n$-два нормальных молекулярных замкнутых пространства размерности m и n. Пусть Ф есть некоторое гомеоморфное преобразование пространства $G^m$. Тогда прямая сумма $G^m \oplus H^n$ гомеоморфна прямой сумме $\Phi(G^m) \oplus H^n$, $G^m \oplus H^n \approx \Phi(G^m) \oplus H^n$.*

Доказательство.

Пусть $G^m$ есть некоторое молекулярное нормальное пространство, $v_1$ и $v_2$ есть две точки этого пространства со связью $(v_1v_2)$. Обозначим $O(v_1v_2)= A^{m-2}$. Заменим связь $(v_1v_2)$ на точку v гомеоморфным образом. При этом окаем точки v имеет вид $O(v)=S^0(v_1,v_2)\oplus A^{m-2}$ и $G^m$ переходит в $\Phi G^m$. В пространстве $\Phi(G^m)\oplus H^n$ окаем точки v имеет вид $O(v)= S^0(v_1,v_2)\oplus(A^{m-2}\oplus H^n)$. Рассмотрим $G^m\oplus H^n$. Выберем те же две точки $v_1$ и $v_2$ из подпространства $G^m$ со связью $(v_1v_2)$ Очевидно, что в $G^m\oplus H^n$ $O(v_1v_2)=A^{m-2}\oplus H^n$. Гомеоморфным образом удалим связь $(v_1v_2)$ и соединим точку v со всеми точками из $O(v_1v_2)$, а также с $v_1$ и $v_2$. При этом окаем точки v имеет вид $O(v)=S^0(v_1,v_2) \oplus(A^{m-2}\oplus H^n)$ в $\Phi(G^m \oplus H^n)$. Так как шары точки v в $\Phi(G^m)\oplus H^n$ и $\Phi(G^m\oplus H^n)$ совпадают, окаемы всех остальных точек также совпадают, это означает, что $\Phi(G^m)\oplus H^n=\Phi(G^m\oplus H^n)$. Так как $G^m\oplus H^n\approx\Phi(G^m\oplus H^n)= \Phi(G^m)\oplus H^n$, то $G^m\oplus H^n\approx\Phi(G^m)\oplus H^n$. Для преобразования замена связи на точку теорема доказана.

Пусть $G^m$ есть то же самое молекулярное нормальное пространство, $v_1$ и $v_2$ есть две точки этого пространства без связи и для точки v выполняются условия $O(v) = S^0(v_1,v_2)\oplus O(v_1v_2v)$ и $O(v_1v_2) = v\oplus O(v_1v_2v)$, где $A^{m-2} =O(v_1v_2v)$ есть нормальное замкнутое молекулярное пространство размерности (m-2). Заменим теперь точку v на связь $(v_1v_2)$ гомеоморфным

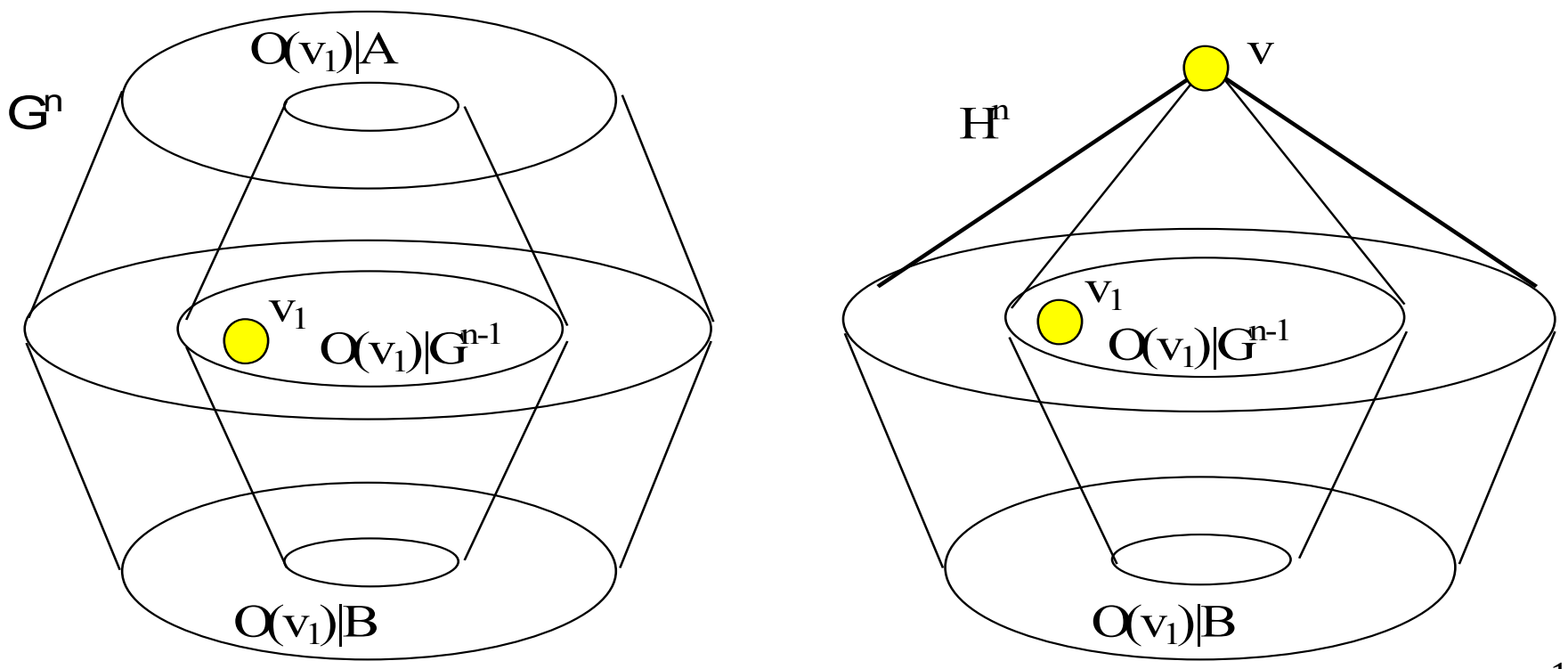

Рис. 58 Отбрасывание A и соединение точки v со всеми точками из $G^{n-1}$ не меняет размерность пространства. Пространство $O(v_1)$ имеет такую же структуру как $G^n$.

образом. При этом $G^m$ переходит в $\Phi g^m$. Рассмотрим $G^m\oplus H^n$. Выберем ту же точку v те же две точки $v_1$ и $v_2$ без связи из подпространства $G^m$. При этом $O(v) = S^0(v_1,v_2)\oplus(A^{m-2}\oplus H^n)$ и $O(v_1v_2) = v\oplus(A^{m-2}\oplus H^n)$, где $(A^{m-2}\oplus H^n)= O(v_1v_2v)$ в $G^m\oplus H^n$. Следовательно в $G^m\oplus H^n$ можно гомеоморфно отбросить точку v и установить связь между точками $v_1$ и $v_2$. При этом пространство $G^m\oplus H^n$ гомеоморфно переходит в $\Phi(G^m\oplus H^n)$. С другой стороны пространство $\Phi(G^m\oplus H^n)$ по-прежнему является прямой суммой двух подпространств, одно из которых, $H^n$, осталось без изменений, поскольку все преобразования происходили над точками подпространства $G^m$ . То есть $\Phi(G^m\oplus H^n)$ можно представить как прямую сумму $(\Phi G^m)\oplus H^n$. Отсюда $G^m\oplus H^n\approx\Phi(G^m\oplus H^n)=(\Phi G^m)\oplus H^n$. Следовательно, $G^m \oplus H^n \approx (\Phi G^m) \oplus H^n$. Для преобразования замена точки на связь теорема доказана. Это завершает доказательство теоремы. □

Таким образом, мы видим, что любое гомеоморфное преобразование на $G^m$ (или $H^n$) автоматически индуцирует соответствующее гомеоморфное преобразование на $G^m \oplus H^n$. С другой стороны не всякое гомеоморфное преобразование на $G^m \oplus H^n$ индуцируется каким-либо преобразованием на $G^m$ или $H^n$. Это происходит только в том случае, если точки $v_1$, $v_2$ и v, участвующие в преобразовании, принадлежат одному из подпространств, или $G^m$, или $H^n$, одновременно. Для этого случая справедливо соотношение $(\Phi G^m) \oplus H^n = \Phi(G^m \oplus H^n)$.

В практической работе часто возникает необходимость уменьшить объем молекулярного пространства, чтобы снизить время, затрачиваемое на обработку этого пространства для различных целей. Поэтому встает вопрос, можно ли отбрасывать точки не по одной, а целыми группами, так чтобы сохранялись определенные характеристики пространства. Частично ответить на этот вопрос позволяет следующая теорема, которая сохраняет, по крайней мере, размерность и нормальность пространства. В других главах мы еще вернемся к этому вопросу. В следующей теореме мы рассмотрим специальный вид молекулярных пространств представленных на Рис. 57.

Теорема 37

*Пусть молекулярное нормальное n-мерное пространство $G^n$ может быть представлено как объединение трех непустых подпространств A, B и $G^{n-1}$, где между точками пространств A и B нет связей, окаемы подпространств A и B есть замкнутое нормальное (n-1)-мерное пространство $G^{n-1}$, $G^n = A \cup G^{n-1} \cup B$, $A \cap O(B) = O(A) \cap B = \varnothing$, $O(A) = O(B) = G^{n-1}$ (Рис. 57). Отбросим все точки подпространства A и соединим изолированную точку v со всеми точками подпространства $G^{n-1}$. Полученное пространство $H^n$ является нормальным молекулярным (n-1)-мерным пространством.* □

Доказательство.

Используем индукцию. Для размерности n=1 теорема проверяется непосредственно. Предположим, что теорема верна для размерности $n < k$. Пусть $n = k$. Рассмотрим точку $v_1$ пространства $G^n$, лежащую в $G^{n-1}$. Окаем $O(v_1)$ этой точки в $G^n$ является нормальным молекулярным замкнутым (n-1)-мерным пространством, окаем $O(v_1)|G^{n-1}$ этой точки в $G^{n-1}$ является нормальным молекулярным замкнутым (n-2)-мерным пространством, между точками подпространств $O(v_1)|A$ и $O(v_1)|B$ нет связей. Таким образом $O(v_1)$ имеет ту же структуру, что и $G^n$. Следовательно, по индукционному предположению при отбрасывании $O(v_1)|A$ и приклеивании точки v к $O(v_1)|G^{n-1}$ $O(v_1)$ перейдет в нормальное замкнутое (n-1)-мерное пространство. Это означает, что окаемы всех точек в $G^{n-1}$ останутся замкнутыми нормальными (n-1)-мерными пространствами. Окаемы точек из B не изменились. Окаем точки $v_1$ также является нормальным молекулярным замкнутым (n-1)-мерным пространством $G^{n-1}$. Отсюда пространство $H^n$ является нормальным молекулярным замкнутым n-мерным пространством. Теорема доказана. □

Рассмотрим на нескольких конкретных примерах как работает эта теорема. На Рис. 59 вверху изображена замена подпространства A на точку v в окружности $S^1{}_1$, состоящей из 6-и точек. При этом полученное пространство $S^1{}_2$ так же является окружностью на 4-х точках. Внизу

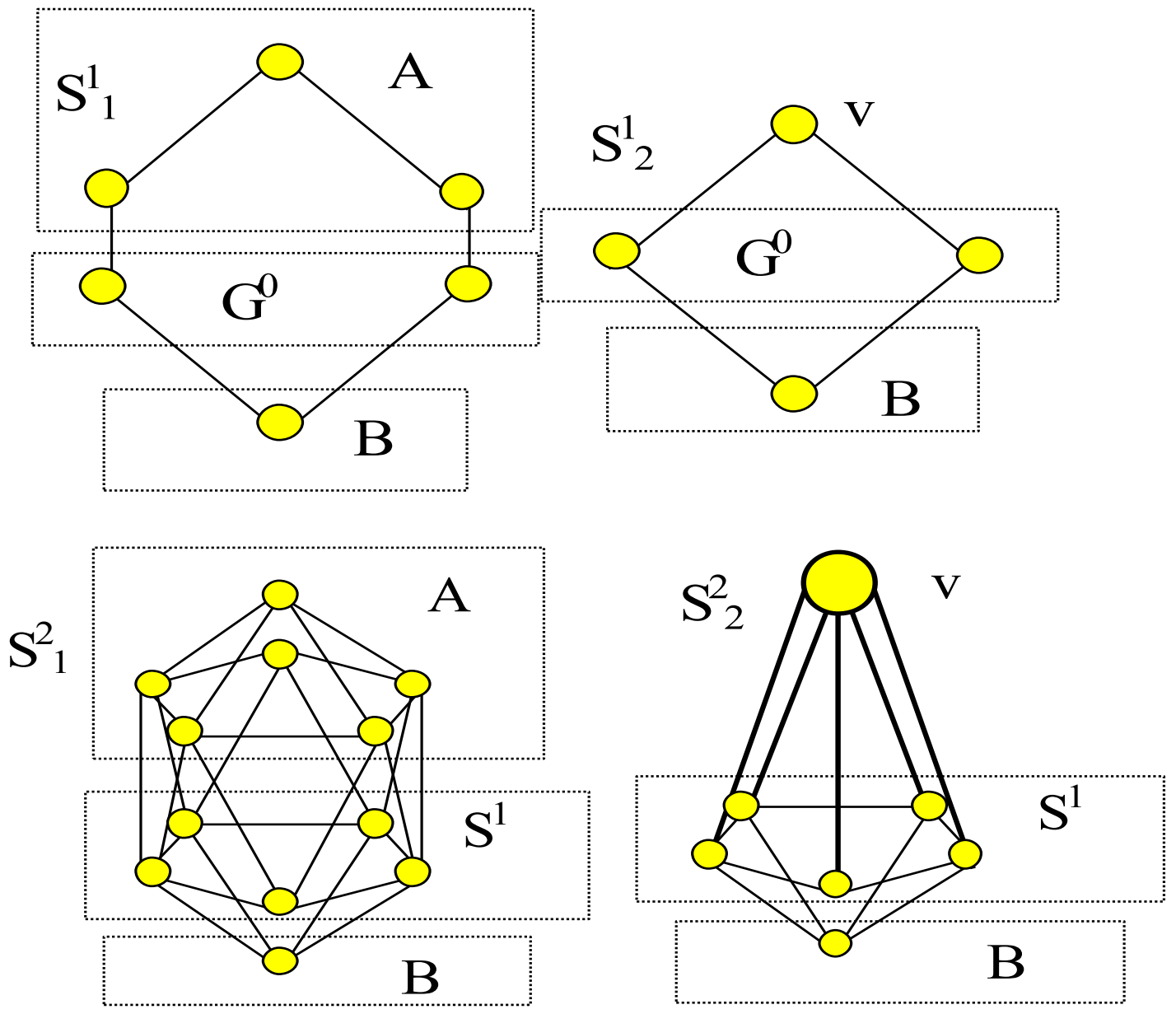

Рис. 59 Замена подпространства A на точку v в окружности и сфере.

изображена сфера $S^2{}_1$, состоящая из 12-и точек. Пространство A состоит из 6-и точек, пространство B-из 1-й точки. При замене A на точку v полученная сфера $S^2{}_2$ состоит из 6-и точек.

Особенно успешно работает эта теорема в тех случаях, когда молекулярное пространство является дигитальным образом одно, двух или трехмерной непрерывной поверхности, то есть многообразием. Как правило, такие образы состоят из значительного количества точек, что увеличивает время обработки таких поверхностей в компьютерах и занимает большой объем машинной памяти. Поэтому удобнее вначале уменьшить по возможности количество точек дигитальной модели. В таких случаях удобно разделять пространство на полные n-окаемы некоторой точки v, $O^1(v)$, $O^2(v)$, $O^3(v)$,...

$O^{p-1}(v)$, $O^{p}(v)$, B (Рис. 60). Последовательно проверяя, являются ли они (n-1)-мерными замкнутыми молекулярными пространствами, мы выбираем последний в серии, который обозначаем как $G^{p-1}$. Объединение всех предыдущих полных окаемов вместе с точкой v будет составлять подпространство A, подпространство, образованное точками, не входящими в A и $G^{n-1}$ является подпространством B. Заменяя A точкой в соответствии с теоремой 3, мы получаем пространство меньшим количеством точек, гомеоморфное первоначальному.

В одном из предыдущих разделов мы ввели точечные преобразования молекулярных пространств. Эти преобразования применимы к произвольным молекулярным пространствам. Покажем, что эти гомеоморыные преобразования являются последовательностью точечных преобразований и применимы к любым пространствам, включая нормальные.

Теорема 38

*Гомеоморфное преобразование замена связи на точку является последовательностью двух точечных преобразований: точечного приклеивания точки и точечного отбрасывания связи в молекулярном пространстве $G^n$.*

Доказательство.

Используем определение замены связи на точку. Вначале присоединим точку v к $G^n$ таким образом, что $O(v)=v_1\oplus v_2\oplus O(v_1v_2)$. O(v), очевидно, является точечным подпространством и, следовательно, присоединение точки v является точечным приклеиванием. Рассмотрим, теперь, общий окаем $O(v_1v_2)$ в $G^n\cup v$. Он имеет вид $O(v_1v_2)|(G^n\cup v)=v\oplus O(v_1v_2)|G^n$. Этот окаем есть конус, и он является точечным пространством. Следовательно, связь $O(v_1v_2)$ может быть отброшена. Теорема доказана. □

Аналогично доказывается теорема о замене точки на связь.

Теорема 39

*Гомеоморфное преобразование замена точки на связь является последовательностью двух точечных преобразований: точечного приклеивания связи и точечного отбрасывания точки в молекулярном пространстве $G^n$.*

Доказательство.

Используем определение замены точки на связь. Введем связь $v_1v_2$ между точками $v_1$ и $v_2$, которые несмежны между собой. Это есть точечное преобразование, так как $O(v_1v_2) = v\oplus O(v_1v_2v)$ является конусом. В полученном пространстве $G^n\cup(v_1v_2)$ окаем точки v есть O(v) $=v_1 \oplus v_2 \oplus O(v_1 v_2 v)$, то есть точечное пространство. Следовательно, точка v может быть отброшена. Теорема доказана. □

## Список литературы к главе 4.

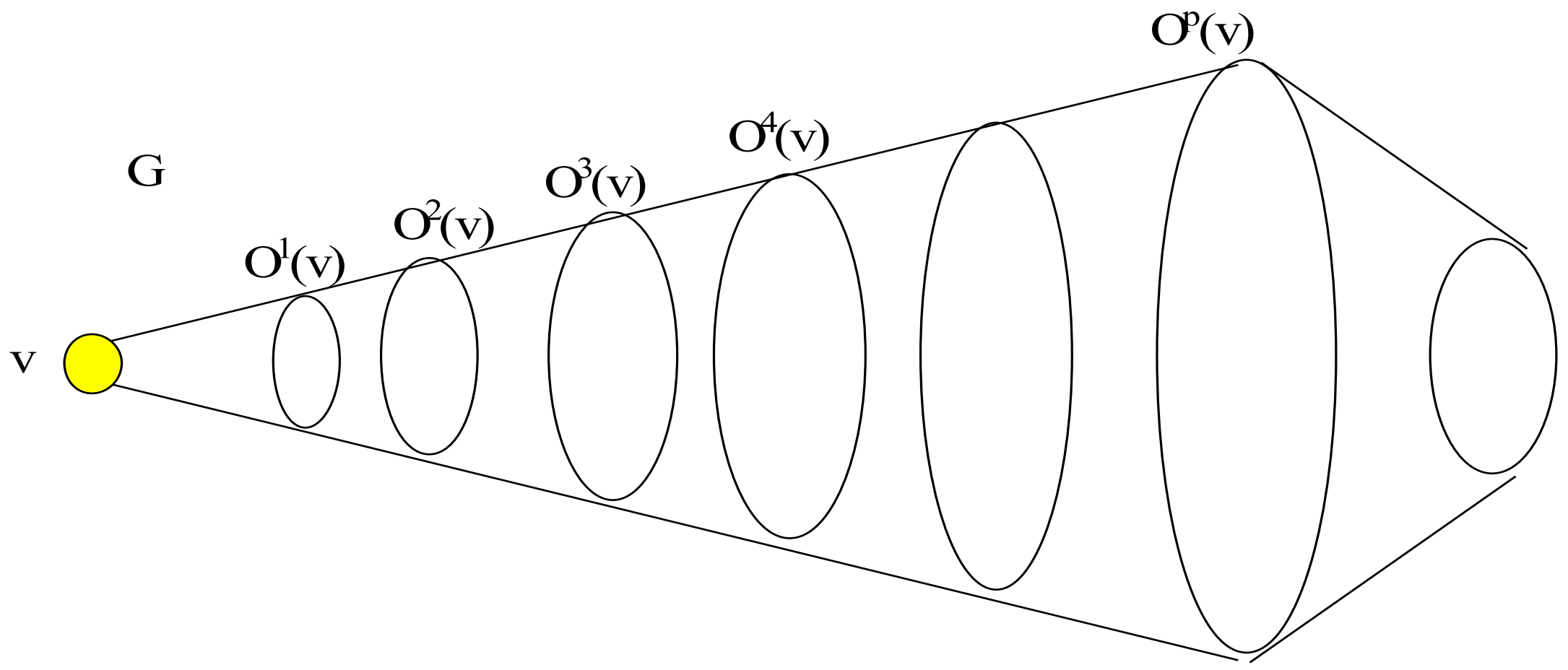

Рис. 60 Пространство G образовано последовательностью полных окаемов $O^1(v)$. $O^2(v)$, $O^3(v)$.... $O^p(v)$ являющихся (n-1)-мерными нормальными замкнутыми пространствами. Объединение всех полных окаемов до (p-1)-полного окаема включительно заменяется на точку v, соединенную со всеми точками пространства $O^n(v)$.

# НОРМАЛЬНЫЕ МОЛЕКУЛЯРНЫЕ ПРОСТРАНСТВА С КРАЕМ И НОРМАЛЬНЫЕ МОЛЕКУЛЯРНЫЕ МНОГООБРАЗИЯ

We study some properties of normal molecular spaces with an adge and normal manyfolds.

Некоторые результаты этой главы имеются в работах [43].

## *НОРМАЛЬНЫЕ МОЛЕКУЛЯРНЫЕ ПРОСТРАНСТВА С КРАЕМ И ИХ СВОЙСТВА*

В предыдущих разделах мы рассматривали в основном пространства, в которых каждая точка является нормальной и n-мерной. Эти пространства являются или конечными, замкнутыми, или же бесконечными нормальными пространствами. Рассмотрим теперь нормальные пространства с краем. В непрерывном пространстве размерности n край может быть связным или состоящим из нескольких связных компонент пространством размерности (n-1) (Рис. 61).

Сходным образом определим молекулярное нормальное пространство с

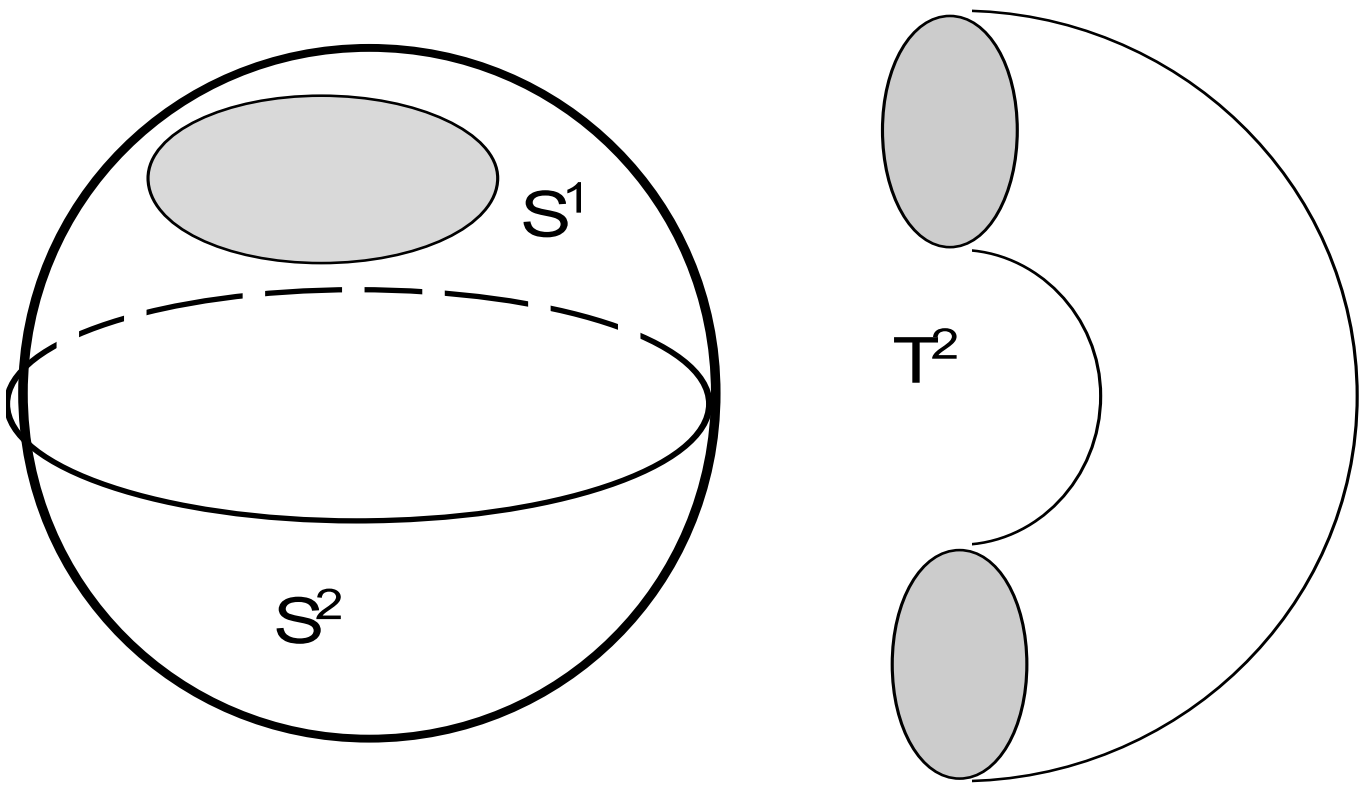

Рис. 61 Пространства с краем. Двумерная сфера $S^2$ с отверстием $S^1$, которое является одномерной окружностью. Цилиндр $T^2$ с краем, который является несвязным пространством, состоящим из двух окружностей.

краем. Если пространство G состоит из пространства H, и точки v, смежной со всеми точками подпространства P, P⊆H, мы будем использовать для краткости следующее обозначение: G=H∪(v⊕P). Например G=(G-v)∪(v⊕O(v)). Поскольку объединение пространств является прежде всего объединением точек, образующих пространство, то также можно использовать обозначение G=G∪(v⊕O(v)). Предварительно сделаем одно замечание. Во избежание недоразумений в дальнейшем мы будем рассматривать пространства с краем размерности более 0. Тем не менее этого ограничения можно избежать, если учесть, что мы уже определили, что единственное 0-мерное нормальное замкнутое пространство есть 0-мерная сфера, а единственное 0-мерное пространство с краем есть уединенная точка.

Определение пространства со связным краем

Пусть $G^n$ является нормальным n-мерным пространством, и точка v принадлежит $G^n$. Тогда пространство $H^n=G^n-v$ называется нормальным n-мерным пространством со связным краем B=O(v) (Рис. 62).

Любое нормальное n-мерное пространство $H^n$ с связным краем B можно представить в виде объединения подпространства A, $A=G^n-B$, и края B,

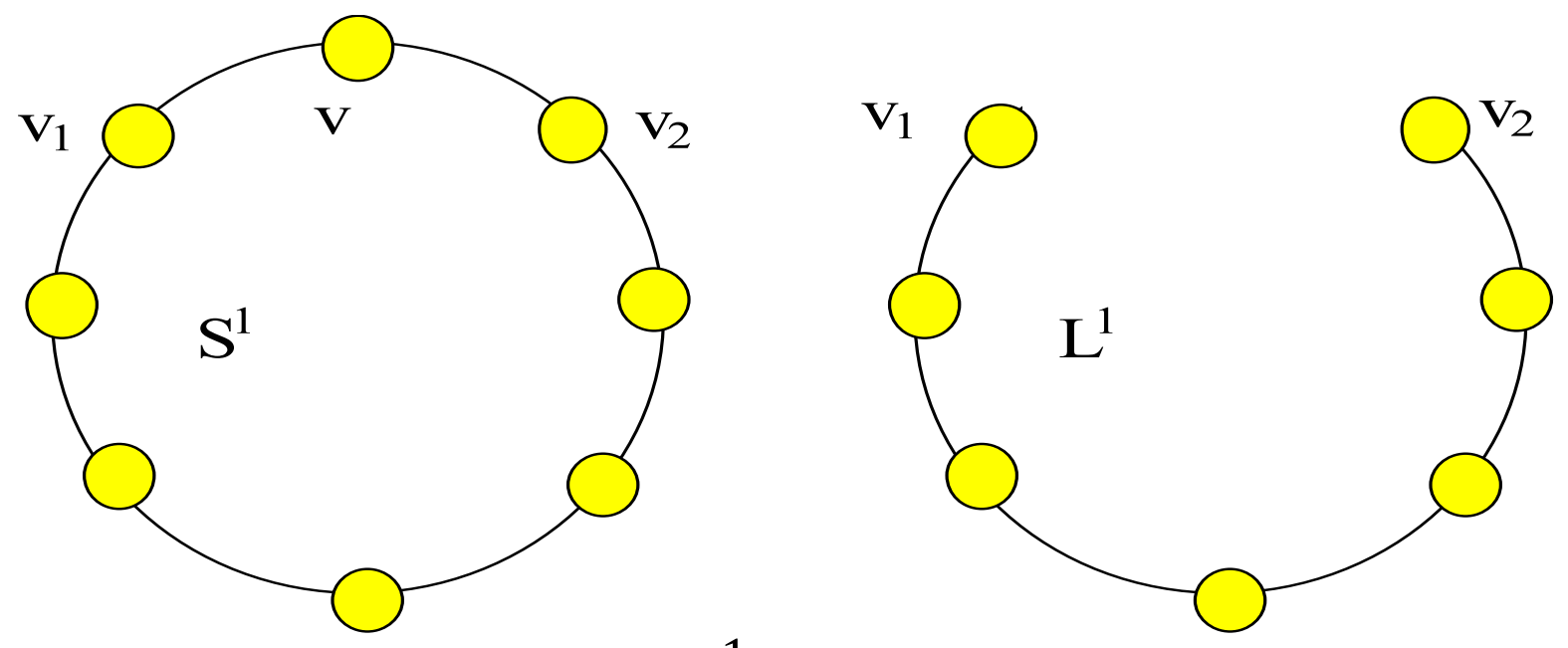

Рис. 62 В нормальной окружности $S^1$ точка v отбрасывается, и получившееся пространство $L_1$ становится пространством с краем $v_1$ и $v_2$. $L_1$ есть нормальный отрезок, край B есть нормальное 0-мерное пространство, состоящее из 2-х несвязных точек $v_1$ и $v_2$. Все остальные точки пространства $L_1$ являются одномерными.

$H^n=A\cup B=(A,B)$. Подпространство A, состоящее только из n-мерных нормальных точек, называется внутренностью или ядром пространства $H^n$. Таким образом $H^n$ есть объединение его ядра A и края B, $H^n=A\cup B=(A,B)$. Точки, принадлежащие ядру A, называются внутренними точками пространства $H^n$, точки, принадлежащие краю B, называются краевыми или граничными точками $H^n$.

Простейшими нормальными n-мерными пространствами с краем являются шары точек n-мерного нормального пространства. Действительно, шар нормальной n-мерной точки состоит из самой точки и ее окаема, который по определению есть нормальное замкнутое (n-1)-мерное пространство. На Рис. 63 показаны эти простейшие пространства с краем $D^1$, $D^2$, $D^3$ и $D^4$ размерности 1, 2, 3 и 4, являющиеся шарами 1, 2, 3 и 4-мерных нормальных точек. Нормальное пространство со связным краем возникает не только при отбрасывании точки, но и при отбрасывании связи. Докажем, что отбрасывание связи также приводит к нормальному n-мерному пространству со связным краем.

Теорема 40

*Пусть $H^n$ является нормальным n-мерным замкнутым пространством, и $(v_1v_2)$ есть связь между точками $v_1$ и $v_2$. Тогда пространство $G^n=H^n-(v_1v_2)$ является нормальным n-мерным пространством со связным краем.*

Доказательство.

Отбросим связь $(v_1v_2)$ из $H^n$. Окаемы точек $v_1$, $v_2$ и точек из $O(v_1v_2)$. изменились. Приклеим точку v к $H^n-(v_1v_2)$ по $v_1 \oplus v_2 \oplus O(v_1v_2)$. По теореме о гомеоморфной замене связи на точку полученное пространство $G^n \cup v=(H^n-(v_1v_2)) \cup v$ является нормальным n-мерным пространством. Следовательно, отбрасывание точки v превращает пространство $G^n \cup v$ в нормальное n-

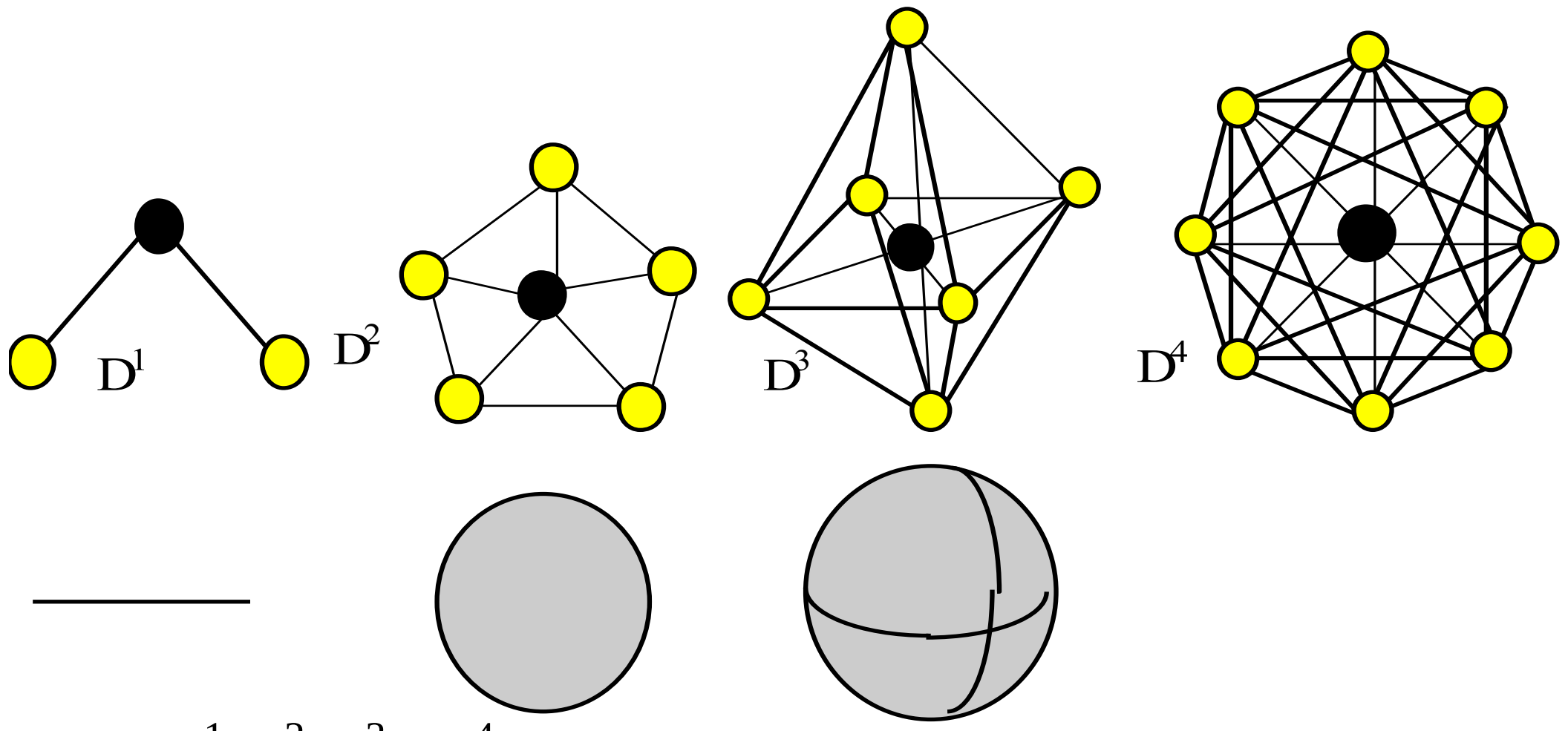

Рис. 63 $D^1$, $D^2$, $D^3$ и $D^4$ являются простейшими пространствами с краем. Они являются шарами 1, 2, 3 и 4-мерных нормальных точек и моделируют 1, 2, 3 и 4-мерные непрерывные диски.

мерное пространство со связным краем. Поскольку $(G^n \cup v)-v=((H^n-(v_1v_2)) \cup v)-v=H^n-(v_1v_2)$, то $H^n-(v_1v_2)$ также есть нормальное n-мерное пространство со связным краем. Теорема доказана.□

Свойства пространства с краем следуют из определения нормального n-мерного пространства.

Теорема 41

*Пусть $H^n$ является нормальным n-мерным пространством со связным краем B. Тогда край B будет нормальным замкнутым (n-1)-мерным пространством.*□

Доказательство.

Приклеим точку v к пространству $H^n$ по краю B. Полученное пространство

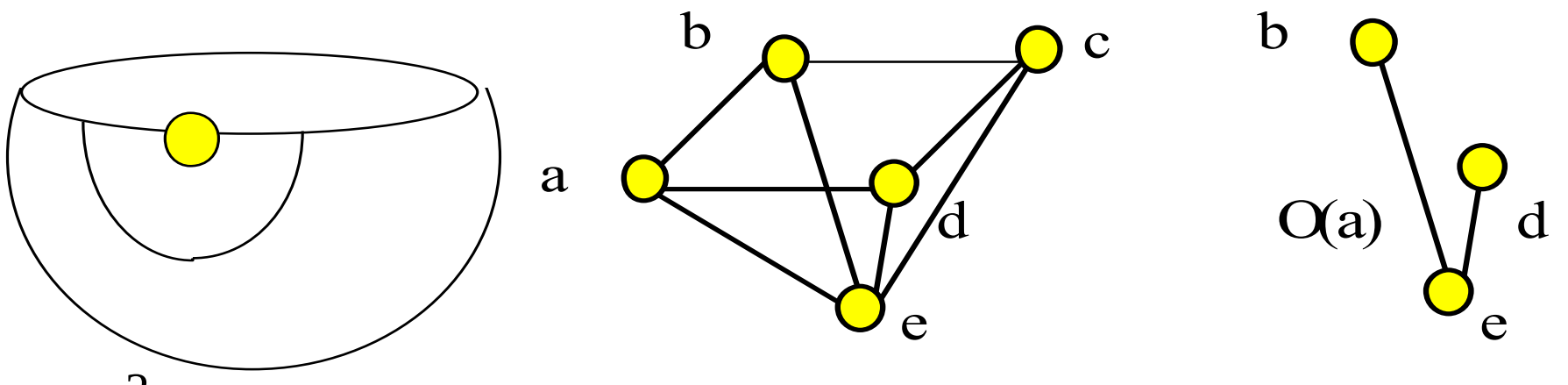

Рис. 64 $D^2$ является простейшим 2-мерным нормальным пространствам с краем (a, b, c, d). Окаем точки a, в свою очередь, есть 1-мерное нормальное пространство с краем (b, d).

$G^n=H^n\cup(v\oplus B)$. является нормальным n-мерным пространством. Следовательно, O(v)=B есть нормальное замкнутое (n-1)-мерное пространство. O(v)=B является связным, если (n-1)>0. Теорема доказана.□
На Рис. 63 диски являются пространствами с краями, края, в свою очередь, представляют нормальные сферы.

Теорема 42

*Пусть $H^n$ является нормальным n-мерным пространством со связным краем B, $H^n=A\cup B=(A,B)$. Тогда подпространство A не пусто, и все точки подпространства A являются n-мерными нормальными точками.* .

Доказательство.

Приклеим точку v к пространству $H^n$ по краю B. Полученное пространство $G^n=A\cup(v\oplus B)$. является нормальным n-мерным пространством. Следовательно, по следствию к предыдущей теореме существует по крайней мере одна точка, не смежная с v. Очевидно эта точка принадлежит A. Пусть точка $v_1$ принадлежит A. $O(v_1)$ есть нормальное замкнутое (n-1)-мерное пространство в $G^n$. Отбрасывание точки v из $G^n$ не меняет $O(v_1)$. Отсюда все точки $H^n$, не принадлежащие краю B, являются нормальными n-мерными точками. Теорема доказана.□

На Рис. 62 край B нормального отрезка $L_1=A\cup B$ состоит из 2-х несвязных точек и является нормальным 0-мерным пространством. Все точки подпространства A являются нормальными одномерными точками.

Теорема 43

*Пусть $H^n$ является нормальным n-мерным пространством со связным краем B, $H^n=A\cup B=(A,B)$. Тогда окаем O(v) любой точки v принадлежащей B есть нормальное (n-1)-мерное пространство со связным краем $B\cap O(v)$.*

Доказательство.

Приклеим точку v к пространству $H^n$ по краю B. Полученное пространство $G^n=A\cup(v\oplus B)$. является нормальным n-мерным пространством. Тогда в $G^n$

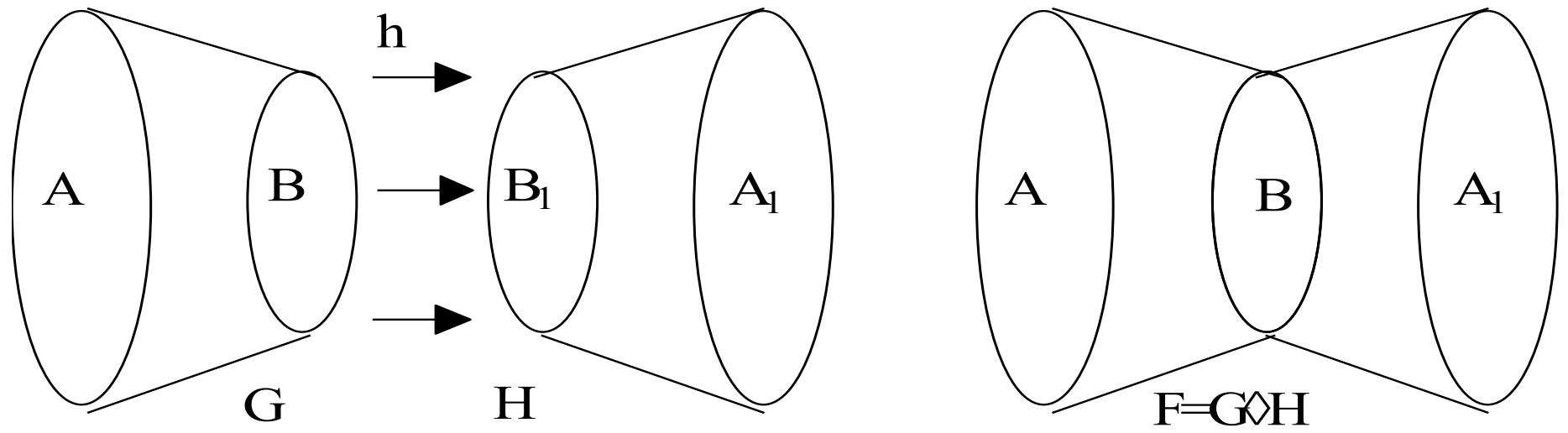

Рис. 66 Пространства G и H склеиваются в пространство F по подпространству B.

окаем $O(v_1)$ любой точки $v_1$, принадлежащей B, есть нормальное (n-1)-мерное пространство. Общий окаем $O(vv_1)$ есть нормальное (n-2)-мерное

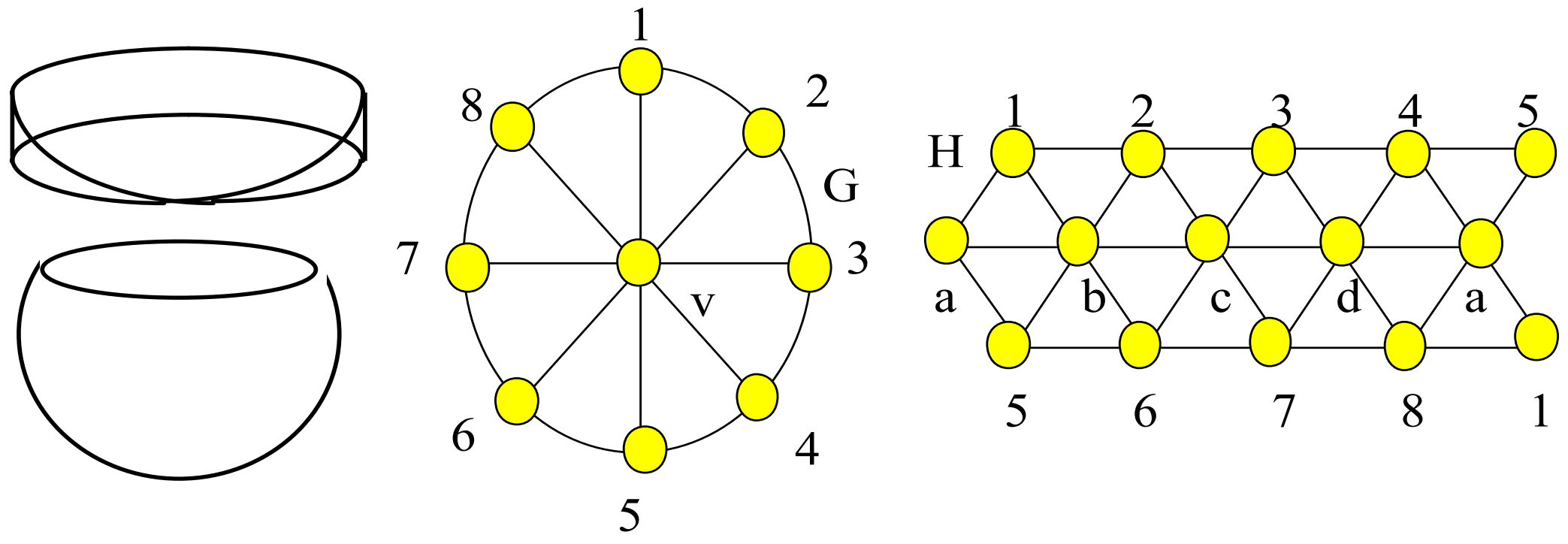

Рис. 65 Слева непрерывная сфера с дыркой заклеивается лентой Мебиуса. Справа молекулярная нормальная сфера с дыркой G (2-мерный диск) заклеивается молекулярной нормальной лентой Мебиуса H. G и H с краем B склеиваются в пространство F по подпространству B, являющемуся окружностью, состоящей из 8 точек 1,2,3,4,5,6,7,8. Склейка F=G◊H будет проективной плоскостью.

пространство. При отбрасывании точки v $O(v_1)$ переходит в $O(v_1)$-v, то есть становится нормальным (n-1)-мерным пространством со связным краем $B\cap O(v_1)$. Теорема доказана.□

На Рис. 64 $D^2$ является простейшим 2-мерным нормальным пространствам с краем (a, b, c, d), который есть окружность. Точка е будет двумерной точкой. Окаем точки a, в свою очередь, есть 1-мерное нормальное пространство с краем (b, d). Напомним определение склейки пространств.

Определение склейки двух пространств

Пусть $G=A\cup B$ и $H=A_1\cup B_1$ есть два молекулярных пространства и подпространство B изоморфно $B_1$, h: $B\to B_1$. Тогда склейкой G и H по изоморфизму h называется пространство $F=G\Diamond H=(G\cup H)/B=A\cup B\cup A_1$, в котором любая точка $v\in B$ отождествлена с $h(v)\in B_1$ (Рис. 66, Рис. 65).

При склейке нормальных пространств точки, принадлежащие краям, становятся нормальными n-мерными точками, точки, принадлежащие ядрам, не меняются. В дальнейшем если между B и $B_1$ уже установлен определенный изоморфизм, пространство $B_1$ мы будем обозначать как B. Очевидно, что в пространстве F шары $U(A)=G$, $U(A_1)=H$, и $U(A)\cap A_1=U(A_1)\cap A=\varnothing$.

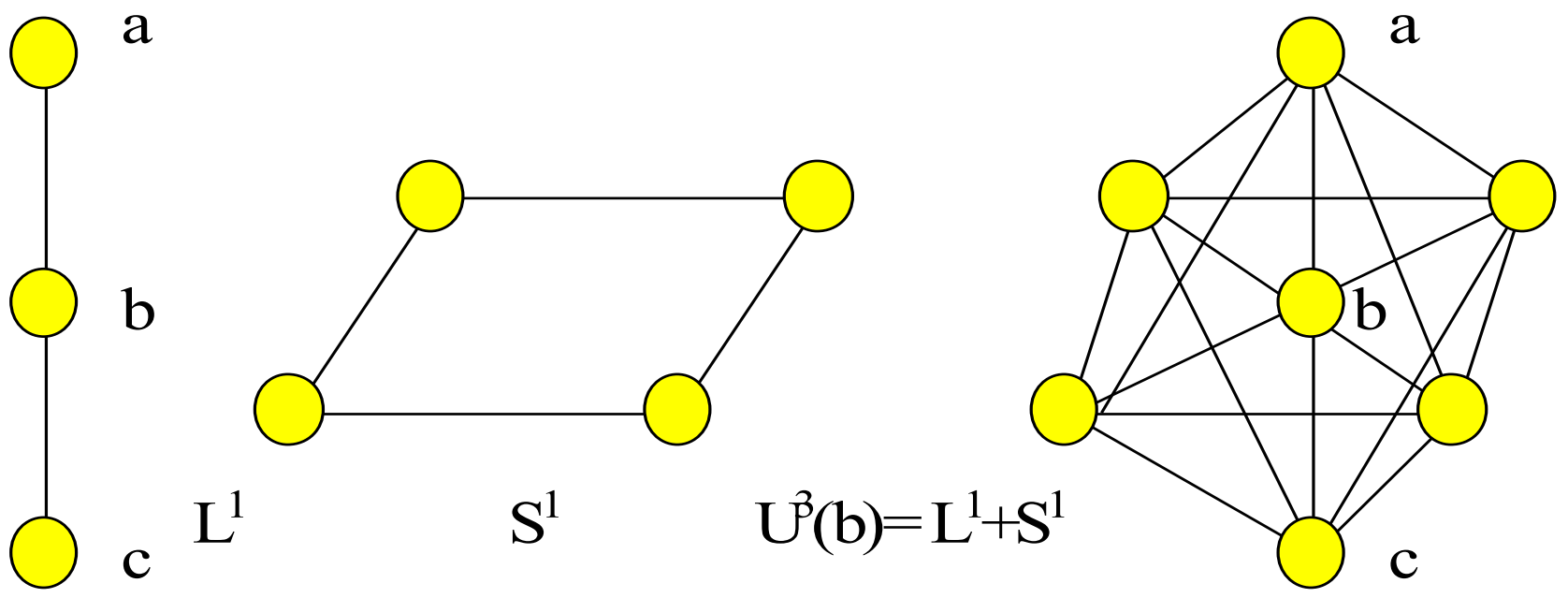

Рис. 67 Прямая сумма одномерного отрезка $L^1$, состоящего из 3-х точек, и окружности $S^1$ является трехмерным нормальным пространством $U^3(b)=L^1\oplus S^1$ с краем $S^2$ и шаром нормальной трехмерной точки.

Теорема 44

*Пусть $G^n=A\cup B$ и $H^n=A_1\cup B_1$ есть два нормальных пространства с краями и подпространство B изоморфно $B_1$, h: $B\to B_1$. Тогда склейка $F=G^n\Diamond H^n=A\cup B\cup A_1$ пространств $G^n$ и $H^n$ по изоморфизму h, является нормальным n-мерным пространством без края, (Рис. 66, Рис. 65).*

Доказательство.

Используем индукцию. Для n=1 теорема проверяется непосредственно. Предположим, что теорема верна для $G^n$ и $H^n$, где $n \le p$. Пусть $n=p+1$. Выберем произвольную точку $v_1 \in B$ и рассмотрим окаем этой точки в $F^n$. Очевидно, что $O(v_1)|F^n=O(v_1)|A \cup O(v_1)|B \cup O(v_1)|A_1$. Так как $O(v_1)|A_1 \cup O(v_1)|B$ и $O(v_1)|A \cup O(v_1)|B$ являются нормальными

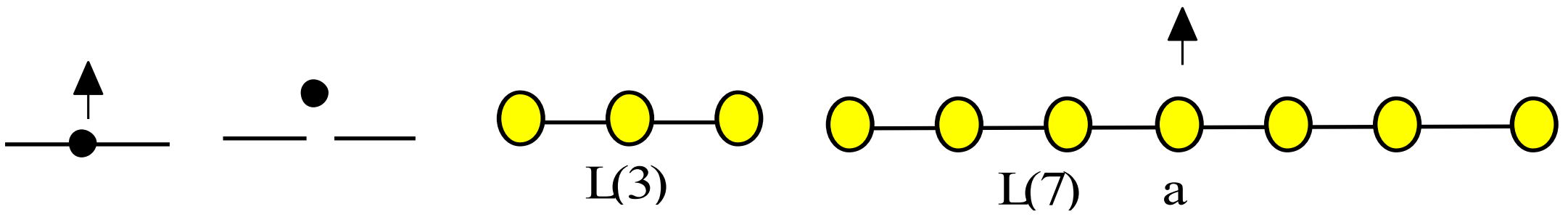

Рис. 68 Непрерывный отрезок может быть разбит на два одномерных отрезка удалением любой внутренней точки. Из молекулярного отрезка L(3) не может быть выброшена ни одна точка. Только отрезок L(7), состоящий из не менее чем 7 точек может быть разбит на два одномерных отрезка удалением точки а.

пространствами со связным краем размерности n-1, то $O(v_1)|F^n=O(v_1)|A_1 \cup O(v_1)|B \cup O(v_1)|A$ является их склейкой и, по предположению индукции, есть нормальное замкнутое (n-1)-мерное пространство. Следовательно все точки из B, а также из $A_1$ и A являются в $F^n$ нормальными n-мерными. Это означает, что $F^n$ есть нормальное n-мерное пространство. Теорема доказана. □

На Рис. 65 изображены два нормальных 2-мерных пространства G и H с краем B. Пространство G есть диск, пространство H есть лента Мебиуса. Край B является окружностью, состоящей из 8 точек. После склеивания G и H по краю B образуется проективная плоскость или склеенная сфера, уже рассматриваемая нами ранее. В непрерывном случае проективная плоскость получается точно так же, приклеиванием ленты Мебиуса к диску.

Теорема 45

*Пусть $G^n$, $n>0$, является нормальным n-мерным замкнутым пространством. Тогда конус $v \oplus G^n$ будет нормальным (n+1)-мерным пространством со связным краем $G^n$.*

Доказательство.

Рассмотрим прямую сумму 0-мерной нормальной сферы $S^0(v,u)$, состоящей из двух изолированных точек, и $G^n$. Тогда $T^{n+1}=S^0 \oplus G^n$ будет (n+1)-мерным нормальным замкнутым пространством. Отбросим из $T^{n+1}$ точку u. Получившееся пространство будет нормальным (n+1)-мерным пространством с краем $G^n$, а также конусом $v \oplus G^n$. Теорема доказана. □

Легко видеть, что в конусе $v\oplus G^n$ только одна точка v является (n+1)-мерной. Простейшими конусами такого рода являются шары нормальных точек, изображенные на Рис. 63.

Теорема 46

*Пусть $G^n$ является нормальным n-мерным пространством на конечном множестве точек с ядром A и со связным краем B, $G^n=A\cup B$. Пусть $H^m$ является m-мерным нормальным замкнутым пространством. Тогда их прямая сумма $G^n\oplus H^m$ будет нормальным (n+m+1)-мерным пространством с ядром A и со связным краем $B\oplus H^m$, являющимся (n+m)-мерным нормальным замкнутым пространством, (Рис. 67).*

Доказательство.

Приклеим точку v к $G^n$ по B, O(v)=B. Полученное пространство $F^n=G^n\cup v=A\cup(B\oplus v)$ будет нормальным замкнутым n-мерным

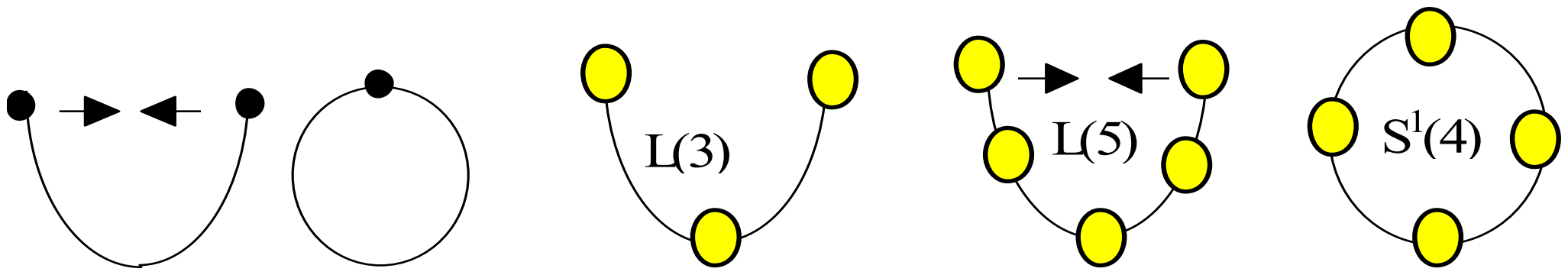

Рис. 69 Непрерывный отрезок всегда может быть склеен своими концами в окружность. Молекулярный отрезок L(3) не может быть склеен в окружность. Только отрезок L(5), состоящий из не менее чем 5 точек может быть склеен в минимальную окружность. $S^1(4)$.

пространством. Тогда прямая сумма $P^{n+m+1}=F^n\oplus H^m$ есть нормальное замкнутое (n+m+1)-мерное пространство. Отбросим из $P^{n+m+1}$ точку v. $P^{n+m+1}-v=G^n\oplus H^m$ будет нормальным (n+m+1)-мерным пространством со связным краем $B\oplus H^m$, являющимся (n+m)-мерным нормальным замкнутым пространством. Теорема доказана. □

В качестве примера рассмотрим прямую сумму одномерного отрезка $L^1$, состоящего из 3-х точек, и окружности $S^1$ (Рис. 67). Полученное пространство $U^3(b)= L^1\oplus S^1$ является трехмерным нормальным пространством с краем и шаром нормальной трехмерной точки. Несложно доказать также, что прямая сумма двух пространств с краем будет пространством без ядра, все точки которого являются краевыми.

З а м е ч а н и е.

Пусть $G^n$ является нормальным n-мерным пространством на конечном множестве точек с ядром A и со связным краем B, $G^n=A\cup B$. Пусть $H^m$ является также m-мерным пространством на конечном множестве точек с ядром $A_1$ и со связным краем $B_1$, $H^n=A_1\cup B_1$. Тогда их прямая сумма $G^n\oplus H^m$ будет связным пространством без ядра.

Рассмотрим теперь нормальное пространство с несвязным краем, который является объединением нескольких связных нормальных пространств (Рис. 61).

Операция склейки непрерывных пространств с краем по краю хорошо известна в классической топологии. Два пространства G и H с гомеоморфными краями могут быть склеены по гомеоморфизму краев. Эта операция позволяет получать новые пространства на основе уже известных. Например все замкнутые двумерные многообразия получаются при заклеивании сферы с дырками лентами Мебиуса или ручками.

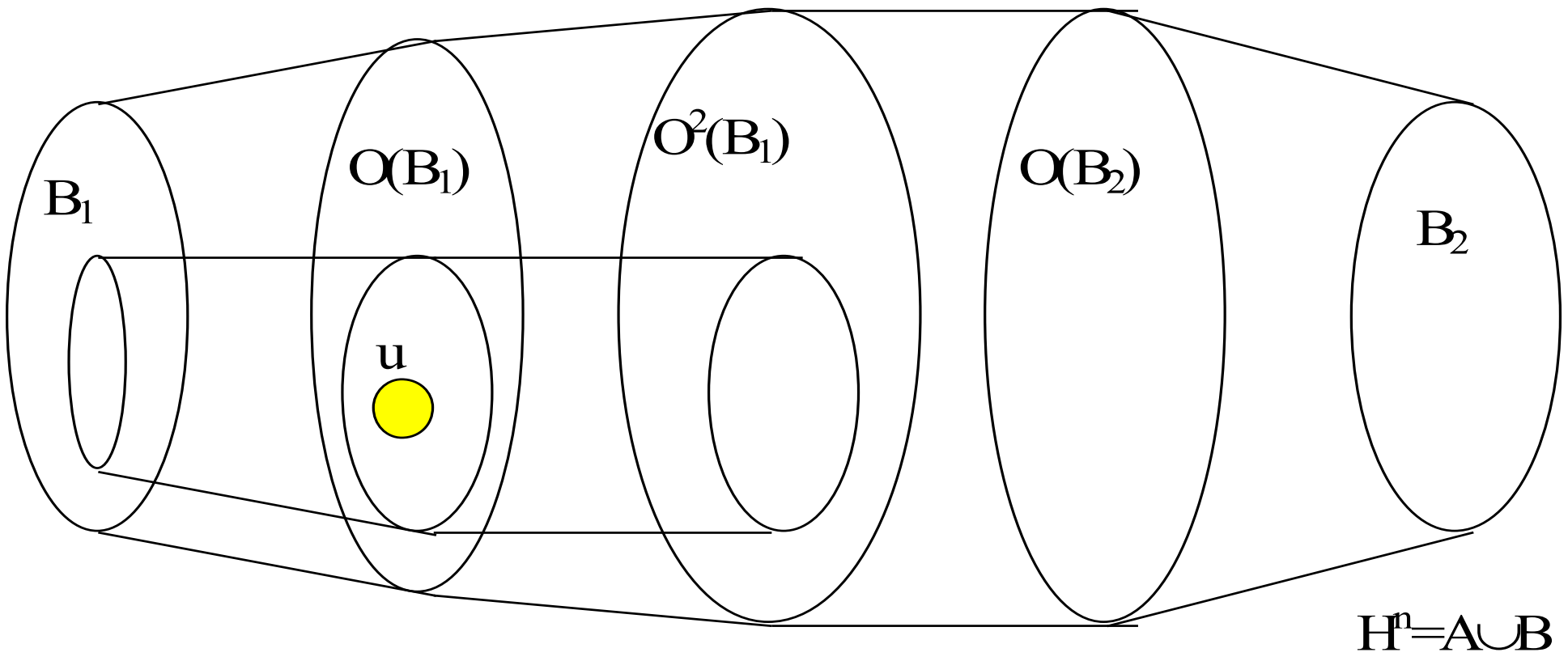

Рис. 70 При склейке в пространстве $H^n$ двух несвязных изоморфных краев $B_1$ и $B_2$ окаем внутренней точки u не меняется, если $U^2(B_1)\cap U(B_2)=\varnothing$.

При работе с многомерными компьютерными образами также возникают операции, аналогичные склеиванию непрерывных пространств, особенно в таких областях как компьютерная графика, графическая обработка результатов измерений, медицина. Поэтому важным является точно определить как саму операцию склеивания, так и границы, в которых мы можем применять эту операцию. При этом важно придерживаться классического топологического подхода, используя везде, где это возможно, хорошо развитый аппарат классической геометрии и топологии. При определении пространства со связным краем как нормального n-мерного пространства, из которого выброшена одна точка, никаких особых трудностей не возникает. Иное дело когда мы выбрасываем две или более точек. Первый вопрос, который возникает при этом, как близко могут быть расположены выбрасываемые точки. В непрерывном случае такой проблемы просто не существует. Из непрерывного единичного отрезка мы можем выбросить любое конечное число точек, при этом каждый из оставшихся кусков будет сохранять все свойства отрезка и содержать бесконечно много точек. В молекулярном пространстве любой конечный отрезок состоит всегда из конечного числа точек, и может так случиться, что при выбрасывании некоторого их числа отрезок перестанет

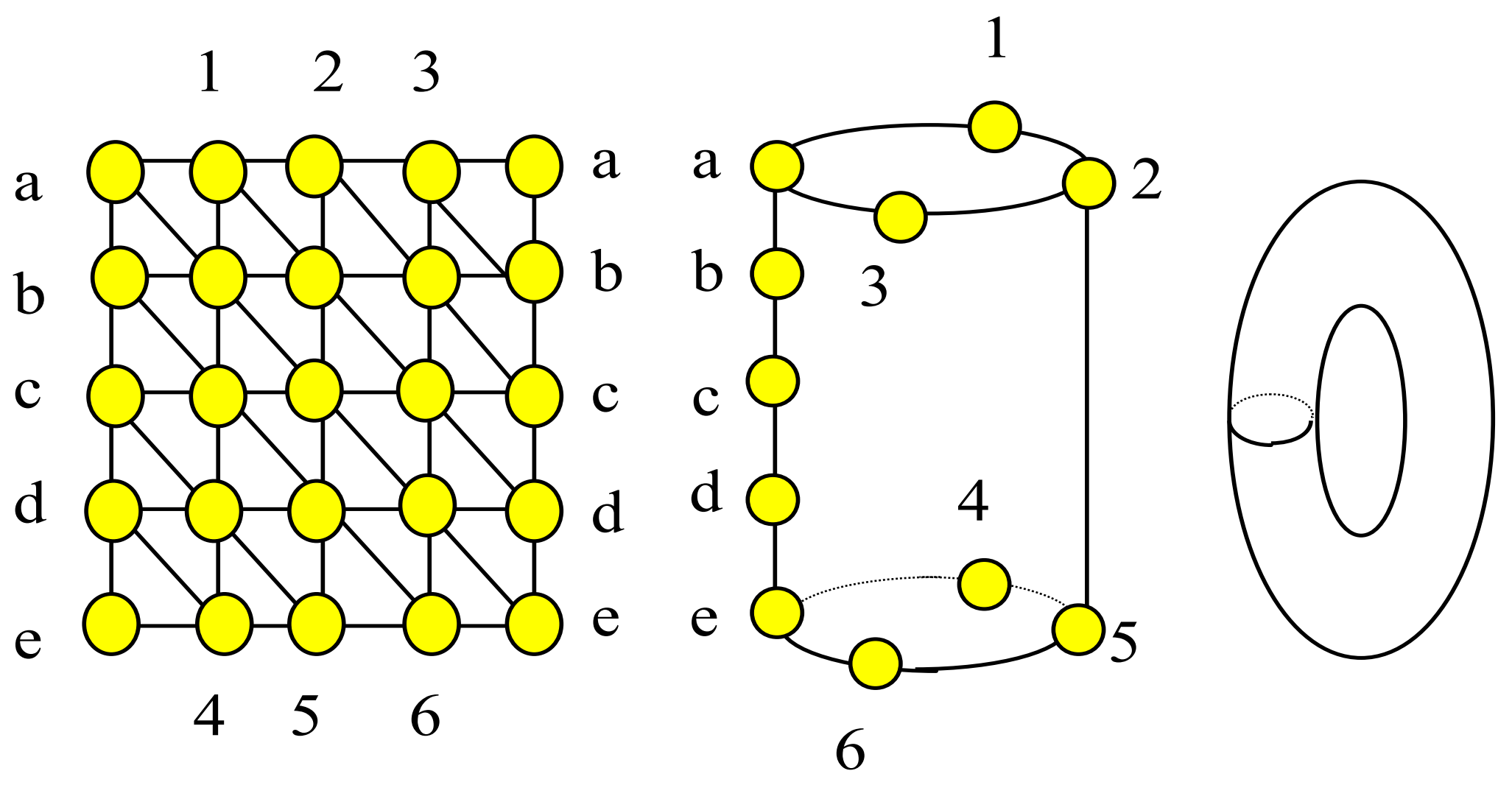

Рис. 71 Цилиндр $C^2$ является 2-мерным нормальным пространством с краем, состоящим из двух несвязных окружностей (1,2,3,a) и (4,5,6,a). При склеивании краев по изоморфизму (1,4), (2,5), (3,6), (a,e) получается тор $T^2$, изображенный справа. При склеивании краев по изоморфизму (1,6), (2,5), (3,4), (a,e) получается бутылка Клейна.

существовать как одномерное пространство. На Рис. 68 из отрезка L(3) мы не можем выбросить ни одной точки, а из отрезка L(7) можно выбросить точку a, при этом каждый из оставшихся отрезков будет по-прежнему одномерным нормальным пространством с краем.

То же самое относится к склеиваниям (Рис. 69). Непрерывный отрезок любой длины превращается в окружность при склеивании его концов. В случае молекулярного пространства только отрезок, состоящий из 5-и и более точек переходит в окружность при склеивании его концов. С другой стороны два отрезка L(3) как одномерные пространства с краем могут быть склеены по краю, образуя одномерную окружность $S^1(4)$.

Цель этих простых примеров показать, что операция склеивания в молекулярном пространстве является более сложной, чем при склеиваниях непрерывных пространств. Общий принцип должен быть следующий: при любых склеиваниях пространства с краем внутренние точки не должны быть затронуты, то есть топология внутренних точек меняться не должна.

Поскольку топология точки есть не что иное, как ее окаем, то окаем меняться не должен. Следует заметить, что если речь идет о пространстве со связным краем, то это условие выполняется автоматически.

Определение пространства с несвязным краем

Пусть $G^n$ является нормальным n-мерным пространством, и точки $v_1$, $v_2$, $v_3$,... $v_k$, принадлежащие $G^n$, удовлетворяют условию $U^2(v_s)\cap U(v_p)=\varnothing$, для любых $s\neq p$, s, p=1,2,...k. Тогда пространство $H^n=G^n-v_1-v_2-v_3-...v_k$, полученное из $G^n$ выбрасыванием этих точек, называется нормальным n-мерным пространством с несвязным краем $B=B_1\cup B_2\cup B_3\cup ...B_k$, где $B_s=O(v_s)$ s=1,2,...k, состоящим из k связных компонент.□

Край B состоит из k компонент связности, Каждая компонента связности $B_s$ окружена внутренними n-мерными нормальными точками, $O(B_s)\subseteq A$, где A есть ядро пространства $H^n$. Рассмотрим условия склеивания двух изоморфных компонент связности пространства $H^n$.

Теорема 47

*Пусть пространство $H^n=A\cup B$ является нормальным n-мерным пространством с несвязным краем $B=B_1\cup B_2\cup B_3\cup ...B_k$, и пусть подпространство $B_1$ изоморфно $B_2$, $h: B_1\to B_2$. Если $U^2(B_1)\cap U(B_2)=\varnothing$, то $B_1$ и $B_2$ могут быть склеены по изоморфизму h.* □

Доказательство.

Склеим $B_1$ и $B_2$ по изоморфизму h и выберем некоторую точку $v\in B_1$. По одной из предыдущих теорем окаем этой точки в получившемся пространстве становится нормальным (n-1)-мерным пространством. Нам также необходимо показать, что окаемы точек, принадлежащих ядру, не меняются. Пусть точка u принадлежит ядру A и $O(u)\cap B_1\neq\varnothing$ непусто. Следовательно $O(u)\in U^2(B_1)=B_1\cup O(B_1)\cup O^2(B_1)$ (Рис. 70), и при отождествлении $B_1$ и $B_2$ по изоморфизму h O(u) не меняется, то есть остается нормальным (n-1)-мерным пространством. Точно также не меняются окаемы всех остальных точек из ядра A. Теорема доказана. □

Образовавшееся пространство $F^n$ имеет ядро, к которому добавились точки из $B_1$ и край, в котором отсутствуют компоненты $B_1$ и $B_2$, после склейки перешедшие в ядро. Очевидно, что если край состоял только из этих двух компонент связности, то $F^n$ является пространством без края.

На Рис. 71 показано склеивание цилиндра $C^2$, который есть нормальное 2-мерное пространство с краем, состоящим из двух несвязных окружностей, в нормальное замкнутое пространство без края. Для наглядности слева цилиндр изображен в разрезе, где точки 1, 2, 3, 4 и 5 изображены дважды. Склеивание происходит по краям (1,2,3,a) и (4,5,6,e). При двух типах склеиваний получаются или тор, или бутылка Клейна.

## *НОРМАЛЬНЫЕ МОЛЕКУЛЯРНЫЕ МНОГООБРАЗИЯ И ИХ СВОЙСТВА*

Переходим к описанию свойств многообразий. Как уже отмечалось ранее большинство многомерных объектов, обрабатываемых компьютерами, являются многообразиями. Это относится ко всем образам, возникающим в компьютерных играх, компьютерной графике, медицинских анализах и

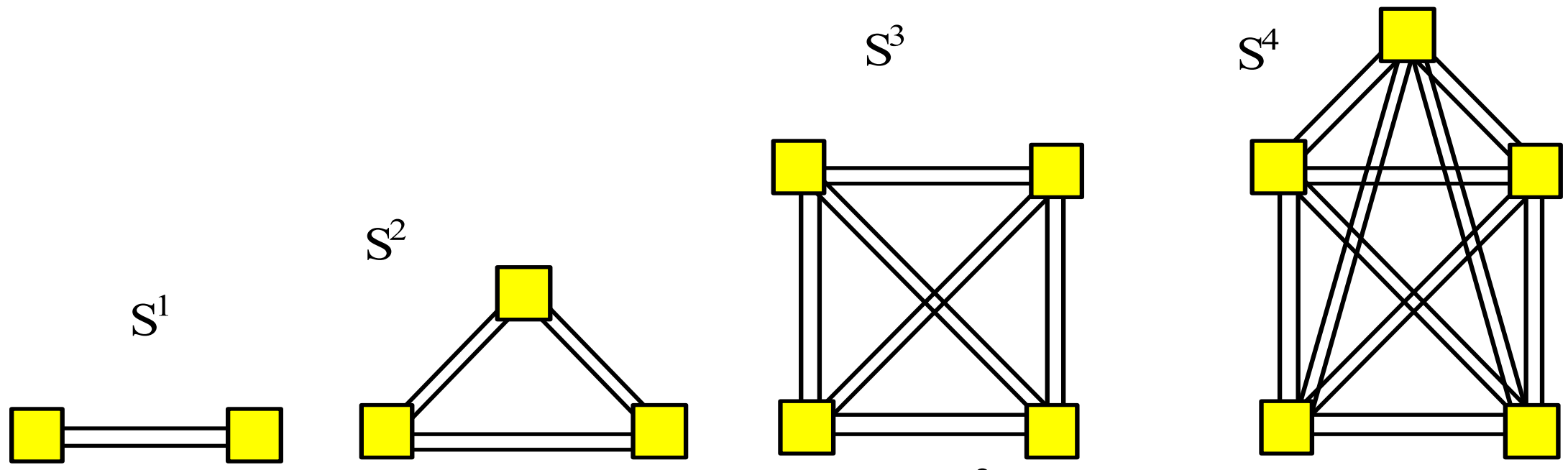

Рис. 72 Квадраты изображают 0-мерную сферу $S^0$. Каждая точка одной сферы соединена со всеми точками всех остальных сфер. Прямые суммы двух, трех, четырех и пяти 0-мерных сфер образуют минимальные сферы размерности 1, 2, 3 и 4.

исследованиях и во многих других областях. Поэтому изучение свойств многообразий составляет важное направление в дигитальной топологии. Для сохранения полноты изложения мы повторим некоторые из определений и рисунков, которые были введены ранее.

Рассмотрим, вначале, несколько более подробно n-мерные сферы $S^n$, важный класс нормальных молекулярных пространств. Важность этого класса определяется прежде всего тем, что в многообразиях окаем каждой точки есть сфера.

Определение минимальной нормальной сферы

Минимальной нормальной n-мерной сферой называется прямая сумма (n+1) копий 0-мерной сферы $S^0$, состоящей из двух уединенных точек.

$$S^n = S^0 \oplus S^0 \oplus S^0 \ldots \oplus S^0 = \sum_{k=0}^{k=n} S_k^0, \quad S_k^0 = S^0, \quad k = 0,1,2,\ldots,n.$$

На Рис. 72 условно изображены некоторые минимальные сферы. Многие свойства нормальной n-мерной сферы почти очевидны и не нуждаются в доказательствах. Другие свойства будут доказаны позднее. Перечислим основные свойства $S^n$.

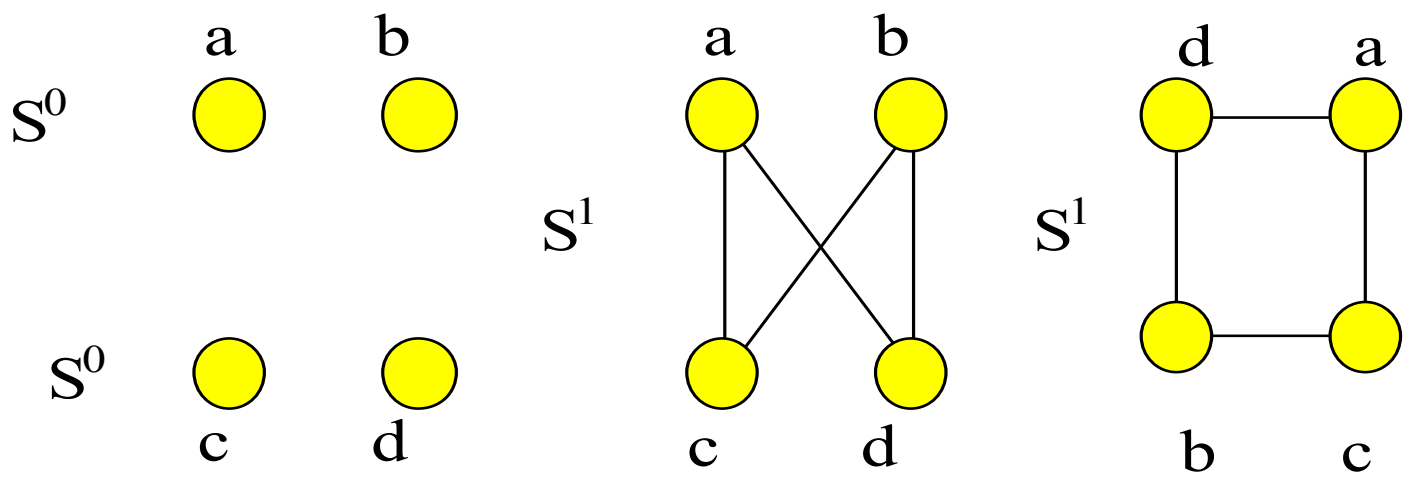

Рис. 73 Сумма двух 0-мерных сфер является одномерной окружностью $S^1=S^0 \oplus S^0$

### Свойства минимальной нормальной n-мерной сферы

1. $S^n$ есть нормальное замкнутое n-мерное пространство.
2. Минимальная нормальная n-мерная сфера есть однородное пространство, то есть окаемы всех точек изоморфны.
3. Число точек минимальной нормальной n-мерной сферы равно 2n+2, $|S^n| = 2n+2$, (Рис. 73, *Рис. 76*, Рис. 74, Рис. 75).
4. Минимальная сфера $S^n$ есть полное (n+1)-дольное пространство, $S^n = К(2,2,...2)$.

Свойство 1 является прямым следствием теоремы о прямой сумме замкнутых молекулярных пространств.

Свойство 2 также очевидно. Окаем любой точки n-мерной минимальной

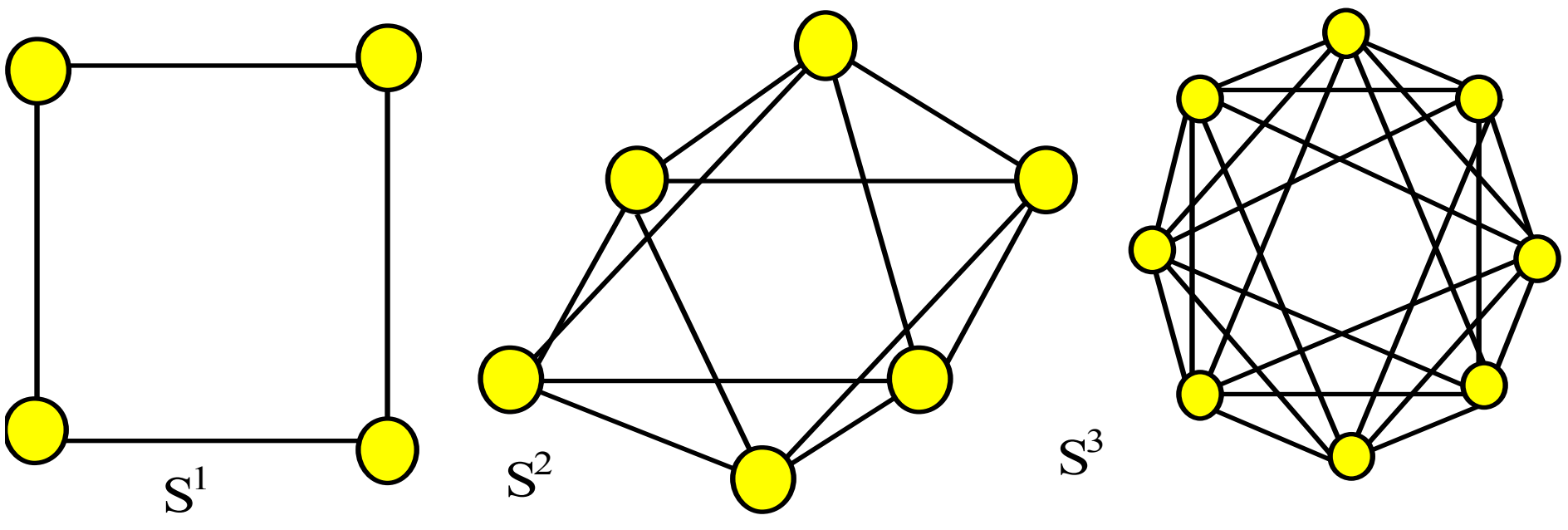

Рис. 74 Сферы $S^1$, $S^2$ и $S^3$ являются единственными замкнутыми n-мерными пространствами, для которых в размерностях 1, 2 и 3 выполняется соотношение $|G^n|=2n+1$.

сферы есть (n-1)-минимальная сфера. Следовательно, $S^n$ однородна. Свойство 3 есть следствие теоремы о наименьшем количестве точек замкнутого нормального n-мерного пространства, которая гласит, что 2n+2 есть минимальное число точек такого пространства. Поскольку $|S^n|=2n+2$, то $S^n$ минимальна по числу точек (Рис. 73, *Рис. 76*, Рис. 74, Рис. 75).
Свойство 4 есть просто определение полного (n+1) дольного пространства К(2, 2,..2).
Обратимся к примерам.

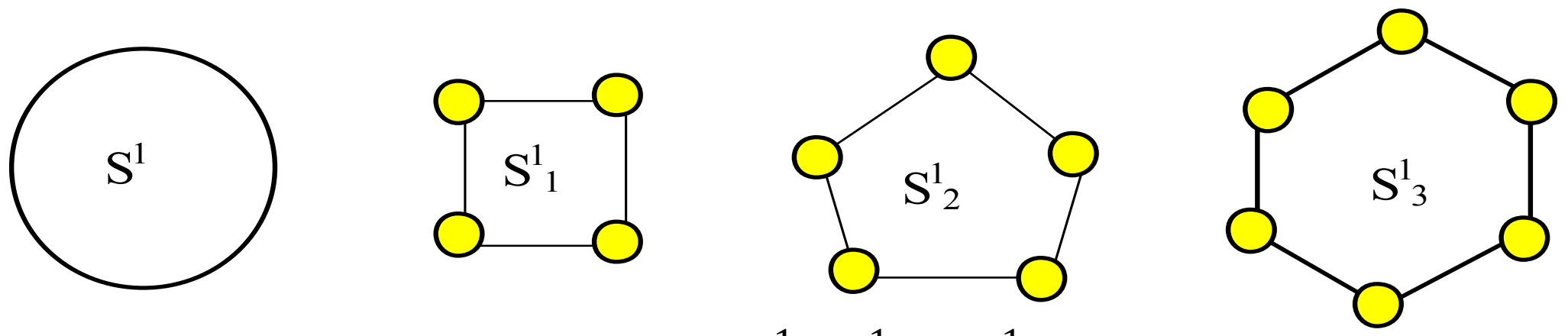

Рис. 75 Изображены три окружности $S^1{}_1$, $S^1{}_2$, и $S^1{}_3$, состоящие из 4, 5 и 6 точек. Сфера $S^1{}_1$, состоящая из 4-х точек, минимальна.

Пусть имеется две 0-мерных сферы $S^0{}_1$ и $S^0{}_2$ с точками (a,b) и (c,d) соответственно, Рис. 73. Соединяем связями каждую точку $S^0{}_1$ со всеми точками $S^0{}_2$. Из рисунка легко видеть, что мы получаем минимальную одномерную сферу $S^1$-окружность $S^1 = S^0 \oplus S^0$.
Теперь возьмем сумму 0-мерной и минимальной одномерной сфер $S^0$ и $S^1$, как это показано на *Рис. 76*. Соединяя все точки одной сферы со всеми точками другой сферы, очевидно, мы получаем минимальную двумерную сферу $S^2=S^1 \oplus S^0=S^0 \oplus S^0 \oplus S^0$.
Легко проверить непосредственно, что для трехмерной сферы (Рис. 74). $S^3=S^0 \oplus S^0 \oplus S^0 \oplus S^0=S^1 \oplus S^1=S^0 \oplus S^2$.
На примерах в предыдущих разделах мы видели, что существует

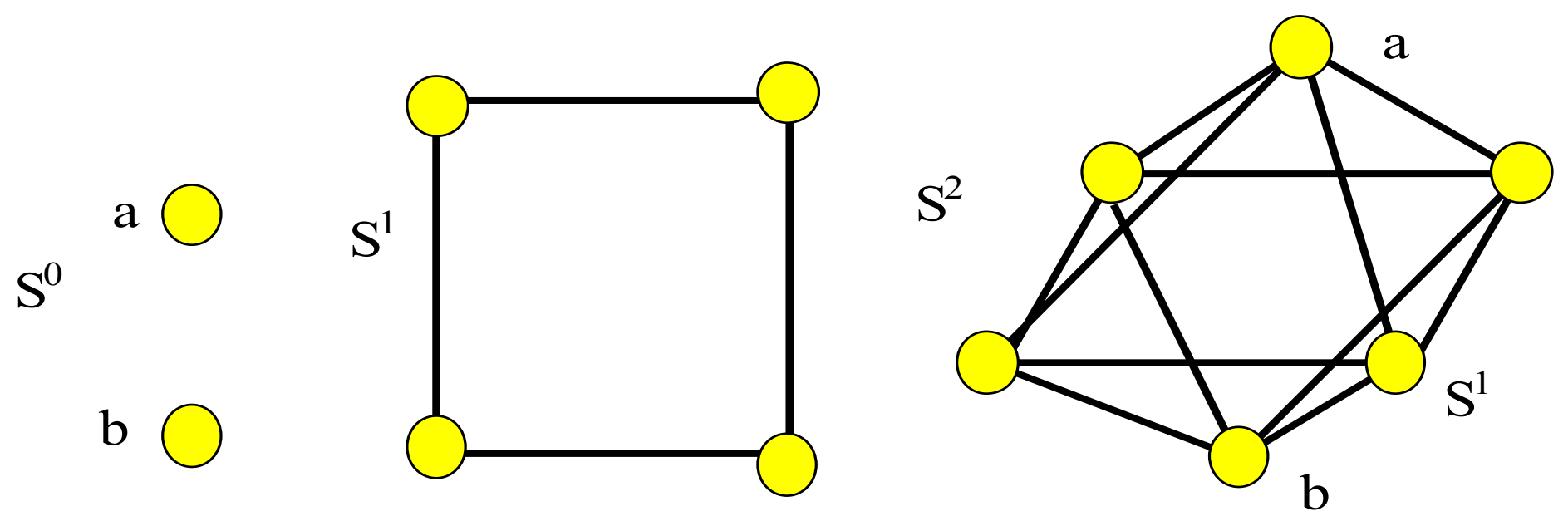

Рис. 76 Сумма 0-сферы $S^0$ и 1-сферы $S^1$ является 2-мерной сферой. $S^2 = S^0 \oplus S^1= S^0 \oplus S^0 \oplus S^0$.

несколько молекулярных окружностей и 2-мерных сфер, отличающихся числом точек и связями между ними. На Рис. 75 и Рис. 77 мы воспроизводим окружности и 2-мерные сферы.
Связь между двумя топологически гомеоморфными молекулярными пространствами обеспечивается гомеоморфными преобразованиями. На основе минимальной n-мерной сферы и гомеоморфных преобразований дадим теперь общее определение нормальной n-мерной сферы.

Определение нормальной n-мерной сферы

Нормальной n-мерной сферой $S^n$ является МП, полученное из минимальной n-мерной сферы при помощи гомеоморфных преобразований замена связи на точку и замена точки на связь, $S^m=\Phi_s\Phi_{s-1}\Phi_{s-2}...\Phi_1 S^m_{min}$.

Из определения следует, что нормальной n-мерной сферой $S^n$ является МП, сводимое к минимальной n-мерной сфере при помощи гомеоморфных преобразований замена связи на точку и замена точки на связь. На Рис. 78 точка минимальная окружность из 4-х точек заменяется на окружность из 5-и точек заменой связи на точку; 6-и точечная окружность переводится в пятиточечную заменой точки на связь.

Прямая сумма двух нормальных молекулярных замкнутых пространств также будет нормальным замкнутым пространством. Аналогичное утверждение может быть доказано для сфер.

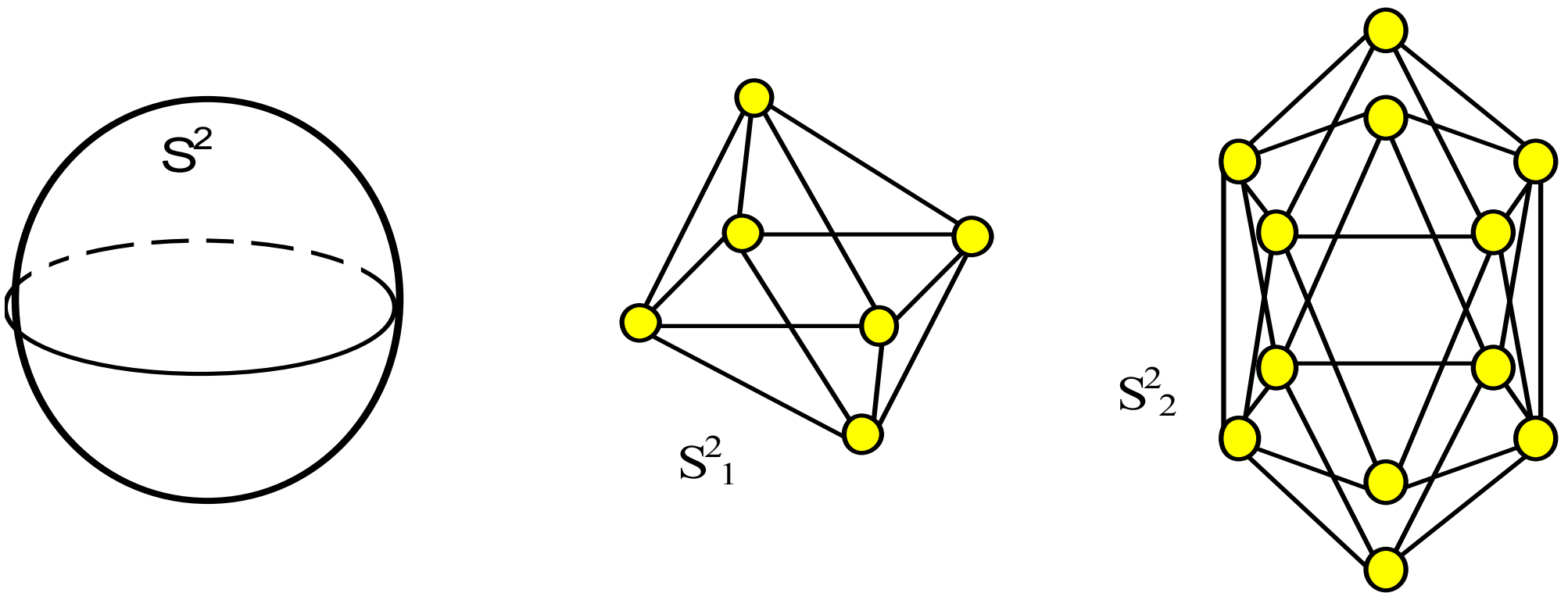

Рис. 77 Двумерные сферы. Окаем каждой точки двумерной сферы есть окружность. Двумерная сфера $S^2_1$ минимальна, она состоит из 6 точек- наименьшего числа точек, необходимых для образования двумерной сферы. Сфера $S^2_2$ состоит из 12 точек. Обе сферы однородны.

Теорема 48

*Пусть $S^m$ и $S^n$-нормальные сферы размерности m и n. Тогда их сумма $W^{m+n+1} = S^m \oplus S^n$ является нормальной сферой размерности (m+n+1). $S^{m+n+1} = S^m \oplus S^n$.*

Доказательство.

По сути дела эта теорема является следствием теоремы о гомеоморфизме прямой суммы пространств. Сфера $S^m$ по определению получена из минимальной m-мерной сферы $S^m_{min}$ при помощи ряда гомеоморфных преобразований $\Phi_1, \Phi_2, \Phi_3,\ldots \Phi_s$. $S^m=\Phi_s\Phi_{s-1}\Phi_{s-2}\ldots \Phi_1 S^m_{min}$. Согласно теореме о гомеоморфизме прямой суммы пространств при каждом преобразовании $\Phi_k$ или $\Phi_k^{-1}$ выполняется соотношение

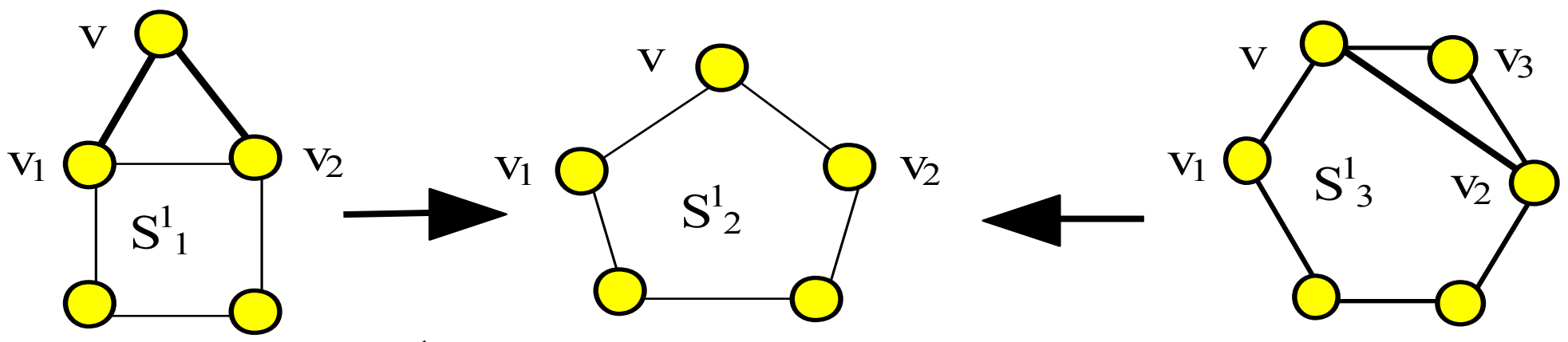

Рис. 78 Окружность $S^1_1$, состоящая из 4-х точек переходит в окружность $S^1_2$, состоящую из 5-и точек, заменой связи ($v_1v_2$) на точку v. Окружность $S^1_3$, состоящая из 6-х точек также переходит в окружность $S^1_2$, состоящую из 5-и точек, заменой точки $v_3$ на связь ($vv_3$).

$(\Phi_k S^m)\oplus S^n \approx \Phi_k(S^m \oplus S^n) \approx S^m \oplus S^n$. Применяя обратные преобразования в обратной последовательности к $S^m$ и, следовательно, к $S^m\oplus S^n$ мы получаем: $S^m\oplus S^n \approx \Phi_s^{-1} S^m\oplus S^n \approx \Phi_{s-1}^{-1}\Phi_s^{-1} S^m\oplus S^n \approx\ldots\approx \Phi_1^{-1}\Phi_2^{-1}\ldots \Phi_{s-1}^{-1}\Phi_s^{-1} S^m\oplus S^n \approx S^m_{min}\oplus S^n$. То есть $S^m\oplus S^n \approx S^m_{min}\oplus S^n$. Такое же рассуждение применимо к $S^n$, $S^m_{min}\oplus S^n \approx S^m_{min}\oplus S^n_{min}$. Это означает, что $S^m\oplus S^n \approx S^m_{min}\oplus S^n_{min} = S^{m+n+1}_{min}$. Таким образом прямая сумма двух (и более) сфер $S^m\oplus S^n$ является (m+n+1)-мерной сферой $S^{m+n+1}$, $S^m\oplus S^n = S^{m+n+1}$. Теорема доказана. □

Перейдем от сфер к многообразиям общего вида.

Определение нормального многообразия

Молекулярное нормальное n-мерное пространство $M^n$ называется многообразием, если окаем каждой точки является (n-1)-мерной нормальной молекулярной сферой $S^{n-1}$.

Иными словами все точки замкнутого многообразия являются сферическими. Примерами многообразий являются все n-мерные сферы $S^n$ и n-мерные плоские пространства, двумерные тор $T^2$, проективная плоскость $P^2$ и бутылка Клейна $K^2$. Какова локальная структура многообразия?

Теорема 49

*Пусть $G^n$ является молекулярным (не обязательно замкнутым) нормальным n-мерным многообразием, где K(m) есть связка (подпространство из m попарно связных точек $v_1$, $v_2$, $v_3$,... $v_m$. Тогда общий окаем этих точек есть нормальная (n-m)-мерная сфера $O(v_1,v_2,v_3,...v_m)=O(v_1)\cap O(v_2)\cap O(v_3)\cap ...\cap O(v_m)=S^{n-m}$.*

Доказательство.

По определению нормального молекулярного n-мерного многообразия окаем $O(v_1) = S^{n-1}$ произвольной точки пространства $v_1$ является молекулярной замкнутой нормальной (n-1) мерной сферой. Пусть $v_1$ и $v_2$ смежны, то есть между этими точками имеется связь. Это означает, что $v_2 \in S^{n-1} = O(v_1)$. Следовательно, $O(v_2)|S^{n-1} = O(v_2v_1)$ является нормальной (n-2)-мерной сферой $S^{n-2}$. Пусть $v_3$ смежна с $v_2$ и $v_1$, то есть $v_3 \in O(v_2v_1) = S^{n-2}$. Следовательно, $O(v_3)|S^{n-2} = O(v_3v_2v_1)$ является нормальной (n-3)-мерной

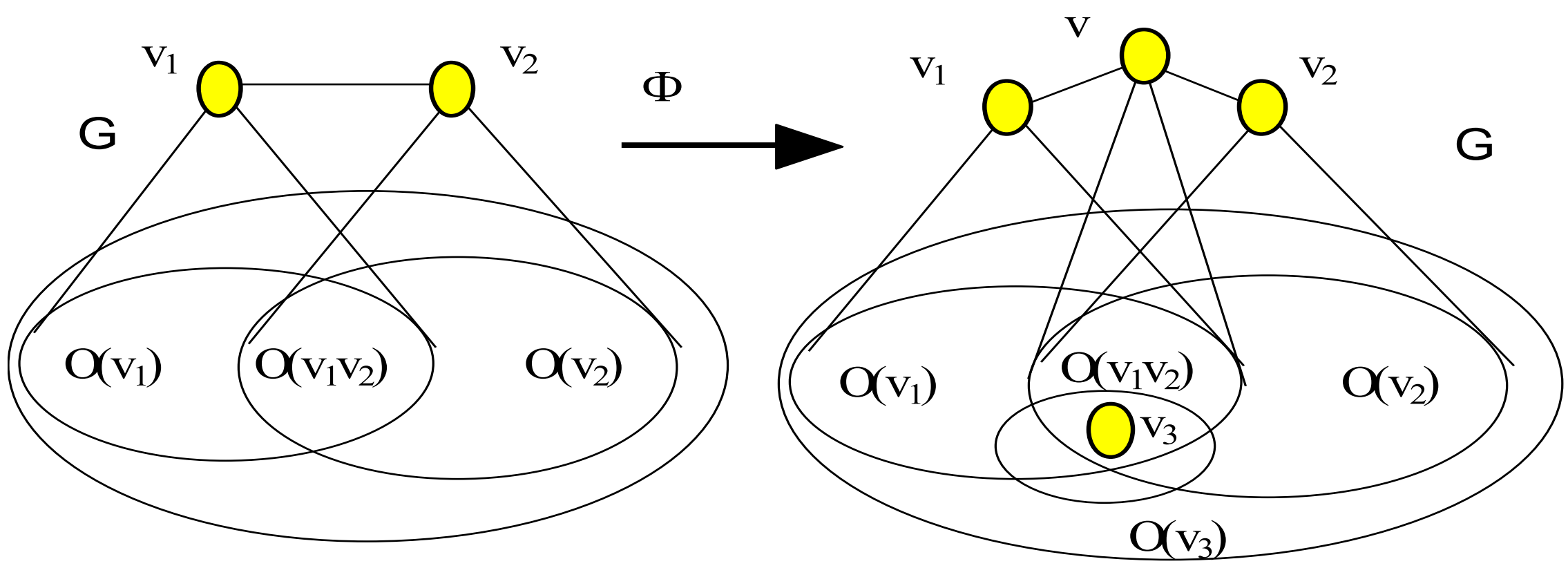

Рис. 79 Гомеоморфные преобразования многообразия. $O(v_1v_2)=S^{n-2}$. Связь $(v_1v_2)$ меняется на точку v, для которой $O(v) = S^0(v_1,v_2)\oplus O(v_1v_2)$.

сферой. Повторяя этот процесс для точек $v_1$, $v_2$, $v_3$,... $v_m$ мы получаем, что $O(v_1,v_2,v_3,....v_m)=O(v_1)\cap O(v_2)\cap O(v_3)\cap ...\cap O(v_m)=S^{n-m}$. Теорема доказана.□

Эта теорема говорит о том, что общий окаем любого количества попарно смежных точек является всегда сферой, в отличие от нормальных пространств общего вида, где окаем может быть замкнутым n-мерным пространством любой топологии. Важным является тот факт, что гомеоморфные преобразования сохраняют многообразия, то есть при гомеморфных преобразованиях окаемы точек остаются сферами. Как мы уже отмечали, в произвольном нормальном молекулярном пространстве окаемы могут быть любыми замкнутыми нормальными пространствами, сферами, торами и так далее. При гомеоморфных преобразованиях окаем данной точки может изменить свою топологию. И только в одном случае топология окаема сохраняется, это когда пространство является многообразием.

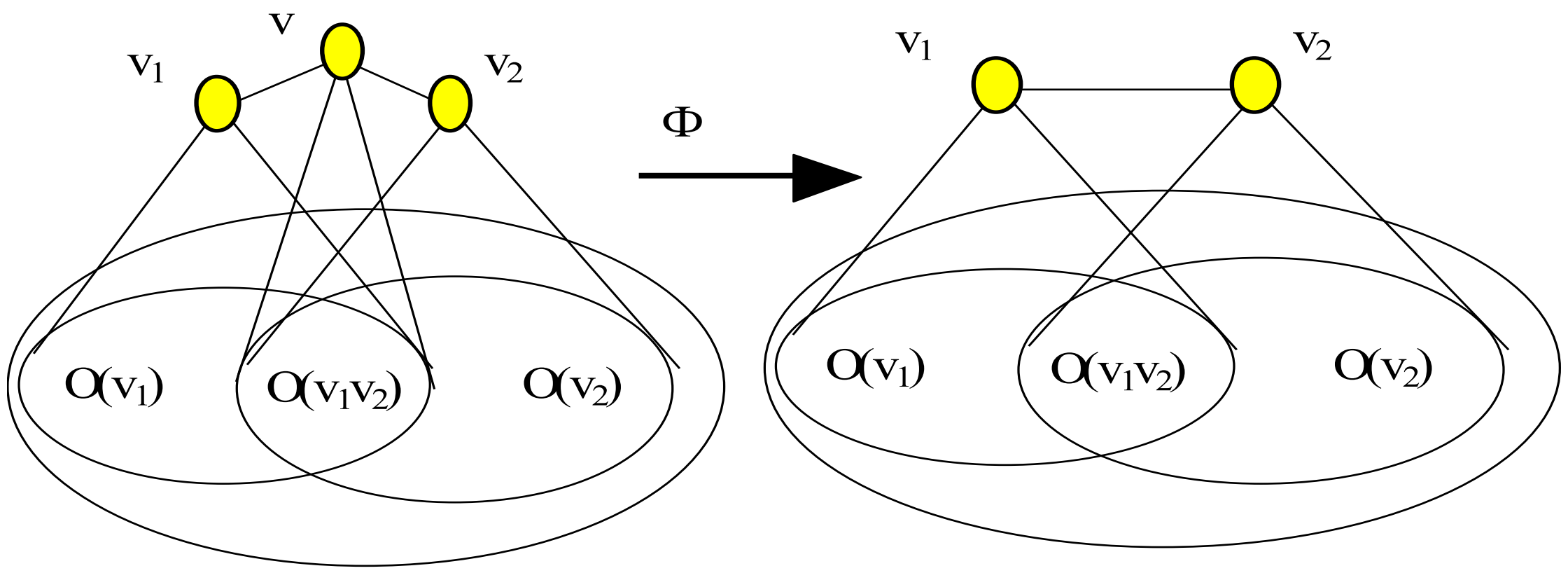

Рис. 80 Гомеоморфные преобразования молекулярного многообразия. Точка v меняется на связь ($v_1v_2$). При этом должны соблюдаться условия: $O(v) = S^0(v_1,v_2)\oplus O(v_1v_2v)$, $O(v_1v_2) = v\oplus O(v_1v_2v)$, где $S^0$ состоит из точек $v_1$ и $v_2$.

Теорема 50

*Гомеоморфное преобразование Φ замена связи на точку переводит многообразие $G^n$ в многообразие $H^n=\Phi G^n$.*

Доказательство.

Доказательство теоремы сходно с доказательством аналогичной теоремы для нормальных пространств. Для малых размерностей теорема проверяется непосредственно. Предположим, что теорема верна для размерности $n < k$. Пусть $n = k$.

Пусть $G^n$ есть некоторое многообразие, $v_1$ и $v_2$ есть две точки этого пространства со связью ($v_1v_2$) (Рис. 79). Согласно определению O($v_1v_2$) есть нормальная (n-2)-мерная сфера. Удалим связь ($v_1v_2$) и соединим точку v с точками со всеми точками из O($v_1v_2$), а также с $v_1$ и $v_2$. Легко видеть, что окаем O($v_1$) точки $v_1$, не изменился кроме того, что точка $v_2$ в этом окаеме заменилась на точку v. Сходным образом не изменился окаем O($v_2$), точки $v_2$. Следовательно, окаемы этих точек остались (n-1)-мерными нормальными сферами после замены связи на точку.

Окаем O(v) точки v имеет структуру суммы 0--мерной сферы, состоящей из двух несвязных точек $v_1$ и $v_2$ и (n-2)-мерной сферы O($v_1v_2$), $O(v)=S^0\oplus O(v_1v_2)$. Следовательно, O(v) есть нормальная (n-1)-сфера. Пусть точка $v_3$ лежит в O($v_1v_2$). O($v_3$) есть нормальная сфера размерностью (n-1).

В нем происходит точно такое же преобразование замена связи на точку. Согласно предположению индукции после этого преобразования $O(v_3)$ остается нормальной сферой той же размерности. Все остальные точки пространства $G^n$ своих окаемов не меняют. Таким образом окаем каждой точки остается нормальной (n-1)-мерной сферой. Теорема доказана. □

Теорема 51

*Гомеоморфное преобразование Ф замена точки на связь переводит многообразие $G^n$ в многообразие $H^n=ФG^n$.*

Доказательство.

Доказательство теоремы сходно с доказательством аналогичной теоремы для нормальных пространств. Используем индукцию. Для размерности n=1 теорема проверяется непосредственно. Предположим, что теорема верна для размерности $n < k$. Пусть $n = k$. Пусть $G^n$ есть некоторое молекулярное пространство, v есть какая-то точка этого пространства, для которой $O(v) = S^0(v_1,v_2)\oplus O(v_1v_2v)$ и $O(v_1v_2) = v\oplus O(v_1v_2v)$, где $S^0$ есть две несвязные точки $v_1$ и $v_2$ из U(v) (Рис. 80). Удалим точку v и установим связь $(v_1v_2)$. Легко видеть, что окаемы всех точек, не лежащих в $O(v_1v_2)$, не изменились. Следовательно, окаемы этих точек остались (n-1)-мерными нормальными сферами после замены точки на связь. Рассмотрим окаем $O(v_3)$ точки $v_3$ лежащей в $O(v_1v_2)$. Этот окаем является нормальной замкнутой (n-1)-мерной сферой и имеет структуру исходного пространства $G^n$, то есть $O(v)|O(v_3) = S^0(v_1,v_2)\oplus O(v_1v_2v)|O(v_3)$ и $O(v_1v_2)|O(v_3) = v\oplus O(v_1v_2v)|O(v_3)$, где $S^0$ есть две несвязные точки $v_1$ и $v_2$ из $O(v)|O(v_3)$. В $O(v_3)$ происходит точно такое же преобразование замена связи на точку. Согласно предположению индукции после этого преобразования подпространство $O(v_3)$ остается нормальной сферой той же размерности (n-1). Таким образом окаем каждой точки остается замкнутой нормальным (n-1)-мерной сферой. Теорема доказана. □

Эти две теоремы говорят, что многообразие сохраняется при любых гомеоморфных преобразованиях.

Следствие.

Гомеоморфные преобразования переводят многообразие в многообразие той же размерности.

## *КЛАССИФИКАЦИЯ МОЛЕКУЛЯРНЫХ ПРОСТРАНСТВ*

Здесь уместно еще раз подчеркнуть, что дигитальная топология как наука, которой посвящена эта книга, возникла всего несколько лет назад, и исключительно благодаря требований практики, то есть когда возникла практическая необходимость формализовать многомерные образы, с которыми работает компьютер. Поразительно, но сто лет усилий величайших математических и физических умов, включая Б. Римана, А. Эйнштейна, А. Эддингтона не привели к сколь-нибудь значительным результатам в этой области. Единственное объяснение этому то, что усилия не были велики, так как ученые до недавнего времени не видели никакой существенной необходимости в геометрии и топологии на конечном или счетном множестве точек.

В подавляющем большинстве случаев программист имеет дело с дигитальными моделями непрерывных многообразий, то есть с молекулярными многообразиями. Однако многообразия составляют только один, и сравнительно узкий, класс молекулярных пространств. Многообразия принадлежат более широкому классу нормальных пространств. Отличительные особенности класса нормальных пространств

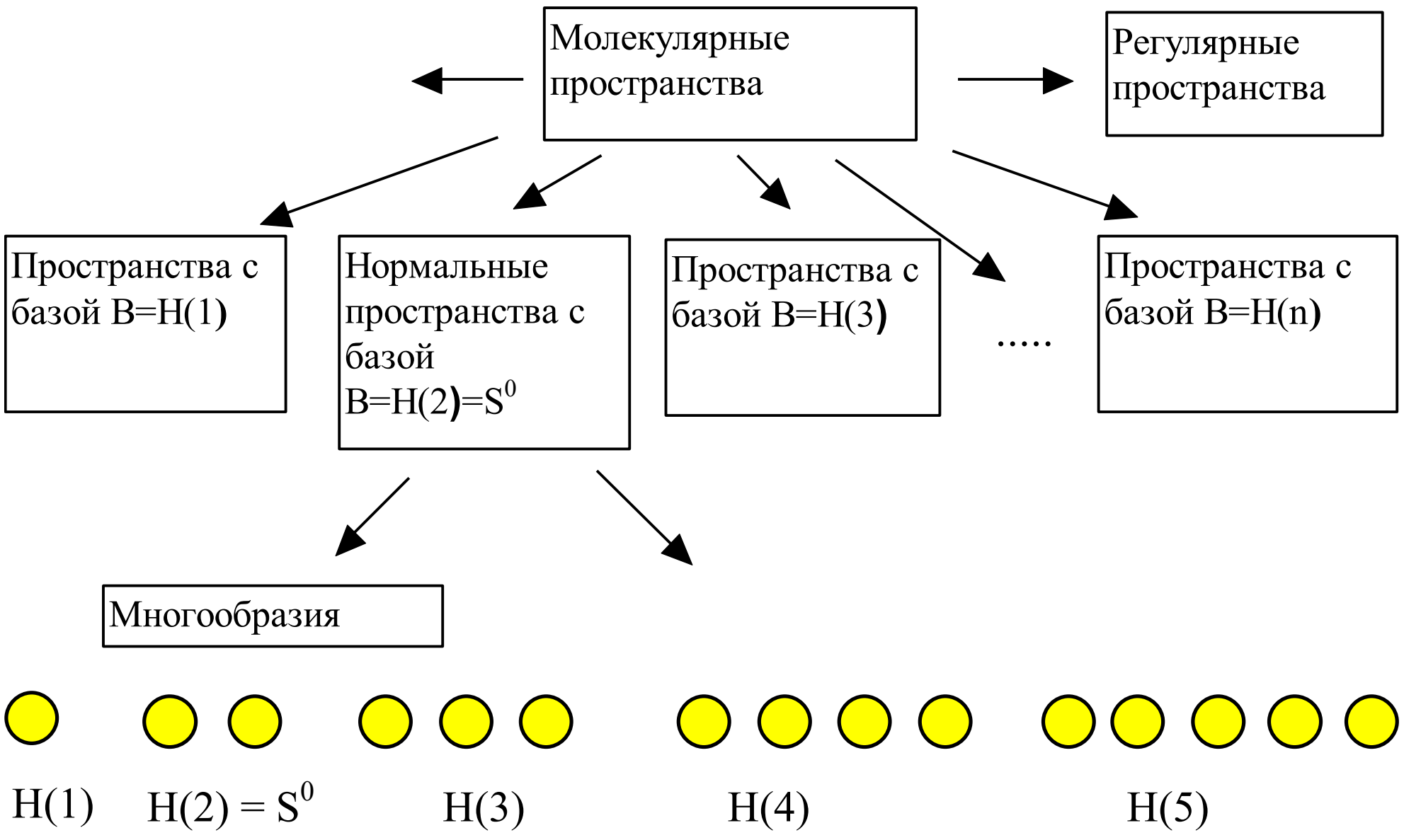

Рис. 81 Классификация молекулярных пространств. Аксиоматическое построение классов пространств на основе базового элемента H(n) 0-мерные базовые пространства H(1), H(2) = $S^0$, H(3), H(4) и H(5), содержащие 1, 2, 3, 4 и 5 точек. Пространство H(2) является базовым элементом для класса нормальных пространств.

видны из способа построения n-мерных точек и множества 0-мерных пространств. Напомним, что в классе нормальных пространств существует только два 0-мерных нормальных пространства: 0-мерная точка и 0-мерная сфера, состоящая из двух изолированных точек. Замкнутые нормальные n-мерные пространства и нормальные n-мерные точки определяются по индукции. Причем 0-мерная сфера играет роль базового элемента B, определяющего всю остальную структуру нормальных пространств. При более общем подходе определение соответствующего класса молекулярных пространств почти не зависит от того, что выбирается в качестве базового элемента. В нормальных пространствах мы выбрали сферу $S^0$ как пространство, состоящее из двух изолированных точек только потому, что это соответствует классическому определению и дает разумные и ожидаемые результаты при сравнении с непрерывными пространствами. Сразу же возникает вопрос, что произойдет, если в качестве 0-мерной сферы, то есть базового элемента, взять пространство из 3-х изолированных точек H(3).

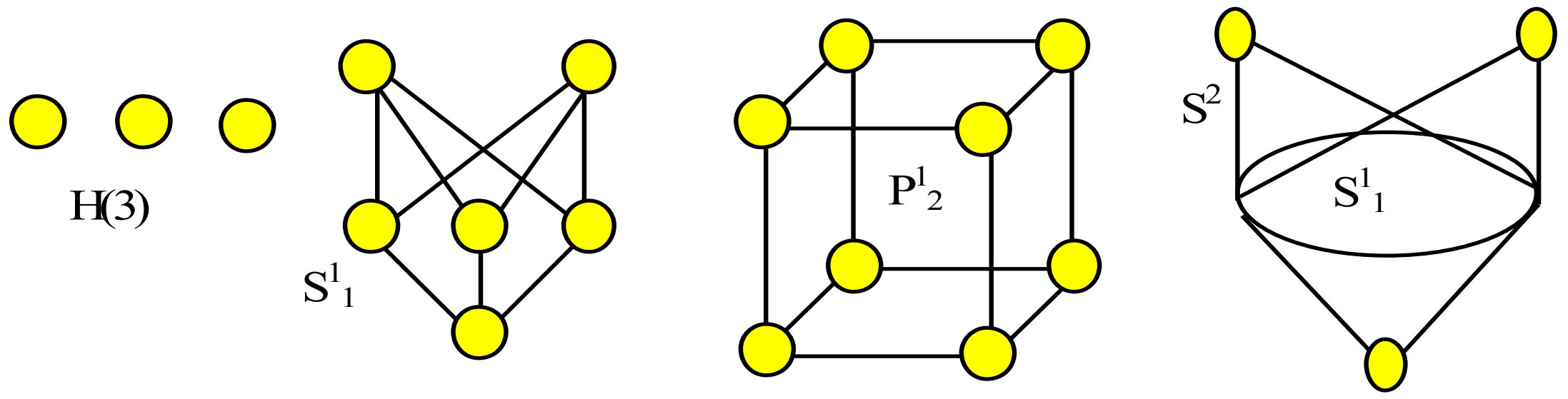

Рис. 82 0-мерное пространство H(3) используется как базовый элемент построения аналога одномерной сферы $S^1_1$=H(3)⊕H(3). $S^1_1$ также есть одномерное замкнутое пространство. $S^2$ является двумерным пространством, $S^2$=H(3)⊕$S^1_1$.

Очевидно получится другой класс пространств со своими свойствами. Исходя из можно построить грубую классификацию пространств в соответствии с базовым элементом B (Рис. 81.). Пространства с базой H(1) в какой-то мере сходны с симплициальными комплексами, изучаемыми в комбинаторной топологии. Пространства с базой H(2) образуют класс нормальных пространств. Пространства с другими базами еще не исследованы.

Важность нормальных пространств обусловлена тем, что как в реальной окружающей нас действительности, так и в науке, мы сталкивались и сталкиваемся только с многообразиями. В n-мерном многообразии каждая точка имеет окрестность, гомеоморфную открытому шару в n-мерном евклидовом пространстве. Следовательно, в молекулярном пространстве, моделирующем это многообразие, окаем каждой точки должен быть (n-1)-мерной молекулярной сферой. Если это не выполняется, то данное молекулярное пространство не может быть моделью многообразия. До недавнего времени в физике признавались только хаусдорфовы многообразия как пространства, математически моделирующие физические реалии. Однако многообразия составляют только небольшую часть молекулярных пространств. Большая часть пространств не является многообразиями. Утверждая это мы имеем ввиду следующее. Если мы построим все возможные замкнутые молекулярные пространства, содержащие n и менее точек, то окажется, что многообразия составляют очень небольшую часть общего числа различных пространств.

В ТМП замечательно то, что не-многообразия легко строятся и наглядны. Простота построения и наглядность не-многообразий в ТМП резко расширяют наши возможности в изучении о реальных объектов,

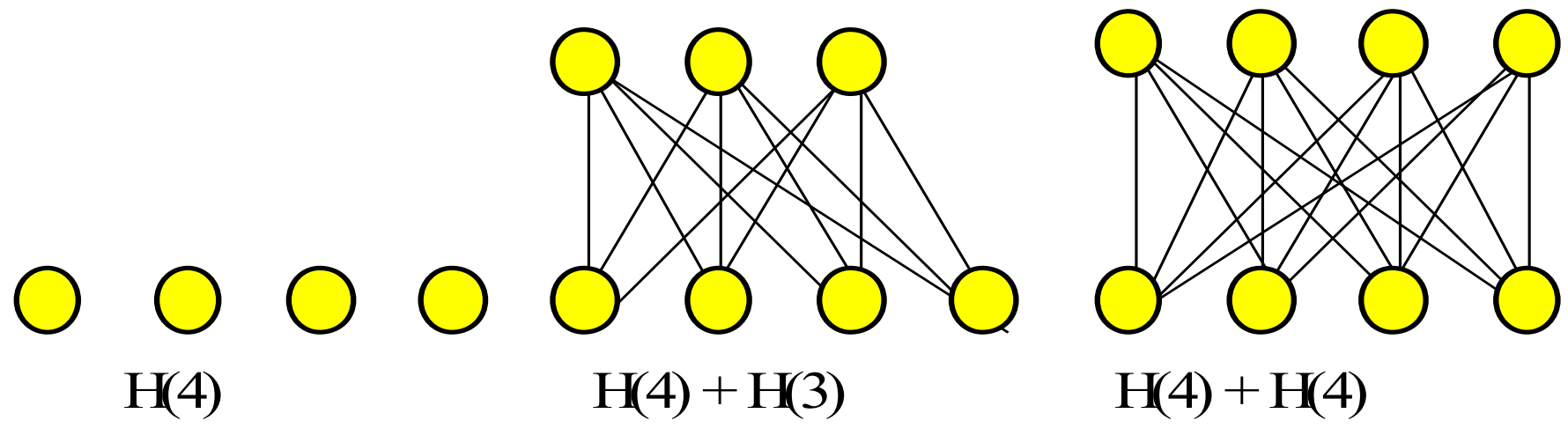

Рис. 83 Одномерные замкнутые пространства с базой H(n), n>2, содержащие 6, 7 и 8 точек.

существующих или возможно существующих в окружающем нас физическом мире. В качестве очевидного примера таких пространств можно рассмотреть так называемые структурные формулы молекул, в которых не учитываются кратность и длина связей в молекуле. Большинство структурных формул являются молекулярными пространствами, не принадлежащими к классу нормальных пространств и тем более многообразий.

Поэтому одно из привлекательных направлений для исследования-изучить свойства пространств, в которых в качестве базового элемента берется не 0-мерная сфера $S^0$=H(2), а, например, H(3) или H(4). При выборе в качестве базового элемента пространства H(3), состоящего из трех несвязных точек, получаются одномерные замкнутые пространства, изображенные на Рис. 82. В традиционном ключе базовое пространство H(3) может считаться 0-мерной сферой, пространства и $S^1{}_1$=H(3)⊕H(3) и $P^1{}_2$ являются одномерными замкнутыми пространствами, пространство $S^2$ будет аналогом двумерной сферы. На Рис. 83 пространство H(4) является базовым пространством. Пространство H(3)⊕H(4) состоит из 7 точек, три из которых имеют окаем H(3), а остальные четыре-окаем H(4). Одномерное пространство H(4)⊕H(4) состоит из 8 точек, каждая из которых имеет окаем H(4).

Среди нормальных пространств также можно выделить пространства, не являющиеся многообразиями. Они возникают, когда мы берем прямую сумму замкнутых нормальных пространств, в том числе и многообразий, одно из которых не есть сфера. Как пример такого не-многообразия, нормального пространства, полученного из двух многообразий, рассмотрим сумму нульмерной сферы $S^0$, состоящей из точек a и b, и двумерного тора $T^2$, $V^3 = S^0 \oplus T^2$.

Схема построения изображена на Рис. 84. Легко видеть, что окаемами точек a и b в $V^3$ является тор, и эти точки тороидальны. Следовательно $O(a) = O(b) = T^2$, и это не-многообразие. Окаемом каждой точки тора является двумерная сфера. Таким образом данное не-многообразие имеет две не-сферические, тороидальные точки; все остальные точки есть сферические.

Более сложные примеры-пятимерные замкнутые молекулярные нормальные пространства $T^2 \oplus T^2$, $T^2 \oplus P^2$, $T^2 \oplus K^2$ и так далее. Ни одна точка в этих пространствах не является сферической. Пока трудно сказать, каково может быть математическое или прикладное применение таких пространств. Можно только предположить, что найдутся среди химических или биологических структур объекты, являющиеся такими

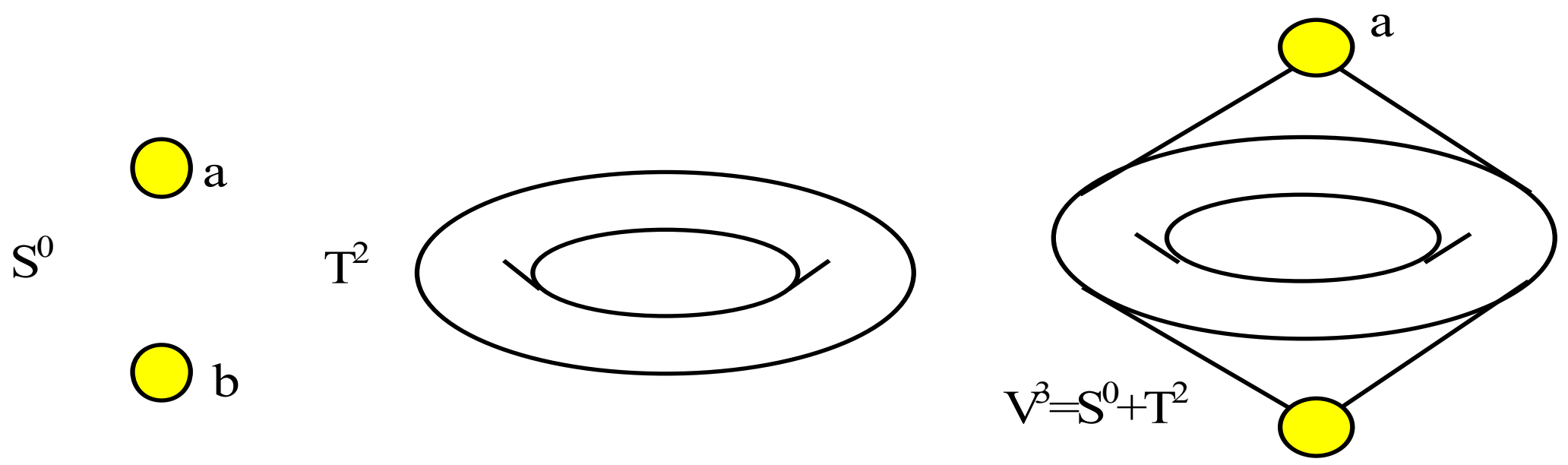

Рис. 84 Трехмерное нормальное пространство не-многообразие, полученное как сумма двух многообразий, двумерного тора и 0-мерной сферы.

экзотическими (пока) пространствами.

З а д а ч а .

Найти физические, химические, биологические или какие-либо другие объекты, являющиеся моделями не-многообразий.

Перейдем теперь к регулярным пространствам. Регулярное пространство содержит нормальное пространство, как свое подпространство. Кроме этого регулярное пространство содержит некоторое дополнительное число точек, которые не меняют топологию пространства в целом, Любую из этих дополнительных точек можно отбросить, при этом основные топологические характеристики пространства не изменятся. Пользуясь аналогиями можно сказать, что нормальное подпространство, содержащееся в регулярном пространстве, является скелетом, к которому присоединяется некоторое количество точек, позволяющее получить более полное соответствие математической или компьютерной модели реальным объектам. При общем подходе можно рассматривать регулярное пространство, в которой скелетом является пространство с произвольной вазой H(n), n>1. Однако исходя из целей данной книги, мы ограничимся регулярными пространствами, в которых скелетами являются нормальные пространства с базой H(2).

Список литературы к главе 5.

# ПРЯМАЯ СУММА МОЛЕКУЛЯРНЫХ ПРОСТРАНСТВ И ЕЕ СВОЙСТВА

In this section we study properties of the direct sum of two or more molecular spaces. We prove that the direct sum of a contractible and any space is a contractible space. The direct sum of two noncontractible spaces is a noncontactible space. The direct sum of two normal closed spaces is a normal closed space. Contractible transformations in one of the spaces are at the same time contractible transformations in the direct sum.

Результаты этой главы частично изложены в работах [7,43 ].

## *ОБЩИЕ СВОЙСТВА ПРЯМОЙ СУММЫ МОЛЕКУЛЯРНЫХ ПРОСТРАНСТВ*

Прямая сумма молекулярных пространств позволяет получать новые пространства на основе уже известных. Повторим уже сделанные определения.

<u>Определение прямой суммы пространств G⊕H.</u>

Пусть имеется два пространства: G с набором точек $V=(v_1, v_2, v_3, .... v_n)$ и набором связей W, и H с набором точек $R=(r_1, r_2, r_3, .... r_m)$ и набором связей S. Прямой суммой (суммой) двух пространств G и H назовем пространство G⊕H состоящее из точек V и R, связей W и S, и всех связей, соединяющих точки V и R

Как видно, при этой операции число точек пространства G⊕H равно сумме точек пространств G и H, |G⊕H|=|G|+|H|, число связей равно |W|+|S|+nm.

<u>Определение конуса G⊕v.</u>

Если одно из пространств, например, второе в предыдущем определении является точкой v, то сумма G⊕v называется конусом (cone) пространства G (Рис. 85).

<u>Определение прямой суммы подпространств G⊕H.</u>

Пусть в пространстве A имеется два подпространства: G с набором точек $V=(v_1, v_2, v_3, .... v_n)$, и H с набором точек $R=(r_1, r_2, r_3, .... r_m)$, причем множества V и H не пересекаются. Если каждая точка из V смежна с каждой точкой из R, то подпространство G∪H=G⊕H, состоящее из точек V и R, называется прямой суммой (суммой) двух подпространств G и H.

Например, шар точки есть подпространство, являющееся прямой суммой точки и ее окаема.

Теорема 52

*Пусть G и H есть точечное и произвольное МП. Тогда прямая сумма*

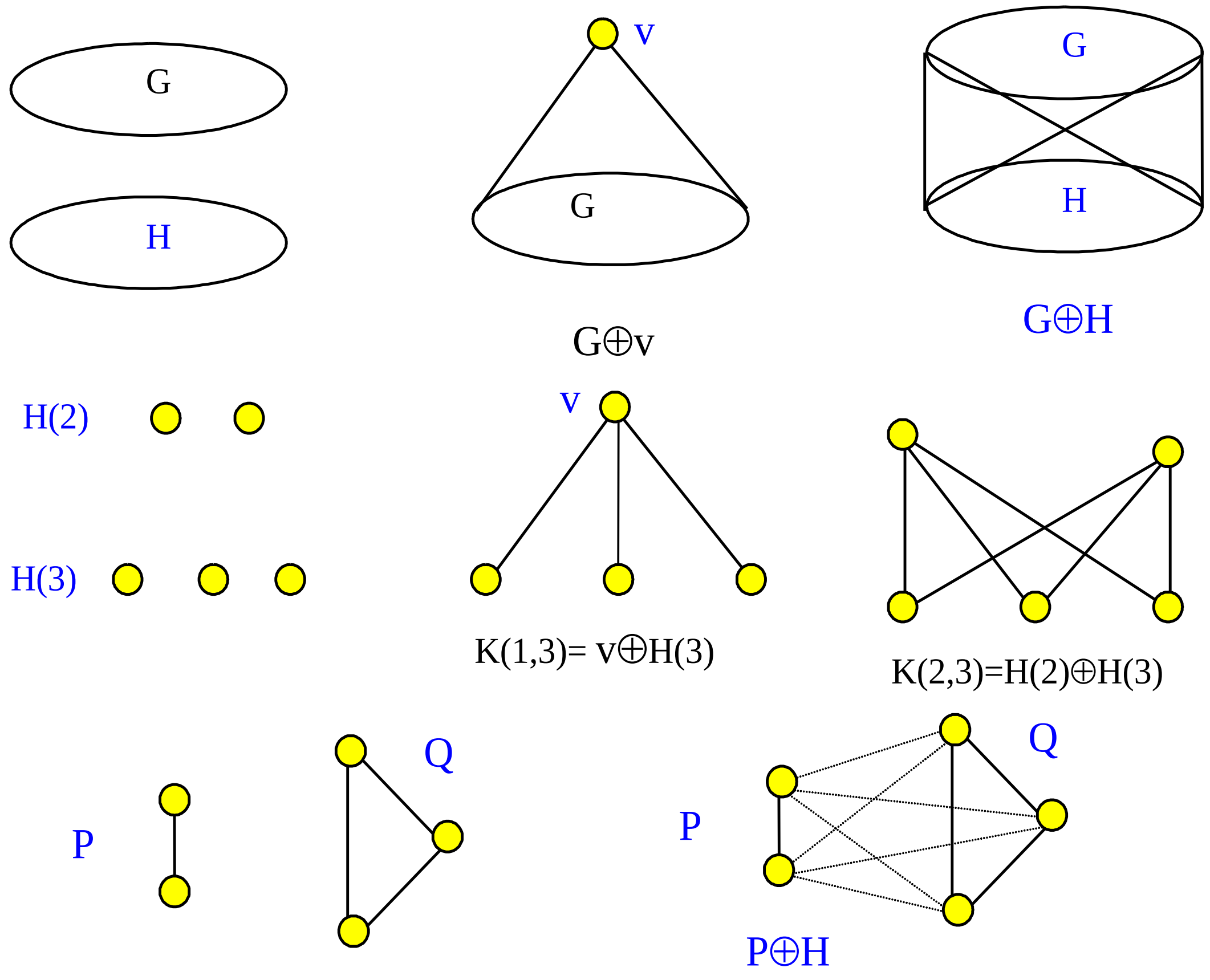

Рис. 85 Прямые суммы пространств G⊕v, G⊕H, v⊕H(3), H(2)⊕H(3) и P⊕Q.

*G⊕H есть точечное МП.*

Доказательство.

Проведем доказательство по индукции. Для МП G с малым числом точек теорема проверяется непосредственно. Предположим, что теорема верна для |G|<n. Пусть |G|=n. Рассмотрим W=G⊕H. Так как G есть точечное МП, то он имеет по крайней мере одну точечную точку $v_1$, $O(v_1)\in T$. В МП W окаем $OW(v_1)$ точки $v_1$ определяется выражением $OW(v_1)=O(v_1)\oplus H$. Так как $|O(v_1)|<n$, то согласно индукции OW $(v_1)$ будет точечным, и точка $v_1$ может быть отброшена в W. $W-v_1$ снова есть прямая сумма точечного $G-v_1$ и H, где $|G-v_1|<n$. Согласно предположению индукции $W-v_1=(G-v_1)\oplus H$ является точечным. Доказательство закончено.□

Теорема 53

*Пусть G и H есть два неточечных пространства, в каждом из которых нет точек с точечными окаемами. Тогда прямая сумма*

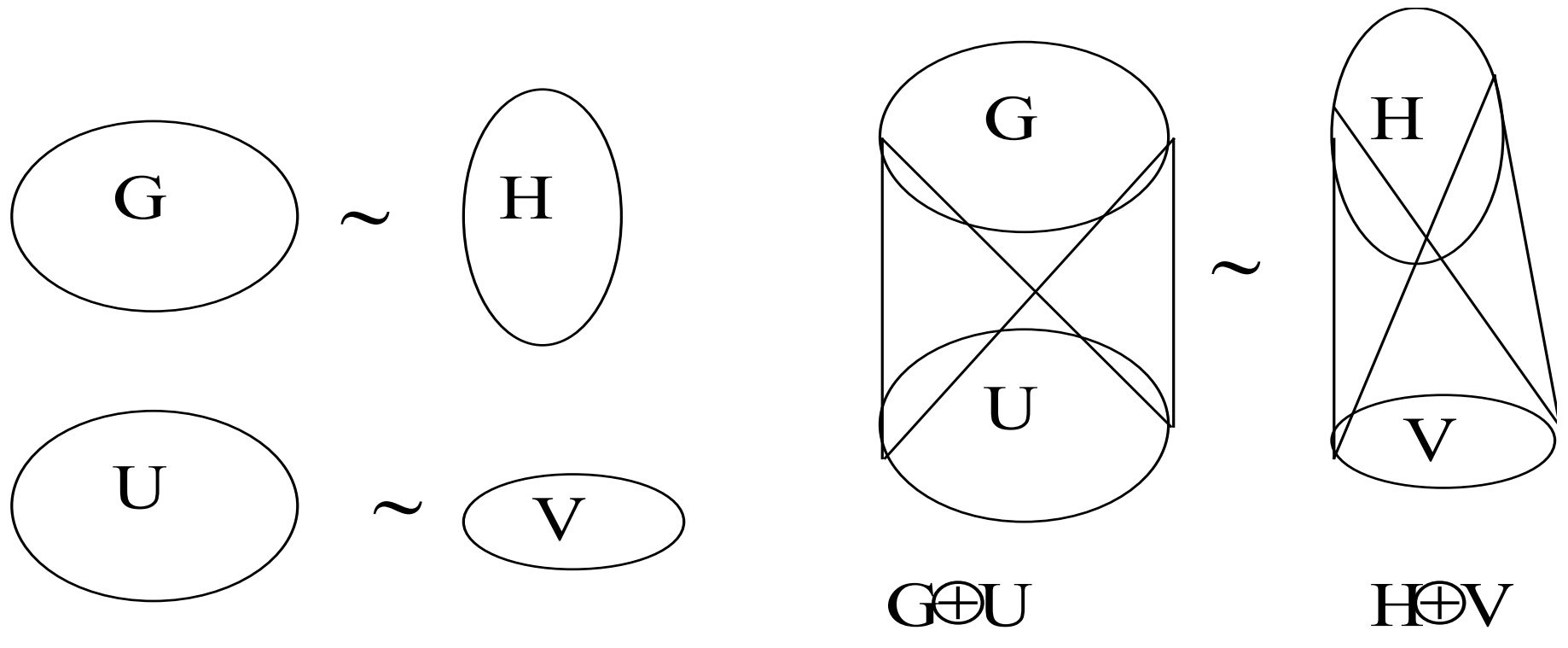

Рис. 86 G гомотопно H, U готопно V. Тогда G⊕U гомотопно H⊕V.

*G⊕H есть неточечное МП., в котором нет точек с точечными окаемами.*

Доказательство.

Проведем доказательство по индукции. Для МП G с малым числом точек теорема проверяется непосредственно. Предположим, что теорема верна для |G⊕H|<n. Пусть |G⊕H|=n. В W=G⊕H OW(v)=O(v)⊕H, где v∈G. Следовательно |OW(v)|<n, и из индукции O(v) не является точечным. То же самое справедливо для любой точки, принадлежащей H. Теорема доказана. □

Прямая иллюстрация этой теоремы содержится в прямой сумме нормальных замкнутых пространств.

Теорема 54

*Пусть G и H есть два молекулярных пространства. Точечное преобразование в G является точечным преобразованием в G⊕H.*

Доказательство.

Пусть v∈G и O(v) есть точечное пространство в G. Тогда $O_1(v)$=O(v)⊕H в G⊕H и, следовательно, является также точечным пространством. Это означает, что точка v может быть отброшена (или присоединена) из G⊕H. Пусть $v_1$∈G и $v_2$∈G и $O(v_1v_2)$ есть точечное пространство в G. Тогда $O_1(v_1v_2)=O(v_1v_2)$⊕H в G⊕H и, следовательно, является также точечным пространством. Это означает, что связь между точками может быть отброшена (или присоединена) из G⊕H. Теорема доказана. □

Теорема 55

*Пусть G гомотопно H и U гомотопно V. Тогда G⊕U гомотопно H⊕V.*

Доказательство.

Так как G~U, H~V, то существует две последовательности точечных преобразований, переводящих G в U и H в V. Эти же преобразования с соответствующей заменой окаемов являются точечными в G⊕U. Применив эти преобразования к G⊕U получим H⊕V, G⊕U~H⊕V (Рис. 86). Теорема доказана. □

Определение подпространства [A⊕B] прямой суммы G⊕H.

Пусть A есть подпространство пространства G, B есть подпространство пространства H, тогда [A⊕B] есть подпространство пространства G⊕H, содержащее точки $v_s$ и $u_t$,, где $v_s$ принадлежит A, $u_t$ принадлежит B.

Теорема 56

*Пусть A есть подпространство пространства G, B есть подпространство пространства H, тогда подпространство [A⊕B], где [A⊕B]⊆G⊕H изоморфно пространству A⊕B, [A⊕B]≈A⊕B.*

Доказательство

Нужно показать, что между точками U=[A⊕B] и V=A⊕B можно установить взаимно-однозначное соответствие, сохраняющее смежность точек. Установим соответствие естественным образом, то есть точка $v_k$=[$v_k$] принадлежащая A в [A⊕B], соответствует точке $v_k$, принадлежащей A в A⊕B, $u_p$=[$u_p$] принадлежащая B в [A⊕B], соответствует точке $u_p$, принадлежащей B A⊕B. Окаем OU([$v_k$])=(O([$v_k$])∩A)⊕B в [A⊕B]. Точно так же OV($v_k$)=(O($v_k$)∩A)⊕B в A⊕B. То же самое применимо к точкам из B. Следовательно, взаимно-однозначное соответствие между точками из [A⊕B] и A⊕B, сохраняет смежность точек. Теорема доказана.□

В дальнейшем, если это не оговорено специально, мы не будем различать [A⊕B] и A⊕B, и для [A⊕B] будем использовать обозначение A⊕B.

Теорема 57

*Пусть G и H есть неточечные пространства, и U и V являются их вложенными гомотопными подпространствами соответственно.*

*Тогда U⊕V является гомотопным и вложенным подпространством пространства G⊕H.*

Доказательство.

Так как G и H стягиваются точечными отбрасываниями точек к U и V, то G⊕H стягивается путем отбрасывания точечных точек в U⊕V. Следовательно U⊕V есть вложенное гомотопное подпространство пространства G⊕H. Теорема доказана. □

Сформулируем еще одну теорему о свойствах суммы, которую доказывать

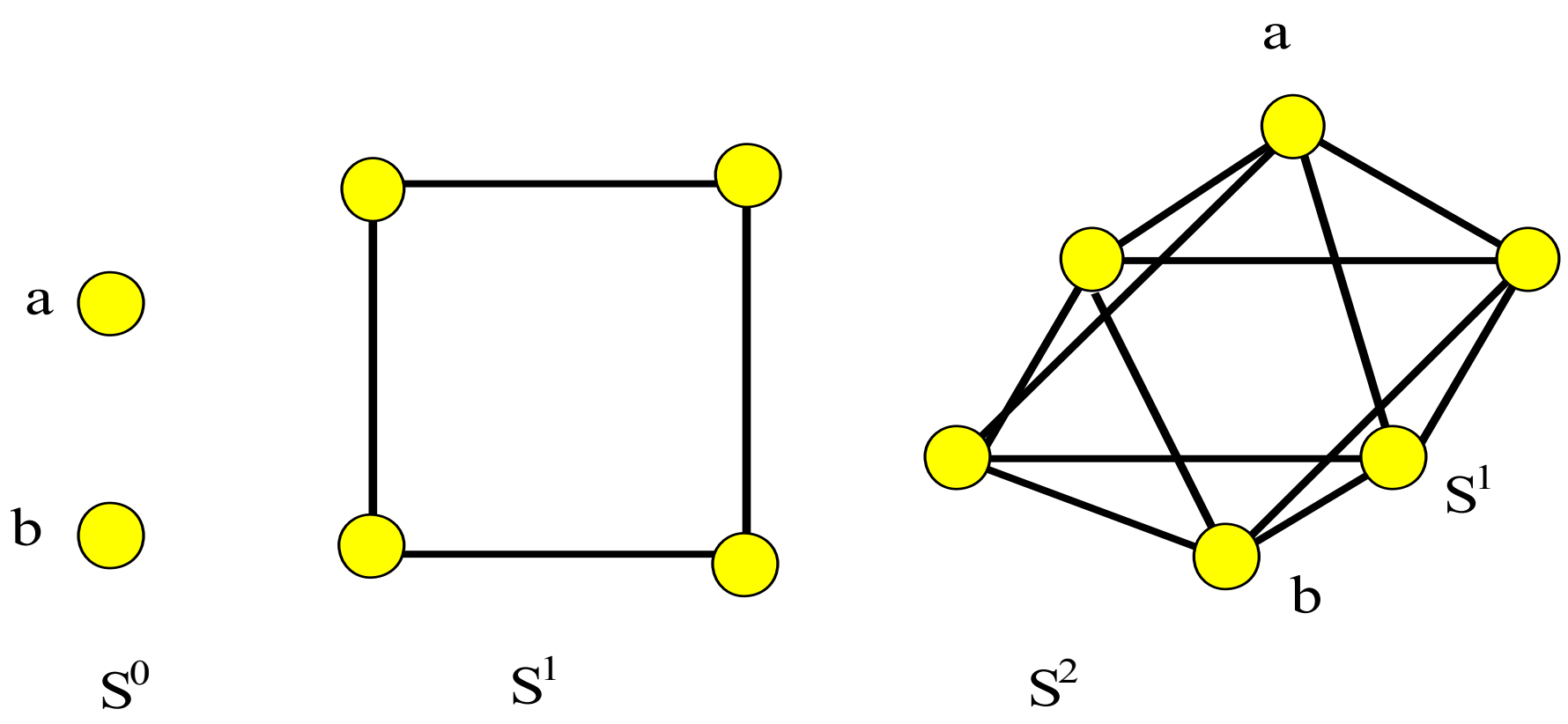

Рис. 87 Сумма 0-сферы $S^0$ и 1-сферы $S^1$ является 2-мерной сферой. $S^2=S^0\oplus S^1=S^1\oplus S^0$.

не будем в силу ее очевидности.

Теорема 58

*Сумма молекулярных пространств коммутативна и ассоциативна. W=G⊕H=H⊕G. W=(G⊕H)⊕F=G⊕(H⊕F).*

Из Рис. 87 видно, что сумма одномерной и нульмерной сфер является двумерной сферой. Следует отметить то обстоятельство, что сумма минимальных сфер также дает минимальную сферу. Не исключено, что этим свойством обладает прямая сумма любых минимальных пространств. В связи с этим сформулируем в виде теоремы следующее предположение.

Теорема 59

*Пусть G и H-два нормальных молекулярных замкнутых минимальных пространства. Тогда их сумма W=G⊕H также будет минимальным пространством.*

Доказательство.

Самостоятельно доказать или опровергнуть.□

Из определения прямой суммы следует, что пространства-сомножители должны быть пространствами конечного объема. В противном случае окаемы некоторых точек будут пространствами бесконечного объема.

## *СВОЙСТВА ПРЯМОЙ СУММЫ НОРМАЛЬНЫХ МОЛЕКУЛЯРНЫХ ПРОСТРАНСТВ*

Рассмотрим свойства прямой суммы нормальных пространств. Эта операция позволяет получить молекулярные пространства с новыми топологическими свойствами. Размерное свойство такой операции аналогично свойству соответствующей операции в классической топологии.

Теорема 60

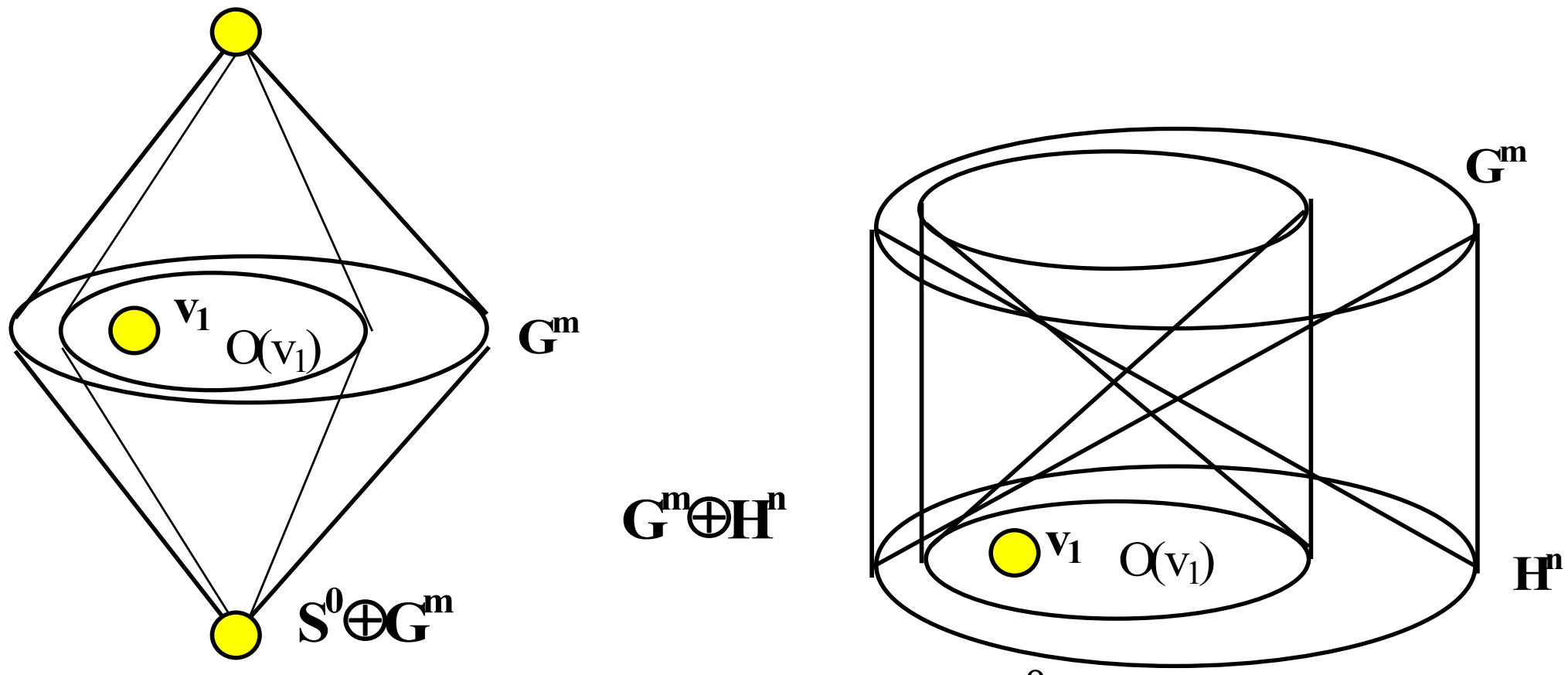

Рис. 88 Прямая сумма 0-мерной сферы $S^0$ и замкнутого n-мерного пространства $G^n$ имеет размерность (n+1). Прямая сумма замкнутого n-пространства $G^n$ и замкнутого m-пространства $H^m$ имеет размерность (n+m+1).

*Пусть $G^m$ и $S^0$-нормальное молекулярное замкнутое пространство размерности m и 0-мерная сфера с набором точек $V=(v_1,v_2, v_3,....v_q)$ и $R=(a,b)$ соответственно. Тогда их прямая сумма $W^{m+1}=G^m \oplus S^0$ является нормальным замкнутым молекулярным пространством размерности (m+1).*

Д о к а з а т е л ь с т в о .

Используем индукцию.

Для m=0, 1 теорема проверяется непосредственно. Предположим, что теорема верна для $G^m$, где $m \le k$. Пусть m=k+1. Возьмем любую точку из $G^m$, скажем, $v_1$. Окаем этой точки в пространстве $W^{m+1}$ есть сумма окаема этой точки в пространстве $G^m$ и пространства $S^0$, то есть $O(v_1)|W^{m+1}=(O(v_1)|G^m)\oplus S^0$. Так как окаем этой точки в пространстве $G^m$ является к-мерным пространством, то, согласно индукции, $O(v_1)|W^{m+1}$ будет нормальным к+1 мерным пространством. Окаемы точек a и b также являются нормальным к+1 мерным пространством $G^m$. Следовательно, $W^{m+1}=G^m\oplus S^0$ есть нормальное к+2 мерное пространство (*Рис. 88*). Теорема доказана.□

На *Рис. 89* иллюстрируется сумма двух 0-мерных сфер, дающая одномерную сферу, то есть окружность, $S^1=S^0 \oplus S^0$. Следующая теорема обобщает предыдущую на сумму произвольных нормальных замкнутых пространств.

Теорема 61

*Пусть $G^m$ и $H^n$-два нормальных молекулярных замкнутых пространства размерности m и n с набором точек $V=(v_1, v_2, v_3,....$*

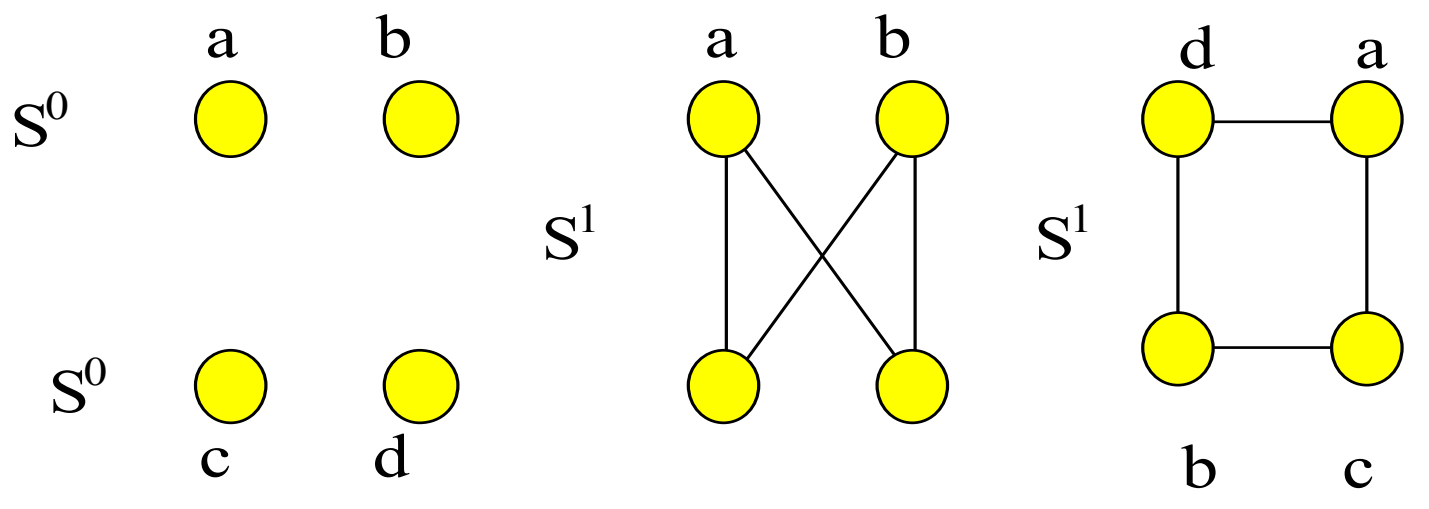

Рис. 89 Сумма двух 0-мерных сфер является одномерной окружностью, $S^1=S^0\oplus S^0$.

*$v_q)$ и $R=(r_1, r_2, r_3,.... r_t)$ соответственно. Тогда их сумма $W^{m+n+1}=G^m\oplus H^n$ является нормальным замкнутым молекулярным пространством размерности m+n+1.*

$$W^{m+n+1}=G^m\oplus H^n.$$

Доказательство.

Используем индукцию.

Для n, m=0, 1 теорема проверяется непосредственно.

Предположим, что теорема верна для $G^m$ и $H^n$, где $m+n \le p$. Пусть $n+m=p+1$. Так как теорема уже доказана для m=0 или n=0, рассмотрим случай, когда m>0 или n>0.

Выберем произвольную точку $v_1 \in G^m$ и рассмотрим окаем этой точки в $W^{m+n+1}=G^m\oplus H^n$. Очевидно, что $O(v_1)|W^{m+n+1}=(O(v_1)|G^m)\oplus H^n$, и раз$(O(v_1)|G^m)=m-1$, раз$(H^n)=n$. Следовательно, по предположению, их сумма является нормальным замкнутым молекулярным пространством размерности m+n, раз($O(v_1)|W^{m+n+1}$)=раз($O(v_1)|G^m$)+раз($H^n$)+1= (m-1)+n+1=m+n. Это означает, что размерность пространства $W^{m+n+1}=G^m\oplus H^n$ равна m+n+1. Теорема доказана.□

Чтобы сделать изложение законченным, напомним те свойства прямого произведения нормальных пространств, которые были уже рассмотрены ранее.

Теорема 62

*Пусть $G^n$, n>0, является нормальным n-мерным замкнутым пространством. Тогда конус $v\oplus G^n$ будет нормальным (n+1)-мерным пространством со связным краем $G^n$.*

Теорема 63

*Пусть $G^n$ является нормальным n-мерным пространством на конечном множестве точек с ядром A и со связным краем B, $G^n=A\cup B$. Пусть $H^m$ является m-мерным нормальным замкнутым пространством. Тогда их прямая сумма $G^n\oplus H^m$ будет нормальным (n+m+1)-мерным пространством с ядром A и со связным краем $B\oplus H^m$, являющимся (n+m)-мерным нормальным замкнутым пространством.*

Список литературы к главе 6.

# ПРЯМОЕ ПРОИЗВЕДЕНИЕ МОЛЕКУЛЯРНЫХ ПРОСТРАНСТВ И ЕГО СВОЙСТВА

In this chapter we study properties of the direct product of two or more molecular spaces. We prove that the direct product of a noncontractible and any space is a noncontractible space. The direct product of two contractible spaces is a contactible space. Contractible transformations in one of the spaces generate at the same time contractible transformations in the direct product.

Результаты этой главы частично изложены в работе [43].

## *ОПРЕДЕЛЕНИЕ И ОБЩИЕ СВОЙСТВА ПРЯМОГО ПРОИЗВЕДЕНИЯ МОЛЕКУЛЯРНЫХ ПРОСТРАНСТВ*

Топологическое произведение двух пространств используется в математике довольно широко [53]. Это понятие позволяет рассматривать плоскость как произведение двух прямых, а тор как произведение двух окружностей. Аналогичная операция уже была нами введена для молекулярных пространств. Для полноты изложения повторим еще раз

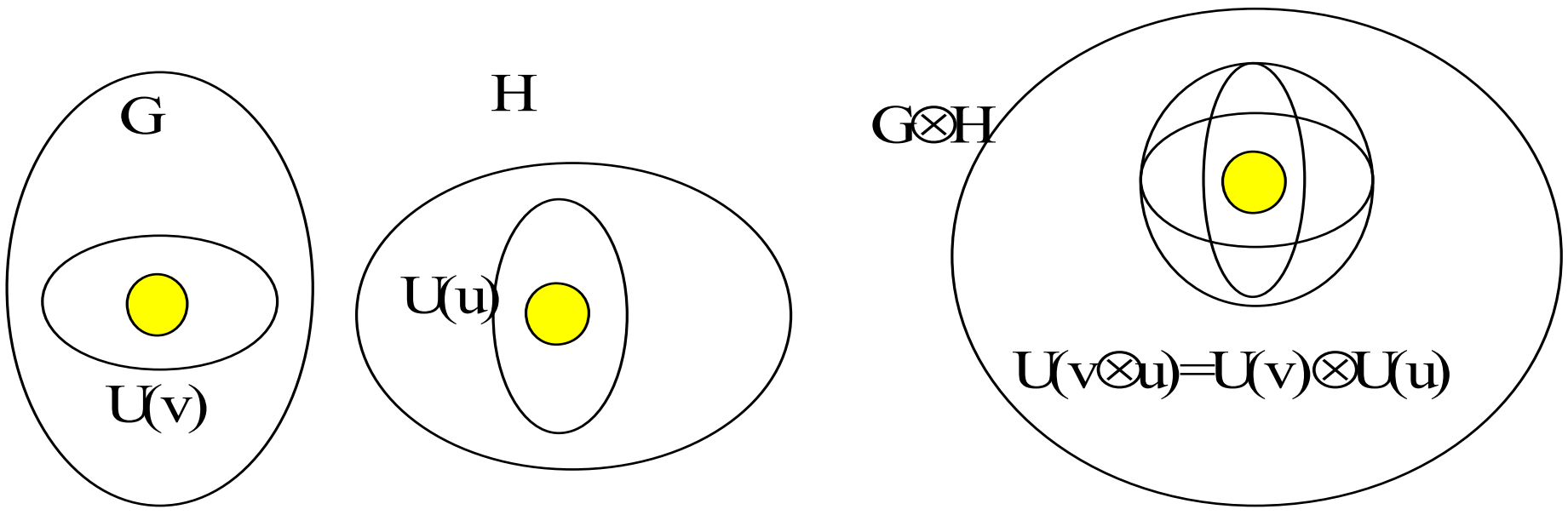

Рис. 90 В прямом произведении G⊗H шар U(v⊗u) точки v⊗u равен прямому произведению U(v)⊗U(u) шаров U(v) и U(u). U(v⊗u)=U(v)⊗U(u).

определение прямого произведения молекулярных пространств.

Определение прямого произведения G⊗H.

> Пусть имеется два пространства: G с набором точек ($v_1$, $v_2$, $v_3$,.... $v_m$...) и набором связей W, и H с набором точек ($u_1$, $u_2$, $u_3$,.... $u_n$...) и набором связей S. Прямым произведением (произведением) двух пространств G и H назовем пространство G⊗H, состоящее из точек $w_{kp}=v_k\otimes u_p$, в котором шар $U(v_k\otimes u_p)$ каждой точки $w_{kp}=v_k\otimes u_p$ состоит из точек $w_{st}=v_s\otimes u_t$,, где $v_s$ принадлежит $U(v_k)\subseteq$,G, $u_t$ принадлежит $U(u_p)\subseteq H$.

Следует напомнить, что по определению шара точки точка $v_k\otimes u_p$ является смежной со всеми остальными точками, принадлежащими этому шару.

Легко видеть из этого определения, что топологические и иные свойства прямого произведения двух пространств в окрестности каждой точки определяются через локальные свойства каждого из пространств. Рассмотрим более подробно структуру окрестности точки в прямом произведении двух пространств. Прежде всего докажем, что прямое произведение обладает неким наследственным свойством, а именно, шар любой точки в прямом произведении сам является прямым произведением шаров точек.

Теорема 64

*Шар $U(v_k \otimes u_p)$ любой точки $w_{kp}=v_k \otimes u_p$ прямого произведения двух пространств G с набором точек ($v_1$, $v_2$, $v_3$,.... $v_n$..) и H с набором точек ($u_1$, $u_2$, $u_3$,.... $u_m$...) есть прямое произведение шаров $U(u_k)$ и $U(v_p)$, $U(v_k \otimes u_p)=U(v_k) \otimes U(u_p)$ (Рис. 90).*

Доказательство

Рассмотрим подпространство $U(v_k \otimes u_p)$. В соответствии с определением для точки $v_k \otimes u_p$ шар этой точки состоит точек $v_s \otimes u_t \in U(v_k \otimes u_p)$, где $v_s \in U(v_k)$, $u_t, \in U(u_p)$. Пусть точка $v_s \otimes u_t$ принадлежит $U(v_k \otimes u_p)$. Ее шар состоит из точек вида $v_m \otimes u_n \in U(v_s \otimes u_t)$, где $v_m \in U(v_s)$, $u_n, \in U(u_t)$. Следовательно, в пространстве $U(v_k \otimes u_p)$ шар $U(v_s \otimes u_t)$ точки $v_s \otimes u_t$ состоит из точек $v_m \otimes u_n$, для которых $v_m \in U(v_s) \cap U(v_k)=U(v_s v_k)$, $u_n \in U(u_t) \cap U(u_p)=U(u_t u_p)$. Подпространства $U(v_s) \cap U(v_k)=U(v_s v_k)$ и

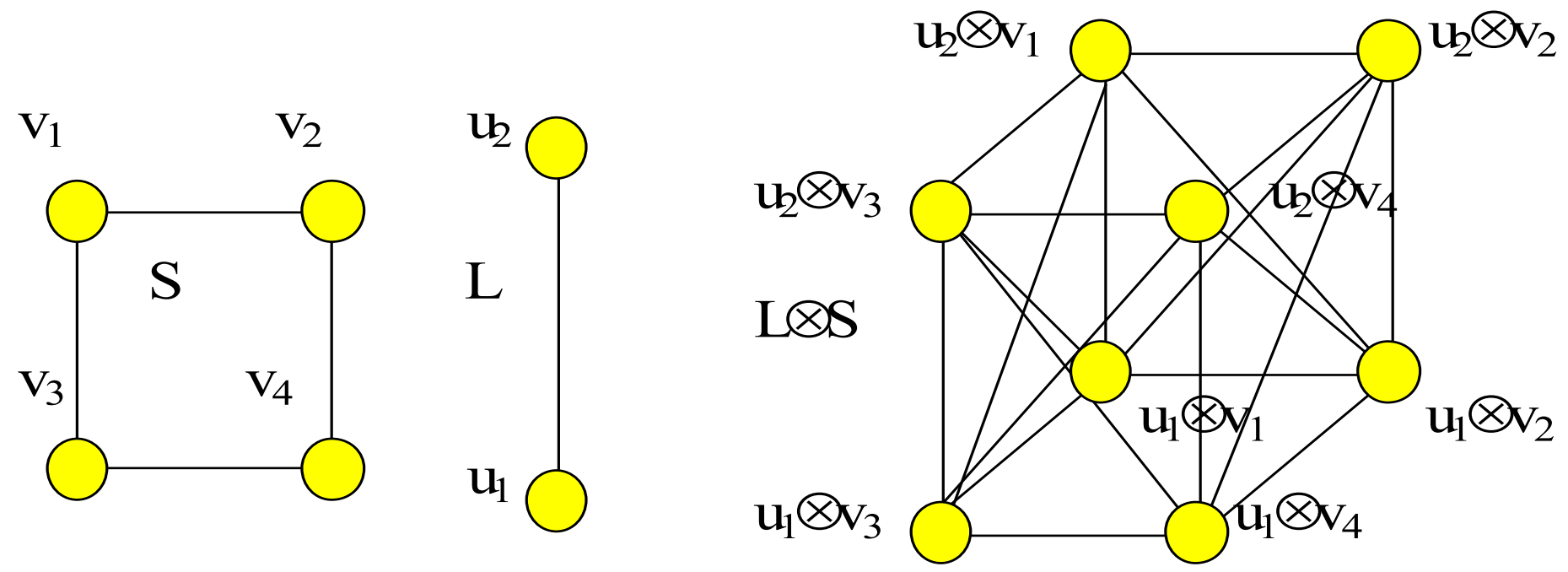

Рис. 91 Прямое произведение L⊗S связки L и окружности S.

$U(u_t) \cap U(u_p)=U(u_t u_p)$ есть не что иное, как шары точек $v_s$ и $v_t$ в $U(v_k)$ и $U(v_k)$. Это, по определению, означает, что $U(v_k \otimes u_p)=U(v_k) \otimes U(u_p)$. Теорема доказана.□

Эта теорема, очевидно, обобщается на любое число сомножителей.

Определение координат точки и слоя в прямом произведении G⊗H.

Пусть имеется прямое произведение $G_1 \otimes G_2 \otimes ... \otimes G_k$ пространств $G_1$, $G_2$,... $G_k$, содержащих точки ($v_1$, $v_2$, $v_3$,.... $v_m$...), ($u_1$, $u_2$, $u_3$,.... $u_m$...),... ($w_1$, $w_2$, $w_3$,.... $w_m$...) соответственно. Тогда для любой точки $v_k \otimes u_p ... \otimes w_s \in G_1 \otimes G_2 \otimes ... \otimes G_k$ символы $v_k$,$u_p$,... $w_s$ называются ее координатами в подпространствах $G_1$, $G_2$,... $G_k$ соответственно. Подпространство $S_1 = G_1 \otimes u_p ... \otimes w_s \subseteq G_1 \otimes G_2 \otimes ... \otimes G_k$, определяемое точками $v_k \otimes u_p ... \otimes w_s \in G_1 \otimes G_2 \otimes ... \otimes G_k$, где координаты $v_k$ могут меняться, а остальные координаты фиксированы, называется слоем в $G_1 \otimes G_2 \otimes ... \otimes G_k$. Точно так же определяются любые другие слои в прямом произведении.

Теорема 65

*Слой $G_1 \otimes u_p ... \otimes w_s \subseteq G_1 \otimes G_2 \otimes ... \otimes G_k$ изоморфен пространству $G_1$.*

Доказательство

Очевидно взаимно-однозначное соответствие между точками слоя и $G_1$. □

Теорема 66

*Пусть имеется прямое произведение $G_1 \otimes G_2 \otimes ... \otimes G_k$ пространств $G_1$, $G_2$,... $G_k$ содержащих точки ($v_1$, $v_2$, $v_3$,.... $v_m$...), ($u_1$, $u_2$, $u_3$,.... $u_m$...),... ($w_1$, $w_2$, $w_3$,.... $w_m$...) соответственно. Тогда шар любой точки из прямого произведения определяется выражением $U(v_k \otimes u_p ... \otimes w_s) = U(v_k) \otimes U(u_p) ... \otimes U(w_s)$.*

Доказательство

Теорема сразу же доказывается по индукции..□

На Рис. 91 показано прямое произведение связки L, состоящей из двух точек и окружности S, состоящей из 4-х точек. Рассмотрим несколько общих свойств прямого произведения пространств, используемых в дальнейшем изучении.

Определение подпространства [A⊗B] прямого произведения G⊗H.

Пусть A есть подпространство пространства G, B есть подпространство пространства H, тогда [A⊗B] есть подпространство пространства G⊗H, содержащее точки vs⊗ut,, где vs принадлежит A, ut принадлежит B.

Теорема 67

*Пусть A есть подпространство пространства G, B есть подпространство пространства H, тогда подпространство [A⊗B], где $[A \otimes B] \subseteq G \otimes H$ изоморфно пространству A⊗B, $[A \otimes B] \approx A \otimes B$.*

Доказательство

Нужно показать, что между точками [A⊗B] и A⊗B можно установить взаимно-однозначное соответствие, сохраняющее смежность точек. Установим соответствие естественным образом, то есть точка $v_k \otimes u_p$, принадлежащая [A⊗B], соответствует точке $v_k \otimes u_p$, принадлежащей A⊗B. Пусть точка $v_k \otimes u_p$ смежна с точкой $v_r \otimes u_s$, в [A⊗B]. Это означает, что $v_r \otimes u_s \in U(v_k \otimes u_p)$. Следовательно, $v_r \in U(v_k) \cap A$, $u_s \in U(u_p) \cap B$. Так как $U(v_k) \cap A$ и $U(u_p) \cap B$ являютя шарами точек $v_k$ и $u_p$ в A и B соответственно, то $v_r \otimes u_s \in U(v_k \otimes u_p)$ в A⊗B, и эти две точки смежны в A⊗B. Точно так же доказывается, что если пара точек смежна в A⊗B, то она будет смежной в [A⊗B]. Теорема доказана.□

В дальнейшем, если это не оговорено специально, мы не будем различать [A⊗B] и A⊗B, и для [A⊗B] будем использовать обозначение A⊗B. Рассмотрим свойства объединения двух подпространств.

Теорема 68

*Пусть $A_1$ и $A_2$ есть подпространства пространства G, B есть подпространство пространства H, тогда $(A_1 \cup A_2) \otimes B = (A_1 \otimes B) \cup (A_2 \otimes B)$.*

Доказательство

Согласно определению $(A_1 \cup A_2)$ есть подпространство пространства G, содержащее точки как из $A_1$, так и из $A_2$. Согласно предшествующей теореме $(A_1 \cup A_2) \otimes B$ есть подпространство пространства G⊗H. Подпространство ($A_1$(A2)(B ñîñòîèò èç òî÷åê v(u, ãäå v ïðèíàäëåæèò èëè A1 èëè A2, è u ïðèíàäëåæèò B. Ïîäïðîñòðàíñòâî (A1(B)((A2(B) òàêæå åñòü ïîäïðîñòðàíñòâî ïðîñòðàíñòâà G(H è ñîñòîèò èç òî÷åê v(u, ïðèíàäëåæàùèõ èëè (A1(B) èëè (A2(B). Ïóñòü v(u((A1(B)((A2(B). Òîãäà v(u((A1(B) èëè v(u((A2(B). Îòñþäà ñëåäóåò, ÷òî, v(A1 èëè v(A2, è u(B. Ýòî îçíà÷àåò, ÷òî v((A1(A2). Ñëåäîâàòåëüíî, ìíîæåñòâî òî÷åê â îáîèõ ïîäïðîñòðàíñòâàõ ñîâïàäàåò. Íåò íåîáõîäèìîñòè ðàññìàòðèâàòü ìíîæåñòâî ñâÿçåé, òàê êàê ýòè ïîäïðîñòðàíñòâà ïðèíàäëåæàò îáà ïðîñòðàíñòâó G(H, ñîäåðæàò îäíî è òî æå ìíîæåñòâî òî÷åê è, ñëåäîâàòåëüíî, ñîâïàäàþò. Òåîðåìà äîêàçàíà.(

Òåîðåìà 69

*Пусть $A_1$ и $A_2$ есть подпространства пространства G, B есть подпространство пространства H, тогда $(A_1 \cap A_2) \otimes B = (A_1 \otimes B) \cap (A_2 \otimes B)$.*

Доказательство

Согласно предшествующей теореме $(A_1 \cap A_2) \otimes B$ есть подпространство пространства G⊗H. Подпространство $(A_1 \cap A_2) \otimes B$ состоит из точек $v \otimes u$, где v принадлежит как $A_1$ так и $A_2$, u принадлежит B. Подпространство $(A_1 \otimes B) \cap (A_2 \otimes B)$ также есть подпространство пространства $G \otimes H$ и состоит из точек $v \otimes u$, принадлежащих одновременно как $(A_1 \otimes B)$ так и $(A_2 \otimes B)$. Пусть та же точка $v \otimes u \in (A_1 \otimes B) \cap (A_2 \otimes B)$. Тогда $v \otimes u \in (A_1 \otimes B)$ и $v \otimes u \in (A_2 \otimes B)$. Отсюда следует, что, $v \in A_1$, $v \in A_2$, $u \in B$. Это означает, что $v \in (A_1 \cap A_2)$. Следовательно, множество точек в обоих подпространствах совпадает. Нет необходимости рассматривать множество связей, так как эти подпространства принадлежат оба пространству $G \otimes H$, содержат одно и то же множество точек и, следовательно, совпадают. Теорема доказана.□
Несложно обобщить предыдущие теоремы.

Теорема 70

*Пусть $A_1$ и $A_2$ есть подпространства пространства G, $B_1$ и $B_2$ есть подпространства пространства H. Тогда в пространстве $G \otimes H$ выполняется соотношение*

$$(A_1 \cap A_2) \otimes (B_1 \cap B_2) = (A_1 \otimes B_1) \cap (A_2 \otimes B_2) = (A_2 \otimes B_1) \cap (A_1 \otimes B_2).$$

$$(A_1 \cup A_2) \otimes (B_1 \cup B_2) = (A_1 \otimes B_1) \cup (A_2 \otimes B_2) \cup (A_2 \otimes B_1) \cup (A_1 \otimes B_2).$$

Доказательство

Очевидно, что $(A_1 \cap A_2) \otimes (B_1 \cap B_2)$, $(A_1 \otimes B_1) \cap (A_2 \otimes B_2)$ и $(A_2 \otimes B_1) \cap (A_1 \otimes B_2)$ являются подпространствами пространства $G \otimes H$. Покажем, что множества точек в них совпадают. Пусть $v \otimes u$ принадлежит подпространству $A = (A_1 \cap A_2) \otimes (B_1 \cap B_2)$. Это означает, что v принадлежит $A_1 \cap A_2$, u принадлежит $B_1 \cap B_2$. Отсюда, $v \in A_1$, $u \in B_1$ и $v \in A_2$, $u \in B_2$. Следовательно, $v \otimes u$ принадлежит $(A_1 \otimes B_1)$ и $(A_2 \otimes B_2)$ одновременно, то есть $v \otimes u$ принадлежит $(A_1 \otimes B_1) \cap (A_2 \otimes B_2)$. Аналогично $v \otimes u$ принадлежит $(A_2 \otimes B_1) \cap (A_1 \otimes B_2)$. Точно также проводится доказательство в обратном порядке. Вторая формула есть очевидное следствие из предшествующего соотношения для объединения. Теорема доказана.□

Замечание

Следует отметить, что, аналогично непрерывному случаю, если $A_1$, $A_1$ и $A_3$ есть подпространства пространства G, $B_1$, $B_1$ и $B_2$ есть подпространства пространства H, то

$$(A_1 \otimes B_1) \cap (A_2 \otimes B_2) \otimes (A_3 \otimes B_3) = (A_1 \cap A_2 \cap A_3) \otimes (B_1 \cap B_2 \cap B_3),$$

$$(A_1 \otimes B_1) \cup (A_2 \otimes B_2) \neq (A_1 \cup A_2) \otimes (B_1 \cup B_2).$$

В общем случае получаются следующие формулы.

Следствие

Пусть $A_1$, $A_2$,...$A_n$ есть подпространства пространства G, $B_1$, $B_2$,...$B_m$ есть подпространства пространства H. Тогда в пространстве $G \otimes H$ выполняется соотношения:

$$(A_1 \cap A_2 ... \cap A_n) \otimes (B_1 \cap B_2 ... \cap B_m) = \bigcap_{k=1}^{n} \bigcap_{K \le S}^{m} (A_k \otimes B_s)$$

$$(A_1 \cup A_2 ... \cup A_n) \otimes (B_1 \cup B_2 ... \cup B_m) = \bigcup_{k=1}^{n} \bigcup_{s=1}^{m} (A_k \otimes B_s)$$

Формула для пересечения может быть записана и в других формах. Важным результатом полученных формул является пересечение прямого произведения шаров.

Теорема 71

*Пусть $v_k \in G$, $v_s \in G$, $u_p \in H$, $u_r \in H$, тогда*
*$U(v_k \otimes u_p) \cap U(v_s \otimes u_r) = (U(v_k) \cap U(v_s)) \otimes (U(u_p) \cap U(u_r)) = U(v_k v_s) \otimes U(u_p u_r) = U(v_k \otimes u_r) \cap U(v_s \otimes u_p)$.*

Доказательство

В соответствии с предшествующей теоремой $(U(v_k) \cap U(v_s)) \otimes (U(u_p) \cap U(u_r)) = (U(v_k) \otimes U(u_p)) \cap (U(v_s) \otimes U(u_r)) = U(v_k \otimes u_p) \cap U(v_s \otimes u_r)$. Теорема доказана.□

Геометрически прямое произведение двух пространств $G \otimes H$ есть усложнение структуры пространства G (или H) путем замены каждой его точки на другое пространство H, и установления дополнительных связей между точками, принадлежащими смежным копиям H (Рис. 92). Пусть G состоит из точек ($v_1$, $v_2$, $v_3$,.... $v_m$), H состоит из точек ($u_1$, $u_2$, $u_3$,.... $u_n$). Каждую точку $v_k \in G$ заменяем на $H_k = H$ и получаем набор одинаковых пространств ($H_1$, $H_2$, $H_3$,.... $H_m$). Назовем $H_k$, и $H_p$, смежными, если точки $v_k$, и $v_p$, смежны. Устанавливаем связи между точками, принадлежащими различным $H_k$, следующим образом: если $H_k$, и $H_p$, смежны, то точка $u_s \in H_k$ соединяется со всеми точками шара $U(u_s) \subseteq H_p$. На Рис. 92 показан процесс построения прямого произведения окружности S, состоящей из 4-х точек ($v_1$,$v_2$,$v_3$,$v_4$) и связки K(3) из трех точек ($u_1$,$u_2$,$u_3$). Напомним, что если $A_1$ и $A_2$ есть подпространства пространства G, не имеющие общих точек, и каждая точка из $A_1$ является смежной с каждой точкой из $A_2$, то подпространство $A_1 \cup A_2$ называется прямой суммой подпространств $A_1$ и $A_2$ и обозначается $A_1 \oplus A_2$. Рассмотрим, как ведет себя $A_1 \oplus A_2$ в прямом произведении.

Теорема 72

*Пусть подпространство $A_1 \oplus A_2$ пространства G есть прямая сумма $A_1$ и $A_2$, и пусть B есть подпространство пространства H, являющееся связкой. Тогда $(A_1 \oplus A_2) \otimes B = (A_1 \otimes B) \oplus (A_2 \otimes B)$.*

Доказательство

Согласно одной из предшествующих теорем $(A_1 \oplus A_2) \otimes B = (A_1 \otimes B) \cup (A_2 \otimes B)$. Пространство $(A_1 \oplus A_2) \otimes B$ состоит из точек вида $v \otimes u$, где v принадлежит либо $A_1$ либо $A_2$, u принадлежит B. Такое же множество точек образует подпространство $(A_1 \otimes B) \cup (A_2 \otimes B)$. Рассмотрим, как соотносятся связи в этих подпространствах. Пусть точка $v \otimes u$ принадлежит $(A_1 \oplus A_2) \otimes B$, где v принадлежит $A_1$. Тогда точки, смежные этой точке в подпространстве $(A_1 \oplus A_2) \otimes B$, определяются объединением подпространств $v \otimes (O(u) \cap B)$, $((O(v) \cap A_1) \oplus A_2) \otimes u$, $((O(v) \cap A_1) \oplus A_2) \otimes (O(u) \cap B)$. Следовательно, точка $v \otimes u$ смежна со всеми точками, принадлежащими подпространству $(A_2 \otimes u) \cup (A_2 \otimes (O(u) \cap B))$. Так как $O(u) \cap B = B - u$, то $(A_2 \otimes u) \cup (A_2 \otimes (O(u) \cap B)) = (A_2 \otimes u) \cup (A_2 \otimes (B - u)) = (A_2 \otimes B)$. Из этого следует, что любая точка из $(A_1 \otimes B)$ смежна со всеми точками из $(A_2 \otimes B)$. Теорема

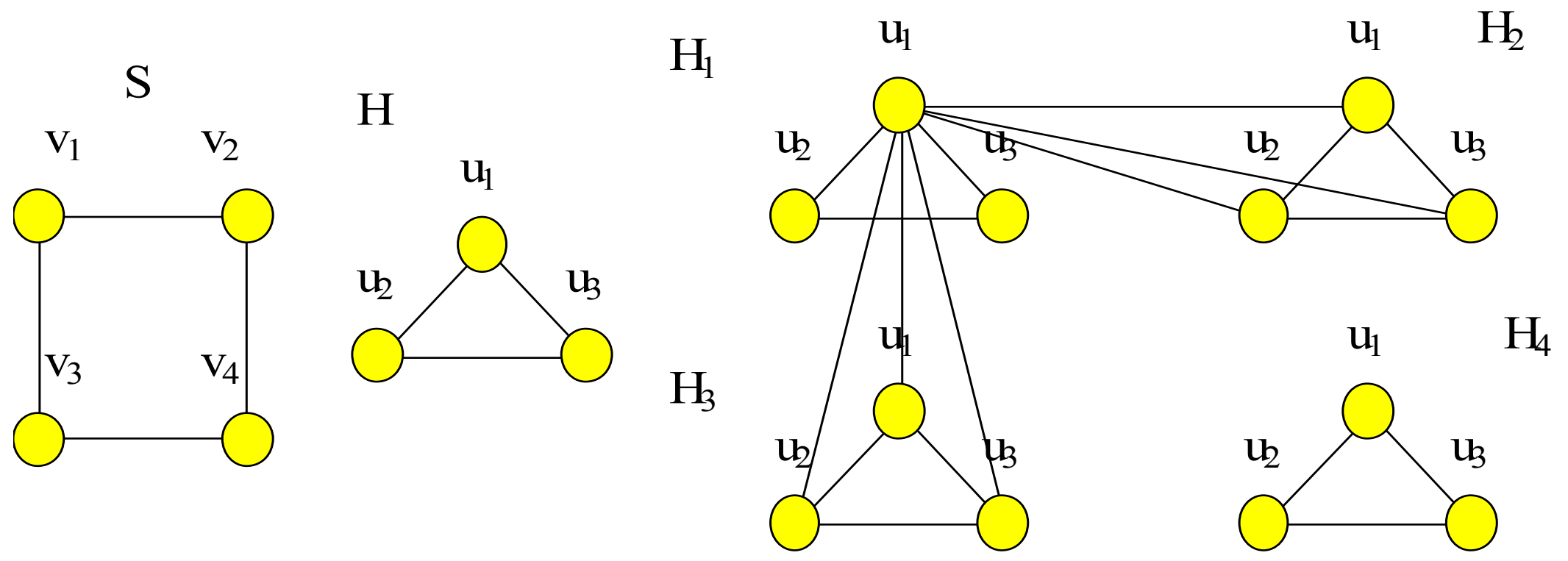

Рис. 92 Построение прямого произведения $S \otimes H$, где $H = K(3)$.Каждая точка окружности S заменяется на K(3). Точка $u_1 \in H_1$ соединяется с точками $u_1$, $u_2$ и $u_3$ из $H_2$ и $H_3$.

доказана.□

Следует еще раз отметить, что для того, чтобы на подмножество точек образовывало подпространство, необходимо, чтобы все связи между этими точками входили в это подпространство. При различного рода операциях это требование может нарушаться. Поэтому важно проверять выполнение этого условия в каждом отдельном случае. Сведем полученные соотношения в таблицу для удобства их использования вдальнейшем.

| |
|---|
| $U(v_k \otimes u_p) = U(v_k) \otimes U(u_p)$ |
| $U(v_k \otimes u_p ... \otimes w_s) = U(v_k) \otimes U(u_p) ... \otimes U(w_s)$. |
| $A \subseteq G$, $B \subseteq H$, $[A \otimes B] \subseteq G \otimes H \Rightarrow [A \otimes B] \approx A \otimes B$. |
| $A_1 \subseteq G$, $A_2 \subseteq G$, $B \subseteq H$, $\Rightarrow (A_1 \cap A_2) \otimes B = (A_1 \otimes B) \cap (A_2 \otimes B)$ в $G \otimes H$. |
| $A_1 \subseteq G$, $A_2 \subseteq G$, $B \subseteq H$, $\Rightarrow (A_1 \cup A_2) \otimes B = (A_1 \otimes B) \cup (A_2 \otimes B)$ в $G \otimes H$. |
| $A_1 \subseteq G$, $A_2 \subseteq G$, $B_1 \subseteq H$, $B_2 \subseteq H$, $\Rightarrow$ $(A_1 \cap A_2) \otimes (B_1 \cap B_2) = (A_1 \otimes B_1) \cap (A_2 \otimes B_2) = (A_2 \otimes B_1) \cap (A_1 \otimes B_2)$ в $G \otimes H$. |
| $v_k \in G$, $v_s \in G$, $u_p \in H$, $u_r \in H$, $\Rightarrow$ $U(v_k \otimes u_p) \cap U(v_s \otimes u_r) = (U(v_k) \cap U(v_s)) \otimes (U(u_p) \cap U(u_r)) = U(v_k v_s) \otimes U(u_p u_r) =$ $U(v_k \otimes u_r) \cap U(v_s \otimes u_p)$ в $G \otimes H$. |
| Пусть подпространство $A_1 \oplus A_2$ пространства G есть прямая сумма $A_1$ и $A_2$, и пусть B есть подпространство пространства H, являющееся связкой. Тогда $(A_1 \oplus A_2) \otimes B = (A_1 \otimes B) \oplus (A_2 \otimes B)$. |

Следующая теорема очевидна и не нуждаются в доказательстве.

Теорема 73

*Объем $|G \otimes H|$ прямого произведения пространств G и H равен произведению объемов $|G|$ и $|H|$, $|G \otimes H| = |G| * |H|$.*

## ЛОКАЛЬНЫЕ СВОЙСТВА ПРЯМОГО ПРОИЗВЕДЕНИЯ МОЛЕКУЛЯРНЫХ ПРОСТРАНСТВ

Рассмотрим более подробно структуру окрестности точки прямого произведения.

Теорема 74

*Пусть имеется два пространства: G с набором точек $(v_1, v_2, v_3, \ldots v_m)$ и набором связей W, и H с набором точек $(u_1, u_2, u_3, \ldots u_n)$ и набором связей S. Тогда в прямом произведении этих пространств $G \otimes H$ окрестность точки $v_k \otimes u_p$, определяется следующими условиями:*

*$v_k \otimes u_p$ смежна с $v_k \otimes u_s$, если $u_s$ принадлежит $O(u_p)$.*

*$v_k \otimes u_p$ смежна с $v_s \otimes u_p$, если $v_s$ принадлежит $O(v_k)$.*

*$v_k \otimes u_p$ смежна с $v_r \otimes u_s$, если $k \neq r$, $p \neq s$, $v_r$ принадлежит $O(v_k)$, $u_s$ принадлежит $O(u_p)$ одновременно.*

*Кроме того*

$$U(v_k \otimes u_p) = (v_k \otimes u_p) \oplus [(v_k \otimes O(u_p)) \oplus (O(v_k) \otimes u_p)) \cup (O(v_k) \otimes O(u_p)],$$

$$O(v_k \otimes u_p) = (v_k \otimes O(u_p) \oplus O(v_k) \otimes u_p)) \cup (O(v_k) \otimes O(u_p)).$$

Доказательство

Все свойства следуют из следующих выражений (Рис. 93).

$U(v_k \otimes u_p) = U(v_k) \otimes U(u_p) = (v_k \oplus (O(v_k)) \otimes (u_p \oplus O(u_p)) =$

$(v_k \otimes u_p) \cup (v_k \otimes O(u_p)) \cup (O(v_k) \otimes u_p)) \cup (O(v_k) \otimes O(u_p)) =$

$(v_k \otimes u_p) \oplus [(v_k \otimes O(u_p)) \oplus (O(v_k) \otimes u_p)) \cup (O(v_k) \otimes O(u_p)],$

$O(v_k \otimes u_p)=U(v_k \otimes u_p)-(v_k \otimes u_p)=(v_k \otimes O(u_p) \oplus O(v_k) \otimes u_p)) \cup (O(v_k) \otimes O(u_p))$.
Теорема доказана.□

Теорема 75

*В прямом произведении $G\otimes H$ окаем $O(v\otimes u)$ любой точки $v\otimes u$, где $v\in G$, $u\in H$, гомотопен прямой сумме $O(v)\oplus O(u)$ окаемов $O(v)$ и $O(u)$, $O(v\otimes u)\sim O(v)\oplus O(u)$. Кроме того $O(v\otimes u)$ стягивается к $O(v)\oplus O(u)$ с помощью точечных отбрасываний всех точек из $(O(v)\otimes O(u))$, принадлежащих $O(v\otimes u)$.*

Доказательство

Рассмотрим $O(v_1\otimes u_1)=[(v_1\otimes O(u_1))\oplus(O(v_1)\otimes u_1)]\cup(O(v_1)\otimes O(u_1))$. Возьмем произвольную точку $v_2\otimes u_2$ из $O(v_1)\otimes O(u_1)$. Легко видеть, что ее окаем в $O(v_1\otimes u_1)$ является конусом с вершинами $v_1\otimes u_2$ и $v_2\otimes u_1$. Следовательно, любая точка из $O(v_1)\otimes O(u_1)\subseteq O(v_1\otimes u_1)$ имеет точечный окаем, являющийся конусом, и может быть отброшена. При таких отбрасываниях структура окаемов других точек, являющихся конусами, не меняется. Это означает, пространство $O(v_1\otimes u_1)$ гомотопно $(v_1\otimes O(u_1))\oplus(O(v_1)\otimes u_1)$. Кроме того, $O(v_1\otimes u_1)$ стягивается к $(v_1\otimes O(u_1))\oplus(O(v_1)\otimes u_1)\approx O(u_1)\oplus O(v_1)$. Теорема доказана.□

На Рис. 94 показано, что окаем прямого произведения двух одномерных шаров является двумерным шаром, в котором окаем точки $v\otimes u$ есть регулярная окружность $O(v\otimes u)$, приводимая точечными отбрасываниями точек к нормальной минимальной окружности $O(v)\oplus O(u)=S^1=S^0\oplus S^0$, являющейся прямой суммой двух 0-мерных сфер.

В классической топологии прямое произведение единичного отрезка и окружности есть цилиндр, который гомотопен окружности. В ТМП имеет место подобное свойство.

Как следствие этой теоремы рассмотрим прямое произведение двух конусов.

Теорема 76

*Прямое произведение $(v\oplus G)\otimes(u\oplus H)$ конусов $(v\oplus G)$ и $(u\oplus H)$ пространств $G$ и $H$ определяется выражением*
*$(v\oplus G)\otimes(u\oplus H)=(v\otimes u)\oplus[((v\otimes H)\oplus(G\otimes u))\cup(G\otimes H)]$, где*
*$O(v\otimes u)=((v\otimes H)\oplus(G\otimes u))\cup(G\otimes H)\sim H\oplus G$ (Рис. 93).*

Доказательство

Очевидно, что $(v\oplus G)=U(v)$, $G=O(v)$ и $(u\oplus H)=U(u)$, $H=O(u)$ в пространствах $(v\oplus G)$ и $(u\oplus H)$ соответственно. Дальнейшее есть следствие предшествующей теоремы. Теорема доказана.□

Следующая теорема показывает, что применение точечного отбрасывания (или приклеивания) точки в одном из сомножителей в прямом произведении молекулярных пространств, порождает точечное отбрасывание точек в целом слое прямого произведения пространств и переводит прямое произведение в пространство, ему гомотопное.

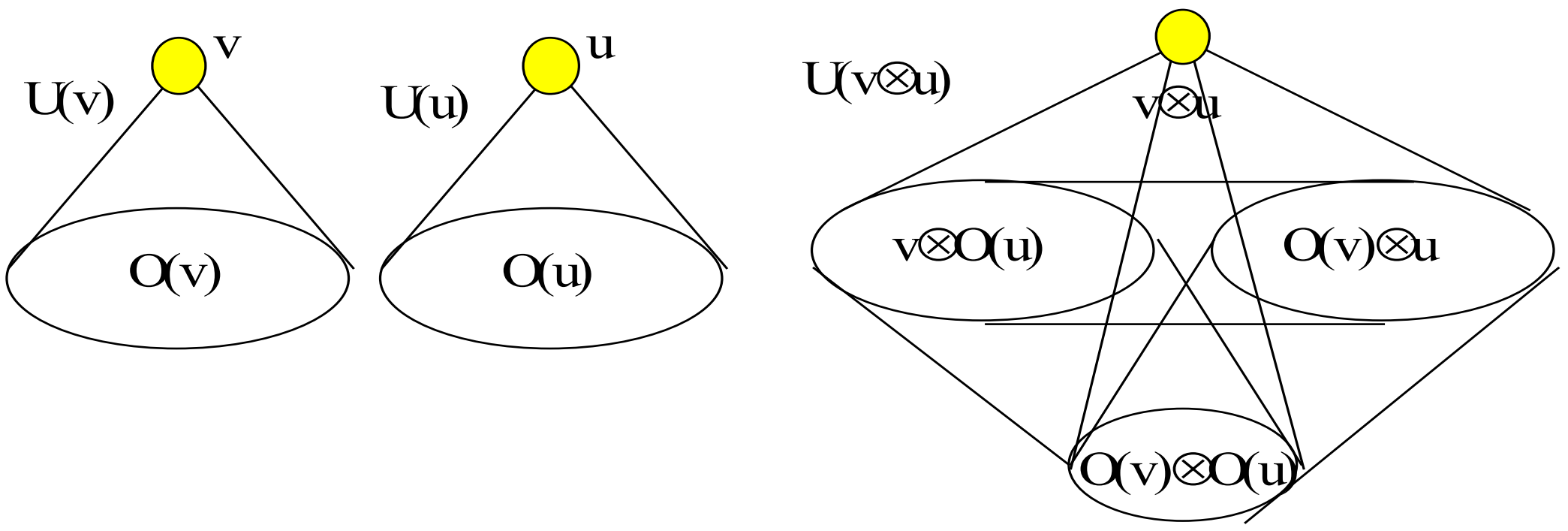

Рис. 93 Шар $U(v\otimes u)$ и окаем $O(v\otimes u)$ точки $v\otimes u$ в прямом произведении пространств $G\otimes H$ определяются выражениями:
$U(v\otimes u)=(v\otimes u)\oplus((v\otimes O(u))\oplus(O(v)\otimes u))\cup(O(v)\otimes O(u))$,
$O(v\otimes u)=((v\otimes O(u))\oplus(O(v)\otimes u))\cup(O(v)\otimes O(u))$.

Очевидно, что при этом достаточно использовать в доказательстве только отбрасывание точки.

Теорема 77

*Пусть u есть точка, $K(n)=(u_1,u_2,...u_n)$ есть связка из n точек, G есть произвольное пространство. Тогда $u\otimes G$ изоморфно G, $K(n)\otimes G$ гомотопно G, $u\otimes G\approx G$, $K(n)\otimes G\sim G$.*

Доказательство

Изоморфизм $u\otimes G$ и G очевиден. Пусть G состоит из точек $(v_1,v_2,...v_n)$.
Рассмотрим любую точку $u_k\otimes v_s$, $s\neq 1$, из $K(n)\otimes G$. Ее окаем имеет вид
$O(u_k\otimes v_s)=$
$((u_k\otimes O(v_s))\oplus(O(u_k)\otimes v_s))\cup(O(u_k)\otimes O(v_s))=$
$(u_1\otimes v_s)\oplus)[(u_k\otimes O(v_s))\oplus((O(u_k)-u_1)\otimes v_s)]\cup(O(u_k)\otimes O(v_s))]$ и является всегда конусом с вершиной $(u_1\otimes v_s)$. Следовательно эта точка может быть отброшена и при отбрасывании этой точки окаемы остальных точек остаются конусами. После отбрасывания всех таких точек пространство есть $u_1\otimes G$, изоморфное G. Теорема доказана.□

Эта теорема обобщается на точечный окаем и произвольное точечное пространство.

Теорема 78

*Пусть $v\in G$ и ее окаем $O(v)$ есть точечное подпространство. Тогда $G\otimes H$ гомотопно $(G-v)\otimes H$, $G\otimes H\sim(G-v)\otimes H$.*

Доказательство

Пусть для точки v ее окаем O(v) является точечным (Рис. 95). Тогда для любой точки $v\otimes u_1\in G\otimes H$ окаем $O(v\otimes u_1)$ этой точки определяется выражением $O(v\otimes u_1)=[(v\otimes O(u_1))\oplus(O(v)\otimes u_1)]\cup(O(v)\otimes O(u_1))$. Как показано выше, это подпространство гомотопно $O(u_1)\oplus O(v)$. Поскольку O(v) точечно, то $O(u_1)\oplus O(v)$ также точечно. Следовательно, $O(v\otimes u_1)$ гомотопно точечному пространству, то есть точечно. Это означает, что точка $v\otimes u_1$ может быть отброшена. Пусть точка $u_2$ смежна с $u_1$. Окаем этой точки определяется выражением $O(v\otimes u_2)=[(v\otimes O(u_2))\oplus(O(v)\otimes u_2)]\cup(O(v)\otimes O(u_2))-v\otimes u_1=[((v\otimes O(u_2))-v\otimes u_1)\oplus(O(v)\otimes u_2)]\cup(O(v)\otimes O(u_2))$. По-прежнему любую точку $v_k\otimes u_p$ из $O(v)\otimes O(u_2)$ можно отбросить, как имеющую точечный окаем, являющийся конусом с вершиной $v_k\otimes u_2$ из $O(v)\otimes u_2$. Пространство $((v\otimes O(u_2))-v\otimes u_1)\oplus(O(v)\otimes u_2)$ является точечным как прямая сумма некоторого пространства, в данном случае $((v\otimes O(u_2))-v\otimes u_1)$, и точечного пространства $O(v)\otimes u_2$. Следовательно, точка $v\otimes u_2$ имеет точечный окаем и может быть отброшена. То же самое рассуждение применимо ко всем точкам типа $v\otimes u_k$, где $u_k\in H$. Все эти точки могут быть отброшены. Оставшееся пространство есть прямое произведение $(G-v)\otimes H$. Теорема доказана.□

Следует еще раз отметить, что отбрасывание точки в G равносильно отбрасыванию целого слоя, подпространства H, в прямом произведении $G\otimes H$.

В классической топологии прямое произведение единичного отрезка и окружности есть цилиндр, который гомотопен окружности. В ТМП имеет место подобное свойство, которое доказывается в следующей теореме.

Теорема 79

*Если G есть точечное пространство, то $G\otimes H$ гомотопно H, $G\otimes H\sim H$.*

Доказательство

Так как пространство G точечно, то оно имеет по крайней мере две точки с точечными окаемами, если $|G|>1$. Пусть для точки v ее окаем O(v) является точечным. Тогда в соответствии с предшествующей теоремой все эти точки вида типа $v\otimes u_k$, где $u_k\in H$, могут быть отброшены. Оставшееся пространство есть прямое произведение $(G-v)\otimes H$, где G-v также является точечным и, следовательно, имеет точку $v_1$ с точечным окаемом. Тогда, согласно предшествующему, $(G-v)\otimes H$ гомотопна $(G-v-v_1)\otimes H$. Так как при отбрасывании точек G приводится к одноточечному пространству v, G~v, то G⊗H приводится к v⊗H≈H. Таким образом G⊗H~v⊗H≈H. Теорема доказана.□

Следствие

Если G и H есть точечные пространства, то G⊗H также является точечным.

Теперь можно более подробно определить структуру шаров и окаемов различных точек пространства G⊗H. Как и в классической топологии, в ТМП именно локальные свойства определяют такие характеристики пространства, как размерность, связность и некоторые другие. Выясним,

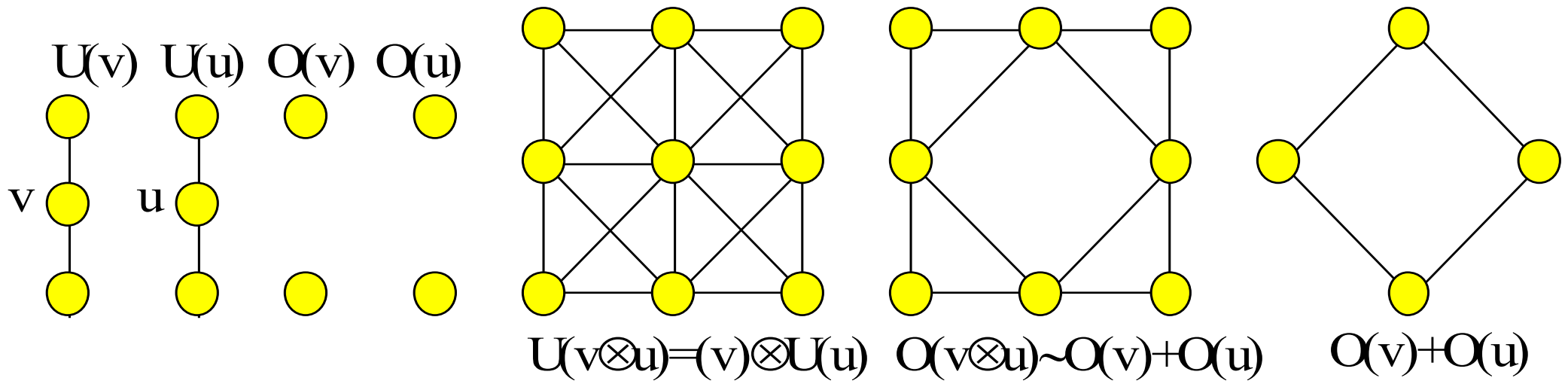

Рис. 94 Прямое произведение двух одномерных нормальных шаров является двумерным шаром, окаем O(v⊗u) есть регулярная окружность, гомотопная нормальной окружности O(v⊗u)~O(v)⊕O(u).

какова структура окаема точки общего окаема двух точек в прямом произведении пространств.

Теорема 80

*Для любых точек $v_1 \in G$, $v_2 \in G$, $u_1 \in H$ и $u_2 \in H$ общий окаем $O((v_1 \otimes u_1)(v_2 \otimes u_2))$ в $G \otimes H$ определяется следующими соотношениями:*

- *Если $v_1=v_2$ и $u_1$ и $u_2$ смежны, то*
$O((v_1 \otimes u_1)(v_1 \otimes u_2))=(v_1 \otimes O(u_1u_2)) \cup (O(v_1) \otimes u_1 \oplus u_2) \cup (O(v_1) \otimes O(u_1u_2))=$
$(v_1 \otimes O(u_1u_2)) \oplus (O(v_1) \otimes (u_1 \oplus u_2)) \cup (O(v_1) \otimes O(u_1u_2))=$
$(U(v_1) \otimes O(u_1u_2)) \cup (O(v_1) \otimes (u_1 \otimes u_2)) \sim O(v_1) \oplus O(u_1u_2)$.
- *Если $v_1=v_2$, $u_1$ и $u_2$ несмежны и $O(u_1u_2)$ непусто, то*
$O((v_1 \otimes u_1)(v_1 \otimes u_2))=(v_1 \otimes O(u_1u_2)) \cup (O(v_1) \otimes O(u_1u_2))=$
$U(v_1) \otimes O(u_1u_2) \sim O(u_1u_2))$.
- *Если $v_1$ и $v_2$ смежны, $u_1$ и $u_2$ смежны. то*
$O((v_1 \otimes u_1)(v_2 \otimes u_2))=(v_1 \otimes u_2) \oplus (v_2 \otimes u_1) \oplus [((v_1 \oplus v_2) \otimes O(u_1u_2)) \oplus (O(v_1v_2) \otimes$
$(u_1 \oplus u_2)) \cup O(v_1v_2) \otimes O(u_1u_2)] \sim (v_1 \otimes u_2) \oplus (v_2 \otimes u_1) \oplus [((v_1 \oplus v_2) \otimes O(u_1u_2)) \oplus$

*$(O(v_1v_2)\otimes(u_1\oplus u_2))]$. Выражение в квадратных скобках гомотопно $(v_1\oplus v_2)\otimes O(u_1u_2))\oplus(O(v_1v_2)\otimes(u_1\oplus u_2))\sim O(u_1u_2)\oplus(O(v_1v_2)$.*

*• Если $v_1$ и $v_2$ смежны, $u_1$ и $u_2$ несмежны и $O(u_1u_2)$ непусто, то $O((v_1\otimes u_1)(v_2\otimes u_2))=(v_1\oplus v_2)\otimes O(u_1u_2)\cup O(v_1v_2)\otimes O(u_1u_2)=$ $(v_1\oplus v_2\oplus O(v_1v_2))\otimes O(u_1u_2)\sim O(u_1u_2)$.*

*• Если $v_1$ и $v_2$ несмежны, $u_1$ и $u_2$ несмежны, $O(v_1v_2)$ непусто, $O(u_1u_2)$ непусто, то $O((v_1\otimes u_1)(v_2\otimes u_2))=O(v_1v_2)\otimes O(u_1u_2)$.*

*Доказательство*

Рассмотрим, вначале чисто формально, для произвольных точек в соответствии с предшествующей теоремой общий окаем двух точек.
$O((v_1\otimes u_1)(v_2\otimes u_2))=O(v_1\otimes u_1)\cap O(v_2\otimes u_2)=U(v_1\otimes u_1)\cap U(v_2\otimes u_2))-(v_1\otimes u_1)-(v_2\otimes u_2)=U(v_1v_2)\otimes U(u_1u_2)-(v_1\otimes u_1)-(v_2\otimes u_2)$.

Если $v_1=v_2$ и $u_1$ и $u_2$ смежны, то $U(v_1v_2)=U(v_1)=v_1\oplus O(v_1)$, $U(u_1u_2)=(u_1\oplus u_2\oplus O(u_1u_2))$. Учитывая это получаем
$O((v_1\otimes u_1)(v_1\otimes u_2))=(v_1\otimes O(u_1u_2))\cup(O(v_1)\otimes u_1\oplus u_2)\cup(O(v_1)\otimes O(u_1u_2))=$ $(v_1\otimes O(u_1u_2))\oplus(O(v_1)\otimes(u_1\oplus u_2))\cup(O(v_1)\otimes O(u_1u_2))=$ $(U(v_1)\otimes O(u_1u_2))\cup(O(v_1)\otimes(u_1\otimes u_2))\sim O(v_1)\oplus O(u_1u_2)$.

Если $v_1=v_2$, $u_1$ и $u_2$ несмежны и $O(u_1u_2)$ непусто, то $U(v_1v_2)=U(v_1)=v_1\oplus O(v_1)$, $U(u_1u_2)=O(u_1u_2)$. Учитывая это получаем $O((v_1\otimes u_1)(v_1\otimes u_2))=(v_1\otimes O(u_1u_2))\cup(O(v_1)\otimes O(u_1u_2))=U(v_1)\otimes O(u_1u_2)\sim O(u_1u_2))$.

Если $v_1$ и $v_2$ смежны, $u_1$ и $u_2$ смежны. то $U(v_1v_2)=(v_1\oplus v_2\oplus O(v_1v_2))$, $U(u_1u_2)=(u_1\oplus u_2\oplus O(u_1u_2))$. Учитывая это и приводя к удобной форме получаем
$O((v_1\otimes u_1)(v_2\otimes u_2))=(v_1\otimes u_2)\oplus(v_2\otimes u_1)\oplus[((v_1\oplus v_2)\otimes O(u_1u_2))\oplus(O(v_1v_2)\otimes(u_1\oplus u_2))\cup O(v_1v_2)\otimes O(u_1u_2)]\sim(v_1\otimes u_2)\oplus(v_2\otimes u_1)$ $\oplus[((v_1\oplus v_2)\otimes O(u_1u_2))\oplus(O(v_1v_2)\otimes(u_1\oplus u_2))]$. Выражение в квадратных скобках гомотопно
$((v_1\oplus v_2)\otimes O(u_1u_2))\oplus(O(v_1v_2)\otimes(u_1\oplus u_2))\sim O(u_1u_2)\oplus(O(v_1v_2)$.

Если $v_1$ и $v_2$ смежны, $u_1$ и $u_2$ несмежны и $O(u_1u_2)$ непусто, то $U(v_1v_2)=(v_1\oplus v_2\oplus O(v_1v_2))$, $U(u_1u_2)=O(u_1u_2)$. Следовательно, $O((v_1\otimes u_1)(v_2\otimes u_2))=(v_1\oplus v_2)\otimes O(u_1u_2)\cup O(v_1v_2)\otimes O(u_1u_2)=$ $(v_1\oplus v_2\oplus O(v_1v_2))\otimes O(u_1u_2)\sim O(u_1u_2)$.

Если $v_1$ и $v_2$ несмежны, $u_1$ и $u_2$ несмежны, $O(v_1v_2)$ непусто, $O(u_1u_2)$ непусто, то $U(v_1v_2)=O(v_1v_2)$, $U(u_1u_2)=O(u_1u_2)$. Следовательно, $O((v_1\otimes u_1)(v_2\otimes u_2))=O(v_1v_2)\otimes O(u_1u_2)$.

Теорема доказана.□

Теперь рассмотрим, как влияет точечное отбрасывание связи в одном из сомножителей на прямое произведение пространств. Мы хотим показать,

что так же, как и в случае отбрасывания точки с точечным окаемом, прямое произведение остается в том же гомотопном классе.

Теорема 81

*Пусть точки $u_1$ и $u_2$ в пространстве H смежны и их общий окаем $O(u_1u_2)$ есть точечное подпространство. Пусть пространство $H-(u_1u_2)$ получено из H путем отбрасывания точечной связи $(u_1u_2)$. Тогда $G\otimes(H-(u_1u_2))$ гомотопно $G\otimes H$, $G\otimes(H-(u_1u_2))\sim G\otimes H$.*

Доказательство

Отбрасывание связи $(u_1u_2)$ в H равносильно отбрасыванию серии связей между точками $v_k\otimes u_1$ и $v_k\otimes u_2$, где $v_k$ есть любая точка из G, а также между точками $v_k\otimes u_1$ и $v_p\otimes u_2$, где $v_k$ и $v_p$ есть любая пара смежных точек из G. Рассмотрим общие окаемы этих точек.

$O((v_k\otimes u_1)(v_p\otimes u_2))=(v_k\otimes u_2\oplus v_p\otimes u_1)\oplus[(v_k\otimes O(u_1u_2)\oplus O(v_kv_p)\otimes u_1)\cup(v_p\otimes O(u_1u_2)\oplus O(v_kv_p)\otimes u_2)]\cup O(v_kv_p)\otimes O(u_1u_2)$.

Выделим подпространство

$[(v_k\otimes O(u_1u_2)\oplus O(v_kv_p)\otimes u_1)\cup(v_p\otimes O(u_1u_2)\oplus O(v_kv_p)\otimes u_2)]\cup O(v_kv_p)\otimes O(u_1u_2)$.

Точки, принадлежащие подпространству $O(v_kv_p)\otimes O(u_1u_2)$, имеют точечные окаемы, являющиеся конусами с вершинами в подпространстве $(v_k\otimes O(u_1u_2)\oplus O(v_kv_p)\otimes u_1)\cup(v_p\otimes O(u_1u_2)\oplus O(v_kv_p)\otimes u_2)$. Следовательно, они могут быть отброшены. В свою очередь подпространство $(v_k\otimes O(u_1u_2)\oplus O(v_kv_p)\otimes u_1)\cup(v_p\otimes O(u_1u_2)\oplus O(v_kv_p)\otimes u_2)=(v_k\oplus v_p)\otimes O(u_1u_2)\cup(O(v_kv_p)\otimes(u_1\oplus u_2)=(v_k\oplus v_p)\otimes O(u_1u_2)\oplus(O(v_kv_p)\otimes(u_1\oplus u_2)$ является прямой суммой точечного подпространства $(v_k\oplus v_p)\otimes O(u_1u_2)$ и подпространства $(O(v_kv_p)\otimes(u_1\oplus u_2)$, и, следовательно, есть точечное пространство вне зависимости от отсутствия или наличия связи между точками $v_k\otimes u_2$ и $v_p\otimes u_1$. Таким образом, поскольку $O(u_1u_2)$ есть точечное подпространство, все между точками $v_k\otimes u_1$ и $v_p\otimes u_2$, где $v_k$ и $v_p$ есть любая пара смежных точек из G, могут быть отброшены. Рассмотрим, теперь, после этих отбрасывания связей, общий окаем точек $v_k\otimes u_1$ и $v_k\otimes u_2$, где $v_k$ есть любая точка из G. Он имеет вид

$O((v_k\otimes u_1)(v_k\otimes u_2))=O(v_k)\otimes O(u_1u_2)\oplus O(v_k)\otimes(u_1\oplus u_2)$.

Поскольку $O(u_1u_2)$ является точечным подпространством, $O((v_k\otimes u_1)(v_k\otimes u_2))$ также будет точечным. Из вида общих окаемов следует, что отбрасывание части связей не влияет на структуру общих окаемов оставшихся точек. Таким образом переход от $G\otimes H$ к $(G-(v_1v_2))\otimes H$ осуществляется путем точечного отбрасывания связей. Теорема доказана.□

Как следствие двух предшествующих теорем сформулируем теорему, определяющую как влияют гомотопные преобразования на прямое произведение пространств.

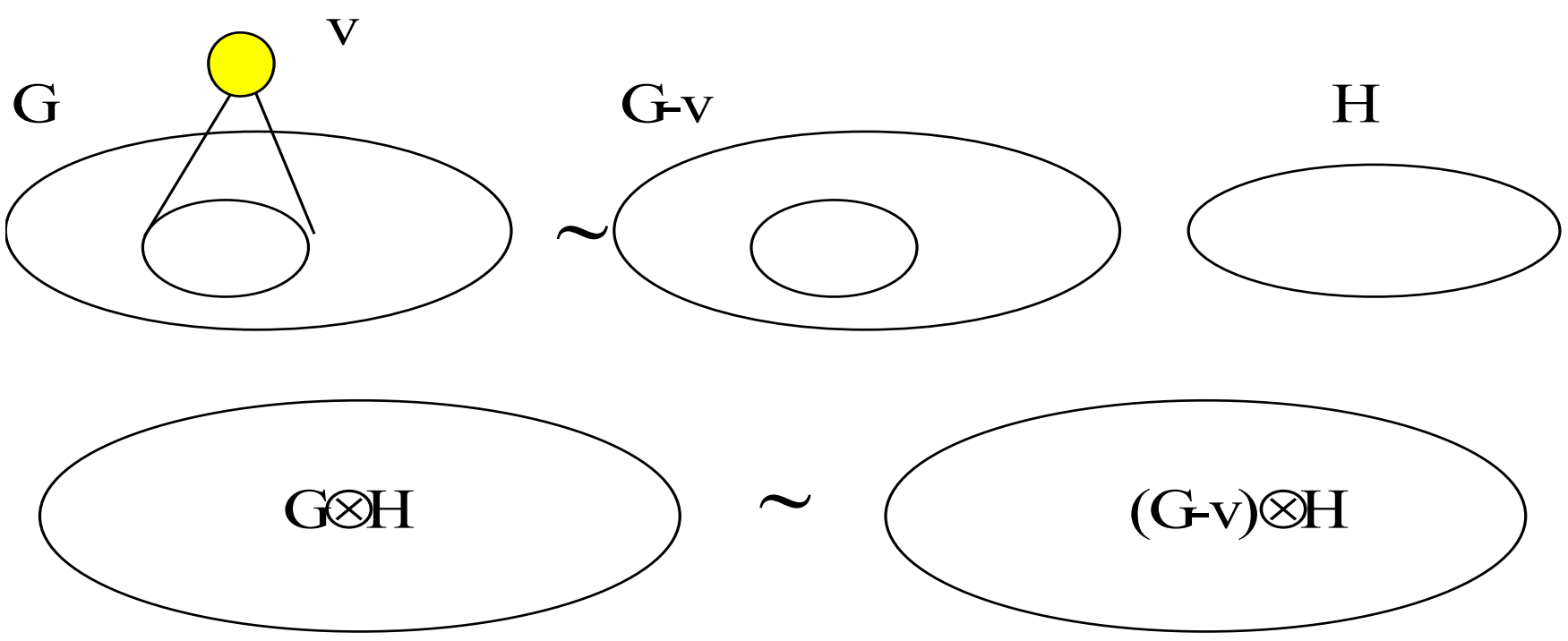

Рис. 95 Если G и G-v гомотопны, то есть O(v) является точечным пространством, то (G-v)⊗H гомотопно G⊗H.

Теорема 82

*Если G~H, то для любого пространства P прямое произведение G⊗P гомотопно H⊗P, G⊗P~H⊗P.*

Доказательство

Переход от G к H состоит из последовательных точечных отбрасываний и приклеиваний точек и связей. Как показано выше, при каждом таком преобразовании прямое произведение переходит в гомотопное. Следовательно, G⊗P~H⊗P (Рис. 96). Теорема доказана.□

Теорема 83

*Пусть v∈G, u∈H и их окаемы O(v) и O(u) есть неточечные подпространства. Тогда окаем точки v⊗u∈G⊗H является неточечным.*

Доказательство

Окаем точки v⊗u определяется выражением

O(v⊗u)=[(v⊗O(u))⊕(O(v)⊗u)]∪(O(v)⊗O(u))~(v⊗O(u))⊕(O(v)⊗u)≈O(u)⊕O(v). Прямая сумма двух неточечных пространств есть неточечное пространство. Теорема доказана.□

Рассмотрим, какова топологическая структура окаемов точек прямого произведения в слоях. Связь может быть неточечной в пространстве, но быть точечной в прямом произведении, где это пространство является слоем. Это представляет значительный интерес в прямом произведении нормальных пространств.

Теорема 84

*Пусть $G_1$ состоит из точек $x_1$, $G_2$ состоит из точек $x_2$. Пусть две смежные точки $v_1 \oplus v_2$ и $u_1 \oplus v_2$ принадлежат слою $G_1$. Тогда в прямом произведении $G_1 \otimes G_2$ :*

*если $O(v_1u_1)$ и $O(v_2)$ неточечные, то $O((v_1 \oplus v_2)(u_1 \oplus v_2))$ также неточечный,*

*если $O(v_1u_1)$ точечный и/или $O(v_2)$ точечный, то $O((v_1 \oplus v_2)(u_1 \oplus v_2))$ точечный.*

*Пусть две несмежные точки $v_1 \oplus v_2$ и $u_1 \oplus v_2$ принадлежат слою $G_1$. Тогда:*

*$O((v_1 \oplus v_2)(u_1 \oplus v_2))$ неточечный тогда и только тогда, когда $O(v_1u_1)$ неточечный.*

Доказательство

Согласно предшествующей теореме если $v_1 \oplus v_2$ и $u_1 \oplus v_2$ смежные, то $O((v_1 \oplus v_2)(u_1 \oplus v_2)) \sim O(v_1u_1)) \oplus (O(v_2)$. По свойству прямой суммы если оба слагаемых неточечные, то прямая сумма неточечная. Если хотя бы одно из слагаемых точечное, то прямая сумма также точечная. Если $v_1 \oplus v_2$ и $u_1 \oplus v_2$ несмежные, то $O((v_1 \oplus v_2)(u_1 \oplus v_2)) \sim O(v_1u_1) \otimes U(v_2)) \sim O(v_1u_1)$. Следовательно, $O((v_1 \oplus v_2)(u_1 \oplus v_2))$ и $O(v_1u_1)$ точечные или неточечные одновременно. Теорема доказана.□

Следствие для молекулярных пространств

Прямое и важное следствие из этих теорем для пространств с неточечными окаемами и, в частности, нормальных пространств состоит в том, в прямом произведении нормальных пространств $G_1 \otimes G_2$ слой $G_1$ как подпространство в $G_1 \otimes G_2$ обладает той же топологической структурой, что и $G_1$ само по себе. В слое можно ввести или отбросить связь тогда и только тогда, когда это можно сделать в $G_1$.

Рассмотрим, как меняются окаемы точек и само прямое произведение при отбрасываниях связей.

Теорема 85

*Пусть имеется два произвольных молекулярных пространства: пространство A с набором точек ($v_1$, $v_2$, $v_3$,.... $v_m$...) и набором связей W, и пространство B с набором точек ($u_1$, $u_2$, $u_3$,.... $u_n$...) и набором связей S. Тогда в их прямом произведении $A \otimes B$ отбрасывание одной связи из каждой пары связей вида $((v_i \otimes u_k)(v_p \otimes u_s))$ и $((v_i \otimes u_s)(v_p \otimes u_k))$, $i \neq p$, $k \neq s$, является точечным преобразованием. При этом окаем* *любой точки пространства всегда остается гомотопным исходному окаему этой точки.*

Доказательство

Выберем произвольную точку $v_1 \otimes u_1$ в $A \otimes B$. Ее окаем определяется выражением $O(v_1 \otimes u_1)=(v_1 \otimes O(u_1) \oplus O(v_1) \otimes u_1) \cup O(v_1) \otimes O(u_1)$. Рассмотрим общий окаем двух смежных точек $v_1 \otimes u_1$ и $v_2 \otimes u_2$ в $A \otimes B$ при условии, что $v_1 \neq v_2$ и $u_1 \neq u_2$.

$O((v_1 \otimes u_1)(v_2 \otimes u_2))=(v_1 \otimes u_2 \oplus v_2 \otimes u_1) \oplus [(v_1 \otimes O(u_1 u_2) \oplus O(v_1 v_2) \otimes u_1) \cup (v_2 \otimes O(u_1 u_2) \oplus O(v_1 v_2) \otimes u_2)] \cup O(v_1 v_2) \otimes O(u_1 u_2)=v_1 \otimes u_2 \oplus v_2 \otimes u_1 \oplus W$. Очевидно, что при этом точки $v_1 \otimes u_2$ и $v_2 \otimes u_1$ также будут смежны, и их общий окаем будет определяться аналогичным выражением $O((v_1 \otimes u_2)(v_2 \otimes u_1))= v_1 \otimes u_1 \oplus v_2 \otimes u_2 \oplus [(v_1 \otimes O(u_1 u_2) \oplus O(v_1 v_2) \otimes u_1) \cup (v_2 \otimes O(u_1 u_2) \oplus O(v_1 v_2) \otimes u_2)] \cup O(v_1 v_2) \otimes O(u_1 u_2)=v_1 \otimes u_1 \oplus v_2 \otimes u_2 \oplus W$. Из этих выражений видно, что вне зависимости от вида пространства W, оба общих окаема являютя точечным пространствами, а именно, конусами с двумя вершинами. Следовательно, одна из связей, или $((v_1 \otimes u_1)(v_2 \otimes u_2))$, или $((v_1 \otimes u_2)(v_2 \otimes u_1))$, может быть отброшена. Точно по той же причине может быть отброшена любая связь из каждой пары связей вида $((v_i \otimes u_k)(v_p \otimes u_s))$ и $((v_i \otimes u_s)(v_p \otimes u_k))$, $i \neq p$, $k \neq s$. Это завершает доказательство первой части теоремы.

Покажем, что отбрасывание связи в $A \otimes B$ соответствует точечному отбрасыванию точки в $O(v_1 \otimes u_1)$. Рассмотрим окаем точки $v_1 \otimes u_1$ как самостоятельное пространство $O(v_1 \otimes u_1)=(v_1 \otimes O(u_1) \oplus O(v_1) \otimes u_1) \cup O(v_1) \otimes O(u_1)$. Рассмотрим окаем смежной точки $v_2 \otimes u_2$ в пространстве $O(v_1 \otimes u_1)$ при условии, что $v_1 \neq v_2$ и $u_1 \neq u_2$. Очевидно, что $O(v_2 \otimes u_2)|O(v_1 \otimes u_1)= O((v_1 \otimes u_1)(v_2 \otimes u_2))=v_1 \otimes u_2 \oplus v_2 \otimes u_1 \oplus W$. При этом точки $v_1 \otimes u_2$ и $v_2 \otimes u_1$ также будут смежны, и, вне зависимости от вида пространства W, $O(v_2 \otimes u_2)|O(v_1 \otimes u_1)$ является точечным пространствам, а именно, конусом с двумя смежными вершинами $v_1 \otimes u_2$ и $v_2 \otimes u_1$. Следовательно, точка $v_2 \otimes u_2$ может быть отброшена из пространства $O(v_1 \otimes u_1)$ точечным преобразованием. Выберем ту же произвольную точку $v_1 \otimes u_1$ в $A \otimes B$. Ее окаем имеет вид $O(v_1 \otimes u_1)=(v_1 \otimes O(u_1) \oplus O(v_1) \otimes u_1) \cup O(v_1) \otimes O(u_1)$. Покажем, что отбрасывание связи в $A \otimes B$ соответствует точечному отбрасыванию связи в $O(v_1 \otimes u_1)$. Рассмотрим общий окаем смежных точек $v_1 \otimes u_2$ и $v_2 \otimes u_1$ в пространстве $O(v_1 \otimes u_1)$. Очевидно, что он определяется выражением $O((v_1 \otimes u_2)(v_2 \otimes u_1))|O(v_1 \otimes u_1)=v_2 \otimes u_2 \oplus W$. Вне зависимости от вида пространства W, $O((v_1 \otimes u_2)(v_2 \otimes u_1))|O(v_1 \otimes u_1)$ является точечным пространствам, а именно, конусом с вершиной $v_2 \otimes u_2$. Следовательно, связь точка $((v_1 \otimes u_2)(v_2 \otimes u_1))$ может быть отброшена из пространства $O(v_1\otimes u_1)$ точечным преобразованием. Отсюда следует, что при любом отбрасывании связи из условия теоремы окаем любой точки пространства остается гомотопным исходному. Теорема доказана. □

Эта теорема легко обобщается на некоторое подпространство пространства $G\otimes H$.

Теорема 86

*Пусть имеется два произвольных молекулярных пространства: пространство A с набором точек ($v_1, v_2, v_3,\ldots. v_m\ldots$) и набором связей W, и пространство B с набором точек ($u_1, u_2, u_3,\ldots. u_n\ldots$) и набором связей S. Пусть C является подпространством прямого произведения $A\otimes B$, для каждой пары смежных точек которого выполняется условие: если ($v_i\otimes u_k$) смежна с ($v_p\otimes u_s$) в $A\otimes B$, то $i\neq p$, $k\neq s$. Тогда точечным отбрасыванием связей в пространстве $A\otimes B$ окаем любой точки $v_p\otimes u_s$ из C может быть приведен к виду $O(v_p\otimes u_s)=v_p\otimes O(u_s)\oplus O(v_p)\otimes u_s$, где $O(v_p)$ есть окаем точки $v_p$ в A, из $O(u_s)$ есть окаем точки $u_s$ в B (Рис. 94).*

Доказательство

Выберем произвольную точку $v_1\otimes u_1$ в C. Ее окаем определяется выражением $O(v_1\otimes u_1)=(v_1\otimes O(u_1)\oplus O(v_1)\otimes u_1)\cup O(v_1)\otimes O(u_1)$. Этот окаем не содержит никаких других точек из C. В соответствии с предшествующей

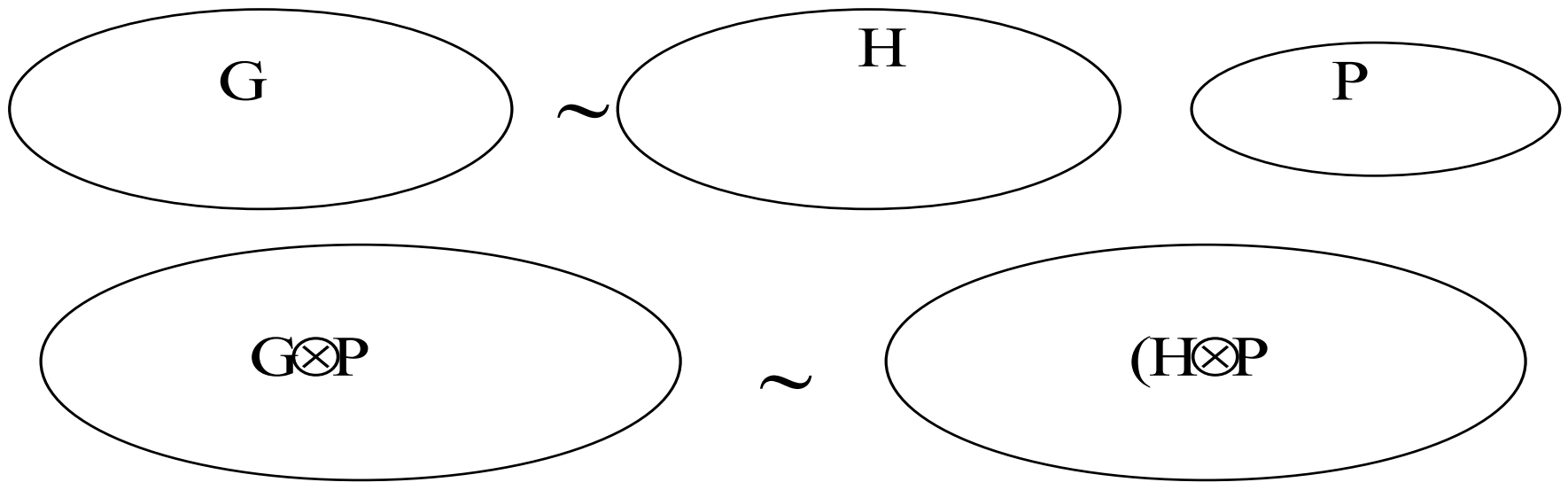

Рис. 96 Если G и H гомотопны, G~H, то при любом пространстве P $G\otimes P$ гомотопно $H\otimes P$, $G\otimes P\sim H\otimes P$.

теоремой отбрасываем все связи вида $((v_1\otimes u_1)(v_p\otimes u_s))$, где $v_p\otimes u_s$ есть точки, смежные $v_1\otimes u_1$ при условии, что $v_1\neq v_p$ и $u_1\neq u_s$. Эти преобразования являются точечными как в $A\otimes B$, так и в $O(v_1\otimes u_1)$. При этом $O(v_1\otimes u_1)=(v_1\otimes O(u_1)\oplus O(v_1)\otimes u_1)\cup O(v_1)\otimes O(u_1)$ переходит в $O(v_1\otimes u_1)=v_1\otimes O(u_1)\oplus O(v_1)\otimes u_1$. Очевидно, что эти преобразования с тем же результатом можно применить к любой другой точке из C. Проделывая подобные точечные отбрасывания связей для всех точек из подпространства C получаем, что окаем любой точки вида $((v_p\otimes u_s)$ из C примет вид $O(v_p \otimes u_s)=v_p \otimes O(u_s) \oplus O(v_p) \otimes u_s$. Теорема доказана. □

На Рис. 97 изображено прямое произведение $L^1{}_1 \otimes L^1{}_2$ двух одномерных пространств, в которых окаем каждой точки есть 0-мерная сфера. Из рисунка видно, что окаем каждой точки пространства C приводится к прямой сумме двух 0-мерных сфер $S^0 \oplus S^0 = S^1$, которая является одномерной минимальной окружностью.

В заключение этого раздела рассмотрим свойства прямого произведения нескольких молекулярных пространств. Легко видеть, что свойства прямого произведения двух пространств обобщаются на прямое произведение любого конечного числа пространств. Заново формулировать все из этих свойств мы не будем, докажем только наиболее важные из них, касающиеся нормальных пространств. Хотя результаты этой теоремы скорее технические, они важны тем, что позволяют доказать многие принципиальные теоремы о размерных пространствах.

## *СВОЙСТВА ПРЯМОГО ПРОИЗВЕДЕНИЯ НОРМАЛЬНЫХ МОЛЕКУЛЯРНЫХ ПРОСТРАНСТВ*

Все предшествующие теоремы были доказаны также и для того, чтобы выяснить структуру прямого произведения размерных пространств. Метрические свойства прямого произведения представляют большой интерес, поскольку именно они важны в компьютерных приложениях. Оказывается, например, что прямое произведение нормальных пространств, хотя само по себе и не является нормальным пространством, может быть стянуто к нормальному пространству только точечными отбрасываниями связей. Этот результат важен потому, что он в точности соответствует аналогичному результату в классической топологии. Так как в прямом произведении структура пространства определяется локально, достаточно изучить свойства прямого произведения двух шаров.

Теорема 87

*Пусть шары U(v) и U(u) не содержат подпространств, гомотопных n-мерному и m-мерному замкнутым нормальным пространствам соответственно. Тогда их прямое произведение $U(v \otimes u)$ не содержит никакого подпространства, гомотопного нормальному замкнутому (n+m+1)-мерному пространству.*

Доказательство

Предположим противное. Пусть шар $U(v \otimes u)=U(v) \otimes U(u)$ содержит подпространство, гомотопное нормальному замкнутому (n+m+1)-мерному пространству. Оно должно содержаться в $v \otimes O(u) \oplus O(v) \otimes u \approx O(u) \oplus O(v)$, поскольку точки из $O(v) \otimes O(u)$ имеют точечные окаемы и могут быть отброшены. По одному из свойств прямой суммы пространства O(v) и O(u) должны быть гомотопны нормальным замкнутым s-мерному и p-мерному пространствам, где s+p=m+n. Это означает, что или s≥m, или p≥n. По условию теоремы это невозможно. Теорема доказана. □

Докажем сейчас одну из основных теорем, определяющую размерные свойства прямого произведения двух нормальных пространств. Если U(v) и U(u) являются n-мерным и m-мерным нормальными шарами, то окаемы точек v и u есть (n-1)-мерное и (m-1)-мерное нормальные замкнутые пространства. Мы покажем, что U(v⊗u) является (n+m)-мерным шаром, то есть окаем точки (v⊗u) содержит нормальное (n+m-1)-мерное замкнутое подпространство, и не содержит никакого подпространства, гомотопного (n+m)-мерному замкнутому нормальному пространству. Аналогичное свойство имеется в классической топологии [40].

Теорема 88

*Пусть шары U(v) и U(u) являются n-мерным и m-мерным нормальными шарами соответственно. Тогда их прямое произведение U(v⊗u) является (n+m)-мерным шаром, то есть содержит (n+m-1)-мерное нормальное подпространство в своем окаеме и не содержит подпространства, гомотопного нормальному замкнутому (n+m)-мерному пространству.*

Доказательство

Так как шары U(v) и U(u) являются n-мерным и m-мерным нормальными шарами соответственно, то окаемы O(v) и O(u) являются (n-1)-мерным и (m-1)-мерным нормальными замкнутыми подпространствами. Как известно, O(v⊗u)=[(v⊗O(u))⊕(O(v)⊗u)]∪(O(v)⊗O(u)). Это подпространство содержит (v⊗O(u))⊕(O(v)⊗u)≈O(u))⊕O(v), которое является прямой суммой окаемов и, следовательно, есть (m+n-1)-мерное замкнутое нормальное пространство. В пространствах O(v) и O(u) шары точек являются (n-1) и (m-1)-мерными нормальными замкнутыми шарами, следовательно, каждая точка пространства O(v)⊗O(u)) содержит в своем окаеме (m+n-3)-мерное замкнутое нормальное пространство. Все точки из O(v)⊗O(u)) имеют точечные окаемы в O(v⊗u) и могут быть отброшены. Кроме того, согласно предшествующей теореме, U(v⊗u) не содержит никакого подпространства, гомотопного нормальному замкнутому (n+m)-мерному пространству. Теорема доказана. □

Теорема 89

*Пусть имеется два нормальных пространства: n-мерное пространство $G^n$ с набором точек ($v_1$, $v_2$, $v_3$,.... $v_m$...) и набором связей W, и m-мерное пространство $H^m$ с набором точек ($u_1$, $u_2$, $u_3$,.... $u_n$...) и набором связей S. Тогда их прямое произведение $G^n \otimes H^m$ обладает следующими свойствами:*

- *Шар $U(v\otimes u)$ каждой точки $v\otimes u$ является $(n+m)$-мерным шаром, то есть содержит $(n+m-1)$-мерное нормальное подпространство в своем окаеме и не содержит подпространства, гомотопного нормальному замкнутому $(n+m)$-мерному пространству.*
- *$G^m\otimes H^n$ гомотопно замкнутому нормальному пространству размерности $(m+n)$, и приводится к нему точечным отбрасыванием связей.*
- *Точечным отбрасыванием одной связи из каждой пары связей вида $((v_i\otimes u_k)(v_p\otimes u_s))$ и $((v_i\otimes u_s)(v_p\otimes u_k))$, $i\neq p$, $k\neq s$, пространство $G^n\otimes H^m$ приводится к нормальному $(n+m)$-мерному пространству.*

Доказательство

Первое свойство доказано в предшествующей теореме. Дальнейшее доказательство состоит из нескольких этапов.

Первый этап доказательства.

На первом этапе покажем, что окаем произвольной точки приводится точечными отбрасываниями связей к нормальному замкнутому пространству. Выберем произвольную точку $v_1\otimes u_1$ в $G^m\otimes H^n$. Ее окаем определяется выражением $O(v_1\otimes u_1)=(v_1\otimes O(u_1)\oplus O(v_1)\otimes u_1)\cup O(v_1)\otimes O(u_1)$. Рассмотрим общий окаем $O((v_1\otimes u_1)(v_p\otimes u_s))=v_1\otimes u_s\oplus v_p\otimes u_1\oplus[[(v_1\otimes O(u_1u_s)\oplus O(v_1v_p)\otimes u_1)\cup(v_p\otimes O(u_1u_s)\oplus O(v_1v_p)\otimes u_s)]\cup O(v_1v_p)\otimes O(u_1u_s)]$ двух смежных точек $v_1\otimes u_1$ и $v_p\otimes u_s$ в $G^m\otimes H^n$ при условии, что $v_1\neq v_p$ и $u_1\neq u_s$. Легко видеть, что этот окаем является точечным подпространством, и связь $((v_1\otimes u_1)(v_p\otimes u_s))$ может быть отброшена. Точно по той же причине могут быть отброшены все связи между $(v_1\otimes O(u_1)\oplus O(v_1)\otimes u_1)$ и $O(v_1)\otimes O(u_1)$. Следовательно, после отбрасывания связей окаем $O(v_1\otimes u_1)=v_1\otimes O(u_1)\oplus O(v_1)\otimes u_1$ точки $v_1\otimes u_1$

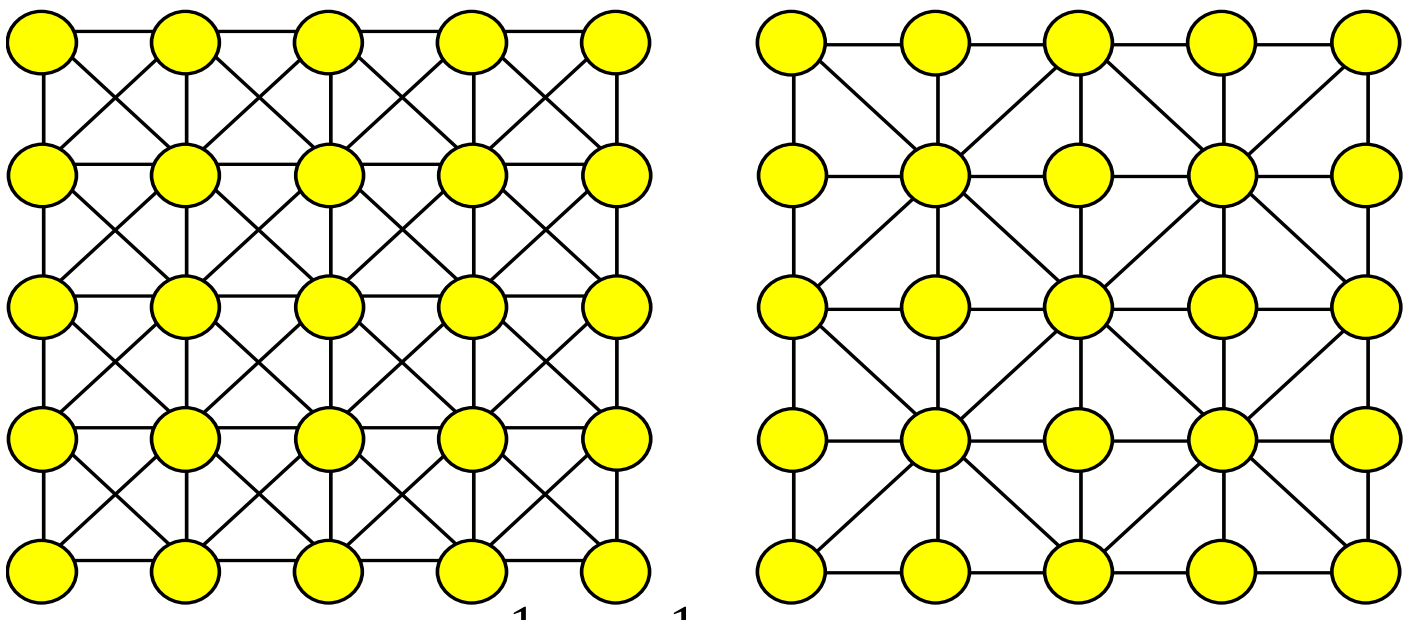

Рис. 97 Прямое произведение $L^1{}_1\otimes L^1{}_2$ двух одномерных прямых $L^1{}_1$ и $L^1{}_2$ содержит подпространство C. Окаемы точек из C точечными отбрасываниями связей приводятся к виду $v_p\otimes O(u_s)\oplus O(v_p)\otimes u_s$.

есть прямое произведение двух нормальных замкнутых пространств $v_1\otimes O(u_1)$ и $O(v_1)\otimes u_1$, то есть является нормальным замкнутым $(n+m-1)$- мерным пространством, поскольку $O(v_1)$ и $O(u_1)$ являются нормальными замкнутыми (n-1)-мерным и (m-1)-мерным пространствами соответственно.

Второй этап доказательства.

Докажем, что любую связь из вида $((v_1\otimes u_s)(v_p\otimes u_1))\in(v_1\otimes O(u_1)\oplus O(v_1)\otimes u_1)$ из полученного окаема можно заменить на связь $((v_1\otimes u_1)(v_p\otimes u_s))$. Проделаем это следующим образом: установим сначала связь $((v_1\otimes u_1)(v_p\otimes u_s))$, а затем отбросим связь $((v_1\otimes u_s)(v_p\otimes u_1))$. Связь между точками $(v_1\otimes u_1)$ и $(v_p\otimes u_s)$ можно установить, так как их общий окаем имеет вид конуса с двумя смежными вершинами $v_1\otimes u_s$ и $v_p\otimes u_1$, то есть $O((v_1\otimes u_1)(v_p\otimes u_s))=v_1\otimes u_s\oplus v_p\otimes u_1\oplus(v_1\otimes O(u_1u_s)\oplus O(v_1v_p)\otimes u_1)$. После установления этой связи общий окаем вершин $v_1\otimes u_s$ и $v_p\otimes u_1$ также станет конусом с двумя общими вершинами $v_1\otimes u_1$ и $v_p\otimes u_s$,
$O((v_1\otimes u_s)(v_p\otimes u_1))=v_1\otimes u_1\oplus v_p\otimes u_s\oplus$
$[(v_1\otimes O(u_1u_s)\oplus O(v_1v_p)\otimes u_s)\cup(v_p\otimes O(u_1u_s)\oplus O(v_1v_p)\otimes u_1)\cup O(v_1v_p)\otimes O(u_1 u_s)]$.
Разумеется, что подпространство в квадратных скобках может отличаться от соответствующего подпространства до начала преобразований. Поскольку общий окаем является точечным пространством, связь $((v_1\otimes u_s)(v_p\otimes u_1))$ может быть отброшена.

Третий этап доказательства.

Мы покажем, что при применении этих двух преобразований окаем точки $v_1\otimes u_1$ остается нормальным замкнутым (n+m-1)-мерным пространством. Рассмотрим $O(v_1\otimes u_1)=v_1\otimes O(u_1)\oplus O(v_1)\otimes u_1$ отдельно от остального пространства. Предшествующая пара преобразований, то есть установление связи $((v_1\otimes u_1)(v_p\otimes u_s))$ и отбрасывание связи $((v_1\otimes u_s)(v_p\otimes u_1))$ в $G^m\otimes H^n$, означает просто гомеоморфное преобразование замена связи на точку в нормальном пространстве $O(v_1\otimes u_1)=v_1\otimes O(u_1)\oplus O(v_1)\otimes u_1$. Как известно, при таком преобразовании пространство остается нормальным, и его размерность не меняется. Полезно напомнить это преобразование. Пусть G есть некоторое нормальное пространство, $v_1$ и $v_2$ есть две точки этого пространства со связью $(v_1v_2)$. Преобразование, состоящее из отбрасывания связи $(v_1v_2)$ и присоединения точки v называется гомеоморфным, если окаем O(v) точки v O(v) содержит точки $v_1$, $v_2$ и $O(v_1v_2)$, то есть $O(v)=S^0(v_1,v_2)\oplus O(v_1v_2)$, где сфера $S^0(v_1,v_2)$ состоит из двух точек $v_1$, $v_2$. В нашем случае с одной стороны общий окаем точек $(v_1\otimes u_s)$ и $(v_p\otimes u_1)$ в пространстве $A=v_1\otimes O(u_1)\oplus O(v_1)\otimes u_1$ есть подпространство $v_1\otimes O(u_1u_s)\oplus O(v_1v_p)\otimes u_1$. С другой стороны, точка $v_p\otimes u_s$ приклеивается к A по подпространству $B=(v_1\otimes u_s\cup v_p\otimes u_1)\oplus(v_1\otimes O(u_1u_s)\oplus O(v_1v_p)\otimes u_1)$, как это следует из вышеизложенного, причем между точками $(v_1\otimes u_s)$ и $(v_p\otimes u_1)$ нет связи. Это означает, что $B=S^0((v_1\otimes u_s),(v_p\otimes u_1))\oplus(v_1\otimes O(u_1u_s)\oplus O(v_1v_p)\otimes u_1)$. Таким образом окаем точки $v_1\otimes u_1$ остается нормальным (n+m-1)-мерным замкнутым пространством при отбрасывании связи $((v_1\otimes u_s),(v_p\otimes u_1))$ и установлении связи $((v_1\otimes u_1)(v_p\otimes u_s))$.

Четвертый этап доказательства.

Очевидно, что такую же операцию можно проделать для любой связи $((v_1\otimes u_s)(v_p\otimes u_1))$ из $O(v_1\otimes u_1)$. Иными словами, связь $((v_1\otimes u_s)(v_p\otimes u_1))$ можно отбросить, введя при этом связь $((v_1\otimes u_1)(v_p\otimes u_s))$. Окаем $O(v_1\otimes u_1)$ точки $v_1\otimes u_1$ все время остается нормальным замкнутым (n+m-1)-мерным пространством.

Пятый этап доказательства.

Из предшествующего ясно что подобную процедуру можно провести для

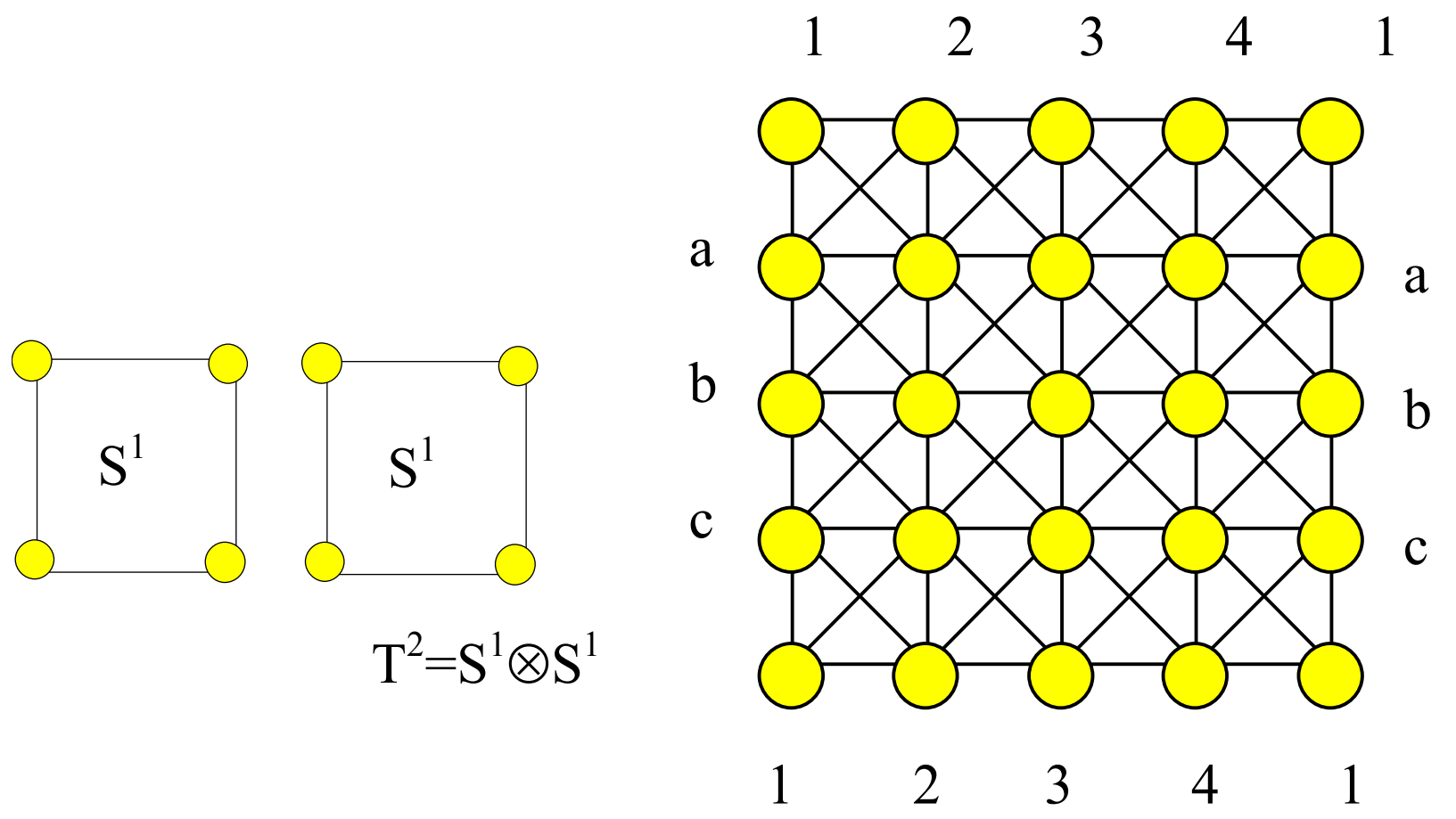

Рис. 98 В теории молекулярных пространств прямое произведение двух окружностей $S^1\otimes S^1$ дает двумерный тор $T^2$, так же как и в классической топологии.

любой точки из пространства $G^m\otimes H^n$. При этом из каждой пары связей вида $((v_i\otimes u_k)(v_p\otimes u_s))$ и $((v_i\otimes u_s)(v_p\otimes u_k))$, $i\neq p$, $k\neq s$, остается только одна связь, а вторая отбрасывается. После всех отбрасываний окаем любой точки получившегося пространства будет нормальным замкнутым (n+m-1)-мерным пространством. Следовательно $G^m\otimes H^n$, подвергшееся таким точечным преобразованиям, перейдет в нормальное замкнутое (n+m)-мерное пространство. Теорема доказана. □

Прямое произведение $L^1{}_1\otimes L^1{}_2$ двух нормальных одноименных прямых $L^1{}_1$и $L^1{}_2$ является двумерным пространством, но не является нормальным пространством. Отбрасывание точечных связей превращает это пространство в нормальное двумерное пространство $nob(L^1{}_1\otimes L^1{}_2)$ (Рис. 97). В качестве еще одной иллюстрации к этой теореме рассмотрим прямое произведение двух нормальных окружностей. В классической топологии это есть двумерный тор. Точно также в ТМП это будет двумерным тором, но уже содержащим некоторый избыток связей. В качестве окружностей возьмем их минимальный случай из 4-х точек. Прямое произведение двух окружностей $S^1 \otimes S^1$ дает двумерный тор $T^2$, что видно на Рис. 98, Рис. 99. Пространство $T^2=S^1 \otimes S^1$ стягивается точечными отбрасываниями наклонных вправо связей к нормальному двумерному молекулярному тору.

Рассмотрим, теперь, прямое произведение трех и более нормальных пространств.

<u>Определение нормального образа прямого произведения $G_1 \otimes G_2 \otimes ... \otimes G_k$ трех или более нормальных пространств.</u>

Пусть имеется набор нормальных пространств $G_1$, $G_2$,... $G_k$ размерности $n_1$, $n_2$,...$n_k$ соответственно. Нормальным образом прямого произведения $G_1 \otimes G_2 \otimes ... \otimes G_k$ этих пространств (НОБ) называется нормальное пространство $nob(G_1 \otimes G_2 \otimes ... \otimes G_k)$ размерности $n_1+n_2+...+n_k$, полученное из $G_1 \otimes G_2 \otimes ... \otimes G_k$ точечным отбрасыванием или установлением связей между точками пространства.

НОБ может быть получен различными способами, и все они гомотопны один другому, окаемы данной точки в различных НОБ также гомотопны. Однако, в большинстве случаев окаемы данной точки не изоморфны. Рассмотрим один из способов получения нормального образа прямого произведения $G_1 \otimes G_2 \otimes ... \otimes G_k$, который легко может быть алгоритмизирован и использован при компьютерной обработке.

Теорема 90

*Для любого набора нормальных пространств $G_1$, $G_2$,... $G_k$ размерности $n_1$, $n_2$,...$n_k$ соответственно всегда существует нормальное пространство $nob(G_1 \otimes G_2 \otimes ... \otimes G_k)$ размерности $n_1+n_2+...+n_k$, гомотопное прямому произведению $G_1 \otimes G_2 \otimes ... \otimes G_k$, и полученное из $G_1 \otimes G_2 \otimes ... \otimes G_k$ точечным отбрасыванием или установлением связей между точками пространства.*

Доказательство

Пусть имеется набор нормальных пространств $G_1$, $G_2$,... $G_k$ размерности $n_1$, $n_2$,...$n_k$ соответственно. Построим прямое произведение $W=G_1 \otimes G_2 \otimes ... \otimes G_k$ этих пространств.

Первый этап.

Возьмем прямое произведение двух пространств $G_1 \otimes G_2$ и в соответствии с одной из предшествующих теорем точечным отбрасыванием (установлением) связей стянем $G_1 \otimes G_2$ в $nob(G_1 \otimes G_2)$. Поскольку $G_1 \otimes G_2 \sim nob(G_1 \otimes G_2)$ то, опять таки в соответствии с одной из предшествующих теорем, пространство $G_1 \otimes G_2 \otimes G_3$ стягивается к пространству $nob(G_1 \otimes G_2) \otimes G_3$ точечными отбрасываниями связей.

Второй этап.

Рассмотрим прямое произведение двух нормальных пространств $nob(G_1 \otimes G_2) \otimes G_3$, и в соответствии с одной из предшествующих теорем точечным отбрасыванием связей стянем $nob(G_1 \otimes G_2) \otimes G_3$ в $nob(nob(G_1 \otimes G_2) \otimes G_3)$, как и в предыдущем случае.

Третий этап.

Повторяя подобные рассуждения мы получаем, что $nob(...nob(nob(G_1 \otimes G_2) \otimes G_3) \otimes ...G_k) \sim G_1 \otimes G_2 \otimes G_3 \otimes ...G_k$. Теорема доказана.□

Естественно, отбрасывание связей можно осуществить многими способами, и все полученные нормальные пространства имеют одинаковую размерность и гомотопны друг другу.

Алгоритм получения нормального образа прямого произведения $G_1 \otimes G_2 \otimes ... \otimes G_k$ нормальных пространств.

Рассмотренный в теореме метод можно назвать послойным или покоординатным приведением к нормальному виду-ППНВ. Вначале мы приводим к нормальному виду подпространство $G_1 \otimes G_2$, затем подпространство $G_1 \otimes G_2 \otimes G_3$ и так далее вплоть до $G_1 \otimes G_2 \otimes G_3 \otimes ...G_k$.

Применяя эту теорему к конкретным преобразованиям выбираем каждую пару смежных точек $v_s$, $v_t$ из $G_1$ и каждую пару смежных точек $u_r$, $u_q$ из $G_2$ и в $G_1 \otimes G_2 \otimes G_3 \otimes ...G_k$ отбрасываем связи между точками вида $v_s \otimes u_r \otimes ... \otimes ...$, сохраняя связи между точками вида $v_t, \otimes u_q \otimes ... \otimes ...$. При этом подпространство $G_1 \otimes G_2$ в $G_1 \otimes G_2 \otimes G_3 \otimes ...G_k$ приводится к нормальному виду, то есть переходит в пространство $nob(G_1 \otimes G_2)$, а $G_1 \otimes G_2 \otimes G_3 \otimes ...G_k$ переходит в $nob(G_1 \otimes G_2) \otimes G_3 \otimes ...G_k$. Затем выбираем каждую пару смежных точек w(s,t), w(r,q) из $nob(G_1 \otimes G_2)$ и каждую пару смежных точек $u_b$, $u_c$ из $G_3$ и в $G_1 \otimes G_2 \otimes G_3 \otimes ...G_k$ отбрасываем связи между точками вида $w(s,t) \otimes u_b \otimes ... \otimes ...$, сохраняя связи между точками вида $w(r,q) \otimes u_c \otimes ... \otimes ...$.В этом случае слой из трех подпространств приведен к нормальному виду, то есть $G_1 \otimes G_2 \otimes G_3$ в $G_1 \otimes G_2 \otimes G_3 \otimes ...G_k$ приводится к нормальному виду, то есть переходит в пространство $nob(G_1 \otimes G_2 \otimes G_3)$, а $G_1 \otimes G_2 \otimes G_3 \otimes ...G_k$ переходит в $nob(G_1 \otimes G_2 \otimes G_3) \otimes G_3 \otimes ...G_k$.

Таким образом переходя от слоя к слою мы превращаем $G_1 \otimes G_2 \otimes G_3 \otimes ...G_k$ в нормальное пространство $nob(G_1 \otimes G_2 \otimes G_3 \otimes G_3 \otimes ...G_k)$. Очевидно, что каждый слой сам по себе имеет произвольную размерность. В случае, если речь идет об одномерных нормальных пространствах, возникает полная аналогия с прямым произведением непрерывных пространств. Нечто, подобное прямому произведению было использовано в работах Халимского, Конга, Коппермана и других [30,31,34]. Однако, поскольку они использовали по сути дела неоднородное одномерное пространство и специальным образом определенное произведение пространств, уже в двумерном случае возникало возникали точки различной топологии, и причем число топологически различных точек нарастало примерно по экспоненте в зависимости от числа сомножителей в произведении. При этом подходе является проблемой определить даже простую окружность, не говоря уже о замкнутых пространствах более высокой размерности. В теории молекулярных пространств эти вопросы не возникают вообще. Рассмотрим, теперь, применение аксиомы о гомотопии нормальных n-мерных пространств к прямому произведению. Напомним эту аксиому. Если пространство не имеет точечных точек и связей и гомотопно нормальному n-мерному пространству, то оно само является нормальным n-мерным пространством. В применении к нашему случаю сформулируем теорему.

Теорема 91

*Прямое произведение $G_1 \otimes G_2 \otimes ... \otimes G_k$ любого набора нормальных пространств $G_1$, $G_2$,... $G_k$ размерности $n_1$, $n_2$,...$n_k$ соответственно всегда приводится к нормальному пространству $nob(G_1 \otimes G_2 \otimes ... \otimes G_k)$ размерности $n_1+n_2+...+n_k$, точечным отбрасыванием связей между смежными точками пространства в произвольном порядке и последующем отбрасыванием висячих точек, имеющих точечный окаем, состоящий из единственной точки.*

Доказательство

Пусть имеется набор нормальных пространств $G_1$, $G_2$,... $G_k$ размерности $n_1$, $n_2$,...$n_k$ соответственно. Построим прямое произведение $W=G_1 \otimes G_2 \otimes ... \otimes G_k$ этих пространств. Согласно предыдущему оно гомотопно нормальному молекулярному пространству размерности $n_1+n_2+...+n_k$. Кроме того, точечное отбрасывание и установление связей между точками приводит к тому, что этих связей не остается вообще, при этом окаемы некоторых точек становятся неточечными, а окаемы некоторых других становятся точечными. Точки, имеющие точечные окаемы, очевидно, являются висячими, то есть имеющими только одну смежную точку. Если их отбросить, то согласно аксиоме о гомотопности нормальных пространств получившееся пространство будет нормальным пространством размерности $n_1+n_2+...+n_k$, поскольку не содержит точечных связей и точек. Теорема доказана.□

Прямое произведение молекулярных пространств имеет то преимущество по сравнению с прямой суммой, что это хорошо изученная и используемая операция в классической математике. Во многих случаях свойства непрерыевных и молекулярных пространств, являющихся их моделями, совпадают. Хотя этот факт кажется очевидным, тем не менее он будет нуждаться в доказательстве в каждом конкретном случае, прежде, что будет применяться в практике. Также.как и прямая сумма, прямое

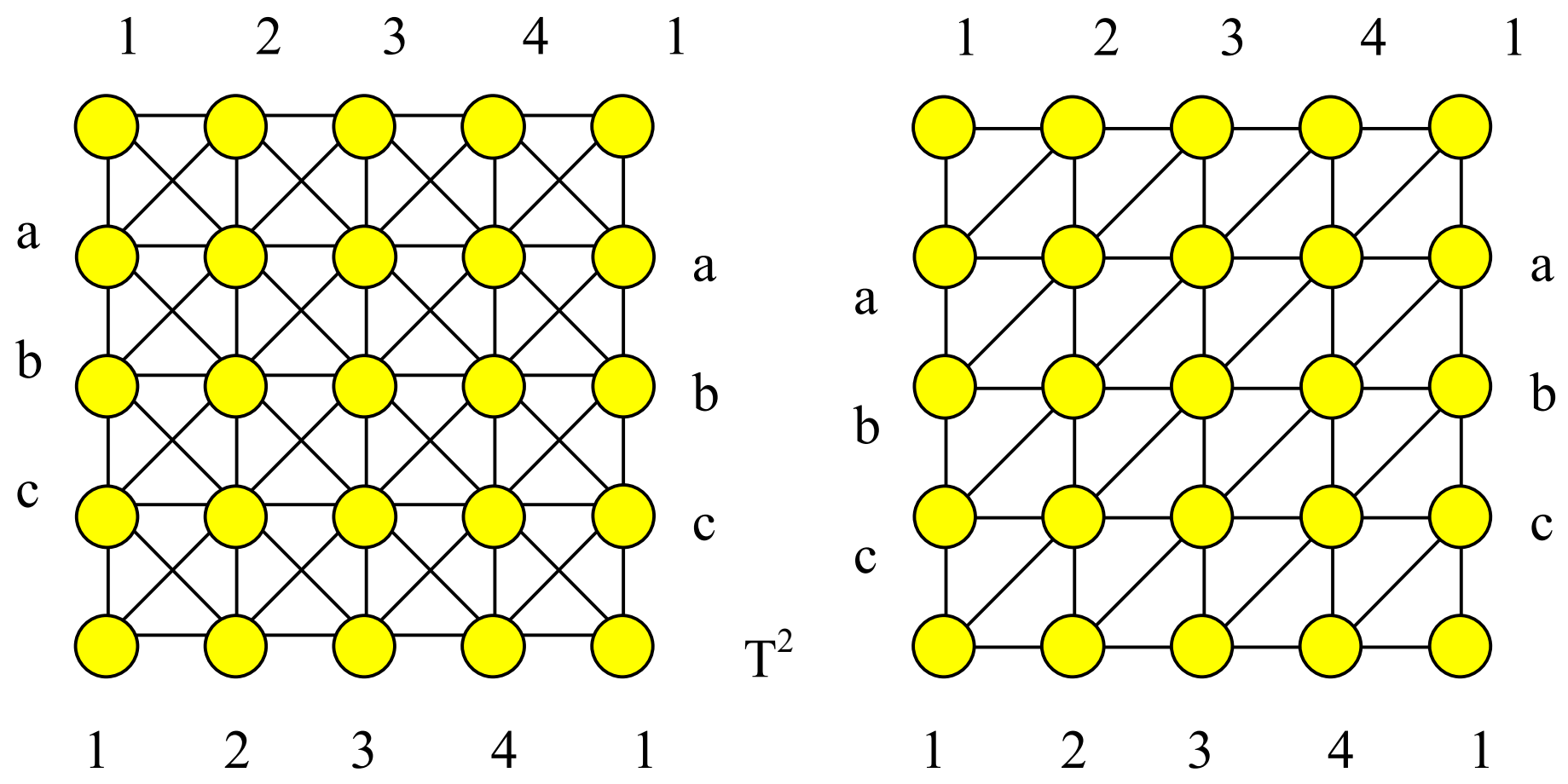

Рис. 99 Пространство $T^2=S^1\otimes S^1$ стягивается точечными отбрасываниями наклонных вправо связей к нормальному двумерному молекулярному тору.

произведение дает многообразия только в тем случаях, где сомножителями являются многообразия, то есть пространства, в которых окаемы являются сферами. Если же хотя бы один из сомножителей не является многообразием, то прямое произведение также не будет многообразием, поскольку в нем окаем точки $v\otimes u$ гомотопен прямой сумме $O(v)\oplus O(u)$, и, следовательно, не есть сфера, если хотя бы один из исходных окаемов не является сферой. Примечательным является тот факт, что прямое произведение двух нормальнах пространств уже не является нормальным пространством. Исходя из этого можно сказать, что хотя нормальные пространства очень удобны для работы, не они являются моделями, наиболее близкими по своим к непрерывным пространствам, а регулярные пространства.

Список литературы к главе 7.

# n-МЕРНЫЕ МОЛЕКУЛЯРНЫЕ ПРОСТРАНСТВА НА БАЗЕ НОРМАЛЬНЫХ МОЛЕКУЛЯРНЫХ ПРОСТРАНСТВ

We consider general molecular (regular) n-dimensional spaces and study their properties.

Некоторые результаты этой главы изложены в [43].

Перейдем определению и свойствам n-мерных пространств, общего вида. Как уже неоднократно говорилось, любая модель пространства должна

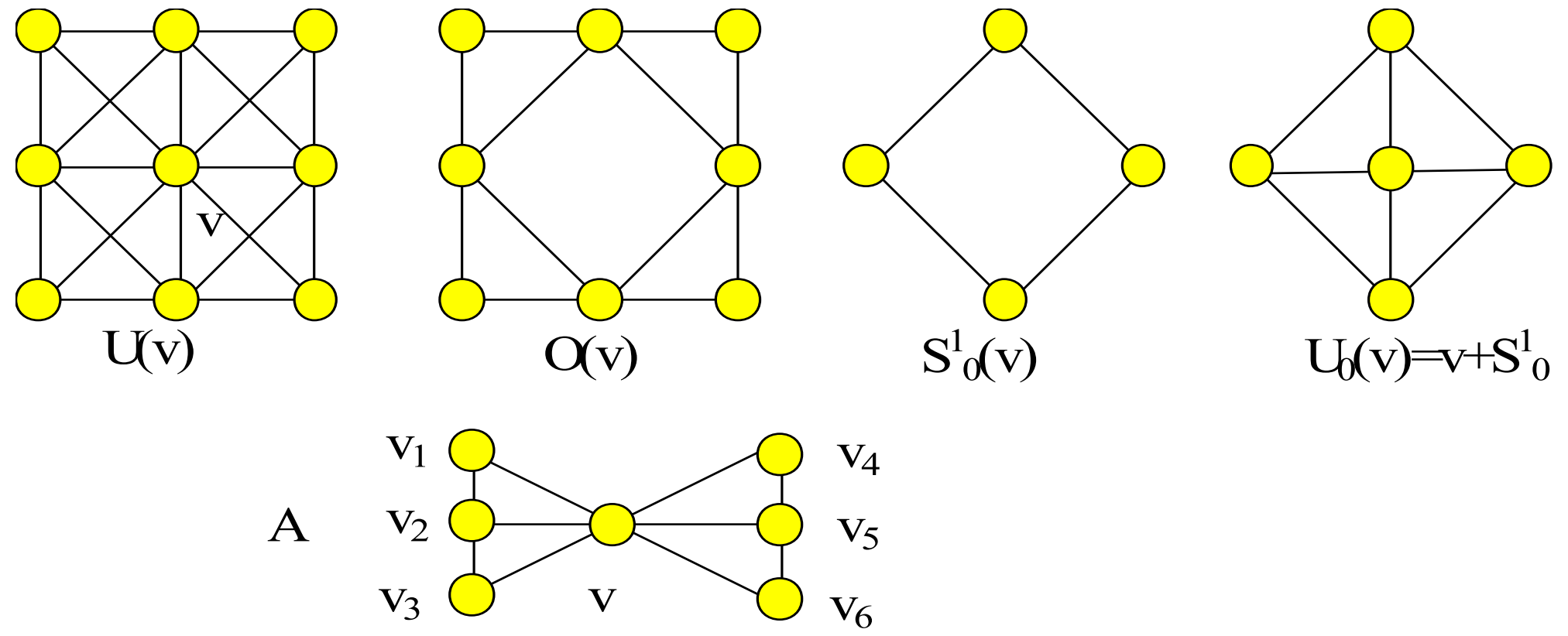

Рис. 100 Окаем O(v) точки v содержит нормальную одномерную сферу $S^1_0$. Остовный шар $U_0(v)$ точки v определяется выражением $U_0(v)=v\oplus S^1_0$. В пространстве A остовные шары точки v есть подпространства вида $(v_k,v,v_p)$, где k=1,2,3, p=4,5,6.

подтверждаться какими-то опытами. Компьютерные эксперименты описаны в одной из следующих глав. Как следствие этих экспериментов был сделан вывод, что молекулярные модели непрерывных пространств имеют ту же размерность, что и их непрерывные прообразы. В качестве молекулярных образов использовались нервы покрытий или разбиений пространств. При этом возникали молекулярные модели более сложного вида, чем нормальные. Появляется необходимость дать определение размерности для такого вида пространств. Так как в покрытиях или разбиениях многое зависит от формы элементов, то необходимо определить n-мерное пространство и n-мерную точку не привязывая эти определения к какому-то определенному виду покрытия.

З а м е ч а н и е

Для определения размерности произвольного пространства мы используем нормальные пространства. Однако, нормальные пространства могут быть заменены на регулярные и, в частности, на квазинормальные.

Определение n-мерной точки

Точка v пространства A называется n-мерной, если окаем O(v) этой точки содержит (n-1)-мерное замкнутое нормальное пространство $G^{n-1}$, и не содержит n-мерного замкнутого нормального пространства $G^n$, $G^{n-1} \subseteq O(v)$, $G^n \not\subseteq O(v)$.

Очевидно, что все точки любого n-мерного нормального пространства без края являются n-мерными точками в смысле определения, данного выше. В n-мерном нормальном пространстве с краем внутренние точки являются n-мерными, а краевые являются (n-1)-мерными.

Определение остовного шара $U_0(v)$ n-мерной точки v

Пусть точка v пространства A является n-мерной, и $G^{n-1} \subseteq O(v)$, есть некоторое нормальное замкнутое (n-1)-мерное пространство. Тогда подпространство $v \oplus G^{n-1} = U_0(v)$, где $v \oplus G^{n-1} = U_0(v) \subseteq U(v) = v \oplus O(v)$ называется остовным шаром этой точки.

На *Рис. 100* показан двумерный шар точки v и остовный шар этой же точки, являющийся прямой суммой точки v и нормальной минимальной окружности $S^1_0$ из 4-х точек, $U_0(v) = v \oplus S^1_0$. Точка может иметь несколько остовных шаров, как это видно из *Рис. 100*.

Определение n-мерного пространства

Пространство G называется n-мерным, если в нем существует по крайней мере одна n-мерная точка, и не существует (n+1)-мерной точки.

Примером двумерного пространства является U(v) на *Рис. 100*. Точка v является двумерной, все остальные точки являются либо одномерными либо нульмерными. Таким образом размерность пространства определяется по наивысшей размерности точек, его составляющих.

Как уже говорилось, точечные преобразования молекулярных пространств вообще говоря не сохраняют размерность пространства. Однако при работе с дигитальными моделями непрерывных пространств естественным требованием во многих случаях является сохранение размерности пространства. Например, в компъютерных играх объекты являются движущимися трехмерными телами и не меняют своей размерности во время работы программы. Поэтому необходимо наложить такие дополнительные ограничения на точечные преобразования, которые позволяют сделать размерность инвариантом, не меняющимся при точечных преобразованиях.

В качестве примера рассмотрим точечное отбрасывание (приклеивание) точки к пространству (*Рис. 101*). Точка u является нормальной трехмерной точкой, пространство A есть шар этой точки. Отбросим точку v из A. Точка u станет двумерной. Отбросим вторую точку w. Размерность точки u останется прежней, точка u будет двумерной. При первом отбрасывании размерность пространства уменьшилась, а при втором не изменилась.

Таким образом возможны такие точечные отбрасывания (приклеивания) точек, при которых размерность пространства не меняется. Несколько более сложная ситуация возникает при отбрасывании и введении связей. Отбрасывание или введение точечной связи может увеличить или уменьшить размерность точки v на 1 (*Рис. 103*). Докажем вначале несколько простых теорем, которые будут использоваться ниже. Еще раз

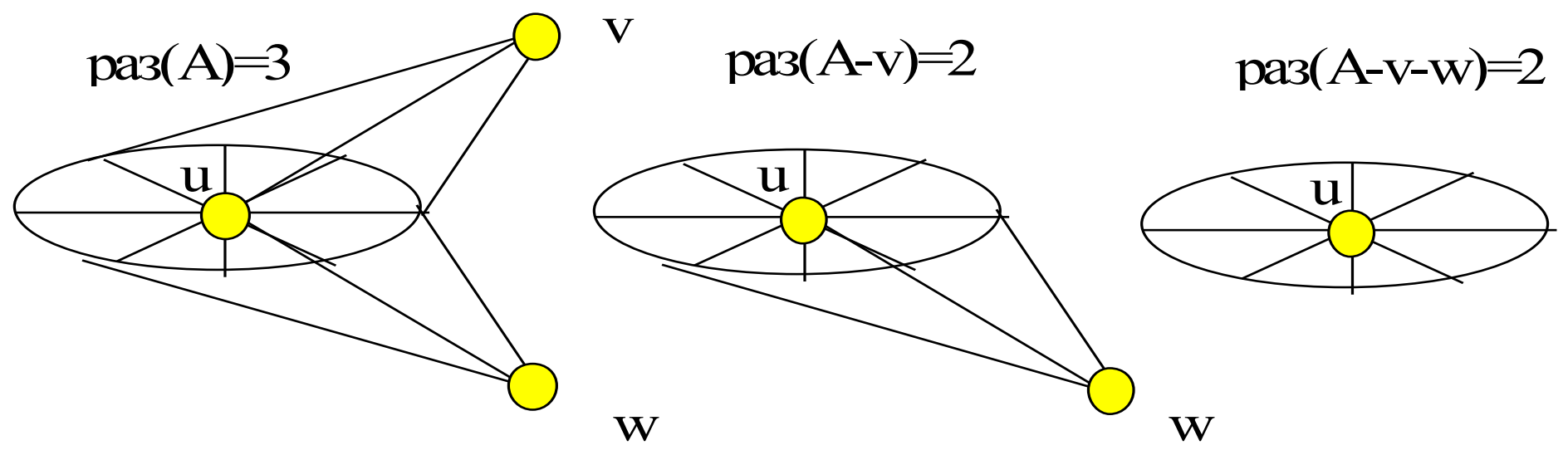

Рис. 101 При отбрасывании точки v из трехмерного пространства A размерность пространства уменьшается на 1. При отбрасывании точки w из двумерного пространства A-v размерность пространства не меняется. Обратная операция приклеивания точки или же не меняет размерность пространства, ил увеличивает размерность на 1.

отметим, что могут использоваться квазинормальные пространства, но это усложнит доказательства. Тем не менее все последующие результаты справедливы и для n-мерных пространств общего вида.

Теорема 92

*При любом отбрасывании точки v из n-мерного пространства A размерность точек пространства или не меняется или уменьшается на 1.*

Доказательство

Смысл теоремы ясен из *Рис. 101*. Размерность точек, не принадлежащих O(v) не меняется. Рассмотрим как меняется размерность n-мерной (или (n-k)-мерной) точки u, принадлежащей O(v). Пусть существует остовный n-мерный шар $U_0(u)$ в A, не содержащий v. Тогда размерность точки u не меняется при отбрасывании точки v.

Пусть не существует никакого остовного n-мерного шара $U_0(u)$ в A, не содержащего v. То есть точка v содержится в любом нормальном замкнутом (n-1)-мерном подпространстве $G^{n-1}$, принадлежащем O(v). Тогда при отбрасывании точки v любое $G^{n-1}$ переходит в пространство со связным краем $G^{n-1}$-v. Окаем O(v) в A-v уже не содержит никакого нормального замкнутого (n-1)-мерного подпространства $G^{n-1}$, но при этом содержит нормальное замкнутое (n-2)-мерное подпространство $G^{n-2}$. Следовательно, размерность точки u уменьшается на 1(*Рис. 101*). Теорема доказана.□

Используя эту теорему мы можем сказать, как ведет себя размерность пространства в целом при отбрасывании точки.

С л е д с т в и е

При любом отбрасывании точки из n-мерного пространства A размерность пространства или не меняется или уменьшается на 1.

Следствие почти очевидно. Если точка v является единственной n-мерной точкой в пространстве A, то размерность пространства однозначно уменьшается на 1. С другой стороны размерность пространства A не может уменьшиться более чем на 1, поскольку O(v) содержит нормальное замкнутое (n-1)-мерное подпространство $G^{n-1}$. Если отбрасываемая точка v не является n-мерной, то размерность n-мерной точки, содержащейся в пространстве, не может уменьшиться более чем на 1 согласно предыдущей теореме. Теорема доказана.□

Смысл этой теоремы в том, что любое отбрасывание точки может или уменьшить размерность пространства на 1, или не же не изменить ее, но при таком преобразовании невозможно увеличение размерности пространства. Здесь возникает одно интересное продолжение. Предположим, что пространство A имеет только одну n-мерную точку v. Может ли при этом быть окаем O(v) точки v точечным пространством. Эта задача не решена. Сформулируем ее отдельно.

З а д а ч а

Доказать или опровергнуть следующее утверждение: Если n-мерное пространство A имеет единственную n-мерную точку v, то окаем этой точки O(v) не может быть точечным пространством.

Теорема 93

*При любом приклеивании точки v к n-мерному пространству A размерность точек пространства или не меняется или увеличивается на 1.*

Д о к а з а т е л ь с т в о

Теорема также ясна из *Рис. 101*. Для любой смежной с v точки u к ее окаему O(u) также приклеивается v. Размерность пространства O(u) не превышает (n-1) по условию теоремы. Предположим, что точка u есть n-мерная точка, Это означает, что ее окаем содержит некоторое замкнутое (n-1)-мерное нормальное пространство, но не содержит никакого замкнутого нормального n-мерного пространства, Предположим, также, что ее окаем содержит некоторое n-мерное нормальное пространство $D^n$ с краем $B=G^{n-1}$, являющимся замкнутым нормальным (n-1)-мерным пространством. Пусть B содержится в O(v). Тогда приклеивание точки v к A есть в то же самое время приклеивание точки v к $(u\oplus B)$. При этом $D^n\cup v$ становится замкнутым нормальным n-мерным пространством, принадлежащим O(u). Это означает, что размерность точки u увеличивается на 1 и становится равной (n+1). То же самое рассуждение применимо, если размерность точки u в A равна (n-k), где $0<k<n$. Если вышеприведенные условия не выполняются, то размерность точки u не меняется (*Рис. 101*). Теорема доказана.□

С л е д с т в и е

При любом приклеивании точки v к n-мерному пространству A размерность пространства или не меняется или увеличивается на 1.

Пусть к n-мерному пространству A каким-то образом приклеивается некоторая точка v. Если O(v) содержит нормальное замкнутое n-мерное пространство $G^n$, то точка v становится (n+1)-мерной, и пространство становится (n+1)-мерным. Точки, принадлежащие пространству A являются не более чем n-мерными, и не могут увеличить свою размерность

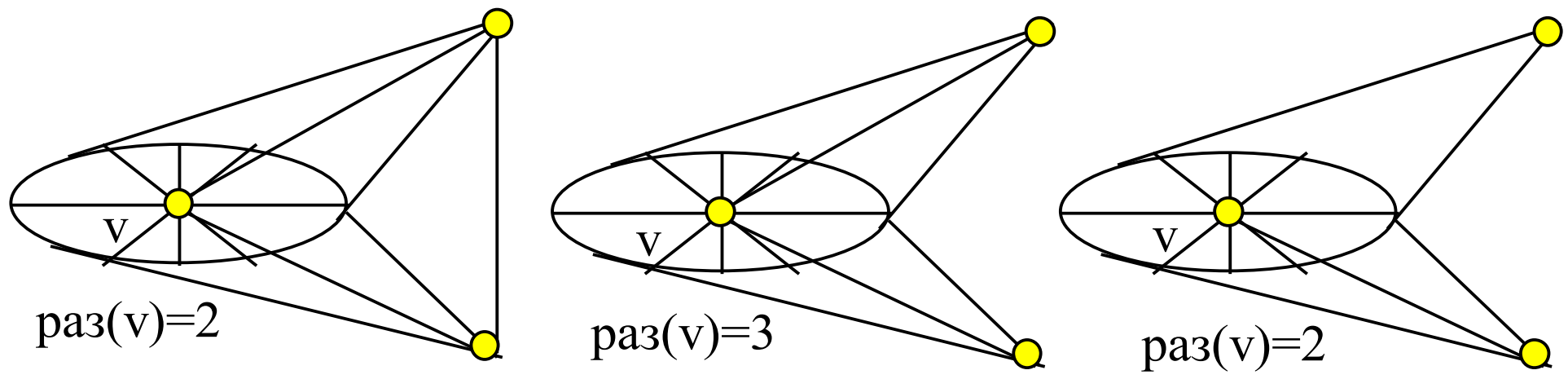

*Рис. 102 При отбрасывании связи размерность точки может как увеличиваться так и уменьшаться на 1.*

не более чем на 1. Таким образом размерность пространства может увеличиться не более чем на 1 (*Рис. 101*).

Сформулируем задачу, решение которой может помочь в составлении программ по обработке дигитальных моделей n-мерных пространств.

З а д а ч а

Доказать или опровергнуть следующее утверждение: Если n-мерное пространство A содержит нормальное замкнутое n-мерное подпространство $G^n$, то A не является точечным пространством.

Рассмотрим теоремы, касающиеся отбрасывания и приклеивания связи к пространству.

Теорема 94

*При отбрасывании любой связи в n-мерном пространстве A размерность точек может увеличиваться или уменьшаться на 1, или же сохраняться (*Рис. 102*).*

Д о к а з а т е л ь с т в о

При отбрасывании связи могут менять свою размерность точки $v_1$, $v_2$ и точки, принадлежащие $O(v_1v_2)$. Точки $v_1$ и $v_2$ или не меняют или уменьшают свою размерность согласно предыдущей теореме.

Пусть точка v, принадлежащая $O(v_1v_2)$, является n-мерной, и хотя бы один ее остовный шар $U^n{}_0(v)$ также содержится в $O(v_1v_2)$. Тогда после отбрасывания связи остовный шар точки v примет вид $U_1(v)=S^0(v_1v_2)\oplus U^n_0(v)=S^0(v_1v_2)\oplus v\oplus G^{n-1}=v\oplus G^n$. Следовательно, точка v станет (n+1)-мерной точкой.

Пусть точка v, принадлежащая $O(v_1v_2)$, является n-мерной, и все ее остовные шары $U^n_0(v)$ содержат точки $v_1$ и $v_2$. То есть каждый остовный

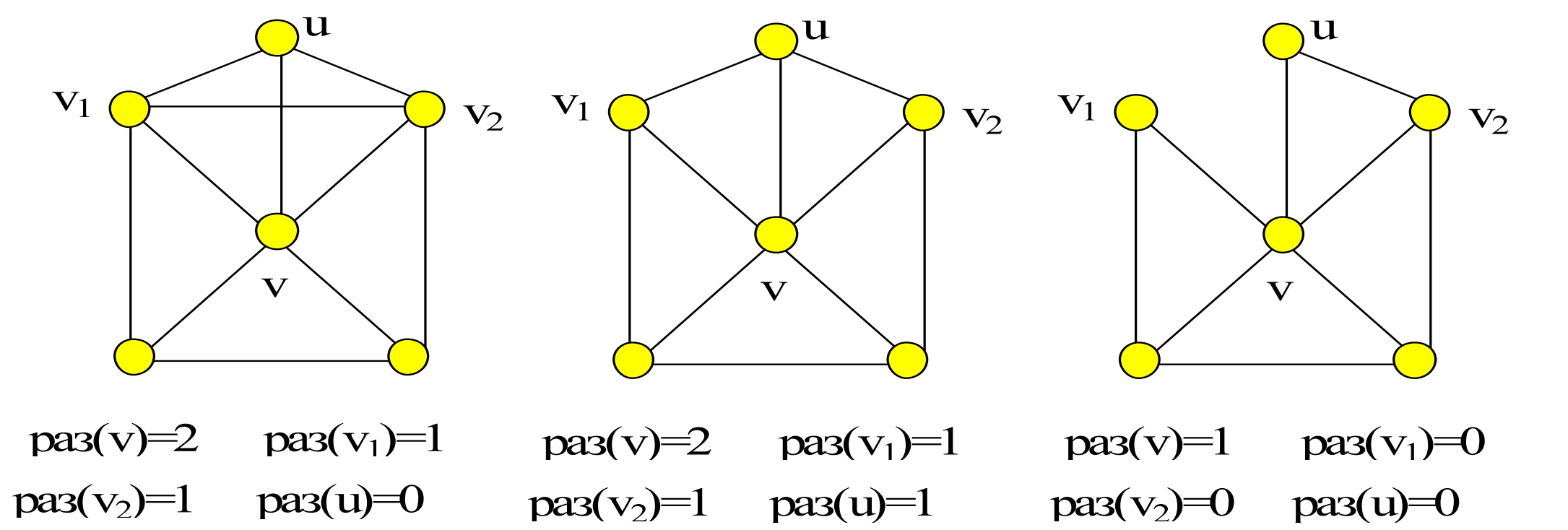

Рис. 103 При отбрасывании связей и размерность точек v и может сохраняться или уменьшаться на 1, размерность точек, принадлежащих общим окаемам и может увеличиваться или уменьшаться на 1,или сохраняться.

шар включает в себя точки $v_1$ и $v_2$. Тогда после отбрасывания связи любой остовный шар точки v примет вид $U_1(v)=U^n_0(v)-(v_1v_2)=v\oplus(G^{n-1}-(v_1v_2))$. Пространство $(G^{n-1}-(v_1v_2))$ уже не является нормальным замкнутым (n-1)-мерным пространством, становится нормальным (n-1)-мерным пространством со связным краем, хотя и содержит (n-2)-мерное нормальное замкнутое подпространство. Следовательно, точка v станет (n-1)-мерной точкой.

Пусть точка v, принадлежащая $O(v_1v_2)$, является n-мерной, ни один ее остовный шар $U^n_0(v)$ не содержится в $O(v_1v_2)$, и существует остовный шар $U^n_0(v)$, не содержащийся в $v_1\oplus v_2\oplus O(v_1v_2)$, Тогда после отбрасывания связи размерность точки v не увеличится и не уменьшится, то есть точка v останется n-мерной точкой.. Теорема доказана.□

Как следствие этой теоремы выделим более подробно, как может меняться размерность точек при отбрасывании связей.

Следствие

Пусть в n-мерном пространстве A отбрасывается связь $(v_1v_2)$ между смежными точками $v_1$ и $v_2$,. Тогда:

- размерность точек $v_1$ и $v_2$ или не меняется или уменьшается на 1
- размерность точек v, принадлежащих общему окаему $O(v_1v_2)$ может увеличиваться или уменьшаться на 1, или же сохраняться.

Эта теорема иллюстрируется на *Рис. 102* и *Рис. 103*. Отбрасывание связи между двумя точками является более сложным процессом, чем отбрасывание точки, поскольку фактически вовлекает также и отбрасывание точки. С одной стороны для точек $v_1$, и $v_2$ это означает отбрасывание точек $v_2$, и $v_1$ из пространств $U(v_1)$ и $U(v_2)$, то есть уменьшение размерности пространства A или ее сохранение. С другой стороны для точек v, принадлежащих $O(v_1v_2)$, эта же операция ведет к увеличению, уменьшению или сохранению размерности. В конечном итоге

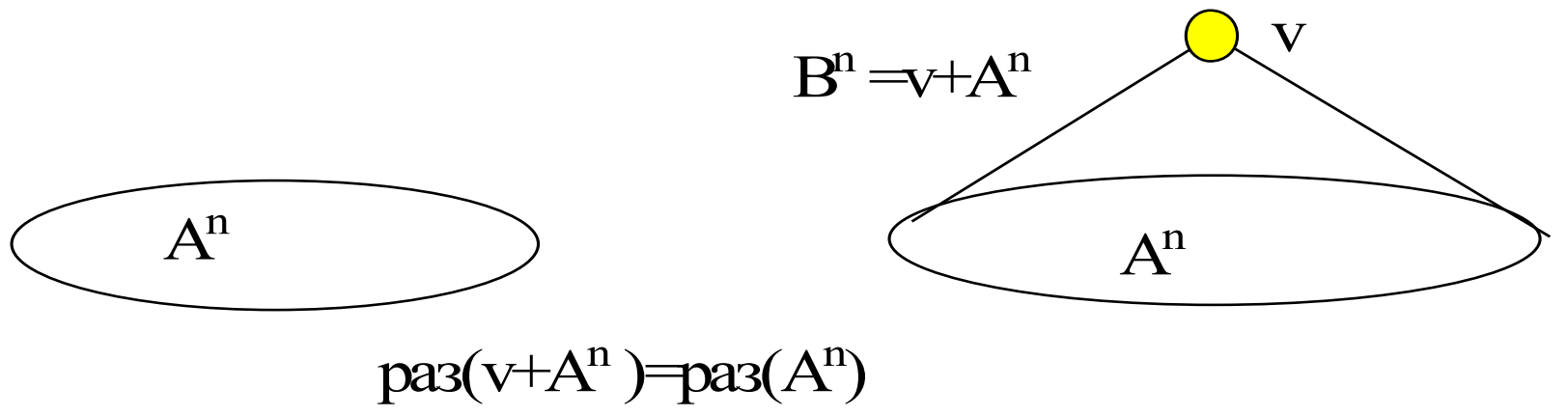

Рис. 104 Размерность пространства $A^n$ равна размерности пространства $v\oplus A^n$, если $A^n$ не содержит нормальное замкнутое пространство $G^n$ как подпространство.

возможно увеличение размерности некоторых точек на 1, и уменьшение размерности некоторых других также на 1. Следовательно, размерность пространства также может увеличиться или уменьшиться на 1. Аналогичный результат получается при введении связи между точками пространства.

Теорема 95

*При введении любой связи в n-мерном пространстве A размерность точек может увеличиваться или уменьшаться на 1, или же сохраняться(*Рис. 102, Рис. 103*).*

Доказательство

При установлении связи могут менять свою размерность точки $v_1$, $v_2$ и точки, принадлежащие $O(v_1v_2)$. Точки $v_1$ и $v_2$ или не меняют или увеличивают на 1 свою размерность согласно предыдущей теореме.

Пусть точка v, принадлежащая $O(v_1v_2)$, является n-мерной, и все ее остовные шары $U^n{}_0(v)$ содержатся в $S^0(v_1v_2)\oplus O(v_1v_2)$, но не содержатся в $O(v_1v_2)$. То есть каждый остовный шар включает в себя точки $v_1$ и $v_2$. Очевидно, что при этом $U^n{}_0(v)\cap O(v_1v_2)=v\oplus G^{n-1}$. Тогда после введения связи $(v_1v_2)$ любой остовный шар точки v примет вид $U_1(v)=U^n{}_0(v)-v_1-v_2=v\oplus G^{n-1}$. Следовательно, точка v станет (n-1)-мерной точкой.

Пусть точка v, принадлежащая $O(v_1v_2)$, является n-мерной, и существует n-мерное нормальное пространство $D^n$ со связным краем, принадлежащее O(v), и содержащее точки $v_1$ и $v_2$, такое, что введение связи $(v_1v_2)$ переводит $D^n$ в n-мерное нормальное замкнутое пространство $G^n=D^n\oplus(v_1v_2)$. Тогда после введения связи окаем O(v) будет содержать $G^n=D^n\oplus(v_1v_2)$. Следовательно, точка v станет (n+1)-мерной точкой.

Пусть точка v, принадлежащая $O(v_1v_2)$, является n-мерной, не существует n-мерное нормальное пространство $D^n$ со связным краем, принадлежащее O(v), и содержащее точки $v_1$ и $v_2$, такое, что установление связи $(v_1v_2)$

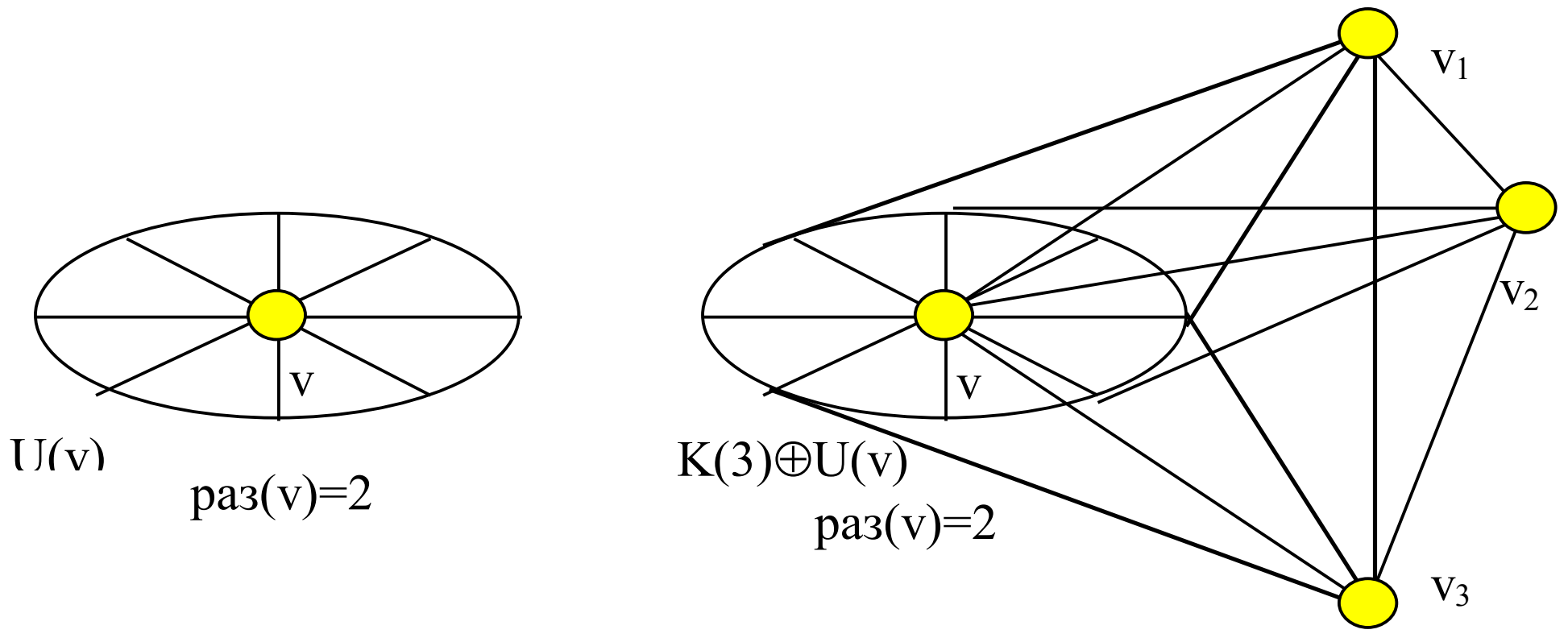

Рис. 105 Прямая сумма двумерного шара U(v) и связки K(3) является также двумерным шаром для точки v, а также точек $v_1$, $v_2$ и $v_3$.

переводит $D^n$ в n-мерное нормальное замкнутое пространство $G^n=D^n\oplus(v_1v_2)$, и существует ее остовный шар $U^n{}_0(v)$ не содержащий по крайней мере одной из точек $v_1$ или $v_2$. Тогда после введения связи размерность точки v не увеличится и не уменьшится, то есть точка v останется n-мерной точкой (*Рис. 102*, *Рис. 103*). Теорема доказана.□

Подчеркнем более подробно, как могут изменяться окаемы различных точек.

Следствие

Пусть в n-мерном пространстве A устанавливается связь $(v_1v_2)$ между несмежными точками $v_1$ и $v_2$,. Тогда:

- размерность точек $v_1$ и $v_2$ или не меняется или увеличивается на 1
- размерность точек v, принадлежащих общему окаему $O(v_1v_2)$ может увеличиваться или уменьшаться на 1, или же сохраняться (*Рис. 102*, *Рис. 103*).

Как результат предыдущих теорем мы можем сформулировать следующее общее следствие.

Следствие

Любое отбрасывание или приклеивание связи или точки меняет размерность пространства не более чем на 1.

Рассмотрим теперь точечные преобразования, не меняющие размерность пространства.

Теорема 96

*Если n-мерное пространство A не содержит никакого замкнутого n-мерного нормального подпространства $G^n$, то конус $v\oplus A$ также*

*является n-мерным пространством и не содержит никакого замкнутого n-мерного нормального подпространства $G^n$.*

Доказательство

Очевидно, что v⊕A содержит $G^{n-1}$. Предположим, что v⊕A содержит $G^n$. Так как A не содержит $G^n$, то присоединенная точка v должна входить в

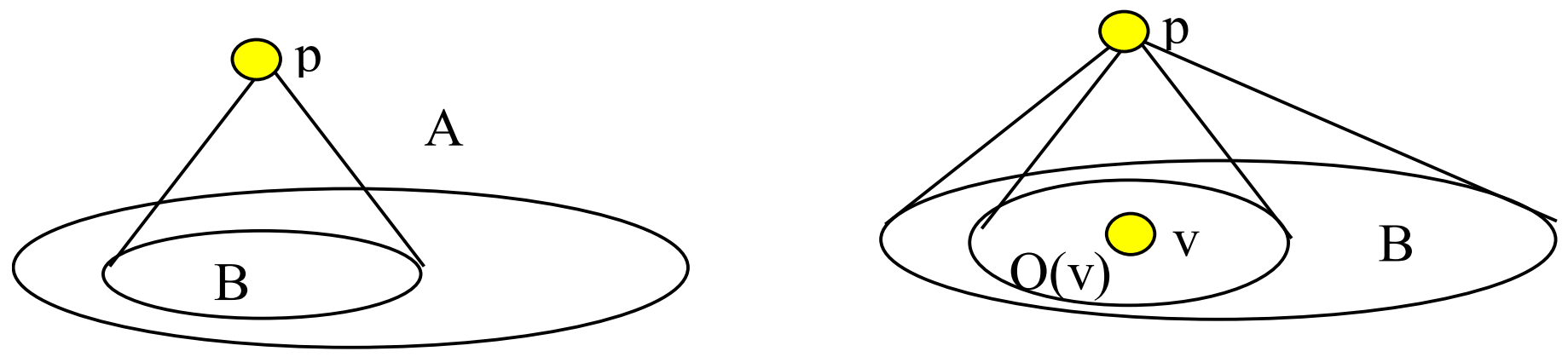

Рис. 106 Приклеивание точки p к A по B. Если остовный шар n-мерной точки v лежит в B, то весь шар точки v также лежит в B.

$G^n$, а это означает, что $G^n$=v⊕B, B⊆A. Однако, конус вообще не может быть нормальным пространством. Следовательно, предположение неверно, и v⊕A не содержит $G^n$ (*Рис. 104*). Теорема доказана.□

Эта теорема важна потому, что с одной стороны шар точки является всегда конусом, с другой стороны шар точки есть базовая конструкция для дальнейших построений.

Эту теорему можно обобщить, если взять вместо точки v связку K(m), то есть пространство из m попарно связанных точек.

Следствие

Если n-мерное пространство A не содержит никакого замкнутого n-мерного нормального подпространства $G^n$, то конус K(m)⊕A также является n-мерным пространством и не содержит никакого замкнутого n-мерного нормального подпространства $G^n$. K(m) является связкой из m точек ($v_1$,$v_2$,...$v_m$) (*Рис. 105*).

Важно подчеркнуть еще одно следствие этой теоремы применительно к шару точки.

Следствие

Если пространство $U^n$(v) есть n-мерный шар точки v, то K(m)⊕$U^n$(v), где K(m) является связкой из m точек ($v_1$,$v_2$,...$v_m$), также есть n-мерный шар точки v, а также точек ($v_1$,$v_2$,...$v_m$) (*Рис. 105*).

Рассмотрим, какие условия должны выполняться, чтобы размерность пространства сохранялась при точечных приклеиваниях к отбрасываниях точек и связей в более общем случае. При этом необходимо предварительно определить, что мы будем называть сохранением размерности пространства, поскольку от этих требований зависит ход дальнейших построений. Возможны несколько подходов. Мы будем использовать подход, который предполагает, чтобы при любом преобразовании n-мерного пространства пространство всегда имело хотя бы одну n-мерную точку. Сформулируем эти условия в виде определения.

Определение точечных преобразований с сохранением размерности

Точечные преобразования Ф n-мерного пространства A называются точечными преобразованиями с сохранением размерности, если в процессе преобразований пространство не меняет свою размерность, раз(ФA)=раз(A).

Рассмотрим более подробно, какие условия необходимо наложить на точечные преобразования, чтобы они не меняли размерность пространства.

Теорема 97

*Пусть A есть n-мерное пространство, B, $B \subseteq A$ есть его подпространство, удовлетворяющее следующим условиям:*

- *B является точечным пространством.*
- *B не содержит никакого нормального n-мерного подпространства $G^n$, $G^n \not\subset B$.*
- *если для некоторой n-мерной точки v существует ее остовный шар $U_0(v)$, целиком принадлежащий B, то шар U(v) этой точки также принадлежит B, $U_0(v) \subseteq B \Rightarrow U(v) \subseteq B$.*

*Тогда точечное приклеивание точки p к A по B не меняет размерность пространства, то есть $A \sim C = A \cup (p \oplus B)$, раз(A)=раз(C).*

Доказательство

При приклеивании точки размерность пространства не уменьшается. Покажем, что она также не увеличивается. Так как O(p)=B не содержит $G^n$, то p не может быть (n+1)-мерной точкой. Если v есть n-мерная точка и $U_0(v) \subseteq B$, то $U(v) \subseteq B$. После приклеивания точки p шар точки v перейдет в $p \oplus U(v)$. Согласно предыдущему раз(U(v))=раз($p \oplus U(v)$). Следовательно, точка v остается n-мерной точкой в $C = A \cup (p \oplus B)$ (*Рис. 106*). Предположим, что $v \in B$, и не существует (n-1)-мерного нормального подпространства из O(v), целиком лежащего в B. Тогда окрестность O(v)|C точки v в пространстве C $O(v)|C = O(v) \cup (p \oplus (O(v) \cap B))$. Если O(v)|C содержит n-мерное замкнутое подпространство $G^n$, то точка p обязана принадлежать $G^n$, отсюда $O(p)|G^n \supseteq G^{n-1}$ и $G^{n-1} \subseteq O(v) \cap B$,, что противоречит предположению. Следовательно, размерность этой точки также не превышает n в C. Предположим, что точка v не принадлежит B. Тогда окаем этой точки не меняется при переходе от пространства A к C. Следовательно, присоединение точки p к A не меняет размерность пространства. Поскольку B=O(p) является точечным пространством, то A и C гомотопны. Теорема доказана.□

Соблюдая условия теоремы мы можем приклеить неограниченное число точек к A без изменения размерности пространства A. Рассмотрим теперь отбрасывание точки, принадлежащей n-мерному пространству A. В этом случае здесь меньше ограничений, поскольку размерность пространства может только уменьшиться.

Теорема 98

*Точка v, имеющая точечный окаем O(v), может быть удалена из n-мерного пространства A, если A-v содержит хотя бы одну n-мерную точку.*

Доказательство

Согласно предыдущему при отбрасывании точки размерность оставшихся в пространстве точек увеличиться не может, то есть пространство или сохраняет свою размерность, или же размерность пространства уменьшается на 1. Поскольку A-v также содержит n-мерную точку, то размерность пространства не меняется. Теорема доказана.□

Более сложная ситуация возникает при установлении или отбрасывании

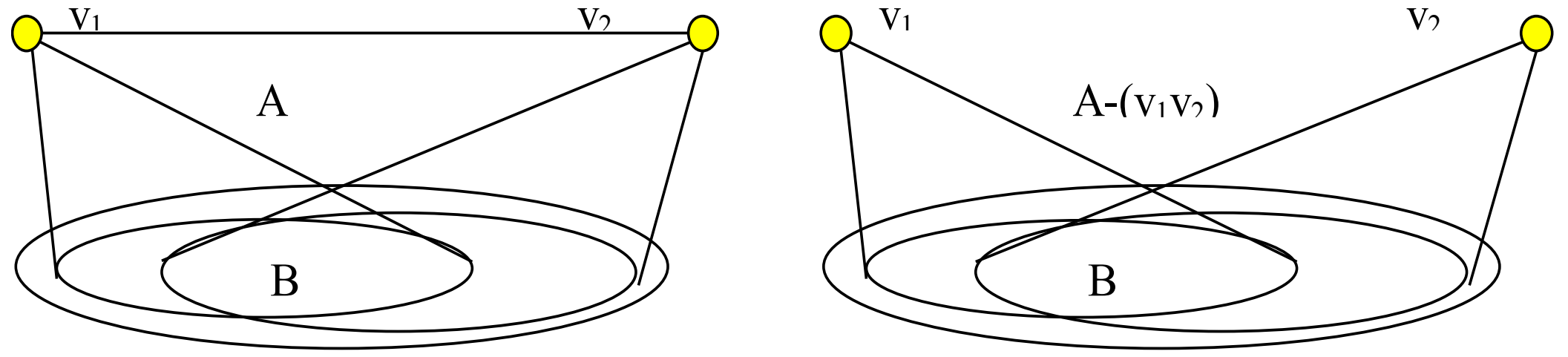

Рис. 107 При отбрасывании или введении связи размерность n-мерного пространства A не меняется, если выполнены условия соответствующих теорем.

связи между точками. Это вызвано тем, что при этом размерность точек может как увеличиваться так и уменьшаться. Кроме того в процесс вовлекаются два типа точек: точки, между которыми устанавливается или отбрасывается связь, а также точки, принадлежащие их общему окаему. Возникает достаточно большое число различных вариантов. Например, пусть отбрасывается связь между точками $v_1$ и $v_2$, и v принадлежит их общему окаему $O(v_1v_2)$. В этом случае размерность точек $v_1$ и $v_2$, может уменьшиться, а размерность точки v может уменьшиться или увеличиться. Если $v_1$ есть единственная n-мерная точка, и ее размерность уменьшается, а v есть (n-1)-мерная точка, и ее размерность увеличивается, то в целом пространство не меняет свою размерность.

Теорема 99

*Пусть A есть n-мерное пространство, B есть общий окаем двух смежных точек $v_1$ и $v_2$, $O(v_1v_2)=B$, и выполняются следующие условия:*

- *B является точечным пространством.*
- *Если точка $v_1$ является n-мерной, то существует некоторое замкнутое нормальное (n-1)-мерное пространство $G^{n-1}$, принадлежащее $O(v_1)$ и не имеющее общих точек с B.*
- *Если точка $v_2$ является n-мерной, то существует некоторое замкнутое нормальное (n-1)-мерное пространство $G^{n-1}$,*

*принадлежащее $O(v_2)$ и не имеющее общих точек с B.*

- *Если точка v, принадлежащая $O(v_1v_2)=B$, является n-мерной, то не существует никакого ее остовного n-мерного шара $U^n{}_0(v)$, лежащего в B и существует ее остовный шар, не лежащий целиком в $v_1 \oplus v_2 \oplus O(v_1v_2)$, или*

*Тогда отбрасывание связи между точками $v_1$ и $v_2$ не меняет размерность каждой точки и пространства вцелом, то есть $A\sim C=A-(v_1v_2)$, раз(A)=раз(C).*

Доказательство

В соответствии с предыдущими теоремами ни одна n-мерная точка пространства A не меняет своей размерности. Следовательно, размерность пространства A при отбрасывании связи не меняется. Теорема доказана.□

Рассмотрим условия, при которых размерность пространства не меняется при введении связи между точками.

Теорема 100

*Пусть A есть n-мерное пространство, B есть общий окаем двух несмежных точек $v_1$ и $v_2$, $O(v_1v_2)=B$, и выполняются следующие условия:*

- *B является точечным пространством.*
- *Если точка $v_1$ является n-мерной, то ее окаем не содержит никакого n-мерного нормального пространства $D^n$ с краем $G^{n-1}$, являющимся замкнутым нормальным (n-1)-мерным пространством, и принадлежащим целиком B.*
- *Если точка $v_2$ является n-мерной, то ее окаем не содержит никакого n-мерного нормального пространства $D^n$ с краем $G^{n-1}$, являющимся замкнутым нормальным (n-1)-мерным пространством, и принадлежащим целиком B.*
- *Если точка v, принадлежащая $O(v_1v_2)=B$, является n-мерной, то не существует n-мерное нормальное пространство $D^n$ со связным краем, принадлежащее O(v), и содержащее точки $v_1$ и $v_2$, такое, что установление связи $(v_1v_2)$ переводит $D^n$ в n-мерное нормальное замкнутое пространство $G^n=D^n \oplus (v_1v_2)$, и существует ее остовный шар $U^n{}_0(v)$ не содержащий по крайней мере одной из точек $v_1$ или $v_2$.*

*Тогда введение связи между точками $v_1$ и $v_2$ не меняет размерность каждой точки и пространства вцелом, то есть $A\sim C=A-(v_1v_2)$, раз(A)=раз(C).*

Доказательство

В соответствии с предыдущими теоремами ни одна n-мерная точка пространства A не меняет своей размерности. Следовательно, размерность пространства A при отбрасывании связи не меняется. Теорема доказана.□

Конечно, последние теоремы налагают слишком строгие ограничения на отбрасывания и введения связи с сохранением размерности пространства. На практике удобнее применять следующее следствие из вышесказанного.

Следствие

Связь $(v_1v_2)$, между точками $v_1$ и $v_2$, n-мерного пространства A, имеющими точечный общий окаем $O(v_1v_2)$, может быть удалена (приклеена) из пространства A, если $A-(v_1v_2)$ содержит хотя бы одну n-мерную точку ($A\oplus(v_1v_2)$ не содержит (n+1)-мерных точек).

Докажем еще несколько теорем о свойствах n-мерных пространств. Ранее мы уже изучали свойства прямой суммы и прямого произведения молекулярных пространств. Рассмотрим, как меняется размерность при таких операциях.

Теорема 101

*Пусть $A^m$ и $B^n$ являются m и n-мерными пространствами соответственно. Если $A^m$ и $B^n$ не содержат замкнутых нормальных n и m-мерных подпространств соответственно, то прямая сумма $A^m\oplus B^n$ этих пространств имеет размерность n+m, раз$(A^m\oplus B^n)=n+m$. Если $A^m$ содержит замкнутое нормальное n-мерное подпространство и/или $B^m$ содержат замкнутое нормальное m-мерное подпространство, то прямая сумма $A^m\oplus B^n$ этих пространств имеет размерность n+m+1, раз$(A^m\oplus B^n)=n+m+1$.*

Доказательство

Докажем, что если $A^m$ и $B^n$ не содержат замкнутых нормальных n и m-мерных подпространств соответственно, то прямая сумма $A^m\oplus B^n$ этих пространств имеет размерность n+m, раз$(A^m\oplus B^n)=n+m$. Пусть $A^m$ с набором точек $V=(v_1, v_2, v_3, \dots vp)$ и $B^n$ с набором точек $R=(r_1, r_2, r_3, \dots r_q)$ не содержат замкнутых нормальных n и m-мерных подпространств соответственно. Рассмотрим точку $v\in A^m$, имеющую размерность m, и точку $u\in B^n$ имеющую размерность n. Это означает, что O(v) содержит нормальное

(m-1)-мерное замкнутое подпространство $G^{v-1}$ и O(u) содержит нормальное (n-1)-мерное замкнутое подпространство $H^{n-1}$. Тогда в $A^m\oplus B^n$ окаем точки v является прямой суммой $O(v)\oplus B^n$ и содержит прямую сумму $G^{v-1}\oplus H^{n-1}$ в качестве подпространства. Так как $G^{v-1}\oplus H^{n-1}$ является нормальным замкнутым (m+n-1)-мерным подпространством, то размерность точки v в

$G^{v-1} \oplus H^{n-1}$ равна (m+n).

Предположим, что размерность некоторой точки v в $A^m \oplus B^n$ есть (m+n+1). Пусть эта точка принадлежит A. Это означает, что существует нормальное (n+m+1)-мерное замкнутое подпространство W, принадлежащее $O(v) \oplus B^n$. Пусть это подпространство состоит из точек $v_1$, $v_2$, $v_3$,.... vp, принадлежащих $A^m$, и всех точек, принадлежащих $B^n$. Тогда $W=G \oplus H$, где G и H являются замкнутыми нормальными подпространствами, принадлежащими $A^m$ и $B^n$ соответственно, причем или размерность G превышает n, или размерность H превышает m. Однако, это противоречит условиям теоремы. Следовательно, наше предположение неверно, и размерность любой точки пространства не может превышать n. Аналогично доказывается, что если $A^m$ содержит замкнутое нормальное n-мерное подпространство и/или $B^m$ содержат замкнутое нормальное m-мерное подпространство, то прямая сумма $A^m \oplus B^n$ этих пространств имеет размерность n+m+1, раз$(A^m \oplus B^n)$=n+m+1. Теорема доказана.□

Теорема 102

> *Пусть $A^m$ и $B^n$ являются m и n-мерными пространствами соответственно. Тогда прямое произведение $A^m \otimes B^n$ этих пространств имеет размерность n+m, раз$(A^m \otimes B^n)=n+m$.*

Доказательство

По сути дела доказательство этой теоремы будет содержиться в доказательстве теоремы о размерности прямого произведения квазинормальных пространств, поскольку размерность точки в прямом произведении определяется локально, через размерность двух точек, определяющих данную.

## Список литературы к главе 8.

# РЕГУЛЯРНЫЕ МОЛЕКУЛЯРНЫЕ ПРОСТРАНСТВА

We consider kvazinormal and block molecular spaces which can be turned into normal spaces by contractible transformations.

Некоторые результаты этой главы изложены в [43].

## *РЕГУЛЯРНЫЕ ПРОСТРАНСТВА*

Молекулярные пространства, помимо всего прочего, являются удобной возможностью представить в наглядном виде многомерные непрерывные поверхности. Например, в непрерывном случае проективная плоскость или бутылка Клейна не могут быть построены в трехмерном пространстве. Для воспроизведения же их молекулярных моделей потребуется лишь набор шариков и стержней, их соединяющих.

Перейдем теперь к более широкому классу молекулярных пространств, который включает в себя предыдущие, нормальные пространства Именно эти пространства будут возникать при создании молекулярных образов непрерывных пространств. В одной из следующих глав мы опишем несколько алгоритмов такого моделирования.

Так как дигитальная топология в некотором смысле экспериментальная наука, то любая модель пространства должна подтверждаться какими-то опытами. Поэтому, прежде чем перейти к определениям, рассмотрим некоторые типичные примеры, позволяющие понять, для чего необходимы и как получаются такие пространства. Предположим, что мы имеем покрытие единичного квадрата $L^2$ небольшими прямоугольниками со сторонами, параллельными сторонам квадрата. Будем считать, что длина наибольшей стороны любого прямоугольника много меньше единицы. Далее, построим нерв этого покрытия и рассмотрим его как молекулярную модель G этого квадрата. Оказывается, что G сохраняет основные свойства квадрата.

- G не содержит трехмерных нормальных точек
- G содержит двумерные нормальные точки
- Точечными отбрасываниями точек G стягивается к нормальному двумерному шару, являющемуся молекулярным образом единичного квадрата

Примеры можно продолжить, обратившись к главе, посвященной компьютерным экспериментам.

Все это приводит к убеждению, что при получении молекулярных образов непрерывных многомерных объектов мы практически всегда получаем МП, отличные от нормальных и должны, если это необходимо, использовать точечные преобразования, чтобы привести эти пространства к нормальной форме. Назовем такие пространства регулярными.

Определение таких пространств основано на результатах экпериментов и может быть осуществлено несколькими способами.

Переходя к обсуждению введения строгих определений поставим условие, что нормальное пространство является как бы остовом или скелетом для регулярного пространства, к которому приклеиваются дополнительные точки так, чтобы топологические характеристики пространства оставались неизменными. Рассмотрим, например, одномерную нормальную окружность, состоящую из 4 точек (Рис. 108). Точечные преобразования позволяют нам приклеивать к ней дополнительные точки и связи, так же как и отклеивать их. Возникает вопрос, какие приклеивания сохраняют пространство как одномерную окружность, а какие меняют ее

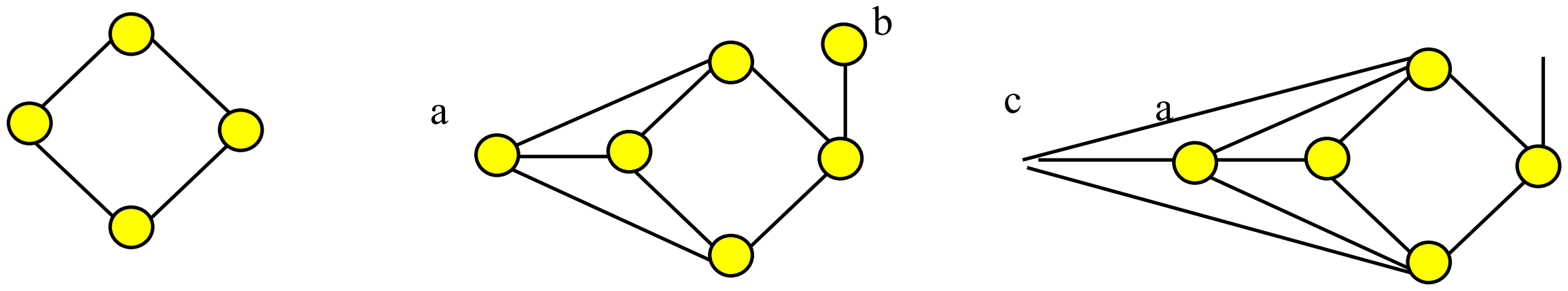

Рис. 108 Приклеивания к окружности точек a и b являются допустимыми, приклеивание точки с недопустимо.

характеристики. Естественно, определение допустимых приклеиваний дается более или менее произвольно, Мы будем приклеивать точки и связи до тех пор, пока размерность точек пространства не превышает 1. Приклеивание, которое увеличивает размерность какой-либо точки пространства до 2 считается недопустимым. На Рис. 108 приклеенные к окружности точки a и b считаются допустимыми, а вот точку c приклеивать уже нельзя, так как при этом точка a становится двумерной, окруженной одномерной окружностью. Если наше объяснение покажется читателю неубедительным, он может придумать свое и ввести свои собственные определения.

<u>Определение 1 регулярного n-мерного замкнутого пространства</u>

Пространство G называется замкнутым регулярным n-мерным, если:

1. G гомотопно замкнутому нормальному n-мерному пространству $H^n$
2. Некоторые точки пространства G являются n-мерными
3. Размерность любой точки пространства G не превышает n
4. Любая точка пространства G является смежной с хотя бы одной n-мерной точкой

Условие 1 определяет топологическую структуру пространства. Условие 2 означает, что пространство, по крайней мере, n-мерно. Условие 3 говорит о том, что пространство не имеет точек высших размерностей. Например, это обеспечивает отсутствие 2-мерных утолщений на молекулярном пространстве, моделирующем 1-мерную прямую. Возможно, что условие 2 являются следствием условия 1. Рассмотрим на простейшем уровне (Рис. 109), какие непрерывные пространства, расположенные в евклидовом непрерывном 3-мерном пространстве, моделируются и не моделируются регулярной одномерной окружностью $S^1_R$. В качестве элементов покрытия выберем небольшие произвольные кубы со сторонами, параллельными координатным осям. В качестве молекулярной модели рассматривается нерв покрытия. Согласно определению $S^1_R$ гомотопна минимальной

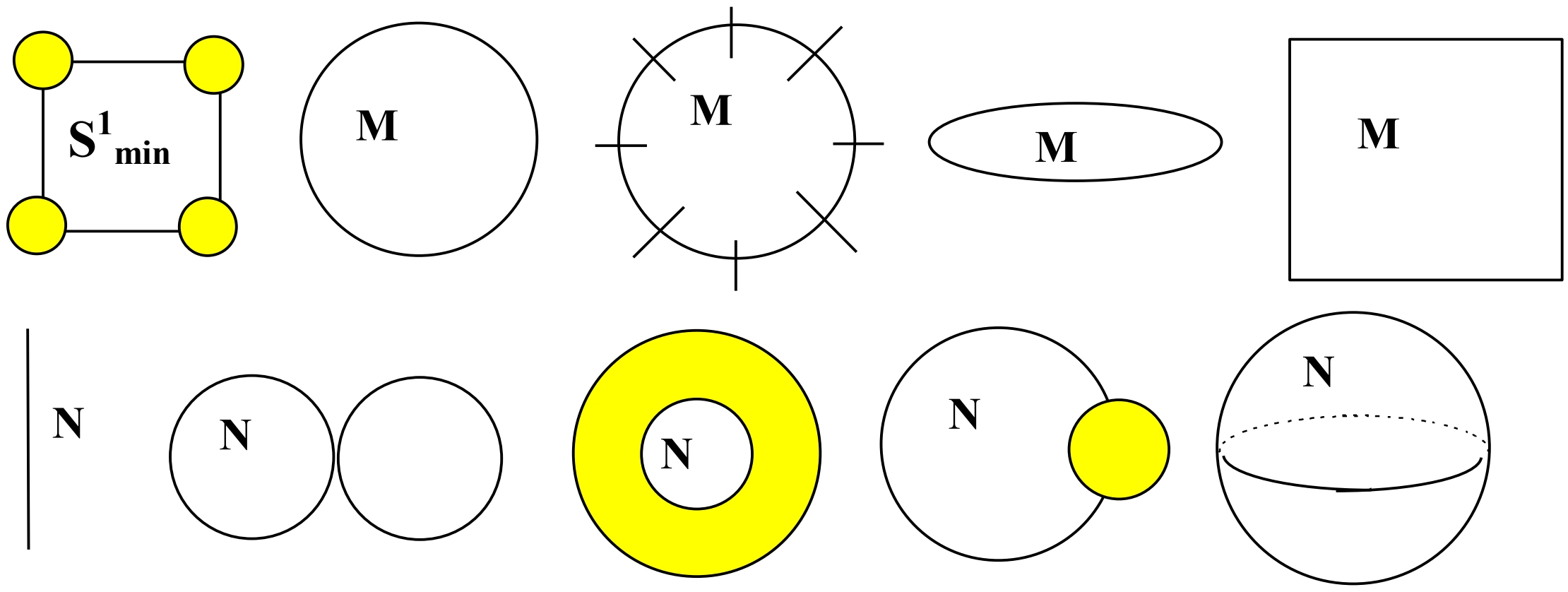

Рис. 109 Пространства M моделируются регулярной окружностью, пространства N не моделируются регулярной окружеостью.

нормальной окружности $S^1_{min}$, изображенной на том же рисунке. Буква M означает, что молекулярная модель данного непрерывного пространства есть $S^1_R$, буква N означает, что молекулярная модель данного непрерывного пространства не есть $S^1_R$.

Введем еще одно определение регулярного замкнутого пространства, являющееся с математической точки зрения более удобным для использования. Вопрос об эквивалентности этих двух определений остается открытым.

<u>Определение 2 замкнутого регулярного пространства</u>

1. Молекулярное пространство G называется регулярным замкнутым регулярным пространством, если
2. G имеет подпространство $H^n$, называемое остовом или скелетом пространства G и являющееся нормальным замкнутым n-мерным пространством
3. G гомотопно $H^n$
4. Размерность любой точки G не превышает n
5. Любая точка из G является смежной хотя бы одной точке из $H^n$.

Перейдем к рассмотреню некоторых простейших видов регулярных пространств.

## *КВАЗИНОРМАЛЬНЫЕ МОЛЕКУЛЯРНЫЕ ПРОСТРАНСТВА*

Для того, чтобы получить молекулярный образ непрерывного пространства, скажем, двумерного тора, мы должны использовать, какой-либо способ дискретизации. Существует несколько таких способов, о них мы будем говорить позднее, которые используются в зависимости от конкретных условий. Например, мы вводим конечное или счетное покрытие непрерывного пространства, и рассматриваем нерв этого покрытия, как его молекулярный образ. При этом довольно сложно подобрать покрытие таким образом, чтобы получившееся молекулярное пространство сразу являлось нормальным пространством. В большинстве случаев такое пространство не является нормальным и приводится к нормальному с помощью точечных преобразований. Поэтому нам необходимо расширить понятие нормального пространства, определив пространства, которые формально не являются нормальными, но близки к

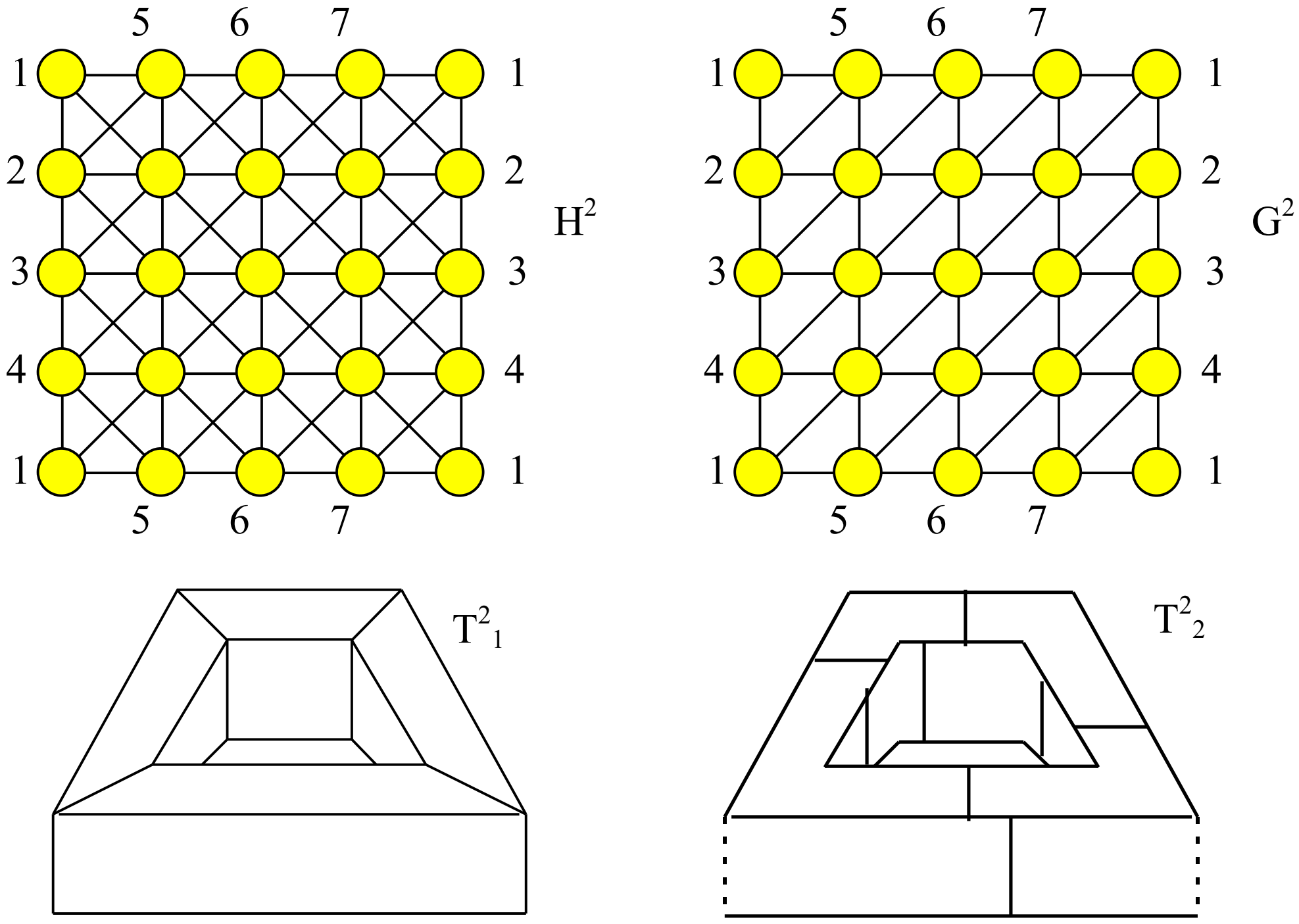

Рис. 110 Квазинормальный (слева) и нормальный (справа) двумерные торы $T^2_1$ и $T^2_2$. В нижней части изображены покрытия торов, нервы которых являются соответствующими молекулярными пространствами.

ним по своим свойствам. В качестве примера рассмотрим покрытия двумерного тора-поверхности бублика, и их молекулярные пространства (*Рис. 110*). Очевидно, что молекулярное пространство $H^2$ не является нормальным. В этом пространстве мы можем удалить точечные связи, как это показано на том же рисунке, и перейти к нормальному пространству $G^2$. Окаем каждой точки пространства $H^2$ содержит нормальную окружность как подпространство, и, кроме того, имеет еще четыре дополнительные точки, не входящие в этот окаем. Такого вида пространства мы назовем квазинормальными.

Определение квазинормального n-мерного замкнутого пространства

Пространство $H^n$ называется квазинормальным n-мерным замкнутым пространством, если выполняются следующие условия:

- Пространство $H^n$ гомотопно нормальному n-мерному замкнутому пространству $G^n$,
- Окаем каждой точки пространства $H^n$ содержит подпространство, гомотопное замкнутому нормальному (n-1)-мерному пространству, но не содержит подпространства, гомотопного замкнутому нормальному n-мерному пространству.

Вторая часть в этом определении является наиболее существенной, поскольку она позволяет определить квазинормальное n-мерное пространство, не являющееся замкнутым.

Определение квазинормального n-мерного бесконечного пространства с бесконечным числом точек

Пространство $H^n$ называется квазинормальным n-мерным пространством с бесконечным числом точек, если выполняются следующие условия:

- Пространство $H^n$ гомотопно нормальному n-мерному пространству $G^n$ с бесконечным числом точек,
- Окаем каждой точки пространства $H^n$ содержит подпространство, гомотопное замкнутому нормальному (n-1)-мерному пространству, но не содержит подпространства, гомотопного замкнутому нормальному n-мерному пространству.

На Рис. 111 изображена квазинормальная одномерная прямая. Окаем каждой точки содержит нуль-мерную нормальную сферу и не содержит никакой одномерной нормальной сферы. Черным цветом обозначены точки нормальной прямой, являющейся подпространством квазинормальной регулярной прямой.

Сразу же возникает вопрос. Пусть окаем каждой точки некоторого пространства $H^n$ на конечном множестве точек содержит подпространство, гомотопное замкнутому нормальному (n-1)-мерному пространству, но не содержит подпространства, гомотопного замкнутому нормальному n-мерному пространству. Всегда ли является это пространство гомотопным замкнутому n-мерному пространству? Эта проблема еще не решена. Поэтому мы ввели первое условие в определение квазинормального n-мерного замкнутого пространства. Введем определение квазинормальной n-мерной точки.

Определение n-мерной квазинормальной точки

Точка v называется n-мерной квазинормальной точкой, если ее окаем содержит подпространство, гомотопное замкнутому нормальному (n-1)-мерному пространству, но не содержит подпространства, гомотопного замкнутому нормальному n-мерному пространству.

Квазинормальные пространства наиболее близки по своим основным свойствам к нормальным пространствам. Большинство теорем о нормальных пространствах легко доказываются для квазинормальных пространств. Легко видеть, что в нуль-мерном случае существуют только нормальные пространства. Отличие между нормальными и квазинормальными пространствами возникает только при размерности 1 и выше. Очевидно, что любое нормальное пространство является также квазинормальным пространством. Выделим это очевидное утверждение в

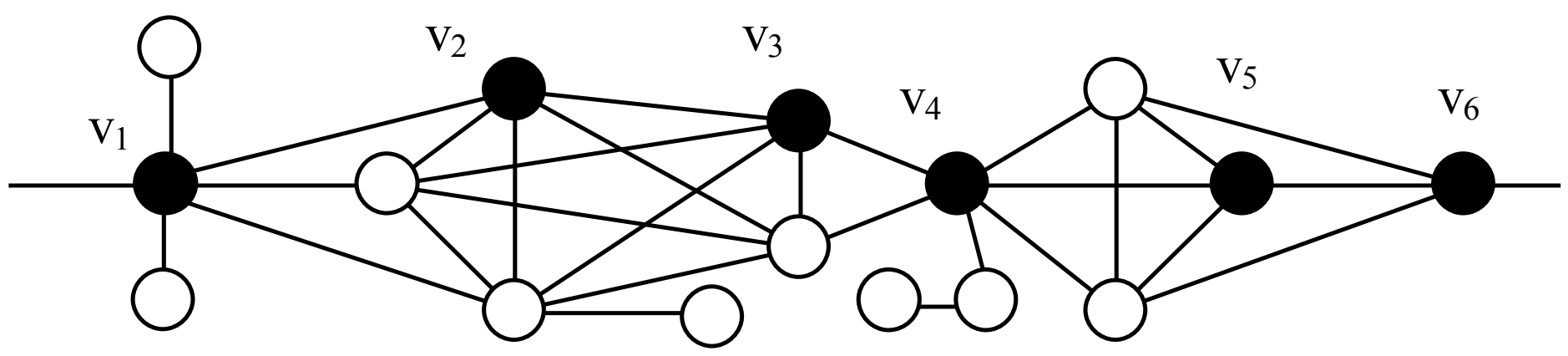

Рис. 111 Квазинормальная одномерная прямая. Окаем каждой точки содержит нуль-мерную нормальную сферу и не содержит никакой одномерной нормальной сферы. Черным цветом обозначены точки нормальной прямой, являющейся подпространством квазинормальной регулярной прямой.

виде теоремы.

Теорема 103

*Любое нормальное пространство есть квазинормальное пространство*

Рассмотрим некоторые свойства квазинормальных пространств.

Теорема 104

*Пусть $G^n$ и $H^m$ есть два квазинормальных замкнутых пространства размерности n и m соответственно. Тогда прямая сумма $G^m \oplus H^m$ есть квазинормальное МП размерности n+m+1.*

Доказательство.

$G^n$ и $H^m$ стягиваются точечными преобразованиями к нормальным замкнутым пространствам $G^n{}_1$ и $H^m{}_1$. Так как эти преобразования являются точечными также в прямой сумме пространств, то $G^n \oplus H^m$ гомотопно $G^n{}_1 \oplus H^m{}_1$. Пусть $v \in G^n$, $O_1(v) \subseteq O(v)$ гомотопен замкнутому нормальному (n-1)-мерному пространству $G^{n-1}{}_2$, $u\in H^m$, $O_1(u)\subseteq O(u)$ гомотопен замкнутому нормальному (n-1)-мерному пространству $H^{m-1}{}_2$. Тогда $O_1(v)|(G^n\oplus H^m)=(O_1(v)|G^n)\oplus H^m$ гомотопен $G^n{}_2\oplus H^m$. Следовательно, окаем точки v в $G^т+H^m$ содержит подпространство, гомотопное нормальному замкнутому (n+m)-мерному пространству $G^n{}_2\oplus H^m$. То же самое справедливо для точки $u\in H^m$. Пусть окаем точки v в $G^n\oplus H^m$ содержит подпространство P, гомотопное нормальному замкнутому (n+m+1)-мерному пространству. Очевидно, что $A\oplus B=P$, $A\subseteq O(v)$ и $B\subseteq H^m$. A не может быть гомотопно нормальному замкнутому n-мерному пространству $F^n$. Тогда $A\oplus B=P$ не может быть гомотопно нормальному замкнутому (n+m+1)-мерному пространству $V^{n+m+1}$. Получили противоречие. Теорема доказана. □

Прямое произведение квазинормальных пространств также обладает аналогичными свойствами. Отметим, что в этом случае уже нет необходимости ограничиваться замкнутыми нормальными пространствами.

Теорема 105

*Пусть $G^n$ и $H^m$ есть два квазинормальных пространства размерности n и m соответственно. Тогда прямое произведение $G^m\otimes H^m$ есть квазинормальное МП размерности n+m.*

Доказательство.

$G^n$ и $H^m$ гомотопны $G^n{}_1$ и $H^m{}_1$. Следовательно, $G^n\otimes H^m$ гомотопно $G^n{}_1\otimes H^m{}_1$. Так как $O(v\otimes u)=[(v\otimes O(u))\oplus(O(v)\otimes u)]\cup(O(v)\otimes O(u))$ гомотопен $O(u)\oplus O(v)$, то, в соответствии с предыдущей теоремой, $O(v\otimes u)$ содержит подпространство, гомотопное замкнутому нормальному (n+m-1)-мерному пространству, и не содержит никакого подпространства, гомотопного замкнутому нормальному (n+m-1)-мерному пространству. Теорема доказана. Теорема доказана. □

Определим, теперь, квазинормальное пространство с краем. Сделать это несколько сложнее, чем в случае нормальных пространств, поскольку отбрасывание точки еще не означает, что квазинормальное пространство может перейти в пространство другого топологического вида. Из Рис. 111 видно, что при отбрасывании точки $v_2$ пространство остается квазинормальным пространством того же типа.

Определение квазинормального пространства со связным краем

Пространство $H^n$ называется квазинормальным n-мерным пространством с связным краем B, если выполняются следующие условия:

- после приклеивания точки v к $H^n$ по подпространству B (O(v)=B), полученное пространство $G^n=H^n \cup v$ является квазинормальным n-мерным пространством.
- Пространство $H^n$ гомотопно некоторому нормальному n-мерному замкнутому пространству с краем.

Исходя из этого определения край B является пространством, которое содержит подпространство, гомотопное замкнутому нормальному (n-1)-мерному пространству, но не содержит подпространства, гомотопного замкнутому нормальному n-мерному пространству. Кроме того, видимо, пространство B само гомотопно замкнутому нормальному (n-1)-мерному пространству. Эта проблема еще не решена.

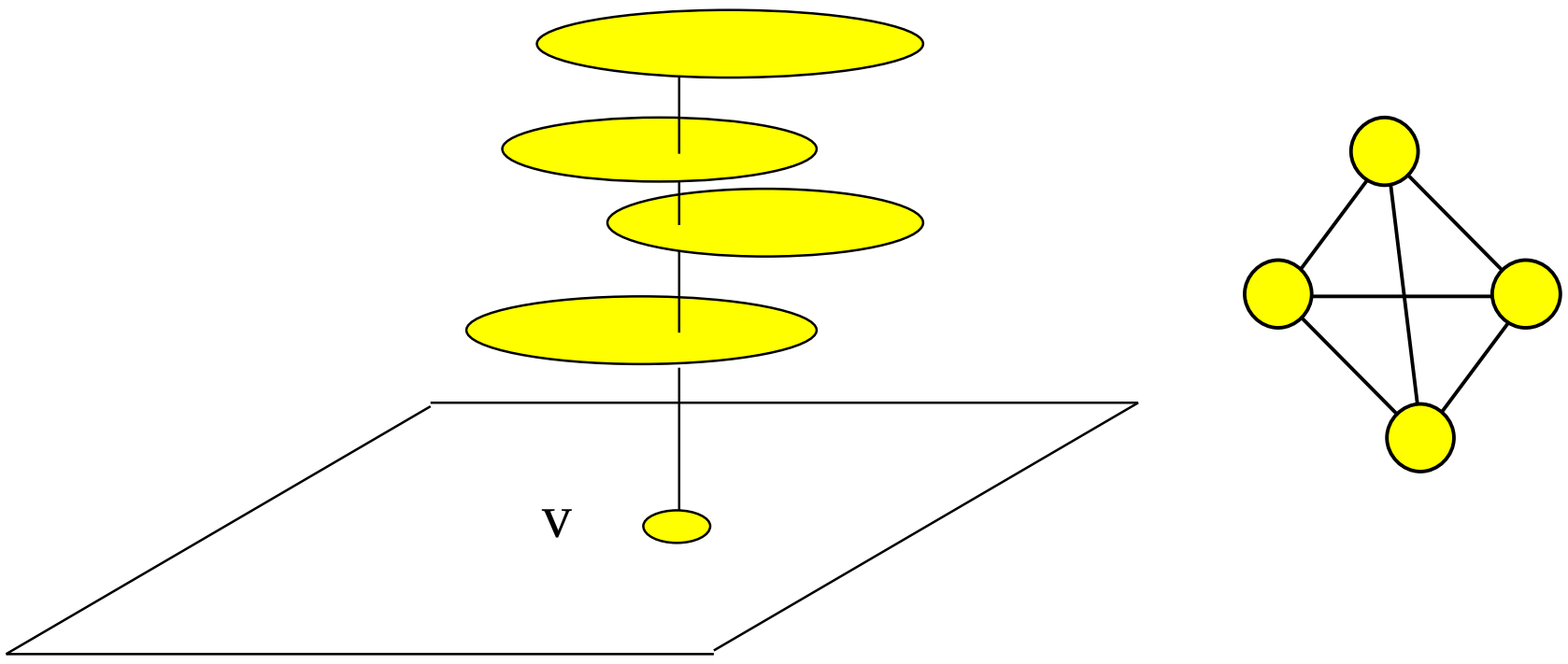

Рис. 112. Точка v принадлежит 4-м открытым дискам. Нерв покрытия является связкой из 4-х точек.

З а д а ч а

Доказать или опровергнуть следующее утверждение: пусть $H^n$ есть квазинормальное пространство размерности n с краем B. Тогда B гомотопен замкнутому нормальному (n-1)-мерному пространству.

## *МОЛЕКУЛЯРНЫЕ БЛОК-ПРОСТРАНСТВА*

Рассмотрим свойства блок-пространств, то есть таких пространств, в которых точки заменены некоторыми пространствами. Перейдем определению и свойствам блок-пространств, которые образуют более широкий класс пространств, чем нормальные, и включают в себя нормальные пространства, как частный случай.

Блок пространство есть не что иное как молекулярное пространство, в котором точки заменены на некоторые пространства. Прежде, чем перейти к определениям, рассмотрим некоторые наглядно-геометрические построения, позволяющие понять, для чего необходимы и как получаются такие пространства. Пусть точка v принадлежит плоскости, которая покрыта некоторым семейством открытых дисков (*Рис. 112*). Поставим в соответствие каждому диску точку молекулярного пространства.

Соединим две точки молекулярного пространства, если соответствующие им диски имеют непустое пересечение. Очевидно, что все точки МП будут соединены между собой.

Этот пример дает основание предположить, что точка непрерывного пространства может представляться связкой из n точек в молекулярном пространстве. Он объясняет путь наших дальнейших действий. В блок пространстве роль точки играет связка точек, то есть несколько точек с полным набором связей между собой. Ключевым, также, является вопрос, как связаны точки между собой в блок пространстве, чтобы сохранилось соответствие нормальным моделям.

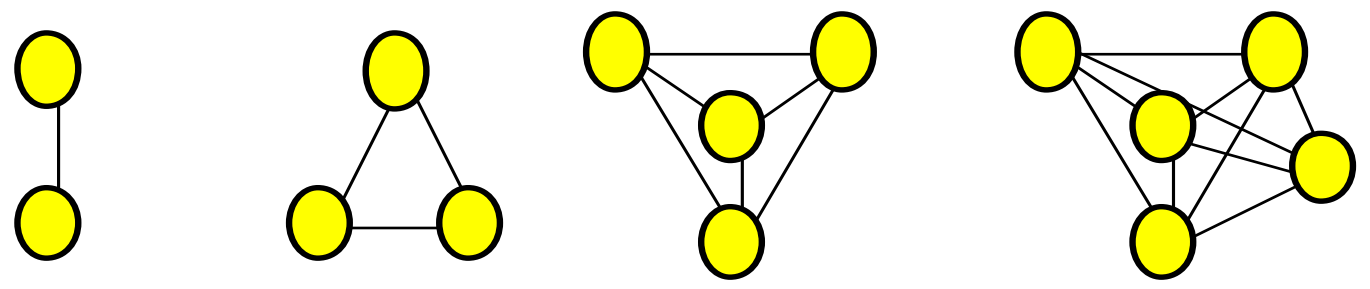

Рис. 113 Блок-точки, являющиеся связками на 2, 3, 4 и 5-и точках.

Рассмотрим еще один пример (*Рис. 114*). Пусть имеется покрытие единичного квадрата пятью базовыми прямоугольниками. Нерв этого покрытия есть минимальный двумерный молекулярный нормальный шар, то есть точка, окаем которой является минимальной окружностью. Добавим к этому покрытию другие элементы, являющиеся небольшими прямоугольниками. При этом каждый добавленный элемент имеет непустое пересечение, по крайней мере, с одним из базовых прямоугольников, кроме того, размерность нерва покрытия всегда должна

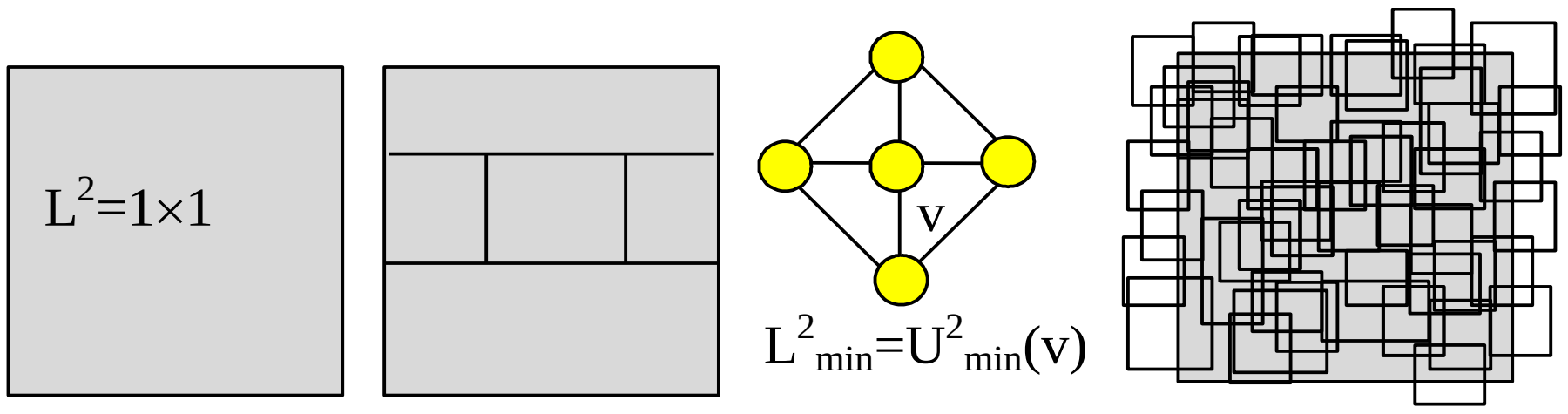

Рис. 114 Единичный квадрат $L^2$ покрыт прямоугольниками. Нерв покрытия является молекулярной моделью квадрата (не изображен). Точечными преобразованиями с сохранением размерности он стягивается к минимальному молекулярному квадрату $L^2_{min}=U^2_{min}(v)$, являющемуся нервом покрытия правого квадрата.

быть равна двум. Геометрически это означает, что к базовому молекулярному пространству приклеиваются по определенным правилам точки при условии, что размерность и другие топологические и геометрические свойства пространства не меняются. Было доказано, что эйлерова характеристика и группы гомологий молекулярного пространства не меняются при точечных преобразованиях пространств. Необходимо на точечные преобразования наложить такие условия, которые позволили бы сохранить размерность пространства.

Исходя из этого определим блок пространство как пространство, состоящее из некоторого базового пространства, в котором точки заменены на некоторые другие пространства, и связи также задаются по определенным. законам, не меняющим основных характеристик пространства. Один вид блок пространств был нами уже подробно разобран. Прямое произведение двух молекулярных пространств G⊗H есть пространство G (или H), в котором каждая точка заменена на пространство H при соответствующем распределении связей. В этом разделе мы рассмотрим другой вид блок пространств, в котором каждая точка

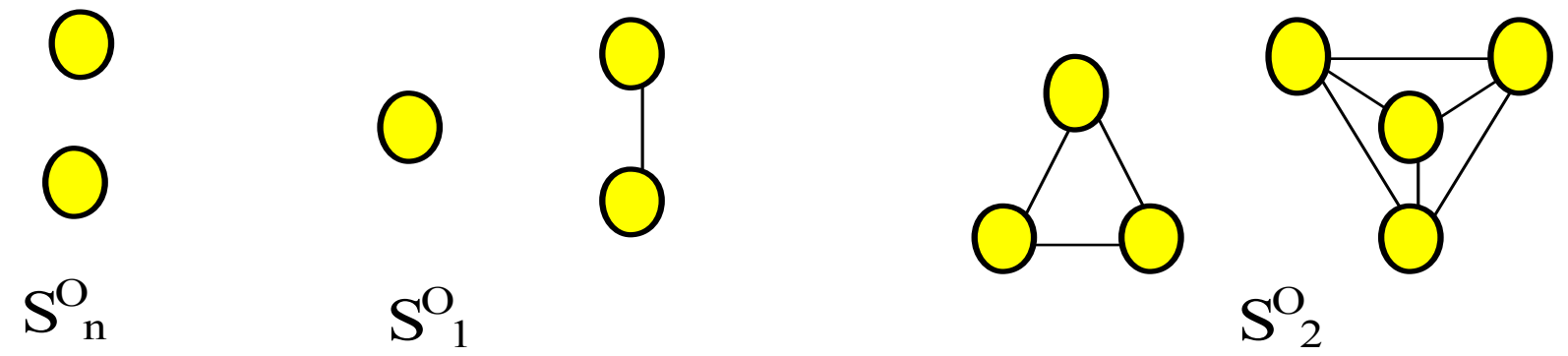

Рис. 115 $S^0_n$ является нормальной нуль-мерной сферой. $S^0_1$ и $S^0_1$ являются регулярными нуль-мерными сферами.

заменяется связкой. Если каждая точка молекулярного пространства G заменяется одной и той же связкой K, то получается прямое произведение G⊗K. Свойства такого пространства изучены достаточно хорошо. Здесь мы рассмотрим пространства, в которых каждая точка заменяется связкой с различным числом точек. Такие пространства более соответствуют реальным дигитальным моделям, получаемым при дискретизации непрерывных пространств.

Определение 0-мерной блок-точки

Связка K(n) из n попарно смежных точек называется 0-мерной блок-точкой

На *Рис. 113* изображено несколько блок-точек. Такое определение хорошо согласовывается с вышеописанным эвристическим подходом, где молекулярное пространство рассматривается, в частности, как нерв покрытия топологического пространства.

Определение 0-мерной регулярной сферы

Регулярной 0-мерной сферой $S^0$ является МП, состоящее из двух изолированных блок-точек. Множество всех регулярных 0-мерных пространств содержит только 0-мерную регулярную сферу (*Рис. 115*).

Определение блок пространства

Пусть A является молекулярным пространством на множестве точек

$(v_1, v_2, \dots v_p \dots)$. Пространство BA называется блок-пространством пространства A если:

- каждая точка $v_k$ заменяется на блок точку $V_k(m_k)$, являющуюся связкой из $m_k$ точек,
- из условия, что $v_s$ и $v_t$ смежны в A следует, что любая точка $v_{kp} \in V_k(m_k)$ смежна любой точке $v_{tr} \in V_t(m_t)$.

Будем говорить, что блок-точки $V_k(m_k)$ и $V_t(m_t)$ смежны в BA, если $v_s$ и $v_t$ смежны в A. Иными словами если в молекулярном пространстве заменить точки на блок-точки и добавить соответствующие связи, то получится регулярное пространство. Любая точка, принадлежащая $V_k(m_k)$, связана со всеми точками из $V_s(m_s)$. На Рис. 116 изображено регулярное пространство $L^1_{reg}$, являющееся одномерной прямой из блок-точек $V_1(1)$, $V_2(3)$, $V_3(2)$, $V_4(1)$, $V_5(3)$, $V_6(1)$.... Выбранные точки $v_1$, $v_2$, $v_3$, $v_4$, $v_5$, $v_6$...обозначены черным цветом и составляют нормальную прямую.

Определение нормального n-мерного блок пространства

Если $G^n$ является нормальным n-мерным пространством, то пространство $BG^n$ называется нормальным n-мерным блок-пространством.

Рассмотрим некоторые свойства блок-пространств. Рассмотрим структуру окаема точки, принадлежащей некоторой блок-точке.

Теорема 106

*Пусть точка $v_1$ принадлежит блок-точке $V_1(m_1)$нормального n-мерного блок-пространства $H^n$. Тогда окаем произвольной точки $v_1$ определяется выражением $O(v_1)=(V_1(m_1)-v_1)\oplus O(V_1(m_1))$, где $(V_1(m_1)-v_1)$ и $O(V_1(m_1))$ являются связкой из $(m_1-1)$ точек и окаемом блок-точки $V_1(m_1)$. При этом $O(V_1(m_1))$ есть замкнутое (n-1)-мерное нормальное блок-пространство.*

Доказательство.

По определению окаем $O(V_1)$ произвольной блок-точки пространства $V_1(m_1)$ является молекулярным замкнутым регулярным (n-1)-мерным пространством. Точка $v_1$ из $V_1(m_1)$ имеет связь со всеми точками из $O(V_1)$, а также с остальными точками из $V_1(m_1)$, то есть с точками $V_1(m_1)-v_1$. В свою очередь точки из $(V_1(m_1)-v_1)$ смежны всем точкам из $O(V_1)$. Следовательно $O(v_1)=(V_1(m_1)-v_1)\oplus O(V_1(m_1))$. Теорема доказана. □

Теорема 107

*Пусть BA есть блок-пространство. Тогда BA приводится к A точечными отбрасываниями точек.*

Доказательство.

Окаем точки $v_k$ есть нормальное замкнутое (n-1)-мерное пространство $O(v_k)$. Из каждой блок-точки $V_k(m_k)$ отбрасываем все точки кроме одной точечными отбрасываниями. Получаем пространство A. Теорема доказана. □

Теорема 108

*Пусть $BG^n$ есть нормальное n-мерное блок-пространство. Тогда окаем любой точки пространства содержит нормальное замкнутое (n-1)-мерное пространство и не содержит нормальное замкнутое n-мерное пространство.*

Доказательство.

Окаем точки $v_k$ есть нормальное замкнутое (n-1)-мерное пространство $O(v_k)$. Пусть точка $v_{kr}$ принадлежит блок-точке $V_k(m_k)$. Выберем по одной

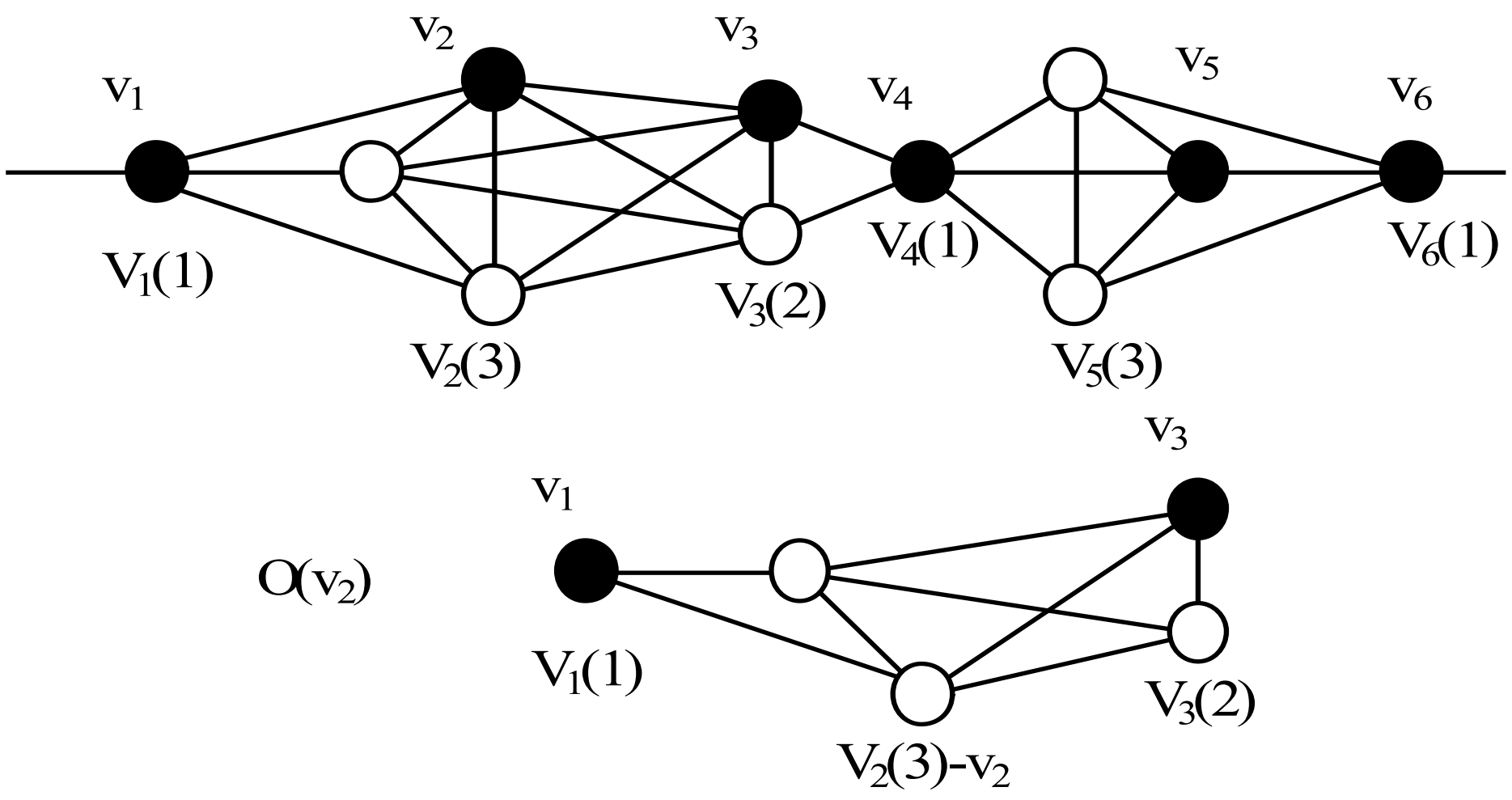

Рис. 116 Регулярная одномерная прямая состоит из блок-точек $V_1(1)$, $V_2(3)$, $V_3(2)$, $V_4(1)$, $V_5(3)$, $V_6(1)$..... Черным цветом обозначены точки нормальной прямой, являющейся подпространством регулярной прямой. Окаем точки $v_2$ есть прямая сумма, $O(v_2)=(V_2(3)-v_2)\oplus O(V_2(3))$..

точке $v_{st}$ из каждой блок-точки $V_s(m_s)$, смежной с $V_k(m_k)$. Очевидно, что выбранное множество точек есть нормальное замкнутое (n-1)-мерное пространство, изоморфное $O(v_k)$. Выберем подпространство, состоящее из блок-точек, смежных с блок-точкой $V_k(m_k)$. Предположим, что это пространство содержит нормальное замкнутое n-мерное подпространство $X^n$. Оно должно содержать некоторые точки, принадлежащие одной и той же блок-точке. Тогда одна из этих точек имеет точечный окаем. Это

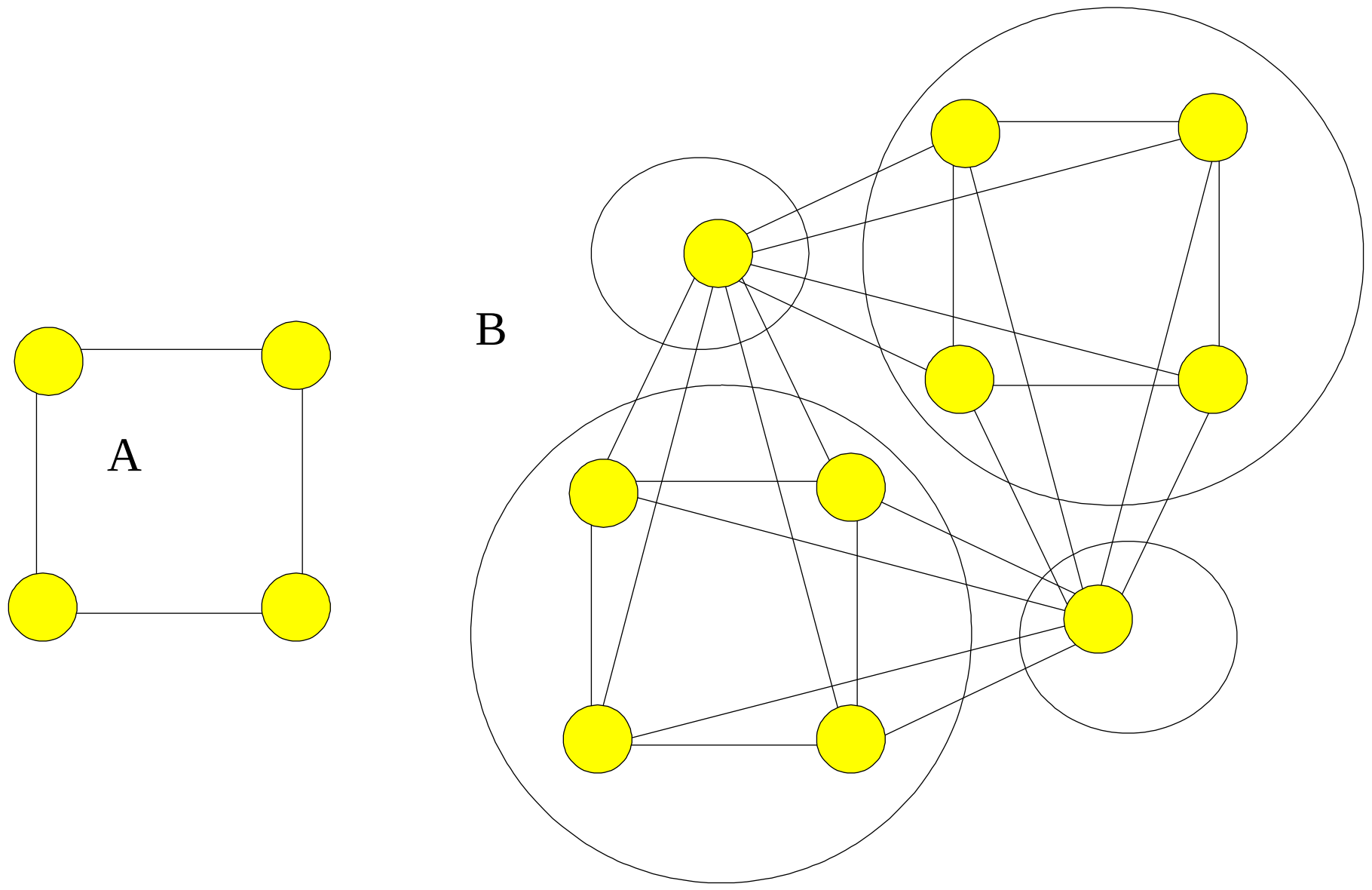

Рис. 117 Нормальная окружность A и блок-пространство B.

противоречит определению нормального замкнутого n-мерного пространства. Таким образом, окаем любой точки пространства не содержит нормальное замкнутое n-мерное пространство. Теорема доказана. □

Следующие теоремы очевидны.

Теорема 109

*Прямая сумма $A_1 \oplus A_2$ блок-пространств $A_1$ и $A_2$, есть блок-пространство*

Теорема 110

*Прямое произведение $A_1 \otimes A_2$ блок-пространств $A_1$ и $A_2$, есть блок-пространство*

Определение остова (нормального n-мерного) блок-пространства

(Нормальное замкнутое n-мерное) пространство $G^n$, являющееся подпространством (n-мерного нормального) блок-пространства $BG^n$ называется (нормальным) остовом или скелетом (n-мерного нормального) блок-пространства $BG^n$.

Очевидно, что блок-пространство может содержать несколько замкнутых нормальных n-мерных пространств, все из которых гомотопны одно другому.

Теорема 111

*Если блок-пространство имеет несколько квазинормальных остовов, то все они гомотопны один другому.*

Доказательство.

Согласно определению любой остов гомотопен исходному блок-пространству. Следовательно, все остовы гомотопны один другому. Теорема доказана.□

Вообще говоря, понятие блок-пространства можно расширить, если в качестве бло-точек брать не связки, а произвольные пространства. Например, в некотором нормальном пространстве часть точек заменить на связки, а оставшиеся точки-на сферы. При этом возникает некоторая внешняя и внутренняя топологии. Внешнее блок- пространство, состоящее из блок-точек как целого, остается нормальным пространством и представляет внешнюю структуру. При этом каждая блок-точка, в свою очередь, состоит из точек, является пространством и имеет внутреннюю топологию. Кроме того, связи между различными блок-точками задают связи между пространствами, составляющими блок-точки, и, следовательно, определяют связь между внутренней структурой каждой блок-точки и внешней общей структурой. Этот путь представляется достаточно интересным и может быть использован в многочисленных физических теориях, микромира, или же в биологических и химических построениях. На Рис. 117 мы показываем такое блок-пространство, в котором внешее пространство, состоящее из блок-точек, является окружностью из 4-х точек, две блок-точки являются обычными точками, две другие блок точки есть нормальные минимальные окружности

Список литературы к главе 9.

# ЭЙЛЕРОВА ХАРАКТЕРИСТИКА МОЛЕКУЛЯРНОГО ПРОСТРАНСТВА

In this chapter we define the Euler characteristic of a molecular space and prove that contractible transformations do not change this characteristic. We study some properties of the Euler characteristic and show by examples that the Euler characteristic of molecular spaces is the same as of their continuous counterparts. The connection is found between the Euler characteristic of a molecular space and the Euler characteristics of rims of its vertices. We apply these results to tilling of two-dimensional closed surfaces: a sphere, torus and projective plane. We investigate in details tilling by one and two different elements. There are a number of theorems describing properties of tilling in general by using the Euler characteristic. For example a projective plane can not have maps containing only quadrangles and pentagons or consisting only of polygons with the number of vertices over 5. Any triangulation of a projective plane by polygons (except triangles) must contain quadrangles and pentagons as well as polygons with the number of vertices more than 5. The definitions of the least and minimal molecular spaces are given and some properties of these MПs are studied.

Полученные в этой главе результаты частично изложены в работах [16,24,25,26,28,43].

## *ОПРЕДЕЛЕНИЕ ЭЙЛЕРОВОЙ ХАРАКТЕРИСТИКИ МОЛЕКУЛЯРНОГО ПРОСТРАНСТВА. ТОЧЕЧНЫЕ ПРЕОБРАЗОВАНИЯ НЕ МЕНЯЮТ ЭЙЛЕРОВУ ХАРАКТЕРИСТИКУ МОЛЕКУЛЯРНОГО ПРОСТРАНСТВА*

Результаты этой главы могут быть выражены в двух предложениях: эйлерова характеристика точечного пространства равна 1, и эйлерова характеристика не меняется при точечных преобразованиях.

Эйлерова характеристика пространства является некоторой числовой функцией пространства как целого. В классической топологии она определяется несколькими способами. Для нас представляет интерес то, что эта функция сохраняет свое значение на топологически эквивалентных поверхностях. Это, например, означает, что как бы мы не изгибали или растягивали двумерную (n-мерную) поверхность, значение эйлеровой характеристики на ней меняться не будет. Поэтому эта характеристика может использоваться для выявления топологического сходства или различия между поверхностями. Важное свойство эйлеровой характеристики в том, что она не меняется при точечных преобразованиях молекулярных пространств. Эта инвариантность эйлеровой характеристики является основным результатом этого раздела. Этот результат свидетельствует о том, что точечные преобразования в теории молекулярных пространств играют ту же роль, что и соответствующие преобразования в классической топологии. Нам предварительно придется определить эйлерову характеристику на пространстве, так как наше определение отличается от стандартного. Тем не менее мы сохраняем название, так как наше определение обобщает традиционное. Известно стандартное определение эйлеровой функции f на графе [22].

$$f = n_1 - n_2 + n_3 \quad (1)$$

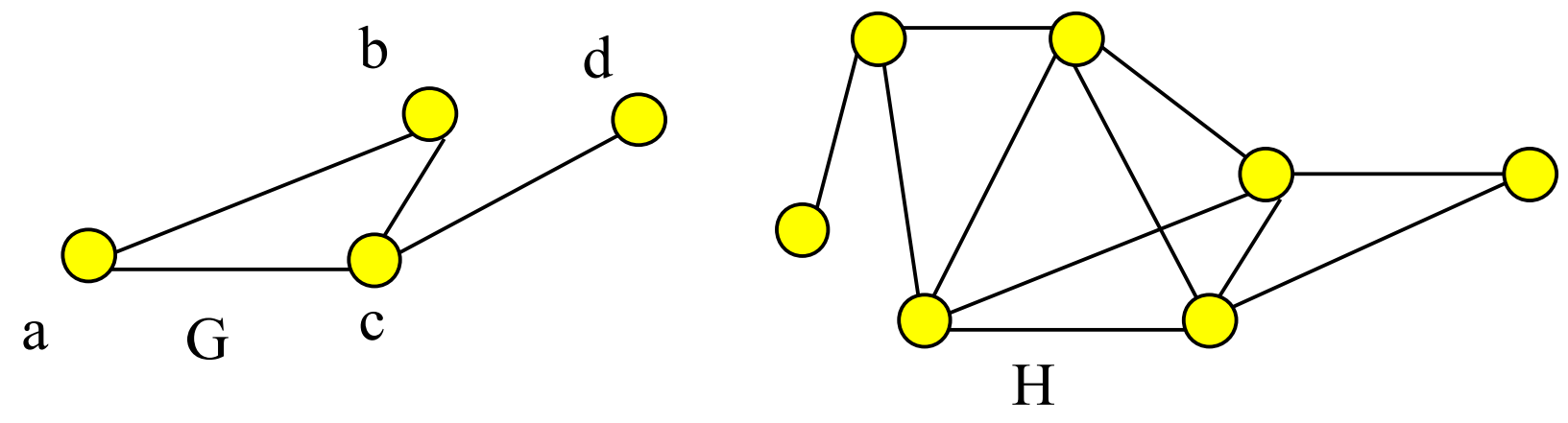

Рис. 118 Эйлерова характеристика пространства G равна 1, эилерова характеристика пространства H равна 1. По формуле Эйлера для графа эйлерова характеристика пространства H равна 2.

где $n_1$, $n_2$, и $n_3$-числа вершин, ребер и треугольников в графе.

Наше определение отличается от данного, оно напоминает определение характеристики Эйлера-Пуанкаре, сформулированное для симплициального комплекса [1,54]. Однако, различие между этими двумя определениями такое же, как различие между симплициальным комплексом и пространством. Поясним это примером.

На Рис. 118 изображена фигура G, которая может быть пространством, а может изображать симплициальный комплекс. Как пространство (или граф), G состоит из одного треугольника (abc) и одного ребра (cd), не входящего в этот треугольник. Рассматривая эту же фигуру как симплициальный комплекс, мы имеем несколько возможных вариантов.

- G есть двумерный симплициальный комплекс, определяемый двумерным симплексом (abc) и одномерным (cd).
- G есть одномерный симплициальный комплекс, определяемый одномерными симлексами (ab), (ac), (bc) и (cd).

Как пространство G имеет одну одномерную точку, остальные точки нульмерны. Этот пример показывает, что пространство не определяет однозначно симплициальный комплекс. Перейдем теперь к определениям.

<u>Определение функционального вектора пространства.</u>

Пусть G и $n_p$ есть конечное пространство и число его связок K(p) на p-вершинах. Тогда набор $(n_1,n_2,....n_s)$ называется функциональным вектором пространства,

$$f(G)=(n_1,n_2,....n_s). \quad (2)$$

В литературе [20,21] в применении к графам этот вектор часто обозначается f(G). На Рис. 118 пространство G имеет 4 вершины, 4 ребра и один треугольник. Следовательно f(G)=(4,4,1). Пространство H на том же рисунке определяется 7 вершинами, 11 ребрами, 6 треугольниками и 1 тетраэдром. Следовательно,

f(H)=(7,11,6,1).

Дадим теперь определение эйлеровой характеристики в применении к молекулярному пространству.

Определение эйлеровой характеристики пространства.

Пусть пространство G имеет функциональный вектор $f(G) = (n_1, n_2,.... n_s)$. Тогда его эйлеровой характеристикой называется функция F(G), определяемая выражением

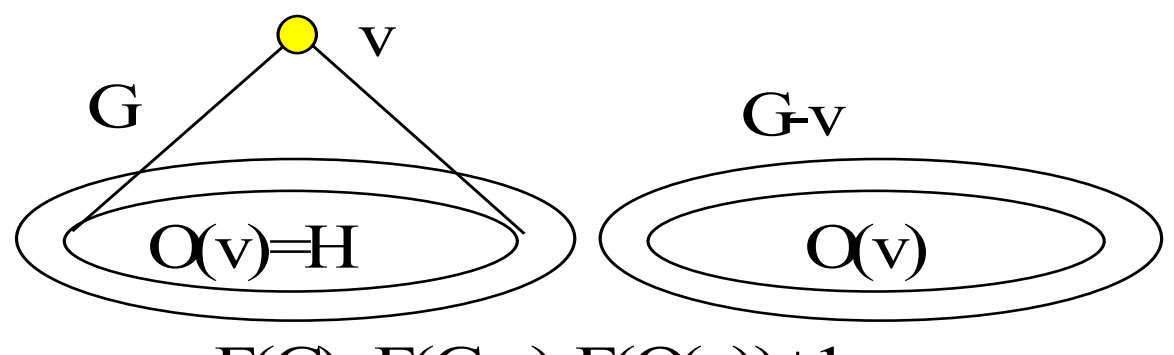

Рис. 119 Показана связь эйлеровой характеристики пространств G и G-v.

$$F(G) = \sum_{k=1}^{s} (-1)^{k+1} n_k . \quad (3)$$

Рассмотрим примеры расчета эйлеровой характеристики пространства H по формуле (1) и нашей формуле (3).

По (1) F(H) = 7-11+6 = 2, по (3) F(H) = 7-11+6-1 = 1.

Легко видеть, что величины различны, и это различие является существенным, о чем говорят результаты этого раздела.

Теорема 112

*Эйлерова характеристика любого конуса v⊕G равна 1,*

$$F(v\oplus G)=1.$$

Д о к а з а т е л ь с т в о

Пусть функциональный вектор пространства G есть $f(G)=(n_1, n_2,... n_s)$. Очевидно, что функциональный вектор пространства v⊕G будет $f(v\oplus G)=(n_1+1, n_2+n_1, n_3+n_2... n_s+n_{s-1}, n_s)$. Тогда $F(v\oplus G)=(n_1+1)-(n_2+n_1)+(n_3+n_2)-...(-1)^{s+1}(n_s+n_{s-1})+(-1)^{k+2}n_s$. После раскрытия скобок получаем 1. $F(v\oplus G)=(n_1+1)-(n_2+n_1)+(n_3+n_2)-...(-1)^{s+1}(n_s+n_{s-1})+(-1)^{k+2}n_s=1$. Теорема доказана.□

Точечные преобразования состоят из точечных отбрасываний и приклеиваний связи и точки. Поэтому нам необходимо определить изменение эйлеровой характеристики при отбрасываниях и приклеиваниях точек и связей.

Теорема 113

*При отбрасывании произвольной точки v из пространства G эйлерова характеристика пространства G-v определяется выражением*

$$F(G)=F(G-v)-F(O(v))+1,$$

*где O(v) есть окаем точки v (Рис. 119).*

Д о к а з а т е л ь с т в о

Пусть функциональный вектор пространства G-v определяется выражением $f(G-v)=(n_1,n_2,n_3,...n_s)$, функциональный вектор пространства O(v) определяется выражением $f(O(v))=(m_1,m_2,m_3,...m_s)$. Легко видеть, что функциональный вектор пространства G есть $f(G)=(n_1+1,n_2+m_1, n_3+m_2,....n_s+m_{s-1},m_s)$. Тогда эйлерова характеристика $F(G)=(n_1+1)-(n_2+m_1)+...+(-1)^{s+1}(n_s+m_{s-1})+(-1)^{s+2}(m_s)$. После простых преобразований мы получаем, что характеристика $F(G)=(n_1-n_2+....+(-1)^{s+1}n_s)-(m_1-m_2+.....+(-1)^s m_{s-1}+(-1)^{s+1}m_s)+1=F(G-v)-F(O(v))+1$. Теорема доказана.□

Теорема 114

*При отбрасывании произвольной связи $v_1v_2$ из пространства G эйлерова характеристика пространства $G-v_1v_2$ определяется выражением*

$$F(G)=F(G-v_1v_2)+F(O(v_1v_2))-1,$$

*где $O(v_1v_2)$ есть общий окаем точек $v_1$ и $v_2$. (Рис. 120).*

Д о к а з а т е л ь с т в о

Пусть функциональный вектор пространства G-v определяется выражением $f(G-v_1v_2)=(n_1, n_2, n_3,...n_s)$, функциональный вектор пространства $O(v_1v_2)$ определяется выражением $f(O((v_1v_2))=(m_1,m_2,m_3,...m_s)$. Легко видеть, что функциональный вектор пространства G есть $f(G)=(n_1,n_2+1,n_3+m_1,n_4+m_2,...n_s+m_{s-2},m_{s-1},m_s)$. Тогда эйлерова характеристика $F(G)=(n_1)-(n_2+1)+(n_3+m_1),...+(-1)^{s+1}(n_s+m_{s-2})+(-1)^{s+2}m_{s-1}+(-1)^{s+3}m_s$. После простых преобразований мы получаем, что характеристика $F(G)=(n_1-n_2+....+(-1)^{s+1}n_s)+(m_1-m_2+.....+(-1)^s m_{s-1}+(-1)^{s+1}m_s)-1=F(G-v)+F(O(v_1v_2))-1$. Теорема доказана.□

Эти две теоремы показывают связь между эйлеровыми характеристиками двух пространств, полученных одно из другого при приклеивании и отбрасывании произвольной точки или связи. Это необходимо для доказательства основной теоремы об инвариантности эйлеровой характеристики при точечных преобразованиях. Докажем еще одну из двух ключевых теорем этого раздела об эйлеровой характеристике точечного пространства.

Теорема 115

*Эйлерова характеристика уединенной точки равна 1.*

*Эйлерова характеристика точечного пространства равна 1.*

Доказательство

Доказательство по индукции. Очевидно, что эйлерова характеристика уединенной точки равна 1. Пусть теорема справедлива для всех пространств с числом вершин $k<n$, числом связей $m<p$. Пусть пространство

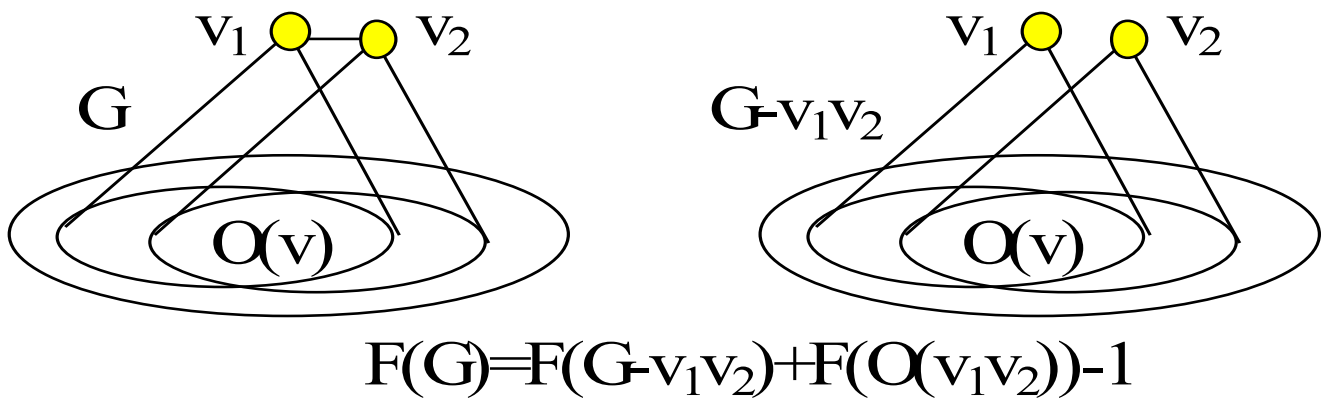

Рис. 120 Показана связь эйлеровой характеристики пространств G и G-$v_1v_2$.

G имеет (n-1) вершин и (m-1) связей. Приклеим к нему точку v по точечному подпространству H, O(v)=H. Обозначим полученное пространство $G_1$. Тогда, согласно предыдущей формуле, $F(G_1)=F(G_1\text{-}v)\text{-}F(O(v))+1$, где $G_1\text{-}v=G$, O(v)=H. Согласно индукции $F(G_1\text{-}v)=F(O(v)=1$. Следовательно, $F(G_1)=1$. Приклеим к G некоторую точечную связь $v_1v_2$, где $O(v_1v_2)$ есть точечный общий окаем, и обозначим полученное пространство $G_1$. Тогда $F(G_1)=F(G_1\text{-}v_1v_2)+F(O(v_1v_2))\text{-}1$. Согласно индукции $F(G_1\text{-}v_1v_2)=F(O(v_1v_2))=1$. Следовательно, $F(G_1)=1$. Теорема доказана.□

Докажем вторую ключевую теорему об инвариантности эйлеровой характеристики при точечных преобразованиях произвольных пространств.

Теорема 116

*Точечные преобразования не меняют эйлерову характеристику пространств.*

Доказательство

Доказательство по индукции. Как и в предыдущем доказательстве $F(G_1)=F(G_1\text{-}v)\text{-}F(O(v))+1$. Согласно индукции $F(O(v)=1$. Следовательно, $F(G_1)=F(G_1\text{-}v)$. Кроме того, $F(G_1)=F(G_1\text{-}v_1v_2)+F(O(v_1v_2))\text{-}1$. Согласно $F(O(v_1v_2))=1$. Следовательно, $F(G_1)=F(G_1\text{-}v_1v_2)$. Теорема доказана.□

С л е д с т в и е .

Эйлерова характеристика любой связки K(n) на n вершинах равна 1, F(K(n)) = 1.

Это очевидно так как любая связка стягивается к тривиальному пространству K(1), состоящему из одной точки, для которого эйлерова характеристика равна 1.

На (Рис. 121) изображены точечные пространства K(1), K (4), K(5) и v x G. Эйлеровы характеристики всех этих пространств равны 1, что легко проверить непосредственным вычислением. Рассмотрим, теперь, чему равна прямая сумма двух молекулярных пространств.

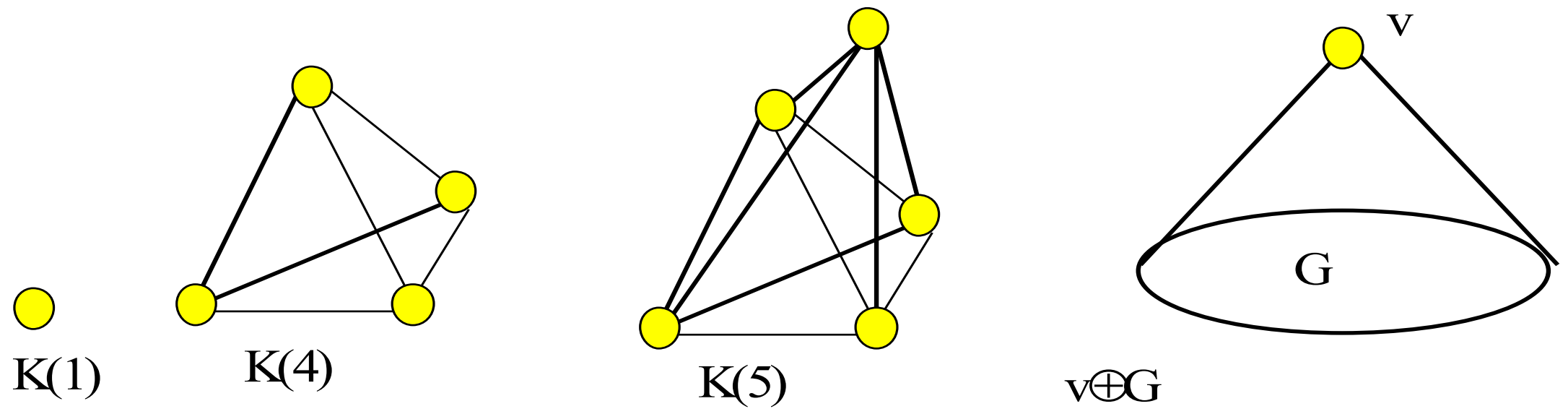

Рис. 121 Эйлеровы характеристики всех пространств равны 1, поскольку все они являются точечными.

Теорема 117

*Эйлерова характеристика прямой суммы двух молекулярных пространств G и H определяется выражением*

$$F(G\oplus H)=F(G)+F(H)-F(G)xF(H).$$

Д о к а з а т е л ь с т в о

Пусть функциональный вектор пространства G определяется выражением $f(G)=(n_1, n_2, n_3,...n_s)$, функциональный вектор пространства H определяется выражением $f(H)=(m_1, m_2, m_3,...m_s)$. В пространстве G⊕H каждая точка из G соединена с каждой точкой из H. Следовательно, все точки каждой связки K(p) из G соединены со всеми точками каждой связки K(r) из H. Следовательно, получается связка K(p+r), и количество таких связок равно $n_p m_r$. Используя этот простой подход получаем функциональный вектор пространства G⊕H, $f(G\oplus H)=((n_1+m_1), (n_2+n_1m_1+m_2), (n_3+n_1m_2+n_2m_1+m_3), (n_4+n_1m_3+n_2m_2+n_3m_1+m_4),...(n_s+n_1m_{s-1},+...+n_{s-1}m_{s-1}+m_s),...(n_sm_s))$. Отсюда получаем, что $F(G\oplus H)=((n_1+m_1)-(n_2+n_1m_1+m_2))+(n_3+n_1m_2+n_2m_1+m_3)-(n_4+n_1m_3+n_2m_2+n_3m_1+m_4)+...+(-1)^{s+1}(n_s+n_1m_{s-1},+...+n_{s-1}m_{s-1}+m_s)+...+(-1)^{2s+1}(n_sm_s))=F(G)+F(H)-F(G)xF(H)$. Теорема доказана.□

Справедливость этой формулы легко проверяется на конкретных примерах.

# *ОБЩИЕ И ЛОКАЛЬНЫЕ СВОЙСТВА ЭЙЛЕРОВОЙ ХАРАКТЕРИСТИКИ*

Эйлерова характеристика является важной функцией, позволяющей, наряду с другими, выяснить топологическую природу молекулярного пространства. Она, прежде, всего характеризует гомотопический тип, и не меняется при точечных стягиваниях.

Перечислим некоторые свойства этой функции и рассмотрим примеры.

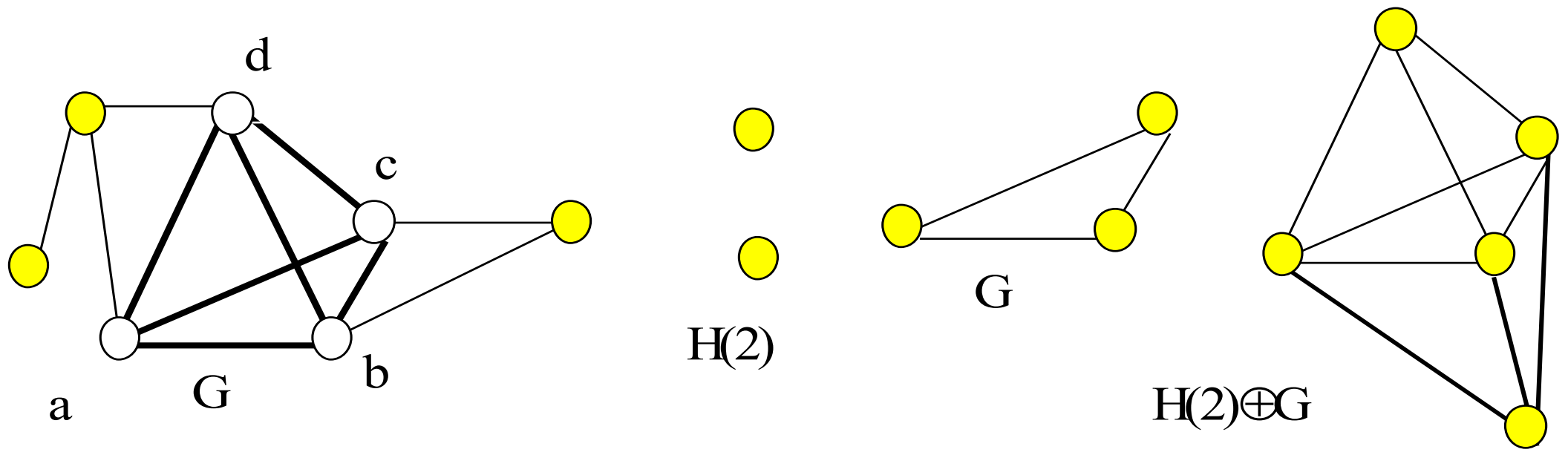

Рис. 122 Подпространство H=K(4) (слева) содержит точки a, b, c и d, и $m_1$=4, $m_2$=6, $m_3$=4, $m_4$=1. Очевидно, что F(H)=1. Кроме того, число вершин, ребер и треугольников, не лежащих целиком в H, определяется числами: $n_1$=3, $n_2$=5, $n_3$=2. Следовательно, F(G)=3-5+2+1=1. Справа F(H(2)⊕G)=2-1=1.

Теорема 118

*Эйлерова характеристика прямой суммы 0-мерной сферы $S^0$ и пространства G определяется выражением (Рис. 122)*

$$(F(S^0 \oplus G)=2-F(G).$$

Доказательство

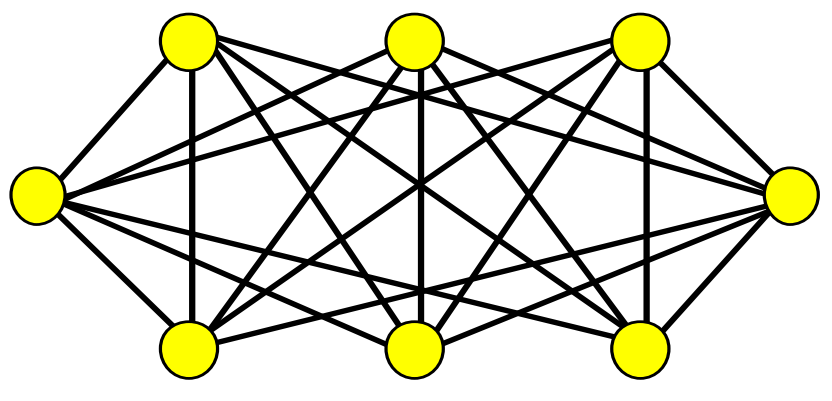

Рис. 123 3-х дольное пространство K(2,3,3). Согласно формуле (3), его эйлерова характеристика определяется выражением F(K(2,3,3))=1-(1-2)(1-3)(1-3)=-3.

Очевидно, что $F(S^0)=2$. Из предыдущего следует, что $F(S^0 \oplus G)=F(S^0)+F(G)-F(S^0)xF(G)$= 2+F(G)-2F(G)=2-F(G). Теорема доказана.□

Теорема 119

*Пусть G и H являются пространством и его подпространством, имеющими функциональные вектора f(G) = ($m_1+n_1$, $m_2+n_2$, $m_3+n_3$,......., $m_s+n_s$),*
*f(H) = ($m_1$, $m_2$, $m_3$,......., $m_s$) соответственно. Тогда*

$$F(G) = \sum_{k=1}^{s} (-1)^{k+1} n_k + F(H) \ ,$$

*где $n_p$ есть число связок K(p) в G, для которых, по крайней мере, одна вершина не принадлежит H, $m_p$ есть число связок K(p) в H.*

Д о к а з а т е л ь с т в о

Пусть G и H являются пространством и его подпространством (Рис. 119), имеющими функциональные вектора f(G) = ($m_1+n_1$, $m_2+n_2$, $m_3+n_3$,......., $m_s+n_s$), f(H) = ($m_1$, $m_2$, $m_3$,......., $m_s$) соответственно. Тогда $m_r$ есть количество связок K(r), целиком лежащих в H, следовательно, $n_r$ есть количество связок K(r), целиком не лежащих в H. Теорема доказана.□

Теорема 120

*Пусть H(n)-вполне несвязное пространство на n вершинах и G-любое пространство. Тогда эйлерова характеристика прямой суммы H⊕G определяется следующим выражением:*
*F(H⊕G) = n-(n-1)F(G).*

Д о к а з а т е л ь с т в о

Доказательство следует из предыдущей теоремы. □

На Рис. 122 F(G)=1. Следовательно, F(H(2)⊕G)=2-(2-1)1=1. Непосредственные вычисления дают тот же результат.

Теорема 121

*Пусть K($n_1$,$n_2$,...$n_p$) является p-дольным пространством. Его эйлерова характеристика определяется по формуле:*

$$F(n_1, n_2, \ldots, n_p) = 1 - \prod_{k=1}^{k=p} (1 - n_k) \ .$$

Д о к а з а т е л ь с т в о

Согласно определению K($n_1$,$n_2$,...$n_p$)=H($n_1$)⊕H($n_2$)⊕H($n_p$). Дальнейшее очевидно. □

Согласно этой формуле, если все n превышают 1, то p-дольный пространство не является точечным, так как его эйлерова характеристика отлична от 1.

На Рис. 123 изображено 3-х дольное пространство K(2,3,3). Согласно формуле его эйлерова характеристика определяется выражением

F(K(2,3,3))=1-(1-2)(1-3)(1-3)=-3. Непосредственный подсчет дает то же самое число.

Рассмотрим, теперь, как эйлерова функция пространства G выражается через некоторые функции окаемов его вершин. В некотором смысле это позволяет описать глобальную функцию на всем пространстве через

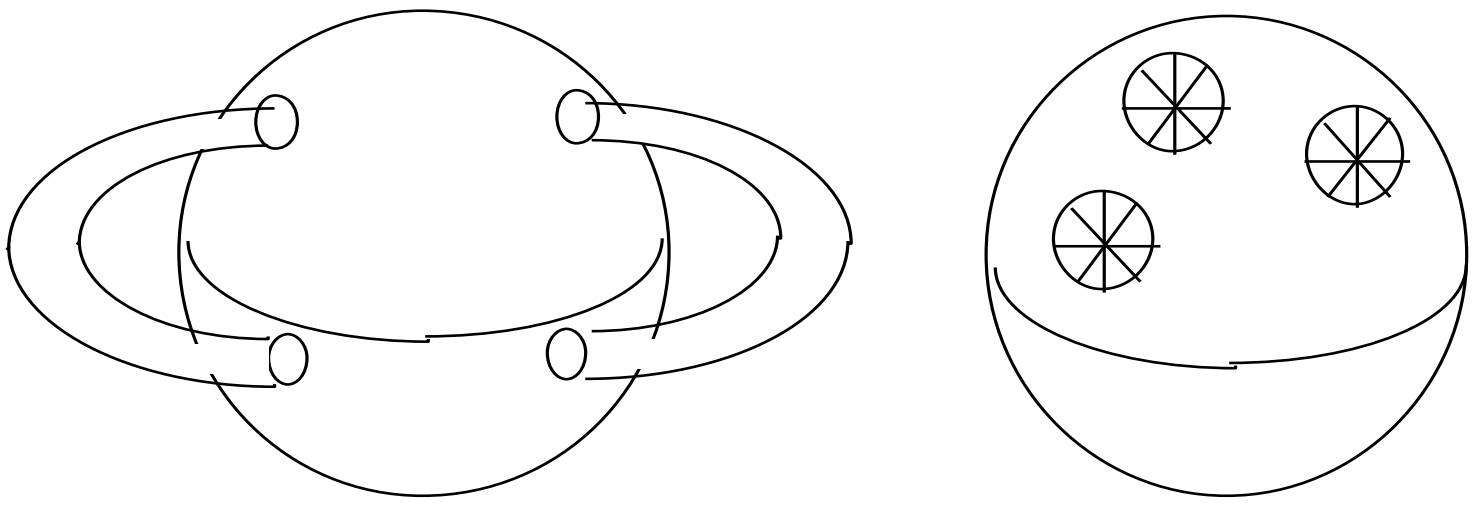

Рис. 124 Слева находится сфера $A_2$ с двумя приклеенными ручками, $F(A_2)=-2$, справа-сфера $B_3$ с тремя лентами Мебиуса, $F(B_3)=-1$.

некоторые функции его подпространств, то есть через локальные свойства пространства. Еще раз напомним понятия, которыми мы будем пользоваться.

1. Пространство G на множестве вершин $(v_1, v_2,... v_t)$

2. f-вектор пространства $f(G) = (n_1, n_2,.... n_s)$, где $n_i$-число связок K(i) пространства G на i вершинах.

3. f-вектор окаема $O(v_i)$ вершины $v_i$ $f(O(v_i)) = (m^i_1, m^i_2,... m^i_s)$

4. эйлеровы функции F(G) и $F(O(v_i))$

5. Логарифмические функции, обозначаемые L, пространства G и окаемов $O(v_i)$, определяемые выражением

$$L(O(v_i)) = \sum_{k=1}^{s} (-1)^{k+1} \frac{m^i_k}{k+1}.$$

6. введем p-степень $d^p(K(r))$ связок K(r) как число связок K(p), каждая из которых содержит K(r) как подпространство. Например, на Рис. 122 слева связка K(2), состоящая точек (db) принадлежит двум треугольникам (adb) и (cbd), а также одному тетраэдру (abcd). Следовательно, $d^3(bd)=2$, $d^4(bd)=1$.

Основная формула, которую мы будем доказывать будет сформулирована нами первой. $F(G) = t - \sum_{i=1}^{t} \sum_{k=1}^{s} (-1)^{k+1} \frac{m^i_k}{k+1} = t - \sum_{i=1}^{t} L(O(v_i))$

Эта формула позволяет изучать структуру и свойства молекулярного пространства в целом по структуре окаемов точек этого пространства. Эта формула будет полезной при анализе однородных покрытий непрерывных поверхностей, таких как сферы, торы, проективная плоскость. Читатель может самостоятельно проверить справедливость этой формулы на

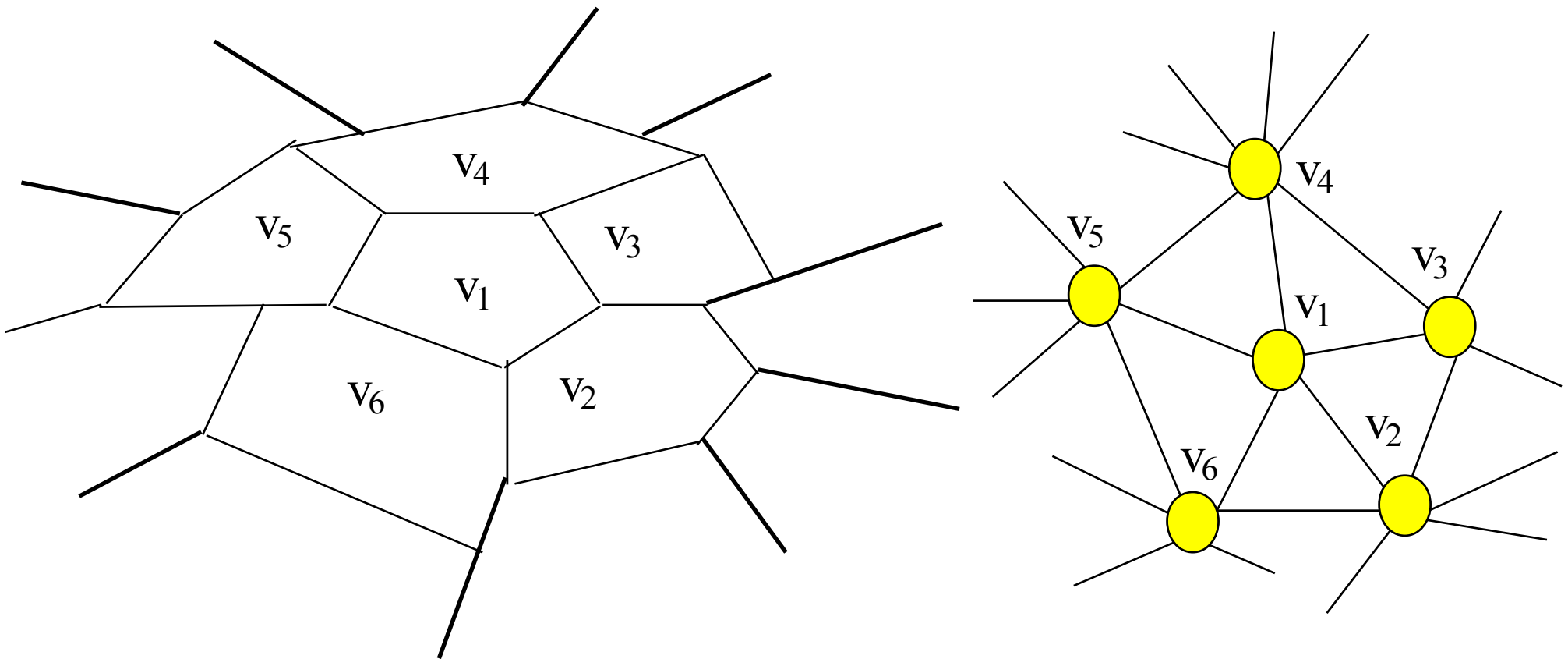

Рис. 125 Разбиение двумерной поверхности и его молекулярное пространство.

пространствах, изображенных на рисунках. Недавно автор обнаружил независимое доказательство сходной формулы в одной из статей по теории графов.

Докажем несколько утверждений, которые позволят нам доказать данную формулу.

Утверждение 1

$$\sum_{K(r)\subseteq G} d^p(K(r)) = C^r{}_p n_p,$$

*где* $m^i{}_0 = 1$.□

Д о к а з а т е л ь с т в о

Эта формула есть, по сути дела, тождество. Пусть функциональный вектор пространства G есть $f(G)=(n_1, n_2,\ldots n_d)$, функциональный вектор пространства $O(v_i)$ есть $f(O(v_i)) = (m^i{}_1, m^i{}_2,\ldots m^i{}_d)$. Перенумеруем все связки в пространстве G. Введем последовательность чисел $a^r{}_s$, определяемых условием: $a^q{}_s(p,r)=1$, если связка $K_s(r)$ принадлежит связке $K_q(p)$, и $a^q{}_s(p,r)=0$ в противном случае. Очевидно, что в каждой связке K(p) содержится $C^r{}_p$ связок K(r). Тогда

$$\sum_{K(r)\subseteq G} d^p(K(r)) = \sum_s \sum_q a_s^q(p,r) = \sum_q \sum_s a_s^q(p,r) = \sum_q C^r{}_p = C^r{}_p n_p$$

Утверждение доказано.□

Утверждение 2

$$\sum_{i=1}^{t} m_p^i = (p+1)n_{p+1},$$

$$m_0^1 \equiv 1.$$

*где $m_0^i = 1$.*

Д о к а з а т е л ь с т в о

Согласно обозначению, $m_p^i$ есть число связок K(p) на p точках в окаеме $O(v_i)$ точки $v_i$. Если присоединить к каждой такой связке саму точку $v_i$, то $m_p^i$ является числом связок K(p+1), содержащих точку $v_i$ в шаре $U(v_i)$. При суммировании по всем точкам пространства G каждая такая связка войдет в сумму (p+1) раз, поскольку она содержит (p+1) точек. Следовательно, $m_p^1+m_p^2+m_p^i+...+m_p^t=(p+1)n_{p+1}$. Утверждение доказано.□

Теорема 122

Пусть молекулярное пространство G содержит t точек $v_p$ и $m_s^p$ есть количество связок K(s) в окаеме точки $v_p$. Тогда эйлерова характеристика G определяется выражением

$$F(G) = t - \sum_{i=1}^{t} \sum_{k=1}^{s} (-1)^{k+1} \frac{m_k^i}{k+1} = t - \sum_{i=1}^{t} L(O(v_i)),$$

где $\sum_{k=1}^{s} (-1)^{k+1} \frac{m_k^i}{k+1} == L(O(v_i))$.

Д о к а з а т е л ь с т в о

Выразим $n_{p+1}$ из предыдущей формулы и подставим в выражение для эйлеровой характеристики.

$$F(G) = \sum_{k=1}^{s} (-1)^{k+1} n_k = \sum_{k=1}^{s} (-1)^{k+1} \frac{1}{k} \sum_{i=1}^{t} m_{k-1}^i = \sum_{i=1}^{t} \sum_{k=1}^{s} (-1)^{k+1} \frac{m_{k-1}^i}{k} =$$

$$t - \sum_{i=1}^{t} \sum_{k=1}^{s} (-1)^{k+1} \frac{m_k^i}{k+1} = t - \sum_{i=1}^{t} L(O(v_i))$$

Теорема доказана.□

Эта формула будет полезной при анализе однородных покрытий непрерывных поверхностей, таких как сферы, торы, проективная плоскость.

## *ПРИМЕНЕНИЕ ЭЙЛЕРОВОЙ ХАРАКТЕРИСТИКИ К ИЗУЧЕНИЮ СВОЙСТВ РАЗБИЕНИЙ И МОЩЕНИЙ НЕКОТОРЫХ ЗАМКНУТЫХ ДВУМЕРНЫХ ПРОСТРАНСТВ*

В этом разделе мы рассмотрим применение полученной формулы на традиционном наборе двумерных пространств: сфере, торе, проективной плоскости и бутылке Клейна, а также на некоторых других пространствах. Существует целое направление в математике, изучающее, как могут быть

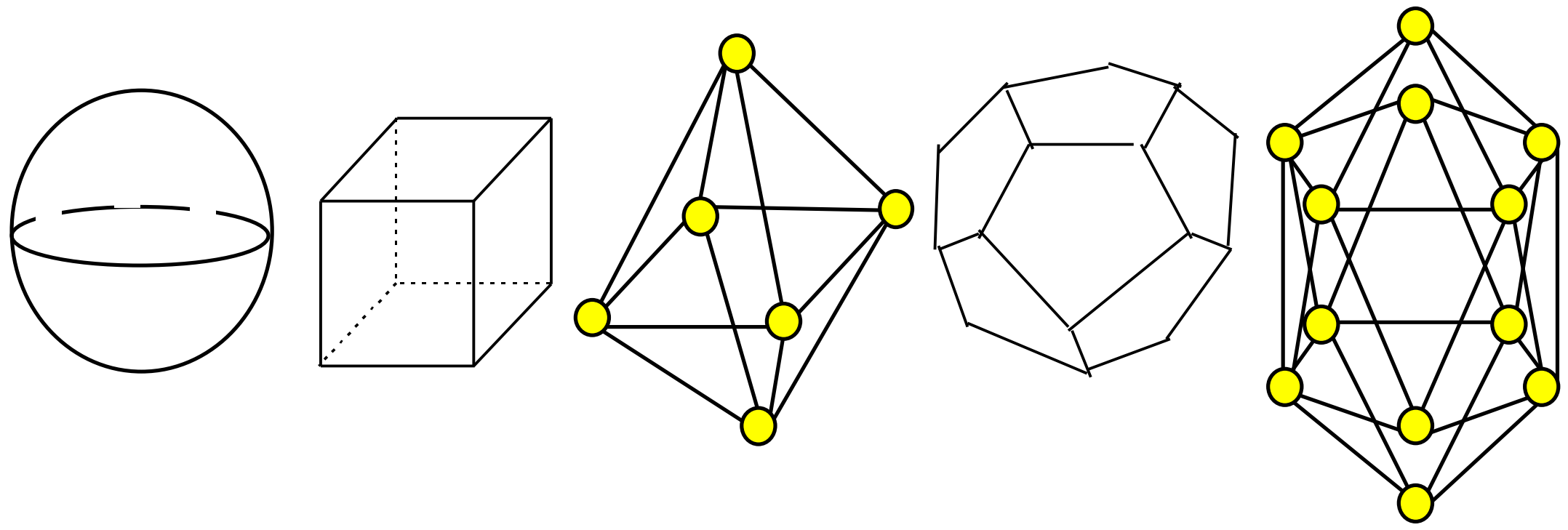

Рис. 126 Сфера и ее однородные покрытия. Существует только два однородных покрытия, нервы которых являются нормальными пространствами. Покрытие четырехугольниками состоит из 6 элементов, покрытие пятиугольниками имеет 12 элементов.

замощены такие поверхности многоугольниками с различным числом вершин. Начало этого изучения традиционно связывается с античными временами, когда были описаны так называемые платоновы тела. Мы не будем погружаться в глубоко в историю, а попробуем применить полученные результаты к анализу покрытий и разбиений двумерных замкнутых поверхностей.

Докажем несколько теорем, которые скорее являются следствиями из предыдущей теоремы. Если пространство однородно, то окаемы всех точек пространства изоморфны.

<u>Определение q-однородного пространства.</u>

> Молекулярное пространство называется q-однородным, если количество различных (неизоморфных) окаемов (или шаров) точек пространства равно q.

Рассмотрим вид эйлеровой характеристики.

Утверждение 3

Если пространство G 1-однородно, то его эйлерова характеристика имеет вид

$$F(G) = t(1 - \sum_{k=1}^{s} (-1)^{k+1} \frac{m_k}{k+1}) = t(1 - L(O(v)).$$

Если пространство G q-однородно, то его эйлерова характеристика имеет вид

$$F(G) = t - \sum_{i=1}^{t} \sum_{k=1}^{s} (-1)^{k+1} \frac{m_k^i}{k+1} = t - t_1 L(O(v_1) - t_2 L(O(v_2) - t_3 L(O(v_3) \ldots . - t_q L(O(v_q).$$

где $\sum_{k=1}^{s} (-1)^{k+1} \frac{m_k^i}{k+1} == L(O(v_i))$, $t_q$ есть число окаемов, изоморфных точке $v_q$.

Доказательство

Доказательство заключается в приведении подобных членов в формуле, доказанной выше..□

Следует напомнить результат из классической топологии, что любая двумерная гладкая компактная замкнутая поверхность гомеоморфна или сфере с 2k дырками, заклеенными k ручками, $A^2_k$, где k есть неотрицательное целое число (род поверхности), или сфере с k дырками, заклеенными k лентами Мебиуса, $B^2_k$, где k есть положительное число (Рис. 124). При этом $F(A^2_k)=2-2k$, $k\geq 0$, $F(B^2_k)=2-k$, $k\geq 1$. Введем, теперь, мощение и разбиение двумерных замкнутых поверхностей на многоугольники и посмотрим, как может использоваться эйлерова характеристика для этих целей. В качестве молекулярного пространства возьмем нерв мощения или разбиения. Напомним, что при этом каждому

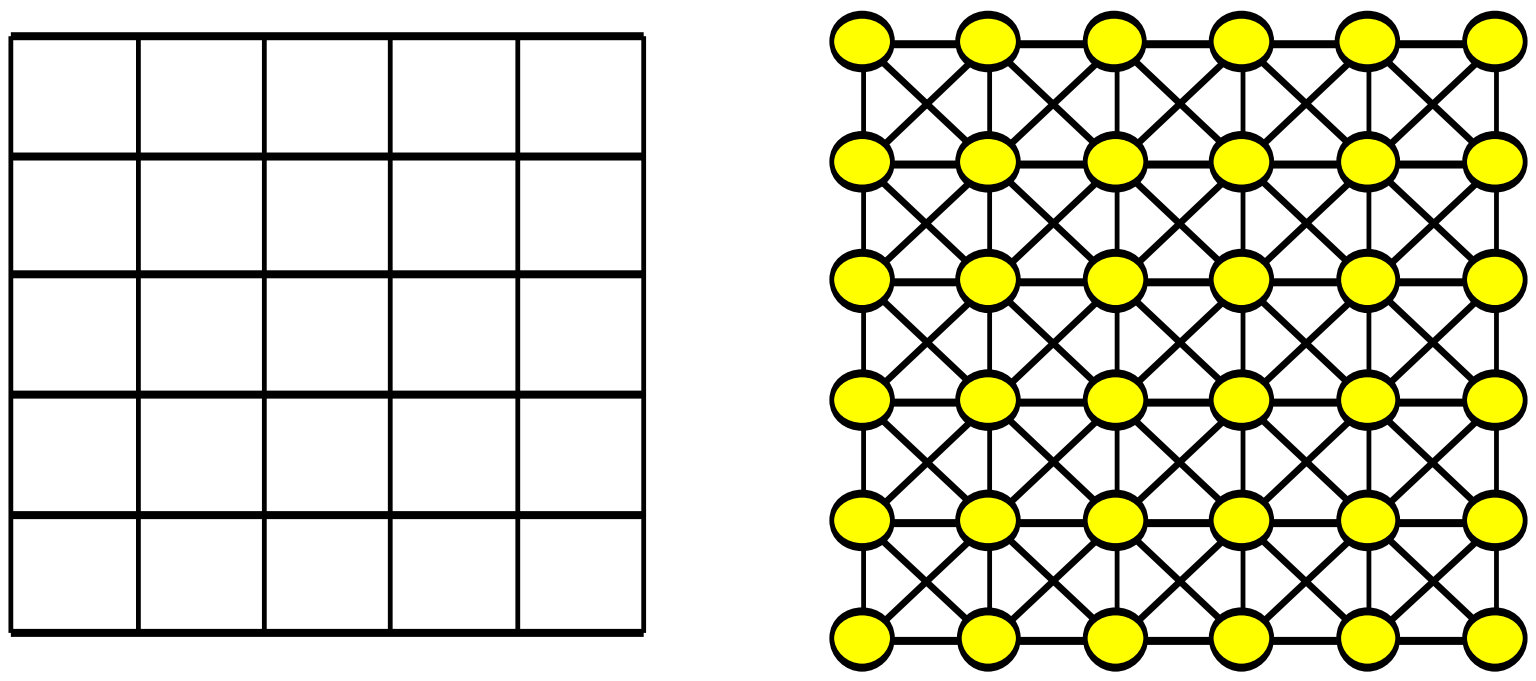

Рис. 127 Покрытие плоскости квадратами, нерв которого является молекулярным пространством, образованным прямым произведением двух одномерных пространств. Сфера не может быть покрыта таким образом, тогда как для тора и бутылки Клейна такое покрытие существует.

многоугольнику мы ставим в соответствие точку молекулярного пространства; если два элемента имеют хотя бы одну общую точку, то соответствующие им точки являются связными. Под одинаковыми многоугольниками (многоугольниками одного типа) мы будем понимать такие многоугольники, окаемы которых в молекулярном пространстве являются изоморфными.

**ДВУМЕРНАЯ СФЕРА $S^2$**

Эйлерова характеристика двумерной сферы $S^2$ равна 2. Рассмотрим, вначале, разбиение двумерной сферы на многоугольники так, чтобы пересечение многоугольников происходило только по их сторонам (Рис. 125). В более общем виде под n-угольником будет пониматься элемент элементы разбиения, имеющий n соседних элементов, причем пересекающийся с каждым из соседних по отрезку ненулевой длины. Рассмотрим покрытие сферы одинаковыми многоугольниками. В этом случае для любой точки $f(O(v_i)) = (m, m, 0, 0,...)$. Тогда, согласно предыдущему,

$$F(S^2) = t(1 - L(O(v\ \ ))) = t(1-(\frac{m}{2}-\frac{m}{3})) = t(1-\frac{m}{6}) = 2$$

Так как t и m должны быть целыми положительными числами и m≥4, мы имеем только два решения: m=4, t=6 и m=5, t=12 (Рис. 126). Из этой формулы также следует, что невозможно покрыть сферу многоугольниками со сторонами больше 5. Этот результат дает еще одно решение одной старинной задачи, существующей со времен Ньютона. Какое количество одинаковых шаров можно расположить на другом шаре так, чтобы ближайшие шары касались друг друга? Геометрическими расчетами было доказано, что это возможно для двенадцати шаров. Долгое время пытались выяснить, возможно ли это для 13 шаров, в конце 19 века было доказано, что это невозможно. Геометрические расчеты в таких случаях очень сложны и трудоемки. Полученный нами результат говорит о том, что невозможно покрыть сферу никаким количеством шаров, кроме 6 и 12.

Рассмотрим, теперь, покрытие сферы двумя типами многоугольников. При этом, естественно, окаем каждой точки должен содержать не менее 4-х точек и $t=t_1+t_2$.

$$F(S^2) = t - t_1L(O(v_1) - t_2L(O(v_2) = t_1(1-L(O(v_1)) + t_2(1-L(O(v_2)) = t_1(1-\frac{m_1}{6}) + t_2(1-\frac{m_2}{6}) = 2,$$

где $m_1$ и $m_2$ количество многоугольников, касающихся каждого из двух данных(или количество сторон в каждом из двух многоугольников). Легко убедиться, что все платоновы тела удовлетворяют этому условию.

Рассмотрим, теперь, покрытие сферы четырехугольниками, образующими молекулярное пространство, которое возникает при прямом произведении двух одномерных отрезков (Рис. 127).

В этом пространстве f-вектор окаема любой точки определяется выражением $f(O(v))=(8, 12, 4)$ и, следовательно, $L(O(v))=1$. Подстановка этого значения в выражение для эйлеровой характеристики обращает эйлерову характеристику в 0. Это означает, что однородное покрытие сферы данным способом невозможно.

**ДВУМЕРНЫЙ ТОР $T^2$**

Как известно, эйлерова характеристика двумерного тора равна нулю. Как и в случае сферы рассмотрим, вначале, разбиение тора на многоугольники одного типа так, чтобы пересечение многоугольников происходило только по их сторонам (Рис. 125). В этом случае для любой точки $f(O(v_i)) = (m, m, 0, 0,...)$. Тогда, согласно предыдущему,

$$F(T^2) = t(1 - L(O(v\ ))) = t(1-(\frac{m}{2}-\frac{m}{3})) = t(1-\frac{m}{6}) = 0$$

Так как t и m должны быть целыми положительными числами и $m \geq 4$, мы имеем только одно решение: $m=6$. Из этой формулы также следует, что тор можно замостить только шестиугольниками, и невозможно покрыть

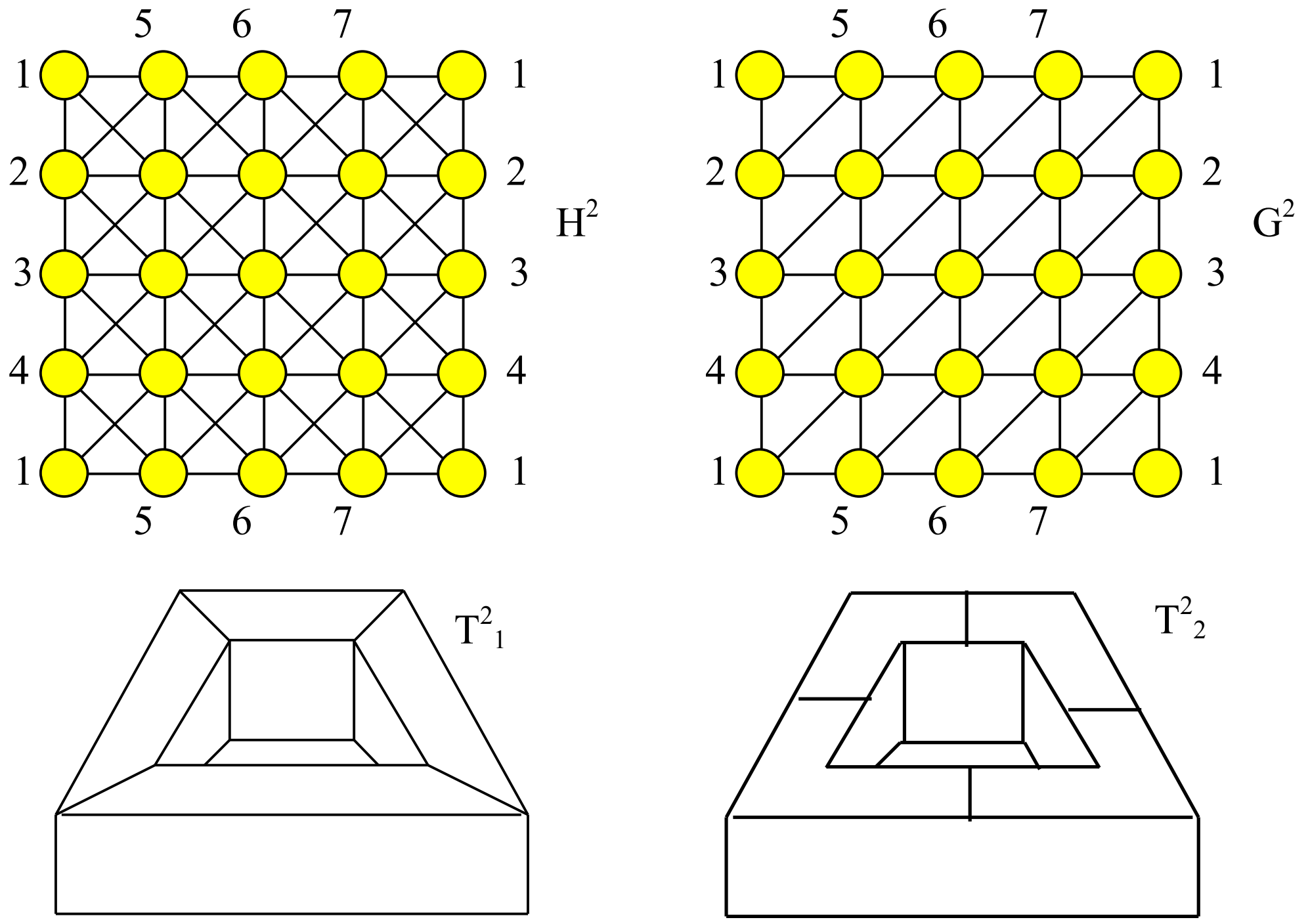

Рис. 128 Справа изображено разбиение тора на 16 элементов таким образом, что каждый элемент разбиения имеет общую сторону с каждым из шести элементов, окружающих его. Слева изображено разбиение тора того же вида, как на Рис. 127. Оба молекулярных пространства состоят из 16 элементов.

многоугольниками со сторонами отличными от 6. При этом количество элементов в разбиении не определяется через эйлерову характеристику. Однако прямая проверка показывает, что 16 шестиугольников составляют наименьшее число элементов разбиения (Рис. 110).

Рассмотрим теперь покрытие тора двумя различными типами многоугольников. При этом, естественно, окаем каждой точки должен содержать не менее 4-х точек. Очевидно, $t=t_1+t_2$.

$$F(T^2) = t - t_1 L(O(v_1) - t_2 L(O(v_2) = t_1(1 - L(O(v_1)) + t_2(1 - L(O(v_2)) = t_1(1 - \frac{m_1}{6}) + t_2(1 - \frac{m_2}{6}) = 0,$$

где $m_1$ и $m_2$ количество многоугольников, касающихся каждого из двух данных (или количество сторон в каждом из двух многоугольников). Из этой формулы следует, что при покрытии тора двумя многоугольниками число сторон одного из них должно быть меньше 6, а другого-больше 6. Рассмотрим в качестве примера $m_1$=4, $m_2$=7. При этом получаем $t_2$=2$t_1$. Рассмотрим, теперь, покрытие тора четырехугольниками, образующими молекулярное пространство, которое возникает при прямом произведении двух одномерных отрезков (Рис. 127). В этом пространстве f-вектор окаема любой точки определяется выражением f(O(v))=(8, 12, 4) и, следовательно, L(O(v))=1. Подстановка этого значения в выражение для эйлеровой характеристики обращает эйлерову характеристику в 0. Такое разбиение и его молекулярное пространство изображены на Рис. 110.

**БУТЫЛКА КЛЕЙНА $K^2$**

Примерно такие же результаты, как и для тора, получаются для бутылки Клейна, для которой эйлерова характеристика равна нулю. Непрерывная бутылка Клейна не может быть реализована в трехмерном пространстве. Как и в случае тора бутылку Клейна можно замостить только шестиугольниками, и невозможно покрыть многоугольниками со сторонами отличными от 6. Возможно, что 16 шестиугольников составляют наименьшее число элементов разбиения (Рис. 129). При покрытии бутылки Клейна двумя многоугольниками число сторон одного из них должно быть меньше 6, а другого-больше 6. Возможно также покрытие бутылки Клейна, как и тора, четырехугольниками, образующими молекулярное пространство, которое возникает при прямом произведении двух одномерных отрезков (Рис. 127).

**ПРОЕКТИВНАЯ ПЛОСКОСТЬ $P^2$**

Эйлерова характеристика проективной плоскости равна единице, $F(P^2)=1$. Проверим, существует ли однородное покрытие проективной плоскости каким-либо n-угольником.

$$F(P^2) = t(1 - \frac{m}{6}) = 1$$

Анализ этой формулы дает m=4, t=3, или m=5, t=6. И в том и в другом случаях число точек недостаточно для построения молекулярной проективной плоскости, в чем можно убедиться простой проверкой. Следовательно, не существует однородной молекулярной проективной плоскости, так же как не существует однородного мощения непрерывной проективной плоскости одной единственной черепицей. В качестве контрпримера мы изобразили минимальную нормальную молекулярную проективную плоскость, состоящую из одиннадцати точек, и покрытие плоскости черепицами, для которого молекулярное пространство является нервом.

Рассмотрим мощение проективной плоскости различными n-угольниками.

$$F(P^2) = t_1\ (1-\frac{m_1}{6}) + t_2\ (1-\frac{m_2}{6}) + ... + t_n\ (1-\frac{m_n}{6}) = 1.$$

Анализ этого равенства дает примерно такие же следствия как и для сферы. Попробуем найти мощение $P^2$ четырех и пятиугольниками. Выберем $m_1=4$, $m_2=5$, все остальные $m_k=0$. Тогда из соотношения следует, что: $2t_1+t_2=6$. Из этого выражения следует, что общее число черепиц не превышает 6. Непосредственная проверка показывает, что не существует

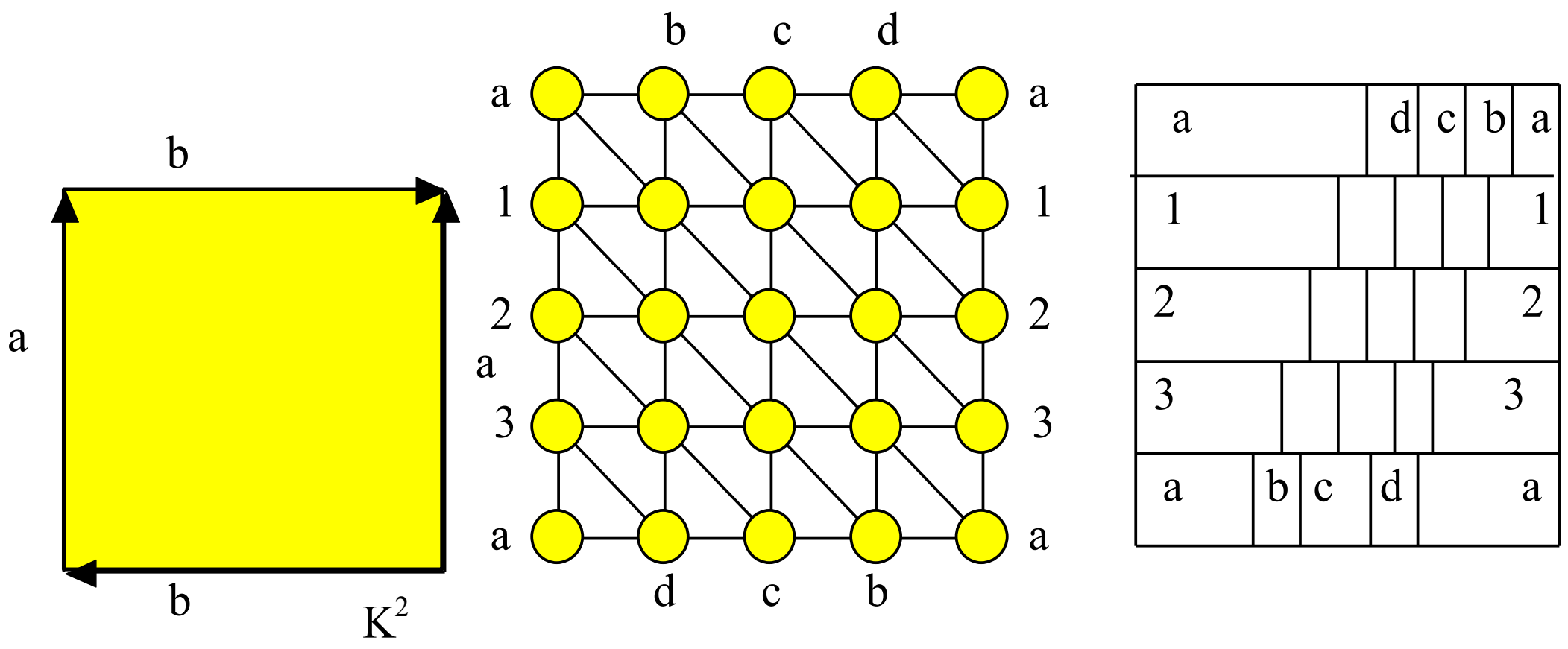

Рис. 129 Слева изображена непрерывная бутылки Клейна, справа дано ее разбиение таким образом, что каждый элемент имеет шесть соседних. В центре изображен нерв этого разбиения. Элементы разбиения могут быть также нарисованы в виде шестиугольников.

такого мощения проективной плоскости. Очевидно, также, что проективная плоскость не может иметь разбиения вида прямого произведения двух одномерных отрезков.

Мы кратко рассмотрели мощения и разбиения четырех замкнутых двумерных поверхностей, которые наиболее часто встречаются в математической литературе. Точно такой-же анализ может быть применен и к другим замкнутым поверхностям. В заключение сформулируем полученные результаты в виде утверждений.

Утверждение 4

- *Существует только два однородных мощения сферы: восемью четырехугольниками и двенадцатью пятиугольниками.*
- *Сфера не может быть вымощена только многоугольниками с числом сторон большим 5. В любом мощении сферы выпуклыми многоугольниками обязательно имеются четырехугольники и/или пятиугольники.*

Утверждение 5

- *Однородное мощение тора состоит только из шестиугольников.*
- *Любое неоднородное мощение тора многоугольниками должно содержать многоугольники как с числом сторон меньшим 6, так и с числом сторон большим 6.*

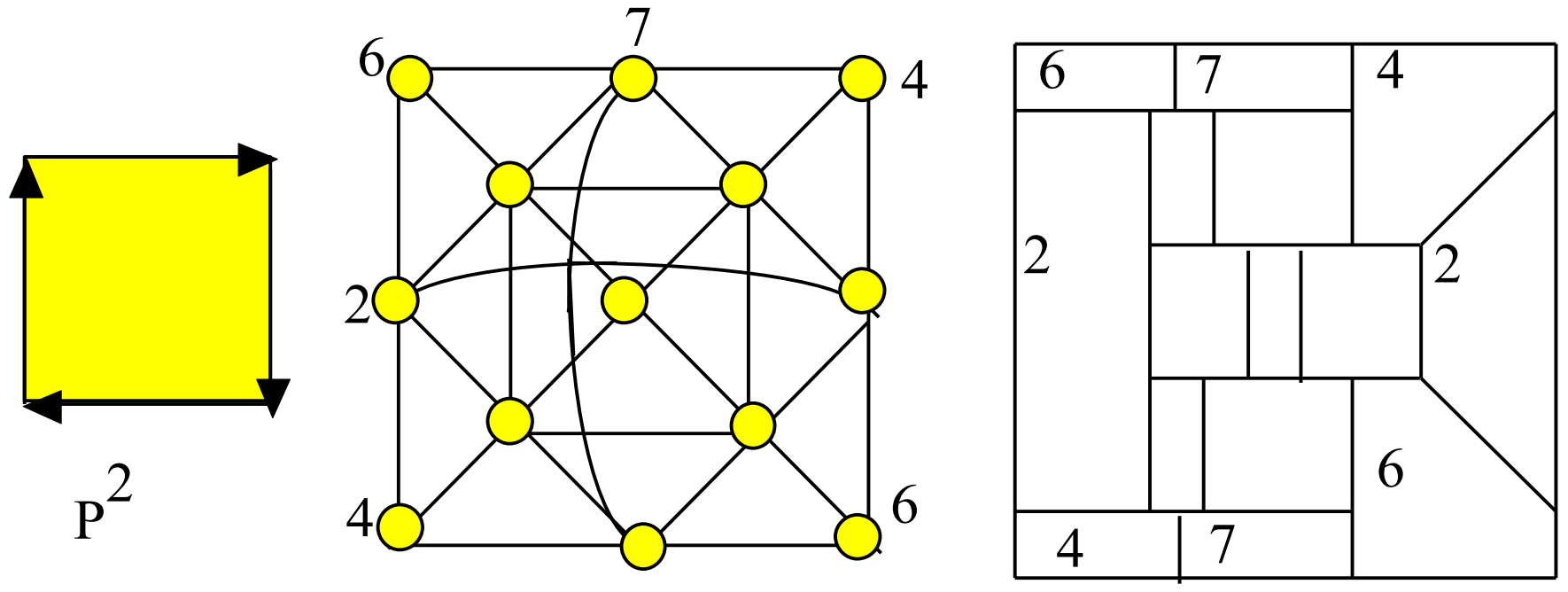

Рис. 130 Разбиение проективной плоскости на 11 элементов и его нерв, являющийся замкнутым двумерным пространством.

Утверждение 6

- *Однородное мощение бутылки Клейна состоит только из шестиугольников.*
- *Любое неоднородное мощение бутылки Клейна многоугольниками должно содержать многоугольники как с числом сторон меньшим 6, так и с числом сторон большим 6.*

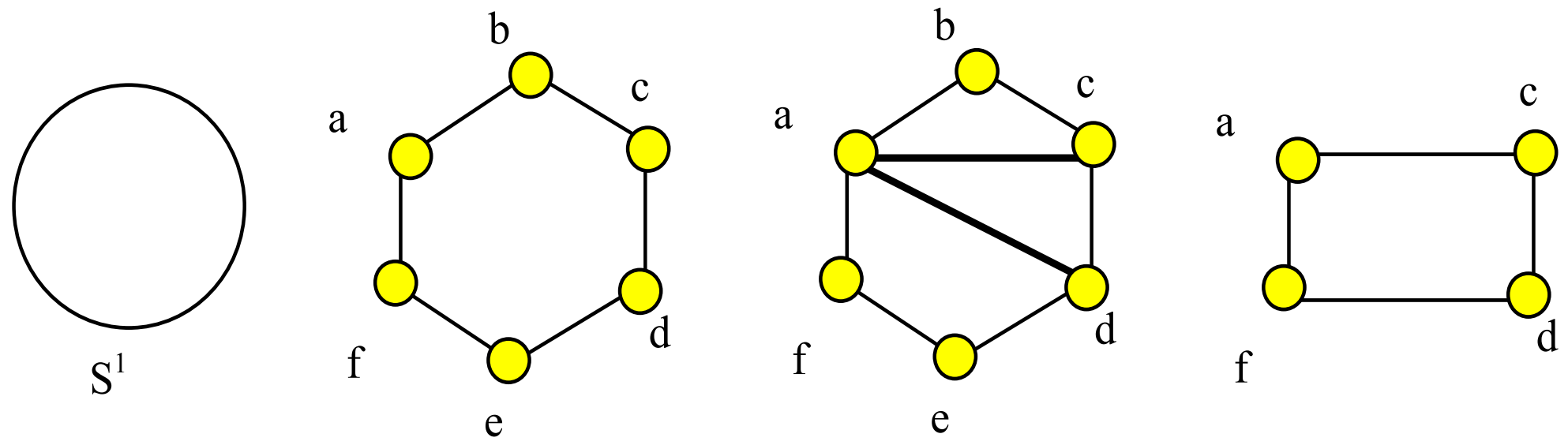

Рис. 131 Приведение окружности к минимальной, состоящей из 4-х точек.

Утверждение 7

- *Не существует однородного мощения проективной плоскости одинаковыми многоугольниками.*
- *Не существует мощения проективной плоскости только 4-х и/или 5-и угольниками. Проективная плоскость не может быть вымощена только многоугольниками с числом сторон большим 5. В любом* *мощении проективной плоскости выпуклыми многоугольниками обязательно многоугольники со сторонами бльше 6 с одной стороны, и имеются четырехугольники и/или пятиугольники с другой стороны.*

## *МИНИМАЛЬНЫЕ МОЛЕКУЛЯРНЫЕ ПРОСТРАНСТВА*

При работе на компьютере часто возникает проблема минимизации хранящейся в памяти машины информации, особенно если идет речь о больших массивах данных. Образы пространственных объектов, широко используются в игровых, научных, промышленных и других программах, и минимизация этих образов позволит освободить часть машинной памяти и уменьшить число операций, необходимых для их обработки. Пространства, с заданными на них математическими структурами, являются дигитальными образами непрерывных многомерных объектов, и их минимизация, таким образом, является составной и важной частью в теории молекулярных пространств. Рассмотрим некоторые аспекты этого направления в рамках теории молекулярных пространств.

Пространство может состоять из различного числа точек. Например, одномерная сфера-окружность содержит четыре и более точек (Рис. 131). Кроме того, при одинаковом числе точек, число связей может также различаться. На Рис. 131 изображены непрерывная одномерная сфера и ее молекулярные образы, являющиеся нормальными и регулярными молекулярными окружностями. Очевидно, что для хранения информации об окружности в машинной памяти желательнее использовать молекулярный образ, состоящий из 4-х вершин и четырех ребер, как занимающий наименьшую память в компьютере. Если окружность представлена шестью точками, то применение точечных преобразований, изображенных на том же Рис. 131 позволяет уменьшить число ребер и точек до 4-х. Однако, прежде чем применять преобразования, мы должны ввести некоторый критерий минимальности молекулярного пространства. Исходя из этого требования, повторим несколько определений и введем критерий минимальности.

<u>Определение. объема МП</u>

Количество точек МП G называется объемом молекулярного пространства G и обозначается |G|.

<u>Определение веса МП.</u>

Весом пространства G называется число ребер P(G) этого пространства.

<u>Определениенаименьшего МП.</u>

Пространство G называется наименьшим, если оно имеет наименьший объем V(G), а, при равных объемах, наименьший вес P(G) среди всех гомотопных ему пространств.

Напомним, что гомотопными называются те пространства, которые могут быть стянуты одно в другое при помощи точечных преобразований. Возникает вопрос, существует ли два гомотопных пространства с различными объемом и весом, одно из которых имеет наименьший вес, а другое-наименьший объем среди всех гомотопных пространств. Анализ

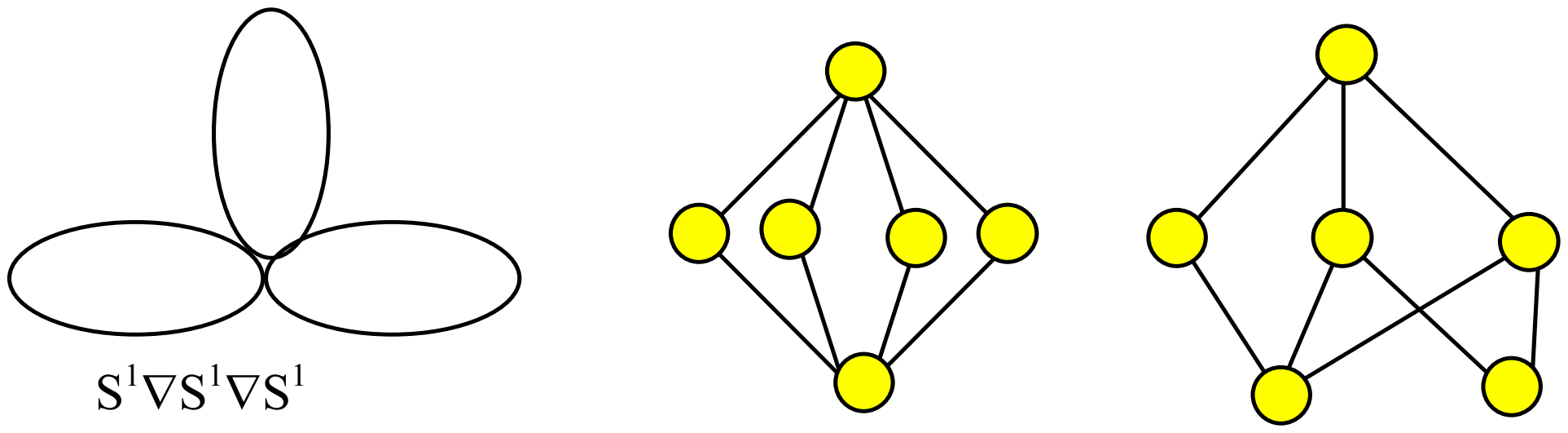

Рис. 132 Минимизация букета из 3-х окружностей.

пространств с небольшим числом точек показывает, что это не так, что пространство, наименьшее по объему, является также наименьшим по весу. Однако в общем случае этот вопрос не решен. Сформулируем это в виде гипотезы.

Г и п о т е з а 1.

> В любом гомотопическом классе существует пространство, наименьшее как по весу, так и по объему.

Во многих случаях для конкретных целей часто достаточно найти пространство со сравнительно небольшим числом точек, начиная с которого поиск наименьшего пространства становится затруднительным или требует простого перебора различных вариантов. Критерии для выбора такого пространства могут быть различными. Мы сформулируем такой критерий в виде определения минимального пространства.

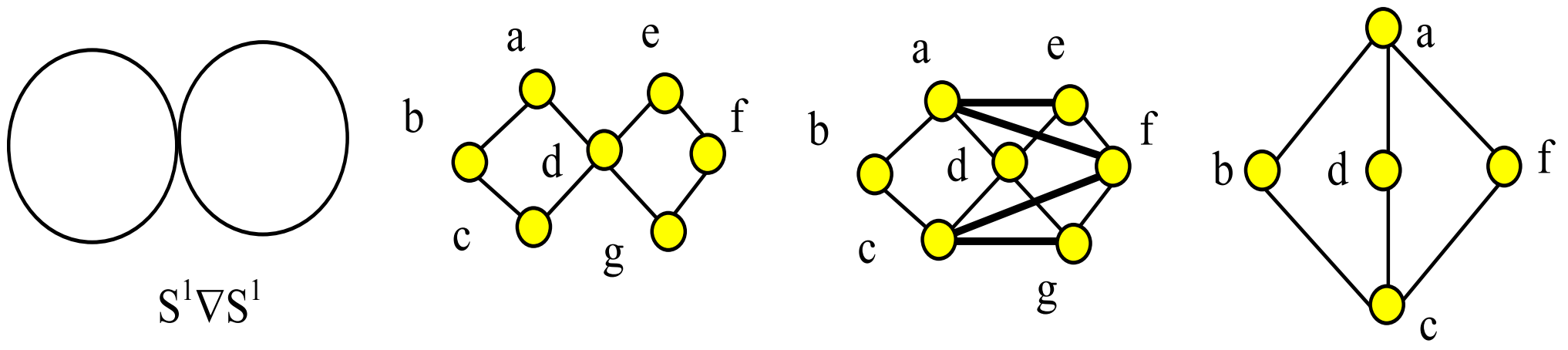

Рис. 133 Минимизация букета из 2-х окружностей.

Минимальное пространство не обязательно является наименьшим из всех гомотопных пространств. Однако, с практической точки зрения, минимальное пространство является часто наилучшей из моделей, поскольку поиск наименьшего пространства может потребовать значительных усилий.

<u>Определение минимального пространства.</u>

Пространство G называется минимальным, если:

1. Окаем любой точки пространства не является точечным.
2. Общий окаем любых двух смежных точек пространства не является точечным.
3. После точечного приклеивания всех возможных связей окаем любой точки пространства по-прежнему не является точечным.

Это определение гарантирует локальную минимальность пространства, а именно, любое точечное преобразование этого пространства может только увеличить его объем или вес.

З а м е ч а н и е .

Если пространство является наименьшим, то, тем самым, оно минимально.

Рассмотрим некоторые свойства минимальных и наименьших пространств.

Теорема 123

*Любое p-дольное пространство* $K(n_1, n_2,\ldots n_p)$ *не является точечным, если число вершин в любой доле не меньше двух,* $n_i > 1$, $i = 1, 2,\ldots p$.

Д о к а з а т е л ь с т в о

Согласно предыдущему, $F(K(n_1,n_2,\ldots,n_p)) = 1 - \prod_{k=1}^{k=p}(1-n_k)$. Следовательно, эйлерова характеристика пространства не равна 1, и пространство не является точечным. Теорема доказана.□

Теорема 124

*Любое p-дольное пространство* $K(n_1, n_2,\ldots n_p)$ *является минимальным, если число вершин в любой доле не меньше двух,* $n_i > 1$, $i = 1, 2,\ldots p$.

Д о к а з а т е л ь с т в о

Мы должны показать, что критерии минимальности выполняются на любом p-дольном пространстве в соответствии с условиями теоремы. Используем индукцию. Для малых p теорема проверяется непосредственно. Предположим, что теорема справедлива для $K(n_1,n_2,\ldots n_s)=H(n_1)\oplus H(n_2)\oplus H(n_s)$, s<p. Пусть s=p. Рассмотрим окаем O(v) любой точки v. Очевидно, что этот окаем является (p-1)-дольным пространством и, согласно предположению, неточечным и минимальным. Точно также общий окаем двух любых точек является (p-2)-дольным пространством и, согласно предположению, неточечным и минимальным. Для проверки критерия 3 рассмотрим две любые несмежные точки пространства. Очевидно, что они принадлежат одной доле, например, $H(n_1)$, и их общий окаем является (p-1)-дольным пространством и, согласно предположению, неточечным и минимальным. Следовательно, никакая связь не может быть установлена в этом пространстве. Теорема доказана.□

Введем частичное отношение порядка на множестве всех молекулярных пространств с конечным множеством точек.

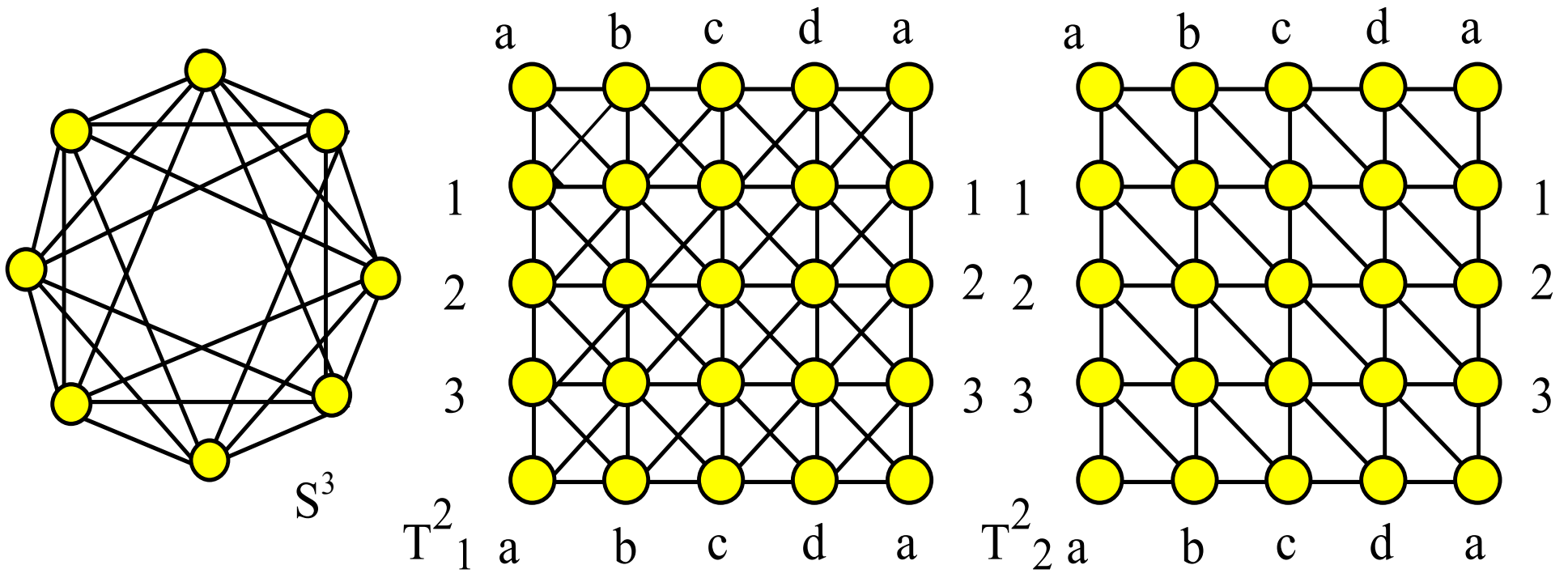

Рис. 134 Минимизация n-мерных сфер и тора.

Обозначение.

Пусть G(a,b) и H(c,d) есть два молекулярных пространства с количеством вершин a и c и ребер b и d соответственно. Обозначим G<H если a<c, или a=c и b<d.

Теорема 125

*Пространство H является наименьшим, если для всех G, G<H выполняется условие F(H)≠F(G).*

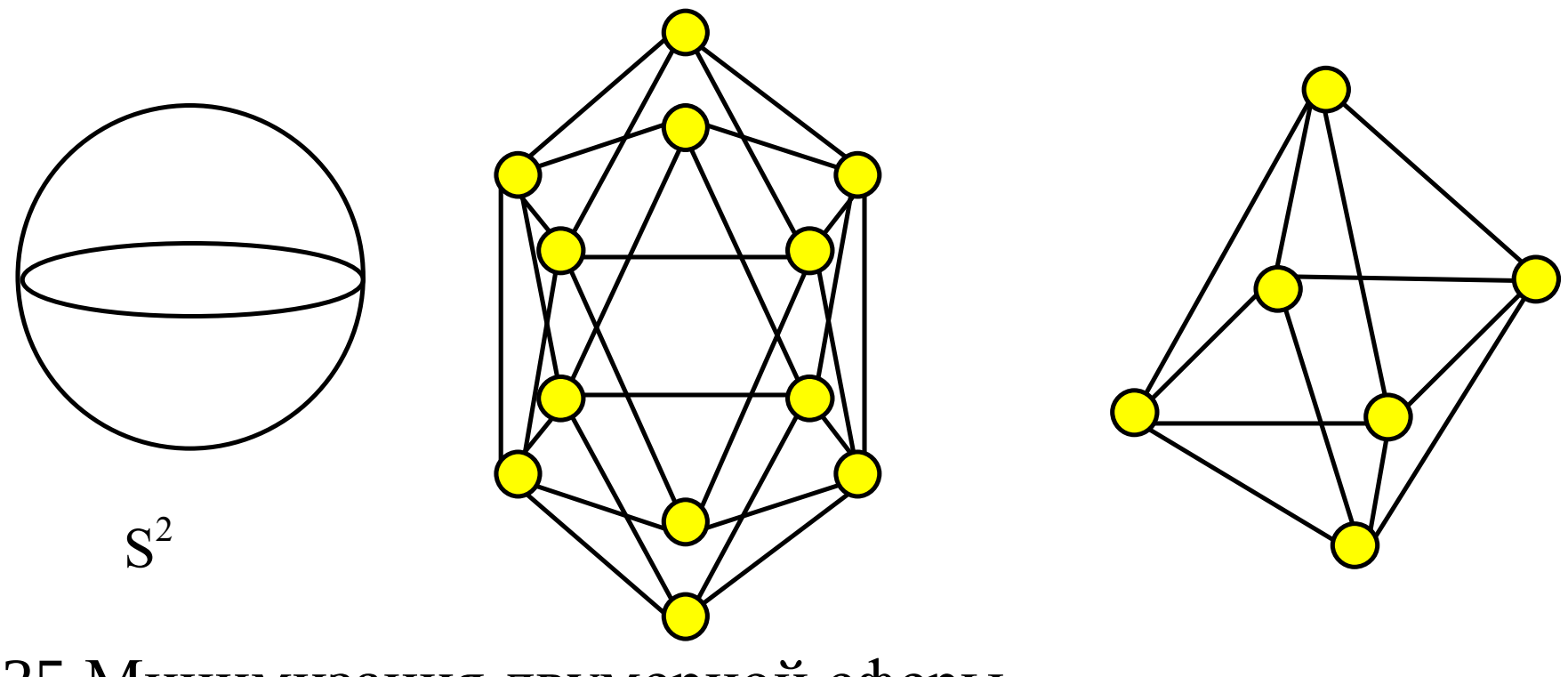

Рис. 135 Минимизация двумерной сферы.

Доказательство

Предположим, что H имеет точечную вершину v. Тогда G=H-v<H, но F(G)=F(H). Так как по условию теоремы это невозможно, то H не содержит точечной вершины. Аналогично доказывается, что H не содержит точечной связи. Предположим, что две несмежные точки имеют общий точечный окаем $O(v_1v_2)$. Введем связь $(v_1v_2)$. Предположим, что окаем некоторой точки v стал точечным после введения этой связи. Тогда эту точку можно отбросить, и возникшее пространство является гомотопным исходному. Это означает, что $(G+(v_1v_2)-v)<G$, $F(G+(v_1v_2)-v)=F(G)$. Это противоречит условию теоремы. Следовательно, после введения точечной связи все точки остаются неточечными. Теорема доказана.□

Рассмотрим некоторые примеры минимальных и наименьших пространств.

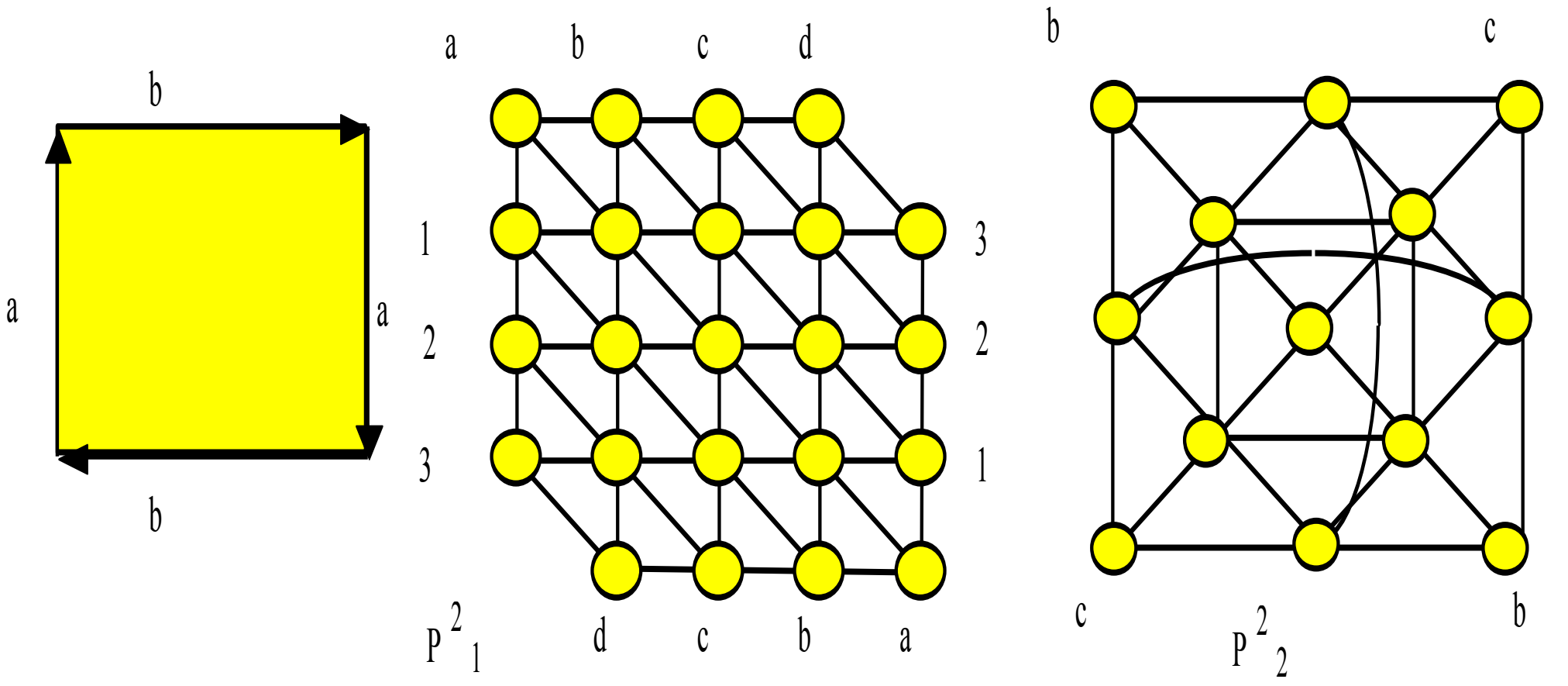

Рис. 136 Минимизация проективной плоскости.

**ОДНОМЕРНАЯ СФЕРА $S^1$ - ОКРУЖНОСТЬ (Рис. 131)**

Окружность образована шестью точками a, b, c, d, e и f. Приклеиваем точечные ребра (ac) и (ad), затем отбрасываем точечные вершины b и c. Получаем окружность из четырех точек a, d, e и f, которая является наименьшей. Она является пространством К(2,2). $F(S^1) = 0$.

**БУКЕТ ДВУХ ОКРУЖНОСТЕЙ $S^1\nabla S^1$ (Рис. 133)**

Букет $S^1\nabla S^1$ состоит из семи точек a, b, c, d, e, f и g. Приклеивая точечные ребра (ae), (af), (cg) и (cf), а затем отбрасывая вершины e и g, получаем наименьшее пространство $S^1\nabla S^1$, которое состоит из точек a, b, c, d, f и является двудольным пространством К(2,3), $F(S^1\nabla S^1) = -1$.

**БУКЕТ ИЗ ТРЕХ ОКРУЖНОСТЕЙ $S^1\nabla S^1\nabla S^1$ (Рис. 132)**

Это пространство является наглядным примером существования двух минимальных и, более того, наименьших неизоморфных пространств в одном гомотопическом классе. Одно из них есть К(2,4), а другое - К(3,3) без одного ребра, то есть К(3,3) - e, где e - одно из ребер в К(3,3). Оба они минимальны, так как ни к одному из них вообще нельзя приклеить ни одного точечного ребра. Очевидно, что переход от одного из них к другому возможен только через увеличение числа вершин, приклеивание точечной вершины и, уже затем, преобразование в другое пространство. Этот пример показывает, что возможно существование нескольких минимальных и даже наименьших неизоморфных пространств в данном гомотопическом классе. $F(S^1\nabla S^1\nabla S^1) = -3$.

**ДВУМЕРНАЯ СФЕРА $S^2$ (Рис. 135)**

Сфера, состоящая из 12 вершин, может быть преобразована путем приклеивания точечных ребер и последующего отбрасывания точечных вершин в минимальную и наименьшую сферу, состоящую из шести точек Эта сфера является трехдольным пространством К(2,2,2). Читатель легко может провести точечные преобразования самостоятельно. $F(S^2) = 2$.

**ТРЕХМЕРНАЯ СФЕРА $S^3$, СФЕРЫ $S^N$ (Рис. 134)**

Трехмерная сфера, изображенная на этом рисунке, состоит из восьми точек и является четырехдольным пространством К(2,2,2,2). Согласно теореме, она является минимальной. Непосредственная проверка показывает, что она является также наименьшей. Любая n-мерная сфера $S^n$, может быть представлена (n+1)-дольным пространством вида К(2,2,...2). В этом случае она минимальна. Не представляет особой сложности показать, что она также наименьшая. Эйлерова характеристика n-мерной сферы $S^n$ определяется выражением:

$F(S^n) = F(K(2,2,...2)) = 0$, если $n = 2к + 1$, $F(S^n) = F(K(2,2,...2)) = 2$, если $n = 2к$.

**ДВУМЕРНЫЙ ТОР $T^2$**

На Рис. 134 изображен двумерный тор. Левая модель не является минимальной, так как мы можем удалить все наклонные влево ребра. Полученная после удаления ребер модель, изображенная в правой части рисунка, является уже минимальной. Можно показать, что эта модель будет также наименьшей. Эйлерова характеристика тора, как это легко проверить, равна нулю. $F(T^2) = 0$.

**ПРОЕКТИВНАЯ ПЛОСКОСТЬ $P^2$**

На Рис. 136 изображены молекулярные модели проективной плоскости. Молекулярное пространство $P^2_1$ не является минимальным. Оно может

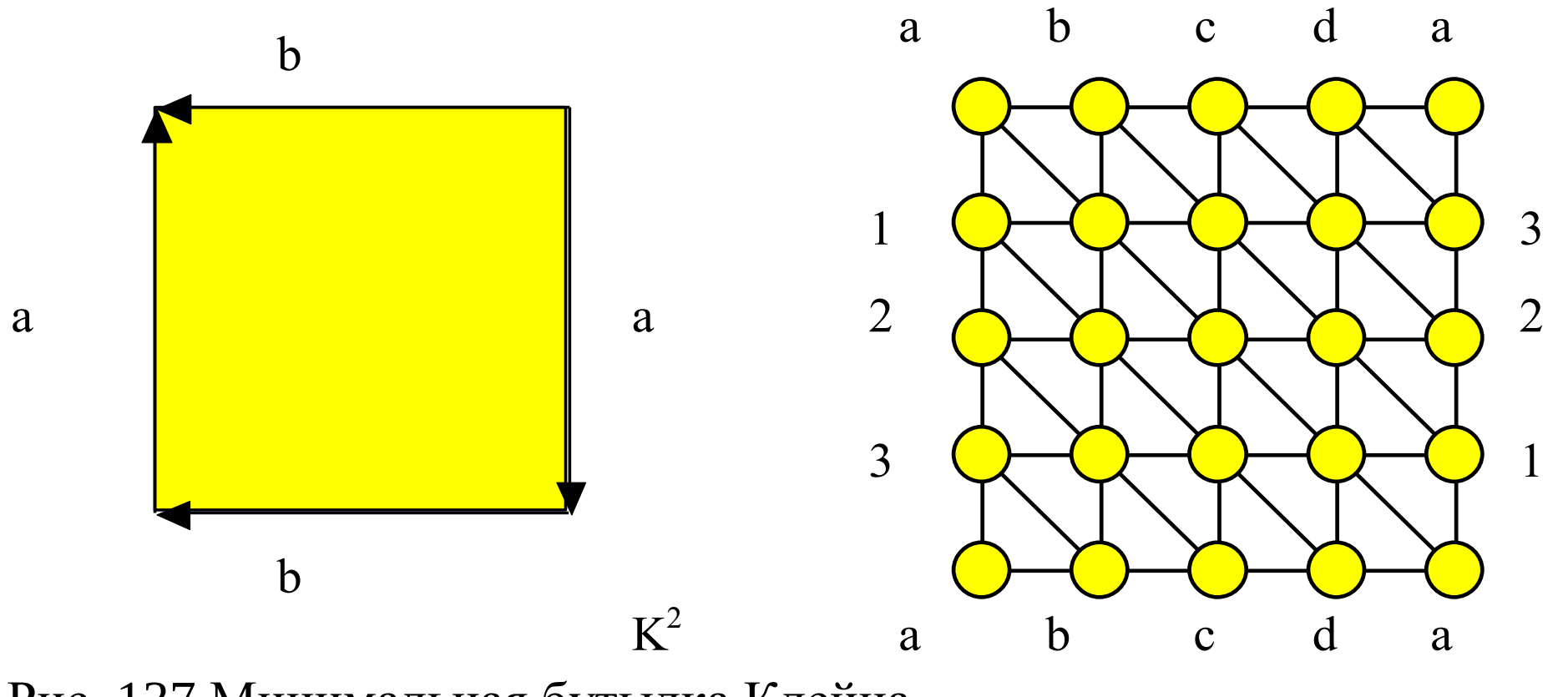

Рис. 137 Минимальная бутылка Клейна.

быть стянуто в пространство $P^2_2$ с меньшим числом точек и связей. Молекулярное пространство $P^2_2$ минимально, однако неизвестно, будет ли оно наименьшим. Легко подсчитать, что эйлерова характеристика проективной плоскости равна 1. $F(P^2) = F(P^2_1) = F(P^2_2)=1$.

З а д а ч а

Доказать, что проективная плоскость $P^2_2$ является наименьшей, или опровергнуть это, найдя наименьшую.

**<u>БУТЫЛКА КЛЕЙНА</u>. $K^2$ (Рис. 137)**

Бутылка Клейна получена, как проективная плоскость и тор, соответствующим склеиванием противоположных сторон квадрата. Изображенная на этом рисунке бутылка Клейна является минимальной. Не исключено, что это молекулярное пространство будет также наименьшим. Эйлерова характеристика дискретной бутылки Клейна равна нулю, так же как и в непрерывном случае. $F(K^2) = 0$.

З а д а ч а

Доказать, что бутылка Клейна $K^2$ является наименьшей, или опровергнуть это, найдя наименьшую.

С п и с о к   л и т е р а т у р ы   к   г л а в е   1 0 .

# ГРУППЫ ГОМОЛОГИЙ МОЛЕКУЛЯРНЫХ ПРОСТРАНСТВ. ТОЧЕЧНЫЕ ПРЕОБРАЗОВАНИЯ НЕ МЕНЯЮТ ГРУППЫ ГОМОЛОГИЙ МОЛЕКУЛЯРНОГО ПРОСТРАНСТВА

We define homology groups on molecular spaces and show that contractible transformations retain homology groups.

Полученные в этой главе результаты частично изложены в работах [16,25,43,47].

Группы гомологий топологических пространств [1,54] являются одной из характеристик, которая не меняется при определенных преобразованиях пространств. В классической комбинаторной топологии группы гомологий вычисляются, в частности, через симплициальное разбиение пространства. Это является прямым указанием на то, что группы гомологий могут быть

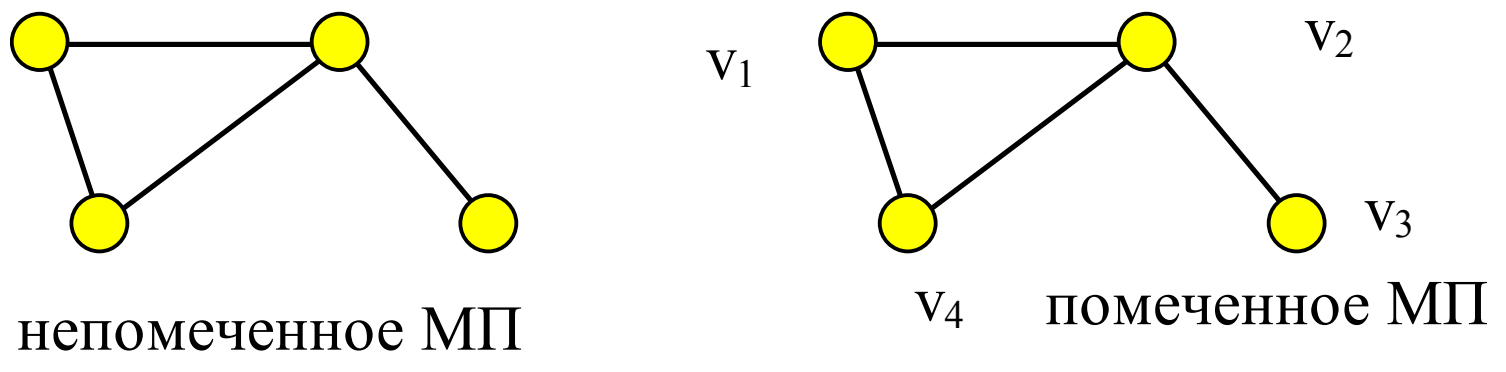

Рис. 138 Помеченное и непомеченное пространства. Пространство имеет 3-цепь $c_3=a_1[v_1v_2v_4]$, 2-цепь $c_2=a_1[v_1v_2]+a_2[v_2v_4]+a_3[v_1v_4]+a_4[v_2v_3]$ и 1-цепь $c_1=a_1v_1+a_2v_2+a_3v_3+a_4v_4$.

разумным образом определены на молекулярных пространствах. Упоминание о циклах и коциклах, которые можно рассматривать как необходимую базу для определения одномерной группы гомологий, имеется в книге Хараре [22] по теории графов. Однако, как ни странно, группы гомологий на графах не были определены и исследованы. Мы восполнили этот пробел. Мы определили группы гомологий на графах и молекулярных пространствах и доказали, что они являются инвариантами точечных преобразований молекулярных пространств, то есть точечные преобразования, примененные к МП, не меняют его групп гомологий. С прикладной точки зрения этот результат дает возможность различить два молекулярных пространства, не прибегая к точечным преобразованиям, которые часто связаны с большим количеством компъютерных операций. Предположим, мы имеем два молекулярных пространства, и нам необходимо выяснить, являются ли они гомотопными. Для этого мы вычисляем группы гомологий каждого из них. Если хотя бы одна из n-мерных групп не совпадает, то эти два пространства заведомо не гомотопны. К сожалению, так же как и в классической топологии, обратное, по-видимому, неверно. Мы не можем сказать, что два пространства гомотопны, если их группы гомологий совпадают. По крайней мере, в теории молекулярных пространств этот вопрос остается открытым.

Перейдем, теперь, к определению групп гомологий на молекулярных

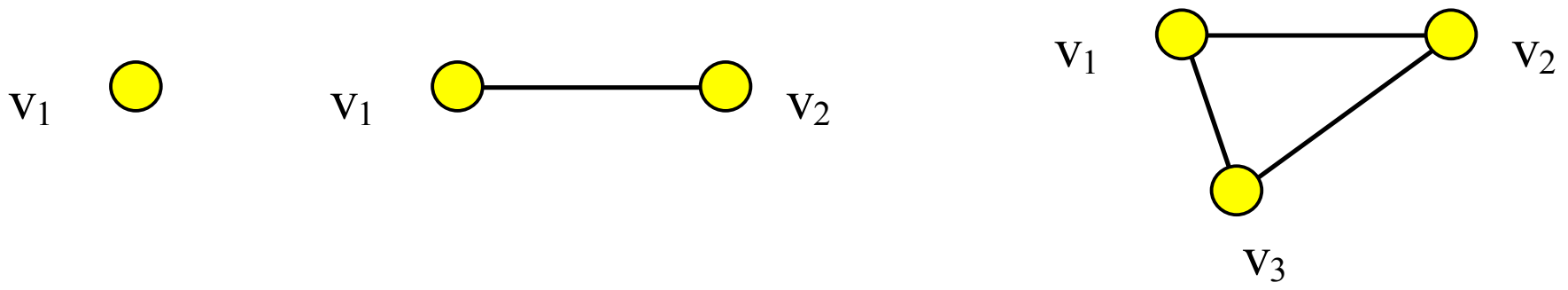

Рис. 139 Ориентированные связки $ф^1=[v_1]$, $ф^2=[v_1, v_2,]$, $ф^3=[v_{1,}, v_{2,}, v_{3,}]=-[v_{2,}, v_{1,}, v_3]$.

пространствах, сделав одно замечание.

Группы гомологий на молекулярном пространстве определяются автоматически, если совместить молекулярное пространство с симплициальным комплексом, а каждую его связку рассматривать как полный симплекс этого комплекса. Однако, как мы уже неоднократно об этом говорили, симплициальный комплекс и молекулярное пространство есть объекты различной природы. Поясним это примером. При стандартном определении группы гомологий для симплициального комплекса [1,54] точка считается 0-мерным объектом, две точки, соединенные ребром, считаются одномерным объектом и так далее. В теории молекулярных пространств такой подход лишен смысла, потому, что любое количество попарно связных точек (связка) является всегда 0-мерным пространством. В связи с этим мы вообще не будем употреблять выражение (n-1)-мерный симплекс, заменив его на связку из n точек, или n-связку.

<u>Определение помеченного пространства</u>

Пусть имеется молекулярное пространство G, содержащее множество V, состоящее из n точек. Этот молекулярное пространство называется помеченным или перенумерованным (labeled), если его точки отличаются одна от другой какими-либо метками, $V=(v_1, v_2, v_3,\dots v_n)$.

На Рис. 138 изображены непомеченное и помеченное молекулярные пространства. В дальнейшем мы рассматриваем помеченные молекулярные пространства.

О б о з н а ч е н и е  с в я з к и  н а  n  т о ч к а х

Связка $ф^n$ на n точках обозначается в виде $ф^n=(v_{1,}, v_{2,}, v_{3,}\dots v_{n,})$. (Напомним, что в связке все точки являются смежными).

Определение ориентированной связки

Связка $ф^n$ называется ориентированной, если выбран порядок точек.

$ф^n=[v_1, v_2, v_3,.... v_n]=-[v_1, v_2, v_3,.... v_n]$.

Ориентированную связку будем обозначать квадратными скобками. На Рис. 139 изображены ориентированная связка $ф^3=[v_1, v_2, v_3]=-[v_2, v_1, v_3]$.

Определение класса эквивалентности ориентированных n-связок

Пусть выбрана некоторая ориентация связки $ф^n=[v_1, v_2, v_3,.... v_n]$, считающаяся положительной. Тогда все ориентированные связки, полученные из данной с помощью четного числа перестановок, образуют класс эквивалентности положительно ориентированных связок, обозначаемый любым его представителем, например, $ф^n=[v_1, v_2, v_3,.... v_n]$. Все ориентированные связки, полученные из данной с помощью нечетного числа перестановок, образуют класс эквивалентности отрицательно ориентированных связок, обозначаемый любым его представителем, например, $-ф^n=[v_2\ v_1, v_3,.... v_n]$.

Определение ориентированного молекулярного пространства

Помеченное молекулярное пространство G называется ориентированным, если задана ориентация каждой из его связок.

Определение ориентации прямой суммы пространств

Пусть имеются два ориентированных молекулярного пространства G и H с точками $V=(v_1, v_2, v_3,.... v_n)$ и $U=(u_1, u_2, u_3,.... u_n)$ соответственно. Тогда ориентация их прямой суммы $G\oplus H$ определяется ориентацией каждого из них и порядком их следования. При этом $v_i<u_k$ для любых $v_i$ и $u_k$. В применении к связке

$ф^n=[v_1, v_2, v_3,....v_n]$

$=v_1\oplus[v_2, v_3,.... v_n]=v_1\oplus v_2\oplus[v_3,.... v_n]$.

Определение ориентированной границы $\partial ф^n$ связки $ф^n$

Ориентированной границей $\partial ф^n$ связки $ф^n$ называется следующая линейная комбинация его подпространств (связок) на (n-1) точках:

$$\partial ф^n=\partial[v_1, v_2, v_3,.... v_n]=\sum_{k=1}^{n}(-1)^{k+1}[v_1 v_2 .. \bar{v}_k ... v_n],$$

где черта над v означает, что соответствующая точка опущена. По определению границей 1-связки является пустое пространство или 0, $\partial v=0$.

Например, $\partial[v_1,v_2,v_3]=[v_1,v_2]-[v_2,v_3]+[v_2,v_3]$.

Теорема 126

*Если пространство G есть связка, то $\partial\partial G=0$*

Доказательство

Теорема проверяется непосредственно.□

Определение. n-цепи молекулярного пространства G и группы n-цепей $C_n(G)$

Пусть задано ориентированное молекулярное пространство G. n-Цепью молекулярного пространства G с целочисленными коэффициентами a (в общем случае с коэффициентами из некоторой абелевой группы A) называется формальная линейная комбинация различных n-связок молекулярного пространства G следующего вида:

$$c_n=\sum a_k\phi^n_k,$$

где ориентация связок определяется порядком ориентации точек молекулярного пространства G, $a_k$ являются произвольными целыми числами. На множестве цепей естественным образом определяется сложение

$$c_n=\sum a_k\phi^n_k,\ d_n=\sum b_k\phi^n_k,\ c_n+d_n=\sum (b_k+b_k)\phi^n_k.$$

Очевидно, что множество всех цепей образует абелеву группу $C_n(G)$ с операцией сложения.

Пространство на Рис. 138 имеет 3-цепь $c_3=a_1[v_1v_2v_4]$, 2-цепь $c_2=a_1[v_1v_2]+a_2[v_2v_4]+a_3[v_1v_4]+a_4[v_2v_3]$ и 1-цепь $c_1=a_1v_1+a_2v_2+a_3v_3+a_4v_4$.

Определение вынесения общего множителя n-цепи молекулярного пространства G

Пусть каждая связка в цепи может быть представлена в виде $\phi^n_k=v_0\oplus\phi^{n-1}_k$, тогда цепь может быть представлена в виде

$$c_n=\sum a_k\phi^n_k=\sum a_k(v_0\oplus\phi^{n-1}_k)=v_0\oplus\sum a_k\phi^{n-1}_k=v_0\oplus c_{n-1},$$

где $c_{n-1}=\sum a_k\phi^{n-1}_k$.

Определение ориентированной границы n-цепи

Ориентированной границей n-цепи $c_n=\sum a_k\phi^n_k$ называется выражение

$$\partial c_n=\sum a_k\partial\phi^n_k.$$

По определению границей 1-цепи (линейной комбинации точек молекулярного пространства) является 0.

Определение n-цикла и n-циклической группы $Z_n(G)$

n-Циклом называется n-цепь $z_n$, граница которой равна 0, $\partial z_n=0$.

Множество всех n-циклов является подгруппой группы всех n-цепей и называется группой n-циклов молекулярного пространства G. Обозначение-$Z_n(G)$.

Определение n-границы и n-граничной группы $B_n(G)$.

n-Границей называется n-цепь $b_n$, являющаяся границей некоторой (n+1)-цепи, $\partial c_{n+1}=b_n$. Множество всех n-границ является подгруппой группы всех n-цепей, называемой группой n-границ молекулярного пространства G.. Обозначение-$B_n(G)$.

Определение n-цикла, гомологичного нулю

Пусть $z_n$ есть цикл, то есть $\partial z_n=0$. Если существует цепь $c_{n+1}$, такая что $\partial c_{n+1}=z_n$, то говорят, что $z_n$ гомологичен нулю, $z_n \sim 0$.

Определение гомологичных n-циклов

Два n-цикла называются гомологичными, если их разность гомологична 0, то есть является границей некоторой цепи. Если $\partial c_{n+1}=z_n-y_n$, то $z_n \sim y_n$.

Определение n-группы гомологий молекулярного пространства G

n-Группой гомологий молекулярного пространства G называется фактор-группа $H_n(G)=Z_{n+1}(G)/B_{n+1}(G)$.

Эта группа корректно определена, так как группа $B_{n+1}(G)$ является подгруппой группы $Z_{n+1}(G)$, поскольку $\partial\partial c_{n+1}=0$.

Замечание

Рассмотрим группу циклов. Легко видеть, что отношение гомологичности есть отношение эквивалентности, разбивающее группу циклов на классы эквивалентности. Эти классы эквивалентности являются элементами группы гомологий.

Для упрощения понимания дальнейших доказательств перечислим несколько свойств групп гомологий, которые мы намереваемся использовать в дальнейшем. В доказательствах мы будем использовать эти свойства без ссылок на них.

Свойства

На любом пространстве G:

$\partial(v_0 \oplus c_n)=c_n-v_0 \oplus \partial c_n$,

$\partial(v_0 \oplus v_1 \oplus c_n)=v_1 \oplus c_n-v_0 \oplus c_n+v_0 \oplus v_1 \oplus \partial c_n$.

$H_{n-1}(G)=0$ тогда и только тогда, когда из условия $\partial c_n=0$ следует, что $c_n=\partial c_{n+1}$.

Если пространство G состоит из одной точки, то $H_0(G)=A$, $H_n(G)=0$, где n>0.

Если пространство G связно, то $H_0(G)=A$.

Если пространство G имеет s компонент связности, то группа $H_0(G)$ равна прямому произведению s групп A, $H_0(G)=A\otimes A\otimes....\otimes A$.

Перейдем, теперь, к теоремам, доказательства которых для графов приведены в указанных ниже работах.

Теорема 127

*Пусть v-точка и G-пространство, тогда группы гомологий прямой суммы $v\oplus G$ есть $H_0(v\oplus G)=A$, $H_n(v\oplus G)=0$, где $n>0$.*

Доказательство

Так как $v\oplus G$ связно, то $H_0(v\oplus G)=A$. Пусть $c_n$ есть некоторая n-цепь в $v\oplus G$, $n>0$. Тогда она может быть представлена в форме $c_n=a_n\text{-}v\oplus a_{n-1}$, где $a_n$ и $a_{n-1}$ цепи в G. Из уравнения $\partial c_n=\partial a_n\text{-}a_{n-1}+v\oplus\partial a_{n-1}=0$ следует, что $\partial a_n\text{-}a_{n-1}=0$, $\partial a_{n-1}=0$. Следовательно, n-цикл может быть представлен в следующей форме: $z_n=a_n\text{-}v\oplus\partial a_n$. Отсюда ясно, что $z_n=a_n\text{-}v\oplus\partial a_n=\partial c_{n+1}=\partial(v\oplus a_n)$. Это означает, что любой n-цикл гомологичен нулю, и n-мерная группа

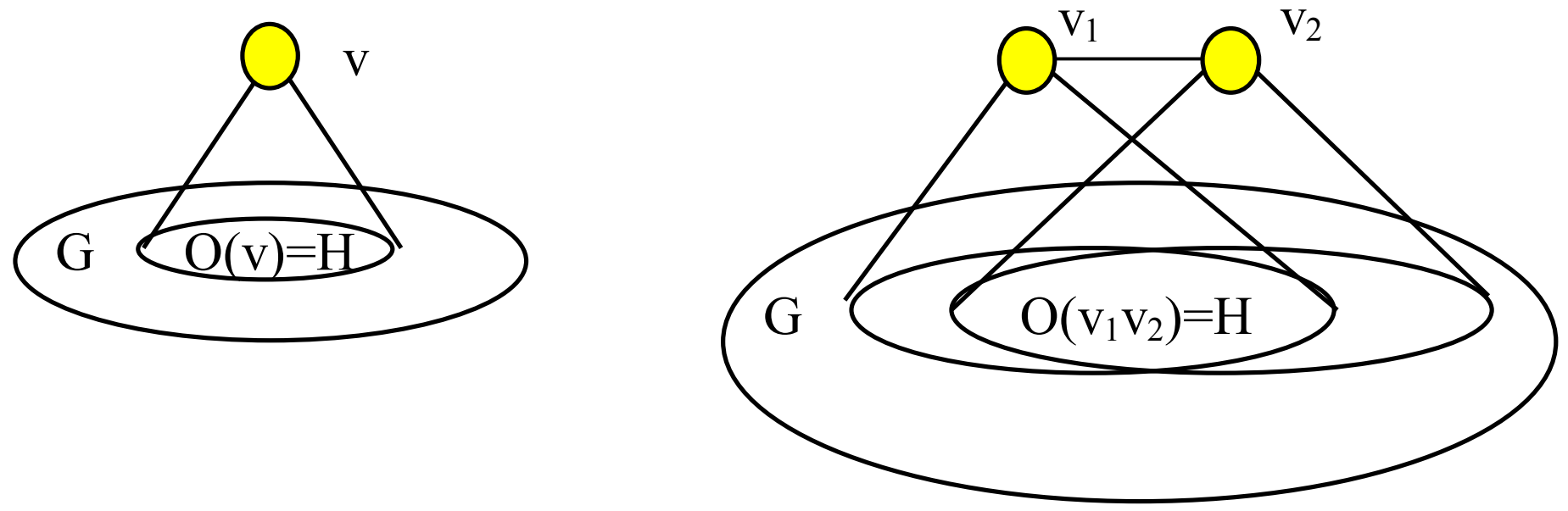

Рис. 140 МП G, окаем O(v)=H точки v (слева) и общий окаем $O(v_1v_2)=H$ двух точек $v_1$ и $v_2$ (справа) имеют группы гомологий точки. Тогда G+v и $G+(v_1v_2)$ также имеют группы гомологий точки.

гомологий равно нулю. Теорема доказана.□

Теорема 128

*Пусть G и H есть пространство и его подпространство, для которых группы гомологий есть $H_0(H)=H_0(G)=A$, $H_n(H)=H_n(G)=0$, где $n>0$. Тогда приклеивание точки v к G по H, O(v)=H, не меняет групп гомологий пространства G. $H_0(G+v)=A$, $H_n(G+v)=0$ (Рис. 140).*

Доказательство

Так как $H_0(G)=A$, то G является связным. При приклеивании точки v к G связность не нарушается, следовательно, $H_0(G+v)=A$.

Пусть G и H являются искомыми пространством и некоторым его подпространством, не совпадающим с G. Приклеим к G по H точку v так, чтобы O(v)=H, U=G+v. Из $c_n \in C_n(U=G+v)$ следует, что $c_n=a_n\text{-}v\oplus b_{n-1}$, где $a_n \in C_n(G)$ и $b_{n-1} \in C_n(H)$. Из $c_n=z_n \in Z_n(U=G+v)$ следует, что $z_n=a_n\text{-}v\oplus b_{n-1}$, $\partial z_n=\partial a_n\text{-}b_{n-1}+v\oplus\partial b_{n-1}=0$. Из $\partial a_n\text{-}b_{n-1}=0$, $v\oplus\partial b_{n-1}=0$ следует, что $\partial b_{n-1}=0$. Из $\partial b_{n-1}=0$ и $H_n(H)=0$, следует, что существует такой $b_n$, что $\partial b_n=b_{n-1}$. Из $\partial a_n\text{-}\partial b_n=0$ и $H_n(G)=0$, следует, что существует такой $a_{n+1}$, что $\partial a_{n+1}=a_n\text{-}b_n$. Из $\partial a_{n+1}=a_n\text{-}b_n$ и $\partial(v_0\oplus c_n)=c_n\text{-}v_0\oplus\partial c_n$, следует, что $z_n=a_n\text{-}v\oplus b_{n-1}=$ $\partial a_{n+1}+b_n\text{-}v\oplus\partial b_n=\partial a_{n+1}+\partial(v\oplus b_n)=\partial(a_{n+1}+v\oplus b_n)$. Следовательно, любой n-цикл в U есть граница, и $H_{n-1}(U)=0$, где n-1>0. Теорема доказана.□

Теорема 129

*Пусть G и H есть пространство и его подпространство, являющееся окаемом точки v, O(v)=H, для которых группы гомологий есть $H_0(H)=H_0(G)=A$, $H_n(H)=H_n(G)=0$, где n>0. Тогда отбрасывание точки v из G не меняет групп гомологий пространства G. $H_0(G\text{-}v)=A$, $H_n(G\text{-}v)=0$ (Рис. 140).*

Доказательство

Отбрасывание точки является операцией, обратной к приклеиванию, рассмотренному в предыдущей теореме. Доказательство теоремы есть взятое в обратном порядке доказательство предыдущей теоремы.□

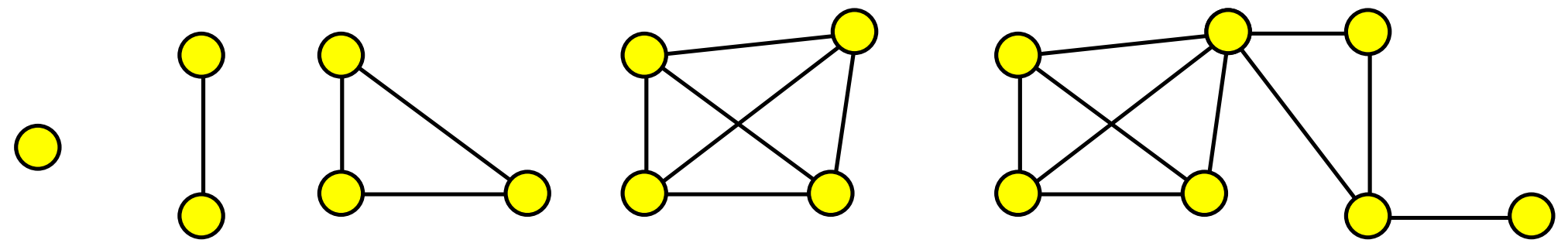

Рис. 141 Группы гомологий данных пространств $H_0(H)=Z$, $H_n(H)=0$, где n>0, поскольку все они являются точечными пространствами.

Теорема 130

*Пусть G и H есть пространство и его подпространство, для которых группы гомологий есть $H_0(G)=H_0(H)=A$, $H_n(G)=H_n(H)=0$, где n>0. Пусть H есть общий окаем двух несмежных точек $v_1$ и $v_2$, $O(v_1v_2)=H$. Тогда приклеивание связи $(v_1v_2)$ к G по H не меняет групп гомологий пространства G. $H_0(G+(v_1v_2))=A$, $H_n(G+(v_1v_2))=0$, n>1 (Рис. 138).*

Доказательство

Так как $H_0(G)=A$, то G является связным. При приклеивании связи $(v_1v_2)$ к G связность не нарушается, следовательно, $H_0(G+(v_1v_2))=A$.

Пусть G и H являются искомыми пространством и некоторым его подпространством, не совпадающим с G. Приклеим к G по H связь $(v_1v_2)$ так, чтобы $O(v_1v_2)=H$. Пусть $U=G+(v_1v_2)$. Из $c_n \in C_n(U=G+(v_1v_2))$ и $n>2$, следует, что $c_n=a_n+v_1\oplus v_2\oplus b_{n-2}$, где $a_n \in C_n(G)$ и $b_{n-2}\in C_{n-2}(H)$. Из $c_n=z_n\in Z_n(U=G+(v_1v_2))$ следует, что $z_n=a_n+v_1\oplus v_2\oplus b_{n-2}$ и $\partial z_n=\partial(a_n+v_1\oplus v_2\oplus b_{n-2})=\partial a_n-v_1\oplus b_{n-2}$ $+v_2\oplus b_{n-2}+v_1\oplus v_2\oplus\partial b_{n-2}=0$. Из $\partial z_n=\partial a_n-v_1\oplus b_{n-2}+v_2\oplus b_{n-2}+v_1\oplus v_2\oplus\partial b_{n-2}=0$ следует, что $\partial a_n-v_1\oplus b_{n-2}+v_2\oplus b_{n-2}=0$, $\partial b_{n-2}=0$. Из $\partial b_{n-2}=0$ и $H_n(G)=0$, $n>1$ следует, что существует такой $b_{n-1}\in C_{n-1}(G)$, что $\partial b_{n-1}=b_{n-2}$. Из $\partial a_n-v_1\oplus b_{n-2}$ +

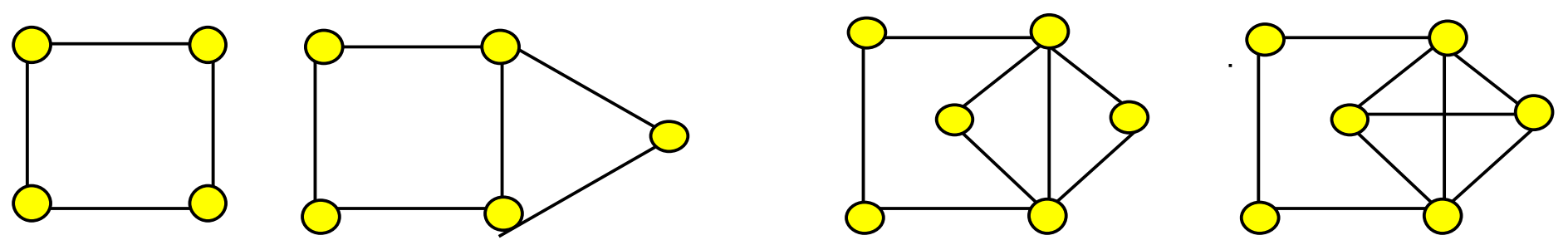

Рис. 142 Все изображенные пространства стягиваются к одномерной окружности и имеют соответствующие группы гомологий: $H_0=H_1=Z$, $H_n=0$, $n > 1$.

$v_2\oplus b_{n-2}=0$, $\partial b_{n-1}=b_{n-2}$, $\partial(v_1\oplus b_{n-1})=b_{n-1}-v_1\oplus\partial b_{n-1}$ и $\partial(v_2\oplus b_{n-1})=b_{n-1}-v_2\oplus\partial b_{n-1}$ следует, что $\partial(a_n-v_2\oplus b_{n-1}+v_1\oplus b_{n-1})=0$. Из $\partial(a_n-v_2\oplus b_{n-1}+v_1\oplus b_{n-1})=0$ и $a_n-v_2\oplus b_{n-1}+v_1\oplus b_{n-1}\in H_n(G+(v_1v_2))=0$, $n>2$, следует, что $a_n-v_2\oplus b_{n-1}+v_1\oplus b_{n-1}=\partial w_n$, где $w_n\in C_n(G)$. Из $a_n-v_2\oplus b_{n-1}+v_1\oplus b_{n-1}=\partial w_n$ следует, $z_n=v_2\oplus b_{n-1}-v_1\oplus b_{n-1}+\partial w_{n+1}+v_1\oplus v_2\oplus\partial b_{n-1}=\partial(v_1\oplus v_2\oplus b_{n-1})+\partial w_{n+1}=\partial(v_1\oplus v_2\oplus b_{n-1}+w_{n+1})$. Следовательно, этот цикл является в то же время границей, и $H_{n-1}(G+(v_1v_2))=0$, $n-1>0$.

Пусть $c_2$ есть некоторая 2-цепь в $C_2(U=G+(v_1v_2))$. Из $c_2\in C_2(U=G+(v_1v_2))$, следует, что $c_2=a_2+av_1\oplus v_2$, где $a_2\in C_2(G)$. Из $c_2=z_2\in Z_2(U=G+(v_1v_2))$ следует, что $z_2=a_2+av_1\oplus v_2$ и $\partial z_2=\partial(a_2+av_1\oplus v_2)=0$. Из уравнения $\partial z_2=\partial a_2-av_1+av_2=0$ следует, что $\partial a_2=av_1-av_2$. Выберем цепь $b_2=a(v_2\oplus v_3+v_3\oplus v_1)\in C_2(G)$ где $v_3$ есть некоторая точка из H. Так как $\partial b_2-\partial a_2=\partial(b_2-a_2)=0$, $b_2-a_2\in C_2(G)$ и $H_n(G)=0$, $n>1$, то $a_2=b_2+\partial w_3$, где $w_3$ есть некоторая 3-цепь, $w_3\in C_3(G)$. Из $z_2=a_2+av_1\oplus v_2$ и $a_2=b_2+\partial w_3$ следует, что $z_2=a(v_2\oplus v_3+v_3\oplus v_1)+av_1\oplus v_2+\partial w_3=\partial(av_1\oplus v_2\oplus v_3)+\partial w_3=\partial(av_1\oplus v_2\oplus v_3+w_3)$. Следовательно, этот цикл является в то же время границей, и $H_1(G+(v_1v_2))=0$. Теорема доказана.□

Теорема 131

*Пусть G и H есть пространство и его подпространство, для которых группы гомологий есть $H_0(G)=H_0(H)=A$, $H_n(G)=H_n(H)=0$, где $n>0$. Пусть H есть общий окаем двух смежных точек $v_1$ и $v_2$, $O(v_1v_2)=H$. Тогда отбрасывание связи $(v_1v_2)$ в G не меняет групп гомологий пространства G, $H_n(G-(v_1v_2))=H_n(G)$.*

Доказательство

Доказательство теоремы состоит из доказательства предыдущей теоремы, взятого в обратном порядке. □

Теорема 132

*Пусть G есть точечное пространство. Тогда группы гомологий пространства G есть группы гомологий точки, $H_0(G)=A$, $H_n(G)=0$, $n>1$ (Рис. 141).*

Доказательство

По определению точечного пространства любое точечное пространство получено из уединенной точки путем точечных приклеиваний и отбрасываний точек и связей, при которых группы гомологий не меняются. Следовательно, группы гомологий точечного пространства G есть группы гомологий точки, $H_0(G)=A$, $H_n(G)=0$, $n>1$. Теорема доказана.□

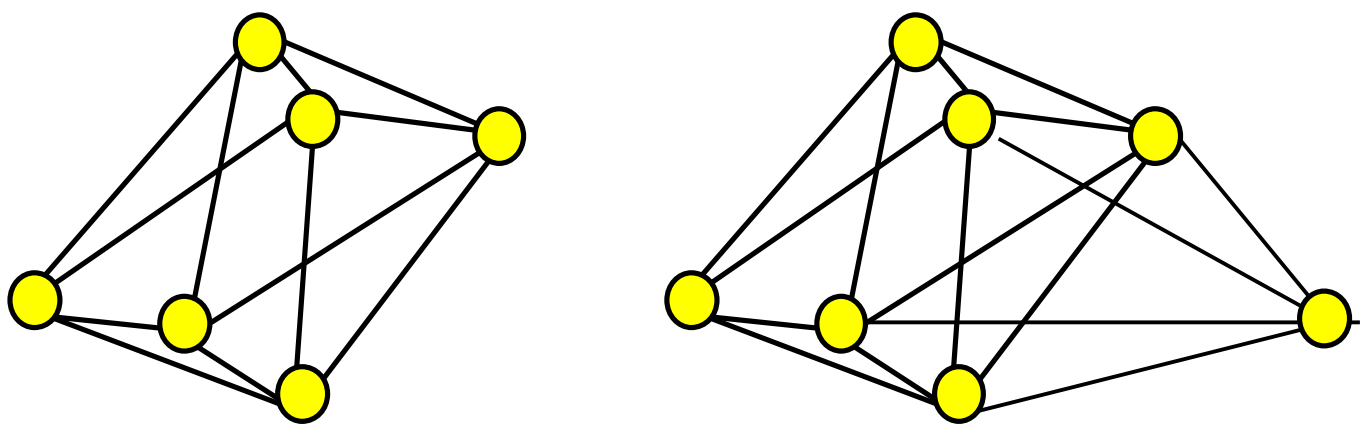

Рис. 143 Двумерная сфера и гомотопное ему пространство, имеющие группы гомологий $H_0=H_2=Z$, $H_1=H_n=0$, $n > 2$.

Теорема 133

*Пусть G есть пространство и H есть его подпространство, для которого группы гомологий есть $H_0(H)=A$, $H_n(H)=0$, где $n>0$. Тогда приклеивание (отбрасывание) точки v к G по H, $O(v)=H$, не меняет групп гомологий пространства G. $H_n(G+v)=H_n(G)$.*

Доказательство

При приклеивании точки v к G число связных компонент не меняется. Следовательно, $H_0(G+v)=H_0(G)$.

Пусть G и H являются искомыми пространством и некоторым его подпространством, не совпадающим с G. Приклеим к G по H точку v так, чтобы O(v)=H, U=G+v. Используя доказательство предыдущей теоремы мы получаем, что любой цикл $z_n$ в U=G+v имеет вид $z_n=w_n+\partial(v\oplus b_n)$, где $w_n$

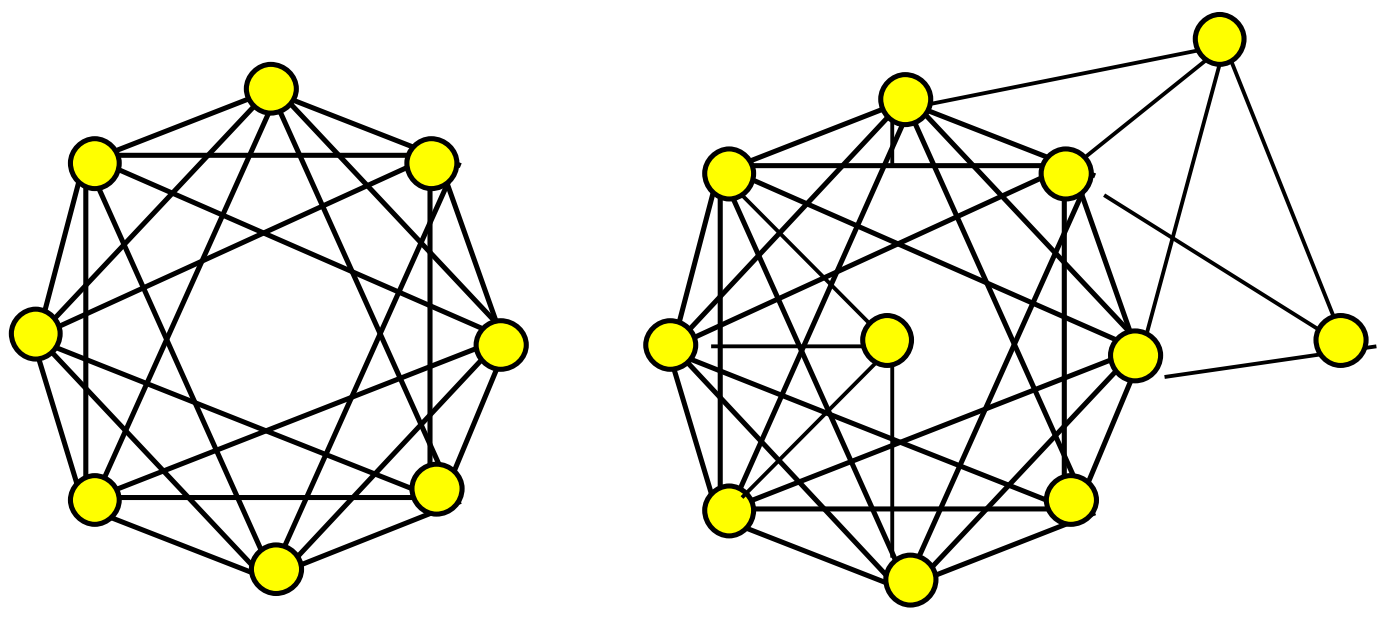

Рис. 144 Трехмерная сфера и гомотопное ему пространство, имеющие группы гомологий $H_0=H_3=Z$, $H_1=H_2=H_n=0$, $n > 3$.

есть некоторый цикл в G. Введем гомоморфизм f: $Z_n(U)\to Z_n(G)$, $f(z_n)=f(w_n+\partial(v\oplus b_n))=w_n$.

Очевидно, что любой цикл на U переходит в цикл на G, и для любого цикла на G существует прообраз, являющийся циклом на U. Очевидно, также, что любая граница из U переходит в границу из G, и для любой границы из G существует прообраз, являющийся границей в U.

$f(\partial(a_{n+1}+v\oplus b_n)=f(\partial a_{n+1}+\partial(v\oplus b_n))=\partial a_{n+1}$.

Следовательно, f является изоморфным отображением $H_n(G+v)$ на $H_n(G)$, то есть $H_n(G+v)$ и $H_n(G)$ изоморфны. Теорема доказана.□

Докажем аналогичную теорему для приклеивания и отбрасывания точечной связи в произвольном пространстве G.

Теорема 134

*Пусть G есть пространство и H есть его подпространство, для которого группы гомологий есть $H_0(H)=A$, $H_n(H)=0$, где $n>0$. Пусть H есть общий окаем двух несмежных точек $v_1$ и $v_2$, $O(v_1v_2)=H$. Тогда приклеивание связи $(v_1v_2)$ к G по H не меняет групп гомологий пространства G, $H_n(G)=H_n(G+(v_1v_2))$.*

Доказательство

При приклеивании связи $(v_1v_2)$ к G количество связных компонент не меняется. Следовательно, $H_0(G)=H_0(G+(v_1v_2))$.

Пусть G и H являются искомыми пространством и некоторым его подпространством, не совпадающим с G. Приклеим к G по H связь $(v_1v_2)$ так, чтобы $O(v_1v_2)=H$. Пусть $U=G+(v_1v_2)$. Из $c_n\in C_n(U=G+(v_1v_2))$ и $n>2$, следует, что $c_n=a_n+v_1\oplus v_2\oplus b_{n-2}$, где $a_n\in C_n(G)$ и $b_{n-2}\in C_{n-2}(H)$. Используя доказательство предыдущей теоремы мы получаем, что любой цикл $z_n$ в U может быть представлен в виде $z_n=\partial(v_1\oplus v_2\oplus b_{n-1})+w_n$, где $w_n$ есть цикл в G. Введем гомоморфизм f: $Z_n(U)\rightarrow Z_n(G)$,

$f(z_n)=f(w_n+\partial(v_1\oplus v_2\oplus b_{n-1}))=w_n$.

Очевидно, что любой цикл на U переходит в цикл на G, и для любого цикла на G существует прообраз, являющийся циклом на U. Очевидно, также, что любая граница из U переходит в границу из G, и для любой границы из G существует прообраз, являющийся границей в U, как это видно из следующего:

$f(\partial(a_{n+1}+v_1\oplus v_2\oplus b_{n-1})=f(\partial a_{n+1}+\partial(v_1\oplus v_2\oplus b_{n-1}))=\partial a_{n+1}$.

Следовательно, f является изоморфным отображением $H_n(U)$ на $H_n(G)$, то есть $H_n(U)$ и $H_n(G)$ изоморфны, n>1.

Пусть $c_2$ есть некоторая 2-цепь в $C_2(U=G+(v_1v_2))$. Из $c_2\in C_2(U=G+(v_1v_2))$, следует, что $c_2=a_2+av_1\oplus v_2$, где $a_2\in C_2(G)$. Из предыдущей теоремы следует, что $z_2=\partial(av_1\oplus v_2\oplus v_3)+w_2$, где $w_2$ есть цикл в G. Введем гомоморфизм f: $Z_2(U)\rightarrow Z_2(G)$,

$f(z_2)=f(\partial(av_1\oplus v_2\oplus v_3)+w_2)=w_2$.

Очевидно, что любой 2-цикл на U переходит в 2-цикл на G, и для любого

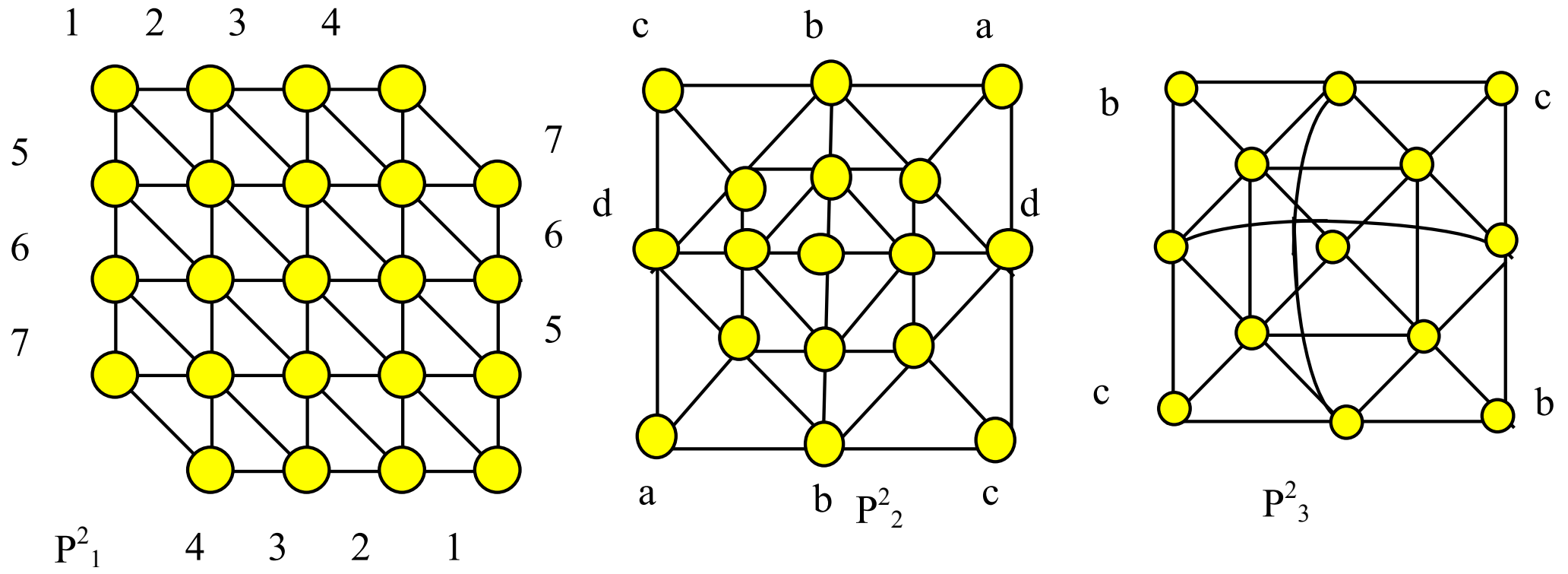

Рис. 145 Молекулярные модели проективной плоскости. Группы гомологий $H_0=Z$, $H_1=Z_2$, $H_n=0$, n >1.

2-цикла на G существует прообраз, являющийся 2-циклом на U. Очевидно, также, что любая 2-граница из U переходит в 2-границу из G, и для любой 2-границы из G существует прообраз, являющийся 2-границей в U, как это видно из следующего:

$f(\partial(a_2+av_1\oplus v_2\oplus b_3)=f(\partial a_3+\partial(av_1\oplus v_2\oplus b_3))=\partial a_3$.

Следовательно, f является изоморфным отображением $H_1(U)$ на $H_1(G)$, то есть $H_1(U)$ и $H_1(G)$ изоморфны. Теорема доказана.□

Сформулируем в заключение основную теорему этого раздела, доказывать которую нет необходимости, поскольку она является результатом предыдущих теорем.

Теорема 135

*Точечные преобразования пространства G не меняют групп гомологий этого пространства.*

Рассмотрим примеры групп гомологий некоторых молекулярных

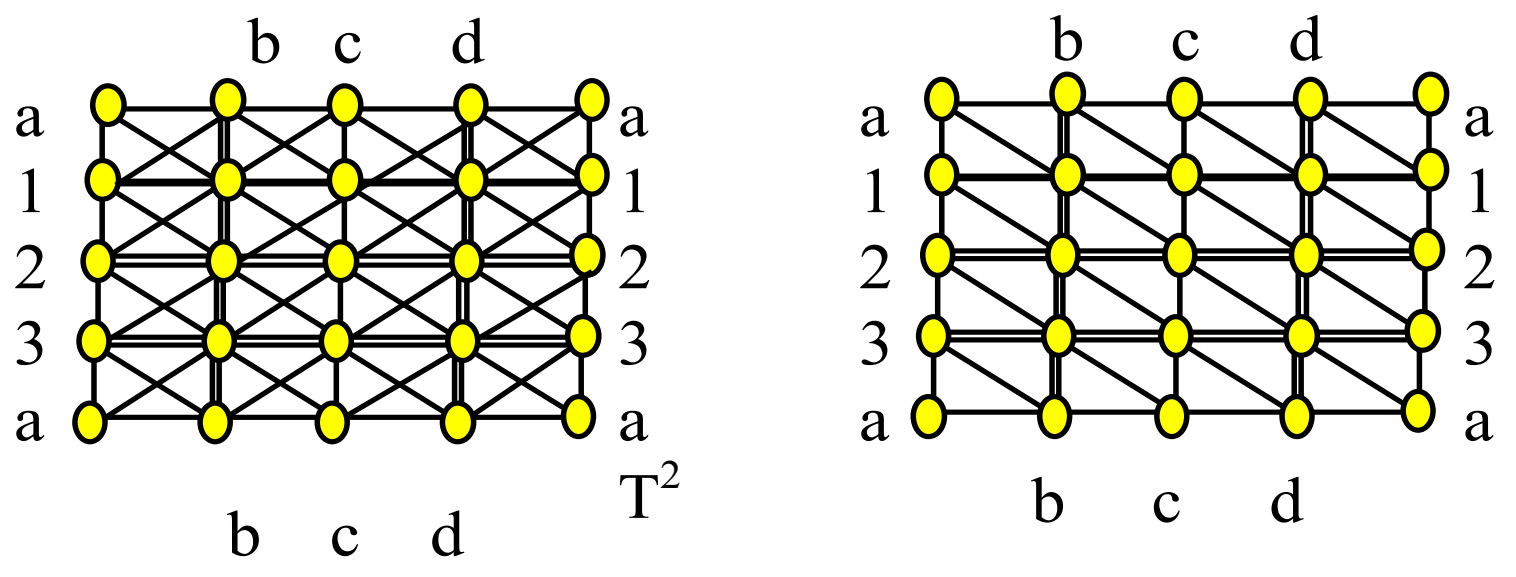

Рис. 146 Пространства, являющиеся тором, $T^2$ с группами гомологий $H_0$=Z, $H_1=Z\otimes Z$, $H_2$=Z, $H_n$=0, n > 2.

пространств.

Примеры

Пусть группа коэффициентов является группой Z целых чисел, Z=(0, ± 1, ± 2, ± 3,..)

**Точечные МП**

Все точечные молекулярные пространства стягиваются к одной точке и имеют соответствующие группы гомологий: $H_0$=Z, $H_n$=0, n > 0. На Рис. 141 изображены некоторые из точечных МП.

**Одномерная окружность $S^1$.**

Все пространства, изображенные на Рис. 142, стягиваются к одномерной окружности и имеют соответствующие группы гомологий: $H_0=H_1$=Z, $H_n$=0, n > 1.

**Двумерная сфера $S^2$**

На Рис. 143 изображена двумерная сфера и гомотопное ему пространство. Их группы гомологий определяются выражениями $H_0=H_2$=Z, $H_1=H_n$=0, n > 2.

### Трехмерная сфера $S^3$

На Рис. 144 изображены МП, гомотопные трехмерной сфере. Их группы гомологий определяются выражениями

$H_0=H_3=Z$, $H_1=H_2=H_n=0$, $n > 3$.

### Двумерный тор

Группы гомологий двумерного тора (Рис. 134) и гомотопных ему пространств, так же, как и в непрерывном случае, определяются выражениями $H_0=H_2=Z$, $H_1=Z\otimes Z$, $H_n=0$, $n>3$.

### Проективная плоскость

Аналогично можно убедиться прямым подсчетом, что группы гомологий проективной плоскости определяются соотношениями $H_0=Z$, $H_1=Z_2$, $H_n=0$, $n >1$, (Рис. 145).

### Бутылка Клейна

Рассмотрим группы гомологий бутылки Клейна (Рис. 137). Легко убедиться непосредственным подсчетом, что ее группы гомологий имеют вид $H_0=Z$, $H_1=Z\otimes Z_2$, $H_n=0$, $n >1$.

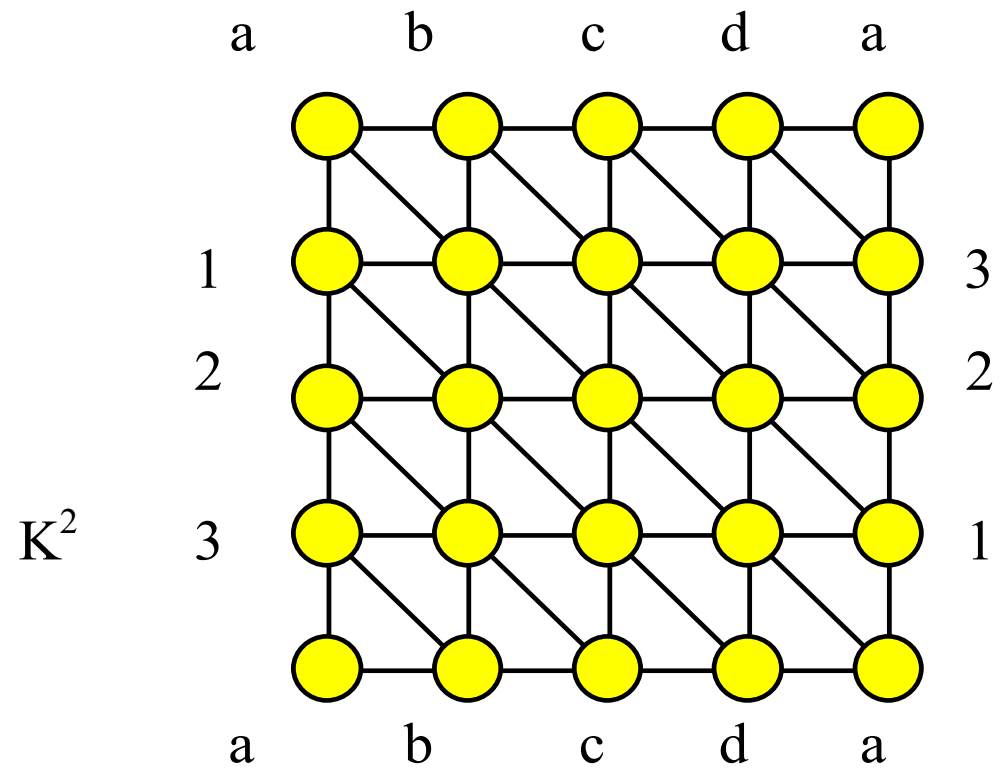

Рис. 147 Группы гомологий бутылки Клейна имеют вид $H_0=Z$, $H_1=Z\otimes Z_2$, $H_n=0$, $n >1$.

Здесь возникает несколько интересных математических направлений для исследований.

З а д а ч а .

Изучить свойства групп гомологий на молекулярных пространствах. Например, как определяются группы гомологий при прямой сумме и прямом произведении молекулярных пространств.

З а д а ч а .

Найти группы гомологий и исследовать их свойства для пространств, не являющихся многообразиями.

З а д а ч а .

Найти метод вычисления групп гомологий и эйлеровой характеристики, по:
Матрице инциденций молекулярного пространства,
Координатной матрице молекулярного пространства.

С п и с о к л и т е р а т у р ы к г л а в е 1 1 .

# ЧАСТЬ 2
# МОЛЕКУЛЯРНОЕ МОДЕЛИРОВАНИЕ

## КИРПИЧНЫЕ (ГЕОМЕТРИЧЕСКИЕ) И МАТРИЧНЫЕ (АЛГЕБРАИЧЕСКИЕ) МОДЕЛИ МОЛЕКУЛЯРНЫХ ПРОСТРАНСТВ

We consider the problem how to find a continuous counterpart for a molecular space. To universalize this procedure we take infinite-dimensional Euclidean space $E^{\infty}$ and define a kirpich which is infinite-dimensional unit cube in $E^{\infty}$. The kirpich is a basic element in this approach. A set of kirpiches forms a universal basic space (UBS). Any finite set of kirpiches in UBS is called a surface S in UBS. The nerve of S is called the molecular space M(S) of S. A matrix that is formed by coordinates of all kirpiches of S is called the coordinate or matrix representation of S. We investigate some properties of S and regard some geometric and visual background for introduction of UBS. The way to obtain a continuous counterpart of a molecular space is the following:
Given a molecular space N we build any S for which N=M(S). Using methods of homotopic topology we can construct a continuous surface in continuous Euclidean $E^n$ that is homotopic to M(S).

Результаты этой главы частично изложены в работах [3,7,14,26,27,43,45,46,49,51].

### *БЕСКОНЕЧНО-МЕРНОЕ ЕВКЛИДОВО ПРОСТРАНСТВО КАК УНИВЕРСАЛЬНОЕ БАЗОВОЕ ПРОСТРАНСТВО*

Некоторые могут посчитать, что молекулярное пространство является достаточно абстрактным инструментом, трудным для геометрической или алгебраической его реализации. На самом деле существует несколько моделей как геометрического так и алгебраического характера, которые являются эквивалентами молекулярных пространств и, возможно, окажутся более удобными для понимания. Геометрическая модель молекулярного пространства также является, по сути дела, непрерывным объектом, который, безусловно, более похож на поверхность, чем само молекулярное пространство. Если мы имеем непрерывную поверхность и строим ее молекулярную модель, мы заведомо огрубляем и упрощаем конструкцию, теряя при этом какую-то часть информации. Обратный путь нахождения непрерывной поверхности по заданному молекулярному пространству не может дать новой информации кроме той, которая уже имеется в самом молекулярном пространстве. С другой стороны, мы стоим перед выбором, какую из непрерывных поверхностей считать достоверной моделью молекулярного пространства. Мы рассмотрим путь, который позволяет свести к минимуму всякого рода неоднозначность и произвол при таких построениях.

Введем понятие универсального кирпичного пространства, к которому мы можем применить гомотопные преобразования для нахождения непрерывной поверхности, гомотопной данному кирпичному пространству. Это достаточно разработанная задача в классической топологии. Так как в молекулярном пространстве все точки обезличены, взаимозаменяемы и одинаковы, в кирпичном пространстве все элементы должны быть идентичны. Переходим к описанию универсального кирпичного пространства и других понятий, связанных с ним.

Определение бесконечномерного евклидова пространства $E^{\infty}$

Бесконечномерным евклидовым пространством $E^{\infty}$ называется множество точек x, $x \in E^{\infty}$, координаты которых являются бесконечной последовательностью действительных чисел $x_k$,

$$x = (x_1, x_2, \ldots x_k, \ldots) = [\, x_k \,],\ k \in N.$$

Определение кирпича.

Кирпичом (kirpich) К в $E^{\infty}$ называется множество точек $x \in K$, координаты x которых удовлетворяют соотношению:

$$n_k \le x_k \le n_k + 1,\ k \in N,\ n_k \text{ - целые числа.}$$

Определение координат кирпича.

Если $x = [x_k]$, $x \in K$, $n_k \le x_k \le n_k + 1$, $k \in N$, $n_k$ - целые числа, то последовательность левых границ неравенств называется координатами кирпича К в $E^{\infty}$, $K = (n_1, n_2, \ldots n_k, \ldots) = [n_k]$, $k \in N$. ¤

Очевидно, что кирпич есть не что иное, как бесконечно-мерный единичный куб в $E^{\infty}$. Термин kirpich уже использовался нами в предыдущих статьях. В силу своего определения. сам по себе кирпич является безразмерным элементом, в отличие от, например, точки, которая является нульмерным элементом в классической математике. Такой подход позволяет использовать кирпич как универсальный строительный элемент для создания пространств различных размерностей. Введем расстояние между кирпичами.

Определение расстояния между кирпичами.

Расстоянием $D(K_1,K_2)$ между кирпичами $K_1 = [n_k]$ и $K_2 = [m_k]$, $k \in N$, называется наибольшая из разностей между их соответствующими координатами, взятая по модулю.

$$D(K_1,K_2) = \max|n_k - m_k|, k \in N. ¤$$

Определение смежных кирпичей.

Два кирпича называются смежными или соседними, если расстояние между ними равно единице. Связью называется любая пара смежных кирпичей.

З а м е ч а н и е

Если кирпичи имеют общие точки, то они всегда смежны. Некоторая неопределенность может возникнуть, если каждая координата одного кирпича отличается от каждой координаты другого кирпича на единицу. Например, $K_1$=(0,0,0,0,...0,...) и $K_2$=(1,1,1,1,...1,...). В этом

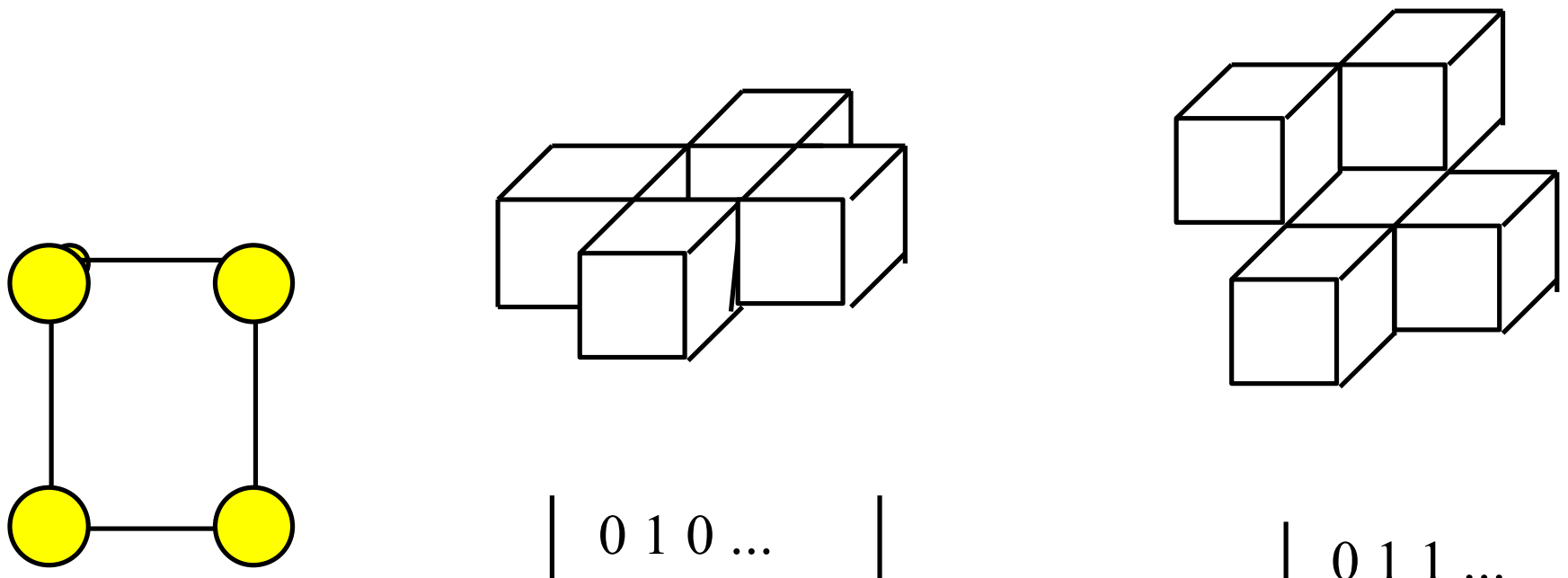

| | | | |
|---|---|---|---|
| 0 | 1 | 0 | ... |
| 1 | 2 | 0 | .... |
| 2 | 1 | 0 | .... |
| 1 | 0 | 0 | .... |
| 1 | 0 | 0 | ... |
| 0 | 1 | 1 | ... |
| 1 | 2 | 1 | .... |
| 2 | 1 | 0 | .... |
| 1 | 0 | 0 | .... |
| 1 | 0 | 0 | ... |

Рис. 148 Молекулярная минимальная окружность, состоящая из четырех точек, представлена кирпичами и координатными матрицами в двумерном и трехмерном пространствах.

случае будем считать по определению, что кирпичи имеют общие точки.

Определение универсального кирпичного пространства $B(E^\infty)$.

Универсальным кирпичным пространством $B(E^\infty)$, сокращенно УБП, называется множество всех кирпичей в $E^\infty$.

Поскольку кирпичное пространство и молекулярное пространство по сути дела есть различные представления дискретного пространства, мы дадим определение универсального молекулярного пространства.

Определение универсального молекулярного пространства $M(E^{\infty})$.

Универсальным молекулярным пространством $M(E^{\infty})$ называется нерв универсального кирпичного пространства $B(E^{\infty})$.

Это означает, что каждому кирпичу из $B(E^{\infty})$ ставится в соответствие точка из $M(E^{\infty})$ с теми же координатами, что и кирпич, и кирпичи в $B(E^{\infty})$ смежны тогда и только тогда, когда смежны точки в $M(E^{\infty})$. При таком подходе различие между кирпичным и молекулярным пространствами исчезает, и все выводы, сделанные относительно кирпичных пространств, применимы к молекулярным пространствам.

Определение кирпичного n-мерного подпространства $B(E^{n})$.

n-Мерным кирпичным подпространством $B(E^{n})$ универсального базового пространства $B(E^{\infty})$ называется множество всех кирпичей в $E^{\infty}$, у которых все координаты, начиная с (n+1) и выше равны 0, $B(E^{n})\subseteq B(E^{\infty})$.

Определение молекулярного n-мерного евклидова подпространства $M(E^{n})$.

n-Мерным молекулярным подпространством $M(E^{n})$ универсального базового пространства $M(E^{\infty})$ называется множество всех точек в $E^{\infty}$, у которых все координаты, начиная с (n+1) и выше равны 0, $M(E^{n})\subseteq M(E^{\infty})$.

Определение кирпичного пространства или поверхности S в $B(E^{\infty})$.

Кирпичным пространством или поверхностью S в $B(E^{\infty})$, называется любое конечное или счетное множество кирпичей.

Для любого кирпичного пространства S строим молекулярное пространство M(S), являющееся нервом множества кирпичей. Каждому кирпичу ставим в соответствие точку молекулярного пространства и если два кирпича являются смежными, то соответствующие им точки молекулярного пространства являются смежными. Это отображение множества кирпичей на множество молекулярных пространств, очевидно является взаимно-однозначным. При этом универсальное базовое пространство $B(E^{\infty})$ также является кирпичным пространством, для которого существует универсальное базовое молекулярное пространство $M(E^{\infty})$. Все выводы сделанные относительно кирпичных пространств, будут справедливы для молекулярных пространств и наоборот.

Определение координатной матрицы A(S) кирпичной поверхности S (и пространства M(S)).

Пусть поверхность S состоит из кирпичей ( $K_1,K_2,\ldots K_p$) с координатами $K_1=(n_{11},n_{12},\ldots n_{1m},\ldots)=[n_{1m}]$, $K_2=(n_{21},n_{22},\ldots n_{2m},\ldots)=[n_{2m}],\ldots K_p=(n_{p1},n_{p2},\ldots n_{pm},\ldots)=[n_{pm}]$, Координатной матрицей A(S) поверхности S и молекулярного пространства M(S) называется матрица $[n_{tm}]$, имеющая p строк и бесконечное множество столбцов (Рис. 148).

Определение вложенной размерности vr(S) поверхности S в $B(E^{\infty})$.

Вложенной размерностью поверхности S в $B(E\infty)$ называется число отличных от нуля координат в множестве кирпичей поверхности.

На Рис. 149 поверхности, состоящие кирпичей, находятся в трехмерном пространстве. Следовательно, их вложенные размерности равны трем. Таким образом с этого момента, как и классическом случае для каждого молекулярного пространства мы имеем внутреннюю размерность пространства, определяемую свойствами окаема, и вложенную размерность, определяемую свойствами пространства вцелом при его

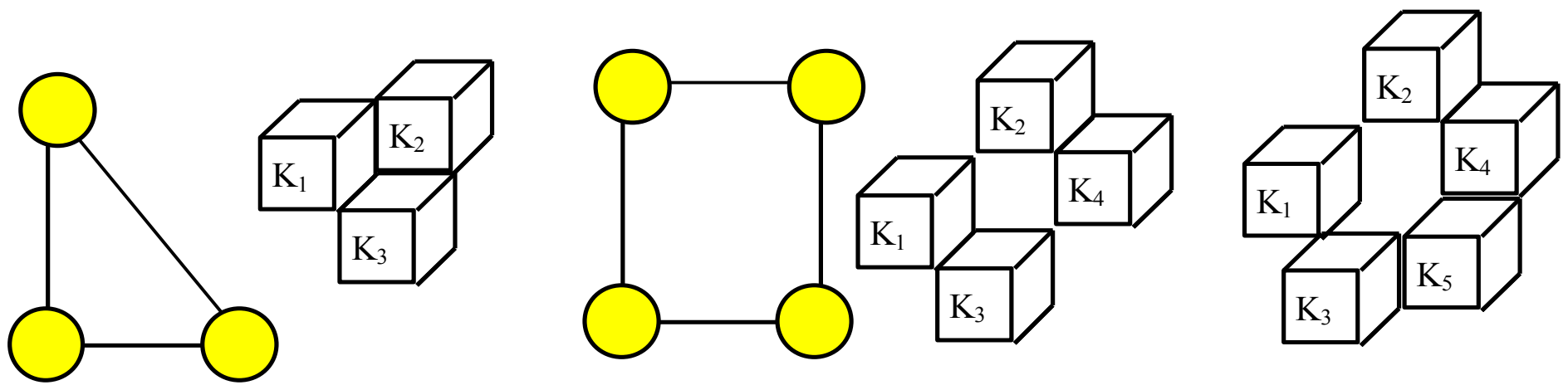

Рис. 149 Невозможно расположить три кирпича таким образом, чтобы между ними была “дырка”. Для четырех и более кирпичей такое возможно,

вложении во внешнее пространство.

Теорема 136

*Для каждого кирпичного пространства S существует изоморфное ему молекулярное пространство M(S).*

*Для каждого кирпичного пространства S в пространстве $B(E^{\infty})$ существует изоморфное ему молекулярное пространство M(S) в $M(E^{\infty})$.*

Доказательство.

Установим взаимно-однозначное соответствие следующим образом: каждому кирпичу кирпичного пространства поставим в соответствие точку молекулярного пространства и две точки будут смежными тогда и только тогда, когда два кирпича являются соседними. Доказательство закончено.□

Использование такого подхода имеет то преимущество, что позволяет геометрически интерпретировать многие свойства молекулярных пространств, а также использовать хорошо развитый координатный аппарат для изучения молекулярных пространств. Например, становится наглядным, почему нормальная молекулярная окружность не может иметь менее четырех точек. Невозможно расположить три кирпича таким образом, чтобы между ними была “дырка”. Для четырех и более кирпичей такое возможно, как это видно на Рис. 149. Очевидно, что для различных поверхностей в $B(E^\infty)$ может быть один и тот же нерв. Это означает, что поверхности эквивалентны и отличаются только расположением в пространстве $B(E^\infty)$ или $E^\infty$. С другой стороны, для данной поверхности существует единственное (с точностью до изоморфизма) молекулярное пространство.

Определение объема V(S) или |S| поверхности S

Объемом V(S) поверхности S (или молекулярного пространства этой поверхности) называется число кирпичей (точек молекулярного пространства), составляющих эту поверхность (молекулярное пространство).

## *СВОЙСТВА ПОВЕРХНОСТЕЙ В УНИВЕРСАЛЬНОМ КИРПИЧНОМ ПРОСТРАНСТВЕ*

Мы ввели матричное координатное представление поверхности S. Это является удобным математическим инструментом, позволяющим исследовать молекулярные пространства и кирпичные поверхности. Следует иметь ввиду, что все, что говорится о поверхности, относится также и к молекулярному пространству этой поверхности. Такое описание поверхности, или молекулярного пространства, позволяет представить эту поверхность в виде прямоугольной матрицы. Рассмотрим пример. Легко проверить, что эта окружность, состоящая из пяти точек, не может быть представлена как семейство единичных квадратов в базовом кирпичном двумерном пространстве $B(E^2)$. Используя трехмерное пространство, мы можем представить эту окружность как семейство единичных кубов в трехмерном пространстве (Рис. 149). Координаты этих кубов образуют матрицу окружности.

Перечислим очевидные свойства координатной матрицы.

Свойства координатной матрицы.

- Если два кирпича смежны, то их одноименные координаты отличаются не более чем на 1. $K_1=(n_{11},n_{12},...n_{1m},...)=[n_{1m}]$, $K_2=(n_{21},n_{22},...n_{2m},...)=[n_{2m}]$, $\Rightarrow |\ n_{1m}-n_{2m}|\le 1$.
- Если два кирпича несмежны, то хотя бы одна пара одноименных координат отличается более чем на 1. $K_1=(n_{11},n_{12},...n_{1m},...)=[n_{1m}]$, $K_2=(n_{21},n_{22},...n_{2m},...)=[n_{2m}]$, $\Rightarrow \exists\ s, |\ n_{1s}-n_{2s}|\ge 2$.

Введем понятие кирпичности. Аналогичные понятия используются в классической математике и, в частности, теории графов. Для всякого молекулярного пространства можно ввести кирпичное пространство различными способами, например, его можно построить так, что его вложенная размерность будет различна. Однако, всегда существует минимальная вложенная размерность. Например, для минимальной нормальной окружности (Рис. 148) можно построить кирпичное

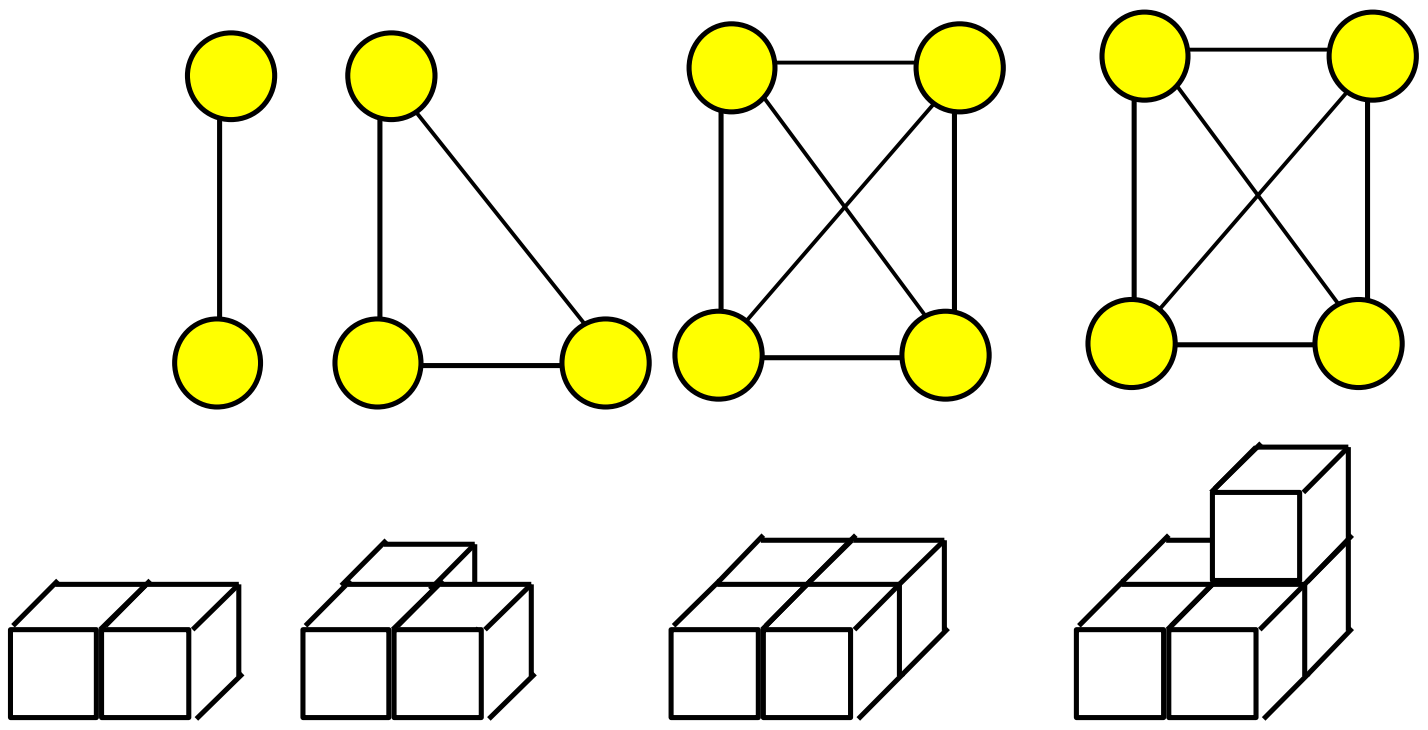

Рис. 150 Укладка кирпичных связок из 2, 3, 4 и 5 кирпичей $kir(K(n))=\lceil \log_2 n\rceil$, где $\lceil \log_2 n\rceil$ есть наименьшее целое число, не меньшее $\log_2 n$.

пространство с вложенной размерностью 2, 3 и выше. Задачи, связанные с минимизацией вложенной размерности возникают, например, в теории графов при рассмотрении биологических задач [20,21].

<u>Определение кирпичности (kirpicity) kir(S) поверхности S в $B(E^\infty)$.</u>

Кирпичностью kir(S) поверхности S в $B(E^\infty)$ называется минимальная вложенная размерность этой поверхности [27,46].

З а м е ч а н и е

Следует отметить, что любое конечномерное пространство можно считать подпространством бесконечномерного пространства $E^\infty$. При этом все отсутствующие координаты можно считать равными нулю. Фактически, мы всегда используем конечномерное пространство. Логическим обоснованием формального введения бесконечномерного пространства, как мы уже говорили, является желание создать некий универсальный элемент, не зависящий от размерности и пригодный для описания элементов различных размерностей: нульмерных точек, одномерных линий, двумерных поверхностей и так далее. С этой точки зрения поверхности любой размерности являются некоторыми семействами одних и тех же элементов, и отличаются только взаиморасположением этих элементов. При таком подходе, в сущности, теряется различие между кирпичным пространством и молекулярным пространством. Можно считать, что молекулярное пространство, кирпичное пространство и координатная матрица есть молекулярное (графическое), геометрическое и алгебраическое (матричное) представления одного и того же объекта - дискретного пространства. При этом молекулярное пространство является наиболее абстрактным представлением, так как не требует вложения в другое пространство. На Рис. 148 изображены кирпичные поверхности и их матрицы для одного и того же молекулярного пространства, являющегося минимальной окружностью. Оценим размерность евклидова

$$A=\begin{vmatrix} 2 & & & 1 \\ B & 2 & \ddots & \\ & & & 2 \\ B & B & & B \end{vmatrix}$$

Рис. 151 Координатная матрица имеет (p-1) столбцов Числа B являются нулями или единицами.

пространства, в котором заведомо можно построить кирпичное пространство для данной поверхности.

Теорема 137

*Пусть S является кирпичным пространством объема V(S)=p. Тогда кирпичность этого пространства n не превышает p-1, kir(S)≤ p-1 (Рис. 150).*

Доказательство.

Обозначим кирпичи поверхности S как $a_1$, $a_2$,...$a_p$. Построим координатную матрицу $a_{ks}$ этой поверхности в которой число столбцов, содержащих ненулевые элементы, определяет размерность вложенного пространства. Выберем координаты $a_{kk}=2$, где k=1,2,...p-1. Выберем координаты $a_{ks}=1$, где $k<s<p$, если $a_k$ и $a_s$ смежны. Выберем координаты $a_{ks}=0$, где $k>s$, $s<p=1,2,...p-1$, если $a_k$ и $a_s$ несмежны. Очевидно, что мы получаем координатную матрицу A(S), имеющую p строк и (p-1) столбцов с различными координатами. Все остальные координаты кирпичей выберем равными 0 (Рис. 150). Следовательно, kir(S)≤ p-1. Доказательство закончено.□

Теорема 138

*В любой конечной поверхности S, два любые смежные кирпичи можно сделать несмежными, не нарушая смежности остальных кирпичей.*

Доказательство.

Обозначим кирпичи поверхности S как $a_1$, $a_2$,...$a_p$. Пусть кирпичи имеют n различных координат, остальные координаты совпадают. Пусть два смежных кирпича $a_r$ и $a_q$ имеют m первых различающихся на 1 координат, где m≤n, и определяются координатами $a_r=(c_1,c_2,...c_n,t,t,..)$, $a_q=(d_1,d_2,...d_n,t,t,..)$. Увеличим (n+1) координату кирпича $a_r$ на 1 и уменьшим (n+1) координату кирпича $a_r$ на 1. То есть $a_r=(c_1,c_2,...c_n,t+1,t,..)$, $a_q=(d_1,d_2,...d_n,t-1,t,..)$. Данные кирпичи станут несмежными, смежность с остальными не изменится. Теорема доказана. □

Теорема 139

*В любой конечной поверхности S два любые несмежных кирпича можно сделать смежными, не нарушая смежности остальных кирпичей.*

Доказательство.

Обозначим кирпичи поверхности S как $a_1$, $a_2$,...$a_p$. Пусть кирпичи имеют n различных координат, остальные координаты совпадают. Пусть два несмежных кирпича $a_r$ и $a_q$ имеют m первых различающихся координат, где m≤n, и определяются координатами $a_r=(c_1,c_2,...c_n,t,t,..)$, $a_q=(d_1,d_2,...d_n,t,t,..)$. Пусть $c_1<d_1$. У кирпичей, имеющих первую координату, меньшую или равную $c_1$, увеличиваем эту координату на 1, одновременно увеличивая первую общую координату t на 1. При этом могут появиться новые смежные кирпичи, которые мы делаем несмежными. Повторяя эти операции нужное число раз делаем несмежные кирпичи смежными, не меняя смежности остальных кирпичей. Теорема доказана. □

Введем определение кирпичной связки, аналогичное молекулярной связке.

<u>Определение кирпичной связки K(n).</u>

Поверхность S, состоящая из n кирпичей, называется кирпичной связкой, если каждая пара кирпичей является смежной.

Рассмотрим кирпичность связки, содержащей n кирпичей.

Теорема 140

*Кирпичность связки K(n) определяется соотношением $\log_2 n \le kir(K(n)) < \log_2 n+1$, или $kir(K(n))=\lceil \log_2 n \rceil$, где $\lceil \log_2 n \rceil$ есть наименьшее целое число, не меньшее $\log_2 n$.*

Доказательство.

Пусть некоторый кирпич связки имеет нулевые координаты $a_1=(0,0,...0,..)$. Тогда координаты всех остальных кирпичей должны состоять из 0 и 1 (или -1). Количество наборов длины s, состоящих из 0 и 1, равно $2^s$. Следовательно, $2^s \le n < 2^{s+1}$. Логарифмируя это выражения мы имеем $s \le \log_2 n < s+1$. Следовательно, $kir(K(n))=\lceil \log_2 n \rceil$, где $\lceil \log_2 n \rceil$ есть наименьшее целое число, не меньшее $\log_2 n$ (Рис. 150). Теорема доказана. □

Докажем еще одну формулу, определяющую кирпичность поверхности.

Теорема 141

*Пусть S является поверхностью, имеющей n кирпичей и к связей. Тогда существует базовое n-мерное кирпичное пространство $B(E^n)$, $S \subseteq B(E^n)$, в пространстве $B(E^\infty)$, такое, что $kir(S) \le C^2_{n-k} + \lceil \log_2 n \rceil$, где $\lceil \log_2 n \rceil$ является наименьшим целым числом, не меньшим, чем $\log_2 n$.*

Доказательство.

В любой кирпичной связке из n кирпичей содержится $C^2_n$ связок из двух кирпичей. При раздвигании двух смежных кирпичей вложенная размерность пространства может увеличиться на 1. Число таких раздвиганий равно ($C^2_{n-k}$). Отсюда, вложенная размерность пространства S не может превышать вложенную размерность кирпичной связки из n кирпичей более чем на ($C^2_n$-k). Следовательно, $kir(S) \le C^2_n - k + \lceil \log_2 n \rceil$. Доказательство закончено.□

Очевидно, что для кирпичного пространства существует два различных вида преобразований. Во-первых, это преобразования, которые не меняют соседства между кирпичами, а меняют только степень этого соседства. При таких преобразованиях два кирпича могут касаться один другого по граням, ребрам или вершинам, это не важно, а важно то, что они всегда остаются соседними (или всегда несоседними). Расположение кирпичной поверхности в евклидовом пространстве меняется, однако, молекулярное пространство остается неизменным. Используя аналогии, такие преобразования напоминают изгибания двумерных поверхностей в пространстве трех и более измерений. Такие преобразования можно назвать гладкими преобразованиями кирпичной поверхности. Второй тип преобразований является прямым отражением точечных преобразований молекулярного пространства. В отличие от первого типа, такие преобразования не могут быть названы гладкими. Приклеивание (или отбрасывание) точки в молекулярном пространстве эквивалентно приклеиванию (или отклеиванию) к кирпичной поверхности кирпича. Приклеивание (или отбрасывание) связи в молекулярном пространстве эквивалентно склеиванию (или расклеиванию) двух кирпичей поверхности. В классической математике такие преобразования, скорее всего, соответствуют гомотопному переходу от одного пространства к другому. Вопрос о геометрии кирпичных пространств как подпространств бесконечно-мерного пространства почти совершенно не исследован, хотя некоторые попытки были сделаны.

Задача

Ввести геометрические характеристики молекулярных и/или кирпичных пространств и исследовать их свойства.

## *УНИВЕРСАЛЬНЫЙ КИРПИЧНЫЙ КУБ И ЕГО СВОЙСТВА*

Рассмотрим еще одно свойство кирпичных поверхностей, называемое компактностью (compactness).

Определение универсального кирпичного n-мерного куба $U^n$ и универсального молекулярного n-мерного куба (Рис. 152).

Универсальным кирпичным n-мерным кубом $U^n$ называется поверхность, состоящая из ($a_1, a_2, \ldots a_p$), $p=3^n$ кирпичей, координатная матрица которой $A(U^n)$ с отличными от нуля столбцами является $3^n$x$n$ матрицей, имеющей $3^n$ строк и n столбцов с координатами, состоящими из 0, 1 и 2,

$$A(U^n) = \begin{vmatrix} 0 & \ldots & 0 & 0 \\ 0 & \ldots 0 \ldots & 0 & 1 \\ 1 & 2 & 2 & 2 \\ 2 & 2 & 2 & 2 \end{vmatrix}.$$

Универсальным молекулярным n-мерным кубом называется нерв (или граф пересечений) универсального кирпичного n-мерного куба $U^n$.

Важной особенностью универсального молекулярного n-куба $U^n$ является то, что он является прямым произведением n копий одномерного нормального отрезка из 3-х точек.

Теорема 142

*Универсальный молекулярный n-куб $U^n$ является прямым произведением n копий одномерного нормального отрезка из 3-х точек,* $U^n = U^1 \otimes U^1 \otimes U^1 \ldots \otimes U^1$.

Доказательство.

Теорема очевидна из Рис. 153. Для произвольной размерности обозначим вершины одномерного универсального куба $U^1$ числами 0, 1 и 2. Рассмотрим $U^n=U^1\otimes U^1\otimes U^1...\otimes U^1$. Легко видеть, что точки в прямом произведении будут смежны только если координаты отличаются не более, чем на 1. Таким образом, координатная матрица будет матрицей

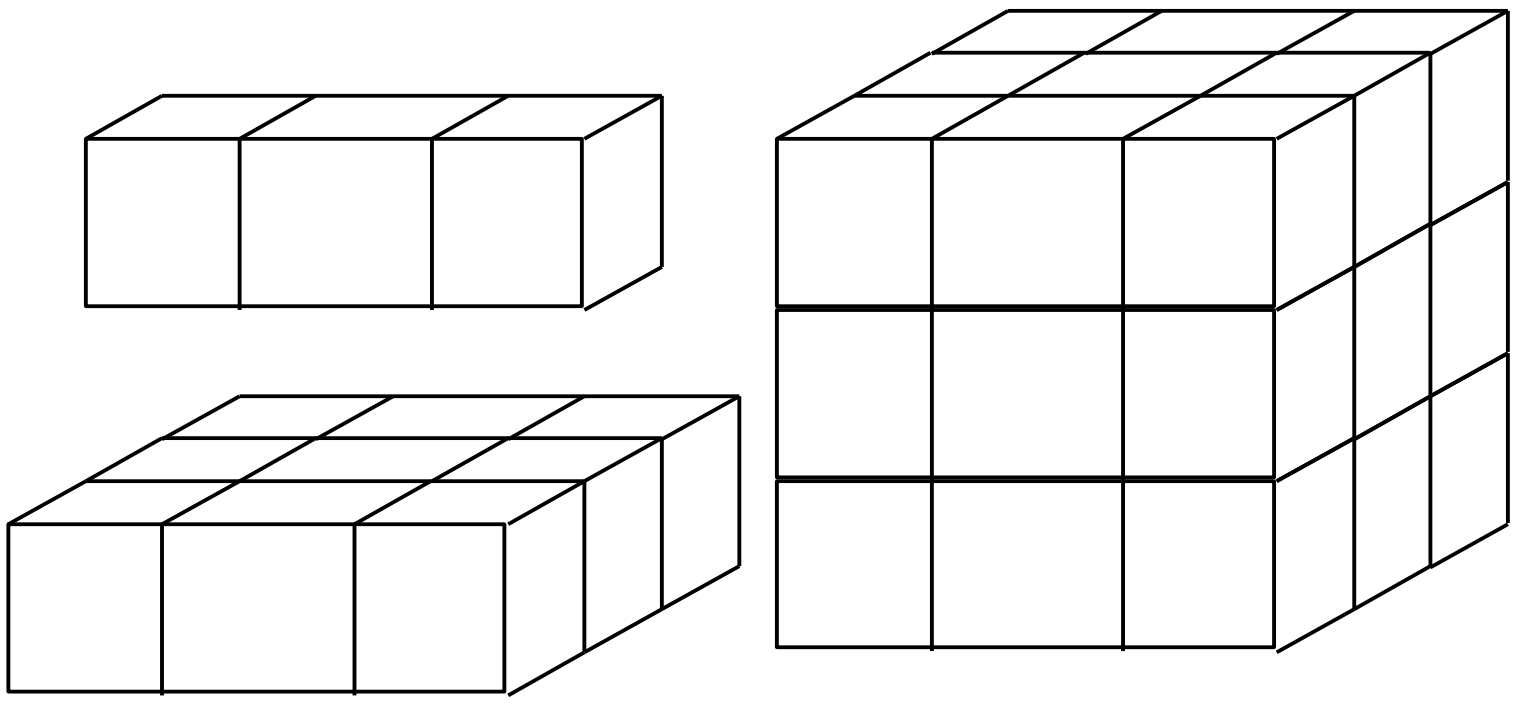

Рис. 152 Одномерный, двумерный и трехмерный универсальные кубы $U^1$, $U^2$ и $U^3$.

универсального n-мерного куба. Теорема доказана..□

Универсальный куб обладает одним важным свойством: он вмещает в себя любое молекулярное пространство. Кроме того, он имеет радиус 1, то есть любая точка универсального куба является смежной центральной точке. Универсальный куб является n-мерным шаром центральной точки v, $U^n=U^n(v)$.

Определение компактности (compactness) com(S) поверхности S в $B(E^\infty)$.

Компактностью com(S) поверхности S в $B(E^\infty)$ называется минимальное целое число n, такое, что поверхность S является подпространством универсального куба $U^n$ [27], $com(S)=min(n, S\subseteq U^n)$.

З а м е ч а н и е

Очевидно, что кирпичность не превышает компактность, $kir(S)\le com(S)$. Например, пространство, состоящее из 4-х несмежных кирпичей имеет кирпичность, равную 1 и компактность, равную 2 (Рис. 155).

Замечание

Компактность подразумевает, что в любой поверхности, вложенной в универсальный куб, расстояние между двумя несмежными кирпичами всегда равно 2. Иными словами, если кирпичи смежны, то расстояние между ними равно 1, если же они несмежны, то расстояние между ними в универсальном кирпичном пространстве равно 2. То есть любая поверхность может быть так изогнута и вложена во внешнее пространство, что любые две несмежные точки находятся на минимальном возможном расстоянии одна от другой. Хотя нечто подобное можно рассматривать и в непрерывных пространствах, там это не имеет большого смысла, поскольку не существует универсального понятия объема поверхности, пригодного для поверхностей любой размерности. Например, двумерная сфера в трехмерном пространстве имеет объем, равный 0, и, следовательно, любые две точки могут находиться на бесконечно близком расстоянии одна от другой. В этом плане объекты классической математики весьма далеки от повседневной практики и являются некими идеализированными объектами, что отчетливо обнажилось в компьютерной практике.

Теорема 143

*Пусть S является кирпичным пространством объема V(S)=p. Тогда компактность этого пространства n не превышает p-1, com(S)≤ p-1.*

Доказательство.

Доказательство слово в слово повторяет доказательство соответствующей теоремы о кирпичности. Координатная матрица A(S) поверхности S

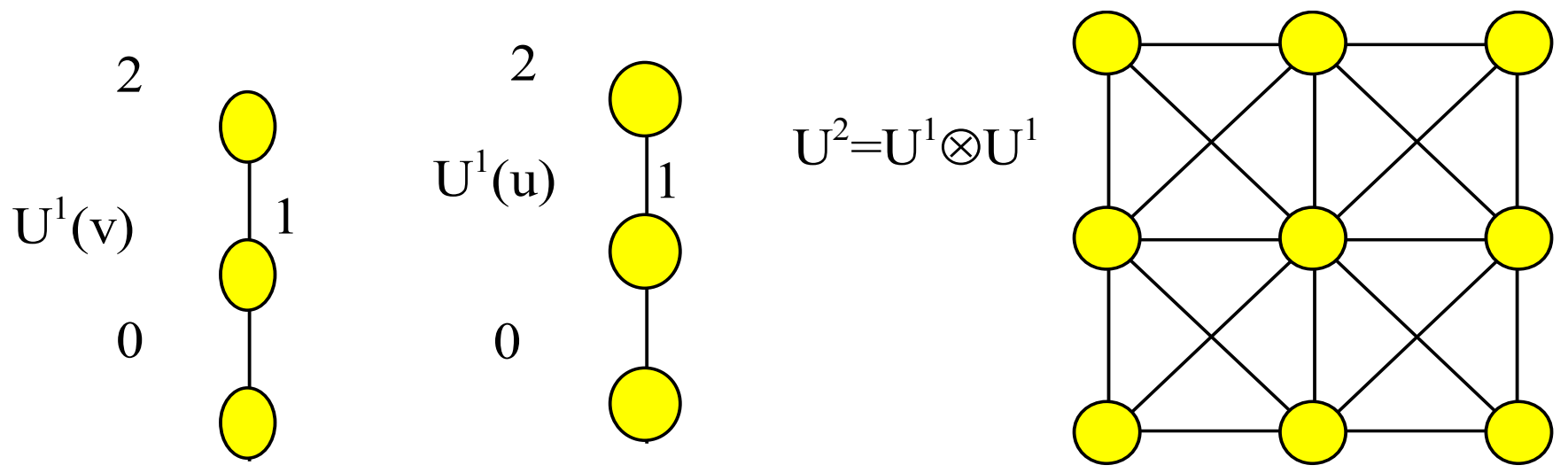

Рис. 153 Прямое произведение двух одномерных нормальных шаров является двумерным шаром и, одновременно, универсальным двумерным кубом.

состоит из 0, 1 и 2, имеет p строк и (p-1) столбцов с различными координатами и является подматрицей матрицы универсального куба $U^{p-1}$. Следовательно, com(S)≤ p-1. Доказательство закончено.□

Теорема 144

*Компактность связки K(n) определяется соотношением $\log_2 n \leq com(K(n)) < \log_2 n+1$, или $com(K(n))=\lceil \log_2 n \rceil$, где $\lceil \log_2 n \rceil$ есть наименьшее целое число, не меньшее $\log_2 n$. Кроме того, com(K(n))= kir(K(n)) (Рис. 150).*

Доказательство.

Доказательство такое же как для кирпичности. Пусть некоторый кирпич связки имеет нулевые координаты $a_1$=(0, 0,...0,..). Тогда координаты всех остальных кирпичей должны состоять из 0 и 1 (или -1). Количество наборов длины s, состоящих из 0 и 1, равно $2^s$. Следовательно, $2^s \leq n < 2^{s+1}$. Логарифмируя это выражения мы имеем $s \leq \log_2 n < s+1$. Следовательно,

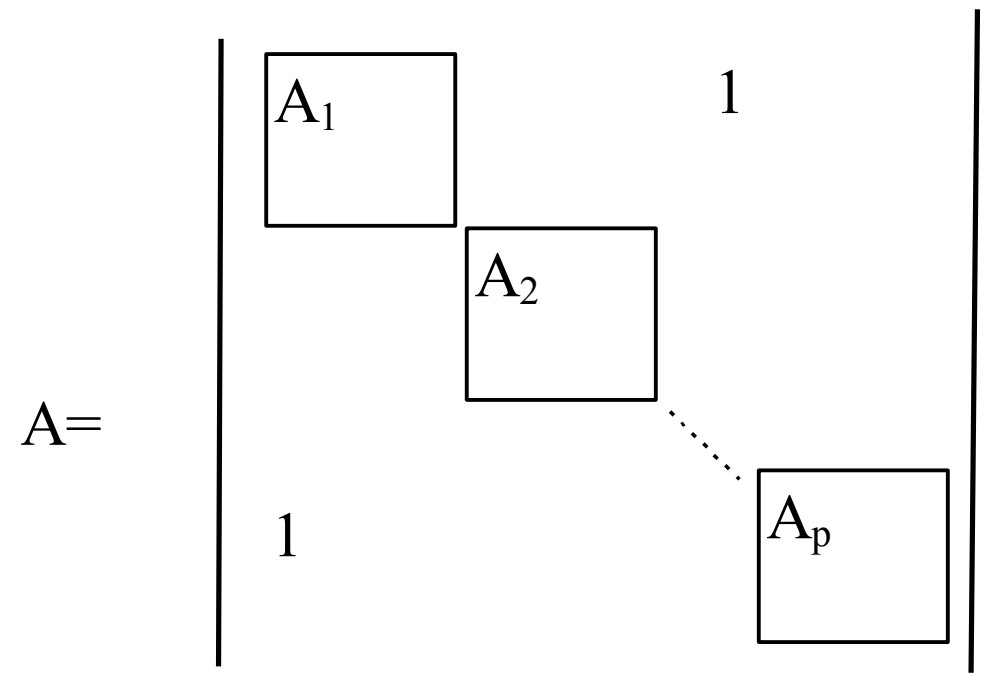

Рис. 154 Построение матрицы p-дольного пространства. Матрицы $A_k$ состоят из нулей.

$com(K(n))=\lceil \log_2 n \rceil$, где $\lceil \log_2 n \rceil$ есть наименьшее целое число, не меньшее $\log_2 n$. Теорема доказана. □

Теорема 145

*Компактность абсолютно несвязного пространства H(n), состоящего из n несвязных кирпичей (Рис. 155), определяется соотношением $\log_2 n \leq com(H(n)) < \log_2 n+1$, или $com(H(n))=\lceil \log_2 n \rceil$, где $\lceil \log_2 n \rceil$ есть наименьшее целое число, не меньшее $\log_2 n$.*

Доказательство.

Доказательство такое же как в предыдущей теореме, если все 1 заменить на 2. Теорема доказана. □

Следствие.

kir(H(n))<com(H(n)) для n>2.

Доказательство

Очевидно, что kir(H(n))=1 для n>1. Следовательно, kir(H(n))=com(H(n)) только для n=1, 2. □

Теорема 146

*Пусть S есть кирпичное пространство, $S_1$ есть его подпространство. Тогда com(S)=com($S_1$)+|S|-|$S_1$|, где |S| и |$S_1$| есть число кирпичей в S и $S_1$.*

Доказательство.

Пусть пространство S имеет кирпичи ($a_1$, $a_2$,...$a_n$, $b_1$, $b_2$,...$b_s$), $S_1$ имеет кирпичи ($a_1$, $a_2$,...$a_n$). Условие com($S_1$)=r означает, что координатная матрица A($S_1$) имеет r ненулевых столбцов. Обозначим первую координату $b_1$ как 1, (r+1) координату как 2. Для каждого кирпича $a_p$ обозначим координату $a_{p(r+1)}$ как 1, если $a_p$ и $b_1$ смежны, и как 0, если $a_p$ и $b_1$ несмежны. Используя эти операции мы присоединяем $b_1$, $b_2$,...$b_s$ к $S_1$ и получаем S. При этом число столбцов в матрице A($S_1$) увеличивается не более чем на s. Следовательно, com(S)≤r+s, где s=|S|-|$S_1$|. Теорема доказана. □

Теорема 147

*Пусть S есть кирпичное пространство, |S|≥4. Тогда com(S)≤|S|-2.*

Доказательство.

Легко проверить непосредственно, что любая поверхность, содержащая 4 кирпича, является подпространством универсального двумерного куба $U^2$, com(S)≤2, если |S|=4. Отсюда, при |S|>4. com(S)≤2+|S|-4=|S|-2. Теорема доказана. □

Напомним определение p-дольного пространства

<u>Определение p-дольного пространства.</u>

Пространство K($n_1$, $n_2$,... $n_p$) называется p-дольным, если его можно представить как сумму p вполне несвязных пространств H($n_1$), H($n_2$)...H($n_p$), K($n_1$,$n_2$,...$n_p$)= H($n_1$)⊕H($n_2$)⊕...⊕H($n_p$).

Теорема 148

Пусть K($n_1$,$n_2$,...$n_p$)= H($n_1$)⊕H($n_2$)⊕...⊕H($n_p$) есть p-дольное пространство, где любое $n_k$>1. Тогда

com(K($n_1$,$n_2$,...$n_p$))=com(H($n_1$)⊕H($n_2$)⊕...⊕H($n_p$))=

$=\lceil \log_2 n_1 \rceil + \lceil \log_2 n_2 \rceil + \ldots + \lceil \log_2 n_p \rceil$.

Доказательство.

Согласно предыдущей теореме пространство H(n) может быть представлено координатной матрицей A(H(n)), имеющей $\lceil \log_2 n_1 \rceil$ столбцов и n строк. Построим матрицу $A(K(n_1,n_2,...n_p))$ состоящей из блоков $A_1=A(H(n_1))$, $A_2=A(H(n_2)),...A_p=A(H(n_p))$, в которой все элементы, не принадлежащие блокам, равны 1 (Рис. 154). Легко видеть, что эта матрица

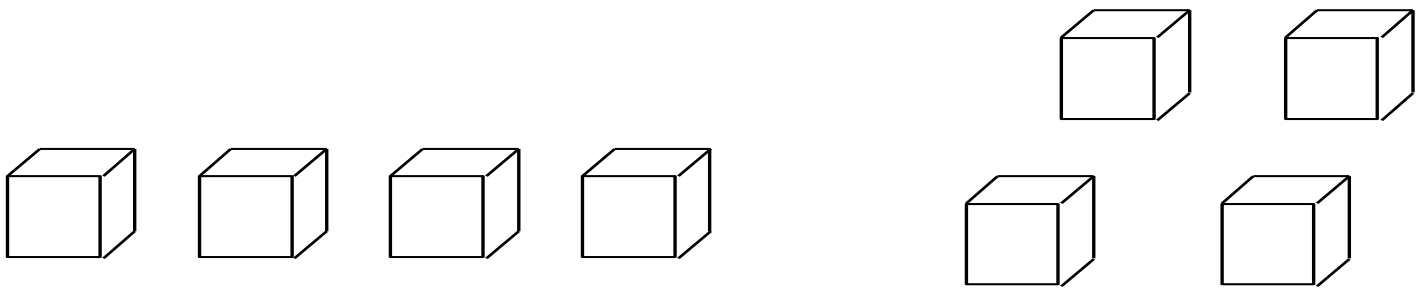

Рис. 155 Кирпичность пространства H(4) равна 1, компактность H(4) равна 2.

имеет $\lceil \log_2 n_1 \rceil + \lceil \log_2 n_2 \rceil + ... + \lceil \log_2 n_p \rceil$ столбцов и является координатной матрицей пространства $K(n_1,n_2,...n_p)= H(n_1)\oplus H(n_2)\oplus ...\oplus H(n_p)$.Теорема доказана. □

## *НЕКОТОРЫЕ НАГЛЯДНЫЕ ОБОСНОВАНИЯ ВВЕДЕНИЯ КИРПИЧЕЙ КАК УНИВЕРСАЛЬНЫХ СТРОИТЕЛЬНЫХ ЭЛЕМЕНТОВ*

Среди читателей, особенно математиков, может возникнуть вопрос: зачем вообще необходимо введение каких-то бесконечно-мерных кирпичей, которые, в принципе, представляют просто некоторое изображение точек молекулярного пространства? Мы здесь изложим нашу концепцию, которая может быть названа физической, или экспериментальной (мы уже описывали компьютерные эксперименты), и которая есть не что иное, как объяснение подхода к дискретным пространствам, и должна быть представлена на страницах этой книги и объяснена как можно более подробно.

Прежде всего, давайте рассмотрим классический математический подход к описанию двумерной поверхности. Математики считают двумерной поверхностью, например, поверхность шара, толщина и объем которой равны нулю. Двумерной поверхностью также является кусок плоскости, толщина плоскости равна нулю, объем куска плоскости также равен нулю. Иными словами, толщина (и объем) двумерной поверхности в математике равны нулю.

Однако, возможна и другая точка зрения. Еще со школьных времен мы знаем, что все тела состоят из молекул. Поверхностью любого материального шара будет внешний слой молекул. Следовательно, эта поверхность имеет ненулевую толщину. Физическая толщина этой поверхности равна одной молекуле. В одной из школьных задач по физике предлагается найти наибольшую площадь нефтяной пленки на поверхности воды. При этом подразумевается, что наименьшая толщина пленки равна также одной молекуле, но не равна нулю. Иными словами, двумерных поверхностей нулевой толщины в природе не существует. Наш подход как раз и отражает эту закономерность и вводит ее в математический и прикладной обиход.

Мы можем считать, что с точностью до ненулевого постоянного коэффициента пропорциональности, минимальная толщина двумерной поверхности равна единице.

С другой стороны, какова наибольшая толщина, при которой данная поверхность еще остается двумерной? Ответ на этот вопрос также достаточно прост и ясен в рамках данного подхода. Чтобы образовалась трехмерное пространство необходимо, чтобы, по крайней мере, некоторые

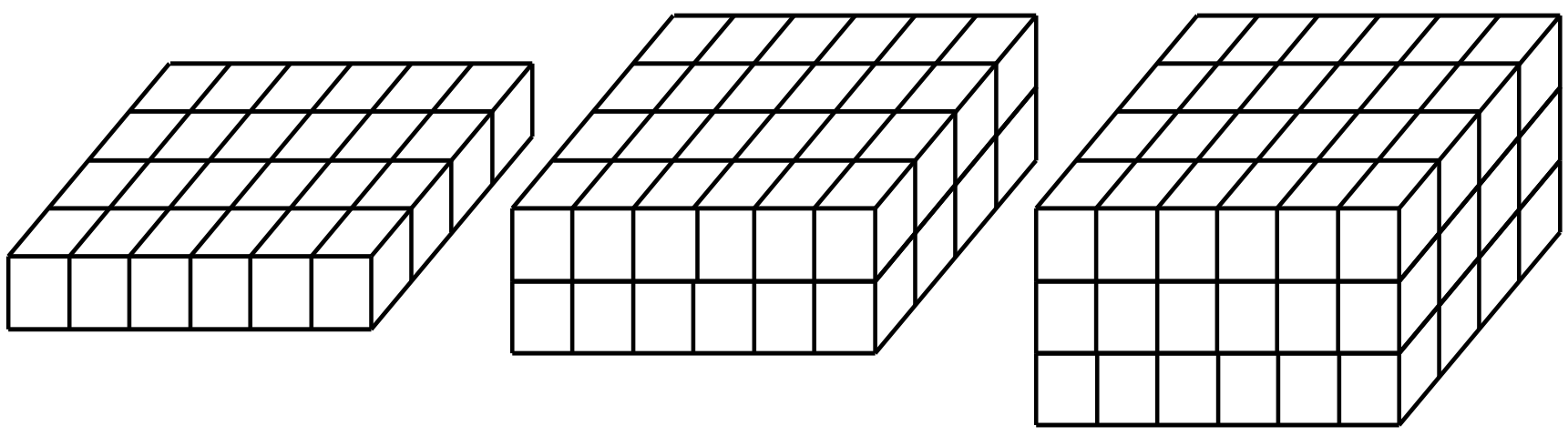

Рис. 156 Двумерная поверхность может иметь толщину один или два и состоять из одного или двух слоев. Три слоя образуют уже трехмерное пространство. Три слоя будут трехмерным пространством, так как каждая точка промежуточного слоя имеет окаем, содержащий двумерную сферу.

его точки были трехмерными. Возьмем несколько слоев двумерной плоскости, образованной кирпичами, и будем поочередно укладывать их один на другой, как это показано на Рис. 156. Выясним, сколько слоев нужно уложить, чтобы получить трехмерные точки.

Легко видеть, что пространство в два слоя, изображенное на Рис. 156, будет все еще двумерным, так как каждая точка имеет в своем окаеме одномерную окружность, но не имеет двумерной сферы. Однако, три слоя, изображенные на Рис. 156, уже будут трехмерным пространством, так как каждая точка промежуточного слоя имеет окаем, содержащий двумерную сферу. Таким образом, двумерная поверхность может иметь толщину один или два и состоять из одного или двух слоев. Три слоя образуют уже трехмерное пространство. Мы можем считать, что с точностью до ненулевого постоянного коэффициента пропорциональности толщина двумерной поверхности равна единице или двум. Минимальная толщина трехмерного (и любого n-мерного) слоя равна трем, причем средний слой состоит из трехмерных (n-мерных) точек, а два внешних слоя содержат двумерные ((n-1)-мерные) точки. Следовательно, можно сформировать n-мерный слой как три слоя (n-1)-мерного пространства, склеенных аналогично предыдущему.

Список литературы к главе 12.

# МОЛЕКУЛЯРНЫЕ МОДЕЛИ N-МЕРНЫХ ЕВКЛИДОВЫХ ПРОСТРАНСТВ

In this chapter we introduce and study molecular models considered as n-dimensional Euclidean spaces $E^n$. Euclidean n-dimensional molecular space is a molecular space whose points are points of $E^n$ with integer coordinates. Connection between these points determines the dimension and other features of the model. We study homogeneous normal, homogeneous regular and other molecular models.

Результаты, полученные в этой главе, частично опубликованы в работах [3,4,7,11,12,13,14,26,42,43].

## *ОДНОРОДНОЕ ПОЛНОЕ ЕВКЛИДОВО МОЛЕКУЛЯРНОЕ ПРОСТРАНСТВО*

Мы рассмотрим молекулярную модель n-мерного евклидова пространства, наиболее близкую по своим свойствам к непрерывному плоскому пространству. Эта модель обладает однородностью, то есть окаемы точек в ней являются изоморфными пространствами. Пространство изотропно и в него можно вложить конечное молекулярное пространство. Свойства такого пространства инвариантны относительно прямого произведения двух и более пространств. Важное свойство такого пространства есть то, что оно допускает координатное представление, аналогичное координатному представлению плоского пространства, что позволяет

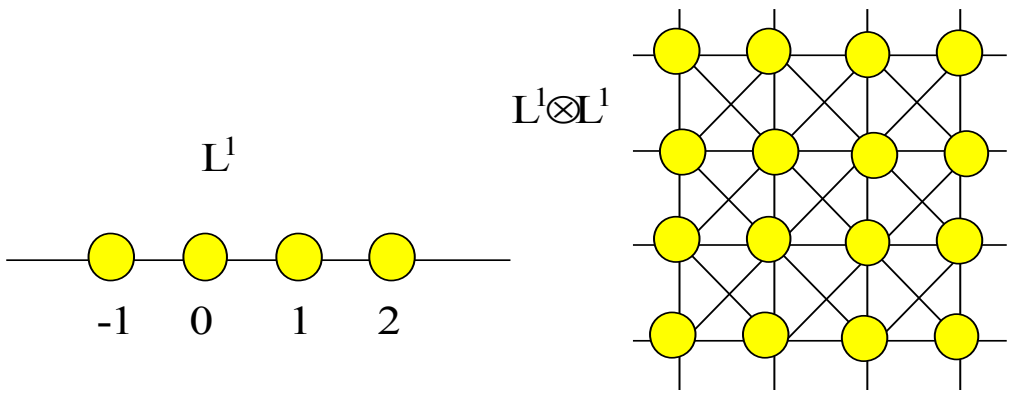

2.

3. *Рис. 1* Прямое произведение двух нормальных прямых есть двумерное молекулярное евклидово пространство.

естественным образом ввести на нем алгебраические и геометрические конструкции, используемые в классической модели. Метрика на этом пространстве, как и на любом молекулярном пространстве, задается естественным образом. Это пространство является подпространством универсального молекулярного пространства, рассмотренного ранее, и к нему применимы многие из свойств универсального молекулярного пространства. Кроме того, из этого пространства путем точечного отбрасывания связей получаются все ранее известные дискретные пространства. Начнем с определений. Как известно, нормальное связное пространство, состоящее из бесконечного множества одномерных точек, будет нормальной бесконечной прямой

Определение молекулярного n-мерного полного евклидова пространства $L^n$

Молекулярным n-мерным полным евклидовым пространством $L^n$ называется прямое произведение n копий одномерной нормальной прямой $L^1$, $L^n=L^1\otimes L^1\otimes\ldots\otimes L^1$. Обозначим точки одномерной прямой $L^1$ целыми числами 0, ±1, ±2, ±3,... Рассмотрим прямое произведение $L^n=L^1\otimes L^1\otimes\ldots\otimes L^1$. Будем считать наборы целых чисел $(s_1, s_2,..s_n.)=[s_k]$, k=0,1,...n, координатами точки $v=s_1\otimes s_2\otimes..\otimes s_n$.

Теорема 149

Евклидово n-мерное пространство $L^n$ есть множество точек $v=(s_1, s_2,..s_n.)=[s_k]$, k=0,1,...n, с целочисленными координатами s в n-мерном евклидовом пространстве $E^n$, причем две точки $v_1$=[sk] и v2=[mk] считаются смежными, если наибольшая из разностей между их соответствующими координатами, взятая по модулю, равна 1,

$$D(v_1,v_2) = \max|s_k - m_k|=1, k=0,1,...n.$$

Доказательство.

Обозначим точки одномерной прямой $L^1$ целыми числами 0, ±1, ±2, ±3,... Теорема очевидна для n=2, как это следует из *Рис. 1* и *Рис. 157*. Рассмотрим прямое произведение $L^n=L^1\otimes L^1\otimes\ldots\otimes L^1$. Оно содержит точки $v=(s_1, s_2,..s_n.)=[s_k]$, k=0,1,...n, с целочисленными координатами. Согласно определению прямого произведения шар любой точки, например $v_1$=(0,

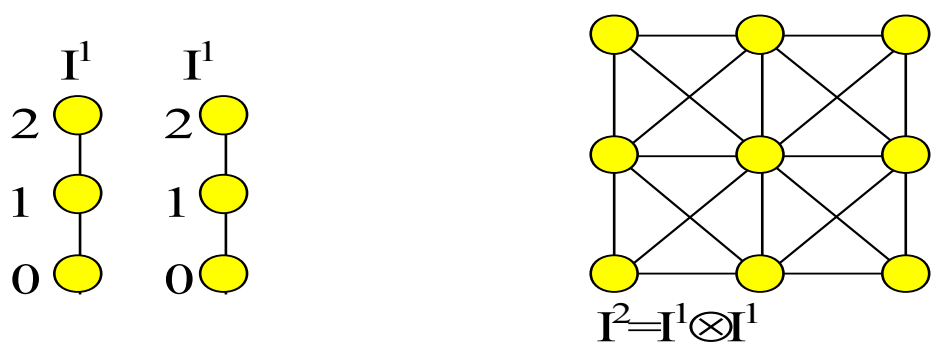

*Рис. 157* Прямое произведение двух одномерных нормальных отрезков $I^1$ является двумерным шаром $I^2$ и, одновременно, универсальным двумерным кубом.

0,..0) состоит из точек, у которых первая координата принадлежит одномерному шару $I^1(s_1)=I^1(0)$ точки 0, вторая координата принадлежит одномерному шару $I^1(s_2)=I^1(0)$ точки 0 так далее, и “n” координата принадлежит одномерному шару $I^1(s_n)=I^1(0)$ точки 0. Очевидно, что шар точки $v_1$=(0, 0,..0) состоит из точек, у которых координаты принимают значения 0, ±1. Это означает, что все точки, смежные точке $v_1$=(0, 0,..0) и только они удовлетворяют условию $D(v_1,v_2) = \max|s_k - m_k|=1$, $k=0,1,...n$. То же самое справедливо для произвольной точки в $L^n$. Теорема доказана..□

На Рис. 1 показано двумерное евклидово пространство. Согласно свойствам прямого произведения шар любой точки из прямого произведения определяется выражением $U(v_k \otimes u_p ... \otimes w_s)=U(v_k) \otimes U(u_p) ... \otimes U(w_s)$, то есть является прямым произведением шаров точек. В применении к данному случаю это означает, что шар любой точки n-мерного евклидова пространства $L^n$ есть прямое произведение n копий одномерного отрезка из трех точек (Рис. 157). Обозначим шар точки v в n-мерном евклидовом пространстве $L^n$ традиционно как $I^n$. Пространство $L^n$ обладает важными особенностями, которые максимально сближают его с непрерывным евклидовым пространством. Попробуем рассмотреть некоторые, наиболее интересные из этих свойств.

Теорема 150

Пространство $L^n$ однородно.

Доказательство.

Однородность пространства следует из его определения, как прямого произведения n-копий нормальной одномерной прямой $L^1$, которая, естественно, однородна. Окаем точки $v_0$ с координатами (0,0,...0) состоит из точек с координатами 0, 1, -1. Окаемы всех точек имеют сходный вид, отличаясь от координат точки $v_0$ наличием постоянного слагаемого по каждой координате. $O(v) \ni (x_1+a_1,x_2+a_2,...x_n+a_n)$ где $O(v_0) \ni (x_1, x_2,...x_n)$. Теорема доказана..□

Теорема 151

Окаем любой точки пространства $L^n$ гомотопен (n-1)-мерной сфере $S^{n-1}$ и содержит минимальную сферу $S_0^{n-1}$ как свое подпространство.

Доказательство.

Топологически структура любого окаема имеет вид $O(v)=(x_1 \otimes x_2 \otimes ... \otimes x_{n-1} \otimes O(x_n)) \oplus (x_1 \otimes x_2 \otimes ... \otimes O(x_{n-1}) \otimes x_n) \oplus ... \oplus (O(x_1) \otimes x_2 \otimes ... \otimes x_{n-1} \otimes x_n) \cup B$, где все точки B могут быть отброшены путем точечных преобразований.

Подпространство $(x_1 \otimes x_2 \otimes ... \otimes x_{n-1} \otimes O(x_n)) \oplus (x_1 \otimes x_2 \otimes ... \otimes O(x_{n-1}) \otimes x_n) \oplus ... \oplus (O(x_1) \otimes x_2 \otimes ... \otimes x_{n-1} \otimes x_n)$ определяется следующим образом

$(x_1 \otimes x_2 \otimes ... \otimes x_{n-1} \otimes O(x_n)) \oplus (x_1 \otimes x_2 \otimes ... \otimes O(x_{n-1}) \otimes x_n) \oplus ... \oplus (O(x_1) \otimes x_2 \otimes ... \otimes x_{n-1} \otimes x_n) \approx O(x_n) \oplus O(x_{n-1}) \oplus ... \oplus O(x_1)=S^0 \oplus S^0 \oplus ... \oplus S^0=S_0^{n-1}$.

Легко видеть, что точки, принадлежащие $S_0^{n-1}$, содержащейся в окаеме точки $v_0$= (0,0,...0) определяются координатной матрицей

$$\begin{vmatrix} 1 & 0 & \dots & 0 & 0 \\ 0 & 1 & \dots & 0 & 0 \\ \dots & \dots & \dots & \dots & \dots \\ 0 & 0 & \dots & 1 & 0 \\ 0 & 0 & \dots & 0 & 1 \end{vmatrix}$$

Теорема доказана..□

Теорема 1

Точечными преобразованиями пространство $L^n$ может быть стянуто к нормальному n-мерному евклидову пространству.

Доказательство.

Доказательство является следствием из соответствующей теоремы о прямом произведении произвольных нормальных пространств.□□

Отметим еще несколько очевидных свойств, уже рассмотренных ранее.

Замечание

- Шар $I^n(v)$ точки v в n-мерном евклидовом пространстве $L^n$ есть

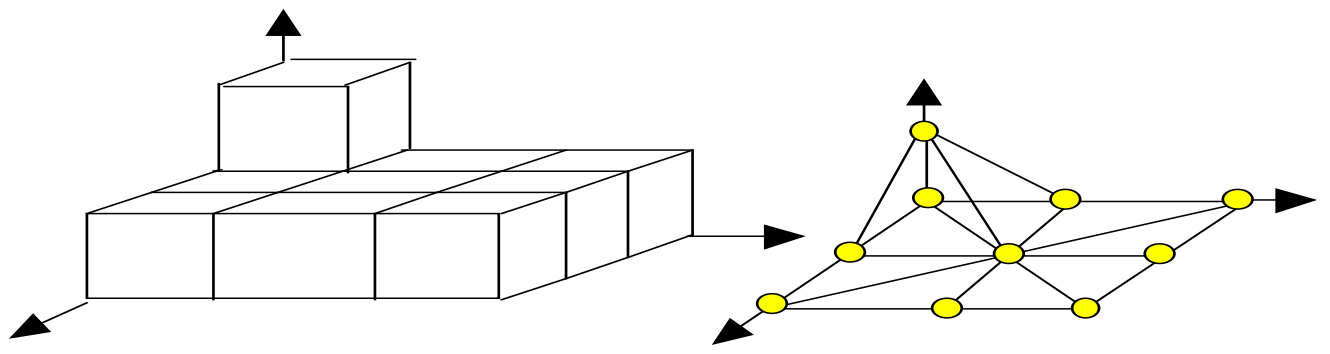

*Рис. 2* Изоморфизм между кирпичным и молекулярным универсальными пространствами.

прямое произведение n копий одномерного отрезка из трех точек (Рис. 157). $I^n(v)=I^1\otimes I^1 \dots \otimes I^1$.
- Шар любой точки v в n-мерном евклидовом пространстве $L^n$ $I^n(v)=I^1\otimes I^1 \dots \otimes I^1$ является универсальным n-кубом, в который можно вписать любое молекулярное пространство, состоящее из (n-2) точек. Шар $I^n(v)$ точки v в n-мерном евклидовом пространстве $L^n$ есть универсальный молекулярный n-куб $U^n$.
- Окаем любой точки пространства $L^n$ не гомотопен никакому n-мерному нормальному замкнутому пространству.

Это свойства следуют из теорем, доказанных ранее.

На Рис. 2 показан изоморфизм между $L^n$, являющимся молекулярным пространством $M(E^n)$, и кирпичным пространством $B(E^n)$.

Теорема 152

$L^n$ есть подпространство универсального бесконечномерного точечного пространства.

Доказательство.

Это теорема является следствием свойств универсального координатного бесконечномерного пространства. □

Векторная алгебра на $L^n$.

На $L^n$ естественным образом задается векторная алгебра, как множество векторов с целочисленными координатами. Это позволяет определить все стандартные структуры, связанные с векторами и существующие в непрерывном случае. Естественным путем задаются такие стандартные операции на векторах, как длина вектора, определяемая как наибольшая по модулю координата вектора, скалярное произведение и тому подобное.

## *ОДНОРОДНОЕ НОРМАЛЬНОЕ ЕВКЛИДОВО МОЛЕКУЛЯРНОЕ ПРОСТРАНСТВО*

Вообще говоря, существует множество молекулярных моделей евклидовых пространств, как нормальных, так и не являющихся нормальными. Мы опишем те из них, которые либо обладают специфическими свойствами, либо уже изучались в научной литературе, а также дадим алгоритм получения молекулярных моделей.

Рассмотрим молекулярное нормальное евклидово пространство, в котором окаемы любых двух точек являются изоморфными подпространствами. Это пространство было введено нами как дискретная нормальная однородная модель n-мерного евклидова пространства. Следует заметить,

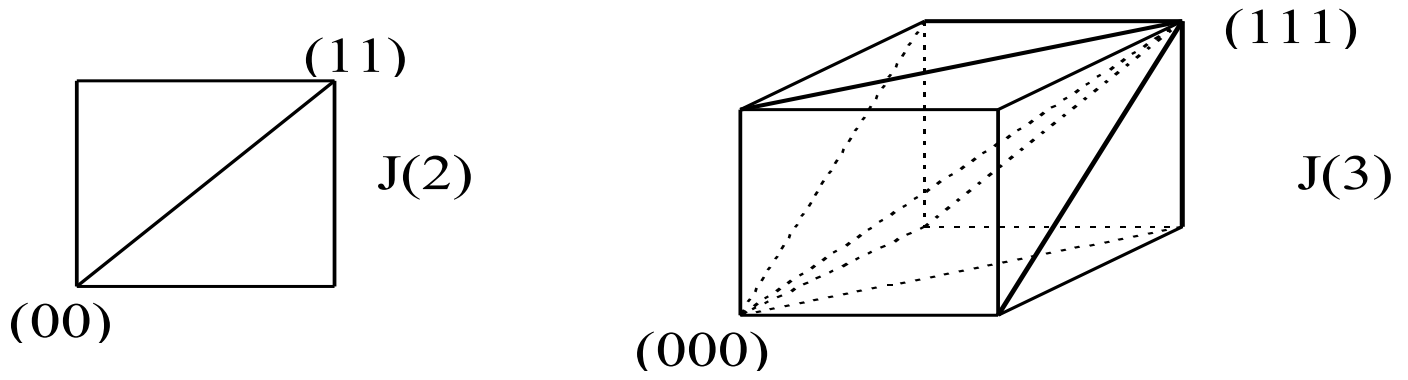

*Рис. 158* Структурные единичные блоки J(2) и J(3). Их диагонали являются 0-мерной и 1-мерной сферами.

что пространство $L^n$, рассмотренное ранее, не является нормальным.

Прежде всего рассмотрим свойства некоторого специального пространства, полученного отбрасыванием точечных связей из прямого произведения n связок K(2) из двух точек. Дадим определение этого пространства.

Определение структурного единичного блока J(n) (СЕБ-1).

Пусть K(2), является связкой из двух смежных точек, обозначенных 0 и 1. Введем прямое произведение $(K(2))^n=K(2)\otimes K(2)...\otimes K(2)$ и отбросим все связи кроме связей между точками $v=(v_1,v_2,...v_n)$ и $u=(u_1,u_2,...u_n)$, удовлетворяющими условиям $v_k \le u_k$, $k=1,2,..n$, или $v_k \geq u_k$, k=1,2,..n. Полученное пространство называется структурным единичным блоком J(n) (СЕБ-1). Будем считать, что $v \leq u$, если $v_k \leq u_k$,. На Рис. 158 показаны J(2) и J(3). Из этого рисунка видно, что если отбросить две точки, имеющие все нулевые координаты и все единичные координаты, то полученное пространство является нуль-мерной и одномерной нормальными сферами Мы покажем, что это свойство выполняется при любом n.

Определение диагонали структурного единичного блока J(n) (СЕБ-1).

Диагональю D(n) структурного единичного блока J(n) называется

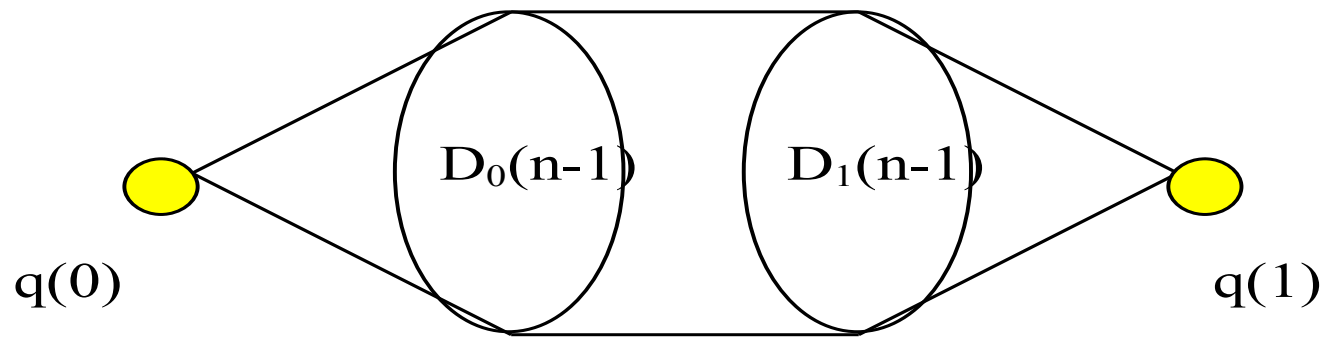

*Рис. 3* Конструкция диагонали D(n), являющейся (n-2)-мерной нормальной сферой.

пространство, полученное из J(n) отбрасыванием двух точек со всеми нулевыми и со всеми единичными координатами. $D(n)=J(n)-p_0(0,0,..0)-p_1(1,1,..1)$.

Теорема 2

1. $|J(n)|=2^n$.
2. $|D(n)|=2^n-2$.
3. J(n) является точечным пространством, а именно, конусом с вершинами $p_0(1,1,...1)$ и $p_1(0,0,...0)$.
4. Если $v \leq v_1$ и $v_1 \leq v_2$, то $v \leq v_2$ ($v \geq v_1$ и $v_1 \geq v_2$, то $v \geq v_2$).

Доказательство.

Первые два свойства очевидны. Третье свойство следует из того, что для любая точка v из J(n) удовлетворяет соотношению $p_0(0,0,..0) \leq v \leq p_1(1,1,..1)$. Следовательно, $p_0$ и $p_1$ будут смежными всем остальным точкам пространства. Для доказательства четвертого свойства возьмем точку $v_1=(a_1,a_2,...a_s,b_1,b_2,...b_t)$, где все $a_k=0$, все $b_k$. Если $v \leq v_1$, то $v=(a_1,a_2,...a_s,c_1,c_2,...c_t)$, если $v_1 \leq v_2$, то $v_2=(c_1,c_2,...c_s,b_1,b_2,...b_t)$, где все $c_k=0,1$. Отсюда видно, что $v \leq v_2$. □

Теорема 153

Диагональ D(n), n>1, структурного единичного блока J(n) есть нормальная (n-2)-мерная сфера $S^{n-2}$.

Доказательство.

Для n=2,3 эта теорема является очевидной (Рис. 3). Используем индукцию. Пусть теорема выполняется для D(k), k<n. Рассмотрим D(n). выберем произвольную точку из D(n), $v_1=(a_1,a_2,...a_s,b_1,b_2,...b_t)$, где все $a_s=0$, все $b_k=1$, s>1, t>1. Тогда O(v) состоит из точек $v\leq v_1$, $v=(a_1,a_2,...a_s,c_1,c_2,...c_t)$, и $v_1\leq v_2$, $v_2=(c_1,c_2,...c_s,b_1,b_2,...b_t)$, где все $c_k=0,1$. Очевидно, что все точки v образуют диагональ $D(s)=S^{s-2}$, все точки $v_2$ образуют диагональ $D(t)=S^{t-2}$. Так как все точки v смежны со всеми точками $v_2$, то $O(v)=D(s)\oplus D(t)=S^{s-2}\oplus S^{t-2}=S^{s+t-3}= S^{n-3}$. То же самое относится к любой точке, содержащей s нулей и t единиц в своих координатах. Выберем точку $v_1$ из D(n), $v_1=(a_1,b_2,...b_n)$, где $a_1=0$, все $b_k=1$. Тогда O(v) состоит из точек $v\leq v_1$, $v=(a_1,c_2,...c_n)$, где все $c_k=0,1$. Очевидно, что все точки v образуют диагональ $D(n-1)=S^{n-3}$. То же самое относится ко всем точкам, содержащим только один нуль или одну единицу в своих координатах. Следовательно, окаем каждой точки из D(n) есть нормальная сфера $=S^{n-3}$, то есть D(n) является нормальным (n-2)-мерным замкнутым пространством.

Докажем, сейчас, что D(n) является нормальной (n-2)-мерной сферой. Выделим из D(n) две точки q(0)=(0,1,1,...1) и q(1)=(1,0,0,...0). Точки v, смежные с точкой q(0), определяются координатами $v=(0,c_2,...c_n)$, где все $c_k=0,1$, $0<c_2+c_3+...c_n<n-1$. Точки $v_2$, смежные с точкой q(1), определяются координатами $v_2=(1,c_2,...c_n)$, где все $c_k=0,1$, $0<c_2+c_3+...c_n<n-1$. Точки v образуют диагональ $D_0(n-1)$, все точки которой смежны с q(0) и не смежны с q(1). Аналогично точки $v_2$ образуют диагональ $D_1(n-1)$, все точки которой смежны с q(1) и не смежны с q(0). Полученная конструкция изображена на Рис. 3.

Покажем, теперь, что точка q(0) может быть соединена точечными преобразованиями со всеми точками подпространства $D_1(n-1)$, а все точки подпространства $D_0(n-1)$ после этого могут быть отброшены. Рассмотрим общий окаем точки $O(q(0)v_2)$ q(0) и некоторой точки $v_2=(1,a_1,a_2,...a_s,b_1,b_2,...b_t)$ из $D_1(n-1)$, где все $a_k=0$, все $b_k=1$, s>0, t>0. Очевидно, что $O(q(0)v_2)$ есть подпространство, принадлежащее $D_0(n-1)$, и состоящее из точек $v_1\leq q(0)$, $v_1\leq v_2$. Следовательно, $O(q(0)v_2)\ni v_1(0,a_1,a_2,...a_s,c_1,c_2,...c_t)$, где все $c_k=0,1$, $0<c_2+c_3+...c_t$. Это подпространство является точечным подпространством, а именно, конусом с вершиной $v_0(0,a_1,a_2,...a_s,b_1,b_2,...b_t)$. На основании этого мы можем соединить все точки из $D_1(n-1)$, имеющие только одну нулевую координату с точкой q(0). Затем соединяем все точки из $D_1(n-1)$, имеющие только две нулевые координаты, с точкой q(0). Предположим, что мы соединили все точки, имеющие менее s нулевых координат, с точкой q(0).

Возьмем произвольную точку $u(1,a_1,a_2,...a_s,b_1,b_2,...b_t)$ из $D_1(n-1)$, имеющую s нулевых координат и не соединенную с q(0). Общий окаем $O(q(0)u)=A\cup B$, где $A= O(q(0)u)\cap D_0(n-1)$, $B=O(q(0)u)\cap D_1(n-1)$. Очевидно $A\ni v_1(0,a_1,a_2,...a_s,c_1,c_2,...c_t)$, где все $c_k=0,1$, $0<c_2+c_3+...c_t$, $B\ni v_2(1,c_1,c_2,...c_s,b_1,b_2,...b_t)$, где все $c_k=0,1$, $0<c_1+c_2+...+c_s,<s$. При этом точка $v_0(0,a_1,a_2,...a_s,b_1,b_2,...b_t)$ будет смежной со всеми точками из A и B. Следовательно, подпространство O(q(0)u) является точечным подпространством, а именно, конусом с вершиной $v_0(0,a_1,a_2,...a_s,b_1,b_2,...b_t)$, и поэтому можно установить связь между точками q(0) и u. Так как любые две точки в $D_1(n-1)$, имеющие ровно s нулевых координат, не являются смежными, связь точки q(0) со всеми такими точками может быть установлена одновременно. Таким образом, все точки подпространства $D_1(n-1)$ могут соединены с точкой q(0). После этого все точки подпространства $D_0(n-1)$ могут быть отброшены, как имеющие точечные окаемы, являющиеся конусами.

Полученное пространство H является прямой суммой 0-мерной сферы $S^0(q(0),q(1))$, состоящей из двух несмежных точек q(0) и q(1), и подпространства $D_1(n-1)$, являющегося (n-3)-мерной сферой $S^{n-3}$, $H=S^0\oplus S^{n-3}=S^{n-2}$. Следовательно, исходное пространство D(n) есть также (n-2)-мерная нормальная сфера. Теорема доказана. □

Свойства СЕБ-1 изучались нами, потому что они определяют свойства нормального однородного евклидова пространства, определение которого мы сейчас введем.

<u>Определение нормального однородного евклидова пространства</u>

> Молекулярным нормальным однородным n-мерным евклидовым пространством $R^n$ называется пространство, каждая точка v которого определяется набором из n целочисленных координат $v=(s_1, s_2,..s_n.)=[s_k]$, $k=0,1,...n$, и в котором две точки $v_1=(s_1, s_2,..s_n)=[s_k]$ и $v_2=(r_1, r_2,..r_n)=[s_k]$, являются смежными тттк, $|s_k-r_k|\le 1$, $s_k\le r_k$, $k=0,1,...n$. При этом будем считать, что $v_1\le v_2$.

Хотя мы называем это пространство нормальным, нам предстоит доказать, что окаем каждой точки этого пространства является нормальной (n-1)-мерной сферой. На Рис. 159 изображено двумерное пространство $R^2$.

Теорема 154

Окаем каждой точки пространства $R^n$ есть нормальная (n-1)-мерная сфера $S^{n-1}$.

Д о к а з а т е л ь с т в о.

Для n=2 эта теорема является очевидной. Используем индукцию. Пусть теорема выполняется для $R^k$, k<n. Рассмотрим $R^n$. Так как пространство является однородным, достаточно доказать эту теорему для точки $v_0$ с координатами (0,0,...0). Согласно определению пространства, окаем $O(v_0)$ этой точки содержит точки $v=(a_1,a_2,...a_n)$, где все $a_s$=0,1,-1, для которых выполняется условие $v_0 \le v$, или $v \le v_0$. Отсюда видно, что $O(v_0)=A\cup B$, где A состоит из точек $v_1$, $v_0 \le v_1$, B состоит из точек $v_2$, $v_2 \le v_0$. Следовательно $A=\{v_1=(c_1,c_2,...c_n)\}$, где все $c_k$=0,1, $B=\{v_2=(d_1,d_2,...d_n)\}$, где все $d_k$=0,-1. Очевидно, что никакая точка из $O(v_0)$ не содержит одновременно координат 1 и (-1). Докажем, что любая точка из окаема $O(v_0)$ имеет окаем, являющийся нормальной замкнутой (n-2)-мерной сферой.

Пусть $v_1=(a_1,a_2,...a_s,b_1,b_2,...b_t)\in A$, ãде все $a_k$=0, все $b_k$=1, s>1, t>1. Любой другой порядок следования точек всегда может быть приведен к данному простым переобозначением координат без изменения смежности точек. Тоãда

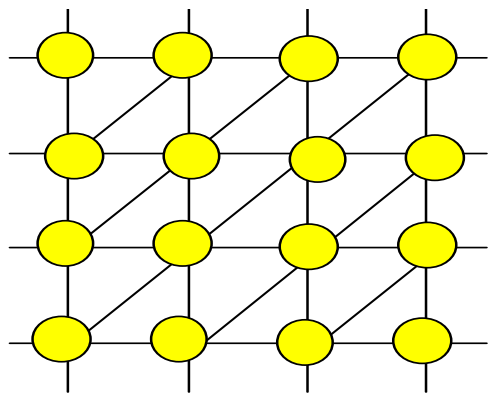

*Рис. 159* Двумерное нормальное однородное евклидово пространство $R^2$.

$O(v_0v_1)=\{v=(d_1,d_2,...d_s,0,0,...0)\}\cup\{v=(c_1,c_2,...c_s,1,1,...1)\}\cup\{v=(0,0,...0,c_1,c_2,...c_t)\}=C\cup D\cup E$, где все $c_k$=0,1, все $d_k$=0,-1. Очевидно, что подпространство, состоящее из точек $C\cup D=\{v=(d_1,d_2,...d_s,0,0,...0)\}\cup\{v=(c_1,c_2,...c_s,1,1,...1)\}$, является окаемом точки с нулевыми координатами в $R^s$. Согласно индукции это подпространство есть нормальная (s-1)-мерная сфера $S^{s-1}$. Подпространство $C=\{v=(0,0,...0,c_1,c_2,...c_t)\}$ является нормальной (t-2)-мерной сферой $S^{t-2}$ как диагональ в J(t). Кроме того каждая точка из $S^{s-1}$ соединена с каждой точкой из $S^{t-2}$. Следовательно, $O(v_0v_1)=S^{s-1}\oplus S^{t-2}=S^{s+t-2}=S^{n-2}$. Рассмотрим, теперь, точку $v_1=(b_1,b_2,...b_n)$, где все $b_k$=1. $O(v_0v_1)=\{v=(c_1,c_2,...c_n)\}$ где все $c_k$=0,1. Это подпространство является нормальной (n-2)-мерной сферой $S^{n-2}$ как диагональ в J(n). Отсюда следует, что каждая точка окаема $O(v_0)$ есть нормальная (n-2)-мерная сфера, то есть $O(v_0)$ является нормальным (n-1)-мерным замкнутым пространством (многообразием).

Докажем, теперь, что $O(v_0)$ является нормальной (n-2)-мерной сферой $S^{n-1}$. Для этого покажем, что точка u=(1,1,...1) может быть соединена точечными преобразованиями со всеми точками из подпространства $B=\{v_2=(d_1,d_2,...d_n)\}$, где все $d_k$=0,-1, а затем все точки из подпространства

A могут быть отброшены как имеющие точечные окаемы. Установим связь точки u со всеми точками из B, имеющими только одну координату (-1). Затем установим связь точки u со всеми точками из B, имеющими точно две координаты (-1). Предположим, что мы уже установили связь точки u со всеми точками из B, имеющими менее s координат (-1). Рассмотрим общий окаем точки u и точки v, имеющей точно s координат (-1), $v=(d_1,d_2,...d_s,0,0,...0)$, где все $d_k=-1$. Тогда $O(uv)=F\cup G=\{v=(d_1,d_2,...d_s,0,0,...0)\}\cup\{v=(0,0,...0,c_1,c_2,...c_t)\}$, где все $d_k=0,-1$, $0<d_1+d_2+...+d_s<s$, все $c_k=0,1$. Очевидно, что F состоит из точек подпространства B, с которыми уже установлена связь. Также ясно, что O(uv) есть точечное пространство, а именно, конус с вершиной $w=(0,0,...0,e_1,e_2,...e_t)$, где все $e_k=1$. Следовательно, можно ввести точечную связь между u и v. Так как любые две точки, имеющие точно s координат (-1) являются несмежными, можно установить связи точки u со всеми такими точками одновременно. Соединим связями точку u со всеми точками из B, кроме точки $w=(-1,-1,...,-1)\in B$. Рассмотрим окаем произвольной точки, кроме u, $v=(c_1,c_2,...c_s,1,1,...1)$, где все $c_k=0,1$, из A. Окаем каждой такой точки есть конус с вершиной $u=(1,1,...1)$. Следовательно, любую точку из A можно отбросить. Получившееся пространство H состоит из несмежных точек u и w, каждая из которых соединена со всеми точками из подпространства $B=\{v_2=(d_1,d_2,...d_n)\}$, где все $d_k=0,-1$, $0<d_1+d_2+...+d_n)<n$. Очевидно, что B является диагональю $D(n)=S^{n-2}$. Следовательно, $H=S^0(u,w)\oplus D(n)= S^0(u,w)\oplus S^{n-2}=S^{n-1}$. Следовательно, исходное пространство $O(v_0)$ есть также (n-1)-мерная нормальная сфера. Теорема доказана. □

Рассмотрим некоторые свойства пространства $R^n$.

Теорема 3

Окаем O(v) любое́ точки пространства $R^n$ имеет обúем $2^{n+1}-2$, $|O(v)|=|S^{n-1}|=2^{n+1}-2$.

Доказательство.

Для n=2 эта теорема является очевидной. Рассмотрим $R^n$. Так как пространство является однородным, достаточно доказать эту теорему для точки $v_0$ с координатами (0,0,...0). Согласно определению пространства, окаем $O(v_0)$ этой точки содержит точки $v=(a_1,a_2,...a_n)$, где все $a_s=0,1,-1$, для которых выполняется условие $v_0\le v$, или $v\le v_0$. Отсюда видно, что $O(v_0)=A\cup B$, где A состоит из точек $v_1$, $v_0\le v_1$, B состоит из точек $v_2$, $v_2\le v_0$. Следовательно $A=\{v_1=(c_1,c_2,...c_n)\}$, где все $c_k=0,1$, $B=\{v_2=(d_1,d_2,...d_n)\}$, где все $d_k=0,-1$. Объем подпространства A равен, очевидно, $2^n-1$, $|A|=2^n-1$.

Таков же объем подпространства B, $|B|=2^n-1$. A и B не имеют общих точек. Следовательно, $|O(v_0)|=|A\cup B|= O(v_0)=|A|+|B|=2^{n+1}-2$. Теорема доказана. □
Докажем еще одну теорему, показывающую связь между полным пространством $L^n$ и нормальным пространством $R^n$.

Теорема 155

Полное евклидово пространство $L^n$ стягивается к нормальному пространству $R^n$ точечными отбрасываниями связей (*Рис. 4*).

Доказательство.

Рассмотрим $L^n$. Выберем в $L^n$ две смежные точки u и v, для топологических координат которых не выполняется ни одно из соотношений $v\le u$, $u\le v$. Для простоты пусть точка v имеет геометрические координаты $v=(a_1,a_2,...a_s,b_1,b_2,...b_t)$, где все $a_k=0$, все $b_k=1$; $u=(b_1,b_2,...b_p,a_1,a_2,...a_q)$, где все $a_k=0$, все $b_k=1$, $p>s$, $q<t$. Общий окаем этих точек в $L^n$ $O(uv)=\{w=(d_1,d_2,...d_s,c_1,c_2,...c_r,e_1,e_2,...e_m)\}$ где все $d_p=0,1$, все $c_k=0,1,2$, все $e_k=0,1$, $r=p-s$, $m=n-p$. Очевидно, что это подпространство является конусом с вершиной $v_0=(1,1,..1)$, причем $w\le v_0$ для любой точки из $O(uv)$. Следовательно, связь между точками u и v можно отбросить. Так как отбрасывание связей не меняет отношений порядка между точками, то на любом шаге для любой пары смежных точек, не удовлетворяющих соотношениям $v\le u$, $u\le v$, их общий окаем будет конусом. Таким образом все связи между смежными точками, не удовлетворяющими соотношениям $v\le u$, $u\le v$, можно отбросить, получив в результате пространство $R^n$.
Теорема доказана.□

## *НЕОДНОРОДНОЕ НОРМАЛЬНОЕ ЕВКЛИДОВО МОЛЕКУЛЯРНОЕ ПРОСТРАНСТВО*

Неоднородных нормальных евклидовых молекулярных пространств существует, по-видимому, бесконечное множество и чтобы выделить специальный тип пространства, необходимо вводить дополнительные ограничения. Как один пример неоднородного пространства представляющий рассмотрим молекулярное пространство, внешне напоминающее в размерности 2 и 3 так называемое топологическое произведение двух и трех связных упорядоченных топологических пространств (COTS). Упорядоченное топологическое пространство (COTS) было введено Халимским как дискретная модель одномерного евклидова пространства. Топологические произведения двух и трех COTS рассматривались как молекулярные модели двух и трехмерного евклидовых пространств и изучались Вилсоном, Мейером, Конгом, Копперманом, Халимским и другими [30,31,32,33,34,35]. В рамках указанной модели был доказан дискретный аналог Жордановой теоремы для размерностей два и три. Недостатком такого подхода является то, что, хотя одномерное пространство COTS является полностью упорядоченным множеством, двух и более мерные модели становятся частично-упорядоченными множествами, и, следовательно, появляется значительное число различного типа соотношений между точками (например, различные виды соседства двух точек). Причем число различных типов возрастает с увеличением размерности. Принципиальное отличие молекулярных пространств от пространств типа COTS заключается в том, что в молекулярных пространствах не существует отношения порядка, тогда как отношение порядка является основой определения пространств типа COTS.

Введем определение и рассмотрим некоторые свойства молекулярного евклидова нормального неоднородного пространства N.

Напомним, что $|s_k|$mod2 равно 0, если целое число $s_k$ является четным, и 1, если $s_k$ является нечетным.

<u>Определение топологических координат точки $v=(s_1, s_2,..s_n.)$</u>

Топологическими координатами $tv=(t_1,t_2,...t_n)$ точки $v=(s_1,s_2,...s_n)$ называются числа $t_k$, равные 0 или 1 и определяемые выражением $t_k=|s_k|$ mod2.

Иными словами, если геометрическая координата является четной, то топологическая координата есть 0, если геометрическая координата является нечетной, то топологическая координата равна 1. Например, v=(2,3,4,-5,6), tv=(0,1,0,1,0).

<u>Определение молекулярного n-мерного евклидова нормального неоднородного пространства $N^n$.</u>

Молекулярное евклидово n-мерное нормальное неоднородное пространство $N^n$ есть множество точек $v=(s_1,s_2,..s_n.)=[s_k]$, k=0,1,...n, с целочисленными координатами s в n-мерном евклидовом непрерывном пространстве $E^n$, причем две точки $v_1=[s_k]$ и $v_2=[m_k]$ считаются смежными, если выполняются следующие условия:

- наибольшая из разностей между их соответствующими координатами, взятая по модулю, равна 1, $\max|s_k - m_k|=1$, k=0,1,...n
- $tv_1 \le tv_2$, или $tv_1 \ge tv_2$ (Рис. 161).

З а м е ч а н и е

В пределах этого раздела будем считать, что точки $v_1=[s_k]$ и $v_2=[m_k]$ связаны соотношением $v_1 \le v_2$, если геометрические координаты этих точек связаны соотношением $\max|s_k - m_k|=1$, k=0,1,...n, а топологические координаты точки $v_1$ не превышают топологические координаты точки $v_2$, $tv_1 \le tv_2$.

В качестве примера рассмотрим три точки: v=(0,0,1,-1,2), u=(1,-1,3,5,1), w=(2,-2,2,4,2),. tv=(0,0,1,1,0), tu=(1,1,1,1,1), w=(0,0,0,0,0). Тогда v и u несмежны, v и w также несмежны, u и w смежны и $w \le u$. Рассмотрим, какое количество точек имеет заведомо неизоморфные окаемы в этом пространстве.

Теорема 156

В пространстве $N^n$ размерности n имеется не менее [n/2] точек с неизоморфными окаемами, где [n/2] есть наибольшее целое число, не превышающее n/2.

Доказательство.

Пусть две точки $v=(s_1,s_2,...s_n)=[s_k]$, k=0,1,...n, и $u=(m_1,m_2,...m_n)=[m_k]$ с геометрическими координатами являются смежными, то есть наибольшая из разностей между их соответствующими координатами, взятая по

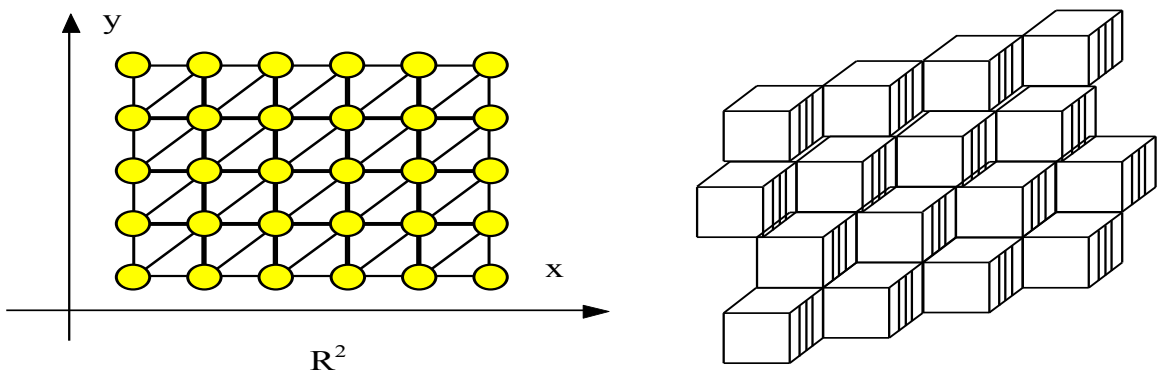

*Рис. 4* Пространство $R^2$ и его модель построенная из кирпичей.

модулю, равна 1, $\max|s_k-m_k|=1$, k=0,1,...n и $|s_k| \bmod 2 \le |m_k| \bmod 2$, k=1,2,...n. Увеличим каждую геометрическую координату на 1. При этом четные координаты перейдут в нечетные, а нечетные станут четными. Выражение $\max|s_k - m_k|=1$, k=0,1,...n не изменится, неравенства $|s_k| \bmod 2 \le |m_k| \bmod 2$, k=1,2,...n, перейдут в неравенства $|s_k| \bmod 2 \ge |m_k| \bmod 2$, k=1,2,...n. Это означает, что точки останутся смежными. Следовательно, окаем точки, имеющей s нулей и t единиц среди топологических координат, s+t=n, является изоморфным окаему точки, имеющей t нулей и s единиц среди топологических координат. В дальнейшем достаточно рассмотреть только точки, имеющие геометрические координаты, состоящие из 0 и 1, поскольку геометрические координаты любой точки могут быть приведены к 0 и 1 простым изменением начала отсчета в евклидовом пространстве. Кроме того, достаточно рассмотреть точки, у которых первые координаты есть нули, а последние координаты состоят из 1. Это осуществляется изменением порядка следования координатных осей.

Связи между точками при этих преобразованиях не меняются. При этом начало координат всегда помещается в точку с четными геометрическими координатами.

Пусть некоторая точка v имеет геометрические (и топологические) координаты $v=(a_1,a_2,...a_s,b_1,b_2,...b_t)$, где все $a_k=0$, все $b_k=1$, $s>0$, $t>0$. Окаем O(v) этой точки состоит из точек u, $u\leq v$, составляющих подпространство A(s) и точек w, $w\geq v$, составляющих подпространство B(t). Тогда $u=(a_1,a_2,...a_s,c_1,c_2,...c_t)$, где все $a_k=0$, все $c_k=0,1,2$, $s>0$, $t>0$, $w=(d_1,d_2,...d_s,b_1,b_2,...b_t)$, где все $d_k=0,1,-1$, все $b_k=1$, $s>0$, $t>0$, $|c_1|+|c_2|+...|c_s|>0$, $|d_1|+|d_2|+...|d_t|>0$. Очевидно, что $u\leq w$. Следовательно, каждая точка из A(s) смежна с каждой точкой из B(t). Тогда $O(v)=A(s)\oplus B(t)$, $s+t=n$. Из этой формулы очевидно, что при различных s окаемы точек могут быть различны. С учетом того, что окаемы точек со всеми четными или нечетными геометрическими координатами являются изоморфными, получаем, что s меняется от 0 до n-1. Таким образом, число полных окаемов, которые могут быть неизоморфными, равно n. (в дальнейшем мы покажем, что все такие окаемы неизоморфны). Теорема доказана. □

Определение чистых и смешанных точек

Точка, пространства $N^n$ у которой все геометрические координаты

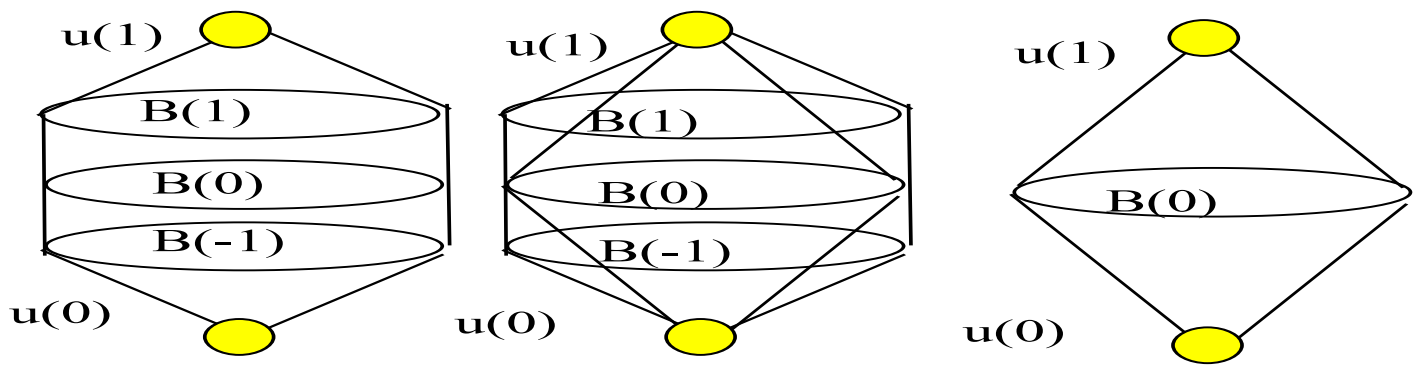

*Рис. 160* Структура окаема чистой точки в пространстве $N^n$ и его преобразование к (n-2)-мерной нормальной сфере.

четны (нечетны) называется чистой точкой. Остальные точки называются смешанными. У чистой точки топологические координаты или все 0, или все 1. Так как окаемы всех чистых точек изоморфны, обозначим окаем чистой точки v $O_p(v)$. Окаем смешанной точки v, имеющей s четных геометрических координат и t нечетных геометрических координат (или t четных геометрических координат и s нечетных геометрических координат), обозначим $O_{st}(v)$.

Теорема 157

*Окаем каждой точки пространства $N^n$ есть нормальная (n-1)-мерная сфера*
*$S^{n-1}$, (Рис. 160).*

Доказательство.

Для n=2 эта теорема является очевидной. Используем индукцию. Пусть теорема выполняется для $N^k$, k<n. Рассмотрим $N^n$. Пусть некоторая смешанная точка v имеет геометрические (и топологические) координаты $v=(a_1,a_2,...a_s,b_1,b_2,...b_t)$, где все $a_k=0$, все $b_k=1$, s>0, t>0. Окаем O(v) этой точки состоит из точек u, u≤v, составляющих подпространство A(s) и точек w, w≥v, составляющих подпространство B(t). Тогда $u=(a_1,a_2,...a_s,c_1,c_2,...c_t)$, где все $a_k=0$, все $c_k=0,1,2$, s>0, t>0, $w=(d_1,d_2,...d_s,b_1,b_2,...b_t)$, где все $d_k=$ 0,1,-1, все $b_k=1$, s>0, t>0, $|c_1|+|c_2|+...|c_s|>0$, $|d_1|+|d_2|+...|d_t|>0$. Очевидно, что u≤w. Следовательно, каждая точка из A(s) смежна с каждой точкой из B(t). Тогда O(v)=A(s)⊕B(t), s+t=n. Пространство A(s) есть окаем чистой точки в пространстве $N^s$. Аналогично пространство B(t) есть окаем чистой точки в пространстве $N^t$. Следовательно, $O_{st}(v)= O_s(v)\oplus O_t(v)$. Так как из индукции $O_s(v)$ и $O_t(v)$ являются нормальными сферами размерности (s-1) и (t-1) соответственно, то $O_{st}(v)$ есть (s+t-1)=(n-1)-мерная нормальная сфера.

Возьмем чистую точку v с геометрическими координатами (0,0,...0) и с теми же топологическими координатами (Рис. 161). Вначале покажем, что любая точка в этом окаеме как подпространстве имеет, в свою очередь, окаем, являющийся (n-2)-мерной сферой. Пусть некоторая смешанная точка u в O(v) имеет геометрические (и топологические) координаты $u=(a_1,a_2,...a_s,b_1,b_2,...b_t)$, где все $a_k=0$, все $b_k=1$, s>0, t>0. Тогда O(vu)=A∪B, где $A=\{w=(a_1,a_2,...a_s,c_1,c_2,...c_t)\}$, где все $a_k=0$, все $c_k=0,1$, $|c_1|+|c_2|+...|c_t|>0$, s>0, t>0, v≤w≤u; $B=\{q=(c_1,c_2,...c_s,b_1,b_2,...b_t)\}$, где все $d_k=0$, все $c_k=0,1,-1$, $|c_1|+|c_2|+...|c_s|>0$, s>0, t>0, v≤q, u≤q. Так как A=D(s) (диагональ в J(s)), $B=O_s(v)$ и w≤q, то $O(vu)=D(s)\oplus O_s(v)=S^{s-2}\oplus S^{t-1}=S^{s+t-2}=S^{n-2}$. Таким образом общий окаем чистой и смешанной смежных точек O(vu) есть нормальная (n-2)-мерная сфера.

Пусть некоторая чистая точка u в O(v) имеет геометрические (и топологические) координаты, равные единице, u=(1,1,...1). Тогда $O(vu)=A=\{w=(c_1,c_2,...c_n)\}$, где все $c_k=0,1$, $|c_1|+|c_2|+...|c_n|>0$, v≤w≤u. Так как A=D(n) (диагональ в J(n)), то $O(vu)=D(n)=S^{n-2}$. Таким образом, общий окаем двух чистых смежных точек O(vu) есть нормальная (n-2)-мерная сфера. Следовательно O(v) есть нормальное замкнутое (n-1)-мерное пространство.

Покажем, теперь, что окаем O(v) может быть приведен точечными преобразованиями к нормальной (n-1)-мерной сфере. Представим O(v) в виде объединения трех подпространств $O(v)=u(1)\cup B(1)\cup B(0)\cup B(-1)\cup u(-1)$, где

$u(1)=(1,0,0,...0)$,

$B(1)=\{w=(1,c_2,c_3,...c_n)\}$, где все $c_k=0,1,-1$, $|c_2|+|c_3|+...|c_n|>0$,

$B(0)=\{w=(0,c_2,c_3,...c_n)\}$, где все $c_k=0,1,-1$, $|c_2|+|c_3|+...|c_n|>0$,

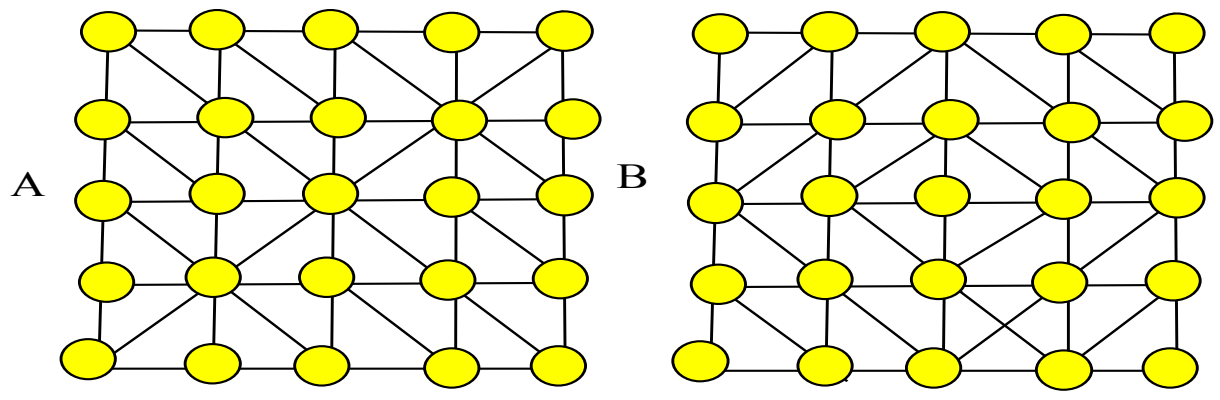

*Рис. 162* Двумерные неоднородные нормальные плоскости.

$B(-1)=\{w=(-1,c_2,c_3,...c_n)\}$, где все $c_k=0,1,-1$, $|c_2|+|c_3|+...|c_n|>0$,

$u(-1)=(-1,0,0,...0)$.

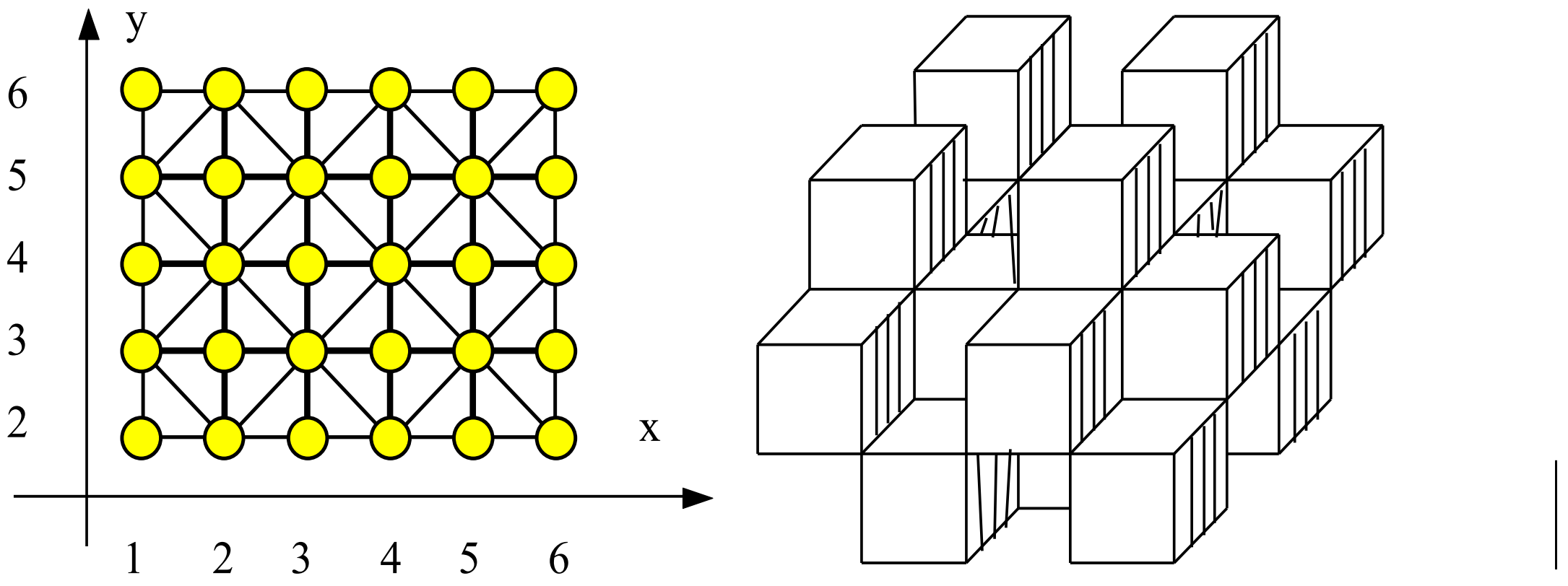

*Рис. 161* Нормальное пространство $N^2$ и его кирпичная модель.

Из индукции следует, что B(0) есть сфера размерности (n-2). Очевидно, что точка u(1) смежна со всеми точками из B(1) и только с ними, точка u(-1) смежна со всеми точками из B(-1) и только с ними, точки из B(1) и B(-1) попарно несмежны. Покажем, что точку u(1) можно соединить точечными связями со всеми точками из B(0), а затем отбросить все точки из B(1), как имеющие точечные окаемы.

Соединим точку u(1) со всеми точками из B(0), имеющими только одну ненулевую координату, затем со всеми точками из B(0), имеющими две ненулевые координаты. Повторяя эту процедуру соединим точку u(1) со всеми точками из B(0), имеющими менее чем s ненулеых координат. Пусть точка w из B(0) имеет точно s ненулевых координат, $w=(0,a_1,a_2,...a_s,0,0,...0)$, где все $|a_k|=1$. Тогда $O(u(1)w)=C\cup D$, где $C\subseteq B(1)$, $D\subseteq B(2)$. Очевидно, что $C=\{q=(1,a_1,a_2,...a_s,d_1,d_2,...d_t)\}$, где все $d_k=0,1,-1$, $|d_1|+|d_2|+...|d_t|\geq 0$, $w\leq q$, $u(1) \le q$. $D=\{p=(0,d_1,d_2,...d_s,0,0,...0)\}$, где все $|d_k|=0,1$, $|d_1|+|d_2|+...|d_s|\ge 0$, поскольку u(1) уже соединена с этими точками. Легко видеть, что точка $r=(1,a_1,a_2,...a_s,0,0,...0)\in B(1)$ является смежной всем точкам из $O(u(1)w)=C\cup D$, следовательно, O(u(1)w) является конусом и можно соединить точку u(1) с точкой w. Так как любые две точки, имеющие ровно s ненулевых координат, являются несмежными, все такие точки можно одновременно соединить с точкой u(1). Используя этот метод мы соединяем точку u(1) с каждой точкой из B(0). После этого любую точку из B(1) можно отбросить, так как окаем этой точки есть конус с вершиной u(1). Точно так же можно соединить точку u(-1) с каждой точкой из B(0) и отбросить все точки из подпространства B(-1). Полученное в результате

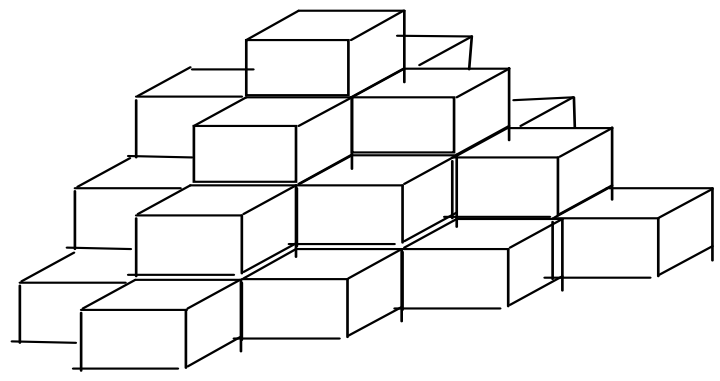

*Рис. 163* Кирпичная модель пространства B.

пространство есть прямая сумма нульмерной сферы, состоящей из двух несмежных точек u(1) и u(-1), и (n-2)-мерной нормальной сферы $B(0)=S^{n-2}$, и является нормальной (n-1)-мерной сферой, $S^0\oplus S^{n-2}=S^{n-1}$. Теорема доказана.□

Рассмотрим, какова связь между полным пространством $L^n$ и нормальным пространством $N^n$, (Рис. 161).

Теорема 158

Полное евклидово пространство $L^n$ стягивается к нормальному пространству $N^n$ точечными отбрасываниями связей.

Д о к а з а т е л ь с т в о.

Доказательство точно такое же, как и в предыдущем разделе. Рассмотрим $L^n$. Выберем в $L^n$ две смежные точки u и v, для топологических координат которых не выполняется ни одно из соотношений $v\le u$, $u\le v$. Для простоты пусть точка v имеет геометрические (и топологические) координаты $v=(a_1,a_2,...a_s,b_1,b_2,...b_t)$, где все $a_k=0$, все $b_k=1$; $u=(b_1,b_2,...b_p,a_1,a_2,...a_q)$, где все $a_k=0$, все $b_k=1$, $p>s$, $q<t$. Общий окаем этих точек в $L^n$ $O(uv)=\{w=(d_1,d_2,...d_s,c_1,c_2,...c_r,e_1,e_2,...e_m)\}$ где все $d_p=0,1$, все $c_k=0,1,2$, все $e_k=0,1$, $r=p-s$, $m=n-p$. Очевидно, что это подпространство является конусом с вершиной $v_0=(1,1,..1)$, причем $w\le v_0$ для любой точки из O(uv). Следовательно, связь между точками u и v можно отбросить. Так как отбрасывание связей не меняет отношений порядка между точками, то на любом шаге для любой пары смежных точек, не удовлетворяющих соотношениям v≤u, u≤v, их общий окаем будет конусом. Таким образом все связи между смежными точками, не удовлетворяющими соотношениям v≤u, u≤v, можно отбросить, получив в результате пространств $N^n$. Теорема доказана.□

В отличие от нормального пространства $R^n$, рассмотренного ранее, пространство $N^n$ неоднородно. Степень неоднородности возрастает с увеличением размерности пространства.

Теорема 4

- Пространство $N^n$ имеет ([n/2]+1) точек с неизоморфными окаемами, где [n/2] есть наибольшее целое число, не превышающее n/2.
- Окаемы всех точек являются нормальными (n-1)-мерными сферами.
- Окаем $O_p(v)=S_p^{n-1}$ чистой точки $v=(a_1,a_2,...a_n)$, имеющей все координаты четные, или все координаты нечетные, содержит $3^n$-1 точек, $|S_p^{n-1}|=3^n-1$.
- Окаем $O_p(v)=S_m^{n-1}$ смешанной точки $v=(a_1,a_2,...a_s,b_1,b_2,...b_t)$, где все $a_k$ четные, все $b_k$ нечетные, или все $b_k$ четные, все $a_k$ нечетные, определяется выражением $O_m(v)=S_m^{n-1}=S_p^{s-1}\oplus S_p^{t-1}$ | $S_m^{n-1}|=3^s+3^t-2$, s+t=n.

Доказательство.

По сути дела эта теорема есть сводка результатов, полученных ранее при доказательствах предыдущих теорем. □

## *3.1 ОБЩИЕ СПОСОБЫ ПОСТРОЕНИЯ ЕВКЛИДОВЫХ МОЛЕКУЛЯРНЫХ ПРОСТРАНСТВ*

Мы рассмотрели два вида евклидовых нормальных пространств, полученных из пространства $L^n$ отбрасыванием точечных связей. Одно пространство, $R^n$, является однородным, второе пространство, рассмотренное потому, что оно сходно с пространством типа COTS, есть неоднородное пространство. Эти пространства иллюстрируют одно из направлений получения нормальных пространств из пространства $L^n$, состоящий в том, что отбрасывание точечных связей дает нормальное пространство. Классификация и изучение таких пространств не входят в наши задачи и требуют отдельного изучения. Дадим еще два нормальных двумерных пространства, обладающих специфическими свойствами: Одно из них является однородным и имеет окаемы, состоящие из 5, 6 и 8 точек, второе обладает центральной симметрией. Оба этих пространства изображены на (Рис. 162). Кирпичная модель пространства B изображена на Рис. 163. Тип пространств можно расширить на произвольное число измерений.

Список литературы к главе 13.

# ПОСТРОЕНИЕ МОЛЕКУЛЯРНЫХ МОДЕЛЕЙ НЕПРЕРЫВНЫХ ПРОСТРАНСТВ

In the first part of this chapter we consider connection between a molecular space and a cover of it by a family of subspaces. We prove a number of theorems which determine when the nerve of a cover is homotopic to the molecular space itself. Then we consider how to build a molecular model for a continuous space. For this purpose we define a special cover of the continuous space which is called complete, then build the nerve of this cover and prove that all nerves of all such covers are homotopic molecular spaces. We expand the class of complete covers and introduce covers with the same properties as complete ones. Finally we describe one more way of obtaining a molecular model of a continuous space, which is a surface in Euclidean space. The method described in this chapter can be used in computers for automatic constructing molecular counterparts. First we define a basic kirpich (and molecular) space as a set of identical cubs (blocks) that form a triangulation of Euclidean space. Given surface S we sort out the set of blocks intersecting S and take the nerve (intersection graph) of this set. In most of cases that are covered by a number of theorems this nerve is a molecular model of S.

Часть результатов, полученных в этой главе, изложена в работах [3,4,10,14,16,17,18,26,27,43,51].

## *ПРАВИЛЬНОЕ ПОКРЫТИЕ МОЛЕКУЛЯРНОГО ПРОСТРАНСТВА СЕМЕЙСТВОМ ПОДПРОСТРАНСТВ И ЕГО СВОЙСТВА*

Для математического обоснования тех методов, которые мы будем использовать при получении молекулярных моделей непрерывных многомерных объектов, нам необходимо рассмотреть связь между пространством и его покрытием подпространствами. Мы хотим установить, когда покрытие обладает теми же свойствами, что и само пространство. Пространство может состоять из большого количества точек, для обработки которых необходима значительная память и время. Покрытие этого же пространства может содержать лишь несколько элементов. Однако, если покрытие и пространство совпадают по своим свойствам, то мы можем работать с элементами покрытия, что существенно может уменьшить объем необходимой памяти компьютера и затраты времени на обработку.
Рассмотрим, вначале чисто формально, покрытие пространства G семейством его подпространств.

<u>Определение покрытия пространства семейством подпространств.</u>

Пусть задано пространство G и семейство $F=(G_1,G_2,...G_n)$ его подпространств. F называется покрытием (cover) пространства G, если любая точка принадлежит по крайней мере одному из подпространств семейства F.

На Рис. 164 изображено пространство G и его покрытие тремя подпространствами $G_1$, $G_2$ и $G_3$. Подпространства, образующие покрытие,

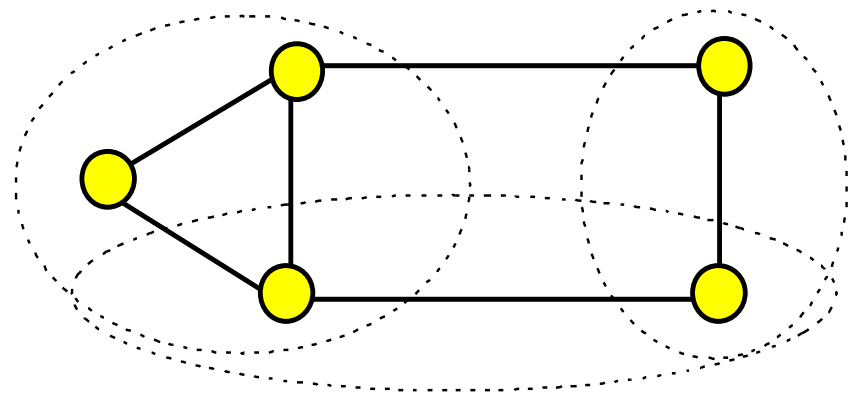

Рис. 164 пространство G и его покрытие тремя подпространствами $G_1$, $G_2$ и $G_3$. Подпространства, образующие покрытие, находятся в соответствующих кружках.

находятся в соответствующих кружках.

Будем считать, что задан пространство G и его покрытие $F=(G_1, G_2,...G_n)$.

<u>.Определение непрерывного семейства подпространств (Рис. 165).</u>

Семейство $F=(G_1, G_2,...G_n)$ называется непрерывным (continuous), если из условия

$G_{p_k} \cap G_{p_m} \neq \varnothing, \forall k, m = 1, 2, \ldots, r$ следует, что

$G_{p_1} \cap G_{p_2} \ldots \cap G_{p_r} \neq \varnothing$.

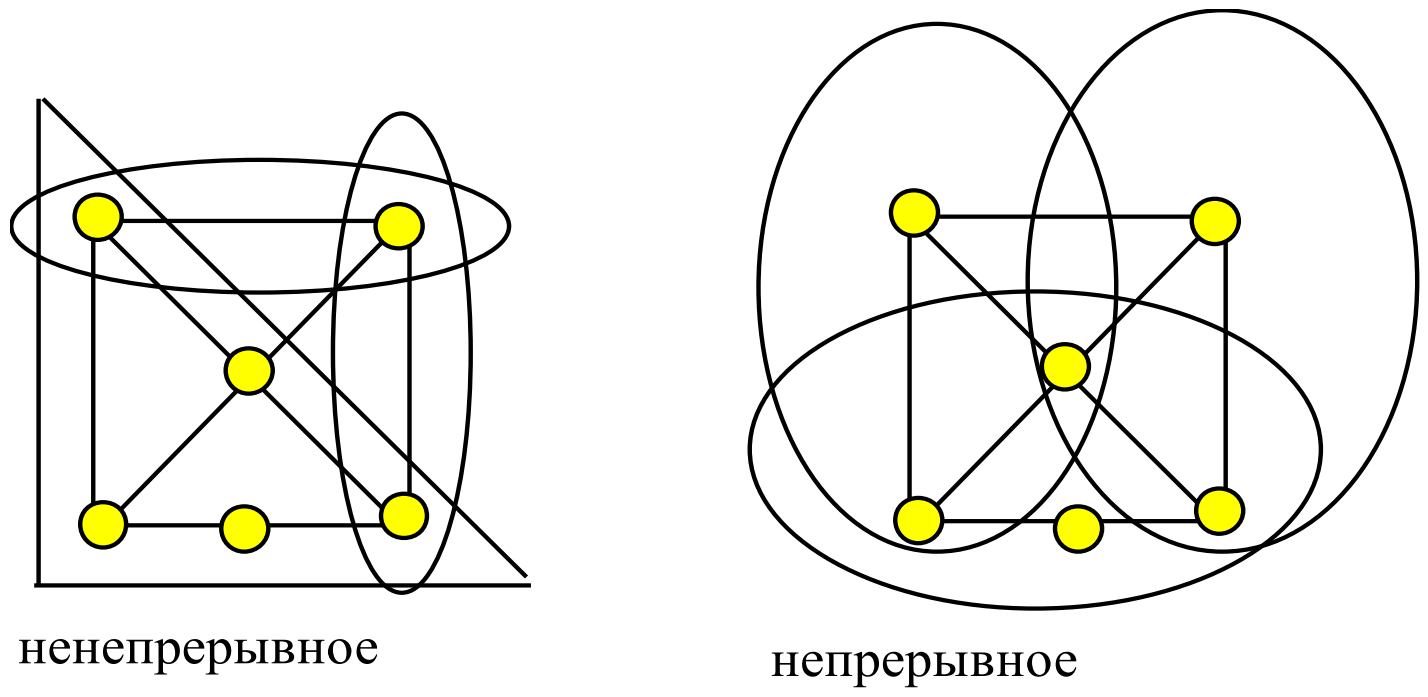

Рис. 165 В непрерывном семействе три подпространства попарно пересекаются, и их общее пересечение также непусто, содержит один элемент. В ненепрерывном семействе три подпространства также попарно пересекаются, но их общее пересечение пусто, не содержит ни одного элемента.

На Рис. 165 изображены непрерывное и ненепрерывное семейства (покрытия). В непрерывном семействе три подпространства попарно пересекаются, и их общее пересечение также непусто, содержит один элемент. В ненепрерывном семействе три подпространства также попарно пересекаются, но их общее пересечение пусто, не содержит ни одного элемента.

Определение стягиваемого семейства подпространств (Рис. 166).

Семейство подпространств(или покрытие) $F=(G_1,G_2,...G_n)$ называется стягиваемым (contractible), если пересечение любого числа

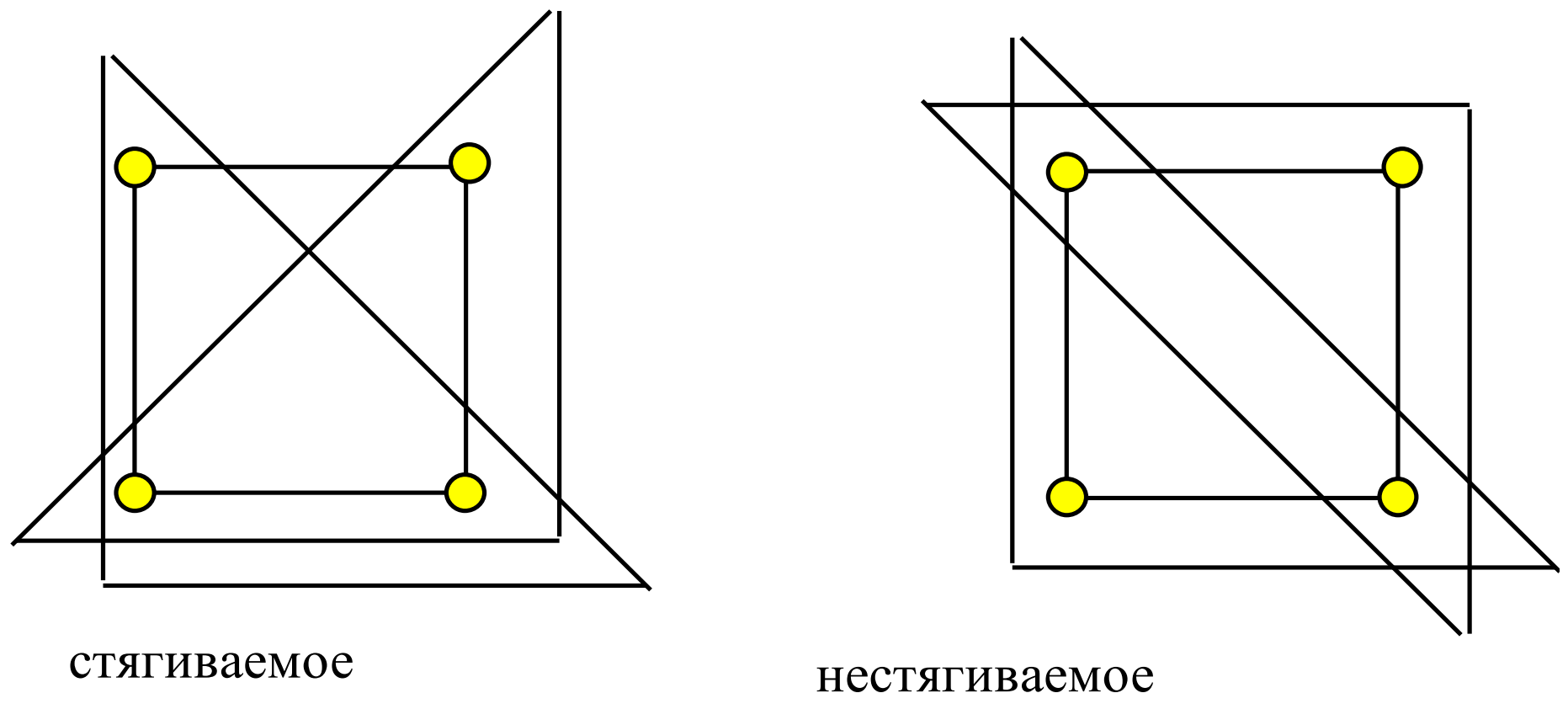

Рис. 166 Семейство содержит два точечных подпространства. В стягиваемом семействе их пересечение также является точечным подпространством, состоящем из двух точек, соединенных связью. В нестягиваемом семействе пересечение является неточечным подпространством, состоящем из двух изолированных точек, то есть 0-мерной.

подпространств покрытия является точечным или пустым.

$G_{k_1} \cap G_p ... \cap G_s = \varnothing$ или точечное подпространство.

Из определения стягиваемости следует, что все подпространства стягиваемого семейства являются точечными. Рассмотрим простые свойства непрерывного и стягиваемого семейств подпространств.

Теорема 159

*Если семейство $F=(G_1,G_2,...G_n)$ является стягиваемым, то*
*Любое подсемейство семейства F будет стягиваемым.*
*Любое семейство $A=(A_1,A_2,...A_p)$, где $A_k$ являются пересечениями элементов из F, будет стягиваемым.*

Теорема 160

*Если семейство $F=(G_1,G_2,...G_n)$ является непрерывным, то*
*Любое подсемейство семейства F будет непрерывным.*
*Любое семейство $A=(A_1,A_2,...A_p)$, где $A_k$ являются пересечениями элементов из F, будет непрерывным.*

Доказательства этих теорем несложны, и мы их приводить не будем..□

На Рис. 166 семейство содержит два точечных подпространства. В стягиваемом семействе их пересечение также является точечным подпространством, состоящем из двух точек, соединенных связью. В нестягиваемом семействе пересечение является неточечным подпространством, состоящем из двух изолированных точек, то есть 0-мерной

Определение регулярного семейства подпространств.

Семейство $F=(G_1, G_2,...G_n)$ называется регулярным (regular), если для шара U(v) любой точки v пространства G найдется подпространство $G_k$ семейства, содержащее пересечение шара этой точки с объединением всех подпространств этого семейства, $(G_1 \cup G_2 \cup ... \cup G_n) \cap U(v) \subseteq G_k$ (Рис. 167).

В применении к покрытию это определение упростится.

Определение регулярного покрытия пространства его подпространствами.

Покрытие $F=(G_1, G_2,...G_n)$ называется регулярным (regular), если для шара любой точки пространства G найдется подпространство покрытия, содержащее этот шар.

На Рис. 167 в регулярном покрытии шар любой точки принадлежит одному из подпространств покрытия. В нерегулярном покрытии шары точек a и b не принадлежат полностью ни одному из подпространств покрытия. Покрытия состоят из двух подпространств. В регулярном покрытии шар любой точки принадлежит одному из подпространств покрытия. В нерегулярном покрытии шары точек a и b не принадлежат полностью ни одному из подпространств покрытия.

Определение правильного семейства подпространств.

Семейство $F=(G_1,G_2,...G_n)$ называется правильным (complete), если оно непрерывно, стягиваемо и регулярно.

Используя семейство подпространств, построим несколько новых пространств, обладающих теми же свойствами, что и исходное пространство.

Определение раздутия пространства G, соответствующего семейству F.

Пусть G, $F=(G_1,G_2,...G_n)$ и D есть пространство с точками $(v_1,v_2,...v_m)$, семейство его подпространств и молекулярное пространство с точками $(d_1,d_2,...d_n)$, являющееся нервом семейства F (каждая вершина $d_k$ соответствует подпространству $G_k$, и две точки $d_k$ и $d_p$ смежны тттк пересечение подпространств $G_k$ и $G_p$ непусто, $G_k \cap G_p \neq \varnothing$). Построим пространство V(G,F,D) с точками $(v_1,v_2,...v_m,d_1,d_2,...d_n)$, содержащее G и D как свои подпространства, и в котором связь между пространствами G и D определяется соотношением $O(d_k) \cap G = G_k$. Пространство V(G,F,D) (или V(G,D)) называется раздутием G, соответствующим семейству F.

Теорема 161

*Пусть G и $F=(G_1,G_2,...G_n)$ являются точечным пространством и стягиваемым м непрерывным семейством его подпространств. Тогда раздутие V(G,F,D) будет точечным молекулярным пространством.*

Доказательство

Используем индукцию. На пространствах с небольшим числом точек эта теорема проверяется непосредственно. Предположим, что теорема верна

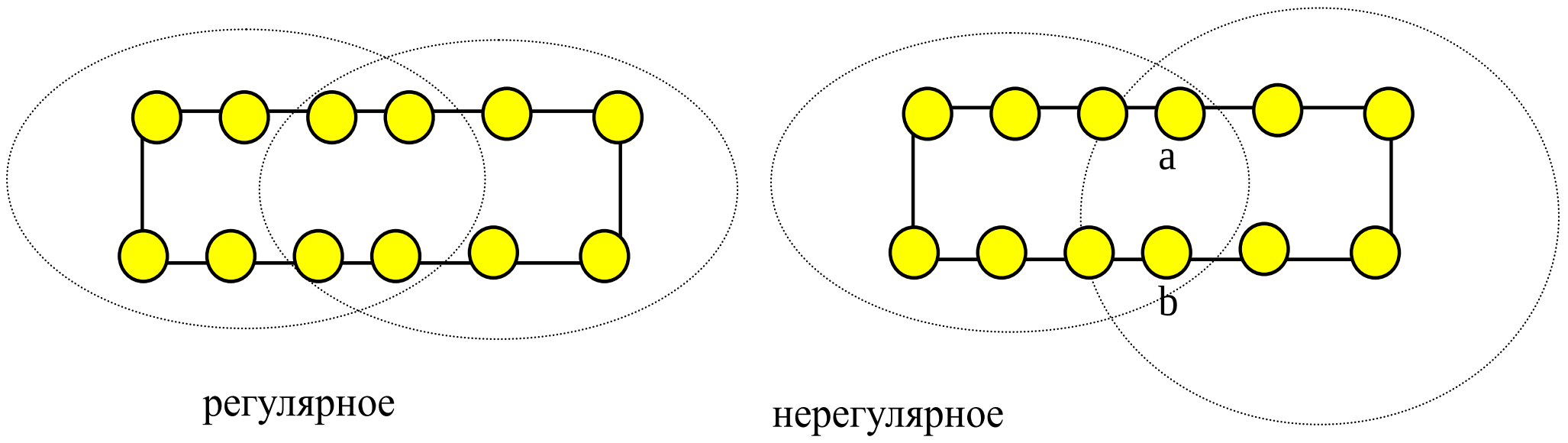

Рис. 167 В регулярном покрытии шар любой точки принадлежит одному из подпространств покрытия. В нерегулярном покрытии шары точек a и b не принадлежат полностью ни одному из подпространств покрытия.

на всех пространствах с числом точек, меньшим чем n. Пусть |G|=n. Приклеим к G точки $d_1,d_2,...d_n$ таким образом, чтобы $O(d_k)=G_k$., k=1,2,...n. Пока все точки $d_k$ несмежны. Полученное пространство гомотопно G. Если $O(d_1d_2)\neq\emptyset$, то этот общий окаем является точечным, и точки $d_1$ и $d_2$ могут быть соединены связью. Предположим, что мы уже установили все возможные связи между точками $d_1,d_2,...d_{s-1}$. Обозначим подпространство на этих точках как $H(s-1)=\{d_1,d_2,...d_{s-1}\}$. Берем точку $d_s$ и устанавливаем необходимые связи этой точки с точками $d_1,d_2,...d_{r-1}\in H(s-1)$. Предположим, что $G_s\cap G_r\neq\emptyset$. Тогда $O(d_sd_r)=V(G_s\cap G_r),H(s-1))$ есть раздутие над $G_s\cap G_r$. Так как $|G_s\cap G_r|<n$, то предположению индукции это пространство является точечным пространством. Следовательно, мы можем соединить точки $d_s$ и $d_r$. Продолжая этот процесс установим все необходимые связи. Так как все преобразования являются точечными, полученное пространство V(G,D) гомотопно исходному. Теорема доказана.

□

Обобщим предыдущую теорему на произвольное пространство G.

Теорема 162

*Пусть G и $F=(G_1,G_2,...G_n)$ являются произвольным пространством и стягиваемым м непрерывным семейством его подпространств. Тогда раздутие V(G,D) будет гомотопно пространству G.*

Д о к а з а т е л ь с т в о

Доказательство повторяет предыдущее. Используем индукцию. На пространствах с небольшим числом точек эта теорема проверяется непосредственно. Предположим, что теорема верна на всех пространствах с числом точек, меньшим чем n. Пусть $|G|=n$. Приклеим к G точки $d_1,d_2,...d_n$ таким образом, чтобы $O(d_k)=G_k$., $k=1,2,...n$. Пока все точки $d_k$ несмежны. Полученное пространство гомотопно G. Если $O(d_1d_2)\neq\emptyset$, то этот общий окаем является точечным, и точки $d_1$ и $d_2$ могут быть соединены связью. Предположим, что мы уже установили все возможные связи между точками $d_1,d_2,...d_{s-1}$. Обозначим подпространство на этих точках как $H(s-1)=\{d_1,d_2,...d_{s-1}\}$. Берем точку $d_s$ и устанавливаем необходимые связи этой точки с точками $d_1,d_2,...d_{r-1}\in H(s-1)$. Предположим, что $G_s\cap G_r\neq\emptyset$ и, следовательно, есть точечное пространство. Тогда $O(d_sd_r)=V(G_s\cap G_r),H(s-1))$ есть раздутие над $G_s\cap G_r$. Так как $|G_s\cap G_r|<n$, то предположению индукции это пространство является точечным пространством. Следовательно, мы можем соединить точки $d_s$ и $d_r$. Продолжая этот процесс установим все необходимые связи. Так как все преобразования являются точечными, полученное пространство V(G,D) гомотопно исходному. Теорема доказана. □

<u>Определение правильного покрытия пространства G.</u>

Покрытие $F=(G_1,G_2,...G_n)$ называется правильным (complete), если оно непрерывно, стягиваемо и регулярно.

<u>Определение накрытия D и накрывающего раздутия COV(G,F,D) пространства G для покрытия F.</u>

Пусть G, $F=(G_1,G_2,...G_n)$ и D есть пространство с точками $(v_1,v_2,...v_m)$, его покрытие и молекулярное пространство с точками $(d_1,d_2,...d_n)$, являющееся нервом покрытия F (каждая вершина $d_k$ соответствует подпространству $G_k$, и две точки $d_k$ и $d_p$ смежны тттк пересечение подпространств $G_k$ и $G_p$ непусто, $G_k\cap G_p\neq\emptyset$). Построим раздутие COV(G,F,D)=V(G,F,D). Пространство D называется накрытием пространства G, пространство COV(G,F,D) накрывающим раздутием (сокращенно накраз) пространства G, соответствующим семейству F, само пространство G называется базой накраза.

На Рис. 168 слева условно изображено покрытие $F=(G_1,G_2,G_3,G_4,G_5)$ некоторого пространства G, состоящее из пяти подпространств, справа изображено накрытие D. Возникает вопрос, когда пространство COV(G,F,D) гомотопно пространству G. Иными словами, когда пространство G может быть заменено другим пространством с меньшим числом точек, но с теми же свойствами. Следующая теорема дает ответ на этот вопрос.

Основная теорема.

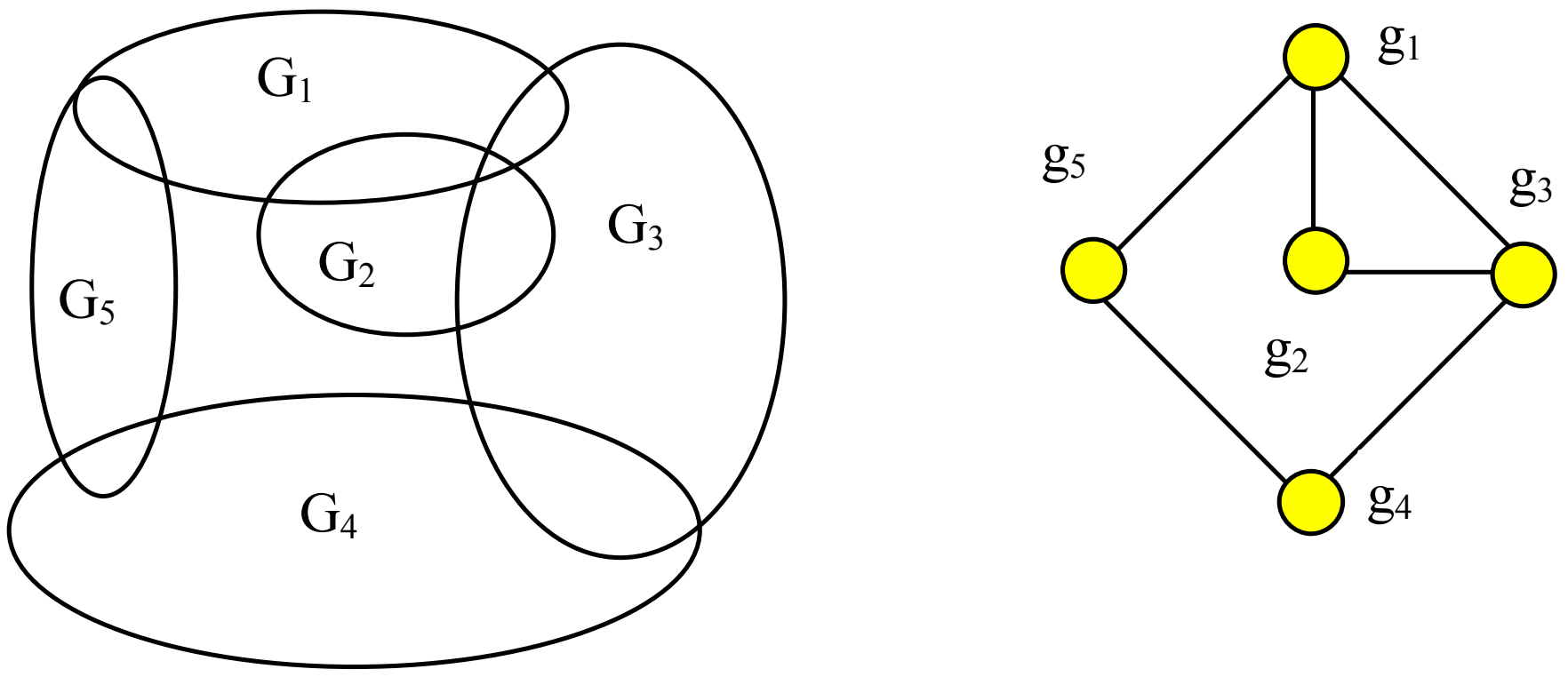

Рис. 168 Слева условно изображено покрытие $F=(G_1,G_2,G_3,G_4,G_5)$ некоторого пространства G, состоящее из пяти подпространств, справа изображено накрытие D

Теорема 163

*Пусть G, $F=(G_1,G_2,...G_n)$ и COV(G,F,D) есть пространство, его правильное покрытие и накрывающее раздутие пространства G. Тогда накрытие D гомотопно G, D~G.*

Доказательство

Согласно предшествующей теореме накраз COV(G,F,D) гомотопен G. Возьмем любую точку v из G и рассмотрим ее окаем O(v) в COV(G,F,D). Очевидно $O(v)=(O(v)\cap G)\cup(O(v)\cap D)$. Так как F регулярное покрытие, то существует такое подпространство $G_k$, что шар $U(v)\in G_k$. Это оэначает, что точка $d_k$ из D является смежной всем точкам из D, принадлежащим O(v), поскольку, если $d_s\in O(v)$, то $v\in G_s$ и, следовательно, $G_s\cup\cap G_k\neq\emptyset$ и $d_s$ смежна с $d_k$. Следовательно, O(v) есть конус в COV(G,F,D), и точка v может быть отброшена. Таким путем мы отбрасываем все точки, принадлежащие подпространству G из COV(G,F,D) и получаем D. Отсюда следует, что COV(G,F,D) гомотопно D. Поскольку COV(G,F,D) гомотопно G, то G гомотопно D, G~COV(G,F,D)~D, G~D. Теорема доказана. □

Важным следствием основной теоремы является еще одна теорема, позволяющая установить связь между покрывающими пространствами пространства G для правильных покрытий.

Теорема 164

*Пусть имеется пространство G и два его правильных покрытия $F_1$ и $F_2$.. Тогда*

- *Накрывающие раздутия $COV(G,F_1,D_1)$ и $COV(G,F_2,D_2)$ гомотопны, $COV(G,F_1,D_1)$~$COV(G,F_2,D_2)$.*
- *Накрытия $D_1$ и $D_2$ гомотопны, $D_1$~$D_2$.*

Доказательство.

G гомотопно как COV(G,F,D) так и COV(G,F,D). Следовательно, $COV(G,F_1,D_1)$~$COV(G,F_2,D_2)$. Кроме того, G гомотопно как $D_1$ так и $D_2$. Следовательно, $D_1$~$D_2$. □

Смысл этой теоремы заключается в том, накрытие и накрывающее

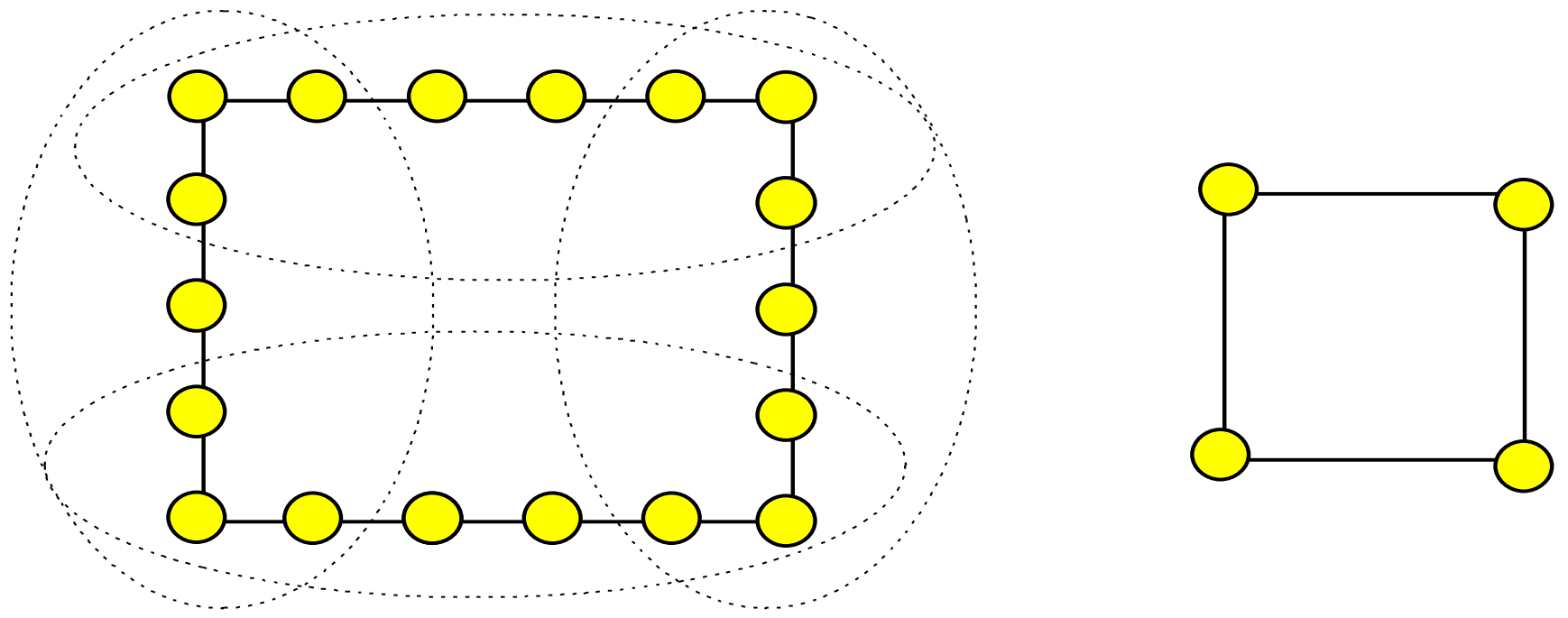

Рис. 169 Исходное пространство является окружностью, состоящей 18-и точек и покрыто четырьмя подпространствами. Данное покрытие является правильным. Легко видеть, что это покрытие дает покрывающее пространство, гомотопное исходному, и состоящее из 4-х точек.

раздутие имеют те же фундаментальные топологические свойства, что и пространство G. С этой точки зрения совершенно безразлично, какое накрытие выбрать вместо исходного.

На Рис. 169 изображены исходное пространство, являющееся окружностью, состоящей 18-и точек и покрывающие его 4-е подпространства. Несложно проверить, что данное покрытие является правильным. Легко видеть, что это покрытие дает покрывающее пространство, гомотопное исходному, и состоящее из 4-х точек. Данная теорема дает достаточно сильные условия для существования покрывающего гомотопного пространства. Легко придумать пример покрытия, не являющегося правильным, тогда как его пространство пересечений гомотопно данному. Очевидно, что достаточно рассматривать покрытие, состоящее из подпространств, которые можно несколько "раздуть", превратив в правильное покрытие.

Теорема 165

*Пусть G и $F_1=(G_1,G_2,...G_n)$ есть пространство и его некоторое его покрытие. Если существует правильное покрытие $F_2=(H_1,H_2,...H_n)$, каждый элемент $H_k$ которого содержит $G_k$, $G_k \subseteq H_k$, $k=1,2,...n$, и $H_k \cap H_p \neq \varnothing$ тттк $G_k \cap G_p \neq \varnothing$, тогда нерв $D_1$ покрытия $F_1$ гомотопен исходному пространству G, $D_1 \sim G$.*

Доказательство.

G гомотопно нерву $D_2$ покрытия $F_2=(H_1,H_2,...H_n)$. С другой стороны пространства $D_1$ и $D_2$ изоморфны, поскольку элементы из $F_1$ и $F_2$ пересекаются или не пересекаются одновременно. Следовательно, $G \sim D_1 \sim D_2$. □

Например, во всех рассмотренных нами конкретных случаях достаточно выполнения следующих условий:

Пусть G, $F=(G_1,G_2,...G_n)$ и COV(G,F,D) есть пространство, его покрытие и накрывающее раздутие пространства G для покрытия F. Если покрытие F:

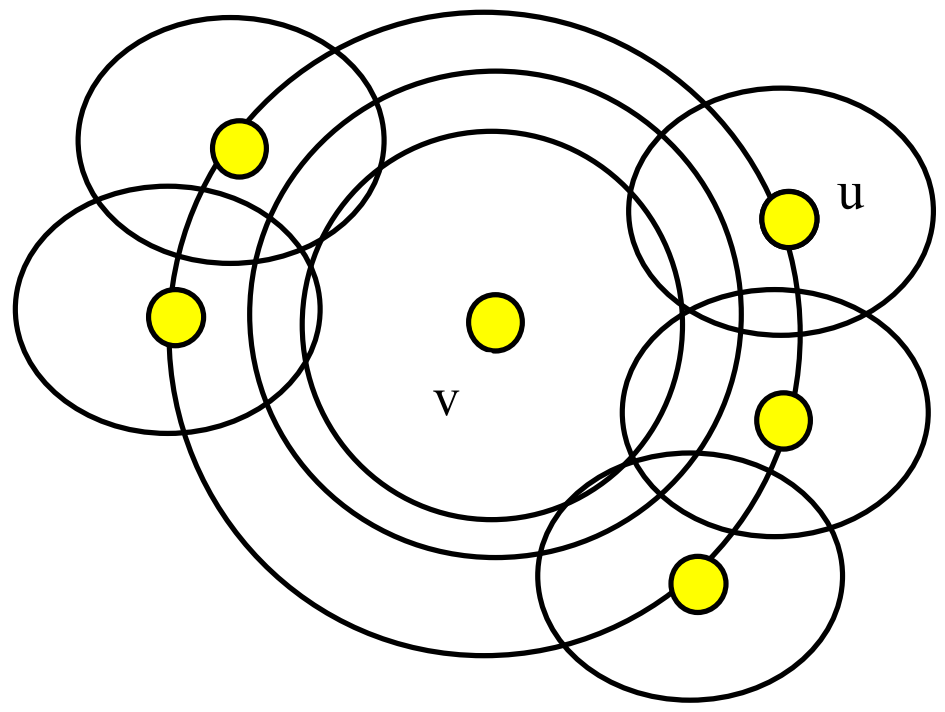

Рис. 170 Алгоритм построения накрытия.

- непрерывно,
- стягиваемо,
- любые две смежные точки принадлежат хотя бы одному подпространству покрытия, то COV(G,F,D) гомотопен G, COV(G,F,D)~G.

Алгоритм построения накрытия.

Опишем теперь простейший алгоритм, позволяющий получить правильное покрытие. Напомним определение 2-шара точки v пространства G. 2-шаром $B^2(v)$ точки v пространства G называется подпространство, содержащий все точки пространства G, соседние точкам шара B(v) точки v, включая точки шара B(v). Аналогично определяется n-шар точки v. Используем 2-шары как основу для построения алгоритма.

- Выбираем произвольную точку v пространства G и ее 2-шар $B^2(v)$.

- Выбираем все точки, находящиеся на расстоянии 3 от v, и выбираем 2-шары этих точек
- Выбираем все точки, находящиеся на расстоянии 6 от v, и выбираем 2-шары этих точек.
- Выбираем все точки, находящиеся на расстоянии 9 от v, и выбираем 2-шары этих точек. Продолжаем этот процесс, пока не покроем пространство подпространствами, являющимися выбранными 2-шарами.
- Проверяем полученное покрытие на правильность.
- Строим накрытие для покрытия.

На Рис. 170 условно изображен первый этап этого алгоритма.

З а д а ч а

Найти другие алгоритмы нахождения правильного покрывающего пространства.

## *ПОСТРОЕНИЕ МОЛЕКУЛЯРНЫХ МОДЕЛЕЙ НЕПРЕРЫВНЫХ ПРОСТРАНСТВ ПРИ ПОМОЩИ ПОКРЫТИЙ*

### 3.1.1СИМПЛИЦИАЛЬНЫЙ КОМПЛЕКС КАК МОЛЕКУЛЯРНАЯ МОДЕЛЬ ПОЛИЭДРА.

Мы рассмотрим способы построения молекулярного пространства по его непрерывному прообразу. Вообще говоря, описываемые здесь методы применимы к любым множествам, определенным в непрерывном n-мерном евклидовом пространстве. Требование, которому должны удовлетворять молекулярные пространства-соответствие свойств непрерывного прообраза и его молекулярного образа. Для того, чтобы было ясно, как обеспечиваются эти требования, рассмотрим идеализированную схему построения молекулярной модели. В математике изучаются так называемые полиэдры, являющиеся некоторым множеством евклидова пространства, которое может быть представлено в виде симплициального комплекса. Попробуем применить симплициальные комплексы для получения молекулярной модели полиэдра. Предположим, что мы имеем некоторое пространство S, которое является симплициальным полиэдром в непрерывном n-мерном евклидовом пространстве $E^n$. Строим для S всевозможные симплициальные комплексы, являющиеся разбиением пространства S. Обозначим семейство всех таких комплексов K(S). Рассматриваем все симплициальные комплексы из K(S) как молекулярные пространства, считая, что любой симплекс является связкой. При достаточно мелком разбиении, что всегда можно сделать при помощи так называемого барицентрического подразделения, получившиеся молекулярные пространства будут гомотопны одно другому и могут рассматриваться как молекулярные модели полиэдра. При этом, конечно, следует иметь ввиду, что не каждый симплициальный комплекс может рассматриваться как молекулярный образ данного пространства, а также что многие понятия, такие как размерность симплекса, неприменимы к молекулярному пространству. Эта схема почти очевидна, выполняется во всех конкретных случаях. Ее доказательство представляет интерес для математического обоснования верности этого подхода. Один путь доказательства-показать, что одно из симплициальных разбиений является правильным покрытием. С практической точки зрения эта схема мало интересна, поскольку построение симплициального комплекса само по себе является достаточно сложной задачей.

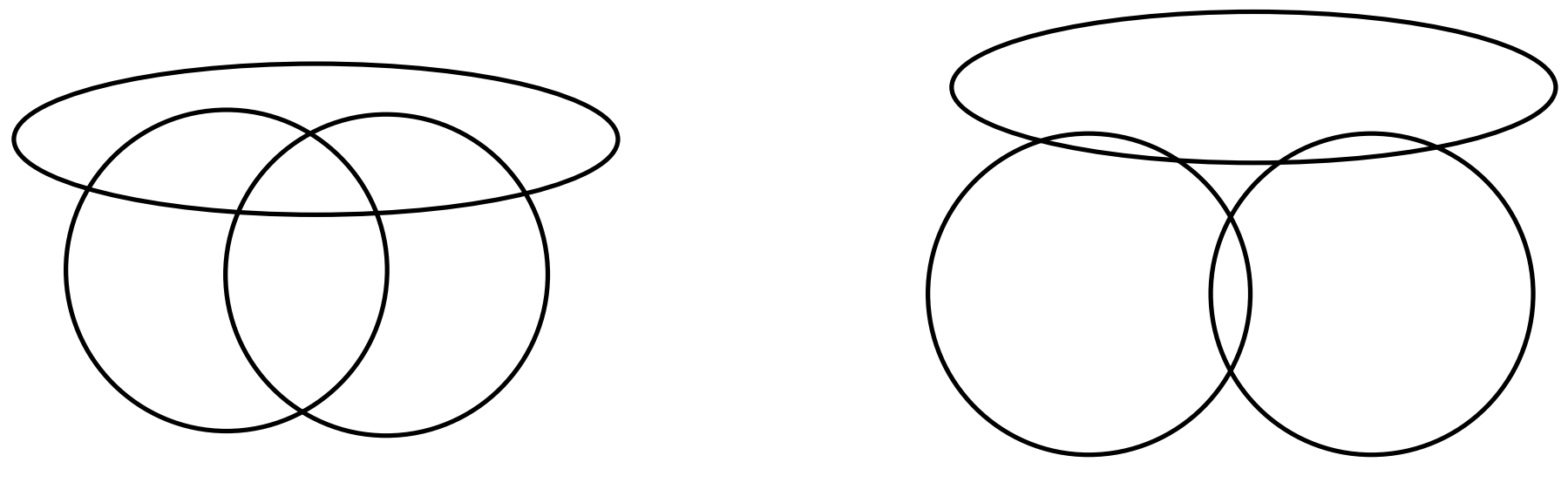

Рис. 171 Непрерывное и ненепрерывное покрытия.

З а д а ч а .

Сформулировать и доказать теорему о верности этой схемы.

## 3.1.2МОЛЕКУЛЯРНЫЕ ПРОСТРАНСТВА КАК НЕРВЫ ПОКРЫТИЙ НЕПРЕРЫВНОГО ПРОСТРАНСТВА.

Мы рассмотрим еще один способ получения молекулярного образа некоторого множества, включая топологическое пространство. Пусть дано множество X и его некоторое его покрытие. Предположим, что мы имеем два конечных покрытия пространства X, $F_1(a_1, a_2, ... a_n)$ и $F_2(b_1, b_2, ... b_n)$.

Возьмем покрытие $F_3(a_1, a_2, ... a_n, b_1, b_2, ... b_n)$, содержащее элементы обоих покрытий. Построим молекулярные пространства $G(F_1)$, $G(F_2)$ и $G(F_3)$, являющиеся нервами покрытий. Очевидно, что можно считать, что $G=G(F_1)$, $D=G(F_2)$ и $COV(G,F,D)=G(F_3)$, или $G=G(F_2)$, $D=G(F_1)$ и $COV(G,F,D)=G(F_3)G(F_2)$. Если $G(F_1)$ можно считать правильным покрытием $G(F_2)$ (или наоборот), то $G(F_1)$ и $G(F_2)$ являются гомотопными. В этом смысле мы можем считать эти покрытия эквивалентными дигитальными моделями пространства X. Заключим эти простые рассуждения в форму теоремы. Важной частью этого подхода является степень соответствия свойств нерва покрытия пространства и самого топологического пространства. Решать эту проблему можно несколькими путями. Можно считать симплициальный комплекс, рассматриваемый как молекулярное пространство, как истинную дигитальную модель топологического пространства. Тогда нервы покрытий сравниваются с

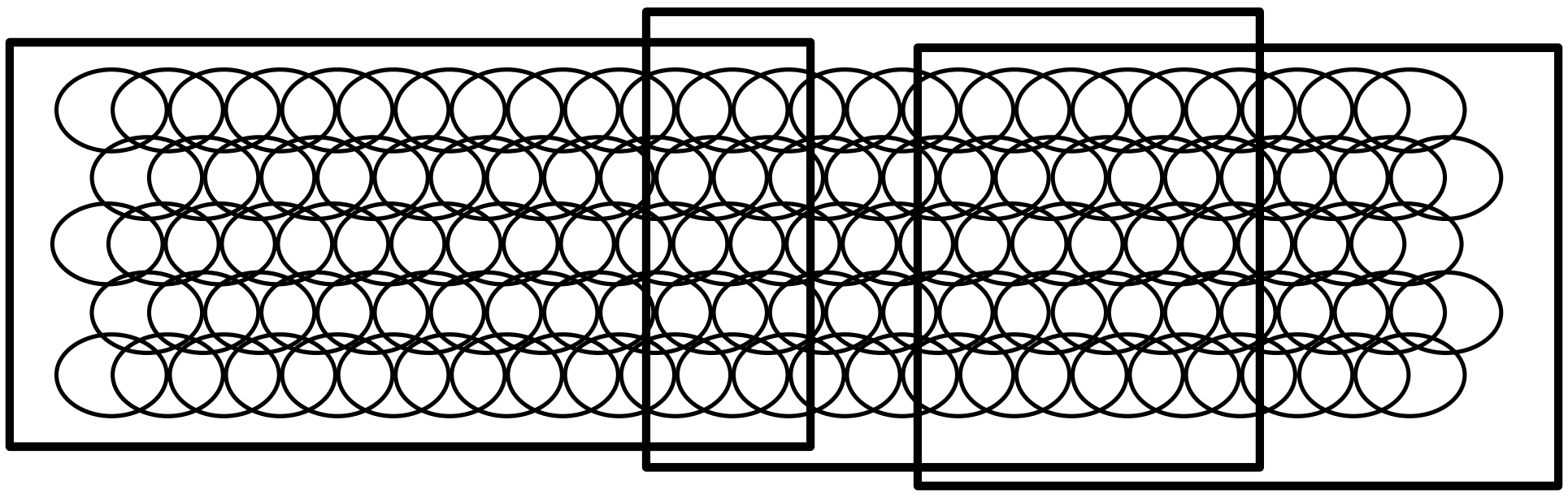

Рис. 172 Элементами покрытия Q являются круги, а элементами покрытия P являются прямоугольники. Легко проверить, что все условия теоремы выполняются, и нерв покрытия P будет правильным накрытием для нерва покрытия Q, следовательно, они гомотопны.

этой моделью. Этот путь по известным причинам, представляет, скорее, интерес для теоретических обоснований, нежели для практической работы. Второй путь заключается в измельчении элементов покрытия и проверке получившихся нервов покрытий на гомотопность. Если при переходе к покрытиям со все более мелкими элементами нервы покрытий остаются гомотопными молекулярными пространствами, то каждый из них можно считать дигитальной моделью топологического пространства. Конечно и в таком подходе имеется много подводных камней. Дело в том, что, например, изменение формы элемента покрытия может существенно изменить смежность точек в нерве покрытия. Поэтому необходимо заранее ввести какое-то ограничения на вид элементов покрытия.

Критерий истинности дигитальной модели топологического пространства

Пусть имеется последовательность $F_a=(G(a)_1,G(a)_2,...G(a)_n,...)$, $a\in A$, покрытий некоторого топологического пространства X, обладающих следующими особенностями:

- Все элементы каждого из покрытий принадлежат некоторому определенному классу K,
- Если $a_1<a_2$, то для любого $G(a_1)_k$ существует $G(a_2)_p$ такой, что $G(a_1)_k\subseteq G(a_2)_p$, и не существует никакого $G(a_2)_t$ такого, что $G(a_2)_t\subseteq G(a_1)_k$.

Тогда, если для всех $a<a_0$ нервы покрытий $F_a=(G(a)_1,G(a)_2,...G(a)_n,...)$ гомотопны один другому, то каждый из этих нервов считается истинной молекулярной моделью топологического пространства X.

Определение гомотопности покрытий топологического пространства X.

Предположим, что мы имеем два конечных покрытия пространства X, $F_1(a_1,a_2,...a_n)$ и $F_2(b_1,b_2,...b_n)$. Они называются гомотопными, если молекулярные пространства $G(F_1)$ и $G(F_2)$, являющиеся нервами этих покрытий, гомотопны.

Возникает вопрос, каким условиям должны удовлетворять покрытия, чтобы быть гомотопными. Введем несколько определений.

Определение непрерывного покрытия множества X.

Покрытие $P=(p_1,p_2,...p_n)$ множества X называется непрерывным (continuous), если из условия $p_{i_k} \cap p_{i_m} \neq \varnothing, \forall\, k, m = 1, 2,...,r$ следует, что $p_{i_1 i_2 ... i_r} = p_{i_1} \cap p_{i_2} ... \cap p_{i_r} \neq \varnothing$.

На Рис. 171 изображены непрерывное и ненепрерывное покрытия.

Теорема 166

*Пусть мы имеем два покрытия $P(p_1,p_2,...p_n)$ и $Q(q_1,q_2,..q_m)$ некоторого множества X со следующими свойствами:*

1. *Покрытие P является непрерывным.*
2. *Элементы из Q, покрывающие любое непустое подмножество $A = p_{i_1 i_2 ... i_r} = p_{i_1} \cap p_{i_2} ... \cap p_{i_r} \neq \varnothing$, имеют нерв, являющийся точечным молекулярным пространством.*
3. *Для любого элемента $q_i$ покрытия Q всегда существует элемент $p_s$ покрытия P такой, что $q_i$ и все пересекающие его элементы покрытия Q целиком принадлежат $p_s$.*

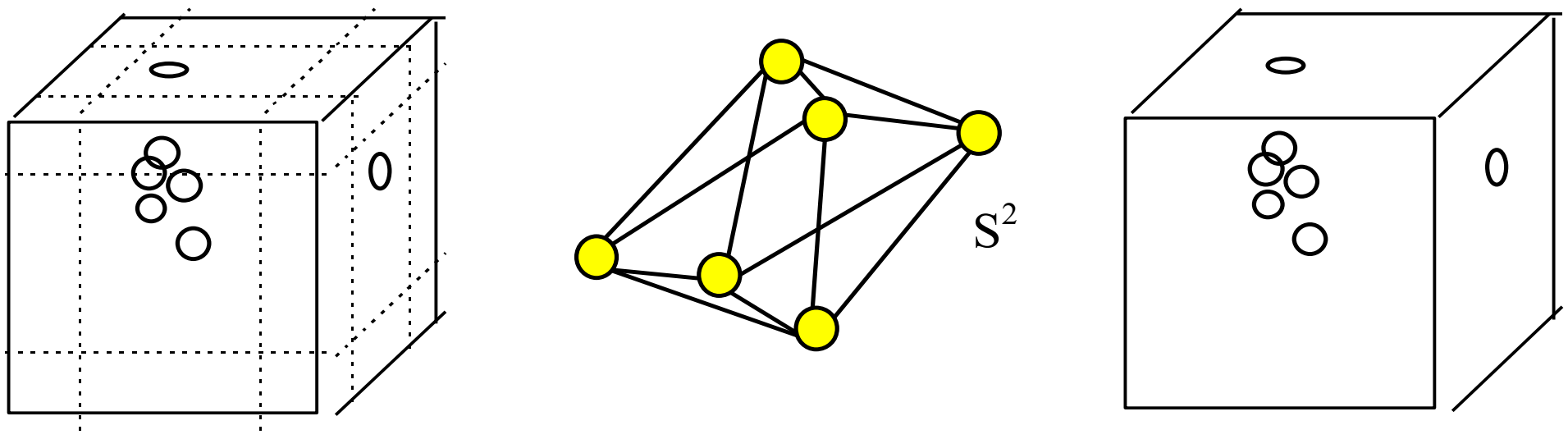

Рис. 173 Покрытие P куба, состоит из шести элементов, каждый из которых содержит свою грань и перекрывает четыре соседних к нему на полосу шириной 2d. Покрытие Q состоит из небольших кругов диаметра d. Тогда в соответствии с теоремой нервы этих двух покрытий гомотопны. Нерв покрытия P является минимальной двумерной сферой.

*Тогда нерв покрытия P гомотопен нерву покрытия Q.*

Доказательство

Рассмотрим покрытие PQ, являющееся объединением покрытий P и Q, и построим нерв G(PQ) этого покрытия (для нервов используем те же обозначения, что и для исходных покрытий). Этот нерв есть молекулярное пространство, в котором молекулярное подпространство G(Q) (нерв покрытия Q) есть база, G(P) (нерв покрытия P) есть накрытие для G(Q), и

G(PQ) есть накрывающее раздутие пространства G(Q). Обозначим $G_k=O(p_k)\cap G(Q)$. Легко видеть, что семейство $F=(G_1,G_2,...G_n)$ является покрытием молекулярного пространства G(Q). Из условия 2 следует, что это покрытие является стягиваемым. Из условия 1 следует, что F является непрерывным. Из условия 3 следует, что F является регулярным. Это означает, что покрытие является правильным и, следовательно, молекулярные пространства G(P) и G(Q) гомотопны. Теорема доказана. □

Сущность этой теоремы очень проста. Имеется два покрытия P и Q множества X. Строим нервы G(P) и G(Q) каждого из этих покрытий. Тогда G(P) будет правильным накрытием для G(Q) и, следовательно, гомотопен G(Q) на основании теоремы, доказанной ранее.

На Рис. 172 элементами покрытия Q являются круги, а элементами покрытия P являются прямоугольники. Легко проверить, что все условия

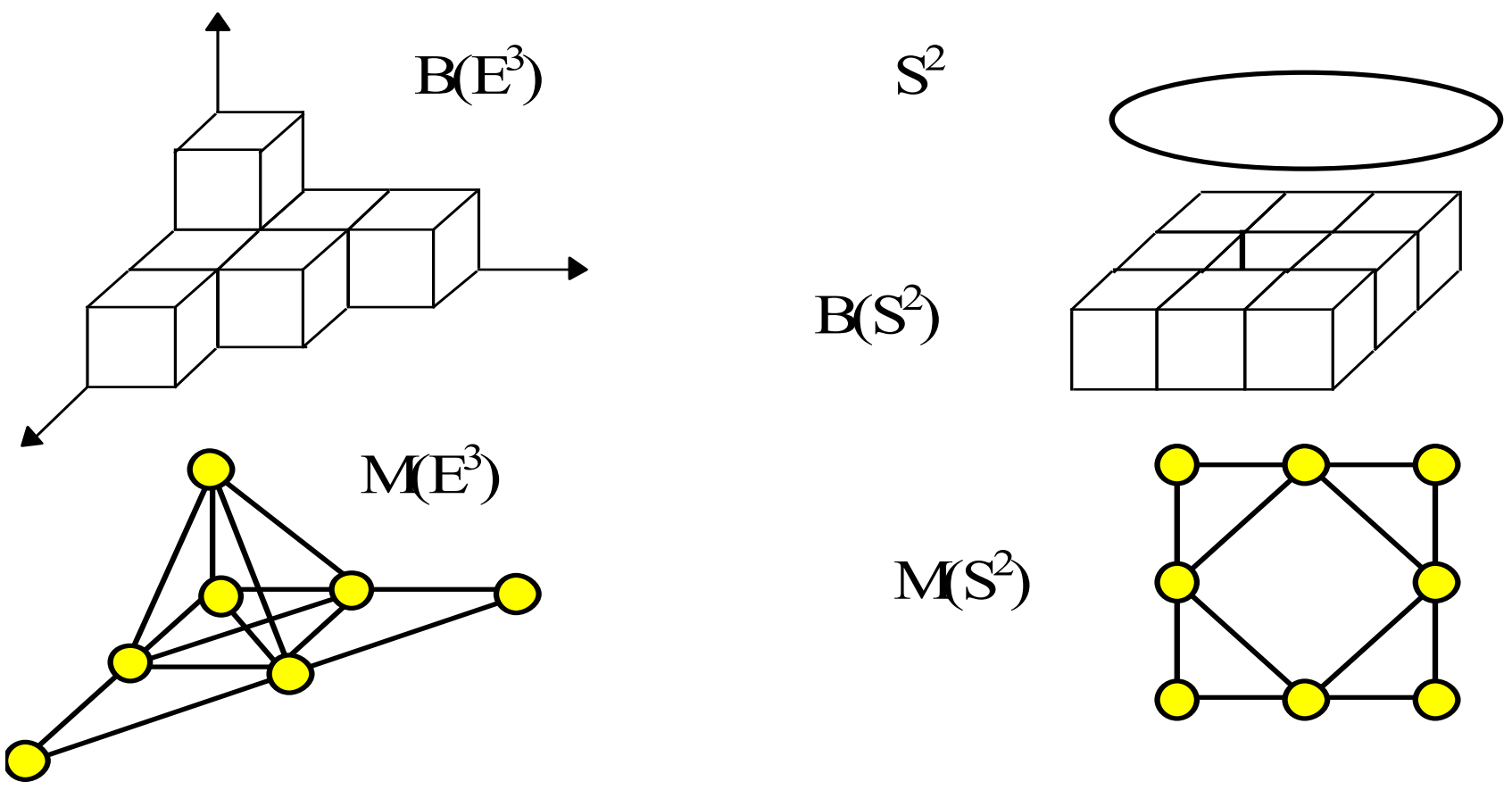

Рис. 174 Базовые кирпичное и молекулярное пространства для трехмерного евклидова пространства $E^3$. Координаты вершин кубов являются числами, кратными $l_0$.

теоремы выполняются, и нерв покрытия P будет правильным накрытием для нерва покрытия Q, следовательно, они гомотопны. Любой из этих нервов может рассматриваться как молекулярная модель множества X, которое в данном случае является частью двумерной плоскости, покрытой кругами. Эта сравнительно простая теорема позволяет установить связь двумя покрытиями в общем виде. Рассмотрим, например, покрытие Q сферы небольшими выпуклыми многоугольниками или кругами. Мы сейчас покажем, что нерв такого покрытия всегда является гомотопным минимальной двумерной сфере. Для этого рассмотрим поверхность единичного куба (Рис. 173), которая гомеоморфна двумерной сфере.

Покроем поверхность куба кругами диаметром меньше d, d<<1. Несложно показать, что покрытие каждой грани куба имеет нерв, который является гомотопным точечному молекулярному пространству. Рассмотрим теперь покрытие P куба, состоящее из шести элементов, каждый из которых содержит свою грань и перекрывает четыре соседних к нему на полосу шириной 2d. Эти элементы напоминают крышку картонной коробки и изображены пунктиром. Тогда в соответствии с теоремой нервы этих двух покрытий гомотопны. Нерв покрытия P является минимальной двумерной сферой. Следовательно, нерв покрытия Q есть также двумерная сфера. Точно такие же рассуждения применимы к любой замкнутой двумерной поверхности.

Рассмотрим одно простое, но важное с практической точки зрения следствие из доказанной выше теоремы, которое позволяет установить связь между двумя покрытиями, которые не удовлетворяют условиям теоремы непосредственно.

Следствие.

Пусть мы имеем два покрытия $R(r_1, r_2,...r_m)$ и $Q(q_1, q_2,...q_m)$ некоторого множества X. Пусть мы имеем покрытие $P(p_1, p_2,...p_n)$ того же X, такое, что каждая пара покрытий $P(p_1, p_2,...p_n)$ и $Q(q_1, q_2,...q_m)$ и $P(p_1, p_2,...p_n)$ и $R(r_1, r_2,...r_m)$ удовлетворяет условиям теоремы, доказанной выше. Тогда нерв G(Q) покрытия Q гомотопен нерву G(R) покрытия R.

## *ПОСТРОЕНИЕ МОЛЕКУЛЯРНЫХ МОДЕЛЕЙ НЕПРЕРЫВНЫХ ПРОСТРАНСТВ ПРИ ПОМОЩИ БАЗОВЫХ ЭЛЕМЕНТОВ*

### 3.1.ЗКРИТЕРИИ ВЫБОРА БАЗОВЫХ ЭЛЕМЕНТОВ, ФОРМИРУЮЩИХ МОЛЕКУЛЯРНЫЙ ОБРАЗ НЕПРЕРЫВНОГО ПРОСТРАНСТВА.

Результаты, полученные выше, и позволяющие строить молекулярное пространство по покрытию топологического пространства, ясны с принципиальной точки зрения, Однако, вопрос конкретного выбора покрытия остается открытым. При этом необходимо учитывать следующие обстоятельства:

- Метод выбора покрытия должен быть универсальным, применимым к любой практической задаче.
- Покрытие должно быть удобным для автоматического его анализа на компьютере.

Поэтому мы должны предложить такой метод, который отвечал бы указанным требованиям, и не был бы слишком сложен для понимания и использования. По сравнению со способом получения молекулярного пространства, рассмотренным в предыдущей главе, предлагаемый способ является более узким, но позволяющим избежать субъективный этап выбора покрытия и осуществлять весь процесс построения молекулярного пространства автоматически на компьютере.

## БАЗОВЫЕ КИРПИЧНОЕ И МОЛЕКУЛЯРНОЕ ПРОСТРАНСТВА.

Как всегда, введем несколько основных определений. Ранее мы ввели кирпичное n-мерное пространство как подпространство универсального базового пространства $B(E^{\infty})$. Несколько изменим масштаб, рассматривая вместо единичных кирпичей произвольные кирпичи с длиной стороны $l_0$.

Определение базового кирпичного n-мерного пространства $B(E^n)$.

Пусть имеется евклидово пространство $E^n$. Рассечем $E^n$ на множество одинаковых кубов (кирпичей) с длиной ребра $l_0$ и с ребрами, параллельными координатным осям. Назовем семейство кубов базовым кирпичным пространством евклидова пространства $E^n$ (под элементами в данном случае понимаются кубы).

Обозначим базовое пространство $B(E^n)$ или просто $B^n$.

Определение базового молекулярного n-мерного пространства $B(E^n)$.

Нерв $M(E^n)$ семейства $B(E^n)$ называется базовым молекулярным n-мерным пространством евклидова пространства $E^n$ и обозначается $M(E^n)$ или $M^n$.

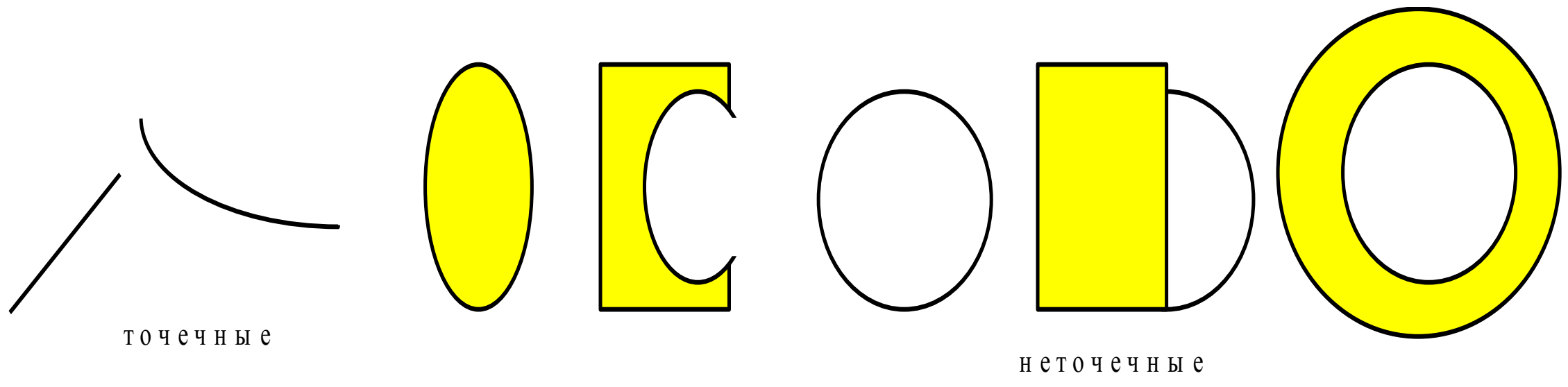

Рис. 175 Множество A называется точечным, если существует δ>0 такое, что для всех $l_0<\delta$ базовое молекулярное пространство M(A) множества A является точечным.

Легко убедиться в том, что $M(E^n)$ является n-мерным молекулярным регулярным пространством, хотя и отсутствует свойство нормальности. Заметим, что два элемента (куба) считаются соседними, если они имеют хотя бы одну общую точку. На Рис. 174 изображены базовые кирпичное и молекулярное пространства для трехмерного евклидова пространства $E^3$. Для простоты будем считать, что координаты вершин кубов являются числами, кратными $l_0$.

Определение координат куба.

Пусть куб содержит точки, задаваемые координатами $(x_1,x_2,...x_n)$ пространства $E^n$, определяемые условиями

$l_0 k_i \le x_i \le l_0(k_i + 1),\ i = 1,2,\dots n,$ $-k_i$ - целые числа. Тогда координатами куба называются целые числа $(k_1, k_2, \dots k_n)$, определяющие левую границу.

Очевидно, что координаты соседних кирпичей отличаются не более чем на единицу. Пусть теперь мы имеем некоторое множество A в $E^n$.

Определение базовой кирпичной поверхности B(A).

Базовой кирпичной поверхностью B(A) поверхности A в $E^n$ называется семейство кубов, пересекаемых этим множеством.

Определение базовой молекулярной поверхности M(A).

Нерв M(A) поверхности B(A) называется базовым молекулярным пространством поверхности A и обозначается M(A).

Рис. 174 содержит окружность в трехмерном пространстве, ее базовое кирпичное пространство и базовое молекулярное пространство этой окружности.

## ТОЧЕЧНОЕ МНОЖЕСТВО В ЕВКЛИДОВОМ ПРОСТРАНСТВЕ.

Определение точечного множества в евклидовом пространстве.

Множество A называется точечным, если существует $\delta>0$ такое, что для всех $l_0<\delta$ базовое молекулярное пространство M(A) множества A является точечным. Число $\delta$ будем называть степенью точечности множества A (Рис. 175).

Согласно этому определению, нерв семейства кубов, пересекающих множество A, стягивается в точку. С классической точки зрения любое точечное множество гомотопно точке.

Определение блока в $E^n$.

Блоком в $E^n$. называется прямой прямоугольный параллелепипед, определяемый координатами $(x_1, x_2, \dots x_n)$ пространства $E^n$, где $a_i \le x_i \le b_i,\ a_i < b_i,\ i = 1,2,\dots n.$ Координатами блока называется последовательность пар вида $((a_1,b_1),(a_2,b_2),\dots(a_n,b_n))$.

З а м е ч а н и е

Предполагается, что некоторые $a_k$ и $b_k$ могут совпадать. В этом случае блок будет иметь размерность, меньшую чем n, но никаких ограничений на ход дальнейшего рассмотрения это не накладывает.

Теорема 167

*Любой блок A в $E^n$ есть точечное пространство.*

Доказательство

Пусть длина ребра кубов в базовом кирпичном пространстве есть $l_0 << \min|b_k - a_k|$. Построим M(A) и рассмотрим все точки, имеющие наибольшую первую координату. Очевидно, что окаем каждой такой точки $(k_1, k_2, \ldots k_n)$ является конусом с вершиной в точке $((k_1-1), k_2, \ldots k_n)$, первая координата которой на 1 меньше первой координаты исходной точки, а все остальные координаты этих точек совпадают. Следовательно, все такие точки можно отбросить. Идя по слоям, отбрасываем все точки из M(A), кроме последней. Теорема доказана. □

Следствие

Точки в верхнем слое блока можно отбрасывать в произвольном порядке. Любая точка в верхнем слое блока, из которого отброшено некоторое количество точек, всегда имеет точечный окаем и может быть отброшена.

Напомним, что в выпуклом множестве A две любые точки, принадлежащие множеству, можно соединить отрезком прямой, также принадлежащим этому множеству. Следовательно, выпуклыми

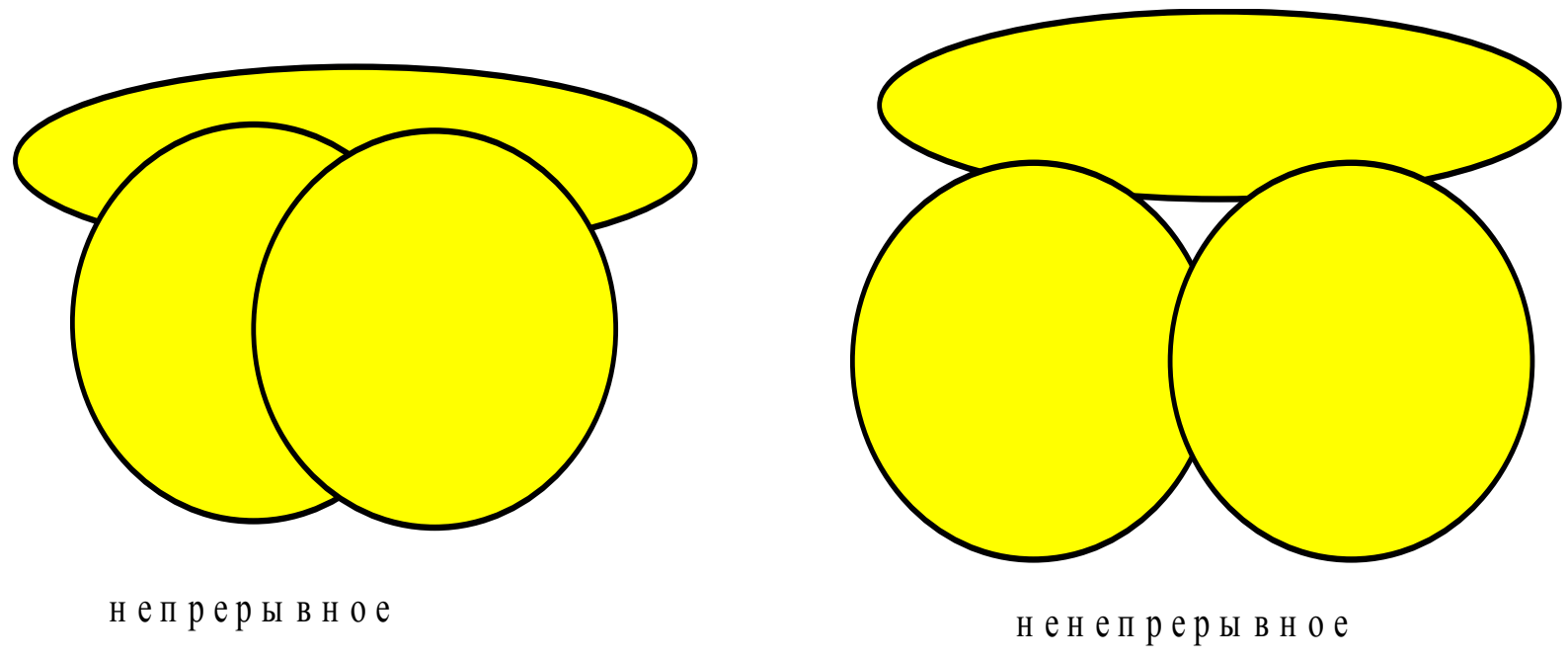

Рис. 176 Непрерывное и ненепрерывное семейства множеств.

множествами будут отрезок прямой, прямоугольник, круг, шар, n-мерный параллелепипед, n-мерный шар и т. д.. На Рис. 175 изображены точечные и неточечные множества.

Теорема 168

*Любое выпуклое множество A в $E^n$ есть точечное пространство.*

Доказательство

Пусть длина диаметра множества A удовлетворяет условию $l_0 << \mathrm{diam}A$, и диаметр множества A параллелен первой координатной оси. Построим M(A) и рассмотрим все точки, имеющие наибольшую первую координату. Очевидно, что окаем каждой такой точки является конусом с вершиной в точке, первая координата которой на 1 меньше первой координаты исходной точки, а все остальные координаты этих точек совпадают. В силу выпуклости множества A это условие выполняется по крайней мере на нескольких слоях, начиная с верхнего. Когда это условие нарушается, выберем точки с наименьшей первой координатой и повторим процесс отбрасывания. Идя, как и в предыдущем случае, по слоям, отбрасываем все точки из M(A), кроме последней. Теорема доказана. □

В общем случае, повидимому, можно доказать следующую теорему.

Теорема 169

*Любое множество A в $E^n$, гомеоморфное некоторому выпуклому множеству, есть точечное пространство.*

Мы не будем доказывать эту теорему, поскольку для стандартных поверхностей теорема очевидна.

## НЕПРЕРЫВНОЕ, СТЯГИВАЕМОЕ, РЕГУЛЯРНОЕ И ПРАВИЛЬНОЕ СЕМЕЙСТВА МНОЖЕСТВ. СВОЙСТВА ПРАВИЛЬНОГО СЕМЕЙСТВА МНОЖЕСТВ.

В дальнейшем будем рассматривать только семейства, состоящие из точечных множеств в евклидовом пространстве $E^n$. Дадим несколько определений, аналогичных определениям различных покрытий молекулярного пространства.

Определение непрерывного семейства множеств.

Семейство множеств $F=(A_1,A_2,..A_n)$ (Рис. 176) называется непрерывным (continuous), если из условия $A_{p_k} \cap A_{p_m} \neq \varnothing,\ \forall\, k, m = 1, 2,\dots,r$ следует, что $A_{p_1} \cap A_{p_2} \dots \cap A_{p_r} \neq \varnothing$.

Определение стягиваемого семейства множеств.

Семейство $F=(A_1,A_2,...A_n)$ называется стягиваемым (contractible) (Рис. 177), если пересечение любого числа множеств является либо точечным либо пустым.

Определение степени точечности семейства множеств.

Степенью точечности δ семейства F называется максимальная степень точечности пересечений множеств семейства.

Из этого определения следует, что все множества семейства являются также точечными. На Рис. 177 слева изображено стягиваемое семейство из 2-х множеств, справа-нестягиваемое, так как пересечение двух множеств семейства состоит из двух несвязных компонент и, следовательно, не является точечным.

Определение регулярного семейства F.

Семейство $F=(A_1,A_2,...A_n)$ называется регулярным (regular) (Рис. 178), если существует ε>0 такое, что для любой точки

$p \in H = A_1 \cup A_2 \cup ... \cup A_n$ и для любого куба K с центром в p и с $l_0 < \varepsilon$ существует множество $A_s$ такое, что $K \cap H \subseteq A_s$. Число $\varepsilon$ называется степенью регулярности семейства.

На Рис. 178 слева показано нерегулярное семейство. Пересечение объединения множеств с квадратом с центром в точке a не принадлежит ни одному из них. Справа показано регулярное семейство из двух множеств. Пересечение объединения множеств с квадратом с центром в точке a принадлежит любому из них. Мы будем рассматривать семейства обладающие всеми тремя свойствами.

Определение правильного семейства F.

Семейство F называется правильным (complete), если оно непрерывно, стягиваемо и регулярно. Характеристикой правильного семейства называются степени точечности и регулярности $(\delta,\varepsilon)$.

Наиболее простым примером правильного семейства множеств в n-мерном пространстве $E^n$ является семейство прямоугольных n-мерных параллелепипедов с ребрами, параллельными координатным осям, в которых пересечение любых двух параллелепипедов или пусто или является также невырожденным прямоугольным n-мерным параллелепипедом, как это показано для двумерного пространства Рис. 179. На этом рисунке также изображено неправильное семейство и базовое кирпичное пространство на двумерном евклидовом пространстве $E^n$. Легко видеть, что в регулярном семействе можно сделать так, чтобы не только куб K принадлежал некоторому множеству семейства, но и все соседние ему принадлежали данному множеству семейства. Для этого нужно взять

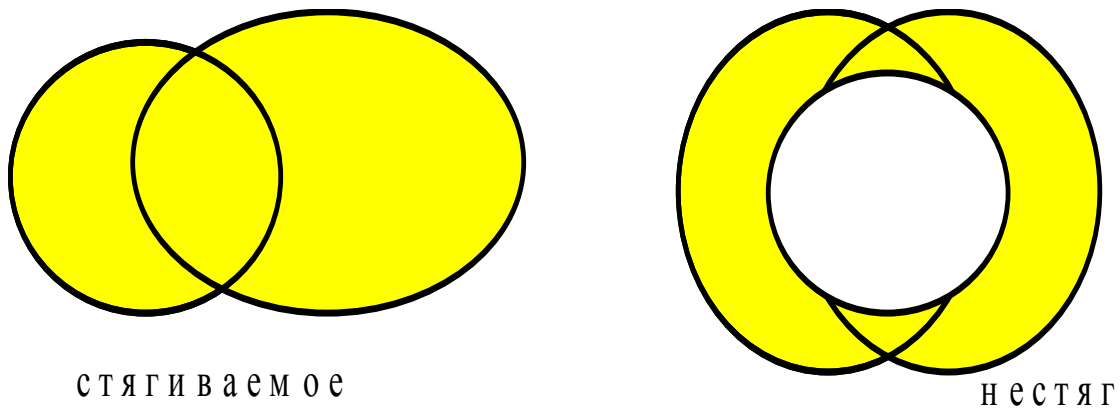

Рис. 177 Стягиваемое и нестягиваемое семейства множеств.

кубы с длиной ребра меньшей $\varepsilon$, $l_0 << \varepsilon$. Например, квадрат, изображенный на Рис. 178 справа можно разбить на 25 одинаковых квадратов с длиной стороны, в 5 раз меньшей. Тогда центральный квадрат в этом разбиении имеет шар, целиком принадлежащий каждому из двух множеств семейства. Под шаром здесь понимается пересечение с объединением всех множеств семейства. Сформулируем это в виде теоремы.

Теорема 170

*Если семейство $F=(A_1,A_2,...A_n)$ является регулярным (regular), то существует $\varepsilon>0$ такое, что для любой точки $p\in H=A_1\cup A_2\cup...\cup A_n$, для любого куба $K$ с центром в $p$ и с $l_0<\varepsilon$ и для шара этого куба $S(K)$ (самого куба и всех кубов, соседних с ним) существует множество $A_s$ такое, что $S(K)\cap H\subseteq A_s$.*

Доказательство.

Пусть семейство F регулярно при $l_0<\varepsilon$. Рассмотрим куб $K_1$ с длиной ребра

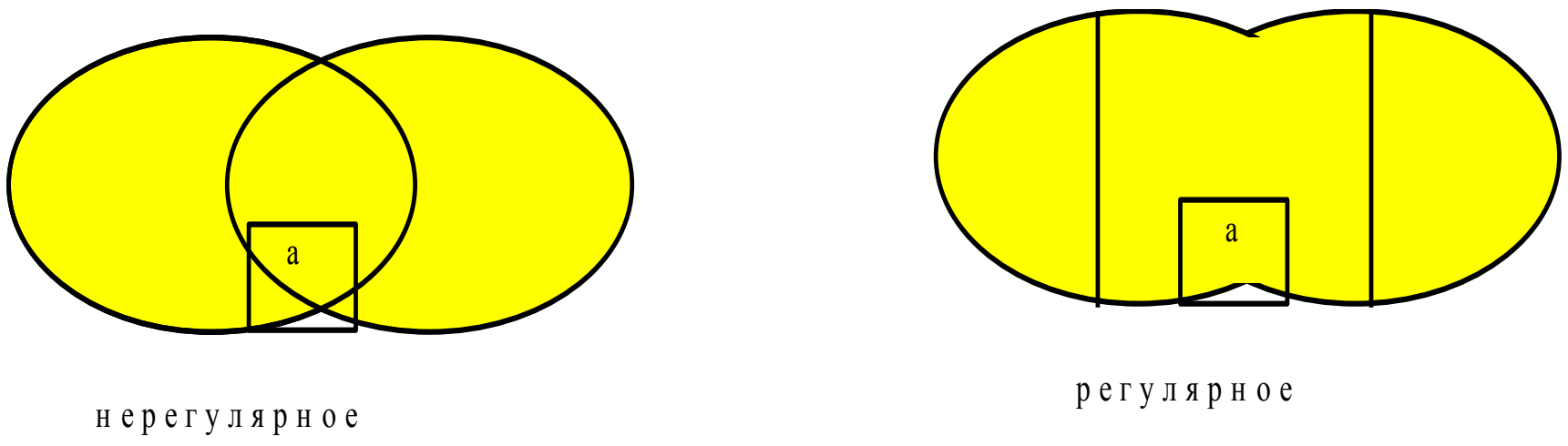

Рис. 178 Регулярное и нерегулярное семейства множеств.

$l_0/8$. Он находится в центре куба К с длиной ребра $l_0$. Следовательно все кубы, имеющие длину ребра $l_0/8$ и соседние $K_1$, будут находиться внутри куба К. Так как $K\cap H\subseteq A_s$, то $S(K)\cap H\subseteq A_s$. Что и требовалось доказать.□

## ПРАВИЛЬНОЕ ПОКРЫТИЕ ПОВЕРХНОСТИ В n-МЕРНОМ ЕВКЛИДОВОМ ПРОСТРАНСТВЕ. НЕРВ ПРАВИЛЬНОГО ПОКРЫТИЯ ГОМОТОПЕН БАЗОВОМУ МОЛЕКУЛЯРНОМУ ПРОСТРАНСТВУ ПОВЕРХНОСТИ

Пусть теперь мы имеем правильное семейство $F=(A_1,A_2,...A_n)$ и $H=A_1\cup A_2\cup...\cup A_n$, в евклидовом пространстве $E^n$.

Пусть $(\delta,\varepsilon)$ является характеристикой этого семейства. Выберем базовое кирпичное пространство для H с достаточно малым ребром $l_0<\min(\delta,\varepsilon/8)$ куба. Построим для $H=A_1\cup A_2\cup...\cup A_n$, нерв G(F) семейства F и базовое молекулярное пространство M(H) для H. Легко убедиться, что эти нервы гомеоморфны, так как базовые молекулярные пространства $M(A_1)$ $M(A_2),...M(A_n)$ являются подпространствами пространства M(H) и образуют его правильное покрытие, и G(F) является нервом этого покрытия. Сформулируем это рассуждение в форме теоремы.

Теорема 171

*Пусть X и $F=(A_1,A_2,...A_n)$ являются некоторым множеством и его правильным покрытием в $E^n$. Тогда существует $\varepsilon>0$ такой, что для всех $l_0<\varepsilon$ базовое молекулярное пространство M(X) множества X гомотопно нерву G(F) пересечений семейства F.*

Доказательство

Если F есть правильное покрытие, то при переходе к M(X) семейство множеств F является правильным покрытием молекулярного пространства M(X) и, следовательно, его нерв гомотопен M(X). Что и требовалось доказать.□

Согласно этой теореме, достаточно найти хотя бы одно правильное покрытие поверхности, или, в общем случае, множества в пространстве $E^n$. Его нерв будет гомотопен базовому молекулярному пространству данного

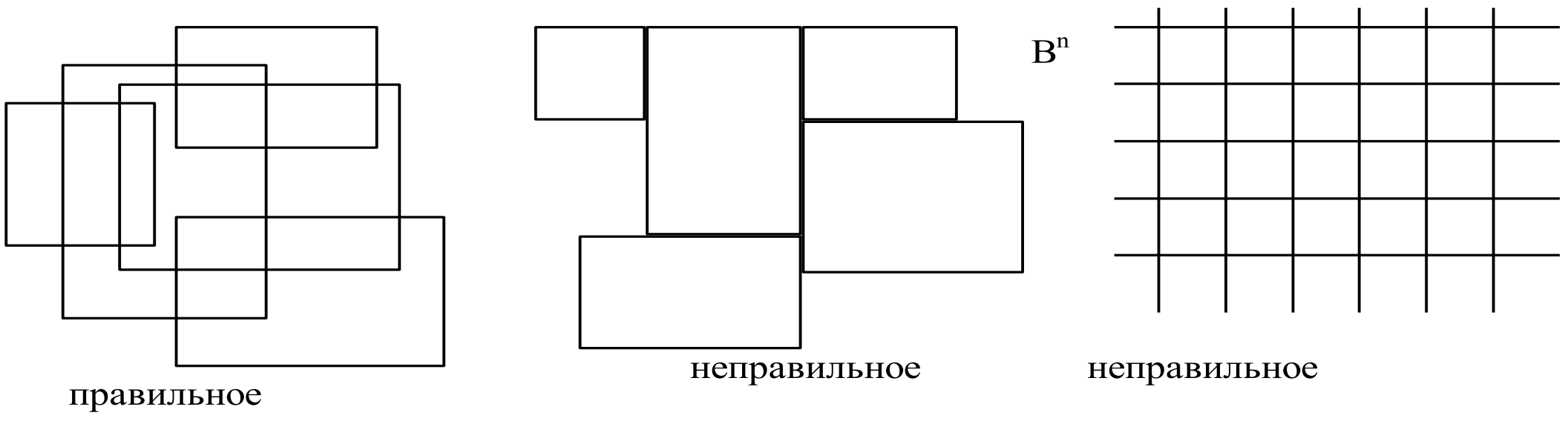

Рис. 179 Правильное и неправильное семейства множеств.

множества и, следовательно, нерву любого другого правильного покрытия. Любое правильное покрытие дает, таким образом, молекулярное пространство, являющееся дискретным образом непрерывного пространства. Однако, можно применять не только правильные покрытия. Если имеется некоторое покрытие, которое может быть преобразовано в правильное небольшими деформациями, то мы также можем его использовать при построении молекулярной модели.

Определение приводимого к правильному покрытия F.

Пусть X и $F=(A_1, A_2,...A_n)$ являются некоторым множеством и его покрытием в $E^n$. Пусть существует покрытие $H=(B_1,B_2,...B_n)$ для X такое, что:

1. H является правильным покрытием множества X,
2. Каждый элемент покрытия F содержится в соответствующем элементе покрытия H, $A_K \subseteq B_K$, k=1,2,...n,
3. Пересечение $A_K$ и $A_P$ непусто тогда и только тогда, когда пересечение $B_K$ и $B_P$ также непусто, $A_K \cap A_P \neq \varnothing \Leftrightarrow B_K \cap B_P \neq \varnothing$, к,р = 1,2,...n.

Тогда покрытие F называется приводимым к правильному.

Очевидно, что нервы пересечений покрытий F и H являются изоморфными. Поэтому мы можем применять приводимое покрытие наравне с правильным.

Следствие

Приводимое покрытие может быть использовано наравне с правильным.

Рассмотрим наиболее простые примеры приводимых и правильных покрытий.

1. Базовое кирпичное пространство.

Наиболее простым и естественным покрытием, близким к правильному, является базовое кирпичное пространство для $E^n$, то есть семейство кубов с длиной ребра $l_0$ и с координатами вершин, кратными $l_0$. Оно является стягиваемым и непрерывным семейством множеств, но не является регулярным. Чтобы сделать его правильным, достаточно "раздуть" каждый куб, увеличив длину каждого его ребра на величину $d < l_0/2$ в обоих направлениях (Рис. 180). Увеличение длины ребра делает базовое кирпичное пространство регулярным семейством, а ограничение означает, что нерв пересечений не меняется. С другой стороны, увеличение размера куба должно быть выбрано достаточно малым, чтобы не привести к увеличению числа кубов, пересекающих данную поверхность, не имеющую особенностей. Это позволяет рассматривать базовое кирпичное пространство как эквивалентное правильному.

2. Пространство из блоков.

Другим примером правильного покрытия в $E^n$ является покрытие n-мерными прямоугольными параллелепипедами с ребрами, параллельными координатным осям. При этом пересечение любых двух параллелепипедов должно быть либо пустым, либо также невырожденным n-мерными

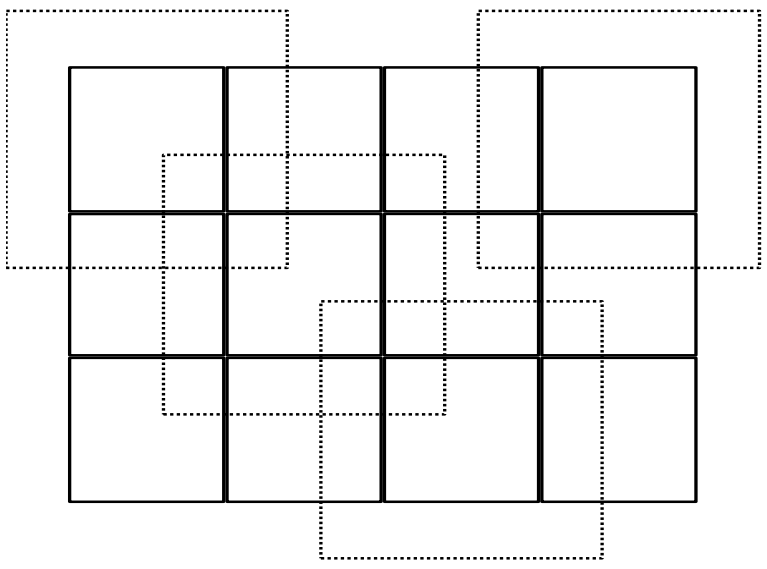

Рис. 180 Раздувание кирпичного пространства и превращение его в правильное покрытие.

прямоугольным параллелепипедом. Однако, в любом случае такое покрытие будет приводимым.

Список литературы к главе 14.

# КОМПЬЮТЕРНЫЕ ЭКСПЕРИМЕНТЫ. МОДЕЛИРОВАНИЕ И ДВИЖЕНИЕ МНОГОМЕРНЫХ ОБЪЕКТОВ В КОМПЬЮТЕРАХ

Here we sketch how computer calculations with finite element method and further computer experiments leeding to the notion of a molecular space and contractible transformations. As an application to computer graphic we describe on a concrete example an algorithm of obtaining molecular models of continuous n-dimensional objects, altering their forms and shapes and moving them in n-dimensional space.

Содержание этой главы частично можно найти в работах [10,14,16,17,18,24,26,43].

## *КОМПЬЮТЕРНЫЕ ЭКСПЕРИМЕНТЫ И ТЕОРИЯ МОЛЕКУЛЯРНЫХ ПРОСТРАНСТВ*

Представляет интерес рассказать, что помогло появлению концепции молекулярных пространств., и описать и компьютерные опыты, которые позволили подвести некоторый экспериментальный фундамент под данную теорию.

В основе компьютерных опытов уже, естественно, лежали некоторые гипотезы, которые подкрепляли убеждение в том, что непрерывный континуум может быть заменен на дискретный, и результаты расчетов не должны сильно зависеть от выбора континуума. По сути дела, дискретный континуум уже используется при работе компьютера, объем памяти которого всегда конечен. Дискретный континуум заменяет непрерывный в

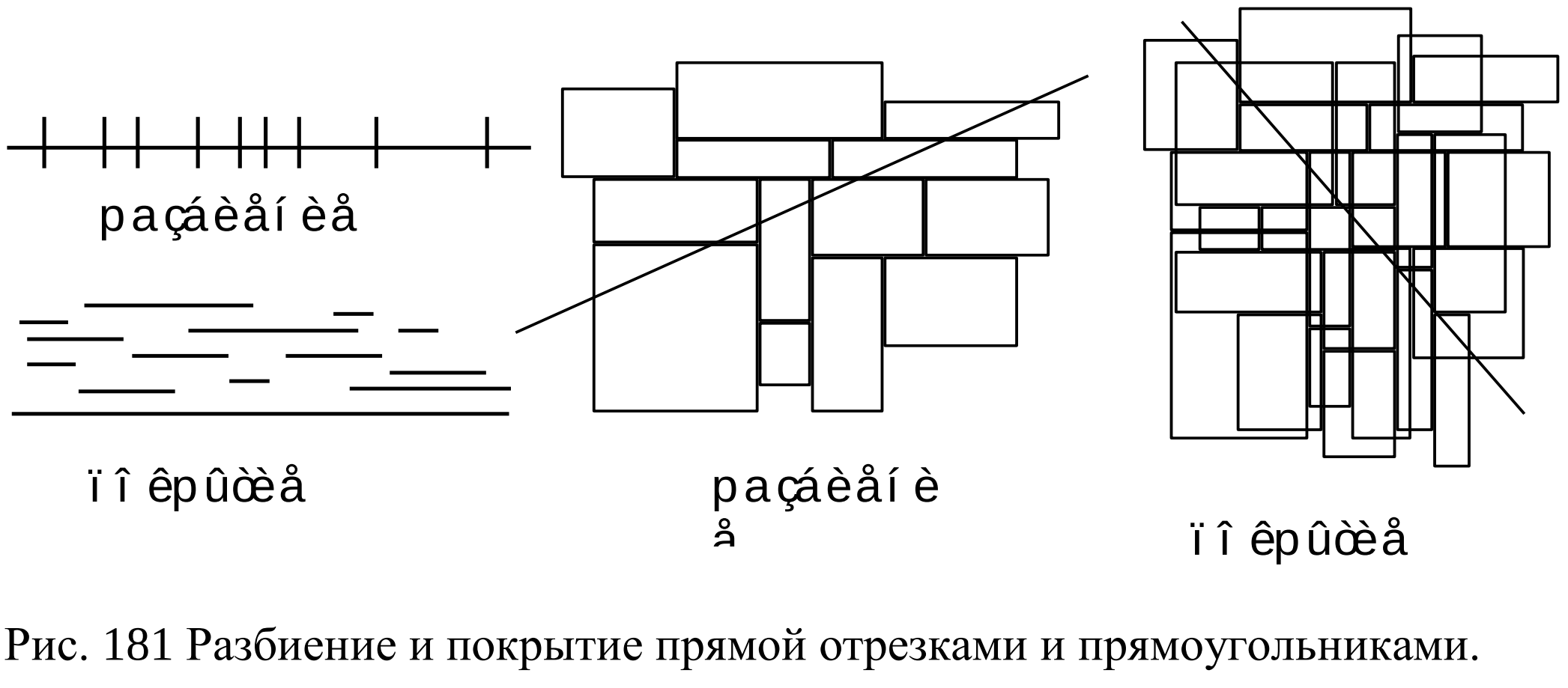

Рис. 181 Разбиение и покрытие прямой отрезками и прямоугольниками.

методе сеток для численного решения дифференциальных уравнений. Сетка есть не что иное, как конечное множество точек, на котором задана и решается система алгебраических уравнений. В математике и смежных науках давно известен метод конечных элементов, который позволяет заменить трудоемкие вычисления по всей заданной области разбиением области на небольшие участки, простыми вычислениями по каждому из участков и последующей сравнительно простой обработкой полученных результатов.

Одним из важных этапов этого процесса является выбор размеров и формы конечных элементов области. Как правило, этот этап не является простым, требует привлечения дополнительных критериев, и для каждого отдельного случая необходимо выбирать свои параметры разбиения. При этом возникает вопрос, как влияют форма, размеры и соседство элементов на результаты вычислений. И, как связанное с этим вопросом, но на более общем уровне, имеется ли связь различных типов разбиений данной поверхности друг с другом. Подобные примеры можно продолжить. Для ответов на подобные вопросы, в значительной степени связанные с однотипными расчетами, мы использовали компьютерные вычисления и соответствующим образом составленные программы.

Переходим теперь к описанию компьютерных экспериментов. Следует заметить, что линейные размеры исходного объекта всегда были значительно больше линейных размеров элементов. Была проделана серия экспериментов по однотипной схеме.

1. Задавалось какое-либо множество, одномерное, двумерное или трехмерное.

2. Выбиралось разбиение или покрытие этого множества некоторыми элементами.

3. Строился нерв (граф пересечений) этого покрытия или разбиения.

4. Полученный нерв анализировался различными методами.

Напоминаем, что нервом (графом пересечений) множества $S=(S_1,S_2,...S_n)$ элементов $S_1,S_2,...S_n$ называется множество $V=(v_1,v_2,...v_n)$ точек, в котором для каждой точки $v_k$ ставится в соответствие элемент $S_k$, и две точки $v_k$ и $v_p$ соединены связью $(v_k v_p)$ тогда и только тогда, когда $S_k$ и $S_p$ имеют непустое пересечение, $S_k \cap S_p \neq 0$.

<u>А. Эксперименты с отрезком прямой</u> (Рис. 181).

1. Выбирался отрезок прямой, разбивался на конечное число небольших отрезков-элементов и строился нерв множества элементов.

2. Выбирался отрезок прямой, покрывался произвольно конечным числом небольших отрезков-элементов и строился нерв множества элементов.

3. Плоскость разбивалась на клетки в виде одинаковых квадратов со сторонами, параллельными координатным осям-сеть квадратов; проводился отрезок прямой в произвольном направлении, и рассматривалось множество квадратов-элементов, пересекающих отрезок хотя бы в одной точке. После этого строился нерв множества таких квадратов.

4. Пространство разбивалось на кубы в виде одинаковых кубов со сторонами, параллельными координатным осям-сеть кубов, проводился отрезок прямой в произвольном направлении, и рассматривалось множество кубов-элементов, пересекающих отрезок хотя бы в одной точке. После этого строился нерв множества элементов.

В. Эксперименты с окружностью.

Эксперименты 1, 2, 3 и 4 проделывались для окружности.

С. Эксперименты с квадратом (Рис. 182).

1. Квадрат разбивался на клетки в виде небольших квадратов-элементов, и строился нерв множества элементов.
2. Квадрат разбивался на клетки в виде неодинаковых небольших прямоугольников-элементов со сторонами, параллельными сторонам квадрата, и строился нерв множества элементов.
3. Квадрат покрывался клетками в виде неодинаковых небольших прямоугольников-элементов со сторонами, параллельными сторонам квадрата, и строился нерв множества элементов.
4. Квадрат покрывался клетками в виде неодинаковых небольших кругов-элементов, и строился нерв множества элементов.
5. Пространство разбивалось на кубы в виде одинаковых кубов со сторонами, параллельными координатным осям-сеть кубов, помещался

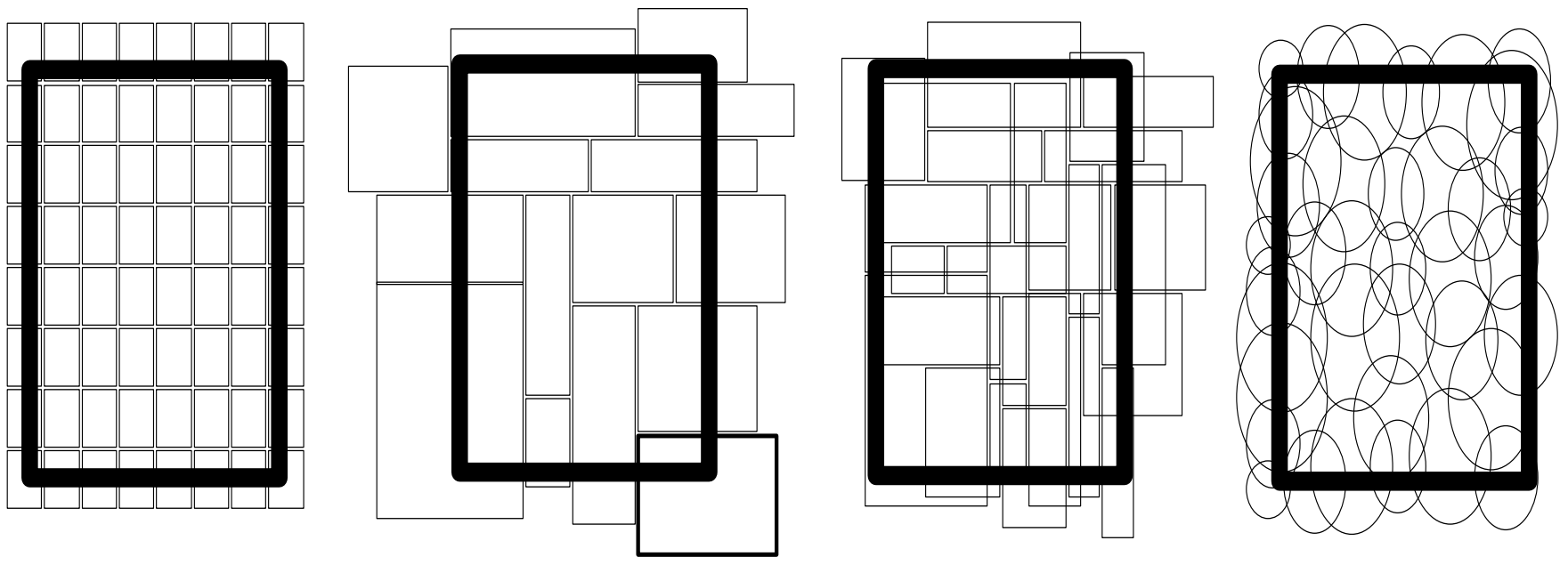

Рис. 182 Разбиения и покрытия квадрата для построения нерва множества.

квадрат произвольной ориентации, и рассматривалось множество кубов-элементов, пересекающих квадрат хотя бы в одной точке. После этого строился нерв множества элементов.

Е. Эксперименты со сферой.

1. Сфера разбивалась на клетки в виде небольших неодинаковых выпуклых многоугольников различной геометрии-элементов, и строился нерв множества элементов.
2. Сфера покрывалась клетками в виде неодинаковых небольших кругов-элементов, и строился нерв множества элементов.
3. Пространство разбивалось на кубы в виде одинаковых кубов со сторонами, параллельными координатным осям-сеть кубов, и

рассматривалось множество кубов-элементов, пересекающих сферу хотя бы в одной точке. После этого строился нерв множества элементов.

4. Пространство разбивалось на небольшие сферы в виде одинаковых сфер и рассматривалось множество этих сфер-элементов, пересекающих исходную сферу хотя бы в одной точке. После этого строился нерв множества элементов.

F. Эксперименты с поверхностью куба.

Эксперименты, аналогичные предыдущим, проделывались на поверхности куба. Была создана компьютерная программа для проведения экспериментов на достаточно сложных поверхностях с автоматическим выбором разбиений или покрытий.

В результате, было получено большое количество нервов, многие из которых выглядели достаточно сложно при их распечатке из памяти компьютера. Однако, некоторые особенности обнаружились довольно быстро.

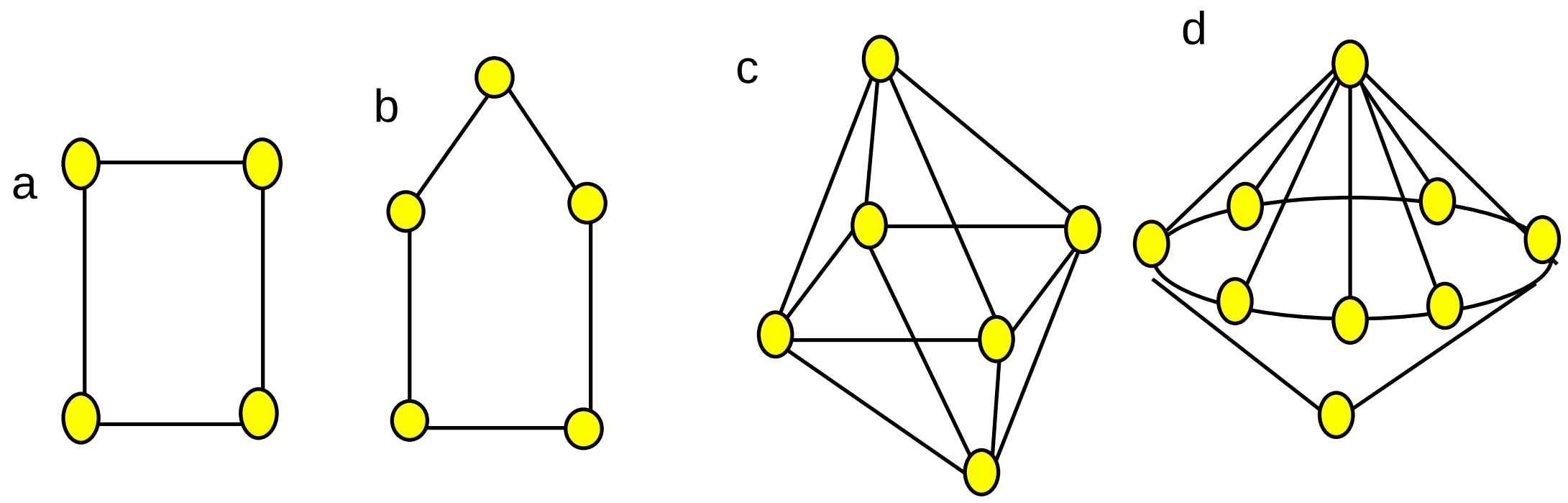

Рис. 183 Нервы покрытий отрезка прямой не содержат подпространств вида a и b. Нервы покрытий плоскости не содержат подпространств вида c и d.

O1. Никакой нерв, построенный для прямолинейного отрезка, криволинейного отрезка или окружности, то есть для одномерной поверхности. не содержит поднервы изображенные на Рис. 183a и b. Как выяснилось, этот результат был доказан несколько лет назад для интервальных графов [20]. Позднее, в работе [10] было доказано, что нерв достаточно длинного отрезка прямой семейством единичных кубов с целочисленными координатами в n-мерном евклидовом пространстве является одномерным молекулярным пространством.

O2. Никакой нерв пересечений, построенный для квадрата или сферы, не содержит поднервы, сходные изображенным на Рис. 183c и d. Этот результат доказывается несложно, но так как мы не встречали в литературе упоминаний о нем, мы сформулируем его в виде теоремы, доказательство которой предоставим читателю.

Теорема.

*Пусть дано произвольное покрытие (разбиение) плоскости выпуклыми замкнутыми множествами, для которых минимальный d и максимальный D диаметры лежат в некотором фиксированном положительном интервале. Тогда нерв этого покрытия не содержит поднервов, являющихся замкнутыми нормальными двумерными молекулярными пространствами.*

Обнаружение еще одной особенности заняло значительно больше времени и дополнительных экспериментов для ее проверки. Суть ее заключалась в следующем. Мы нашли четыре специальных типа преобразований нервов: приклеивание и отбрасывание точек и приклеивание и отбрасывание связей, которые мы назвали точечными. Оказалось, что:

A1. Все нервы отрезка прямой или отрезка кривой могут быть преобразованы один в другой при помощи точечных преобразований.

A2. Все нервы окружности могут быть преобразованы один в другой при помощи таких преобразований.

A3. Ни один нерв отрезка прямой или кривой не может быть преобразован в нерв окружности при помощи точечных преобразований и наоборот.

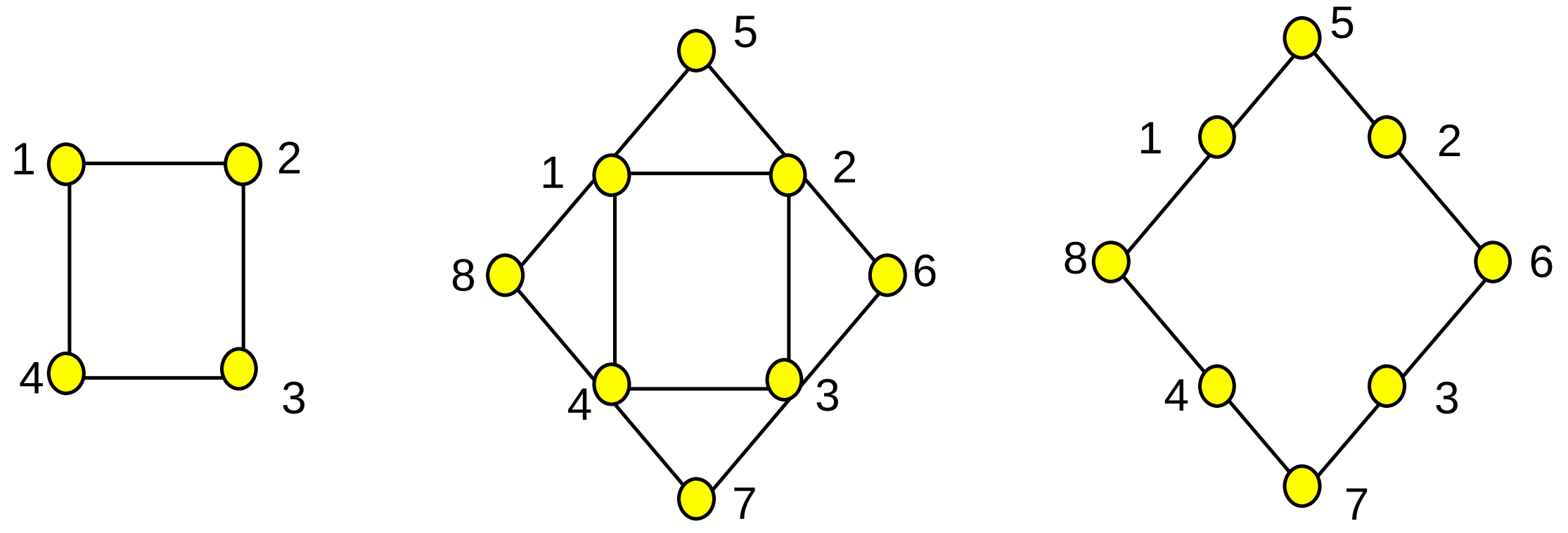

Рис. 184 Увеличение числа точек молекулярной окружности при помощи точечных приклеиваний точек.

B1. Все нервы квадрата могут быть преобразованы один в другой при помощи таких преобразований.

B2. Все нервы сферы могут быть преобразованы один в другой при помощи таких преобразований.

B3. Ни один нерв квадрата не может быть преобразован в нерв сферы при помощи точечных преобразований и наоборот.

Далее мы выяснили другую интересную особенность.

C1. Все нервы отрезка прямой, кривой и квадрата могут быть преобразованы один в другой и в тривиальный нерв, состоящий из одной точки точечными преобразованиями.

C2. Нервы окружности могут быть преобразованы точечными преобразованиями в минимальный нерв $C_4$, изображенный на Рис. 183а.

C3. Нервы сферы могут быть преобразованы точечными преобразованиями в минимальный нерв, изображенный на Рис. 183с.

C4. Нервы окружности, сферы и тривиальный не могут быть преобразованы точечными преобразованиями один в другой.

Мы, естественно, опускаем мелкие детали проведенных экспериментов.

После этого потребовалось еще какое-то время чтобы понять следующее:

1. Мы получили не что иное, как дискретную модель, которая является нервом покрытия или разбиения непрерывного пространства.

2. Отрезок прямой, кривой и квадрат, в топологии известны как пространства, гомотопные одно другому и гомотопные точке. С другой стороны они не гомотопны окружности и сфере, которые, в свою очередь, не гомотопны друг другу. Таким образом, точечные преобразования являются комбинаторными образами отображений, сохраняющих гомотопию. Иными словами, два молекулярных пространства гомотопны, если одно из них можно стянуть в другое точечными преобразованиями.

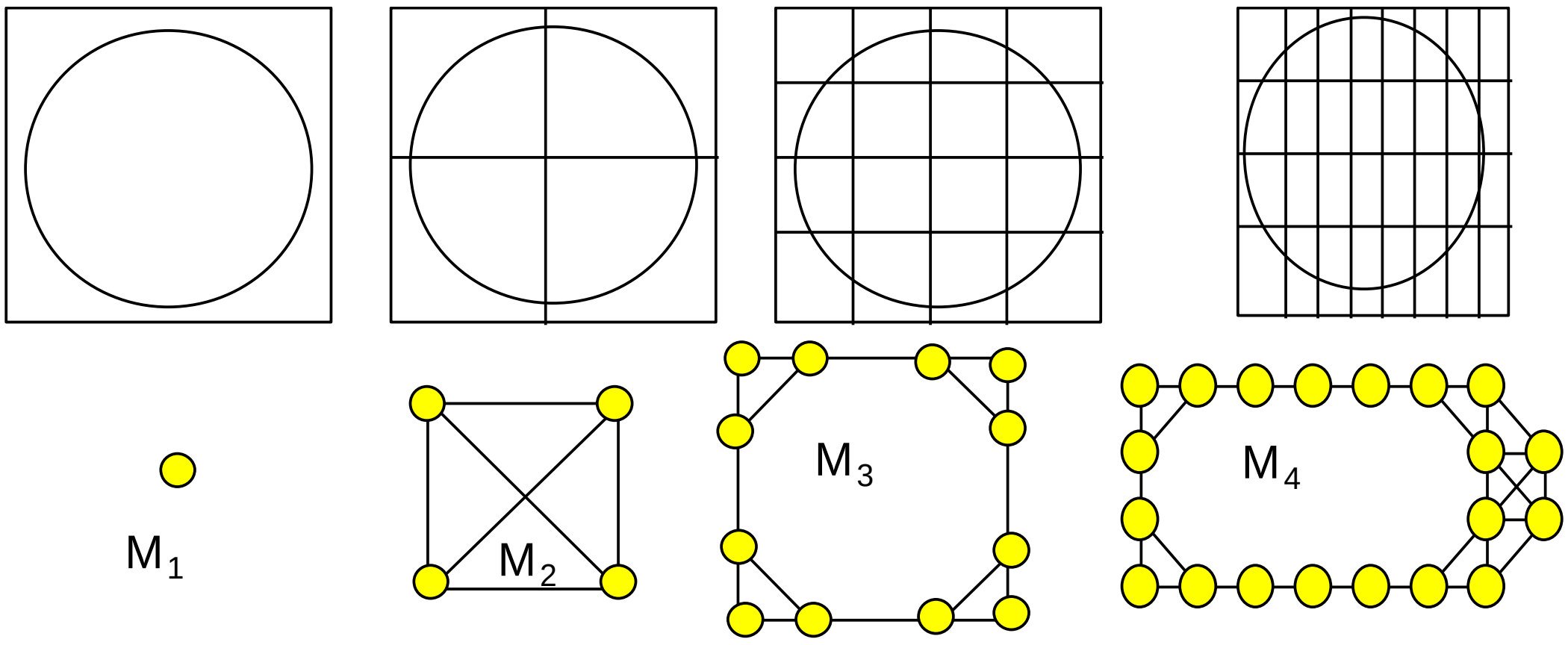

Рис. 185 Процесс нахождения молекулярного пространства для окружности.

Следующий шаг был проверить все еще раз на других поверхностях. Эксперименты на торе сплошном кубе и сплошном шаре и тому подобное подтвердили все результаты.

В заключение заметим, что впервые упоминание о дискретном представлении непрерывного множества согласно [29] встречается у средневековых арабских философов мутакаллимов, с точки зрения которых для образования квадрата (или границы квадрата, то есть окружности) требуется четыре точки.

## *АЛГОРИТМ ПОЛУЧЕНИЯ МОЛЕКУЛЯРНОГО ПРОСТРАНСТВА ДЛЯ НЕПРЕРЫВНОЙ ПОВЕРХНОСТИ*

Пусть поверхность S находится в области ($0 < x_i < d$), i=1,2,...n, пространства $E^n$.

1. Выбираем куб K с длиной ребра d, имеющий точку в начале координат.
2. Разбиваем его на кубы с длиной ребра d/2 и строим базовое молекулярное пространство $M_1$ для S.
3. Разбиваем K на кубы с длиной ребра d/4 и строим базовое молекулярное пространство $M_2$ для S.
4. Проверяем гомотопность $M_1$ и $M_2$.
5. Повторяем этапы 2, 3 и 4, каждый раз уменьшая длину ребер кубов в два раза. Мы получаем последовательность молекулярных пространств $M_1$, $M_2$,....$M_n$...
6. Процесс прекращается, когда, начиная с некоторого номера p, все молекулярные пространства $M_n$ с номерами n>p гомотопны одно другому. В качестве модели поверхности S можно выбрать любое из них.

На Рис. 185 показан процесс нахождения молекулярного пространства для окружности при помощи описанного алгоритма.

Молекулярное пространство $M_1$, состоящее из одной точки, гомотопно $M_2$, но не гомотопно $M_3$. Пространства $M_3$ и $M_4$ гомотопны одно другому и всем последующим, полученным с использованием более мелких разбиений. Легко убедиться, что они являются одномерными замкнутыми пространствами, хотя и не обладают свойством нормальности. Естественно, что они стягиваются точечными преобразованиями к наименьшей нормальной окружности $S^1_{min}$, состоящей из четырех точек.

## *НЕКОТОРЫЕ АЛГОРИТМЫ ДВИЖЕНИЙ МНОГОМЕРНЫХ МОЛЕКУЛЯРНЫХ ОБЪЕКТОВ*

В этом разделе мы рассмотрим алгоритмы движений и преобразований n-мерных объектов, таких как линии, поверхности, тела и тому подобное. В качестве конкретных примеров возьмем одномерные пространства-окружности, преобразования которых легко иллюстрировать рисунками. Обобщения на двумерный, трехмерный и многомерный случаи не представляют трудностей и сводятся лишь к увеличению числа операций.

1. Раздувание окружности.

Пусть имеется минимальная окружность, образованная четырьмя точками 1, 2, 3 и 4, как это показано на Рис. 184.

Для того, чтобы увеличить число точек до восьми, приклеим точку 5 к точкам 1 и 2, точку 6 к точкам 2 и 3, точку 7 к точкам 3 и 4, и точку 8 к точкам 4 и 1. Все эти приклеивания являются разрешенными точечными преобразованиями молекулярного пространства, не меняющими его характеристик. Однако, при этом возникают связи между точками, которые мы можем отбросить.

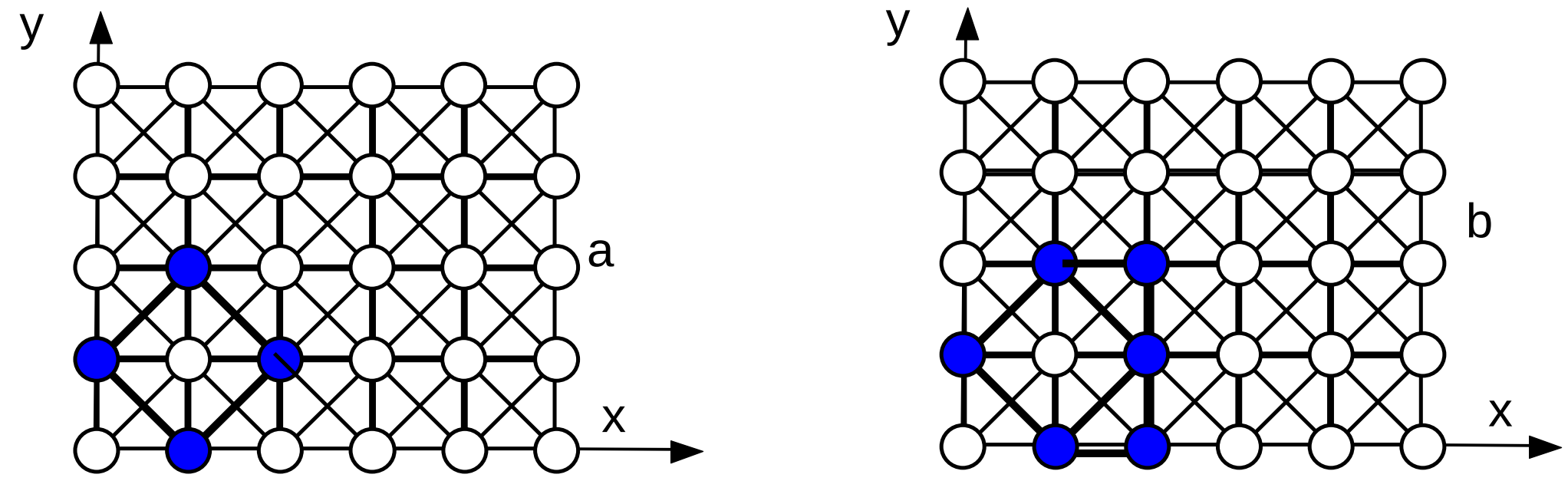

Рис. 186 Сдвиг окружности вправо в евклидовом молекулярном пространстве.

На следующем шаге отбрасываем связи (1,2), (2,3), (3,4) и (4,5). Полученное пространство является одномерной окружностью, состоящей из восьми точек. В этом примере мы не затрагивали перемещения объекта, как целого, в окружающем пространстве.

2. Движение окружности. Для того, чтобы определить такие движения, нам необходимо предварительно ввести систему координат, относительно которой такие движения происходят. Система координат определяется молекулярным евклидовым пространством. Ранее мы подробно описывали способы получения молекулярных моделей евклидовых пространств. Это-отдельная задача, решение которой зависит от конкретных условий в каждом отдельном случае.

В данном случае выберем молекулярную модель двумерной плоскости, в

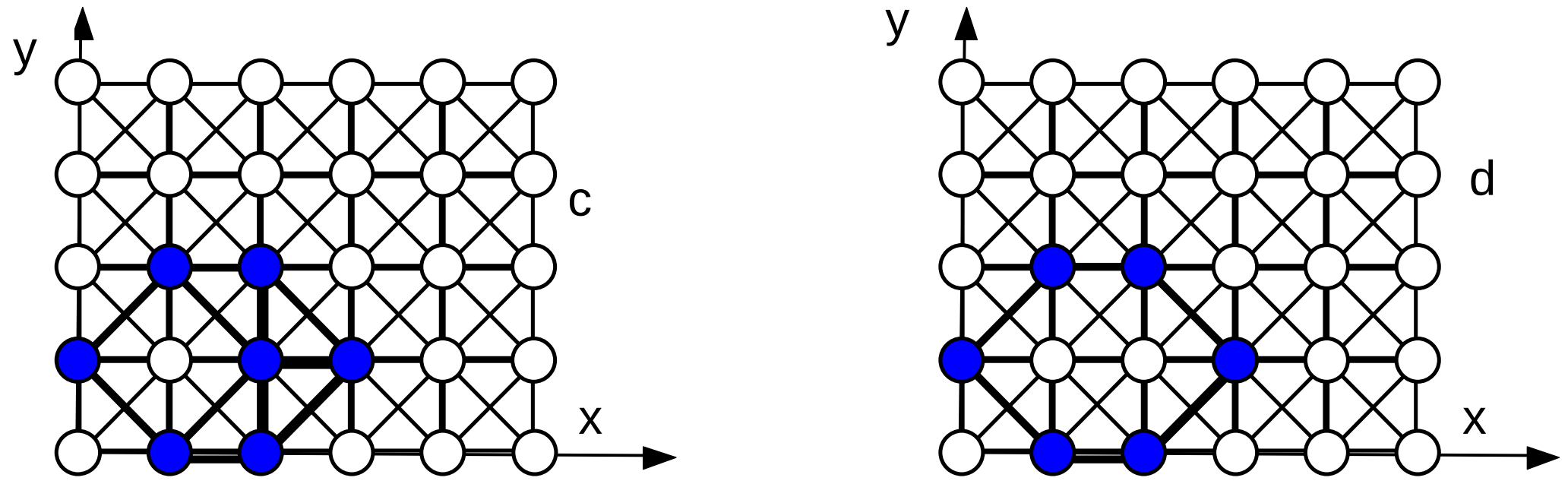

Рис. 187 Сдвиг окружности вправо в евклидовом молекулярном пространстве.

которой точки непрерывного евклидова пространства с целочисленными координатами являются точками молекулярного пространства, и каждая точка соединена связью с ближайшим окружением, как это показано на Рис. 186. Данная модель двумерной плоскости является полным однородным двумерным молекулярным пространством. Окружность состоит из четырех точек темного цвета, соединенных ребрами на Рис. 186.

Для движения вправо приклеим к окружности две точечные точки. Результат приклеивания изображен на Рис. 186, Рис. 187, Рис. 188. Далее к пространству приклеиваем еще одну точку. После этого отбрасываем точку, инцидентную пяти ребрам, и получаем окружность, состоящую из шести точек. Наша задача теперь уменьшить число точек до четырех, сдвинув левую часть окружности вправо. Для этого приклеиваем к окружности еще одну точку, увеличив число точек до семи. После этого мы можем отбросить одну точку, перейдя к окружности из шести точек, а затем еще пару точек, получив, окончательно, окружность из четырех точек, показанную на Рис. 188, и являющуюся первоначальной окружностью, сдвинутой вправо на один шаг. Еще раз подчеркнем, что все

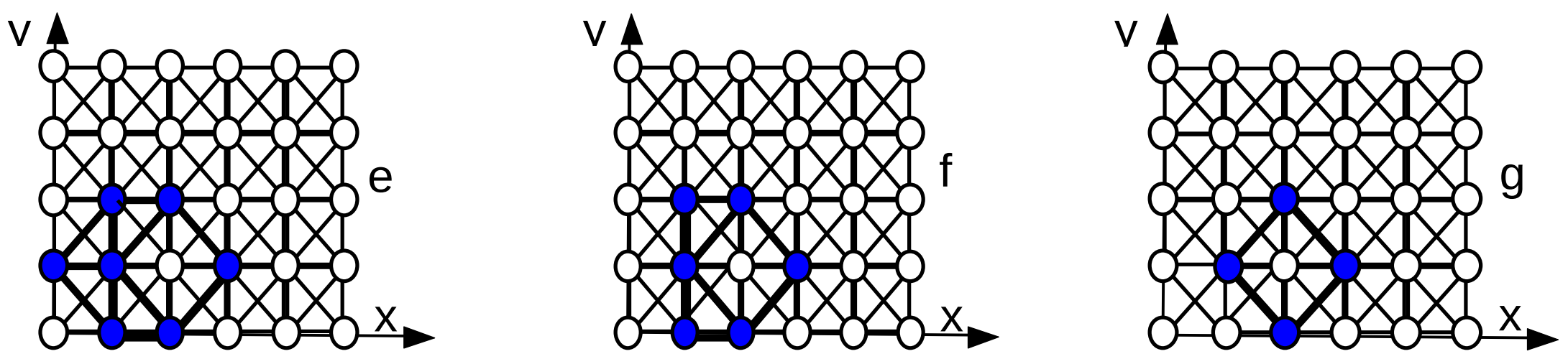

Рис. 188 Сдвиг молекулярной окружности вправо в евклидовом молекулярном пространстве.

преобразования являются точечными, не меняющими характеристик смещаемого пространства. Точно таким же образом мы можем двигать окружность в других направлениях, влево, вверх или вниз, или же изменять ее размеры. Преобразования двумерных трехмерных и вообще произвольных многомерных объектов принципиально ничем не отличаются от преобразований окружности, поскольку необходимое условие, налагаемое на приклеивания и отбрасывания точек и связей единственное-они должны быть точечными. Это условие, как известно, не зависит от размерности, а определяется ближайшим окружением точки или ребра. Преимущество точечных преобразований также и в том, что они локальны, например, чтобы отбросить точку нам нужно знать структуру только ее ближайшего окружения-ее окаем, и при этом нас не интересуют остальные точки пространства.

Список литературы к главе 15.

# ЧАСТЬ 3. ПРИМЕНЕНИЯ ТЕОРИИ МОЛЕКУЛЯРНЫХ ПРОСТРАНСТВ

## ДИФФЕРЕНЦИАЛЬНЫЕ УРАВНЕНИЯ НА МОЛЕКУЛЯРНЫХ ПРОСТРАНСТВАХ

Using well known approximate methods in differential equation theory a dynamic system is defined on a molecular space. This system comprises partial differential equations such as parabolic, elliptic and hyperbolic. Some natural background of definition of differential equations on a molecular space is discussed. Fundamental structural difference between parabolic and hyperbolic equations is found. Supposition is made that differential equations on molecular spaces can be applied to such objects as neuron and electric networks, chemical and biological structures and any objects that can be described by molecular spaces or graphs. As an example numerical solutions of parabolic elliptic and hyperbolic equations on one, two and more dimensional spheres, tory, projective planes, Euclidean spaces and trees are given. A regular equation which is a generalization of given type of differential equations is offered and some properties of the regular equation of the type are studied.

Большинство из полученных результатов содержится в работах [5,6,9,43,48].

### *СХОДСТВО И РАЗЛИЧИЕ ДИФФЕРЕНЦИАЛЬНЫХ УРАВНЕНИЙ НА НЕПРЕРЫВНЫХ И ДИСКРЕТНЫХ ПРОСТРАНСТВАХ*

Прежде всего, мы должны определить дифференциальное уравнение на молекулярном пространстве. Мы сделаем это, пользуясь аналогией с давно уже используемым для численного решения дифференциальных уравнений методом сеток. По сути дела, решение дифференциального уравнения с выбором определенной сетки уже есть его решение на некотором молекулярном пространстве, играющим роль сетки. Отличие состоит в том, что при применении метода сеток мы выбираем сетку, приспособленную к решению конкретной задачи на непрерывном пространстве. В случае молекулярного пространства сетка задается уже самим видом пространства. Однако, это, на первый взгляд незначительное и несущественное отличие, в конце концов привело к тому, что были найдены фундаментальные особенности разных типов дифференциальных уравнений, не зависящие ни от выбора сетки, ни от размерности, ни других математических черт пространства. Эти особенности отражают сущность природных процессов, описываемых уравнениями. Покажем это на простом примере. Параболическое уравнение иногда называют тепловым уравнением или уравнением распространения тепла. Это уравнение показывает, как распространяется тепло в некоторой выделенной области. Помимо распространения тепла, это уравнение описывает множество самых разнообразных физических и природных процессов. Передача тепла происходит локально от точки к точке и никак не связана в своей основе с математическими характеристиками пространства. Следовательно, параболическое уравнение и его решение должны иметь такую форму, которая не содержит параметров, зависящих от математических характеристик пространства, например, от размерности. Тем не менее легко убедиться, что в непрерывном случае как форма записи параболического уравнения, так и вид его решения, напрямую зависят от размерности пространства. Применение уравнений математической физики к молекулярным пространствам как раз и позволяет представить эти уравнения в наиболее общей форме, не содержащей конкретных математических характеристик пространства и определить фундаментальное сходство и отличие параболического уравнения от, например, гиперболического. Об этом мы будем говорить позднее.

С другой стороны, применение аппарата дифференциальных уравнений на дискретных пространствах может дать новые результаты, возможно, не имеющие аналога в непрерывном случае. Особенно это должно касаться пространств с малым объемом (числом точек). Например, распространение тепла (параболическое уравнение), или волн (гиперболическое уравнение) по непрерывной двумерной сфере должно быть аналогично соответствующим процессам на молекулярной модели с достаточно большим числом точек. Однако, в теории молекулярных пространств имеется наименьшая двумерная сфера, число точек которой равно шести, и на которой мы также можем рассматривать процессы передачи тепла и волновые. При этом, в связи с небольшим количеством элементов дискретного пространства, здесь могут возникнуть особенности, не имеющие непрерывного аналога. Другое серьезное отличие связано с существованием естественной наименьшей длины в молекулярном пространстве, определяемой длиной связи, соединяющей две соседние точки пространства. В применении к волновым процессам это означает отсутствие бесконечно коротких волн и бесконечно высоких частот, то есть отсутствие факторов, часто ведущих к появлению расходимостей. Еще одна особенность использования молекулярных пространств состоит в том, что мы можем получать решения дифференциальных уравнений на пространствах, которые не имеют непрерывных аналогов. Например, рассмотрев эллиптическое уравнение на однородном дереве, которое является одномерным регулярным пространством, мы покажем, что вид решения зависит от объема окаема каждой точки. Мы покажем, что можно подобрать такое дерево (то есть молекулярное пространство, не

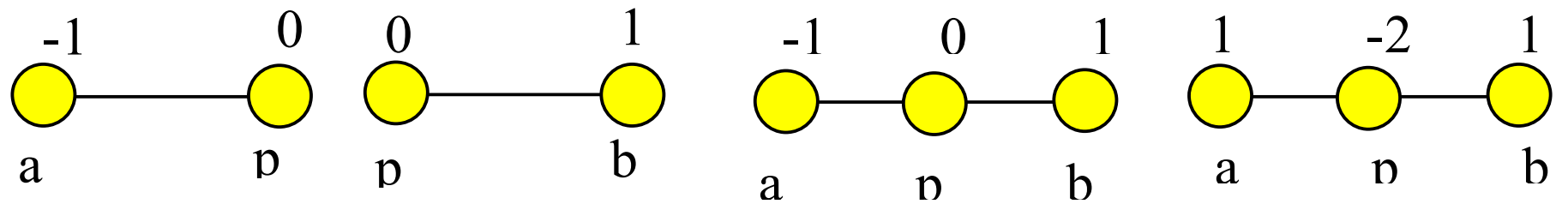

Рис. 189 Шаблоны конечно-разностных схем.

содержащее циклов, или же содержащее очень длинные циклы), на котором решение будет иметь вид, сходный с решением трехмерного эллиптического уравнения. Сходство ставит интересную проблему: отражают ли решения дифференциальных уравнений топологию пространства? Наиболее очевидный ответ-не отражают, однако это-задача более углубленных исследований.

## *ОПРЕДЕЛЕНИЕ УРАВНЕНИЙ МАТЕМАТИЧЕСКОЙ ФИЗИКИ НА МОЛЕКУЛЯРНОМ ПРОСТРАНСТВЕ*

Рассмотрим стандартную аппроксимацию дифференциальных операторов первого и второго порядков разностными операторами соответствующих порядков так, что любая производная будет аппроксимирована соответствующим разностным отношением. Диаграммы таких разностных операторов приведены в [53]. Из диаграмм видно, что первая и вторая частные производные по координате некоторой функции f в точке p в момент t могут быть аппроксимированы следующим образом:

$$\frac{\partial f}{\partial x} = \frac{f_p^t - f_a^t}{h}\ ;\ \frac{f_b^t - f_p^t}{h}\ ;\ \frac{f_b^t - f_a^t}{2h}\ ;$$

$$\frac{\partial^2 f}{\partial x^2} = \frac{f_b^t - 2f_p^t + f_a^t}{h^2}\ ,$$

Здесь h = Δx. Аналогично представим первую и вторую производные по времени в точке p в момент t в соответствии с диаграммами.

$$\frac{\partial f}{\partial t} = \frac{f_p^{t+1} - f_p^t}{\Delta t}; \qquad \frac{\partial^2 f}{\partial t^2} = \frac{f_p^{t+1} - 2f_p^t + f_p^{t-1}}{\Delta t^2}.$$

Рассмотрим дифференциальные уравнения

$$a\frac{\partial f}{\partial t} = b\frac{\partial f}{\partial x}; \quad a\frac{\partial f}{\partial t} = b\frac{\partial^2 f}{\partial x^2}; \quad a\frac{\partial^2 f}{\partial t^2} = b\frac{\partial^2 f}{\partial x^2},$$

и заменим их на конечно-разностные. После соответствующих преобразований каждое из них может быть приведено к виду:

$$f_p^{t+1} = \sum_k c_{pk} f_k^t + \sum_k d_{pk} f_k^{t-1}, \ \ k = a, p, b. \quad (1)$$

где $c_{pk}$ и $d_{pk}$ являются некоторыми коэффициентами, зависящими от типа уравнения.

Из (1) видно, что значение функции в точке p в момент (t+1) определяется значениями функции в соседних с p точках в моменты t и (t-1) (Рис. 190). Величины коэффициентов $c_{pk}$ и $d_{pk}$ определяют, какой вид имеет уравнение: первого или второго оно порядка, гиперболическое, параболическое или эллиптическое. Для параболического уравнения $c_{pa}=e$, $c_{pp}=(1-2e)$, $c_{pb}=e$, $d_{pk}=0$, k=a, p, b,

$$f_p^{t+1} = ef_a^t + ef_b^t + (1-2e)f_p^t.$$

В дальнейшем будем считать, что e достаточно мало, и, следовательно, для параболического уравнения все коэффициенты неотрицательны. Для эллиптического уравнения конечно-разностное уравнение может быть представлено в схожем виде

$$f_p^t = ef_a^t + ef_b^t + (1-2e)f_p^t.$$

Для гиперболического уравнения $c_{pa}=e$, $c_{pp}=(2-2e)$, $c_{pb}=e$, $d_{pk}=0$, k=a, b, $d_{pp}=-1$,

$$f_p^{t+1} = ef_a^t + ef_b^t + (2-2e)f_p^t - f_p^{t-1} = ef_a^t + ef_b^t + (1-2e)f_p^t + f_p^t - f_p^{t-1}.$$

Проанализируем это уравнения с точки зрения его использования в молекулярных пространствах. Первые три слагаемых во всех случаях одинаковы, все коэффициенты перед значениями функций положительны и их сумма равна 1. Гиперболическое конечно-разностное уравнение имеет два дополнительных слагаемых в правой части по сравнению с параболическим и эллиптическим уравнениями. При переходе к двумерному, трехмерному и многомерному случаям легко убедиться непосредственно, что структура всех уравнений не изменится. Каждое из конечно-разностных уравнений может быть представлена в аналогичной форме. Таким образом представим параболическое уравнение в форме

$$f_p^{t+1} = \sum_k c_{pk} f_k^t, \quad \forall c_{pk} \ge 0, \ \sum_p c_{pk} = 1,$$

где в первой сумме суммирование ведется по всем точкам шаблона,

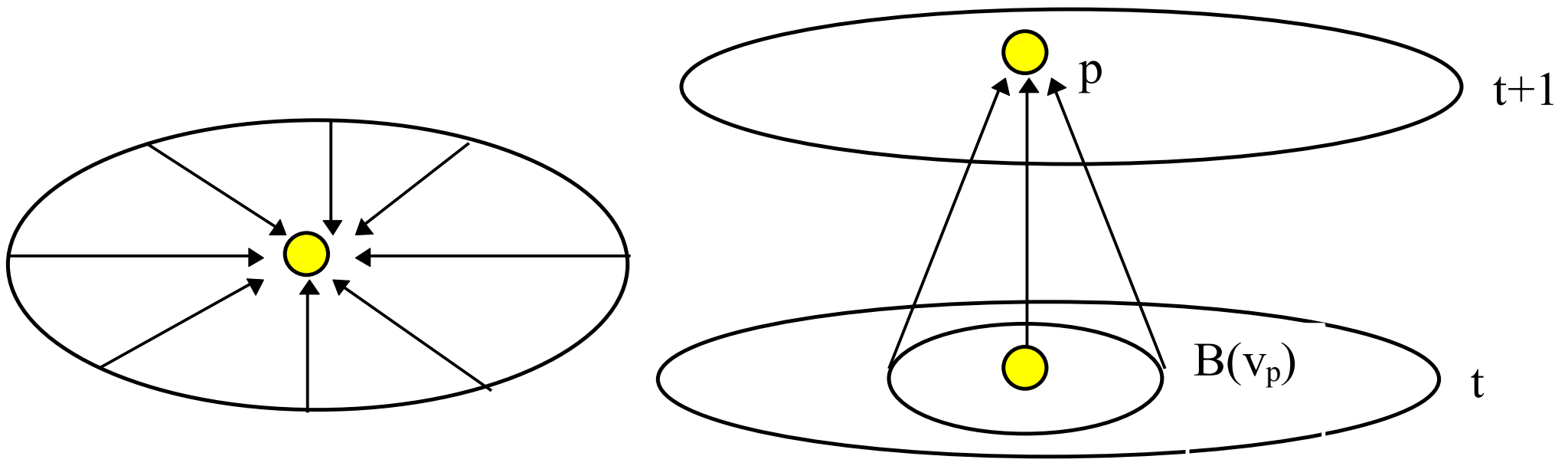

Рис. 190 Количество тепла, перешедшее из точки p в соседнюю точку к, равно $c_{kp}f_p^t$. Очевидно, что общее количество тепла, участвующее в этом действии, должно остаться неизменным.

смежным с точкой p, а во второй сумме по индексу p. В правой части конечно-разностного уравнения имеется только один слой, все коэффициенты в правой части неотрицательны, причем сумма коэффициентов по первому индексу равна единице. Это уравнение можно интерпретировать как, например, уравнение теплопроводности, следующим образом (Рис. 190). Предположим, что $f_p^t$ есть количество тепла в точке p в момент t. В следующий момент (t+1) часть этого тепла переходит в соседние точки, и какая-то часть остается в точке p. Количество тепла, перешедшее из точки p в соседнюю точку к, равно $c_{kp}f_p^t$. Очевидно, что общее количество тепла, участвующее в этом действии, должно остаться неизменным. При этом, для сохранения физической аналогии, будем считать, что $c_{kp}=c_{pk}$, хотя это и не обязательно. Вместо количества тепла можно рассматривать количество частиц, находящихся в данной точке в данный момент времени, и перешедших в соседние точки в следующий момент (уравнение диффузии). Еще одна наглядная иллюстрация-случайное блуждание частицы по пространству, где $c_{pk}$ есть вероятность перехода частицы из точки к в точку p [55].
Эллиптическое уравнение может быть представлено в такой же форме, как и параболическое.

$$f_p^{t+1} = \sum_k c_{pk} f_k^t, \ \forall c_{pk} \geq 0, \ \sum_p c_{pk} = 1,$$

где в первой сумме суммирование ведется по всем точкам шаблона, смежным с точкой p, а во второй сумме по индексу p.
Гиперболическое уравнение представим в виде

$$f_p^{t+1} = \sum_k c_{pk} f_k^t + f_p^t - f_p^{t-1}, \ \forall c_{pk} \geq 0, \ \sum_p c_{pk} = 1,$$

где в первой сумме суммирование ведется по всем точкам шаблона, смежным с точкой p, а во второй сумме по индексу p. Подчеркнем фундаментальный факт связи гиперболического и параболического уравнений

$$(f_p^{t+1})_{hyper} = (f_p^{t+1})_{parab} + f_p^t - f_p^{t-1}.$$

Исходя из этих результатов, определим различные типы уравнений математической физики на произвольном молекулярном пространстве. В сущности, нам остается обобщить уравнение (1) на молекулярное пространство.

Определение динамической системы.

Динамической системой $\{f^t\}$ на молекулярном пространстве G с точками $(v_1, v_2, \ldots v_n, \ldots)$ называется последовательность состояний (значений функции на точках пространства), определяемых номерами t=0,1,2,..., функции f, заданной на множестве точек пространства G.

Значение функции $f^t$ в точке $v_k$ пространства G будет обозначаться $f^t_k = f^t(v_k)$. t можно интерпретировать как, например, время.

Определение однородного линейного дифференциального уравнения второго (или менее) порядка по пространственным и координатам и (m+1) порядка по временной координате.

Дифференциальным линейным однородным уравнением в частных производных порядка не выше второго по пространственным координатам и порядка (m+1) по временной координате в молекулярном пространстве G называется уравнение вида

$$f_p^{t+1} = \sum_{s=t-m}^{s=t} \Big( \sum_{v_k \in B(v_p)} c^s_{pk} f^s_k \Big), \Leftrightarrow a\frac{\partial^{m+1} f}{\partial t^{m+1}} = b_1 \frac{\partial^2 f}{\partial x_1^2} + b_2 \frac{\partial^2 f}{\partial x_2^2} + \ldots b_n \frac{\partial^2 f}{\partial x_n^2}$$

определяющее состояние динамической системы в момент (t+1) через состояния в предыдущие моменты времени. Суммирование ведется по всем точкам шара $B(v_p)$.

Иными словами, значение функции f в точке $v_p$ в момент (t+1) определяется значениями функции в шаре $B(v_p)$ этой точки в предыдущие моменты времени (Рис. 190).

Имеет смысл дать определение неоднородного дифференциального уравнения, которое строится в чисто традиционном плане, добавлением функции к правой части уравнения.

Определение неоднородного линейного дифференциального уравнения второго (или менее) порядка по пространственным и координатам и (m+1) порядка по временной координате.

Дифференциальным линейным неоднородным уравнением в частных производных порядка не выше второго по пространственным координатам и порядка (t+1) по временной координате в молекулярном пространстве G называется уравнение вида

$$f_p^{t+1} = \sum_{s=t-m}^{s=t} \Big( \sum_{v_k \in B(v_p)} c_{pk}^s f_k^s \Big) + h_p^{t+1},$$

определяющее состояние динамической системы в момент (t+1) через состояния в предыдущие моменты времени. Суммирование ведется по всем точкам шара $B(v_p)$.

Как следствие из этого определения, введем основные типы уравнений.

Определение параболического дифференциального уравнениям на молекулярном пространстве G.

Параболическим дифференциальным уравнением на молекулярном пространстве G называется уравнение вида

$$f_p^{t+1} = \sum_{v_k \in B(v_p)} c_{pk}^t f_k^t + h_p^{t+1}, \ \forall c_{pk}^t \geq 0, \quad \sum_{v_k \in B(v_p)} c_{pk}^t = 1.$$

Суммирование ведется по всем точкам шара точки $v_p$, включая саму эту точку. Если все $h^{t+1}{}_p \equiv 0$, то уравнение называется однородным, в противном случае-неоднородным.

Точно также определяется эллиптическое уравнение.

Определение эллиптического дифференциального уравнениям на молекулярном пространстве G.

Эллиптическим дифференциальным уравнением на молекулярном пространстве G называется уравнение вида

$$f_p = \sum_{v_k \in B(v_p)} c_{pk} f_k + h_p, \ \forall c_{pk} \geq 0, \quad \sum_{v_k \in B(v_p)} c_{pk} = 1.$$

Суммирование ведется по всем точкам шара точки $v_p$, включая саму эту точку. Если все $h^{t+1}{}_p \equiv 0$, то уравнение называется однородным, в противном случае-неоднородным.

Перейдем, теперь, к определению гиперболического уравнения.

Рассмотрим гиперболическое уравнение, описывающее волновые процессы.

$$a\frac{\partial^2 f}{\partial t^2} = b_1\frac{\partial^2 f}{\partial x_1^2} + b_2\frac{\partial^2 f}{\partial x_2^2} + \dots b_n\frac{\partial^2 f}{\partial x_n^2},$$

$a, b_1, b_2, \dots b_n > 0.$

Определение гиперболического дифференциального уравнения на молекулярном пространстве G.

Гиперболическим дифференциальным уравнением на молекулярном пространстве G называется уравнение вида

$$f_p^{t+1} = \sum_{v_k \in B(v_p)} c_{pk}^t f_k^t + f_p^t - f_p^{t-1} + h_p^{t+1},\ \forall c_{pk}^t \geq 0,\ \sum_{v_k \in B(v_p)} c_{pk}^t = 1$$

Суммирование ведется по всем точкам шара точки $v_p$, включая саму эту точку. Если все $h^{t+1}{}_p \equiv 0$, то уравнение называется однородным, в противном случае-неоднородным.

Таким образом, правая часть волнового уравнения состоит из параболической части, представленной суммой, плюс значение функции в самой точке p в момент t, минус значение функции в точке p в момент (t-1). В этом фундаментальное отличие волнового уравнения от параболического. Как уже отмечалось, связь между уравнениями определяется следующей формулой.

$$(f_p^{t+1})_{hyper} = (f_p^{t+1})_{parab} + f_p^t - f_p^{t-1}.$$

Полезно рассмотреть, как меняются параболическое и гиперболическое решения от слоя к слою в целом. Введем функцию S, являющуюся суммой значений функции f по всем точкам пространства.

Определение интегральной суммы дифференциального уравнения.

Интегральной суммой дифференциального уравнения на молекулярном пространстве G с точками $(v_1, v_2, \dots v_n, \dots)$ называется сумма значений функции во всех точках пространства в момент времени t.

$$S(f^t) = \sum_{v_p \in G} f_p^t.$$

Легко показать, что для параболического уравнения

$$S(f^{t+1}) = S(f^t).$$

В то время как для гиперболического уравнения

$$S(f^{t+1}) = 2S(f^t) - S(f^{t-1}).$$

Величины $S(f^t)$ в гиперболическом случае образуют арифметическую прогрессию. Если разность прогрессии не равна 0, то $S(f^t)$ бесконечно возрастает по модулю с увеличением t. В конечном счете такое возрастание ведет к появлению расходимостей. Поэтому, мы должны исключить такой случай из рассмотрения. Чтобы это сделать, необходимо правильно определить начальные условия. Для гиперболического уравнения начальные условия определяются значениями исходной функции f на первых двух слоях при t=0 и 1 (или t=-1 и 0). Начальные условия всегда должны задаваться таким образом, чтобы выполнялось

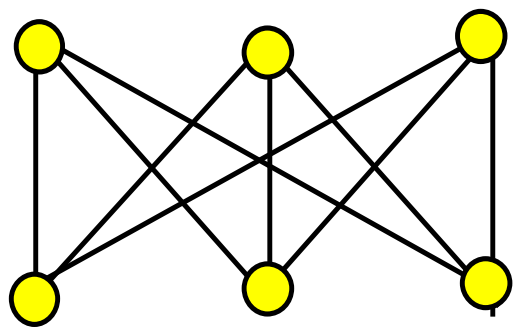

Рис. 191 Пространство не имеет непрерывного аналога и не является многообразием.

соотношение

$$S(f^0)= S(f^1).$$

При выполнении этого соотношения выполняется

$$S(f^t)= S(f^0)= S(f^1).$$

для любого t, и решения на молекулярных пространствах соответствуют решениям на их непрерывных прообразах. Нарушение этого соотношения ведет к появлению расходимостей, что было проверено при непосредственных расчетах на компьютере. Примеры мы дадим ниже.

## *3.2 ДИФФЕРЕНЦИАЛЬНЫЕ УРАВНЕНИЯ НА НЕПРЕРЫВНЫХ ПРОСТРАНСТВАХ ЯВЛЯЮТСЯ ЧАСТНЫМ СЛУЧАЕМ ДИФФЕРЕНЦИАЛЬНЫХ УРАВНЕНИЙ НА МОЛЕКУЛЯРНЫХ ПРОСТРАНСТВАХ?*

Создается впечатление, что дифференциальные уравнения на непрерывных пространствах могут рассматриваться как частные случаи дифференциальных уравнений на молекулярных пространствах. Поясним это следующим рассуждением. Для каждого дифференциального уравнения на непрерывном пространстве можно выбрать конечно-разностную схему и, следовательно, получить его дискретный аналог на молекулярном пространстве.

Однако, не для всякого молекулярного пространства существует его непрерывный образ (или прообраз). Например, пространство на Рис. 191 не имеет стандартного непрерывного прообраза, так как не является многообразием (в одномерном многообразии окаем каждой точки должен содержать две изолированных точки, а не три, как это имеет место в данном случае). Это означает, что дифференциальное уравнение на таком пространстве не имеет своего непрерывного аналога в евклидовом пространстве. Иными словами, дифференциальные уравнения на молекулярных пространствах включают в себя дифференциальные уравнения на непрерывных пространствах как частный случай.

# *ПРИМЕРЫ РЕШЕНИЙ УРАВНЕНИЙ МАТЕМАТИЧЕСКОЙ ФИЗИКИ НА МОЛЕКУЛЯРНЫХ ПРОСТРАНСТВАХ*

## 3.2.1РЕШЕНИЕ ОДНОРОДНОГО ПАРАБОЛИЧЕСКОГО УРАВНЕНИЯ НА ДВУМЕРНОЙ СФЕРЕ. ЗАДАЧА КОШИ

В качестве примера рассмотрим решение однородного параболического уравнения (уравнения распространения тепла) на двумерной сфере, изображенной на Рис. 192.

$$f_p^{t+1} = \sum_{v_k \in B(v_p)} c_{pk}^t f_k^t, \ \forall c_{pk}^t \geq 0, \quad \sum_{v_k \in B(v_p)} c_{pk}^t = 1.$$

Обозначим точки числами 1,2,3,4,5 и 6, начиная с верхней и кончая нижней, и выберем матрицу коэффициентов $c_{pk}$ в уравнении в следующем виде: $c_{kk}=0.7$, $c_{kp}=0.075$, k,p=1,2,...6.

Начальные условия (задача Коши) задаются единицей в верхней вершине и нулями во всех остальных: $f^0_1=1$, $f^0_s=0$, s= 2,3,5,5,6. Зависимость решения f от времени в точках 1 и 6 показана на графике (Рис. 192). Аналогично получается решение на n-мерной сфере.

## 3.2.2РЕШЕНИЕ ОДНОРОДНОГО ПАРАБОЛИЧЕСКОГО УРАВНЕНИЯ НА ДВУМЕРНОЙ СФЕРЕ. КРАЕВАЯ ЗАДАЧА

До сих пор мы рассматривали решения параболических уравнений без краевых условий. В этом разделе мы определим один из видов краевых условий, соответствующий заданию определенных значений функции f в некоторых точках пространства в течение всего процесса. Поясним это следующим примером. Предположим, что в некоторых точках пространства

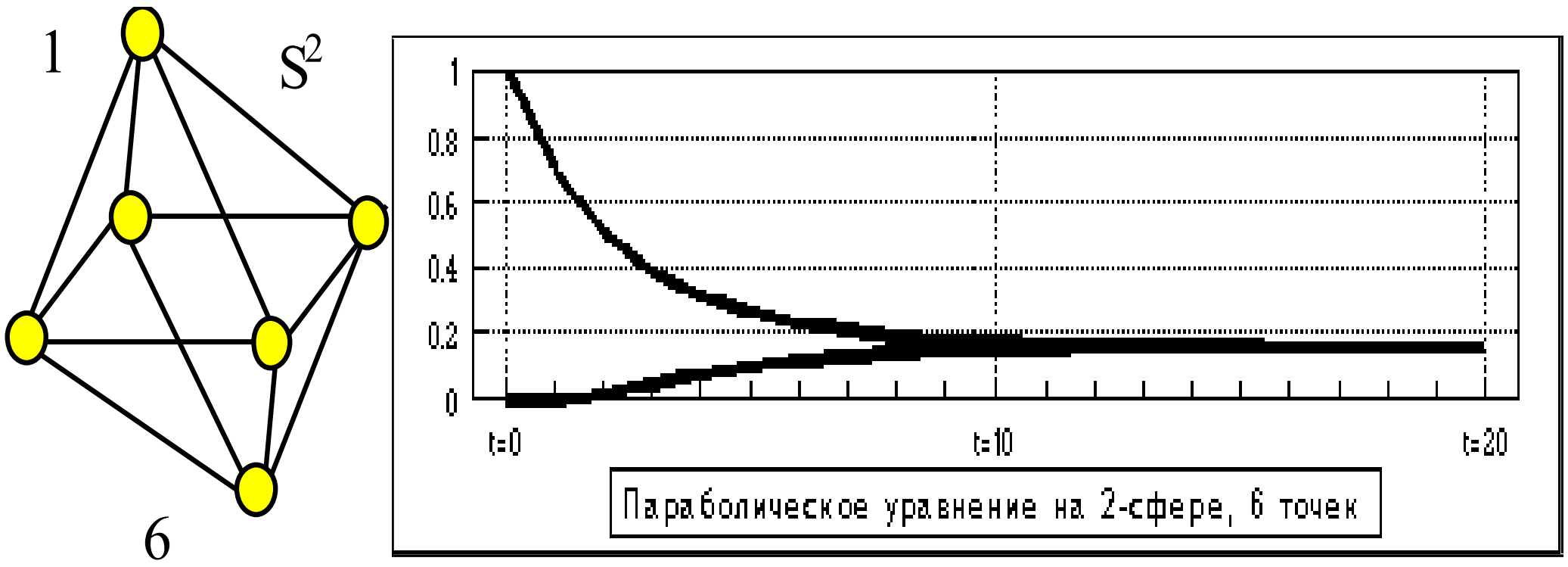

Рис. 192 График решения параболического уравнения на двумерной сфере в противоположных точках 1 и 6.

G задана зависимость температуры, определяемая какими-то внешними условиями, в течение всего процесса. Необходимо найти зависимость температуры от времени в произвольной точке. Возникает типичная краевая задача для параболического уравнения. Здесь нам понадобятся понятия истока и стока.

Определение истока

Точка $v_k$ называется истоком, если все $c_{kp}$ =0 для любого p..

В исток вообще не приходят значения из соседних точек, независимо от их знака. В то же время из истока в соседние точки при каждом шаге могут переходить значения, положительные или отрицательные. Это означает,

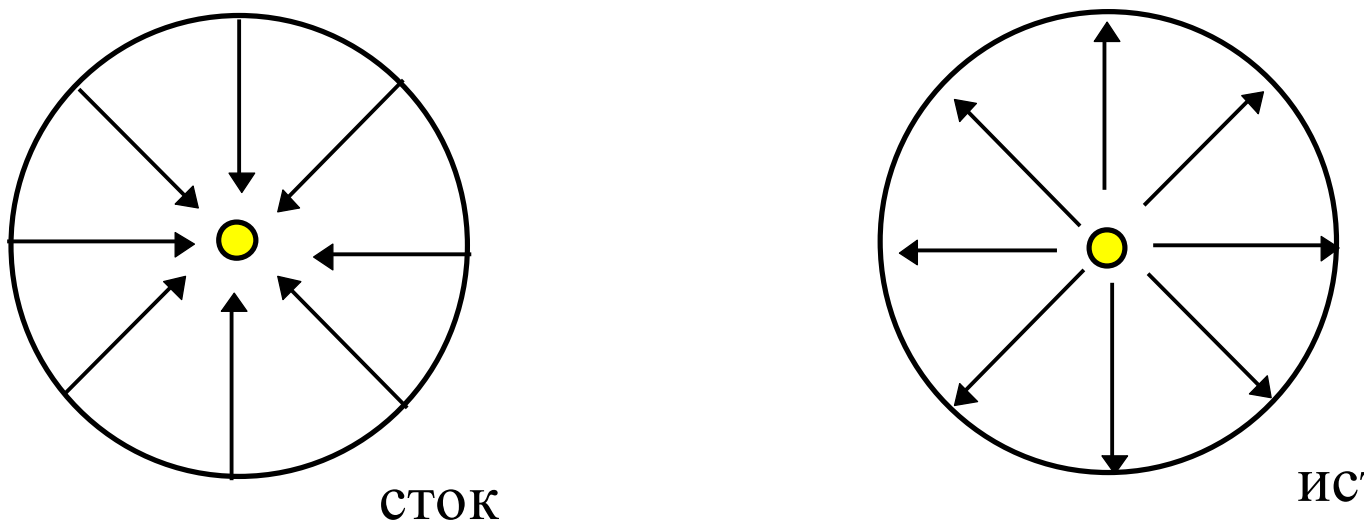

Рис. 193 В стоке p все коэффициенты $c_{kp}$=0, в истоке p все $c_{pk}$=0.

что $c_{pk}$ могут быть отличны от 0.

Определение стока

Точка $v_k$ называется стоком, если все $c_{pk}$=0 для любого p.

В сток могут приходить значения из соседних точек, независимо от их знака. В то же время из стока в соседние точки не могут переходить значения, положительные или отрицательные. В этом случае могут быть

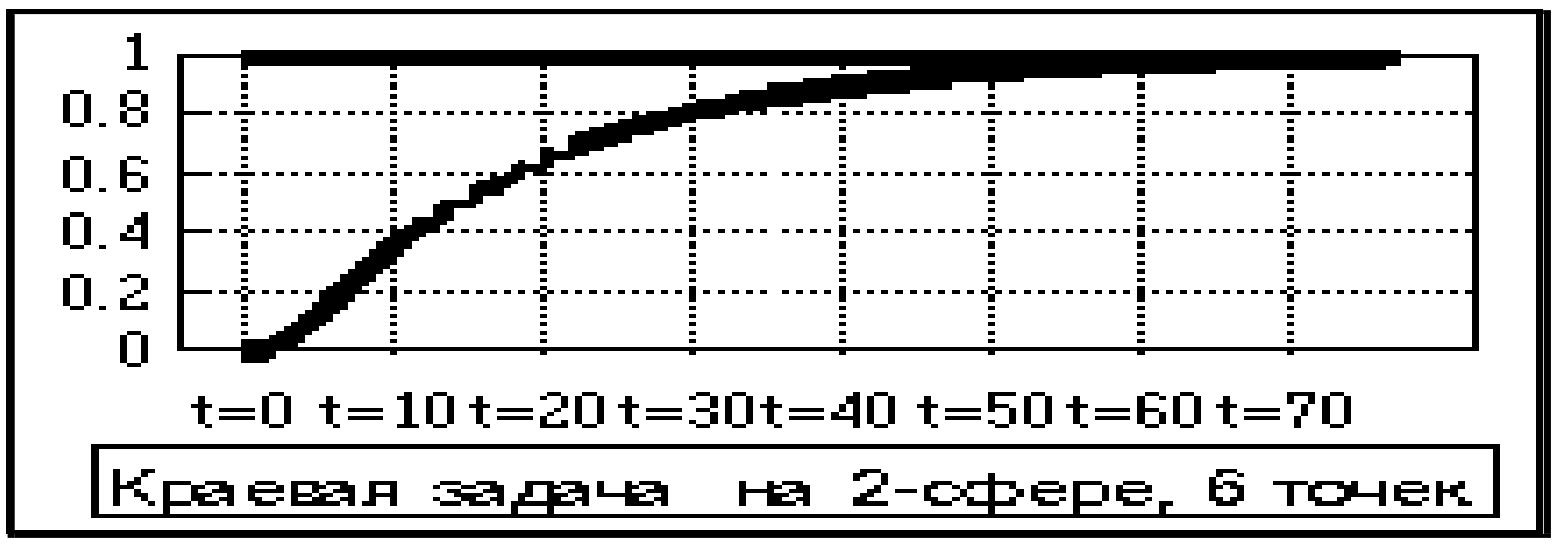

Рис. 194 График решения краевой задачи на двумерной сфере в истоке и противоположной ему точке.

отличными от 0 коэффициенты $c_{kp}$.

Условно сток и исток изображены на Рис. 193.

Введем теперь определение краевой задачи предполагая, что на пространстве заданы истоки и стоки.

Определение краевой задачи типа 1

Краевой задачей типа 1 на пространстве G называется:

1. параболическое уравнение, заданное на этом пространстве,
2. множество истоков, стоков и функций, определенных в истоках.

Рассмотрим пример краевой задачи распространения тепла на двумерной сфере, изображенной на Рис. 192. Обозначим вершины числами 1,2,3,4,5 и 6, начиная с верхней и кончая нижней. Точка 1 на сфере является истоком, в котором задано постоянное значение функции $f_1^t=1$ (температура) для любого момента t. Введем коэффициенты, связанные с истоком и остальными точками следующим образом:

$c_{1p}^t=0$, $c_{p1}^t=0.075$, где p=2,3,4,5; $c_{kk}^t=0.7$, $c_{kp}^t=0.075$, k,p=,2,...6.

Начальные значения функции f во всех точках сферы выберем равными 0, $f^0_s=0$, s= 2,3,5,5,6. Уравнение теплопроводности имеет решения f, зависимость которых от времени в точках 1 и 6 показана на графике (Рис. 194).

В точке 1 функция постоянна и равна 1. В точке 6 значение функции постепенно возрастает от 0 в начальный момент времени и до 1.

## 3.2.ЗРЕШЕНИЯ ПАРАБОЛИЧЕСКОГО УРАВНЕНИЯ НА N-МЕРНОМ ПОЛНОМ ЕВКЛИДОВОМ ПРОСТРАНСТВЕ $L^n$ И МОЛЕКУЛЯРНОМ ДЕРЕВЕ $E^1_r$

В [43] были получены точные решения параболического уравнения на на n-мерном полном евклидовом пространстве $L^n$ и на однородном дереве. Однородное дерево $E^1_r$ вида r есть однородное молекулярное пространство без циклов, в котором каждая точка имеет r соседних точек. На Рис. 195 изображено двумерное однородное регулярное евклидово молекулярное пространство $L^2$ и одномерное дерево $E^1_3$ со степенью 3. Это пространство не является многообразием и не имеет прообраза в виде одномерной непрерывной кривой в евклидовом пространстве. Было показано, что на $E^1_r$ мы можем получить решение, которое сходно с решениями, полученными на n-мерном евклидовом пространстве $L^n$.

Если в $L^n$ имеется точка, являющаяся истоком, то решение

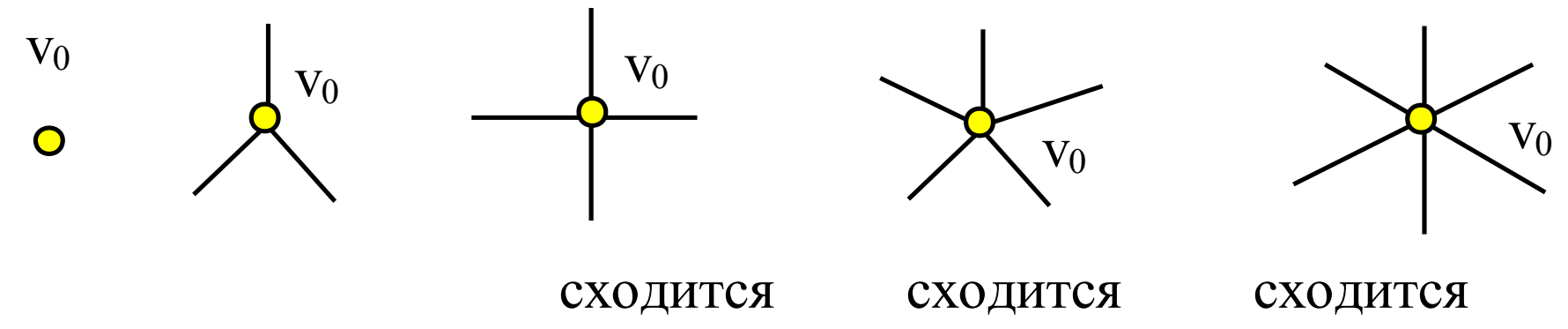

Рис. 195 Параболическое уравнение не имеет бесконечной особенности в истоке, если каждая точка имеет 4 или более смежных точек.

параболического уравнения не имеет расходимости в этой точке при t=∞ при размерности 3 и выше. В протранстве $E^1_r$ также с единственным истоком, расходимость в истоке при t=∞ отсутствует при r≥3. При физической интерпретации в обоих случаях можно считать, что в точке $v_0$ в каждый момент времени t = 0,1,2,... n,... выделяется количество тепла, равное 1, и далее это тепло распространяется по всему пространству. Если $n \geq 3$ в $L^n$ и если $r \geq 4$ в $E^1_r$ (Рис. 196), то количество тепла в $v_0$ и всех

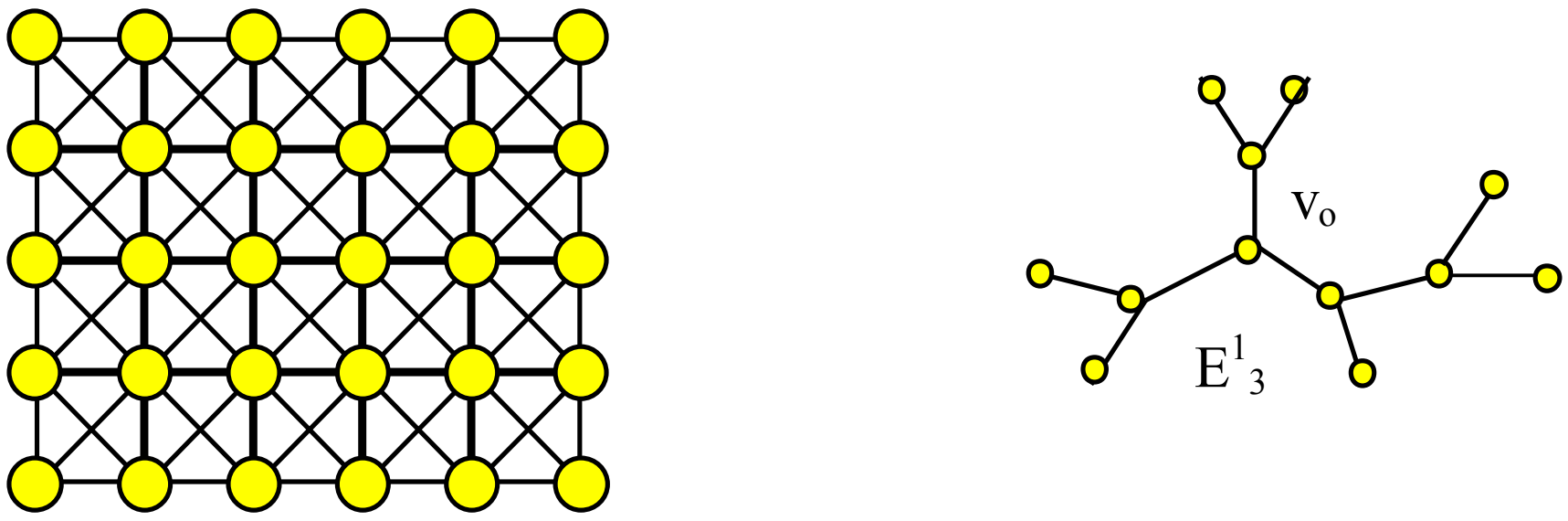

Рис. 196 Справа изображено двумерное полное евклидово пространство, слева-однородное дерево типа 3.

остальных точках в стационарном режиме при $t=\infty$ будет конечной величиной. Топологически пространства $E^1_r$ при $r \geq 4$ и $L^n$ различны. Более того, $E^1_r$ при $r \geq 4$ даже не является многообразием. Однако, в смысле сходимости в источнике тепла, решения на $E^1_r$ при $r \geq 4$ похожи на решения на евклидовых пространствах $L^n$ при $n \geq 3$.

Возникает вопрос, можно ли, опираясь только на вид решения, определить, в евклидовом пространстве или же в одномерном дереве происходит данный процесс. Сходство в виде решений говорит о том, что скорее всего невозможно различить топологию пространства по виду решения. То есть топология не влияет существенно на вид решения. Дифференциальные уравнения являются общими конструкциями, не зависящими от топологии.

### 3.2.4РЕШЕНИЕ ГИПЕРБОЛИЧЕСКОГО УРАВНЕНИЯ НА ОДНОМЕРНОЙ СФЕРЕ

.В качестве примера рассмотрим решение однородного гиперболического уравнения

$$f_p^{t+1} = \sum_{v_k \in B(v_p)} c_{pk} f_k^t + f_p^t - f_p^{t-1}, \text{ где все } c_{pk} \geq 0, \quad \sum_{v_k \in B(v_p)} c_{pk} = 1$$

на одномерной сфере. Пусть одномерная нормальная сфера состоит из 25 точек (Рис. 197). Выберем в уравнении коэффициенты в любой момент времени следующим образом: $c^t_{kk}=0.8$, $c^t_{kp}=0.1$, если $v_k$ и $v_p$ смежны, k,p=1,2,..25. Начальные условия задаются единицей в точке a при t=-1.0, и нулями во всех остальных: $f^0_a = f^{-1}_a = 1$. Физически это означает, что в точке a задано единичное смещение и нулевая скорость. Зависимость решения f от времени в точках a и b, находящейся на 12 позиции (удаленной от a), показана на Рис. 197. Из графика видно, что начальное возмущение из точки a распространяется по точкам окружности. Первое, отличное от нуля значение функции, появляется в точке b в момент t=11.

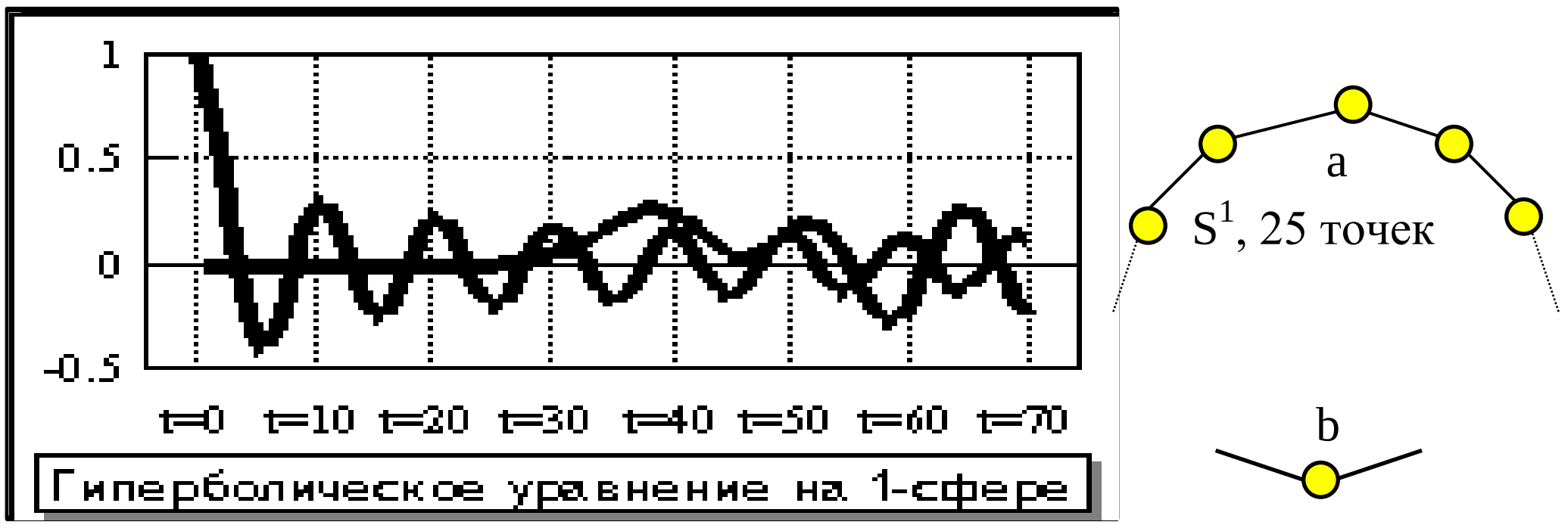

Рис. 197 График решения гиперболического уравнения на одномерной сфере в точках a и b.

### 3.2.5РЕШЕНИЕ ГИПЕРБОЛИЧЕСКОГО УРАВНЕНИЯ ДВУМЕРНОЙ ПРОЕКТИВНОЙ ПЛОСКОСТИ $P^2$

В [43] мы подробно описали решения гиперболического уравнения на двумерной сфере $S^2$, трехмерной сфере $S^3$ и двумерном торе $T^2$. Рассмотрим здесь только проективную плоскость. Если имеется такое g, что норма ||f-g||<e для любого слоя t, то дифференциальное уравнение

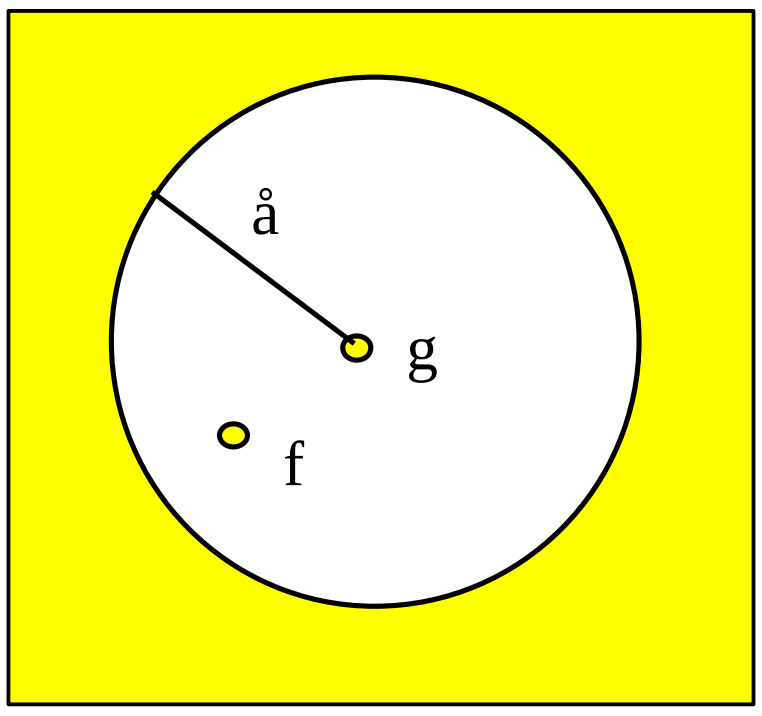

Рис. 198 Если имеется такое g, что норма ||f-g||<e для любого слоя t, то дифференциальное уравнение невырождено.

невырождено (Рис. 193)

Начальные условия определяются выражениями: при t=-1,0, $f^0_a=f^{-1}_a=1$Выберем в уравнении коэффициенты в точке $v_p$ в любой момент времени следующим образом:

$$c^t_{pp} = c^t_{kp}=1/(n+1),$$

где n-число точек, соседних точке $v_p$. Решение представлено соответствующим графиком (Рис. 194).. Одна кривая показывает решение в точке a, другая-в точке b. Для проективной плоскости построены два графика. Первый график получен для начальных условий при t=-1,0, $f^0_a= f^{-1}_a =1$, а второй для начальных условий, $f^0_a=1$, $f^{-1}_a=0$. Напомним, что это означает, нарушается равенство интегральных сумм на слоях.. Легко убедиться, что во втором случае возникает монотонное увеличение значений, на которое накладывается колебательный процесс.

Рис. 199 Графики решений гиперболического уравнения $P^2$.

Если нас интересуют чисто волновые процессы, мы можем исключить такие решения правильным выбором начальных условий..

## 3.3 НЕВЫРОЖДЕННЫЕ ДИФФЕРЕНЦИАЛЬНЫЕ УРАВНЕНИЯ НА МОЛЕКУЛЯРНЫХ ПРОСТРАНСТВАХ

Одним из важных вопросов в теории и применении дифференциальных уравнений является невырожденность их решений. Существует много

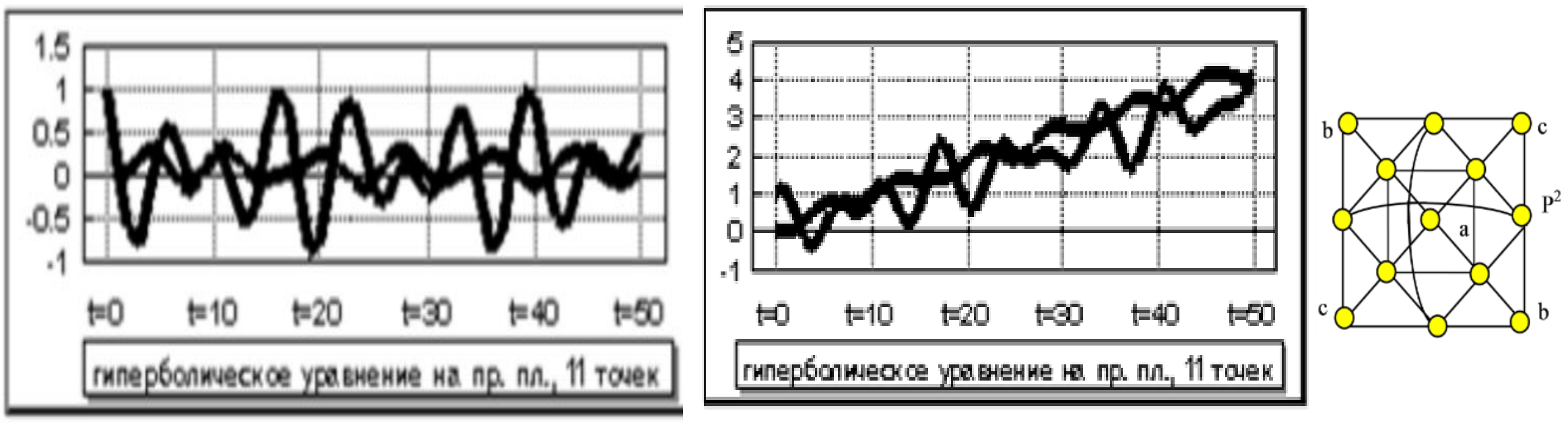

Рис 194. Графическое решение гиперболического уравнения на пр. пл.

определений невырожденности, в том числе и для конечно-разностных схем. Большинство критериев рассчитаны на непрерывную базу, позволяющую выбрать подходящий шаблон. С этой точки зрения все они мало пригодны к применению на молекулярных пространствах, где шаблон определяется самой структурой пространства. Помимо этого, нам необходим критерий, который легко проверяется на компьютере. Поэтому мы введем определение невырожденности дифференциального уравнения, основанное на очевидных предпосылках, и определим условия, при которых оно выполняется.

<u>Определение дифференциального однородного уравнениям первого порядка по времени на молекулярном пространстве G</u>

Дифференциальным однородным уравнением первого порядка по времени на молекулярном пространстве G называется уравнение вида

$$f_p^{t+1} = \sum_{v_k \in B(v_p)} c^t{}_{pk} f_k^t, \quad (1)$$

где $f^t_p = f^t(v_p)$, $c^t_{pk}$ есть произвольные коэффициенты. Как уже говорилось, состоянием $f^t$ решения уравнения (1) называется множество значений функции f на слое t, $f^t = (f^t_1, f^t_2, \ldots f^t_k, \ldots)$.

Множество всех состояний образует динамическую систему f данного уравнения.

Для данного уравнения (1) существует множество динамических систем. По сути дела, динамическая система уравнения полностью задается начальным состоянием $f^0$. В этом параграфе мы будем говорить только о динамических системах, являющихся решениями уравнения (1). Введем понятие нормы состояния динамической системы.

Определение нормы состояния динамической системы.

Нормой состояния t решения дифференциального уравнения (1) называется величина $\|f^t\|$, являющаяся суммой абсолютных величин значений f по всем точкам на данном слое t и определяемая выражением:

$$\|f^t\| = \sum_k |f_k^t|,$$

где суммирование ведется по всем вершинам молекулярного пространства.

Определение невырожденного дифференциального однородного уравнениям первого порядка по времени на молекулярном пространстве G

Пусть для уравнения вида (1) $f_p^{t+1} = \sum_{v_k \in B(v_p)} c_{pk} f_k^t$,.

существует динамическая система g и число e>0 такие, что для любой динамической системы f, удовлетворяющей условию $\|f^0 - g^0\| < e$, выполняется равенство $\|f^t\| = \|f^0\|$ для любого t. Тогда уравнение (1) называется невырожденным.

Смысл невырожденности понятен из Рис. 198. Если f принадлежит кругу радиуса e с центром в g, то, во-первых, норма динамической системы всегда постоянна, и, во-вторых, малые изменения в начальном состоянии вызывают также малые изменения в любом из последующих состояний.

Посмотрим, как ведет себя норма в зависимости от начальных условий, если дифференциальное уравнение является параболическим.

Теорема 172

*Параболическое уравнение*

$$f_p^{t+1} = \sum_{v_k \in B(v_p)} c_{pk} f_k^t, \text{ где } c_{pk} \ge 0, \quad \sum_{v_p \in B(v_k)} c_{pk} = 1 .$$

*есть невырожденное уравнение. Кроме того, с начальными условиями $f^0_p = f^0(v_p) \ge 0$ для любой точки $v_p$ молекулярного*

*пространства норма решения не зависит от времени (состояния) t,* $||f^t||=||f^0||$.

Доказательство.

Рассмотрим начальное состояние $g^0$ системы, при котором значения функции во всех точках равны 1. Рассмотрим семейство $\{e_k\}$, где $|e_k|<1$ для любого k. Пусть $f^0(v_k)=g^0(v_k)+e_k=1+e_k\geq 0$. Тогда $f_p^1 = \sum_{v_k \in B(v_p)} c_{pk} f_k^0 \geq 0$.

Следовательно,

$$||f^1|| = \sum_p |f_p^1| = \sum_p f_k^1 = \sum_p \sum_{v_k \in B(v_p)} c_{pk} f_k^0 = \sum_k f_k^0 \sum_{v_p \in B(v_k)} c_{pk} = \sum_k f_k^0 = ||f^0||.$$

Точно также получается $||f^{t+1}||=||f^t||$ для любого t. Отсюда $||f^t||=||f^0||$. Теорема доказана. □

Теперь мы найдем более широкий класс уравнений вида (1), для которых существуют условия вида $f^0{}_p=f^0(v_p)\geq 0$, сохраняющие норму. Эти уравнения можно назвать квазипараболическими, так как по своим свойствам они напоминают параболические уравнения. Рассмотрим некоторые условия невырожденности дифференциального уравнения. Для этого представим уравнение (1) в матричной форме. При этом мы можем считать, что если две точки $v_p$ и $v_k$ пространства несмежны, то $c_{pk}=c_{kp}=0$.

$$f^{t+1}=Af^t$$

$$f^{t+1} = \begin{bmatrix} f_1^{t+1} \\ f_2^{t+1} \\ f_3^{t+1} \\ * \end{bmatrix}; \quad A = \begin{bmatrix} c_{11} & c_{12} & c_{13} & * \\ c_{21} & c_{21} & c_{23} & * \\ c_{31} & c_{32} & c_{13} & * \\ * & * & * & * \end{bmatrix}; \quad f^t = \begin{bmatrix} f_1^t \\ f_2^t \\ f_n^t \\ * \end{bmatrix};$$

Теорема 173

*Пусть дифференциальное уравнение удовлетворяет следующим условиям:*

*1.* $\sum_{v_k \in B(v_p)} |c_{kp}| = 1$ *для любой вершины* $v_p$..

*2. Матрица A этой системы путем изменения нумерации вершин пространства может быть представлена в виде суммы четырех матриц B(1), B(2), C(1) и C(2), где B(1) и B(2) не имеют отрицательных элементов и общих строк или столбцов, матрицы C(1) и C(2) не имеют положительных элементов и общих строк и столбцов.*

*Тогда уравнение (1) является невырожденным.*

A =
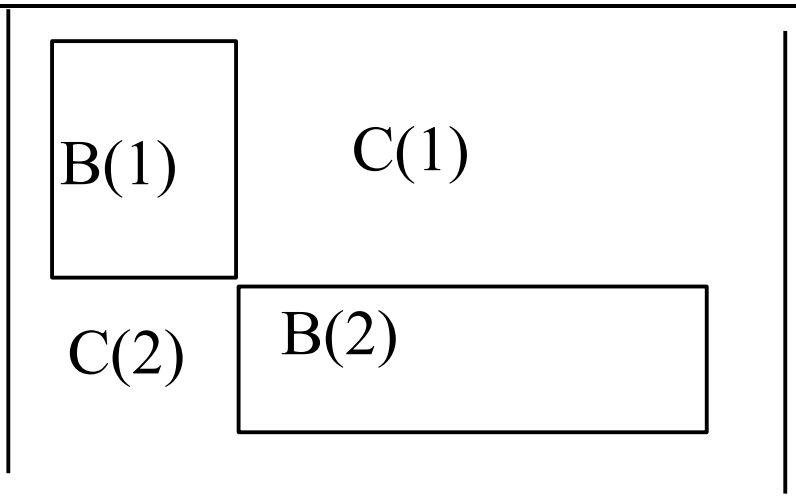

*При этом если B(1) имеет a строк и b столбцов, то в слое $g^0=(g^0_1,g^0_2,..g^0_k,...)$ первые b элементов положительны, остальные отрицательны.*

Доказательство

Любое значение функции на первом слое может быть записано в виде

$$g^1_p=\sum_s c_{pk}g^0_s=\sum_{s\le b} c_{pk}g^0_s+\sum_{s>b} c_{pk}g^0_s.$$

Рассмотрим $g^1_p$, $1\le p\le a$. Тогда $g^1_p=\sum_s c_{pk}g^0_s=\sum_{s\le b} c_{pk}g^0_s+\sum_{s>b} c_{pk}g^0_s\ge 0$, и все слагаемые в сумме положительны. Рассмотрим $g^1_p$, $p>a$. Тогда $g^1_p=\sum_s c_{pk}g^0_s=\sum_{s\le b} c_{pk}g^0_s+\sum_{s>b} c_{pk}g^0_s\le 0$, и все слагаемые в сумме отрицательны. Определим, теперь, норму на первом слое.

$$\|g^1\|=\sum_p |g^1_p|=\sum_{p\le a}|g^1_p|+\sum_{p>a}|g^1_p|=\sum_{p\le a}\sum_s |c_{ps}g^0_s|+\sum_{p>a}\sum_s |c_{ps}g^0_s|=$$

$$\sum_p\sum_s |c_{ps}g^0_s|=\sum_s |g^0_s|\sum_p |c_{ps}|=\sum_s |g^0_s|=\|g^0\|.$$

Так как знаки в точках для слоя 1 такие же как и знаки в тех же точках на слое 0, то $\|g^2\|=\|g^1\|$, и $\|g^{t+1}\|=\|g^t\|$. Следовательно, $\|g^t\|=\|g^0\|$ для любого t. Пусть в слое $g^0$ слое $g^0=(g^0_1, g^0_2,... g^0_k,...)$ первые b элементов равны 1, остальные равны (-1). Тогда для любого набора чисел $e_k$, $|e_k|<1$, k=1,2,..., начальный слой $f^0$ со значениями в точках $f^0_k=g^0_k+e_k$, имеет те же знаки в точках, что и $g^0$. Следовательно, $\|f^t\|=\|f^0\|$ для любого t. Теорема доказана. □

Если матрица A может быть приведена к указанной выше форме, она будет называться невырожденной матрицей. Необходимо отметить, что в невырожденном дифференциальном уравнении произвольный нулевой слой может быть представлен как сумма двух слоев $g^0_1$ и $g^0_2$, у первого слоя первые r элементов положительны, остальные отрицательны, во втором слое $g^0_2$ первые r элементов, наоборот, отрицательны, остальные положительны.

$$f^0=g^0_1+g^0_2$$

Для каждого из этих слоев, как легко видеть, норма не зависит от t.

$$\|g^t_1\|=\|g^0_1\|,\ \|g^t_2\|=\|g^0_2\|.$$

Простой проверкой несложно показать, что норма динамической системы f на любом слое всегда ограничена

$$\|f^t\| \le \|g^0_1\| + \|g^0_2\| = \|f^0\|.$$

Сформулируем это свойство в виде теоремы.

Теорема 174

*В невырожденном дифференциальном уравнении для любой динамической системы* $f^t=(f^t_1, f^t_2,\dots f^t_k,\dots)$ *ее норма на любом слое всегда ограничена, то есть* $\|f^t\| \le \|f^0\|$.

Доказательство.

Пусть нулевой слой представлен как сумма двух слоев $g^0_1$ и $g^0_2$, у первого слоя первые r элементов положительны, остальные отрицательны, во втором слое $g^0_2$ первые r элементов, наоборот, отрицательны, остальные положительны, $f^0=g^0_1+g^0_2$. Тогда во первых $\|g^t_1\|=\|g^0_1\|$, $\|g^t_2\|=\|g^0_2\|$, во-вторых $\|f^t\| \le \|g^t_1\|+\|g^t_2\|$, в третьих $\|f^t\| \le \|g^t_1\|+\|g^t_2\|$. Отсюда $\|f^t\| \le \|g^t_1\|+\|g^t_2\| \le \|g^0_1\|+\|g^0_2\|=\|f^0\|$. Теорема доказана. □

Замечание.

Параболическое уравнение является невырожденным, так как все элементы матрицы A в этом случае являются неотрицательными.

В качестве примера рассмотрим численное решение невырожденного дифференциального уравнения на двумерной сфере. Выберем двумерную

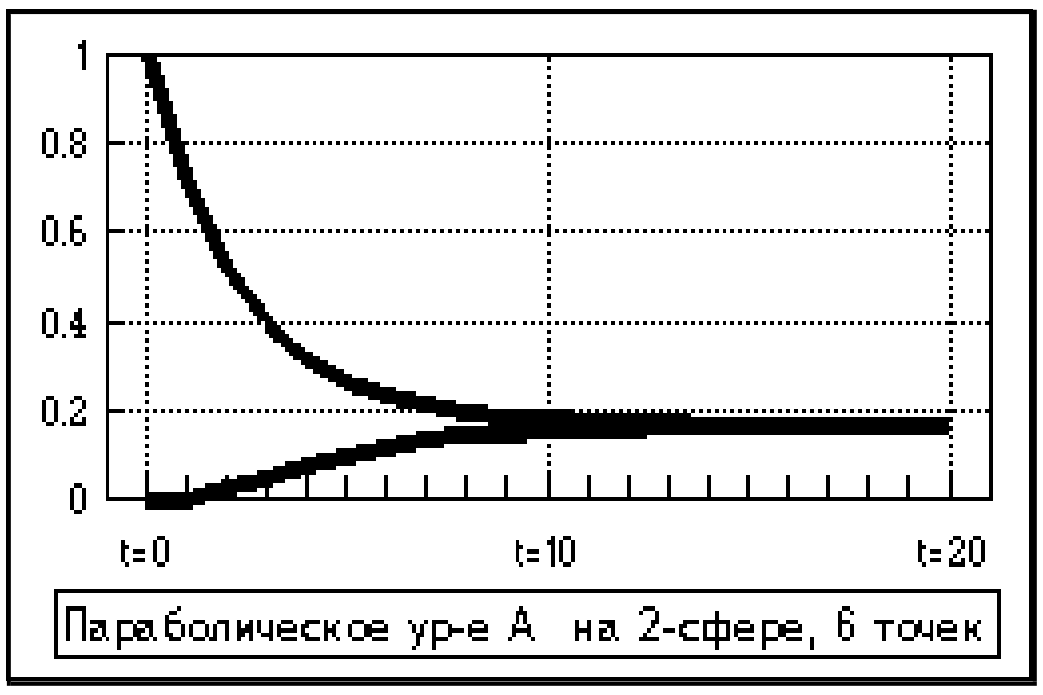

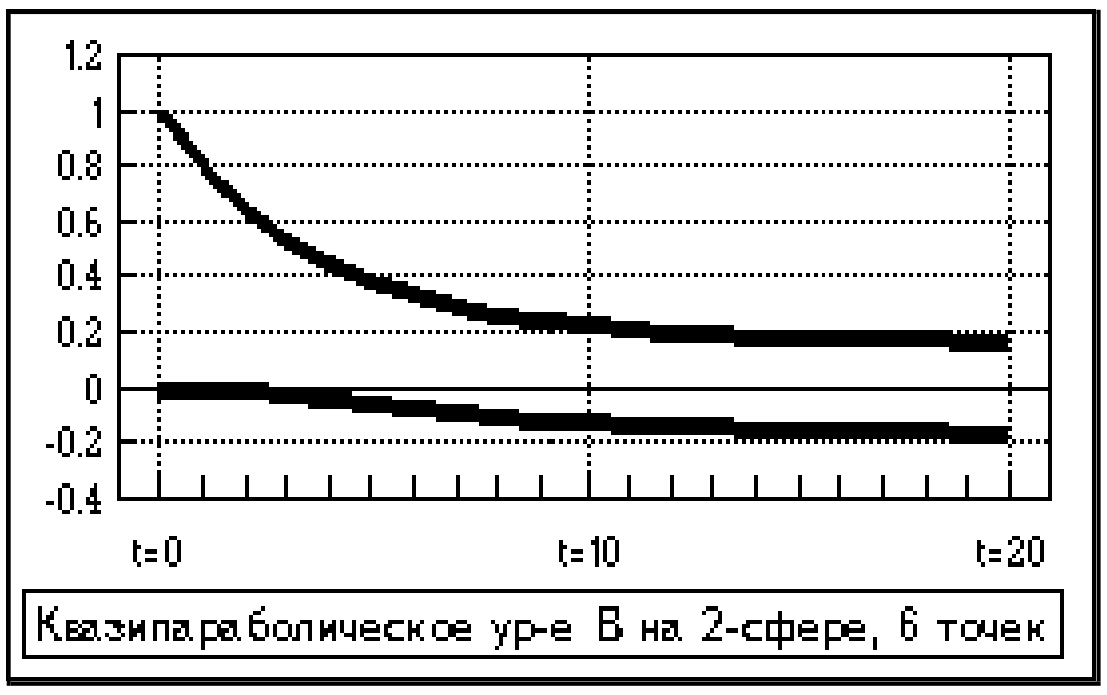

Рис. 200 Графики решений параболического (слева) и квазипараболического (справа) уравнений при одинаковых начальных условиях на двумерной сфере.

сферу содержащую 6 точек (изображена выше), и рассмотрим на ней решения двух уравнений: параболического, определяемого матрицей A, и невырожденного, не являющегося параболическим, определяемого матрицей B. Эти матрицы представлены ниже. Начальный слой одинаков в обоих случаях, $f^0_1=1$, $f^0_k=0$, k=2,3,...6. Решения уравнений для двух точек показаны на Рис. 200, слева для случая A, и справа для случая B.

A=

| | 1 | 2 | 3 | 4 | 5 | 6 |
|---|---|---|---|---|---|---|
| 1 | 0.8 | 0.05 | 0.05 | 0.05 | 0.05 | |
| 2 | 0.05 | 0.8 | 0.05 | | 0.05 | 0.05 |
| 3 | 0.05 | 0.05 | 0.8 | 0.05 | | 0.05 |
| 4 | 0.05 | | 0.05 | 0.8 | 0.05 | 0.05 |
| 5 | 0.05 | 0.05 | | 0.05 | 0.8 | 0.05 |
| 6 | | 0.05 | 0.05 | 0.05 | 0.05 | 0.8 |

B =

| | 1 | 2 | 3 | 4 | 5 | 6 |
|---|---|---|---|---|---|---|
| 1 | 0.8 | 0.05 | 0.05 | 0.05 | -0.05 | |
| 2 | 0.05 | 0.8 | 0.05 | | -0.05 | -0.05 |
| 3 | 0.05 | 0.05 | 0.8 | 0.05 | | -0.05 |
| 4 | 0.05 | | 0.05 | 0.8 | -0.05 | -0.05 |
| 5 | -0.05 | -0.05 | | -0.05 | 0.8 | 0.05 |
| 6 | | -0.05 | -0.05 | -0.05 | 0.05 | 0.8 |

Невырожденные уравнения расширяют класс уравнений, сходных с параболическими по своим свойствам, но не являющимися таковыми. Непрерывный аналог таких уравнений, повидимому, может быть представлен, например, в виде $\frac{\partial f}{\partial t}=\frac{\partial^2 f}{\partial x^2}-\frac{\partial^2 f}{\partial y^2}$.

## Список литературы к главе 16.

# ОБЩИЕ ПРЕОБРАЗОВАНИЯ МОЛЕКУЛЯРНЫХ ПРОСТРАНСТВ И ИХ ИНВАРИАНТЫ

In this chapter we introduce a linear characteristic of a molecular space that generalize the Euler characteristic. We study some properties of this function. In connection with that we define transformations of a molecular space that generalize contractible transformations. We prove that for any given transformation there exists the characteristic function of a molecular space which is not changed by these transformations. In conclusion we discuss possible applications of this approach.

Полученные в этой главе результаты частично изложены в [43].

## *ХАРАКТЕРИСТИЧЕСКИЕ ФУНКЦИИ НА МОЛЕКУЛЯРНЫХ ПРОСТРАНСТВАХ*

В этом разделе мы рассмотрим линейные характеристики молекулярных пространств и покажем связь некоторых из них с различными преобразованиями пространств. Мы обнаружили интересное свойство. Оказывается, что для каждого типа преобразований молекулярного пространства существует класс таких характеристик, причем вполне определенный, которые являются инвариантами этих преобразований. В качестве примера приведем эйлерову характеристику молекулярного пространства, которая не меняется при точечных преобразованиях. Создается впечатление, что для широкого, если не любого класса преобразований существуют такие характеристики. Для точечных преобразований мы имеем классический топологический аналог. Точечные преобразования являются дискретными моделями гомотопных трансформаций топологических пространств. Для других видов преобразований мы не можем пока указать их аналоги в классической математике. Однако, не исключено, что, в компьютерных науках, физике или биологии они описывают какие-то реальные процессы и могут найти применение. Обнаружение одного инварианта преобразования позволяет предположить что существуют другие инварианты, как например, группы. Ранее мы показали, что группы гомологий МП являются инвариантами точечных трансформаций. Возможно, что для других трансформаций также существуют аналогичные группы, но это уже вопрос дальнейших исследований.

Переходим к изложению. Напомним несколько определений, которые мы ввели ранее.

Определение связки K(n).

Пространство K(n) на n точках называется n связкой или полностью связным пространством, если каждая пара его точек смежна.

Определение множества k связок $T_k(G)$.

Пусть имеется молекулярное пространство G. Назовем множеством $T_k(G)$ k связок пространства G множество всех его связок на к точках. Обозначим $n_k=C_k(G)$ мощность $|T_k(G)|$ множества $T_k(G)$ (число различных к-связок), $n_k=C_k(G)=|T_k(G)|$.

З а м е ч а н и е .

Будем считать, что $n_0=C_0(G)=|T_0(G)|=1$ по определению.

Определение функционального вектора пространства.

Пусть G и $n_p$ есть конечное пространство и число его связок К(р) на р точках. Тогда набор всех $(n_1,n_2,....n_s)$ называется функциональным вектором пространства, $n_k=C_k(G)=|T_k(G)|$,

$$f(G)=(C_1(G),C_2(G),...C_s(G))=(n_1,n_2,....n_s).$$

Теорема 175.

*Пусть G и H - пространства и $f(G)=(n_1,n_2,.....n_k)$ и $f(H)=(m_1,m_2,....m_k)$ - их вектора. Тогда для прямой суммы G⊕H справедливы соотношения*

$$C_p(G\oplus H)=C_p(G)+C_p(H)+\sum_{s=1}^{p-1}C_s(G)C_{p-s}(H)=n_p+m_p+\sum_{s=1}^{p-1}n_s m_{p-s}=$$

$$\sum_{s=0}^{p}C_s(G)C_{p-s}(H)=\sum_{s=0}^{p}n_s m_{p-s},\quad n_0=m_0=C_0(G)=C_0(H)=1.$$

Д о к а з а т е л ь с т в о .

Число р связок в G⊕H равно сумме р связок в G, H, а также числу р связок, содержащих точки из G и H одновременно. Так как каждая точка из G соединена с каждой точкой из H, то число таких связок будет определяться входящей в формулу суммой. Теорема доказана. □

Напомним еще несколько теорем, часть из которых была доказана ранее.

Теорема 176.

*Пусть G - пространство и $G_1$ - его подпространство. Пусть $f(G)=(n_1,n_2,.....n_k)$ и $f(G_1)=(m_1,m_2,....m_k)$ - их вектора. Приклеим точку* *v к пространству G по подпространству $G_1$ (установим связи точки v со всеми точками подпространства $G_1$). Обозначим получившийся пространство G+v.*
*Тогда*

$$C_p(G+v)=C_p(G)+C_{p-1}(G_1),\ где\ p=1,2,..\ k+1.$$
$$f(G+v)=(n_1+1,\ n_2+m_1,\ n_3+m_2,...n_k+m_{k-1},\ m_k).$$

Доказательство.
Число точек пространства увеличивается на 1, а каждая p связка порождает (p+1) связку. Теорема доказана. □

Теорема 177.

*Пусть G - пространство, $v_1$ и $v_2$ - его точки и $G_1=O(v_1v_2)$ общий окаем этих точек, не являющихся смежными. Пусть $f(G)=(n_1,n_2,.....n_k)$ и $f(G_1)=(m_1,m_2,....m_k)$ - их вектора. Приклеим связь $v_1v_2$ между точками $v_1$ и $v_2$. Обозначим получившийся пространство $G+v_1v_2$.*
*Тогда*

$$C_p(G+v_1v_2)=C_p(G)+C_{p-2}(G_1),\ где\ p=1,2,..\ k+2.$$
$$f(G+v_1v_2)=(n_1,n_2+1,n_3+m_1,...n_k+m_{k-2},m_{k-1},m_k).$$

Доказательство
Число связей увеличивается на 1. Кроме того каждая p связка из $G_1$ порождает новую (p+2) связку. Теорема доказана. □

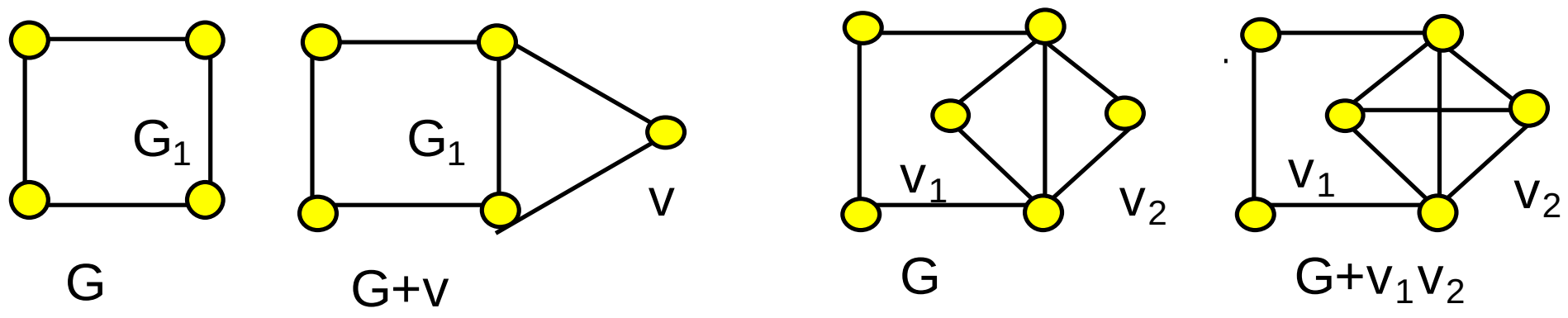

Рис. 201 Приклеивание точки и связи к пространству меняет его векторную функцию в соответствии с

На Рис. 201 слева изображено пространство G к которому приклеена точка v по $G_1$. Легко видеть, что f(G)=(4,4), $f(G_1)$=(2,1), и f(G+v)=(5,6,1). На этом же рисунке справа к точкам $v_1$ и $v_2$ приклеивается связь. В соответствии с формулами f(G)=(6,8,2), $f(G_1)$=(2,1), и $f(G+v_1v_2)$=(6,9,4,1).
Введем определение характеристической функции пространства G.

Определение характеристической функции пространства.

Пусть G есть пространство и $f(G)=(n_1,n_2,.....n_k)$ его вектор. Характеристической функцией F(G) на множестве молекулярных пространств называется функция, заданная на множестве молекулярных пространств и определяемая бесконечной последовательностью коэффициентов $F=(a_1,a_2,...a_k,...)$. значение которой на пространстве G определяется выражением

$$F(G) = \sum_{p=1}^{k} a_p n_p$$

где $a_p$, p=1,2,...к, являются коэффициентами.

Легко видеть, что характеристическая функция определяется, по сути дела, последовательностью коэффициентов $F=(a_1,a_2,...a_k,...)$. называемых вектором коэффициентов. Ее значение на пространстве G можно рассматривать как скалярное произведение векторов f(G) и F.

$$F(G)=\langle f(G)F\rangle$$

Определение степени n характеристической функции.

Пусть задана характеристическая функция, определяемая вектором коэффициентов $F=(a_1,a_2,...a_k,...)$. Степенью n этой функции называется характеристическая функция $F^n$, определяемая коэффициентами

$$F^n=(a_1, a_2,... a_k,...)^n=(a_n,a_{n+1},... a_k,...).$$

Иными словами, если пространство G имеет вектор $f(G)=(n_1,n_2,...n_k,...)$, то

$F(G)=F^1(G)=(a_1n_1+ a_2n_2 +... a_kn_k...)$,

$F^2(G)=(a_2n_1+a_3n_2+...a_{k+1}n_k...)$,

...............................................

$F^p(G)=(a_pn_1+a_{p+1}n_2+... a_{k+p-1}n_k...)$.

Например, для пространства G, состоящего из одной точки v,

$f(v)=(1,0,0,...0,...)$, получаем следующее: $F(v)=F^1(v)=a_1$, $F^2(v)=a_2$,...$F^p(v)=a_p$.

Рассмотрим теперь, как меняется характеристическая функция пространства при различных операциях над пространствами, такими, как приклеивания и отбрасывания точек и связей.

## *СВОЙСТВА ХАРАКТЕРИСТИЧЕСКИХ ФУНКЦИЙ МОЛЕКУЛЯРНОГО ПРОСТРАНСТВА*

Рассмотрим некоторые простые свойства характеристических функций пространства, которые будут использоваться в дальнейшем.

Теорема 178.

*Пусть G и H есть пространства, $f(G)=(n_1,n_2,.....n_s)$ и $f(H)=(m_1,m_2,....m_s)$ есть их вектора., и $F=(a_1,a_2,...a_k,...)$ есть некоторая характеристическая функция.*

*Тогда характеристическая функция на прямой сумме G⊕H определяется выражением*

$$F(G\oplus H)=\sum_{k=1}^{s} a_k n_k+\sum_{k=1}^{s} a_k m_k+\sum_{p=1}^{s} a_p\sum_{k=1}^{p-1} n_k m_{p-k}=F(G)+F(H)+\sum_{k=1}^{s} F^{k+1}(H)n_k.$$

Доказательство.

Пусть функциональный вектор пространства G определяется выражением $f(G)=(n_1,n_2,n_3,...n_s)$, функциональный вектор пространства H определяется выражением $f(H)=(m_1,m_2,m_3,...m_s)$. Как показано выше, функциональный вектор пространства G⊕H имеет вид:

$f(G\oplus H)=((n_1+m_1),(n_2+n_1m_1+m_2),(n_3+n_1m_2+n_2m_1+m_3),(n_4+n_1m_3+n_2m_2+n_3m_1+m_4),...(n_s+n_1m_{s-1},+...+n_{s-1}m_{s-1}+m_s),...(n_sm_s))$. Отсюда получаем, что

$F(G\oplus H)=(a_1(n_1+m_1)+a_2(n_2+n_1m_1+m_2))+a_3(n_3+n_1m_2+n_2m_1+m_3)+a_4(n_4+n_1m_3+n_2m_2+n_3m_1+m_4)+...+a_s(n_s+n_1m_{s-1},+...+n_{s-1}m_{s-1}+m_s)+...+a_{2s}(n_sm_s))=$
$(a_1n_1+a_2n_2+a_3n_3+...a_sn_s)+(a_1m_1+a_2m_2+a_3m_3+...a_sm_s)+n_1(a_2m_1+a_3m_2+...a_{s+1}m_s)+n_2(a_3m_1+a_4m_2+..a_{s+2}m_s)+...n_s((a_{s+1}m_1+a_{s+2}m_2+...a_{2s}m_s)=F(G)+F(H)+F^2(H)n_1+F^3(H)n_2+...+F^{s+1}(H)n_s$. Теорема доказана. □

Следствие.

Из доказательства теоремы следует, что для любых двух пространств G и H справедливо соотношение

$$\sum_{k=1}^{s} F^{k+1}(H)n_k=\sum_{k=1}^{s} F^{k+1}(G)m_k.$$

Теорема 179.

*Для любой характеристической функции справедливы следующие соотношения:*

*1. $F^{n+1}(G)=F^1(F^n(G))$.*

*2. $(F_1+F_2)^n=(F_1)^n+(F_2)^n$, или $(F_1(G)+F_2(H))^n=(F_1(G))^n+(F_2(H))^n$, где пространства G и H могут быть различны.*

Доказательство

Доказательство заключается в проверке. □

Теорема 180.

*Пусть G и H есть пространство и его подпространство. Пусть G+v есть пространство, полученное приклеиванием точки v к G по H. (Точка v соединяется связями со всеми точками в H). Тогда*

$$F(G+v)=F(G)+F^2(H)+F(v)=F(G)+F^2(H)+a_1.$$

Доказательство

Подставим вектор $f(G+v)=(n_1+1,\ n_2+m_1,\ n_3+m_2,...n_k+m_{k-1},\ m_k)$ пространства G+v в F(G+v). После преобразований получим искомое выражение. Теорема доказана. □

Теорема 181.

*Пусть G -пространство, $v_1$ и $v_2$ - две его несмежные точки и H - его подпространство, являющееся общим окаемом точек $v_1$ и $v_2$. Пусть $G+(v_1v_2)$ - пространство, полученное приклеиванием связи $v_1v_2$ к G. Тогда*

$$F(G+(v_1v_2))=F(G)+F^3(H)+F^2(v)=F(G)+F^3(H)+a_2.$$

Доказательство.

Теорема доказывается аналогично предыдущей. Подставим вектор $f(G+v_1v_2)=(n_1,n_2+1,n_3+m_1,...n_k+m_{k-2},m_{k-1},m_k)$ в $F(G+(v_1v_2))$. После несложных преобразований получим исходное выражение. □

## *БАЗИСНЫЕ ПРЕОБРАЗОВАНИЯ ПРОСТРАНСТВ*

Введем теперь определение преобразований пространства, которые обобщают точечные преобразования, и для которых мы найдем характеристическую функцию, являющуюся инвариантом этих преобразований. Следующие два определения должны рассматриваться вместе.

Определение семейства базисных пространств типа B.

Семейство T(B) пространств $B,G_1,G_2,G_3,....G_n,....$,

$T(B)=(B,G_1,G_2,G_3,....G_n,....)$, называется базисным типа B, если:

1. пространство B, называемое базой, принадлежит T(B).
2. Любое пространство семейства T(B) может быть получено из пространства B с помощью базисных преобразований типа B.

Все пространства семейства T(B) называются базисными.

Определение базисных преобразований типа B.

Следующие преобразования называются базисным типа B:

1. Удаление точки v.

Точка v пространства G может быть удалена, если окаем O(v) этой точки является базисным пространством типа B, то есть O(v) принадлежит T(B).

2. Приклеивание точки v.

Точка v может быть приклеена к пространству G по подпространству $G_1$, то есть соединена связями со всеми точками пространства $G_1$, если $G_1$ является базисным пространством типа B, то есть $G_1$ принадлежит T(B).

3. Удаление связи $(v_1v_2)$.

Связь $(v_1v_2)$, соединяющая две смежные точки $v_1$ и $v_2$ пространства G, может быть удалена, если общий окаем $O(v_1v_2)$ этих точек является базисным пространством типа B, то есть $O(v_1v_2)$ принадлежит T(B).

4. Приклеивание связи $(v_1v_2)$.

Две несмежные точки $v_1$ и $v_2$ пространства G могут быть соединены связью $(v_1v_2)$, если общий окаем $O(v_1v_2)$ этих точек является базисным пространством типа B, то есть $O(v_1v_2)$ принадлежит T(B). Семейство T всех базисных пространств может быть получено из

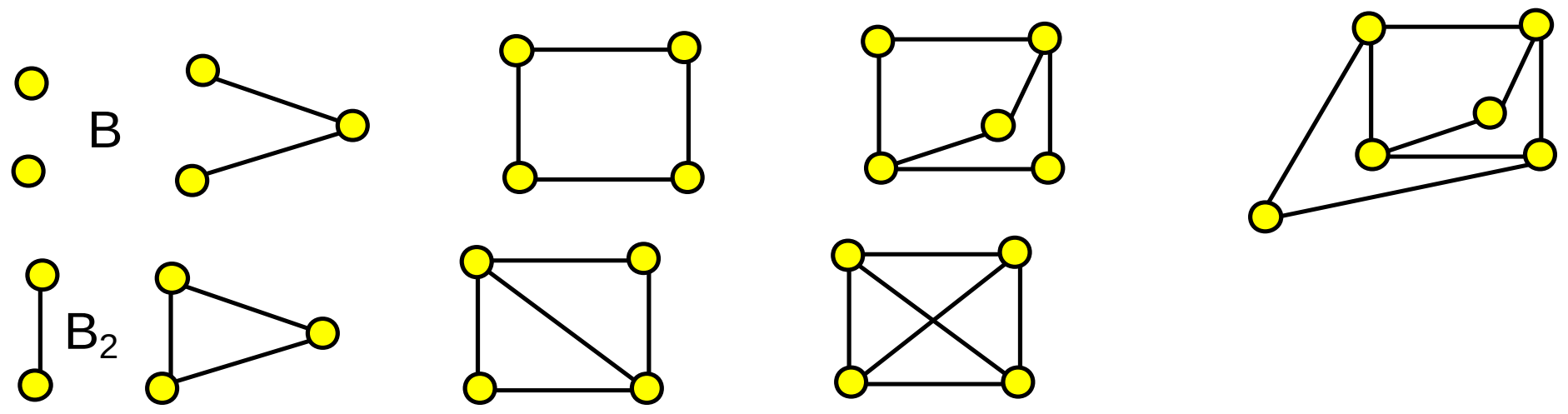

Рис. 202 Два семейства базисных пространств с базами B и $B_2$.

пространства B последовательным применением преобразований 1-4.

З а м е ч а н и е

Если в качестве базы B выбрана уединенная точка, то, очевидно, базисные преобразования являются обычными точечными преобразованиями.

На Рис. 202 изображены два семейства базисных пространств, если базовое пространство состоит из двух несмежных и двух смежных точек соответственно. Сейчас мы докажем основные теоремы этого параграфа, показывающую связь между любым базисным преобразованием и характеристической функцией, не меняющей своего значения на пространстве при базисных преобразованиях.

Определение характеристической функции базисных преобразований типа B.

Пусть база B имеет вектор $f(B)=(n_1,n_2,...n_k)$.

$a_1+a_2n_1+a_3n_2+...+a_{k+1}n_k=0$, или $F^2(B)=-a_1$.
$a_2+a_3n_1+a_4n_2+...+a_{k+2}n_k=0$ или $F^3(B)=-a_2$.
$a_3+a_4n_1+a_5n_2+...+a_{k+3}n_k=0$ или $F^4(B)=-a_3$.
..................................................
$a_{s-1}+a_sn_1+a_{s+1}n_2+...+a_{k+s-1}n_k=0$ или $F^s(B)=-a_{s-1}$.
..................................................

Теорема 182.

*Пусть $F=(a_1,a_2,...a_k,....)$ есть характеристическая функция базы B. Тогда на любом базисном пространстве G*
*$F(G)=F(B)$,*
*$F^{p+1}(G)=F^{p+1}(B)=-a_p$.*

Доказательство

Используя индукцию предположим, что на всех базисных молекулярных пространствах G с числом точек $|G|\le s$ и числом связей $P(G)\le t$ теорема справедлива. Рассмотрим некоторое базисное пространство G, имеющее s точек. Пусть H есть его базисное подпространство, где $|H|<s$. Приклеим к G точку v по H, то есть так, что $O(v)=H$. Тогда, согласно предшествующей теореме $F(G+v)=F(G)+F^2(H)+F(v)$. Так как $F(v)=a_1$ и, согласно предположению индукции, $F^2(H)=-a_1$, то $F(G+v)=F(G)$. Кроме того, $F^p(G+v)=F^p(G)+F^{p+1}(H)+F^p(v)$. Так как $F^p(v)=a_p$ и, согласно предположению индукции, $F^{p+1}(H)=-a_p$, то $F^p(G+v)=F^p(G)$. Таким образом для приклеивания точки теорема справедлива. Рассмотрим приклеивание связи в G. Пусть $v_1$ и $v_2$ - несмежные точки в G и $H=O(v_1v_2)$ - общий окаем этих точек, являющийся базисным пространством. Приклеим связь $v_1v_2$ между точками $v_1$ и $v_2$. Обозначим получившееся пространство $G+v_1v_2$. Из предыдущего $F(G+(v_1v_2))=F(G)+F^3(H)+F^2(v)$. Так как $F^2(v)=a_2$, $F^3(H)=-a_2$, то $F(G+(v_1v_2))=F(G)$. Рассмотрим $F^p(G+(v_1v_2))$. Мы ранее получили, что $F^p(G+(v_1v_2))=F^p(G)+F^{p+2}(H)+F^{p+1}(v)$. Так как $F^{p+1}(v)=a_{p+1}$, $F^{p+2}(H)=-a_{p+1}$, то $F^p(G+(v_1v_2))=F^p(G)$. Таким образом для приклеивания связи теорема справедлива. Аналогично доказывается справедливость теоремы при отбрасывании точек и связей. Теорема доказана. □

Докажем, теперь, что характеристическая функция сохраняется при применении базисных преобразований к произвольному пространству.

Теорема 183.

> *Пусть $F=(a_1,a_2,...a_k,....)$ есть характеристическая функция базы B. Тогда базисные преобразования типа B, примененные к произвольному пространству G не меняют значение характеристической функции на этом пространстве.*

Доказательство

Ход доказательства полностью аналогичен предыдущему. □

Иными словами, характеристическая функция базы B является инвариантом относительно базисных преобразований типа B.

Рассмотрим некоторые примеры.

1. Пусть базовое пространство B состоит из единственной точки v и определяется вектором f(B)=(1,0,...0,...). Формулы дают следующее:

$a_1 + a_2 = 0$

$a_2 + a_3 = 0$

$a_3 + a_4 = 0$

............

$a_{s-1} + a_s = 0$

Из этой системы легко получить, что $a_k=(-1)^{k+1}a_1$. При выборе $a_1=1$ получается эйлерова характеристика молекулярного пространства. При этом базисные преобразования с базой B=v являются точечными преобразованиями. Следовательно, характеристическая функция на пространстве G, определяемом вектором $f(G)=(n_1,n_2,...n_k)$, для этого вида преобразований имеет вид:

$$F = a_1(1,-1,1,-1,1,...(-1)^{k+1}...) \quad F(G) = a_1(n_1 - n_2 + n_3 - ...+(-1)^{k+1}n_k...).$$

1. Пусть базовое пространство $B_1$ (Рис. 202) состоит из двух уединенных точек, и определяется вектором f(B)=(2,0,...0,...). Формулы дают следующее:

$a_1 + a_2 2 = 0$

$a_2 + a_3 2 = 0$

$a_3 + a_4 2 = 0$

..............

$a_{s-1} + a_s 2 = 0$

Из этой системы легко получить, что $a_k=(-1)^{k-1}a_1/2^k$. Следовательно, характеристическая функция на пространстве G, определяемом вектором $f(G)=(n_1,n_2,...n_k)$, для этого вида преобразований имеет вид:

$$F = (a_1, -\frac{a_1}{2}, \frac{a_1}{4}, -\frac{a_1}{8}, ...\frac{(-1)^{k-1}a_1}{2^k}, ....) \quad F(G) = a_1(n_1 - \frac{n_2}{2} + \frac{n_3}{8} - ....+ \frac{(-1)^{k-1}n_k}{2^k} + ....)$$

Легко проверить, что при $a_1=1$ на всех пространствах, полученных из $B_1$ и изображенных в верхней части Рис. 202 характеристическая функция равна 2.

2. Пусть базовое пространство $B_2$ (Рис. 202) состоит из двух смежных точек, и определяется вектором $f(B_2)=(2,1,...0,...)$. Формулы дают следующее:

$a_1+ a_2 2+a_3=0$

$a_2+a_3 2+a_4=0$

$a_3+a_4 2+a_5=0$

....................

$a_{s-1}+a_s 2+a_{s+1}=0$

Из этой системы легко получить, что $a_k=(-1)^k((k-2)a_1+(k-1)a_2)$. Следовательно, характеристическая функция на пространстве G, определяемом вектором $f(G)=(n_1,n_2,...n_k)$, для этого вида преобразований имеет вид:

$$F=(a_1,\ a_2,\ -(a_1+2a_2),\ 2a_1+3a_2,...(-1)^k((k-2)a_1+(k-1)a_2),,...)$$

$$F(G)=(a_1n_1+a_2n_2-(a_1+2a_2)n_3+(2a_1+3a_2)n_4-...+(-1)^k((k-2)a_1+(k-1)a_2)n_k...).$$

Пусть $a_1=a_2=1$. Тогда

$$F(G)=n_1+n_2-3n_3+5n_4-7n_5+...+(-1)^k(2k-3)n_k+...$$

Легко убедиться проверкой, что на всех пространствах, полученных из $B_2$ и изображенных в нижней части Рис. 202 характеристическая функция равна 3. Выберем $a_1=1$, $a_2=-2$. Тогда характеристическая функция примет вид:

$$F(G)=n_1-2n_2+3n_3-4n_4+5n_5-...+(-1)^{k+1}kn_k+...$$

На всех пространствах, полученных из $B_2$ и изображенных в нижней части Рис. 202 характеристическая функция равна 0.

Пусть теперь $a_1=1$, $a_2=-1$. В этом случае характеристическая функция совпадает с эйлеровой характеристикой:

$$F(G)=n_1-n_2+n_3-n_4+n_5-...+(-1)^{k+1}n_k+...$$

На всех пространствах, полученных из $B_2$ и изображенных в нижней части Рис. 202 характеристическая функция равна 1.

Некоторые выводы из полученных результатов.

Мы рассмотрели несколько видов преобразований пространств и нашли функции, не меняющие своего значения при этих преобразованиях. Если основываться на точечных преобразованиях и на инвариантах этих преобразований, как на частном случае из класса базисных преобразований, то, по-видимому, можно ввести группы гомологий для любого базисного преобразования, которые будут инвариантами этого преобразования (аналогично группам гомологий на пространствах введенных ранее). Для этого необходимо определить разумным образом "симплекс", его границу, граничный оператор и другие стандартные математические инструменты, как это мы сделали при введении групп гомологий на пространстве. Кроме того, уже ясно просматриваются некоторые обобщения класса преобразований. Например, можно использовать не одно базовое пространство, а несколько, при сохранении всех остальных операций. Другой возможный путь для изучения-использовать пространства вместо точек, заменив операцию приклеивания точки на произведение пространств.

Список литературы к главе 17.

# МОЛЕКУЛЯРНЫЕ ПРОСТРАНСТВА В ФИЗИКЕ

Тесен мой мир. Он замкнулся в кольцо.
Вечность лишь изредка блещет зарницами.
Время порывисто дует в лицо.
Годы несутся огромными птицами.

Клочья тумана-вблизи...вдалеке...
Быстро текут очертанья.
Лампу Психеи несу я в руке-
Синее пламя познанья.

В безднах скрывается новое дно.
Формы и мысли сместились.
Все мы уж умерли где-то давно...
Все мы еще не родились.

М. Волошин.

We discuss introduction of a discrete background for the space instead of continuous one in Physics. Given the discrete normal 3-dimensional space we prove that (3+1) space-time is a regular 4-dimensional molecular space (not normal). From the discrete point of view the minimal number of points forming closed Universe is equal to eight. We analyze Ehrenfest's conclusion that physical space is three-dimensional based mainly on the solution of a parabolic equation and show that it is not enough to state with certainty that physical space is three-dimensional. We also discuss molecular spaces which have an external (a as the whole) and internal (point) topological structures. In conclusion we discuss several possible applications of molecular spaces to Physics.

Некоторые из полученных здесь результатов можно найти в работах [7,42,43].

## *НЕОБХОДИМОСТЬ ВВЕДЕНИЯ ДИСКРЕТНОГО ПРОСТРАНСТВА В ФИЗИКУ*

Идея введения дискретного пространства, как поля для развития событий, на котором строится физическая теория, начала занимать физиков одновременно с появлением квантовой механики [29,41]. Вначале это были попытки использовать множество точек с целочисленными координатами в непрерывном евклидовом пространстве $E^n$ как n-мерное дискретное пространство. При этом не учитывались топологические и внутренние различия между пространствами различных размерностей. В целом такие модели не принесли успеха, хотя в некоторых случаях, как например, в статистической физике такие модели работали удовлетворительно. Между тем неразрешимые проблемы, связанные с функциональным интегралом на планковской длине, указывают не только на провал классических полевых уравнений, но также и на то, что дифференциальное многообразие должно быть заменено на какую-либо конечную систему. Достаточно сослаться на слова Альберта Эйнштейна, который писал [56]: "Из квантовых явлений, по-видимому, следует, что конечная система с конечной энергией может полностью описываться конечным набором чисел. Это, кажется, нельзя совместить с теорией континуума... ." В другой своей работе-ответе на статью К. Менгера, который высказывал сомнения о справедливости представлений о физическом пространстве как о континууме Эйнштейн говорит следующее: "Я придерживаюсь представлений о континууме не потому, что исхожу из некоторого предрассудка, а потому, что не могу придумать ничего такого, что могло бы органически заменить эти представления. Каким образом следует сохранить наиболее существенные черты четырехмерности, если отказаться от этого представления?" [56]. По поводу применимости понятий римановой геометрии на малых областях пространства он писал: "Предложенная здесь физическая интерпретация геометрии не может быть непосредственно применена к областям пространства субмолекулярных размеров. Только успех может служить оправданием такой попытки приписать физическую реальность основным понятиям римановой геометрии вне области их физического определения. Однако, может оказаться, что подобная экстраполяция имеет не больше оснований, чем распространение понятия температуры на части тела молекулярных размеров" [56].

Среди современных физиков ряд ученых, также как Эйнштейн, считает, что необходимо ввести какую-то конечную структуру взамен дифференциального многообразия. В этом направлении работал Пенроуз, который изобрел спиновую дискретную сеть [36]. Финкельштейн [19] изучал квантовую микроструктуру пространства-времени, Айшем [23] предложил квантовую теорию на множестве всех топологий на данном множестве, и подробно рассмотрел случай конечного множества. Соркин [38] ввел один из путей перехода от локально конечного покрытия пространства к конечному частично-упорядоченному множеству. В своей известной [37] работе Редже предложил использовать симплициальный комплекс, как замену пространства-времени в общей теории относительности, что также можно считать дискретизацией. Его фундаментальными переменными являются множество длин ребер и матрица инцидентности. показывающая связь между ними. В одном подходе длины ребер могут меняться, а матрица инцидентности постоянна. В другом подходе, наоборот, длины ребер посто янны, однако связи в симплициальном комплексе могут меняться. Подход Редже, однако, предполагает, что действие развертывается на непрерывном пространства-времени, что не соответствует пожеланиям Эйнштейна. Кроме того, непонятно, каким образом минимальная (планковская) длина следует из подхода Редже.

Этих недостатков лишены молекулярные пространства, которые могут являться той базой, которая может быть использована физиками как дискретная основа для построения физических теорий.

Укажем на следующие особенности молекулярных пространств, которые делают их удобными для использования в физике:

1. Существование естественной метрики и внутренней геометрии.
2. Существование естественной минимальной длины $l_0$ (планковской длины), являющейся расстоянием между двумя соседними точками молекулярного пространства.
3. Локально-конечное число точек пространства.
4. Естественный переход к непрерывному случаю при увеличении масштабов.

## *ТЕОРЕМА О РЕГУЛЯРНОСТИ ЧЕТЫРЕХМЕРНОГО ПРОСТРАНСТВА-ВРЕМЕНИ*

Рассматривая физическое пространство, как молекулярное трехмерное нормальное, мы показали, что пространство-время теряет свойство нормальности и становится регулярным 4-мерным пространством-временем.

Теорема 184

*Четырехмерное пространство-время является всегда регулярным пространством.*

Доказательство.

Предположим, что пространство нормально и трехмерно (Рис. 203). Источник светового сигнала находится в точке a в момент t. Если происходит испускание светового сигнала, то в следующий момент t+1 световой сигнал может находиться в одной из ближайших к a точек. Процесс повторяется при переходе от слоя t+1 к t+2. Очевидно, что мы должны установить связи (ребра) между точкой a и ближайшими точками как в пространстве, так и во времени. На Рис. 203 изображены в качестве примера одномерное пространство и двумерное пространство-время. Легко видеть, что окаем точки a на слое t+1 является регулярной окружностью $S^1_2$. тогда как нормальной окружностью является $S^1_1$. Следовательно, двумерное пространство-время регулярно, так как не обладает свойством нормальности. Доказательство завершено. □

## *В МОМЕНТ РОЖДЕНИЯ ЗАМКНУТАЯ ВСЕЛЕННАЯ ИМЕЛА ВОСЕМЬ ТОЧЕК*

Один из результатов, следующий непосредственно из теории молекулярных пространств, состоит в том, что 3-мерная сфера (Рис. 203) не может иметь менее чем 8 точек. Как известно, расширяющаяся 3-мерная сфера является физической моделью расширяющейся вселенной. Следовательно, в начальный момент образования вселенной число различных точек пространства равнялось 8, а затем начало увеличиваться.

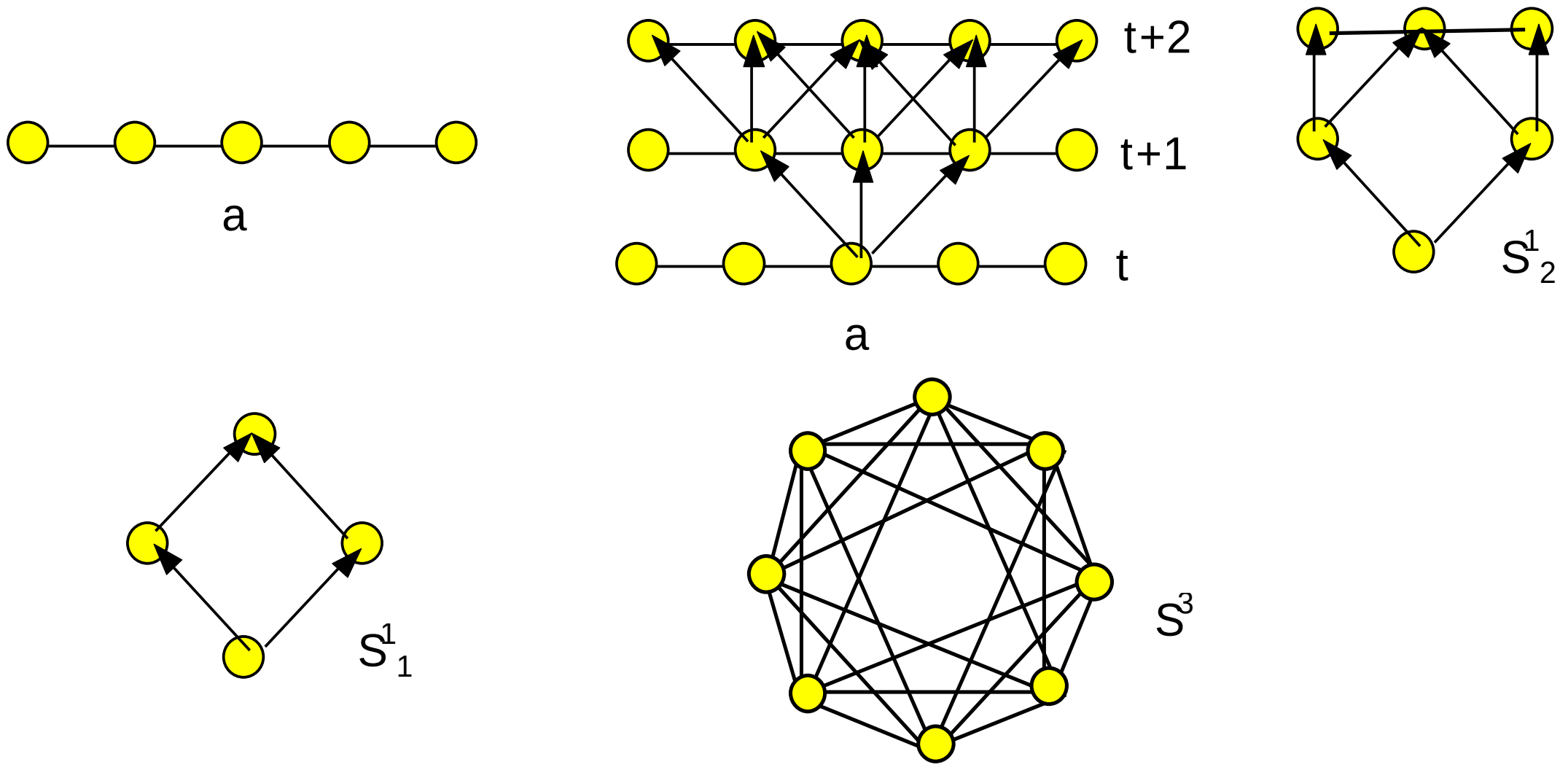

Рис. 203 Вверху. Световой сигнал испускается из точки a в момент t. В следующий момент (t+1) световой сигнал может находиться в одной из ближайших к a точек. Процесс повторяется при переходе от слоя t+1 к t+2. Легко видеть, что окаем точки a на слое t+1 является регулярной окружностью $S^1_2$. тогда как нормальной окружностью является $S^1_1$. Внизу. Минимальная трехмерная сфера, как модель ранней Вселенной, содержит 8 точек.

Возможно, что это связано каким-то образом с теорией элементарных частиц, в которой часто встречается число 8.

## *СУЩЕСТВУЮТ ЛИ ФИЗИЧЕСКИЕ ОПЫТЫ, ДОКАЗЫВАЮЩИЕ, ЧТО РЕАЛЬНОЕ ФИЗИЧЕСКОЕ ПРОСТРАНСТВО ТРЕХМЕРНО?*

Факт трехмерности физического пространства считается надежно обоснованным, по крайней мере в области макроскопических явлений. При этом во многих случаях теоретический анализ физических наблюдений и экспериментов, касающихся фундаментальных свойств пространства, основан на решениях дифференциальных уравнений. Использование молекулярных пространств как поля, на котором развивается физическое действие, расширяет спектр решений уравнений, и, следовательно, возможное толкование физических следствий, которые могут быть сделаны на основе этих решений. Как оказалось, на различных классах пространств дифференциальные уравнения могут давать сходные решения. Ранее мы уже рассматривали решения параболического уравнения на одномерном дереве, не являющемся многообразием и на молекулярном n-мерном евклидовом пространстве и обнаружили сходство в решениях. Это означает, что основываясь только на виде решений мы не можем делать надежных утверждений топологических свойствах пространства, таких как размерность или многообразие оно или нет, и, если реальное физическое пространство дискретно, то необходимо заново проанализировать экспериментальные подтверждения его трехмерности. Возможно, что физическое пространство, в котором мы существуем является не трехмерным, как это принято считать, а, например, одномерным.

Мы показали, [43] что в дискретной модели параболическое дифференциальное уравнение, описывающее физические процессы, дает на одномерном дереве решения, сходные с решениями в трехмерном пространстве. В этом случае мы не имеем реального физического инструмента, позволяющего определить истинную размерность и топологию физического пространства. Остановимся на этом подробнее. Один из аргументов в пользу того, что физическое пространство трехмерно, был предложен Эренфестом [2] в начале века.. Эренфест рассматривал устойчивость планетарных орбит в ньютоновской теории гравитации как экспериментальное доказательство трехмерности физического пространства. В непрерывном n-пространстве ньютоновский потенциал определяется как решение n-мерного уравнения Лапласа

$$\Delta\varphi = 0$$

и определяется выражением

$$\varphi = \frac{c}{r^{n-2}},$$

где n является размерностью пространства. При этом орбиты планет являются устойчивыми только при n=3. На основании этого был сделан вывод, что физическое пространство трехмерно.

В дискретном случае, как это было показано ранее, решение того же уравнения

$$\Delta\varphi = 0$$

дает сходимость на евклидовом молекулярном пространстве $E^n$ при $n \geq 3$, и на (одномерном) дереве $E^1_r$ при $r \geq 4$. Такое сходство в решениях на пространствах различных размерностей означает, что вид решения не зависит, или почти не зависит, от размерности, а определяется структурой связей между точками пространства. Перенося такое заключение на рассмотрение Эренфеста мы можем утверждать, что зависимость потенциала обратно пропорционально расстоянию до центра еще не означает заведомо трехмерности физического пространства. В случае дискретного физического пространства это может означать как трехмерность, так и одномерность. Поэтому необходимо ставить такие эксперименты, или проверять такие результаты, которые явно зависят от топологии физического пространства. Только в этом случае мы можем уверенно говорить о его размерности.

## *НЕСКОЛЬКО ВОЗМОЖНЫХ НАПРАВЛЕНИЙ ИСПОЛЬЗОВАНИЯ МОЛЕКУЛЯРНЫХ ПРОСТРАНСТВ В ФИЗИКЕ*

Использование дискретной модели пространства существенно расширяет как теоретические построения, так и экспериментальные перспективы для физики. Появляется возможность создавать модели вакуума, имеющие различную топологическую молекулярную структуру, которую затем можно проверить в эксперименте. Приведем несколько примеров такого рода пространств.

В предыдущих главах мы описали несколько молекулярных моделей евклидова n-мерного пространства. Именно такое непрерывное пространство в большинстве случаев является полем, в котором строятся физические теории. Использование молекулярных построений позволяет естественным образом ввести внутреннюю структуру в евклидово пространство, заменяя точки на некоторые пространства. При этом сохраняется топология плоского пространства, которую можно считать внешней, и возникает топологическая структура каждой точки, которую можно назвать внутренней. Такие пространства напоминают расслоенные пространства-конструкции, широко используемые как математиками так и физиками. Метод построения таких пространств достаточно прост. Вначале мы строим обычное молекулярное пространство со свойствами, которые мы отнесем к внешним свойствам пространства. Затем каждую точку заменяем каким либо пространством со свойствами, которые мы назовем внутренними. Стандартным примером такого построения является прямое произведение двух пространств, подробно изученное ранее. Рассмотрим еще один пример. Возьмем одномерную окружность, состоящую из 4 точек, а затем каждую точку заменим на двумерную сферу $S^2$ (Рис. 204). При этом сферы, замещающие соседние точки первоначальной окружности, образуют прямую сумму пространств. то есть каждая точка одной сферы соединена со всеми точками другой сферы. Полученное пространство будет иметь внешнюю топологию одномерного многообразия-окружности, в которой сферы играют роль точек.

Внутренняя топология каждой точки двумерной сферы определяется ее окаемом, который является не-многообразием и имеет вид пространства

$$O(v) = S^2{}_1 \oplus S^1{}_1 \cup S^2{}_2 \oplus S^1{}_1 .$$

Размерность каждой такой точки равна 4.

Таким образом внутренняя размерность и внешняя размерность различны.

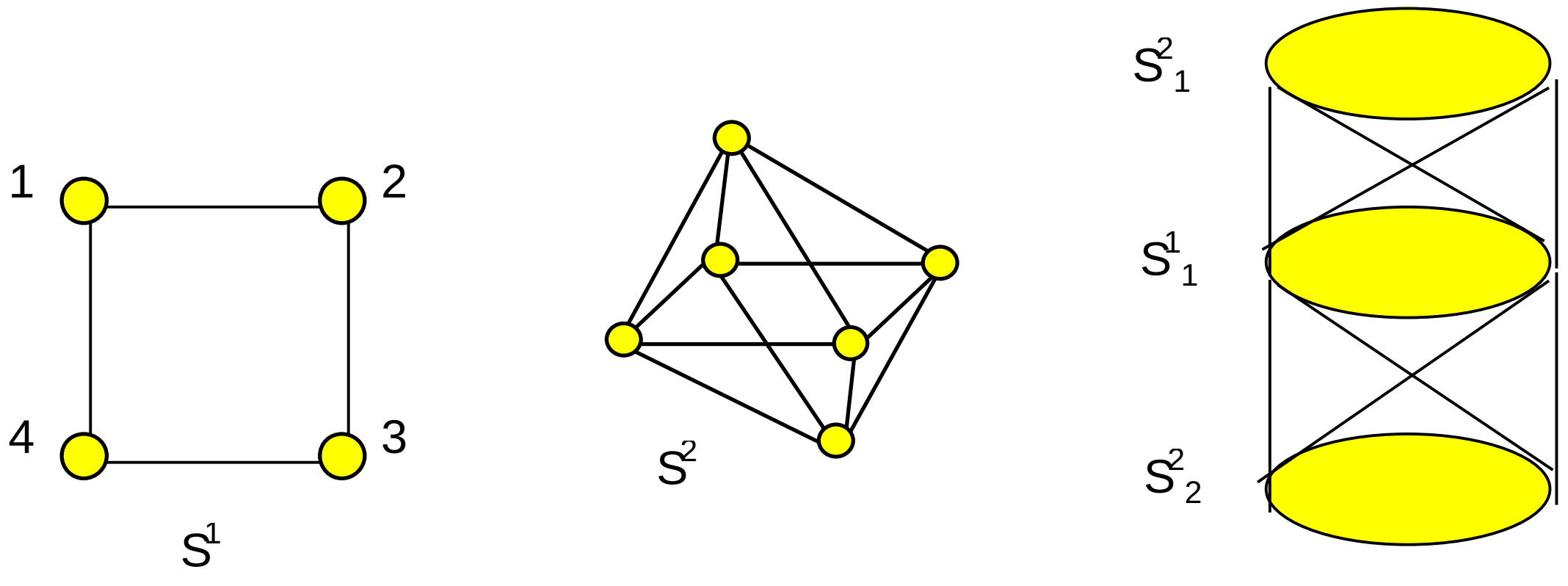

Рис. 204 Каждая точка одномерной окружности заменяется на двумерную сферу $S^2$. Сферы, замещающие соседние точки первоначальной окружности, образуют прямую сумму пространств. Полученное пространство будет иметь внешнюю размерность 1 одномерного многообразия-окружности, в которой сферы играют роль точек. Внутренняя размерность каждой точки двумерной сферы равна 4. $\dim(v)_{int}=4$, $\dim(v)_{ext}=1$.

$$\dim(v)_{int}= 4,\ \dim(v)_{ext}= 1.$$

Дадим точное определение такого вида пространств.

<u>Определение сцепки пространств G и H.</u>

Пусть имеется два пространства: G с набором точек ($v_1$, $v_2$, $v_3$,.... $v_m$,...) и набором связей W, и H с набором точек ($u_1$, $u_2$, $u_3$,.... $u_n$) и набором связей S. Тогда сцепкой двух пространств G и H называется пространство G•H=G(H), состоящее из точек $v_k \bullet u_p$, где окаем точки $v_k \bullet u_p$, определяется следующими условиями:

$v_k \bullet u_p$ смежна с $v_k \bullet u_s$, $s \neq p$, если $u_s$ принадлежит $O(u_p)$.

$v_k \bullet u_p$ смежна с $v_t \bullet u_s$, $t \neq k$, если $v_k$ принадлежит $O(v_t)$.

Сцепка является операцией, совмещающей свойства прямой суммы и прямого произведения пространств. Существенной особенностью является то, что эта операция некоммутативна и неассоциативна, G•H≠H•G, (G•H)•F≠H•G•(H•F). Кроме того, второе пространство H в G•H всегда должно быть конечного объема. Окаем точки $v_k \bullet u_p$ является, как легко видеть, прямой суммой двух пространств $O(v_k \bullet u_p)=(O(v_k)\bullet H)\oplus O(u_p)$.

Помимо упомянутого выше подхода, при котором пространство имеет внутреннюю и внешнюю топологию, можно ввести пространства ячеистого типа, в которых топологии различны в различных участках пространства. Примером такого пространства может служить трехмерное

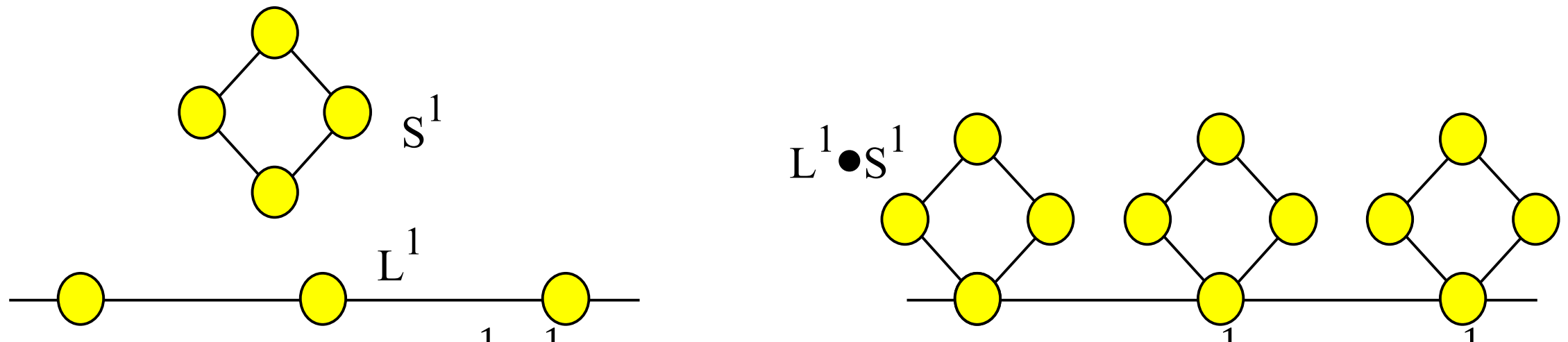

Рис. 205 Пространство $L^1 \bullet S^1$ является подвеской окружности $S^1$ к прямой $L^1$.

полное молекулярное евклидово пространство, в котором каждая точка $v_k$ заменяется на свое пространство $H_k$, а связи между точками устанавливаются как в сцепке двух пространств.

Опишем еще один подход к построению модели физического пространства, который может быть назван, условно, подвеской пространства H к пространству G.

Определение подвески H к G.

Пусть имеется два пространства: G с набором точек ($v_1$, $v_2$, $v_3$,.... $v_m$,...) и H с набором точек ($u_1$, $u_2$, $u_3$,.... $u_n$...). Выберем в H фиксированную точку $u_1$ и приклеим копию пространства H к каждой точке пространства G по точке $u_1$, то есть отождествим каждую точку пространства G с точкой $u_1$ в каждой копии пространства H. В полученном пространстве $W=G\bullet H$ точки можно обозначить как $w_{kp}=v_k\bullet u_p$. Окаем точки $v_k\bullet u_p$, определяется следующими условиями:

$v_k\bullet u_1$ смежна с $v_p\bullet u_1$, $k\neq p$, если $v_p$ принадлежит $O(v_k)$, и с точками $v_k\bullet u_s$, если $u_s$ принадлежит $O(u_1)$.

$v_k\bullet u_p$ смежна с $v_k\bullet u_s$, $p\neq s$, если $u_s$ принадлежит $O(u_p)$.

В W точки $w_{k1}=v_k\bullet u_1$ образуют подпространство, изоморфное G, точки $w_{1p}=v_1\bullet u_p$. $w_{2p}=v_2\bullet u_p$, $w_{3p}=v_3\bullet u_p$,... являются подпространствами, изоморфными H. Окаем точки $w_{k1}=v_k\bullet u_1$ есть, в свою очередь, пространство $O(v_k)\bullet O(u_1)$. Окаем точки $w_{kp}=v_k\bullet u_p$, $p\neq 1$, есть, в свою очередь, пространство $O(u_p)$. На Рис. 205 изображена подвеска окружности к прямой линии. Еще одна конструкция, позволяющая получить внешнюю и внутреннюю топологии в молекулярном пространстве, может быть построена следующим образом.

Определение вставки пространства H в пространство G.

Пусть имеется два пространства: G с набором точек ($v_1$, $v_2$, $v_3$,.... $v_m$,...) и набором связей W, и H с набором точек ($u_1$, $u_2$, $u_3$,.... $u_n$) и набором связей S. Тогда вставкой H в G пространство W=G•H=G(H), состоящее из точек $v_k$•$u_p$, где окрестность точки $v_k$•$u_p$, определяется следующими условиями:

$v_k$•$u_p$ смежна с $v_k$•$u_s$, $s \neq p$, если $u_s$ принадлежит $O(u_p)$.

$v_k$•$u_p$ смежна с $v_t$•$u_p$, $t \neq k$, если $v_k$ принадлежит $O(v_t)$.

На Рис. 206 показана вставка одномерной окружности $S^1$ в одномерную прямую $L^1$.

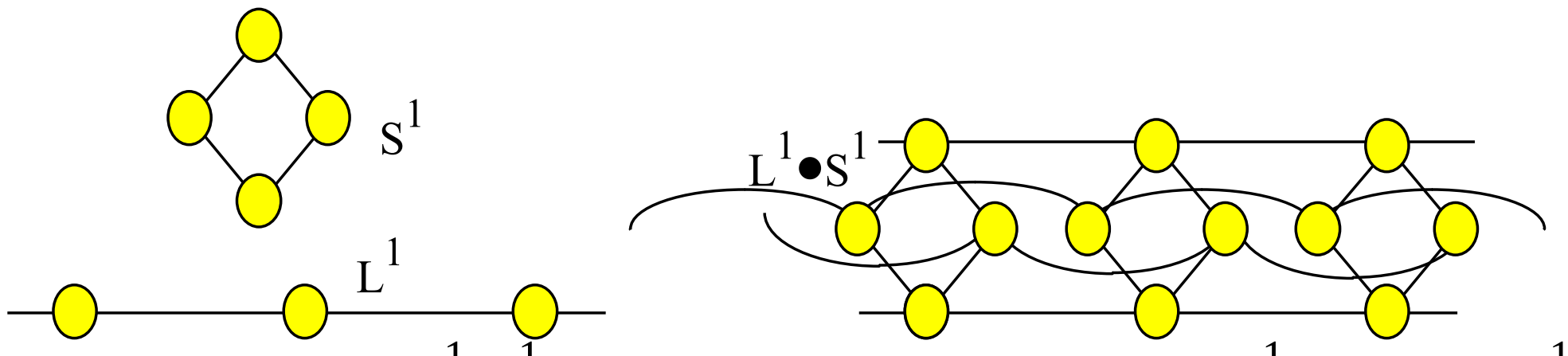

Рис. 206 Пространство $L^1$•$S^1$ является вставкой окружности $S^1$ в прямую $L^1$.

Можно предложить достаточно большое число подобных моделей, однако окончательный выбор должен быть сделан физиками и проверен в эксперименте.

Возможно, что теория молекулярных пространств может быть полезной при анализе геометрической структуры ядра. Ядро, в некотором приближении, можно считать построенным из тождественных частиц, по структуре образующих молекулярное пространство. Если это пространство известно, то не представляет особой сложности найти его характеристики как молекулярного пространства. Исследования в этом направлении являются новыми для физики и могут дать новые практические результаты. Применений молекулярных пространств в физике найдется везде, где физическую структуру можно описать графически при помощи точек и связей между ними, то есть в теории ядра, в теории элементарных частиц. в кристаллографии(см., например, [39]) и многих других областях.

Список литературы к главе 18.

# ПРИЛОЖЕНИЕ. МОЛЕКУЛЯРНОЕ ПРОСТРАНСТВО КАК ДИСКРЕТНАЯ ОСНОВА РЕАЛЬНОГО ФИЗИЧЕСКОГО ПРОСТРАНСТВА

## *ПРЕДПОСЫЛКИ ПОЯВЛЕНИЯ ТЕОРИИ МОЛЕКУЛЯРНЫХ ПРОСТРАНСТВ*

Никто не знает, какова природа реального пространства, в котором мы живем, в частности, непрерывно оно или дискретно. Более того, в современной науке до недавнего времени вообще не существовало концепции структуры пространства в малом, кроме стандартной непрерыной евклидовой модели. Представления о дискретных формах существования пространства имели характер общих рассуждений без всяких конкретных результатов и относились к разделам философии науки. Об этом, в частности, говорит тот факт, что в топологии и геометрии-разделах математики, изучающих различные пространства, не существует аппарата, пригодного для описания пространства, построенного на конечном или счетном множестве точек.

Многочисленные современные физические теории, провозглашаюшие дискретность пространства, на деле используют непрерывную пространственную и временную базу как основу для последующих построений. Между тем, идея дискретности пространства привлекала внимание как выдающихся мыслителей так и простых людей с незапамятных времен.

Дискретность в наиболее простой форме означает, что пространство образовано некоторыми неделимыми элементами, называемыми атомами пространства, между которыми нет ничего, в том числе и пространства. Такая модель напоминает атомную кристаллическую структуру твердого тела, за тем исключением, что между атомами пространства нет самого пространства. Именно отсутствие пространства является по меньшей мере странным для человека с нормальным воображением.

В связи с этим даже выдающиеся ученые совершали элементарные ошибки в трактовке дискретности пространства, в чем можно убедиться, раскрыв наугад почти любую из многих тысяч работ, затрагивающих тему дискретности. Для иллюстрации приведем слова выдающегося немецкого математика Г. Вейля, высссказывающего недоумение о гипотезе дискретности [1].

“Как следует понимать согласно этой идее существующие в пространстве отношения мер длин? Если сделать из “камешков” квадрат, то на диагонали будет лежать столько же “камешков”, сколько их имеется в направлении стороны, таким образом, диагональ должна иметь ту же длину, что и сторона." Вейль наивно применяет непрерывную меру к дискретному пространству, чего делать нельзя. Дискретное расстояние нужно мерить дискретной мерой, то есть числом камешков. С этой точки зрения диагональ действительно имеет ту же длину, что и сторона.

Следует признать, что до недавнего времени все немногочисленные попытки построить дискретное пространство как рабочий инструмент исследователя, не приносили желаемого результата. Между тем, дискретная модель пространства имеет имеет ряд очевидных преимуществ перед непрерывной моделью.

С математической точки зрения исчезают многочисленные особенности и парадоксы, связанные с несоизмеримостью, бесконечно-малыми и бесконечно-большими величинами и тому подобное. Например, исчезают пространства-монстры в топологии и отпадет надобность в иррациональных числах.

С эстетической точки зрения непрерывная модель пространства безлика и несодержательна. Вы уменьшаете размеры, делаете расстояния все меньше, но находитесь все в том же угле, образованном тремя взаимно-перпендикулярными плоскостями, и ничто вокруг вас не меняется. Трудно представить себе что-либо более унылое.

С практической точки зрения дискретная модель пространства более подходит нашему опыту, чем непрерывная модель. Квантовые явления доказывают, что поведение и взаимодействие материи меняются с уменьшением масштабов арены действия. Интуиция и здравый смысл подсказывают, что аналогично должно меняться и пространство. Непрерывная модель не позволяет это сделать, тогда как в дискретной модели изменение пространства происходит естественным образом.

Однако, горздо более интересны те перспективы, которые можно ожидать от применения модели дискретного пространства. Дискретная модель даст нам структуру пространства в малом, то есть покажет, как атомы пространства соединены друг с другом, чтобы образовать двумерное, трехмерное и, в общем случае, многомерное пространство. Знание такой структуры, а также способов изменения связей, позволит менять размерность пространства, создавая в трехмерном пространстве четырехмерные или пятимерные участки, а также совершать переходы между параллельными трехмерными пространствами.

Хотя дискретным пространством в современном понимании ученые начали интересоваться около ста лет назад, реальный путь к решению этой глобальной и в определенной степени мировоззренческой проблемы наметился лишь недавно в сязи с широким применением компьютером и связанной с этим необходимостью хранить в компьютерной памяти, а также преобразовывать многомерные пространственные объекты. Такие задачи возникают в медицине при компьютерном анализе результатов томографических исследований, компьютерной графике (computer graphics), компьютерном моделировании объектов (computer object modeling), диаграммном анализе (pattern analysis), научном представлении объектов в виде визуальных образов (scientific visualization) и многих других областях.

Новое научное направление, занимающееся этими проблемами получило название “дигитальная топология” (digital topology) [2]. Можно сказать, что требования практики породили новую науку.

На наш взгляд название “дигитальная топология” слишком узкое и не совсем верное в математическом плане. Во первых, термин топология не отражает геометрические аспекты указанных задач. Во вторых, дискретное пространство в строгом смысле не есть топологическое пространство и, следовательно, не является объектом топологии как математической науки [3]. В этом смысле более подошло бы название дигитальная геометрия.

В дигитальной топологии имеется несколько альтернативных подходов.

Один из подходов, которому и посвящена предлагаемая читателям статья, называется теорией молекулярных пространств-ТМП. В рамках ТМП строятся дискретные многомерные евклидовы и кривые пространства [4,5], изучаются их деформации, сохраняющие и меняющие пространственные инварианты. Становится ясно, например, как нужно изменить связи между точками, чтобы перейти от трехмерного пространства к четырехмерному, как по непрерывному пространству построить дискретное и наоборот.

## *ГЕОМЕТРИЧЕСКАЯ И ИНТУИТИВНАЯ ОСНОВА ДЛЯ ОПРЕДЕЛЕНИЯ РАЗМЕРНОСТИ НА МОЛЕКУЛЯРНОМ ПРОСТРАНСТВЕ*

Наша задача - построить пространство, состоящее из конечного или счетного числа точек, и обладающее такими свойствами, которые позволят рассматривать его как комбинаторный образ непрерывного n-мерного пространства. Для этого мы выберем в элементарной геометрии те признаки, которые описывают понятия близости, непрерывности и размерности на наиболее фундаментальном уровне, и которые позволят нам в дальнейшем легко перейти к молекулярной модели непрерывного пространства. Картина, которую мы собираемся здесь предложить читателю позволяет понять на интуитивном уровне, что является размерностью в теории молекулярных пространств.

Рассмотрим теперь 3-мерное евклидово пространство $E^3$ [6], в котором выберем некоторую точку p. Пусть $D^3$ будет сплошным шаром радиуса R с центром в точке p. Границей этого шара является, естественно, двумерная сфера $S^2$. Как известно, множество всех таких шаров со всевозможными центрами и радиусами образует базу топологии для $E^3$. Теперь предположим, что пространство молекулярно. Так как в молекулярном пространстве любой конечный объем содержит конечное число элементов, следовательно, шар $D^3$ состоит из конечного числа точек. Для нас здесь важно следующее: поверхность молекулярного шара является молекулярной сферой, которую можно ободрать как внешний слой кожуры с луковицы и получить молекулярный шар меньшего размера. Продолжая этот процесс обдирания, мы получим наименьший молекулярный шар, состоящий из последнего слоя, который еще можно ободрать, и неделимого элемента внутри этого слоя, играющего роль центральной точки. Таким образом, мы можем определить наименьший молекулярный

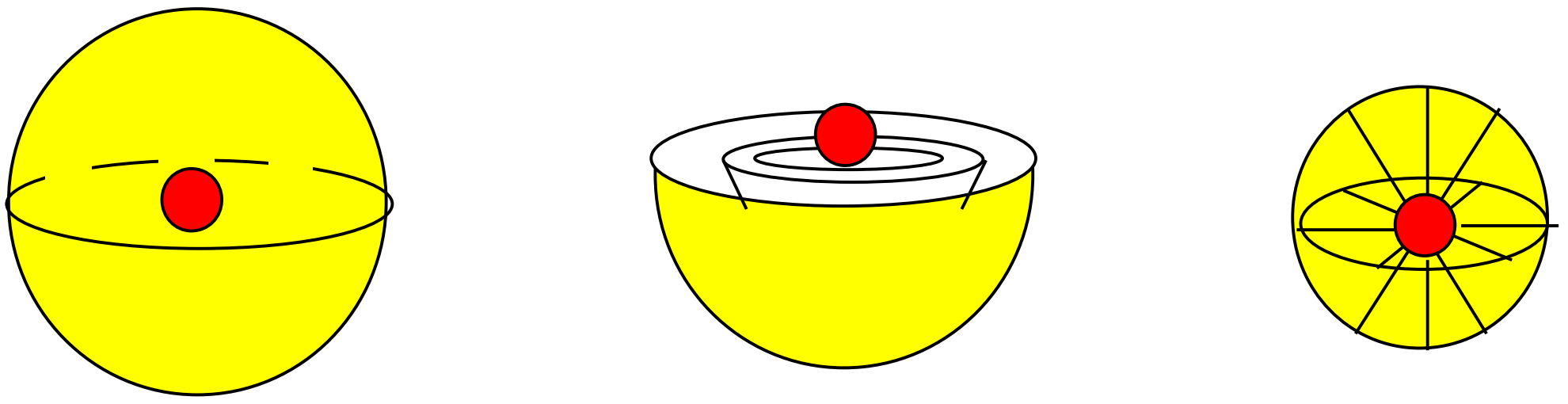

Рис. 207 В дискретном шаре обдирание слоев приводит к наименьшему шару шар как центральную точку, окруженную молекулярной сферой. Точки сферы являются ближайшими к центральной точке (Рис. 207). Для выделения этой близости установим, что центральная точка соединена ребром с каждой точкой ближайшей молекулярной сферы. Наглядной иллюстрацией этого случая является матрешка, где внутри деревянной фигурки находится вторая меньшего размера, внутри которой находится третья и так далее, до последней сплошной куклы, внутри которой уже ничего нет. Можно придумать и растительные образы, например, рассматривать молекулярный шар как кочан капусты или луковицу. Мы уже почти получили определение: молекулярный минимальный трехмерный шар является точкой, окруженной молекулярной двумерной сферой, или, в общем случае, молекулярный минимальный n-мерный шар является точкой, окруженной молекулярной (n-1) - мерной сферой.

## *РАЗМЕРНОСТЬ НОРМАЛЬНОГО МОЛЕКУЛЯРНОГО ПРОСТРАНСТВА. ОПРЕДЕЛЕНИЯ И ПРИМЕРЫ*

Перейдем теперь к определению [5,6] размерности молекулярного пространства, используя вышеприведенный геометрический пример.

В общем виде молекулярное пространство образуется набором точек, являющихся своего рода атомами пространства, наподобие того, как атомы образуют молекулу. Минимальное расстояние между двумя ближайшими точками постоянно и равно L. Расстояние между двумя различными точками всегда кратно L, то есть равно nL, n = 1, 2, 3, …. Из данной точки пространства можно перейти только в ближайшую точку. Если две точки являются ближайшими по отношению друг к другу, то они соединены ребром. Множество всех точек, ближайших к данной, назовем окаемом данной точки. Такая конструкция, состоящая из точек, некоторые из которых соединены линиями, (не обязательно прямыми) называется в математике графом.

Согласно классическому определению, нуль-мерная сфера состоит из двух изолированных точек. Придерживаясь этого определения назовем нормальной нуль-мерной сферой $S^0$ молекулярное пространство, состоящее из двух изолированных точек (атомов пространства) (Рис. 208). Тогда одномерный шар образуется центральной точкой, соединенной с

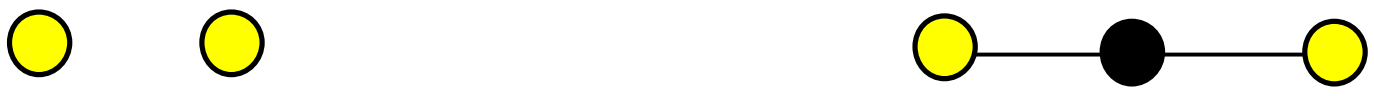

Нульмерная сфера $S^O$

Одномерная точка $v_l$ как центр одномерного шара

Рис. 208 Путь построения одномерной точки обеими точками из $S^0$, при этом центральная точка будет называться одномерной точкой. Замкнутым n-мерным пространством назовем молекулярное пространство, все точки которого n-мерны. Например, в одномерном случае замкнутым пространством будет окружность, состоящая из четырех и более точек. В свою очередь назовем точку n-мерной, если ее окаем является (n-1)-мерным замкнутым пространством. Рис. 209.

Легко убедиться в том, что на Рис. 209 изображены замкнутые одномерные пространства-окружности, а на Рис. 210 двумерные шары, состоящие из центральных точкек, соединенных с окружностями.

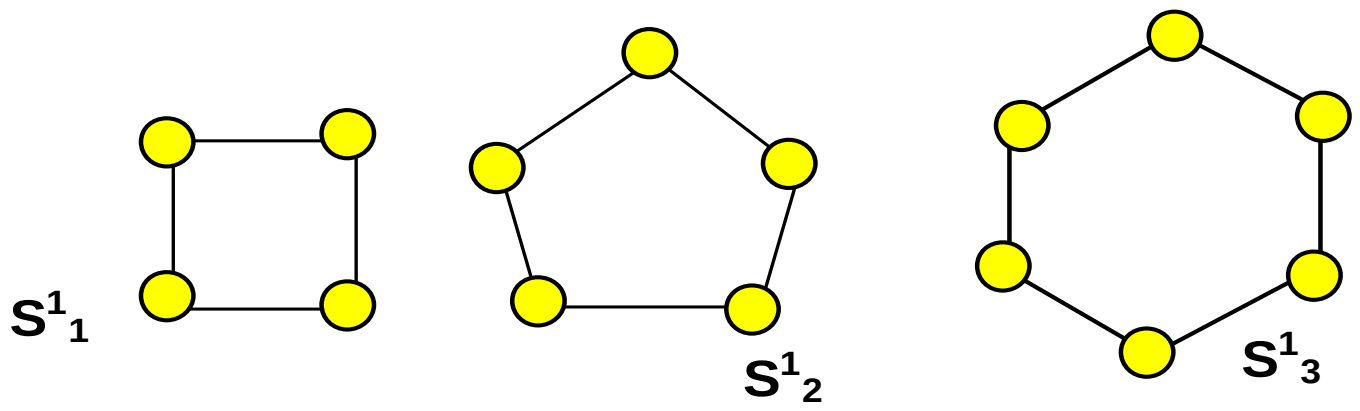

Рис. 209 Одномерные сферы (окружности). Сфера $S^1_1$ является наименьшей возможной.

На Рис. 211 изображены две двумерные молекулярные сферы и проективная плоскость, в непрерывном случае получаемая склеиванием крест-накрест противоположных сторон бумажного листа. Следует подчеркнуть, что проективная плоскость в непрерывном случае не может быть реализована в трехмерном пространстве, для этого необходимо как минимум пространство четырех измерений. Каждая из точек этих двумерных пространств является двумерной точкой одного из типов, изображенных на Рис. 210.

На Рис. 213 показана трехмерная точка, являющаяся центром трехмерного минимального шара, и трехмерная сфера, состоящую из восьми трехмерных точек вида $v_1$. Это минимальное количество точек, необходимое для формирования трехмерной сферы. Как известно, в непрерывном случае трехмерная сфера может быть построена только в

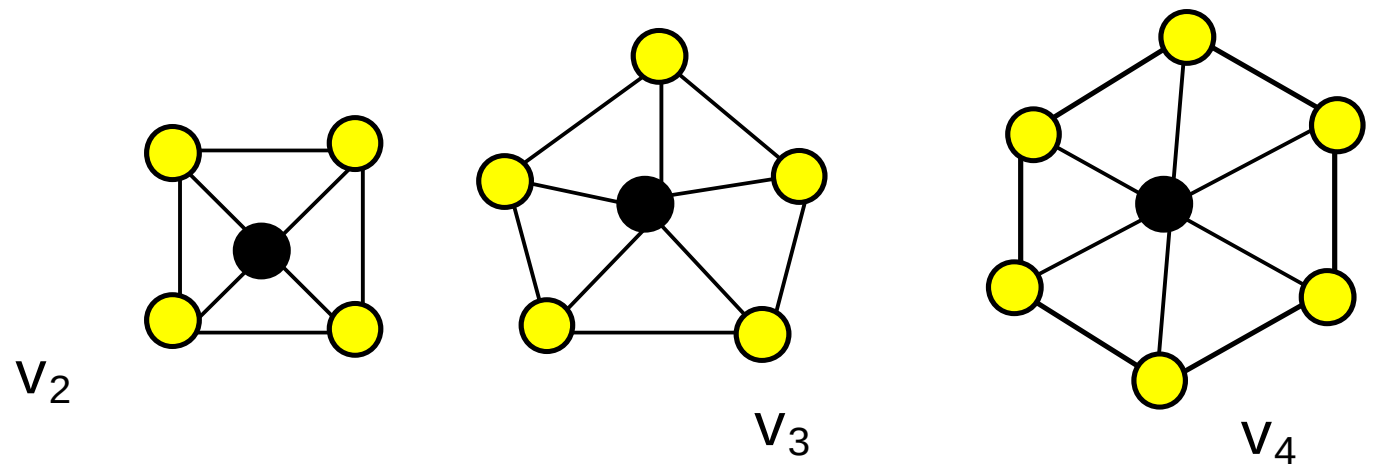

Рис. 210 Двумерные точки, каждая из которых окружена окружностью.

четырехмерном пространстве. Хотя длины ребер на рисунках различны, это сделано для наглядности, на самом деле все ребра имеют одну и ту же длину L. Не представляет сложности построить замкнутое молекулярное пространство любого числа измерений [5,6,7].

Несколько вариантов двумерных молекулярных евклидовых пространств представлены на Рис. 212.

В теории молекулярных пространств предлагаются различные методы моделирования молекулярных пространств по непрерывным и, наоборот, непрерывных пространств по молекулярным [7,8]. Для сравнения с непрерывными пространствами были определены и исследованы некоторые геометрические и топологические характеристики и инварианты молекулярных пространств, аналогичные соответствующим

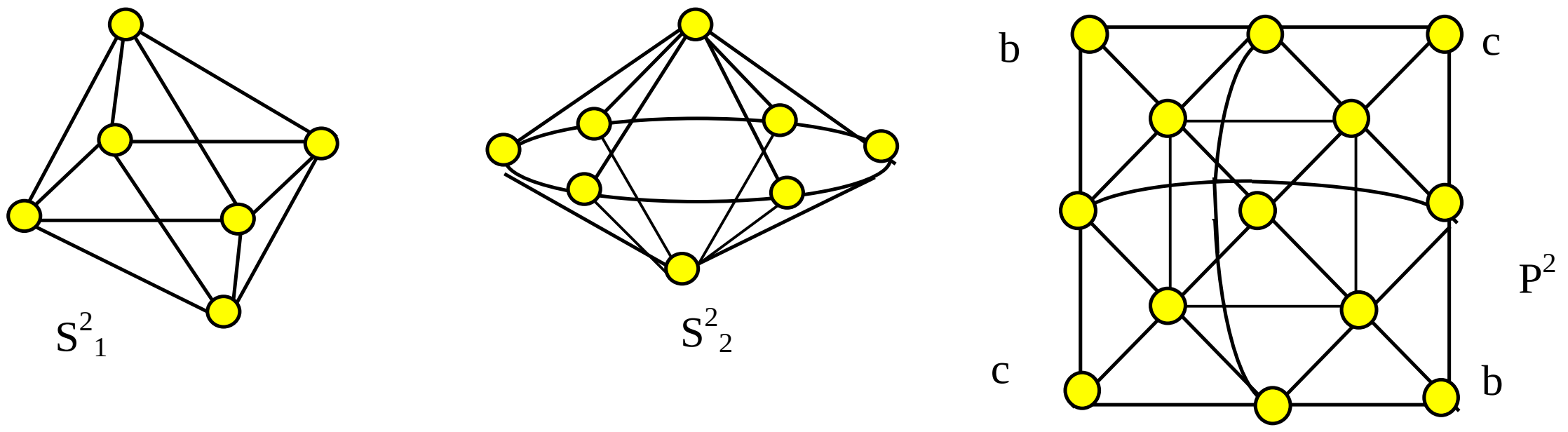

Рис. 211 Две двумерные сферы и проективная плоскость (склеенная сфера).

инвариантам на непрерывных пространствах. В пользу молекулярных пространств говорит тот факт, такие инварианты как эйлерова характеристика и группы гомологий совпадают на непрерывных пространствах и их молекулярных образах. Интересные результаты получены при изучении структуры некоторых дифференциальных уравнений в частных производных на замкнутых и открытых пространствах. Более полное представление об основных направлениях ТМП можно получить, ознакомившись с цитируемой литературой.

## *КАКОВА ТОЛЩИНА ПЛОСКОСТИ?*

Существует несколько способов представить молекулярные пространства алгебраически и геометрически. Алгебраически любое молекулярное пространство может быть описано квадратной матрицей инциденций или прямоугольной матрицей координат. Геометрически молекулярное пространство моделируется наиболее удобным образом как множество кирпичей. Кирпич (kirpich) является бесконечномерным единичным кубом с единичными координатами вершин в бесконечномерном евклидовом пространстве [9]. Два кирпича имеют общую грань, если две соответствующие им точки молекулярного пространства соединены ребром На Рис. 214 изображено построенное из кирпичей молекулярное двумерное евклидово пространство $E^2_2$ (Рис. 212).

Среди читателей, особенно математиков, может возникнуть вопрос: зачем вообще необходимо введение каких-то бесконечно-мерных кирпичей,

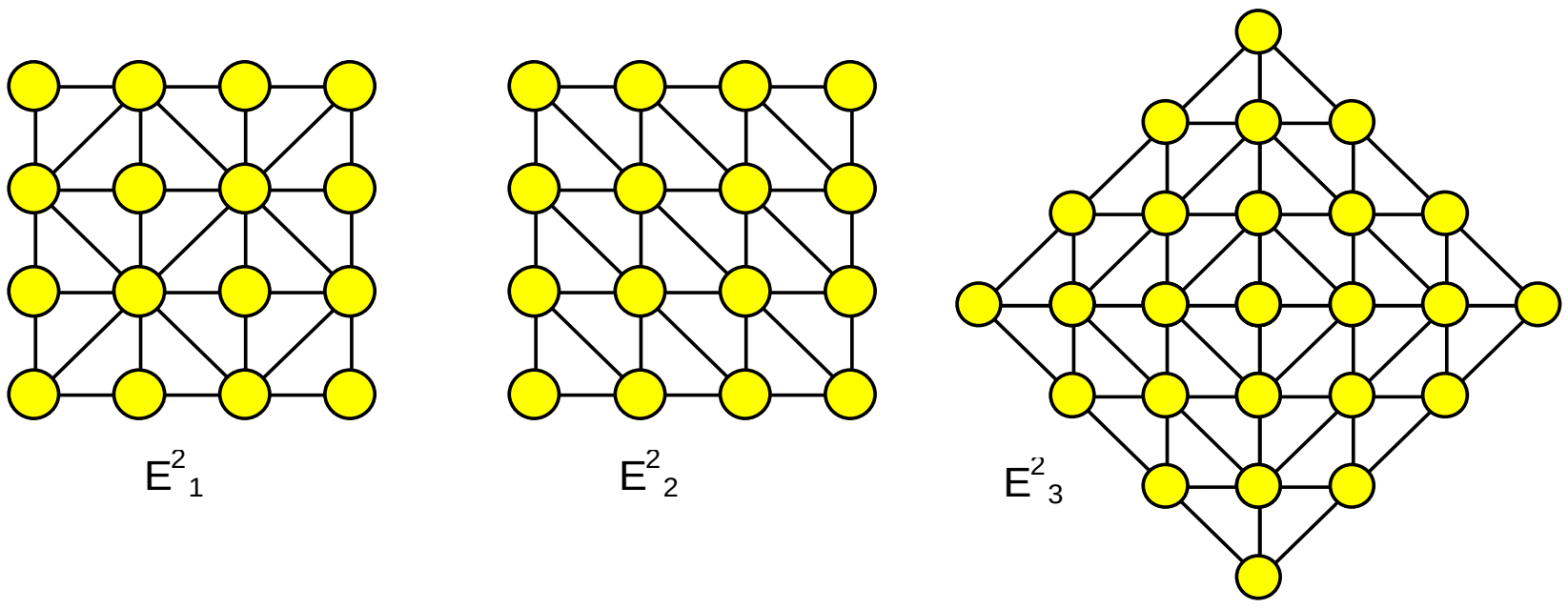

Рис. 212 Различные представления плоского двумерного пространства.

которые, в принципе, представляют просто некоторое изображение вершин молекулярного пространства.

Во первых использование кирпичей является удобным способом получения молекулярных пространств как образов непрерывных пространств, и, наоборот, получения непрерывных моделей молекулярных пространств.

Во вторых, использование кирпичей позволяет ввести понятия классической геометрии в ТМП и рассматривать молекулярное пространство как вложение в непрерывное евклидово пространство с вполне определенными внешними характеристиками. Вообще говоря, непрерывность более соответствует обыденному восприятию окружающего нас мира и использование кирпичей, как элементов деления физического пространства на малые элементы устраняет псмхологическое противоречие между ощущением непрерывности и реальностью дискретности.

Мы здесь изложим нашу концепцию, которая может быть названа физической, или экспериментальной (о компьютерных экспериментах рассказано в [7]), и которая есть не что иное, как объяснение подхода к дискретным пространствам.

Прежде всего, давайте рассмотрим классический математический подход к описанию двумерной поверхности. Математики считают двумерной поверхностью, например, поверхность шара, толщина и объем которой равны нулю. Двумерной поверхностью также является кусок плоскости, толщина плоскости равна нулю, объем куска плоскости также равен нулю. Иными словами, толщина (и объем) двумерной поверхности в математике равны нулю.

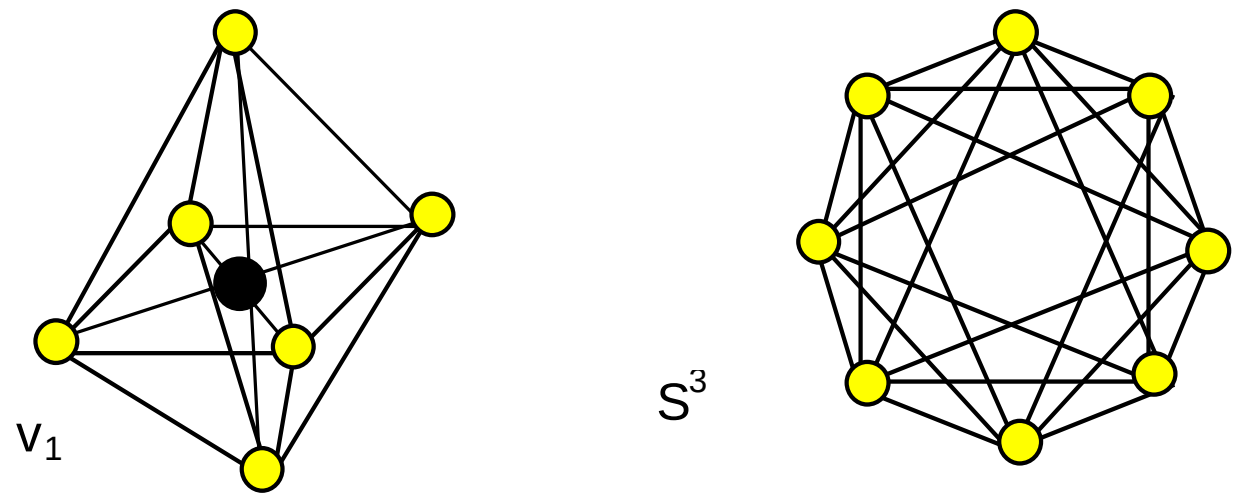

Рис. 213 Трехмерная точка и трехмерная сфера.

Однако, возможна и другая точка зрения.

Еще со школьных времен мы знаем, что все тела состоят из молекул. Поверхностью любого материального шара будет внешний слой молекул.

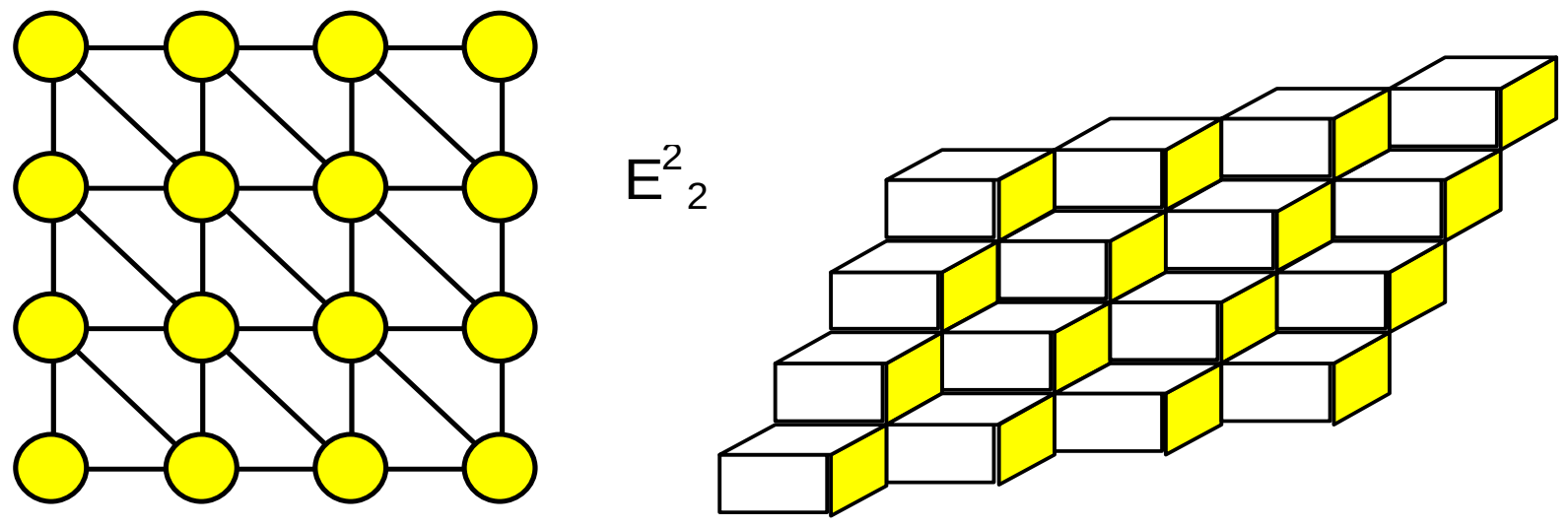

Рис. 214 Молекулярная однородная модель плоскости и ее представление в виде единичных кубов.

Следовательно, эта поверхность имеет ненулевую толщину. Физическая толщина этой поверхности равна одной молекуле. В одной из школьных задач по физике предлагается найти наибольшую площадь нефтяной пленки на поверхности воды. При этом подразумевается, что наименьшая толщина пленки равна также одной молекуле, но не равна нулю. Иными словами, двумерных поверхностей нулевой толщины в природе не существует. Наш подход как раз и отражает эту закономерность и вводит ее в математический и прикладной обиход. Мы можем считать, что с точностью до ненулевого постоянного коэффициента пропорциональности, минимальная толщина двумерной поверхности равна единице. С другой стороны, какова наибольшая толщина, при которой данная поверхность еще остается двумерной? Ответ на этот вопрос также достаточно прост и ясен в рамках данного подхода. Чтобы образовалась трехмерное пространство необходимо, чтобы, по крайней мере, некоторые его точки были трехмерными. Возьмем несколько слоев двумерной плоскости, образованной кирпичами, и будем поочередно укладывать их один на другой, как это показано на Рис. 215. Выясним, сколько слоев нужно уложить, чтобы получить трехмерные точки.

Легко видеть, что пространство в два слоя будет все еще двумерным, так как каждая точка имеет в своем окаеме одномерную окружность, но не имеет двумерной сферы. Однако, три слоя, изображенные на Рис. 215, уже будут трехмерным пространством, так как каждая точка промежуточного слоя имеет окаем, содержащий двумерную сферу. Таким образом, двумерная поверхность может иметь толщину один или два и состоять из одного или двух слоев. Три слоя образуют уже трехмерное пространство. Мы можем считать, что с точностью до ненулевого постоянного коэффициента пропорциональности толщина двумерной поверхности равна единице или двум. Минимальная толщина трехмерного (и любого n-

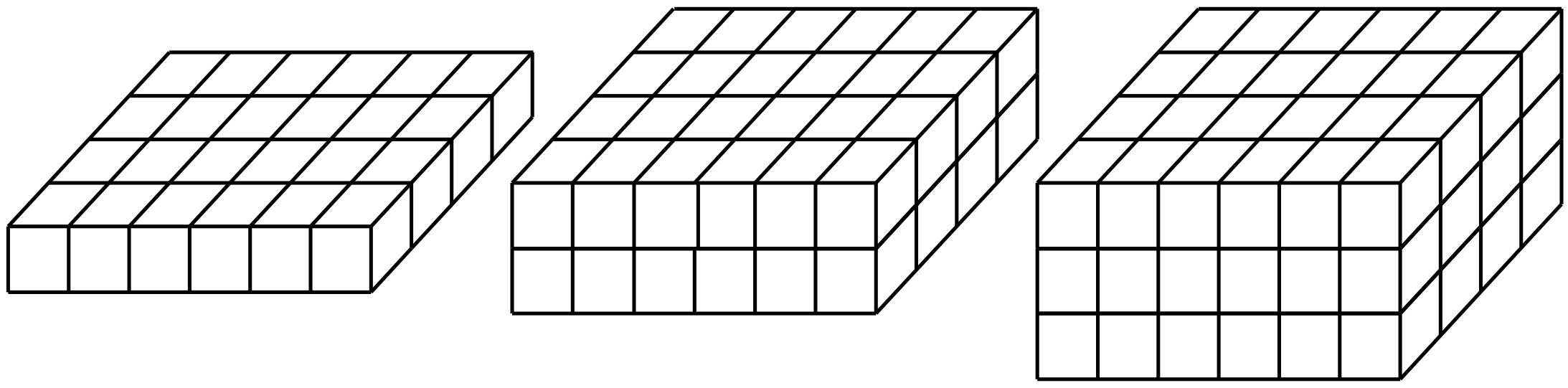

Рис. 215 Один или два слоя кирпичей являются двумерным пространством. Три слоя кирпичей уже образуют трехмерное пространство.

мерного) слоя равна трем, причем средний слой состоит из трехмерных (n-мерных) точек, а два внешних слоя содержат двумерные ((n-1)-мерные) точки. Следовательно, можно сформировать n-мерный слой как три слоя (n-1)-мерного пространства, склеенных аналогично предыдущему.

Точно также линия может иметь ширину 1 или 2, но если линия имеет ширину 3, то это уже часть плоскости.

## *ПЕРЕСТРОЙКА ПРОСТРАНСТВА С ИЗМЕНЕНИЕМ РАЗМЕРНОСТИ*

Одним из возможных применений теории молекулярных пространств является перестройка связей между атомами пространства с целью изменения его размерности. Действительно, как это видно из предыдущего, размерность пространства определяется тем, что является окружением данного атома пространства, в то время как сами атомы не имеют размерности. Одни и те же атомы можно расположить так, что они будут образовывать пространство любого числа измерений. Принципиально не представляет трудностей перейти от одномерного, например, пространства к трехмерному или четырехмерному.

В качестве примера рассмотрим перестройку двумерного шара в трехмерный. Двумерный шар представляется точкой $v_4$, окруженной окружностью, состоящей из шести точек (Рис. 216 слева). Для перехода к трехмерному шару установим последовательно связи в двумерном шаре как это показано на Рис. 216. Правое молекулярное пространство, изображенное на этом рисунке справа, является трехмерным шаром. Возникает вопрос, если реальное физическое пространство дискретно, то как можно осуществить подобную перестройку, создав четырехмерный шар из трехмерного? Ответ на этот вопрос и практическое реализация подобного процесса в первую очередь дело таких наук как физика. Кроме того, результат является, по крайней мере сейчас, непредсказуемым по своим последствиям. Ядерное энергия и связанные с ней военные и экологические катастрофы показали опасность использования научных открытий, даже когда последствия можно предвидеть. Последствия же

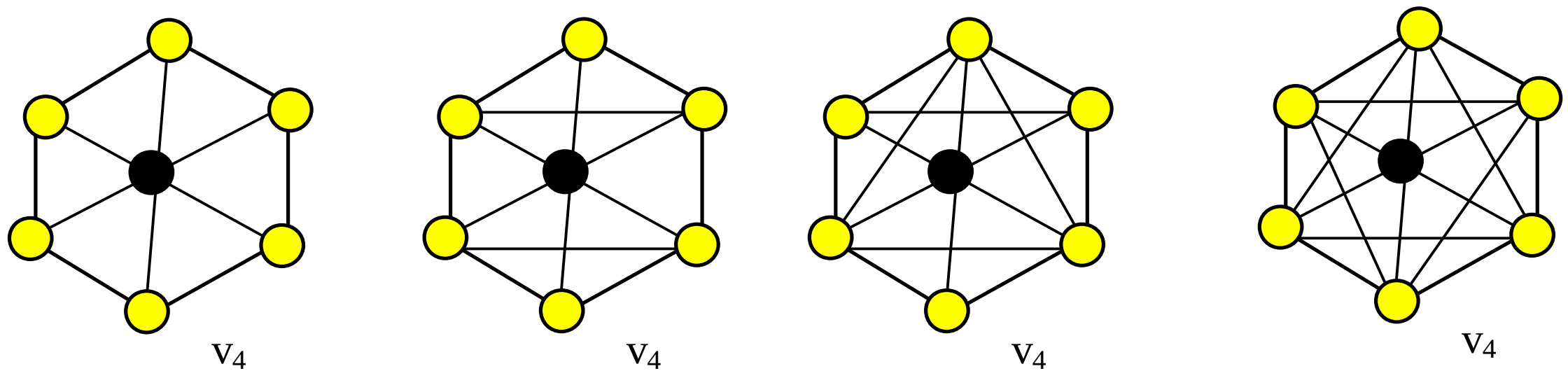

Рис. 216 Перестройка пространства с увеличением его размерности. Слева изображена двумерная точка, справа находится трехмерная точка.

экспериментов с пространством предвидеть в данный момент невозможно. С другой стороны, наука не может остановиться в своем развитии. Поэтому, прежде чем переходить к экспериментам с пространством, необходимо на уровне теории и наблюдений понять, чем нам это может грозить.

## *ДВИЖЕНИЕ В МОЛЕКУЛЯРНОМ ПРОСТРАНСТВЕ*

Представляет определенный интерес рассмотреть, как моделируется движение в молекулярном пространстве.

Следует заранее подчеркнуть, что любое перемещение в реальном пространстве является перемещением из одного слоя престранства в другой, параллельный первому. Рассмотрим, как это осуществляется в двумерном пространстве при переходе из одного одномерного пространства в другое, параллельное ему. Процесс перехода достаточно прост, как это следует из Рис. 217. Предположим, что одномерный стержень, изображенный черным цветом, находится на слое 1 двумерного пространства. Для того, чтобы перевести его на параллельный слой 2, необходимо поочередно сдвинуть все его точки. Это можно сделать, начав, например, с его левого конца. Справа на этом рисунке изображен момент, когда две его левые точки находятся уже на слое 2, а три правые все еще остаются на слое 1. При этом на слое 2 стержень появляется как бы из ничего, а на слое 1 исчезает как бы никуда.

Конечно, это всего лишь теоретическое описание перехода. Однако, для осуществления реального перехода необходимо уметь делать лишь одно: переводить точку стержня из одной точки пространства в другую. Такой процесс осуществляется постоянно при любом движении объектов в реальном трехмерном пространстве. Если предположить, что реальный физический мир есть трехмерный слой в четырехмерном промтранстве, то для перехода в параллельный слой достаточно лишь слегка изменить направление движения объекта.

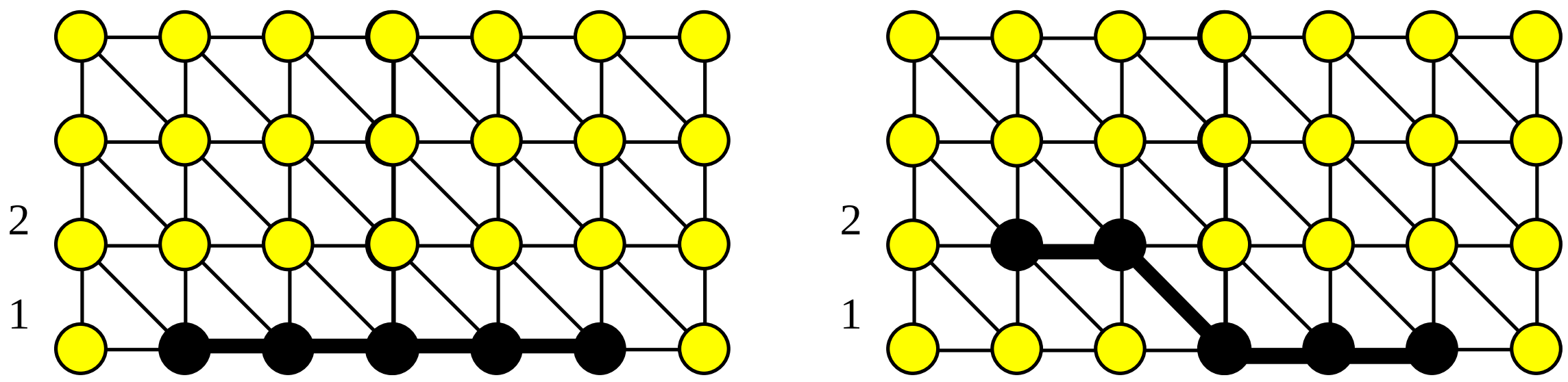

Рис. 217 Процесс перехода из параллельного мира (слоя) 1 в параллельный мир (слой) 2.

Это не представляется очень сложной задачей, особенно по сравнению с перестройками пространств. Решение такой задачи дело, в первую очередь, экспериментальных наук, хотя другие науки также могут внести свой неожиданный вклад проблему переходов между слоями пространства. Не исключено также, что среди того необъятного количества наблюдений и экспериментов, которые уже проделаны, имеются и такие, наиболее естественное объяснение которых лежит в принятии гипотезы о возможности перехода в параллельное пространство.

## 4. НЕКОТОРЫЕ ПРИМЕНЕНИЯ ТМП

В применении к компьютерным наукам ТМП дает фундаментальную математическую базу для принципов формирования и работаты с многомерными объектами при компьютерных вычислениях.

Хотя ТМП еще только формируется как научное направление, уже полученные результаты могут найти применение не только в компьютерных науках, но и в более традиционных областях науки и практики. В применении к физике, в частности, была проанализирована топологическая структура дискретного пространства-времени. Оказалось, что дискретное пространство-время является регулярным, то есть имеет избыток связей между двумя соседними по времени слоями трехмерного пространства[6].

Еще один интересный вывод был сделан относительно начальной стадии замкнутой модели расширяющейся вселенной. Выяснилось, что в момент своего рождения вселенная состояла из не менее чем восьми атомов пространства. [8]. А это может означять, что все многообразие окружающего нас мира, включая элементарные частицы, галактики и самого человека, порождено сочетаниями восьми элементов. Если это так, то такой вывод важен не только для физики, но и для нашего мировоззрения вообще.

Для теоретической физики ТМП дает богатый спектр моделей физического вакуума, особенно если учесть, что в качестве атомов пространства можно, в свою очередь, также брать молекулярные пространства.

Исследователи, работающие в пограничных областях науки, возможно, смогут использовать в своей работе такие процессы, как перестройка отдельных участков пространства с увеличением размерности и смещение в параллельный слой трехмерного пространства.

Хочется надеяться, что после того, как представители различных областей знаний ознакомятся с ТМП и сочтут полезными для себя хотя бы некоторые ее результаты, последуют новые и интересные открытия.

Список литературы.

# СПИСОК ЛИТЕРАТУРЫ

# АЛФАВИТНЫЙ УКАЗАТЕЛЬ

**—Ф—**

**—Ц—**

**—Ш—**

**—Э—**

—Я—

Научное издание

**Александр Вадимович Ивако**

# Введение в теорию молекулярных пространств

*Книга издается в авторской редакции*